\documentclass[11pt]{book}
\usepackage[a4paper,margin=2.6cm]{geometry}
\usepackage{amsmath,amssymb,amsthm,mathtools}
\usepackage{longtable}
\usepackage[colorlinks,linkcolor=blue,citecolor=blue]{hyperref}
\hypersetup{pdftitle={Digital Finite Group Theory},pdfauthor={Vladimir Garc\'{\i}a Morales},pdfsubject={Group theory (math.GR)}}

\theoremstyle{plain}
\newtheorem{theorem}{Theorem}[chapter]
\newtheorem{proposition}[theorem]{Proposition}
\newtheorem{lemma}[theorem]{Lemma}
\newtheorem{corollary}[theorem]{Corollary}
\theoremstyle{definition}
\newtheorem{definition}[theorem]{Definition}
\newtheorem{example}[theorem]{Example}
\theoremstyle{remark}
\newtheorem{remark}[theorem]{Remark}
\numberwithin{equation}{chapter}

\newcommand{\dg}[3]{d_{#1}\!\left(#2,#3\right)}   
\newcommand{\Bf}{\mathcal{B}}                                   
\newcommand{\md}[2]{D_{#1}\!\left(#2\right)}       
\newcommand{\Ap}[1]{\mathcal{A}_{#1}}                          
\newcommand{\Co}{\mathrm{Co}}                                  
\newcommand{\fl}[1]{\left\lfloor #1\right\rfloor}

\newcommand{\Z}{\mathbb{Z}}
\DeclareMathOperator{\Sym}{Sym}
\DeclareMathOperator{\Alt}{Alt}
\DeclareMathOperator{\Aut}{Aut}
\DeclareMathOperator{\PSL}{PSL}

\DeclareMathOperator{\Sp}{Sp}
\DeclareMathOperator{\GL}{GL}
\DeclareMathOperator{\SL}{SL}
\DeclareMathOperator{\sgn}{sgn}
\DeclareMathOperator{\lcm}{lcm}

\title{Digital Finite Group Theory}
\author{Vladimir Garc\'{\i}a Morales}
\date{October 2026}

\usepackage{tikz}
\providecommand{\Aut}{\operatorname{Aut}}
\providecommand{\Inn}{\operatorname{Inn}}
\providecommand{\Out}{\operatorname{Out}}

\usepackage{graphicx}
\makeatletter
\DeclareRobustCommand{\_}{\ifmmode\nfss@text{\textunderscore}\else\textunderscore\allowbreak\fi}
\makeatother
\makeatletter
\renewcommand*\l@section{\@dottedtocline{1}{1.5em}{3em}}
\renewcommand*\l@subsection{\@dottedtocline{2}{4.5em}{3.8em}}
\renewcommand*\l@subsubsection{\@dottedtocline{3}{8.3em}{4.6em}}
\makeatother
\begin{document}
\makeatletter
\begin{titlepage}
\null\vspace{1.5cm}
{\centering
{\Huge\bfseries \@title\par}
\vspace{1.2cm}
{\Large \@author\par}
\vspace{0.4cm}
{\large \@date\par}
\vspace{1.5cm}
{\bfseries Abstract\par}
}
\vspace{0.4cm}
\begin{quote}\small
Let an \emph{alphabet} be the set $\Ap{N}=\{0,1,\dots,N-1\}$ of the first $N$ non-negative integers, the
model of a finite set of size $N$, whose elements serve as labels. We describe every finite group $G$ as an alphabet
endowed with the group operation, that is, as a pair $(\Ap{N},*)$, where $N=|G|$ is the order of the group and $*$ is
the product of $G$ transported to $\Ap{N}$. This is achieved for all finite groups, and it is explained in detail for
all the infinite families of finite simple groups---the cyclic groups of prime order, the alternating groups and the
groups of Lie type---and for the $26$ sporadic groups.
A group of order $N$ becomes then a \emph{digital law}: a labelling $\lambda\colon G\to\Ap{N}$, its inverse, and the
transported product $m*n=\lambda\big(\lambda^{-1}(m)\lambda^{-1}(n)\big)$ written as a closed formula in the digits of
$m$ and $n$.
The book develops this point of view in three directions. The foundations treat codes of
permutations and of $r$-ary operators, and finite magmas. The general theory (construction of laws by diagrams, the
classical theorems, conjugacy classes, stabiliser chains and characters) is restated as identities between digits of
codes. The catalogue gives laws, in one mixed-radix format, for all groups of order at most $32$; for the finite quotients
$G/mL$ of the $17$ plane and $230$ space crystallographic groups, where the law is the point-group label followed by the
digits of a translation and the carry is an integer cocycle; for the symmetric and alternating groups; for the families of
finite simple groups of Lie type; and for all $26$ sporadic groups. For most sporadic groups the law is a chain law
along a permutation action. For the Baby Monster and the Monster, whose smallest actions have about $1.4\cdot10^{10}$ and
$10^{20}$ points, the coset digits are computed structurally from the axes of the $196\,884$-dimensional representation,
and the laws are evaluated, inverted and checked with the \texttt{mmgroup} package. A chapter of worked products, from
$\mathrm{PSL}_2(7)$ to the Monster, multiplies and inverts explicit labels, and every law, table and worked example is
reproduced by the scripts of the accompanying software archive (\href{https://doi.org/10.5281/zenodo.23881253}{doi:10.5281/zenodo.23881253}). The chapters end with exercises, whose solutions are
collected in the last chapter. The structure of these groups is classical; what is new is that they become explicit,
verifiable integer data.
\end{quote}
\vfill
\end{titlepage}
\makeatother
\thispagestyle{empty}\null\newpage
\tableofcontents

\chapter*{Preface}
\phantomsection\addcontentsline{toc}{chapter}{Preface}
\addtocontents{toc}{\protect\setcounter{tocdepth}{0}}

\section*{The idea of the book}
By an \emph{alphabet} we mean the set $\Ap{N}=\{0,1,\dots,N-1\}$ of the first $N$
non-negative integers; it is the model of a finite set of size $N$. Its elements are
\emph{labels}: a finite set $X$ with $N$ elements is identified with $\Ap{N}$ through a
bijection $\lambda\colon X\to\Ap{N}$, its \emph{labelling}. A label is an integer in its
own right and should not be confused with a digit. Each label $m$ has its own digits, in
a base or in a mixed radix adapted to the structure at hand, and it is in these digits
that operations on the alphabet are written.

The idea of the book is to describe every finite group $G$ as an alphabet
endowed with the group operation, that is, as a pair $(\Ap{N},*)$, where
$N=|G|$ is the order of the group and $*$ is the product of $G$ transported
to $\Ap{N}$. This is achieved for all finite groups, and it is explained in
detail for all the infinite families of finite simple groups---the cyclic
groups of prime order, the alternating groups and the groups of Lie
type---and for the 26 sporadic groups.

A finite operation on an
alphabet is read off the digits of a single integer, its \emph{code}. The basic tool is the digit function
\[
\dg{N}{k}{x}=\fl{\frac{x}{N^{k}}}-N\fl{\frac{x}{N^{k+1}}},
\]
which returns the $k$-th digit of $x$ in base $N$. An $r$-ary operation on $\Ap N$ is the integer whose $N^r$ base-$N$ digits are its
values, so a magma, a permutation or a whole Cayley table is one number. A group of order $N$ becomes a \emph{digital law}: a labelling
$\lambda\colon G\to\Ap N$, its inverse, and the transported product $m*n=\lambda\big(\lambda^{-1}(m)\lambda^{-1}(n)\big)$ written as a
closed formula in the digits of $m$ and $n$. For a cyclic group this formula is addition modulo $N$; for a metacyclic group it is the
two-digit rule $m*n=\dg{q}{0}{m+r^{m_1}n}+q\,\dg{s}{0}{m_1+n_1}$; for a finite simple group it packs the coordinates of a normal form
or of a stabiliser chain into mixed-radix digits. Every finite group has such a law, since any bijection transports the product. The
point of the book is that the bijection can be chosen so that the law is explicit, short and computable, and that the structural
properties of the group---associativity, inverses, conjugacy, normal subgroups, characters---become identities between digits of codes.

The book has grown in three directions, and it is worth saying plainly what each contributes. The first is foundational: the digit
function and the $\Bf$-function, codes of permutations and of $r$-ary operators, and the observation that every finite structure is a
single integer. The second is a re-reading of the general theory of finite groups in digital terms: the construction of laws by
diagrams, the classical theorems, conjugacy classes, stabiliser chains and characters. Here the abstract statements are kept, and the
digital form supplies explicit and checkable instances of them. The third is a catalogue: laws for all groups of order at most $32$, for the finite quotients of the plane and space
crystallographic groups, for the symmetric and alternating groups, for the families of finite simple groups of Lie type and
for the sporadic groups, all in one format, so that they can be compared, transmitted and verified as integer arithmetic.
The catalogue is accompanied by worked products in the largest groups and by a software archive that recomputes every
law, table and example of the book. Nothing in the book changes what is known about
the structure of these groups; the classical background is that of the standard texts~\cite{Hall,Roman,Isaacs,Serre,Wilson}. What
changes is that the groups become concrete data.

\section*{Outline of the chapters}
\paragraph*{Foundations (Chapters~\ref{ch:digit} and~\ref{ch:magmas}).}
Chapter~\ref{ch:digit} introduces the digit function and the digital expansion of integers (Section~\ref{sec:digits}), and the
$\Bf$-function, which tests digits and expresses comparisons and carries of positional arithmetic in closed form
(Section~\ref{sec:bfunction}). Chapter~\ref{ch:magmas} begins with permutations as codes---composition, inversion, cycles, conjugacy and
sign (Section~\ref{sec:permutations})---and then encodes arbitrary finite magmas and $r$-ary operators, their homomorphisms and
relabellings, and direct products. Everything that follows uses these two chapters.

\paragraph*{General theory in digital form (Chapters~\ref{ch:diagrammatic}--\ref{ch:characters}).}
Chapter~\ref{ch:diagrammatic} is a method for producing group laws: the digits of a label are arranged in layers, and a diagram of
arrows prescribes how carries couple them. Associativity is built in by construction, and a gallery of recipes recovers many familiar
groups, including simple ones. Chapter~\ref{ch:digitaltheorems} restates the main theorems---Cayley, Lagrange, group actions, conjugacy classes and the counting of colourings (P\'olya--Redfield and de Bruijn), the
isomorphism theorems, the Sylow theorems and the series theorems of Zassenhaus and Schreier---for digital laws, with proofs that run in
the digits. Chapter~\ref{ch:ss} is algorithmic. It presents the Schreier--Sims
algorithm~\cite{sims1970,butler1991,seress2003,holt2005}, whose base and strong generating set is precisely a mixed-radix labelling,
and then computes conjugacy classes (their basic theory is in Section~\ref{ssec:dtclasses}), power maps, rationality and fusion,
ending with a complete example for $M_{12}$.
Chapter~\ref{ch:characters} develops ordinary character theory---Maschke, Schur and Wedderburn, the orthogonality relations,
integrality, induction and restriction---and turns it into a roadmap for computing character tables from class
data~\cite{Isaacs,Serre,FultonHarris,Dixon1967}.

\paragraph*{Families of groups (Chapters~\ref{ch:abelian}--\ref{ch:small}).}
Chapter~\ref{ch:abelian} gives the laws of cyclic groups, additive and multiplicative, and of finite Abelian groups as direct sums, with
explicit isomorphisms for coprime factors. Chapter~\ref{ch:nonabelian} treats semidirect products, dicyclic, dihedral and metacyclic
groups, wreath products and transitive actions; its Section~\ref{sec:metacyclic-chars} computes the characters of metacyclic, dihedral
and dicyclic groups directly from the two digits of their laws. Chapter~\ref{ch:small} classifies the groups of order at most $32$ and
gives a law for each of them~\cite{Holder,BEO}.

\paragraph*{Crystallographic groups (Chapter~\ref{ch:crystal}).}
Chapter~\ref{ch:crystal} gives digital laws for the finite quotients $G/mL$ of the $17$ plane groups and of the $230$ space
groups. An element $\{R_i\mid t_i+24s\}$ receives the label $i\,m^3+s_1m^2+s_2m+s_3$; the product adds the rotated
translation digits, and the carry $c(i,j)$ is an integer cocycle determined by the conventional operations. The quotients
modulo $2T$, $6T$ and $2L$ are identified in the small groups library and separated by their invariants.

\paragraph*{Finite simple groups (Chapters~\ref{ch:simple}--\ref{ch:sporadic}).}
Chapter~\ref{ch:simple} treats $\Sym_n$ and $\Alt_n$: generators and presentations, conjugacy classes, simplicity, maximal subgroups,
automorphisms including the outer automorphism of $\Sym_6$, the characters of $\Sym_n$ and
$\Alt_n$~\cite{James,JamesKerber,Sagan,Macdonald,FRT}, and the Lehmer and rotation laws. Chapter~\ref{ch:lie} builds finite fields in
digital form and then the groups of Lie type, from $\PSL(2,N)$ through the classical families, the Suzuki and Ree groups, $G_2(q)$ and
${}^3D_4(q)$ to $F_4(q)$, ${}^2F_4(q)$ and the groups of type $E$. Their product rules are transcriptions of the Bruhat normal form and
the Chevalley commutator relations~\cite{Dickson,steinberg1967,Carter,Taylor,SpringerVeldkamp,McCrimmon,KLM}. Chapter~\ref{ch:sporadic}
gives laws for the sporadic groups, from the Mathieu groups through the Conway, Fischer and Janko groups to the Monster and the Baby
Monster, relative to standard generators~\cite{ATLAS,AtlasRep,ConwaySloane,GriessLyons}. It ends with canonical labels and with the
recomputation of a character table from a printed law.

\paragraph*{Worked products and software (Chapters~\ref{ch:products} and~\ref{ch:software}).}
Chapter~\ref{ch:products} multiplies explicit labels in $\mathrm{PSL}_2(7)$, $\mathrm{PSL}_3(2)$, $\mathrm{Sz}(8)$,
${}^2G_2(3)$, $M_{24}$, $\Co_0$, $J_4$, $\mathrm{Th}$ and the Monster, and computes the inverses and orders of the factors
and products within the digital laws themselves, with an independent check by decoding. Chapter~\ref{ch:software} describes the
published archives, their installation and layout, the reproduction workflow, and each program and data file, including
those for the crystallographic groups.

Chapter~\ref{ch:solutions} contains solutions to the exercises of all chapters.

\section*{Dependencies and reading paths}
\begin{figure}[htbp]
\centering
\begin{tikzpicture}[x=1.2cm,y=1.05cm,ch/.style={draw,rounded corners,font=\footnotesize,inner sep=3pt,align=center}]
\node[ch] (c1) at (3,0) {1\ Digit function};
\node[ch] (c2) at (3,-1.3) {2\ Permutations, magmas};
\node[ch] (c3) at (0.5,-2.6) {3\ Diagrams};
\node[ch] (c8) at (5.5,-2.6) {7\ Cyclic, Abelian};
\node[ch] (c4) at (0.5,-3.9) {4\ Main theorems};
\node[ch] (c9) at (5.5,-3.9) {8\ Nonabelian};
\node[ch] (c5) at (0.5,-5.2) {5\ Schreier--Sims,\\conjugacy classes};
\node[ch] (c10) at (6.6,-5.2) {9\ Small groups};
\node[ch] (c7) at (0.5,-6.5) {6\ Characters};
\node[ch] (c11) at (1.7,-7.8) {11\ $\Sym_n$, $\Alt_n$};
\node[ch] (c12) at (4.9,-7.8) {12\ Lie type};
\node[ch] (c13) at (3.3,-9.1) {13\ Sporadic};
\node[ch] (ccr) at (6.6,-6.5) {10\ Crystallographic};
\node[ch] (c14) at (3.3,-10.4) {14\ Worked products};
\draw[->] (c1)--(c2);
\draw[->] (c2)--(c3);
\draw[->] (c2)--(c8);
\draw[->] (c3)--(c4);
\draw[->] (c4)--(c5);
\draw[->] (c5)--(c7);
\draw[->] (c8)--(c9);
\draw[->] (c9)--(c10);
\draw[->,dashed] (c7)--(c9);
\draw[->,dashed] (c4)--(c9);
\draw[->] (c7)--(c11);
\draw[->] (c11)--(c12);
\draw[->] (c5.south west) .. controls (-1.0,-7.6) and (0.6,-9.1) .. (c13.west);
\draw[->] (c9)--(c13);
\draw[->] (c12)--(c13);
\draw[->] (c10)--(ccr);
\draw[->] (c9.east) .. controls (8.4,-4.6) and (8.4,-6.0) .. (ccr.east);
\draw[->] (c13)--(c14);
\draw[->] (c12) to[bend left=30] (c14);
\end{tikzpicture}
\caption{Main dependencies between the chapters. An arrow $A\to B$ means that Chapter~$B$ uses results of Chapter~$A$;
Chapter~\ref{ch:digit} and Section~\ref{sec:permutations} are used throughout, and only their first arrows are drawn. The dashed
arrows into Chapter~\ref{ch:nonabelian} mark Sections~\ref{sec:metacyclic-chars} and~\ref{sec:simpleclass}. Chapter~\ref{ch:software}, which
documents the programs of all chapters, is not drawn.}\label{fig:dependencies}
\end{figure}
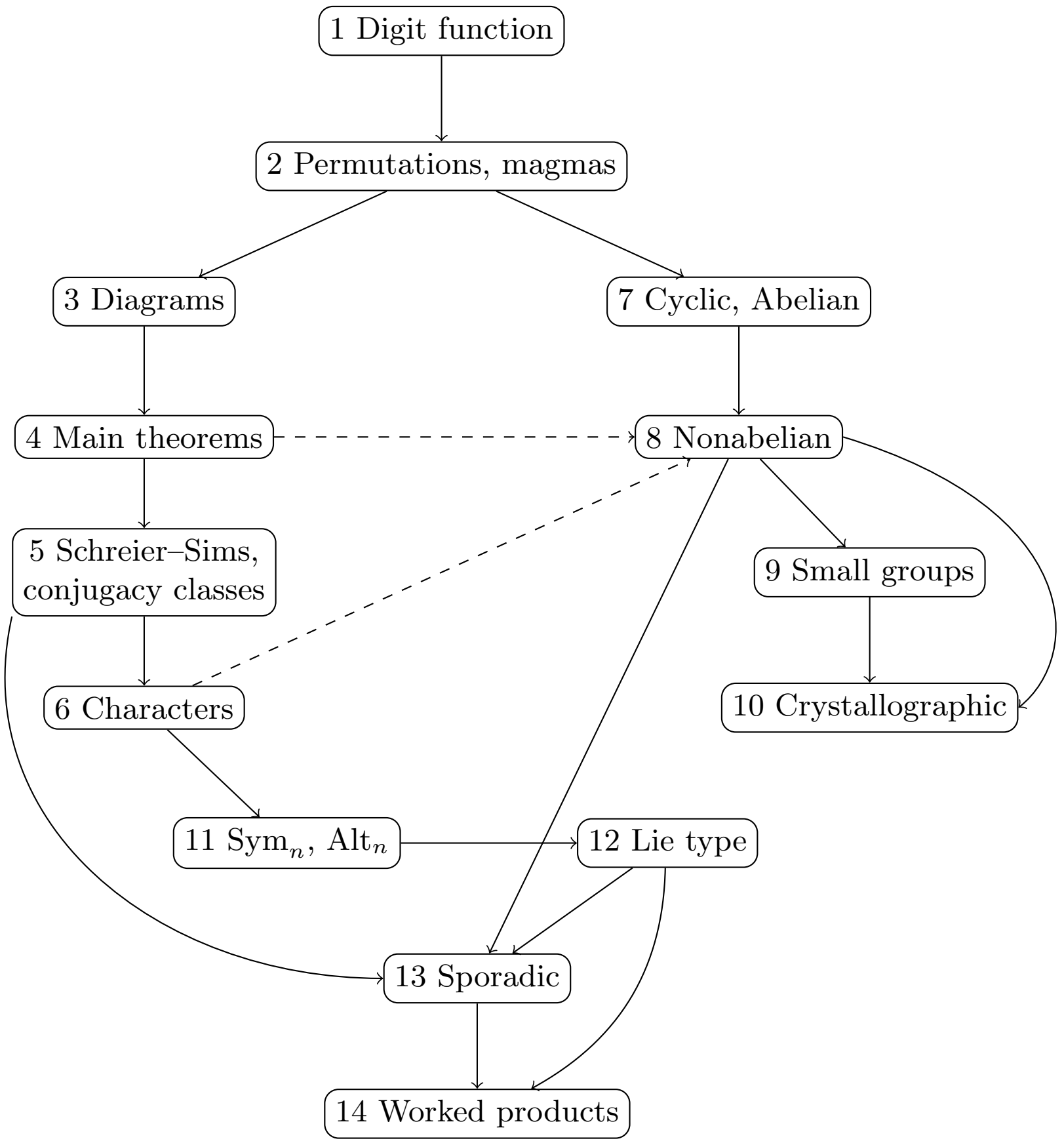

Figure~\ref{fig:dependencies} shows the main dependencies. Chapter~\ref{ch:digit} is used everywhere, and the permutation codes of
Section~\ref{sec:permutations} reappear wherever groups act. Chapter~\ref{ch:digitaltheorems} builds on the diagrammatic construction
of Chapter~\ref{ch:diagrammatic}. Chapter~\ref{ch:ss} combines the Schreier--Sims algorithm with conjugacy classes, and
Chapter~\ref{ch:characters} uses both parts. The families of Chapters~\ref{ch:abelian}--\ref{ch:small} need essentially only
Chapters~\ref{ch:digit} and~\ref{ch:magmas}; the exceptions are Section~\ref{sec:metacyclic-chars}, which uses
Chapter~\ref{ch:characters}, and Section~\ref{sec:simpleclass}, which uses Chapters~\ref{ch:digitaltheorems} and~\ref{ch:characters}. These chapters are placed after the general theory so that Chapters~\ref{ch:diagrammatic}--\ref{ch:characters}
form one block, but a reader who prefers examples first can read them directly after Chapter~\ref{ch:magmas}. Chapters~\ref{ch:simple}--\ref{ch:sporadic}
draw on all of this: Chapter~\ref{ch:simple} on Section~\ref{sec:permutations} and on Chapters~\ref{ch:cc} and~\ref{ch:characters},
Chapter~\ref{ch:lie} on Chapter~\ref{ch:simple}, and Chapter~\ref{ch:sporadic} on the stabiliser chains of Chapter~\ref{ch:ss}, the transitive laws of Theorem~\ref{thm:transitive} and the groups of Lie type.
Chapter~\ref{ch:crystal} uses the semidirect products of Chapter~\ref{ch:nonabelian} and the identification of small groups in
Chapter~\ref{ch:small}. Chapter~\ref{ch:products} draws on Chapters~\ref{ch:lie} and~\ref{ch:sporadic}, and
Chapter~\ref{ch:software} can be consulted at any point alongside the chapter whose computations it reproduces.

Three routes through the book are natural. A first course in finite groups in digital form reads Chapters~\ref{ch:digit}
and~\ref{ch:magmas}, then Chapters~\ref{ch:abelian}--\ref{ch:small}, and then Chapter~\ref{ch:digitaltheorems}. A reader interested
in computation reads Chapters~\ref{ch:digit} and~\ref{ch:magmas}, then Chapters~\ref{ch:ss} and~\ref{ch:characters}. A
reader interested in the simple groups can go from Chapters~\ref{ch:digit} and~\ref{ch:magmas} directly to
Chapters~\ref{ch:simple}--\ref{ch:sporadic}, consulting Chapters~\ref{ch:cc}--\ref{ch:characters} and Chapter~\ref{ch:nonabelian}
when they are cited.

\section*{Verification and reproducibility}
Every law in the book is a finite object, and most of them are small enough to be checked exhaustively. Generated tables---for
instance the character data of Section~\ref{sec:metacyclic-chars} and the tables of Chapter~\ref{ch:simple}---are produced by scripts
kept with the source of the book. They rebuild each group from its printed law, verify the stated identities in exact arithmetic and
compare the results with GAP~\cite{GAP}; the single command \texttt{python3 scripts/sn\_verify.py} regenerates and checks all of them.
For the large groups the laws are tied to standard generators~\cite{ATLAS,AtlasRep} and come with fingerprints, so that they can be reproduced independently of this book. Chapter~\ref{ch:software} lists the published archives and
describes every program and data file; the same scripts reproduce the crystallographic laws and identifications of
Chapter~\ref{ch:crystal} and the worked products and inverses of Chapter~\ref{ch:products}.

\section*{Conventions}
Throughout, $\Ap N=\{0,1,\dots,N-1\}$, $\dg{N}{k}{x}$ is the $k$-th base-$N$ digit of $x$, and $\md{h}{m}$ is the $h$-th digit of $m$
in a mixed radix fixed by the context. Permutations are composed from right to left, as in~\eqref{eq:comp}. Group laws are written
$m*n$ on labels; when a group is given by generators, $x^ay^\alpha$ denotes the element whose label has digits $a$ and~$\alpha$.
The preface closes with a survey of related work and a list of what we believe is new.

\section*{Related work and what is new}\label{sec:related}
\addcontentsline{toc}{section}{Related work and what is new}
\begingroup\setlength{\emergencystretch}{2em}

A digital law is a bijection $\lambda\colon G\to\Ap{|G|}$ together with the transported product
$m*n=\lambda\big(\lambda^{-1}(m)\lambda^{-1}(n)\big)$. Bijections of this kind are not new; what this book adds is a uniform,
explicit and reproducible catalogue of them. We list the known ingredients first.

\paragraph{Stabiliser chains.}
Sims introduced bases and strong generating sets~\cite{sims1970}: for a base $(b_1,\dots,b_k)$ every element of a permutation
group factors uniquely as $s_{1,t_1}\cdots s_{k,t_k}$, with $s_{l,t}$ running over a transversal of the $l$-th basic orbit. Reading
the indices $t_l$ in mixed radix is exactly the chain law of Chapter~\ref{ch:sporadic} and the iterated transitive law of
Theorem~\ref{thm:transitive}. This factorisation is the basis of permutation group algorithms~\cite{butler1991,seress2003,holt2005};
GAP enumerates the elements of a permutation group through it~\cite{GAP}, Knuth uses ``Sims tables'' to generate all elements of
a group~\cite{knuth2011}, and Jerrum studied compact representations built on it~\cite{jerrum1986}. The bijection underlying a chain
law is therefore classical; the contribution here is to treat the resulting integer law as the object of study and to publish it.

\paragraph{Ranking and unranking.}
The permutation codes of Section~\ref{sec:permutations} are relatives of the Lehmer code and the factorial number system~\cite{lehmer1960}. Ranking and
unranking combinatorial objects, including permutations in lexicographic order and in linear time, is a developed subject
\cite{kreher1999,myrvold2001,kagey2022}. Its emphasis is on efficient bijections for whole classes of objects, not on transporting a group
product.

\paragraph{Normal forms in groups of Lie type.}
The laws of Chapter~\ref{ch:lie} pack the coordinates of the Bruhat normal form $G=\bigsqcup_{w}U\,T\,\dot w\,U_w^-$ into one integer.
The normal form and the Chevalley commutator relations that drive the product are classical~\cite{steinberg1967,Carter}, and
algorithms that multiply in these coordinates are implemented in Magma~\cite{cmt2004}. Our product rules for ${}^3D_4(q)$, $F_4(q)$,
${}^2F_4(q)$ and the $E$ series are explicit transcriptions of this theory, checked against matrix models.

\paragraph{Canonical representatives.}
The lexicographically least element of a coset, and more generally the smallest image under a group, is a standard tool in backtrack
search~\cite{butler1991,seress2003,linton2004}. The canonical lexicographic law uses this greedy step; what is new is to use it to fix
a labelling that depends only on the permutation group.

\paragraph{The Monster and the Baby Monster.}
The smallest permutation degrees of $\mathbb B$ and $\mathbb M$, about $1.4\cdot10^{10}$ and $10^{20}$, exclude the orbit enumeration
that works for the other sporadic groups, so the coset digits must be computed from the structure. Our starting point is the
Python package \texttt{mmgroup} of Seysen, which multiplies elements of $\mathbb M$ and reduces them to a normal form, together with
the $196\,884$-dimensional representation, in which the $2A$ involutions correspond to axes. The package is a model of $\mathbb M$, not
a labelling of it: it does not give a canonical coset representative, and we construct one from the orbits of $K=C_{G_{x0}}(\beta^+)$
on the axes.

\paragraph{Crystallographic groups.}
The $17$ plane and $230$ space groups and their symmetry operations are classical and tabulated in the International Tables for
Crystallography~\cite{ITA}. A space group is an extension of its lattice by its point group, described by a cohomology class with values in the
translations, and its finite quotients by sublattices are a standard tool for studying subgroups of finite index. The laws of
Chapter~\ref{ch:crystal} make this extension explicit as integer arithmetic: the label is the point-group index followed by the
translation digits modulo $m$, and the carry is the integer cocycle $c(i,j)$ of the conventional operations. The quotients modulo
$2T$, $6T$ and $2L$ are identified with GAP~\cite{GAP}.

\paragraph{Groups as data.}
The \textsc{Atlas}~\cite{ATLAS} and the online \textsc{Atlas} of group representations~\cite{AtlasRep} fix standard generators,
which is what makes the generator fingerprints of this book reproducible. Integer-labelled group operations have recently become
benchmarks in machine learning~\cite{power2022,chughtai2023}, where they are mostly limited to cyclic and small symmetric groups.

\paragraph{What is new.}
To our knowledge the following are not available elsewhere in one place:
\begin{itemize}
\item a digital law, with explicit $\lambda$, $\lambda^{-1}$ and $*$ on $\Ap{|G|}$, for every family of finite simple groups treated
here, in a single format;
\item explicit coordinate product rules for the exceptional families, verified over every field tested;
\item canonical laws that depend only on the group and its standard generators, with fingerprints that allow independent
reproduction;
\item complete laws for all groups of order at most $32$;
\item digital laws for the finite quotients $G/mL$ of all $17$ plane and $230$ space crystallographic groups, with integer cocycles
as carries, and identifications of the quotients modulo $2T$, $6T$ and $2L$ in the small groups library (Chapter~\ref{ch:crystal});
\item the observation, carried out for $M_{11}$, that the printed laws suffice to recompute classical invariants such as the character
table;
\item a method that produces associative laws layer by layer from diagrams of carries (Chapter~\ref{ch:diagrammatic}), and
digital versions of the main theorems of finite group theory (Chapter~\ref{ch:digitaltheorems});
\item character theory carried out on codes, with closed formulas in the digits of the labels for the characters of metacyclic,
dihedral and dicyclic groups (Section~\ref{sec:metacyclic-chars}) and with generated character data for $\Sym_n$ and $\Alt_n$
(Chapter~\ref{ch:simple});
\item a digital law for $\mathbb M$ whose coset digit indexes the $2A$ involution (equivalently its axis) and whose fibre digit is
a law of $2.\mathbb B$ (Proposition~\ref{prop:monster}), and a digital law for $\mathbb B$ obtained from it (Theorem~\ref{thm:blaw}), in both cases
with canonical coset representatives built from orbits of axes and checked against \texttt{mmgroup} on random elements;
\item worked products, inverses and orders of explicit labels in groups from $\mathrm{PSL}_2(7)$ to the Monster, computed within the
digital laws and checked independently by decoding (Chapter~\ref{ch:products});
\item generated tables whose scripts rebuild each group from its printed law, verify the stated identities exactly and compare them
with GAP~\cite{GAP}, published with all other programs and data in the versioned archives described in
Chapter~\ref{ch:software}.
\end{itemize}
Nothing in this book changes what is known about the structure of these groups. The aim is to make them explicit, computable and
comparable as integer arithmetic.
\par\endgroup

\addtocontents{toc}{\protect\setcounter{tocdepth}{2}}
\chapter{The digit function}\label{ch:digit}

\section{Digits of integers and the digital expansion}\label{sec:digits}

\begin{definition}
Let $N\ge 2$ be an integer. For $x\in\mathbb{R}$ and $k\in\Z$ the \emph{digit function} is
\begin{equation}\label{eq:digit}
\dg{N}{k}{x}=\fl{\frac{x}{N^{k}}}-N\fl{\frac{x}{N^{k+1}}}.
\end{equation}
The set $\Ap{N}=\{0,1,\dots,N-1\}$ is the \emph{digital alphabet} of size $N$.
\end{definition}

\begin{lemma}[Basic properties]\label{lem:basic}
Let $x\in\mathbb{R}$ and $k\in\Z$.
\begin{enumerate}
\item $\dg{N}{k}{x}\in\Ap{N}$ and $\dg{N}{k}{x}=\dg{N}{0}{x/N^{k}}$.
\item If $x\in\Z$, then $\dg{N}{0}{x}=x \bmod N$.
\item If $x\in\Z$ and $0\le x<N^{M}$, then $\dg{N}{k}{x}=0$ for $k\ge M$ and
$x=\sum_{k=0}^{M-1}\dg{N}{k}{x}\,N^{k}$.
\end{enumerate}
\end{lemma}
\begin{proof}
Put $y=x/N^{k}$. Since $\fl{\fl{y}/N}=\fl{y/N}$, Euclidean division gives $\fl{y}=N\fl{y/N}+r$ with $r\in\Ap{N}$, and
$r=\dg{N}{k}{x}$; this proves (1) and (2). For (3) the sum telescopes:
\[
\sum_{k=0}^{M-1}\big(\fl{x/N^{k}}-N\fl{x/N^{k+1}}\big)N^{k}=x-N^{M}\fl{x/N^{M}}=x.
\]
\end{proof}

\begin{lemma}[Modular arithmetic]\label{lem:mod}
For $x,y\in\Z$,
$\dg{N}{0}{x\pm y}=\dg{N}{0}{\dg{N}{0}{x}\pm\dg{N}{0}{y}}$ and
$\dg{N}{0}{xy}=\dg{N}{0}{\dg{N}{0}{x}\,\dg{N}{0}{y}}$.
\end{lemma}
\begin{proof}Immediate from Lemma~\ref{lem:basic}(2).\end{proof}

\begin{lemma}[Block digits]\label{lem:block}
Let $n\ge1$, $x\in\Z_{\ge0}$, $a\in\Z_{\ge 0}$ and $b\in\Ap{n}$. Then
\begin{equation}\label{eq:block}
\dg{N}{an+b}{x}=\dg{N}{b}{\dg{N^{n}}{a}{x}}.
\end{equation}
\end{lemma}
\begin{proof}
\[
\dg{N^n}{a}{x}=\fl{x/N^{an}}\bmod N^{n},\qquad
\dg{N}{an+b}{x}=\fl{x/N^{an+b}}\bmod N.
\]
Because $b<n$, the second expression is the $b$-th base-$N$ digit of the first.
\end{proof}

\begin{definition}[Mixed radix]
Let $p_0,\dots,p_{\ell-1}\ge 2$ and $\Pi_h=\prod_{j<h}p_j$ (so $\Pi_0=1$). The \emph{$h$-th mixed-radix digit} of $x$ is
\begin{equation}\label{eq:mixed}
\md{h}{x}=\dg{p_h}{0}{\frac{x}{\Pi_h}}.
\end{equation}
\end{definition}

\begin{lemma}\label{lem:mixed}
The map $x\mapsto(\md{0}{x},\dots,\md{\ell-1}{x})$ is a bijection
$\Ap{\Pi_\ell}\to\Ap{p_0}\times\dots\times\Ap{p_{\ell-1}}$, with inverse $(x_h)\mapsto\sum_h x_hP_h$.
\end{lemma}
\begin{proof}
Since $\fl{\fl{x/\Pi_h}/p_h}=\fl{x/\Pi_{h+1}}$, we have
\[
\md{h}{x}=\fl{x/\Pi_h}-p_h\fl{x/\Pi_{h+1}}.
\]
Moreover, $\sum_h\md{h}{x}\Pi_h$ telescopes to $x-\Pi_\ell\fl{x/\Pi_\ell}=x$. So the map is injective; both sets have $\Pi_\ell$ elements.
\end{proof}

\begin{theorem}[Coding principle]\label{thm:coding}
Let $r\ge1$. Every map $f:\Ap{N}^{r}\to\Ap{N}$ is of the form
\begin{equation}\label{eq:code}
f(x_0,\dots,x_{r-1})=\dg{N}{\textstyle\sum_{k=0}^{r-1}x_kN^{k}}{R}
\end{equation}
for a unique integer $R\in[0,N^{N^{r}})$, the \emph{code} of $f$, namely
$R=\sum_{x\in\Ap{N}^r}f(x)\,N^{\sum_k x_kN^k}$. Hence there are exactly $N^{N^{r}}$ $r$-ary operators on $\Ap{N}$.
\end{theorem}
\begin{proof}
By Lemma~\ref{lem:basic}(3), $x\mapsto\sum_kx_kN^k$ is a bijection $\Ap{N}^r\to\Ap{N^r}$ and
$R\mapsto(\dg{N}{n}{R})_{n<N^r}$ is a bijection $\Ap{N^{N^r}}\to\Ap{N}^{N^r}$.
\end{proof}

\section{Computing digits in practice}\label{sec:digits:practice}
The digit function~\eqref{eq:digit} is a closed formula: it gives the digit of $x$ in position $k$ directly, for real $x$ and
every integer $k$. This section compares it with the way digits are extracted in computation, and with Lehmer's code of a
permutation~\cite{lehmer1960}, which is the mixed-radix expansion used for $\Sym_n$ in Definition~\ref{def:sn:lehmer}.

\subsection{Sequential extraction}
For a radix $p_h=1$ the right-hand side of~\eqref{eq:digit} still makes sense and vanishes identically; we then read~\eqref{eq:mixed}
as $\md{h}{x}=0$. This is the case of the last radix in Definition~\ref{def:sn:lehmer}.

\begin{proposition}[Sequential digits]\label{prop:seqdigits}
Let $p_0,\dots,p_{\ell-1}\ge1$ be integers, $\Pi_h=\prod_{j<h}p_j$ and $x\in\Z_{\ge0}$. Put $x_0=x$ and $x_{h+1}=\fl{x_h/p_h}$.
Then $x_h=\fl{x/\Pi_h}$ and
\begin{equation}\label{eq:seqdigit}
\md{h}{x}=x_h-p_h\,x_{h+1}=x_h\bmod p_h\qquad(0\le h<\ell).
\end{equation}
In particular $\dg{N}{k}{x}=x_k-N\,x_{k+1}$ for the constant radix $p_h=N$.
\end{proposition}
\begin{proof}
For an integer $p\ge1$ and real $y$ one has $\fl{\fl{y}/p}=\fl{y/p}$. By induction,
$x_{h+1}=\fl{\fl{x/\Pi_h}/p_h}=\fl{x/\Pi_{h+1}}$. By~\eqref{eq:mixed} and~\eqref{eq:digit},
$\md{h}{x}=\fl{x/\Pi_h}-p_h\fl{x/\Pi_{h+1}}=x_h-p_hx_{h+1}$, which is the remainder of $x_h$ modulo $p_h$.
\end{proof}

Proposition~\ref{prop:seqdigits} is the familiar rule ``divide and keep the remainder''. It is formula~\eqref{eq:digit}
evaluated with the quotient carried from one position to the next, instead of being recomputed from $x$ each time.

\begin{remark}[Cost]\label{rem:digitcost}
Exact integers are stored in base $2^{w}$ for a word size $w$, typically $w=30$ or $w=64$. Let $x$ occupy $L$ words and let
every radix fit in one word. A division by a one-word integer costs $O(L)$ operations, so~\eqref{eq:seqdigit} gives all $\ell$
digits in $O(\ell L)$ operations. Evaluating~\eqref{eq:digit} or~\eqref{eq:mixed} separately in every position instead divides
$x$ by the multiword integers $\Pi_h$ and $\Pi_{h+1}$. With schoolbook division, dividing $L$ words by $s$ words costs about
$s(L-s)$ operations; summed over the positions this is of order $L^{3}$ when $\ell$ is proportional to $L$. The closed formula
is therefore the right tool for one digit in a prescribed position and for proofs, and the recurrence is the right tool for a
full expansion.
\end{remark}

\begin{remark}[Divide and conquer]\label{rem:digitdc}
Lemma~\ref{lem:block} gives a faster method for very long expansions. If $0\le x<N^{2n}$, write $x=r+N^{n}q$ with
$r=\dg{N^n}{0}{x}$ and $q=\dg{N^n}{1}{x}$. Then $\dg{N}{b}{x}=\dg{N}{b}{r}$ and $\dg{N}{n+b}{x}=\dg{N}{b}{q}$ for
$b\in\Ap n$. Applied recursively, this halves the problem at each step; with fast multiplication and division the whole
expansion costs only a logarithmic factor more than one product of $L$-word integers. The same splitting works in mixed radix,
with $x=(x\bmod\Pi_h)+\Pi_h\fl{x/\Pi_h}$. Arbitrary-precision libraries convert long binary integers to decimal in this way.
\end{remark}

\begin{remark}[Exact arithmetic]\label{rem:digitexact}
Formula~\eqref{eq:digit} is stated for real $x$, but it must be evaluated with exact floors. In IEEE double precision
every integer is exact only up to $2^{53}$, and $18!<2^{53}<19!$. Hence the Lehmer labels of $\Sym_n$
(Definition~\ref{def:sn:lehmer}) can be handled in double precision only for $n\le18$, and with unsigned $64$-bit integers only
for $n\le20$, since $20!<2^{64}<21!$. Beyond this one needs arbitrary-precision integers, for which~\eqref{eq:seqdigit} is
the native operation. Converting such an integer to decimal for printing is itself the expansion~\eqref{eq:digit} with $N=10$;
for integers with many thousands of decimal digits it is often the most expensive step of a computation.
\end{remark}

\subsection{Lehmer's code and the Lehmer label}
For $R\in\Sym_n$ let $a_i(R)=\#\{j>i:R(j)<R(i)\}\in[0,n-i)$, as in Definition~\ref{def:sn:lehmer}. The sequence
$(a_0(R),\dots,a_{n-1}(R))$ is \emph{Lehmer's code} of $R$. Classically it is read in the factorial number system,
\[
L(R)=\sum_{i=0}^{n-1}a_i(R)\,(n-1-i)!,
\]
whereas this book uses the Lehmer label
\[
\lambda(R)=\sum_{i=0}^{n-1}a_i(R)\,\Pi_i,\qquad \Pi_i=\prod_{j<i}(n-j)=\frac{n!}{(n-i)!}.
\]

\begin{proposition}[Two readings of one code]\label{prop:lehmerreadings}
Both $\lambda$ and $L$ are bijections $\Sym_n\to[0,n!)$, and $\lambda(\mathrm{id})=L(\mathrm{id})=0$.
\begin{enumerate}
\item The mixed-radix digits of $\lambda(R)$ for the radices $(n,n-1,\dots,1)$ are $a_0(R),\dots,a_{n-1}(R)$; thus $a_0(R)$ is
the least significant digit.
\item The mixed-radix digits of $L(R)$ for the radices $(1,2,\dots,n)$ are $a_{n-1}(R),\dots,a_0(R)$; thus $a_0(R)$ is the most
significant digit. Moreover $L(R)$ is the rank of $R$ when $\Sym_n$ is listed in lexicographic order of one-line words.
\end{enumerate}
Hence $\lambda$ and $L$ read the same digit string in opposite directions.
\end{proposition}
\begin{proof}
(1) is Theorem~\ref{thm:sn:lehmer}. For (2), the radices $q_k=k+1$ have place values $\prod_{j<k}(j+1)=k!$, and the digit in
position $k$ is $a_{n-1-k}(R)\in[0,k+1)$, so Lemma~\ref{lem:mixed} gives the bijection. Suppose that $R$ precedes $S$
lexicographically, and let $i$ be the first position with $R(i)\ne S(i)$, so that $R(i)<S(i)$. The values in the positions
before $i$ agree, hence $a_j(R)=a_j(S)$ for $j<i$; and $a_i(R)<a_i(S)$, because $a_i$ is the rank of the value in position $i$
among the values not yet used. Since
$\sum_{j>i}a_j(R)(n-1-j)!\le\sum_{j>i}(n-1-j)\,(n-1-j)!=(n-1-i)!-1$, it follows that $L(R)<L(S)$.
\end{proof}

\begin{table}[ht]
\centering\small
\begin{tabular}{@{}lll@{}}
 & Lehmer label $\lambda$ (this book) & classical rank $L$\\ \hline
place value of $a_i(R)$ & $\Pi_i=n!/(n-i)!$ & $(n-1-i)!$\\
least significant digit & $a_0(R)$ & $a_{n-1}(R)=0$\\
$a_i(R)$ by~\eqref{eq:mixed} & $\dg{n-i}{0}{\lambda/\Pi_i}$ & $\dg{n-i}{0}{L/(n-1-i)!}$\\
$R(0),\dots,R(h-1)$ determined by & $\lambda\bmod\Pi_h$ & $\fl{L/(n-h)!}$\\
transposition $s_i=(i\ i{+}1)$ & $\Pi_i$ & $(n-1-i)!$\\
reversal $i\mapsto n-1-i$ & $n!-1$ & $n!-1$\\
sign & $(-1)^{\sum_i a_i(R)}$ & $(-1)^{\sum_i a_i(R)}$\\
lexicographic order of & $(a_{n-1},\dots,a_0)$ & one-line words\\
\end{tabular}
\caption{Two readings of Lehmer's code of $R\in\Sym_n$.}\label{tab:lehmerreadings}
\end{table}

The difference matters in computation; Table~\ref{tab:lehmerreadings} lists the main points. The low digits of $\lambda(R)$
describe the beginning of the word of $R$: $a_0(R)=\lambda(R)\bmod n$ costs one division by a one-word integer, whereas
$a_0(R)=\fl{L(R)/(n-1)!}$ needs a division by a multiword integer. More generally $\lambda(R)\bmod\Pi_h$ is determined by
$(R(0),\dots,R(h-1))$, so it is constant on each left coset $RH$ of the pointwise stabiliser $H$ of $\{0,\dots,h-1\}$: if $S\in H$,
then $(R\circ S)(i)=R(i)$ for $i<h$. Truncating the label modulo $\Pi_h$ is therefore a step down the point-stabiliser chain.

\begin{example}\label{ex:lehmer4}
Let $n=4$ and $R=(2,0,3,1)$ in one-line notation. Then $(a_0,a_1,a_2,a_3)=(2,0,1,0)$, so
$\lambda(R)=2\cdot1+0\cdot4+1\cdot12+0\cdot24=14$ and $L(R)=2\cdot3!+0\cdot2!+1\cdot1!+0\cdot0!=13$. Indeed twelve words begin
with $0$ or $1$, and $(2,0,1,3)$ is the only word beginning with $2$ that precedes $R$. Proposition~\ref{prop:seqdigits}
recovers the code from $14$ with one small division per digit,
\[
14=3\cdot4+2,\qquad 3=1\cdot3+0,\qquad 1=0\cdot2+1,\qquad 0=0\cdot1+0,
\]
while the closed formula gives a single digit directly, for instance $\md{2}{14}=\fl{14/12}-2\fl{14/24}=1$.
\end{example}

\begin{remark}[Which definition is practical]\label{rem:digitpractical}
Lehmer's code is not a different definition of digits: it is a choice of radices and of the order of the digits, and once these
are fixed its digits are the $\md{h}{x}$ of~\eqref{eq:mixed}. As a definition, the digit function~\eqref{eq:digit} is the more
useful one. It is exact and closed, it applies to real $x$ and to negative $k$ (for instance $\dg{N}{-1}{x}=\fl{Nx}-N\fl{x}$ is the
first digit after the point), and the identities of this chapter are proved from it. As an algorithm it is practical
provided that the floors are exact integer divisions and that a full expansion is computed by~\eqref{eq:seqdigit}, or by
Remark~\ref{rem:digitdc} for very long integers. The accompanying scripts follow exactly this pattern.
\end{remark}

\section{The \texorpdfstring{$\Bf$}{B}-function}\label{sec:bfunction}

The digit function extracts the digits of a code; the \emph{$\Bf$-function} tests them. Together they are
the two basic building blocks of digital mathematics: every finite law, every comparison of labels and every
carry of positional arithmetic is a closed formula in these two functions. This section collects the
properties of the $\Bf$-function needed later, following the treatment in the digital mathematics notes
(Chapter 3, Section 4.8 and Section 11.4 there). We write it as $\Bf$, in calligraphic type, which also keeps it apart from the
digit labels $A,B,C$ used in Chapter~\ref{ch:diagrammatic}.

\subsection{Definition and values}

\begin{definition}\label{def:Bfun}
For real $x,y$ the \emph{$\Bf$-function} is
\[
\Bf(x,y)=\tfrac12\big(\sgn(x+y)-\sgn(x-y)\big),
\]
and we abbreviate $\Bf(x)=\Bf(x,\tfrac12)$.
\end{definition}

\begin{proposition}\label{prop:Bvalues}
\begin{enumerate}
\item For $y>0$, $\Bf(x,y)=1$ if $|x|<y$, $\Bf(x,y)=\tfrac12$ if $|x|=y$, and $\Bf(x,y)=0$ if $|x|>y$.
      Moreover $\Bf(x,0)=0$ and $\Bf(-x,y)=\Bf(x,y)$.
\item For an integer $u$, $\Bf(u)=1$ if $u=0$ and $\Bf(u)=0$ otherwise. In particular
      $\Bf(u)=\dg{n}{u}{1}$ for $u\in\Ap{n}$, the digital Kronecker delta used earlier.
\item For integers $a\le b$ and $x$,
      $\Bf\!\big(x-\tfrac{a+b}{2},\tfrac{b-a+1}{2}\big)=1$ if $a\le x\le b$ and $0$ otherwise.
\end{enumerate}
\end{proposition}
\begin{proof}
(1) If $|x|<y$ then $x+y>0>x-y$ and the bracket is $2$; if $x=\pm y$ exactly one sign vanishes and
the other is $\pm1$ with the right sign; if $|x|>y$ both signs agree. For $y=0$ the two signs coincide, and
$\sgn$ is odd, which gives the symmetry. (2) is (1) with $y=\tfrac12$, since a nonzero integer has $|u|\ge1$.
(3) The condition $|x-\tfrac{a+b}{2}|<\tfrac{b-a+1}{2}$ means $a-\tfrac12<x<b+\tfrac12$, and the
boundary values are half-integers, never attained by $x$.
\end{proof}

Thus $\Bf$ is the indicator of an open interval (with the value $\tfrac12$ on its ends), and on integers it
is an exact \emph{if-then-else}: the map sending $a$ to $b$ and every other integer to $c$ is
$f(n)=c+(b-c)\,\Bf(n-a)$.

\subsection{Identities}

\begin{proposition}[Splitting and block coalescence]\label{prop:Bsplit}
For real $x,y,z,m,\varepsilon$ and $n\ge1$,
\begin{align}
\Bf(x,y+z)&=\Bf(x+z,y)+\Bf(x-y,z),\label{eq:Bsplit}\\
\sum_{k=0}^{n-1}\Bf(m-2k\varepsilon,\varepsilon)&=\Bf\big(m-(n-1)\varepsilon,\,n\varepsilon\big).\label{eq:Bblock}
\end{align}
\end{proposition}
\begin{proof}
The right side of~\eqref{eq:Bsplit} is
$\tfrac12\big(\sgn(x+y+z)-\sgn(x+z-y)+\sgn(x-y+z)-\sgn(x-y-z)\big)$; the middle terms cancel and what
remains is $\Bf(x,y+z)$. For~\eqref{eq:Bblock} use induction on $n$: the case $n=1$ is trivial, and
\eqref{eq:Bsplit} with $x=m-n\varepsilon$, $y=n\varepsilon$, $z=\varepsilon$ gives
$\Bf(m-n\varepsilon,(n+1)\varepsilon)=\Bf(m-(n-1)\varepsilon,n\varepsilon)+\Bf(m-2n\varepsilon,\varepsilon)$.
\end{proof}

Splitting cuts an interval into two adjacent ones; block coalescence glues $n$ adjacent intervals of
width $2\varepsilon$ into one. With $\varepsilon=\tfrac12$ it says that a sum of consecutive Kronecker deltas is a
single $\Bf$-function, which is how carries will be written below.

\begin{proposition}[Integer deltas]\label{prop:Bdelta}
For integers $x,a,b$ and any finite set $S$ of integers containing $a$,
\[
\begin{aligned}
&\Bf(x-a)^2=\Bf(x-a),\qquad \Bf(x-a)\Bf(x-b)=\Bf(x-a)\Bf(a-b),\\
&\Bf(x-a)=\sum_{k\in S}\Bf(x-k)\Bf(k-a).
\end{aligned}
\]
\end{proposition}
\begin{proof}
All three follow from Proposition~\ref{prop:Bvalues}(2): the deltas take values $0,1$; both sides of the
second identity are $1$ exactly when $x=a=b$; in the sum only $k=a$ can contribute, and it contributes
$\Bf(x-a)$.
\end{proof}

The last identity (\emph{propagation}) may be iterated: inserting any number of summed deltas between $x$
and $a$ changes nothing. It is the digital form of composing Kronecker deltas.

The deltas are also the evaluation kernel of every finite operation: a code is turned into a function
by selecting one of its digits.

\begin{proposition}[Evaluation kernel]\label{prop:Bkernel}
Let $R\in[0,N^{N^2})$ be the code of a binary law $*$ on $\Ap{N}$. Then for $a,b\in\Ap{N}$
\[
a*b=\dg{N}{aN+b}{R}=\sum_{n=0}^{N^2-1}\dg{N}{n}{R}\,\Bf(n-aN-b).
\]
More generally, if $R\in[0,N^{N^r})$ is the code of an $r$-ary operator on $\Ap{N}$, then
\[
R(x_0,\dots,x_{r-1})=\sum_{n=0}^{N^r-1}\dg{N}{n}{R}\,\Bf\Big(n-\sum_{j=0}^{r-1}x_jN^{j}\Big).
\]
\end{proposition}
\begin{proof}
By Proposition~\ref{prop:Bvalues}(2) only the term with $n=aN+b$ (respectively $n=\sum_jx_jN^j$) survives, and it
equals the digit selected by the coding of Theorem~\ref{thm:coding}.
\end{proof}

The first form reads one digit of $R$; the second writes the operator as a linear combination of deltas with the
digits of $R$ as coefficients, which is the form used to manipulate tables algebraically (for instance in
Theorem~\ref{thm:Bfactor} below).

\subsection{The \texorpdfstring{$\Bf$}{B}-function and the digits}

\begin{theorem}[Factorization]\label{thm:Bfactor}
For integers $N\ge2$ and $a,b\in[0,N^M)$,
\[
\Bf(a-b)=\prod_{k=0}^{M-1}\Bf\big(\dg{N}{k}{a}-\dg{N}{k}{b}\big).
\]
The same holds for mixed-radix labels: two labels are equal if and only if all their digits are.
\end{theorem}
\begin{proof}
By uniqueness of the radix-$N$ (or mixed-radix) representation, $a=b$ if and only if all digits agree,
and the product is $1$ exactly in that case.
\end{proof}

\begin{corollary}[Boolean form]\label{cor:Bboolean}
For $N=2$, writing $a_k=\dg{2}{k}{a}$ and $b_k=\dg{2}{k}{b}$,
\[
\Bf(a-b)=\prod_{k=0}^{M-1}\big(1-a_k-b_k+2a_kb_k\big).
\]
\end{corollary}
\begin{proof}
By propagation over $S=\{0,1\}$,
$\Bf(a_k-b_k)=\Bf(a_k)\Bf(b_k)+\Bf(a_k-1)\Bf(1-b_k)=(1-a_k)(1-b_k)+a_kb_k$.
\end{proof}

Applied to codes, Theorem~\ref{thm:Bfactor} turns the equality of two magmas into an entrywise comparison of
their tables: for codes $R,R'$ of binary laws on $\Ap{N}$,
\[
\Bf(R-R')=\prod_{m=0}^{N^2-1}\Bf\big(\dg{N}{m}{R}-\dg{N}{m}{R'}\big),\qquad
\sum_{m=0}^{N^2-1}\Big(1-\Bf\big(\dg{N}{m}{R}-\dg{N}{m}{R'}\big)\Big)
\]
being, respectively, the test $R=R'$ and the number of entries in which the two tables differ.

\subsection{Carries}

\begin{proposition}[The carry as a $\Bf$-function]\label{prop:Bcarry}
For $N\ge2$ and $u,v\in\Ap{N}$ let $c_N(u,v)=\fl{(u+v)/N}$. Then
\[
c_N(u,v)=\Bf\Big(u+v-\frac{3N-2}{2},\ \frac{N-1}{2}\Big),\qquad
u+v=\dg{N}{0}{u+v}+N\,c_N(u,v).
\]
\end{proposition}
\begin{proof}
Since $0\le u+v\le2N-2$, the carry is $1$ exactly when $N\le u+v\le2N-2$; apply
Proposition~\ref{prop:Bvalues}(3) with $a=N$, $b=2N-2$. Equivalently, the carry is the sum of the deltas
$\Bf(u+v-s)$ for $s=N,\dots,2N-2$, coalesced by~\eqref{eq:Bblock}. The second identity is the definition of the
lowest digit.
\end{proof}

\begin{proposition}[The carry is a cocycle]\label{prop:Bcocycle}
For $u,v,w\in\Ap{N}$,
\[
c_N(u,v)+c_N\big(\dg{N}{0}{u+v},w\big)=c_N(v,w)+c_N\big(u,\dg{N}{0}{v+w}\big)=\fl{\frac{u+v+w}{N}}.
\]
\end{proposition}
\begin{proof}
By Proposition~\ref{prop:Bcarry} applied twice,
\[
u+v+w=\dg{N}{0}{u+v}+w+N\,c_N(u,v)=\dg{N}{0}{u+v+w}+N\big[c_N(u,v)+c_N(\dg{N}{0}{u+v},w)\big],
\]
and the bracket is $\fl{(u+v+w)/N}$ because the lowest digit is less than $N$. Grouping $v+w$ first gives the other
expression.
\end{proof}

The ordinary addition of integers thus splits into two parts. For $a,b\in[0,N^M)$ define the
\emph{carry-free sum} and the \emph{carry part}
\[
a\oplus_N b=\sum_{k=0}^{M-1}N^k\,\dg{N}{0}{\dg{N}{k}{a}+\dg{N}{k}{b}},\qquad
a\,\hat{+}_N\,b=a+b-a\oplus_N b .
\]
The carry-free sum adds each digit modulo $N$ independently of the others: it is the law of
$C_N\oplus\cdots\oplus C_N$, an operator acting position by position (a graded operator in the sense of
Definition~\ref{def:dggraded}). All coupling between the positions sits in the carry part. For instance, in base $10$,
$56782+36754=93536$ splits as $82436+11100$. Chapter~\ref{ch:diagrammatic} generalises exactly this picture:
a diagram without arrows is a graded operator, a direct sum of its layers, and the carry arrows are
arbitrary admissible replacements of the carry $c_N$ (Section~\ref{sec:dgordinarycarry}).

\begin{remark}
Carries also count. By Kummer's theorem the exponent of $p$ in $\binom{a+b}{a}$ is the number of carries
in the base-$p$ addition $a+b$, that is, the number of positions at which the propagated carry is $1$. This is the
form in which carries enter the proof of the Sylow theorems in Chapter~\ref{ch:digitaltheorems}
(Lemma~\ref{lem:dtlucas} and the discussion after it).
\end{remark}

\section{Finite sets and characteristic functions}\label{sec:sets}

In the rest of the book every finite object is a set of integer codes: subgroups, cosets and classes are subsets of an alphabet $\Ap{n}$. This section records how such subsets are handled with digits and with the $\Bf$-function. For integers $x$ we only use that $\Bf(x)=1$ if $x=0$ and $\Bf(x)=0$ otherwise (Section~\ref{sec:bfunction} and the exercises below).

\begin{definition}
Let $S\subseteq\Ap{n}$. Its \emph{characteristic function} and its \emph{code} are
\begin{equation}\label{eq:charfun}
\mathbf{1}_S(x)=\sum_{s\in S}\Bf(x-s)\quad(x\in\Ap{n}),\qquad R_S=\sum_{s\in S}2^{s}.
\end{equation}
\end{definition}

\begin{proposition}\label{prop:subsetcode}
For $S\subseteq\Ap n$ and $x\in\Ap n$ one has $\mathbf 1_S(x)=\dg{2}{x}{R_S}$ and $|S|=\sum_{x\in\Ap n}\mathbf 1_S(x)$. The map $S\mapsto R_S$ is a bijection from the subsets of $\Ap n$ onto $\Ap{2^n}$; in particular $\Ap n$ has $2^n$ subsets.
\end{proposition}
\begin{proof}
At most one term of the sum in \eqref{eq:charfun} is nonzero, and it equals $1$ exactly when $x\in S$; summing over $x$ counts each element of $S$ once. The binary digits of $R_S$ are $\mathbf 1_S(0),\dots,\mathbf 1_S(n-1)$, so the last claim is Lemma~\ref{lem:basic}(3) with $p=2$.
\end{proof}

\begin{proposition}[Boolean operations]\label{prop:boolean}
For $S,T\subseteq\Ap n$,
\[
\mathbf 1_{S\cap T}=\mathbf 1_S\mathbf 1_T,\quad
\mathbf 1_{S\cup T}=\mathbf 1_S+\mathbf 1_T-\mathbf 1_S\mathbf 1_T,\quad
\mathbf 1_{S\triangle T}=\mathbf 1_S+\mathbf 1_T-2\,\mathbf 1_S\mathbf 1_T,\quad
\mathbf 1_{\Ap n\setminus S}=1-\mathbf 1_S,
\]
and $R_{S\cap T}=\sum_{k<n}\dg2k{R_S}\,\dg2k{R_T}\,2^k$. If $P\subseteq\Ap n$ is the set of elements with a given property and $A\subseteq\Ap n$, the statements ``some $x\in A$ has the property'' and ``every $x\in A$ has the property'' have the truth values
\[
1-\Bf\Big(\sum_{x\in A}\mathbf 1_P(x)\Big)\qquad\text{and}\qquad \prod_{x\in A}\mathbf 1_P(x).
\]
\end{proposition}
\begin{proof}
Each identity is checked on the four possible values of $(\mathbf 1_S(x),\mathbf 1_T(x))$, and the formula for $R_{S\cap T}$ then follows from Proposition~\ref{prop:subsetcode}. The sum in the first truth value is a non-negative integer that vanishes exactly when no $x\in A$ lies in $P$.
\end{proof}

\begin{theorem}[Multiplicities]\label{thm:multiplicity}
Let $S_1,\dots,S_m\subseteq\Ap n$ and let $c(x)=\sum_{j=1}^m\mathbf 1_{S_j}(x)$ be the number of these sets containing $x$. Then $x$ lies in exactly $k$ of the sets if and only if $\Bf\big(k-c(x)\big)=1$, and, for $k\ge1$, in at least one and at most $k$ of them if and only if $\Bf\big(k+1-2c(x),k\big)=1$. In particular
\begin{equation}\label{eq:unionB}
\mathbf 1_{S_1\cup\dots\cup S_m}(x)=1-\Bf\big(c(x)\big)=\Bf\big(m+1-2c(x),m\big).
\end{equation}
\end{theorem}
\begin{proof}
The first claim holds because $k-c(x)$ is an integer. For the second, $|k+1-2c|<k$ holds exactly when $1\le c\le k$, and $|k+1-2c|=k$ is impossible for an integer $c$, so the boundary convention of the two-argument $\Bf$-function plays no role. For \eqref{eq:unionB} take $k=m$ and note that $c(x)\le m$.
\end{proof}

\begin{remark}
When a further set $S$ joins the union, the right-hand side of \eqref{eq:unionB} changes linearly: the first argument gains $1-2\,\mathbf 1_S(x)$ and the second gains $1$. The classical expansion of Corollary~\ref{cor:inclexcl} needs $2^m-1$ terms instead of $m$.
\end{remark}

\begin{corollary}[Inclusion--exclusion]\label{cor:inclexcl}
With the notation of Theorem~\ref{thm:multiplicity}, and with $\bigcap_{j\in\emptyset}S_j=\Ap n$,
\begin{gather*}
\Bf\big(c(x)\big)=\prod_{j=1}^m\big(1-\mathbf 1_{S_j}(x)\big)=\sum_{J\subseteq\{1,\dots,m\}}(-1)^{|J|}\,\mathbf 1_{\bigcap_{j\in J}S_j}(x),
\\
\Big|\Ap n\setminus\bigcup_{j=1}^m S_j\Big|=\sum_{J\subseteq\{1,\dots,m\}}(-1)^{|J|}\,\Big|\bigcap_{j\in J}S_j\Big|.
\end{gather*}
\end{corollary}
\begin{proof}
Both $\Bf(c(x))$ and the product equal $1$ exactly when $x$ lies in no $S_j$, and $0$ otherwise. Expand the product and use $\prod_{j\in J}\mathbf 1_{S_j}=\mathbf 1_{\bigcap_{j\in J}S_j}$ (Proposition~\ref{prop:boolean}); then sum over $x$ (Proposition~\ref{prop:subsetcode}).
\end{proof}

\begin{example}[Generators of $C_{12}$]\label{ex:gen12}
An element of $C_{12}$, coded by $\Ap{12}$, generates the group exactly when it lies in no maximal subgroup. The maximal subgroups are $M_1=\{0,2,\dots,10\}$ and $M_2=\{0,3,6,9\}$, with codes $R_{M_1}=1365$ and $R_{M_2}=585$. Their common binary digits are in positions $0$ and $6$, so $R_{M_1\cap M_2}=65$ and $M_1\cap M_2=\{0,6\}$. Corollary~\ref{cor:inclexcl} gives $12-6-4+2=4$ generators, namely $1,5,7,11$, in agreement with $\varphi(12)=4$.
\end{example}

\paragraph{Relations and quotient sets.} A binary relation on $\Ap n$ is a function $r:\Ap n^2\to\{0,1\}$; via $(x,y)\mapsto x+ny$ it is a subset of $\Ap{n^2}$ and so, by Proposition~\ref{prop:subsetcode}, it is determined by the code $\sum_{x,y}r(x,y)\,2^{x+ny}$. It is an \emph{equivalence relation} if $r(x,x)=1$, $r(x,y)=r(y,x)$ and $r(x,y)\,r(y,z)\le r(x,z)$ for all $x,y,z$. The class of $x$ then has characteristic function $\mathbf 1_{[x]}(y)=r(x,y)$, and the classes form a partition of $\Ap n$.

\begin{proposition}\label{prop:classes}
Let $r$ be an equivalence relation on $\Ap n$. The number of its classes is
\[
\sum_{x\in\Ap n}\frac{1}{\sum_{y\in\Ap n}r(x,y)}=\sum_{x\in\Ap n}\ \prod_{y<x}\big(1-r(x,y)\big).
\]
\end{proposition}
\begin{proof}
In the first sum each class $C$ contributes $|C|$ terms equal to $1/|C|$. In the second, the product equals $1$ exactly when $x$ is the least element of its class, so each class is counted once.
\end{proof}

\begin{example}
For $d\mid n$, congruence modulo $d$ on $\Ap n$ is $r(x,y)=\Bf\big(\dg{d}{0}{x-y}\big)$ by Lemma~\ref{lem:basic}(2). Every class has $n/d$ elements, so Proposition~\ref{prop:classes} gives $n\cdot d/n=d$ classes, with least elements $0,1,\dots,d-1$. Once a group law is given on $\Ap n$, the cosets of a subgroup and the conjugacy classes arise in the same way, and Proposition~\ref{prop:classes} counts them.
\end{example}

\section{Exercises}
\begin{enumerate}
\item Prove directly from~\eqref{eq:digit} that $\dg{N}{k}{x}\in\Ap N$ for every integer $x\ge0$ and that $x=\sum_{k\ge0}\dg{N}{k}{x}\,N^k$.
\item Show that $\dg{N}{k}{x+N^{k+1}y}=\dg{N}{k}{x}$ and $\dg{N}{k}{Nx}=\dg{N}{k-1}{x}$ for all integers $x,y$ and $k\ge1$.
\item Show that $\dg{N^2}{k}{x}=\dg{N}{2k}{x}+N\,\dg{N}{2k+1}{x}$ for every integer $x\ge0$, and state the analogous rule for base $N^{n}$.
\item Using Theorem~\ref{thm:coding}, show that exactly $N^{N(N+1)/2}$ of the codes in $\Ap{N^{N^2}}$ define commutative binary operations on $\Ap N$.
\item Show that for integers $x$ one has $\Bf(x)=1$ if $x=0$ and $\Bf(x)=0$ otherwise, and that for integers $x$ and $n\ge1$,
$\Bf\big(x-\tfrac{n-1}{2},\tfrac n2\big)=1$ if $0\le x<n$ and $0$ otherwise. Write the indicator of $\Ap n$ and the Kronecker delta $\delta_{m,n}$ in this form.
\item Let $a,b\ge0$ and let $c_k=\fl{(a\bmod N^k+b\bmod N^k)/N^k}$. Show that $c_k\in\{0,1\}$ and
$\dg{N}{k}{a+b}\equiv\dg{N}{k}{a}+\dg{N}{k}{b}+c_k\pmod N$, and compare with Proposition~\ref{prop:Bcarry}.
\item Using Kummer's theorem, show that $\binom{2n}{n}$ is even for every $n\ge1$, and that an odd prime $p$ does not divide
$\binom{2n}{n}$ if and only if every base-$p$ digit of $n$ is less than $p/2$.
\item Show that $R_{S\cup T}=R_S+R_T-R_{S\cap T}$ and $R_{S\triangle T}=R_S+R_T-2R_{S\cap T}$ for $S,T\subseteq\Ap n$.
\item In the notation of Theorem~\ref{thm:multiplicity}, show that exactly $\sum_{x\in\Ap n}\Bf\big(k-c(x)\big)$ elements lie in exactly $k$ of the sets, and that $\sum_{k=1}^m\Bf\big(k-c(x)\big)=1-\Bf\big(c(x)\big)$.
\item Use Corollary~\ref{cor:inclexcl} to show that $\Ap n$ has $n!\sum_{k=0}^n(-1)^k/k!$ permutations without fixed points.
\item Apply Corollary~\ref{cor:inclexcl} to the maximal subgroups of $C_n$ to prove $\varphi(n)=n\prod_{p\mid n}(1-1/p)$.
\item Let $f:\Ap n\to\Ap m$. Show that $r(x,y)=\Bf\big(f(x)-f(y)\big)$ is an equivalence relation whose classes are the fibres of $f$, and that it has as many classes as $f$ has values.
\item Show that a relation $r$ on $\Ap n$ is transitive if and only if $\prod_{x,y,z}\big(1-r(x,y)r(y,z)+r(x,y)r(y,z)r(x,z)\big)=1$.
\end{enumerate}

\chapter{Finite magmas and \texorpdfstring{$r$}{r}-ary operators}\label{ch:magmas}

\section{Permutations}\label{sec:permutations}
\subsection{Permutation codes}
\begin{definition}
A \emph{permutation code} of degree $N$ is an integer $R\in[0,N^{N})$ whose $N$ digits
$a_n=\dg{N}{n}{R}$, $n\in\Ap{N}$, are pairwise distinct. It acts by
$R(x)=\dg{N}{x}{R}$. The set of such codes is $\Sym_N$; by Theorem~\ref{thm:coding} (with $r=1$) it is in bijection with the
bijections of $\Ap{N}$, so $|\Sym_N|=N!$. In Cauchy notation,
$R=\left(\begin{smallmatrix}0&1&\cdots&N-1\\ a_0&a_1&\cdots&a_{N-1}\end{smallmatrix}\right)$.
\end{definition}

\begin{proposition}\label{prop:permops}
For $R,S\in\Sym_N$:
\begin{align}
\iota_N&=\sum_{n=0}^{N-1}n\,N^{n}=\frac{N\,(N^{N+1}-2N^{N}+1)}{(N-1)^{2}}\quad\text{(identity)},\label{eq:id}\\
R^{-1}&=\sum_{n=0}^{N-1}n\,N^{\dg{N}{n}{R}}\quad\text{(inverse)},\label{eq:inv}\\
R\circ S&=\sum_{n=0}^{N-1}\dg{N}{\dg{N}{n}{S}}{R}\,N^{n}\quad\text{(composition, right to left)}.\label{eq:comp}
\end{align}
With $\circ$, $\Sym_N$ is a group, the \emph{symmetric group}.
\end{proposition}
\begin{proof}
Each formula is~\eqref{eq:code} for the map in question: $\iota_N$ has digit $n$ at position $n$; $R^{-1}$ has digit $n$ at
position $R(n)$; $R\circ S$ has digit $R(S(n))$ at position $n$. The closed form of $\iota_N$ is the sum
$\sum_{n<N}nx^n$ at $x=N$. Composition of maps is associative.
\end{proof}

\begin{example}
For $N=6$, $\left(\begin{smallmatrix}0&1&2&3&4&5\\3&5&2&1&4&0\end{smallmatrix}\right)$ has code
$3+5\cdot6+2\cdot6^2+1\cdot6^3+4\cdot6^4=5505$, cycle form $(0315)$, and its inverse has code $13055$.
\end{example}

\subsection{Cycles and conjugacy}
\begin{proposition}[Cycle code]\label{prop:cycle}
Suppose $R$ decomposes into disjoint cycles
\[
(j^{(h)}_0\,j^{(h)}_1\cdots j^{(h)}_{s_h-1}),\qquad h=1,\dots,c,
\]
with fixed points included as $1$-cycles. Then
\begin{equation}\label{eq:cyclecode}
R=\sum_{h=1}^{c}\sum_{k=0}^{s_h-1}j^{(h)}_{\dg{s_h}{0}{k+1}}\;N^{\,j^{(h)}_{k}}.
\end{equation}
\end{proposition}
\begin{proof}
The digit of $R$ at position $j^{(h)}_k$ must be its image $j^{(h)}_{k+1}$, indices mod $s_h$; every position occurs exactly once.
\end{proof}

\begin{example}
$(236)\circ(01246)=(0136)(24)$ in $\Sym_7$: by~\eqref{eq:comp}, $0\mapsto1\mapsto1$, $1\mapsto2\mapsto3$, $2\mapsto4\mapsto4$,
$3\mapsto3\mapsto6$, $4\mapsto6\mapsto2$, $6\mapsto0\mapsto0$.
\end{example}

\begin{proposition}[Conjugation]\label{prop:conj}
For $R,S\in\Sym_N$,
\[
S\circ R\circ S^{-1}=\sum_{n=0}^{N-1}\dg{N}{\dg{N}{n}{R}}{S}\,N^{\dg{N}{n}{S}}.
\]
Conjugation replaces each symbol of the cycle form of $R$ by its image under $S$. Hence two permutations are conjugate if and only if
they have the same cycle type, and the conjugacy classes of $\Sym_N$ correspond to the partitions of $N$.
\end{proposition}
\begin{proof}
$SRS^{-1}$ sends $S(n)$ to $S(R(n))$, which is the stated digit placement. Conversely, given two permutations with the same
cycle type, the bijection matching their cycles symbol by symbol is a conjugating $S$.
\end{proof}

\subsection{Sign and the alternating group}
\begin{definition}
The \emph{sign} of $R\in\Sym_N$ is
$\displaystyle\sgn R=\prod_{0\le i<j\le N-1}\frac{\dg{N}{j}{R}-\dg{N}{i}{R}}{j-i}\in\{\pm1\}$.
\end{definition}

\begin{theorem}\label{thm:sign}
$\sgn:\Sym_N\to\{\pm1\}$ is a homomorphism; an $s$-cycle has sign $(-1)^{s-1}$; and $\sgn R=(-1)^{N-c(R)}$, where $c(R)$ is the
number of cycles of $R$ (fixed points included).
\end{theorem}
\begin{proof}
$\sgn(R\circ S)=\prod_{i<j}\frac{R(S(j))-R(S(i))}{S(j)-S(i)}\prod_{i<j}\frac{S(j)-S(i)}{j-i}$, and the first factor is
$\sgn R$ after re-indexing over unordered pairs $\{S(i),S(j)\}$. A transposition $(j\,k)$, $j<k$, reverses exactly the
$2(k-j)-1$ pairs it meets, so has sign $-1$. Moreover
$(j_0\,j_1\cdots j_{s-1})=(j_0\,j_{s-1})\circ\cdots\circ(j_0\,j_2)\circ(j_0\,j_1)$, which gives $(-1)^{s-1}$, and then
$\prod_h(-1)^{s_h-1}=(-1)^{N-c}$.
\end{proof}

\begin{definition}
The \emph{alternating group} is $\Alt_N=\ker\sgn$. For $N\ge2$ it has index $2$ and order $N!/2$.
\end{definition}

\begin{example}[$\Sym_3$]
The six permutation codes of degree $3$ are
\[
\begin{array}{c|cccccc}
R & 5 & 7 & 11 & 15 & 19 & 21\\\hline
(a_0,a_1,a_2) & (2,1,0) & (1,2,0) & (2,0,1) & (0,2,1) & (1,0,2) & (0,1,2)\\
\sgn R & -1 & +1 & +1 & -1 & -1 & +1
\end{array}
\]
so $\Alt_3=\{7,11,21\}$ and $\iota_3=21$.
\end{example}

\section{Magmas and their codes}
\begin{definition}
A \emph{magma} of order $N$ is a map $*:\Ap{N}\times\Ap{N}\to\Ap{N}$. By Theorem~\ref{thm:coding} (with $x_0=b$, $x_1=a$) there is a unique
$R\in[0,N^{N^2})$ such that
\begin{equation}\label{eq:magma}
a*b=\dg{N}{aN+b}{R}.
\end{equation}
We write $(\Ap{N},R)$ for the magma. There are $N^{N^{2}}$ magmas of order $N$.
\end{definition}

\begin{proposition}[Translations are block digits]\label{prop:rows}
Let $\lambda_a(b)=a*b$ and $\rho_b(a)=a*b$. Then $\lambda_a$ has the (not necessarily bijective) code
$L_a=\dg{N^{N}}{a}{R}$, $\rho_b$ has code $C_b=\sum_{a}\dg{N}{aN+b}{R}\,N^{a}$, and
\begin{equation}\label{eq:rowdecomp}
R=\sum_{a=0}^{N-1}L_a\,(N^{N})^{a}.
\end{equation}
\end{proposition}
\begin{proof}
Lemma~\ref{lem:block} with $n=N$ gives $\dg{N}{aN+b}{R}=\dg{N}{b}{L_a}$; \eqref{eq:rowdecomp} is Lemma~\ref{lem:basic}(3) in base $N^N$.
\end{proof}

\begin{definition}[Axioms]
A magma $(\Ap{N},R)$ is
\begin{itemize}
\item a \emph{semigroup} if it is associative, $a*(b*c)=(a*b)*c$;
\item \emph{commutative} if $a*b=b*a$;
\item a \emph{monoid} if it is a semigroup in which $0$ is an identity, $a*0=0*a=a$;
\item a \emph{quasigroup} if every $\lambda_a$ and every $\rho_b$ is bijective (the Cayley table is a Latin square);
\item a \emph{loop} if it is a quasigroup in which $0$ is an identity;
\item a \emph{group} if it is an associative loop. A commutative group is \emph{Abelian}.
\end{itemize}
Placing the identity at $0$ loses no generality, by relabelling (Proposition~\ref{prop:relabel}).
\end{definition}

\begin{proposition}[Digital form of the axioms]\label{prop:axioms}
For a magma $(\Ap{N},R)$:
\begin{enumerate}
\item it is commutative if and only if $R=R^{\mathrm{op}}$, where $R^{\mathrm{op}}=\sum_{a,b}\dg{N}{bN+a}{R}\,N^{aN+b}$;
\item $0$ is an identity if and only if $L_0=C_0=\iota_N$;
\item it is a quasigroup if and only if $L_a,C_b\in\Sym_N$ for all $a,b$;
\item if it is a quasigroup, then
\[
\sum_a\dg{N}{aN+b}{R}=\sum_b\dg{N}{aN+b}{R}=\binom{N}{2},
\qquad
\sum_{a,b}\dg{N}{aN+b}{R}=N\binom{N}{2}.
\]
\end{enumerate}
\end{proposition}
\begin{proof}
(1)--(3) restate the definitions through Proposition~\ref{prop:rows}. For (4), each row and each column contains every element of $\Ap{N}$ once.
\end{proof}

\begin{theorem}[Digital Cayley theorem]\label{thm:cayley}
A loop $(\Ap{N},R)$ is a group if and only if
\begin{equation}\label{eq:cayley}
L_{a*b}=L_a\circ L_b\qquad\text{for all }a,b\in\Ap{N}.
\end{equation}
In that case $a\mapsto L_a=\dg{N^{N}}{a}{R}$ is an injective homomorphism into $\Sym_N$; every group of order $N$ is a subgroup of $\Sym_N$.
\end{theorem}
\begin{proof}
Associativity for all $c$ says exactly $\lambda_a\circ\lambda_b=\lambda_{a*b}$, i.e.\ \eqref{eq:cayley} by~\eqref{eq:comp}.
Injectivity: $L_a(0)=a*0=a$.
\end{proof}

\begin{proposition}[Inverses]
In a loop, each $a$ has a unique right inverse and a unique left inverse, namely $a^{\rho}=L_a^{-1}(0)=\dg{N}{0}{\cdot}$-position of $0$ in $L_a$, i.e.
$a^{\rho}=\dg{N}{0}{\textstyle\sum_n n\,N^{\dg{N}{n}{L_a}}}$, and symmetrically with $C_a$. In a group they coincide.
\end{proposition}
\begin{proof}
$a*x=0$ has the unique solution $x=\lambda_a^{-1}(0)$, whose code is given by~\eqref{eq:inv}. In a group,
$a^{\lambda}=a^{\lambda}*(a*a^{\rho})=(a^{\lambda}*a)*a^{\rho}=a^{\rho}$.
\end{proof}

\begin{definition}
A \emph{Moufang loop} is a loop satisfying the equivalent identities
$z*(x*(z*y))=((z*x)*z)*y$, $x*(z*(y*z))=((x*z)*y)*z$, $(z*x)*(y*z)=(z*(x*y))*z$, $(z*x)*(y*z)=z*((x*y)*z)$.
\end{definition}

\begin{example}
The smallest nonassociative Moufang loop, $M(S_3,2)$, has order $12$:
{\scriptsize
\[
\begin{array}{c|cccccccccccc}
* &0&1&2&3&4&5&6&7&8&9&10&11\\\hline
0&0&1&2&3&4&5&6&7&8&9&10&11\\
1&1&0&4&5&2&3&7&6&10&11&8&9\\
2&2&5&0&4&3&1&8&11&6&10&9&7\\
3&3&4&5&0&1&2&9&10&11&6&7&8\\
4&4&3&1&2&5&0&10&9&7&8&11&6\\
5&5&2&3&1&0&4&11&8&9&7&6&10\\
6&6&7&8&9&11&10&0&1&2&3&5&4\\
7&7&6&10&11&9&8&1&0&5&4&2&3\\
8&8&11&6&10&7&9&2&4&0&5&3&1\\
9&9&10&11&6&8&7&3&5&4&0&1&2\\
10&10&9&7&8&6&11&4&2&3&1&0&5\\
11&11&8&9&7&10&6&5&3&1&2&4&0
\end{array}
\]}
\end{example}

\section{Homomorphisms and relabelling}
\begin{definition}
A \emph{homomorphism} $f:(\Ap{N},R)\to(\Ap{q},S)$ is a map $f:\Ap{N}\to\Ap{q}$ with
\begin{equation}\label{eq:hom}
f\big(\dg{N}{aN+b}{R}\big)=\dg{q}{f(a)q+f(b)}{S}\qquad\text{for all }a,b.
\end{equation}
A bijective homomorphism is an \emph{isomorphism} (then $q=N$); an isomorphism of a magma onto itself is an \emph{automorphism}.
For groups, $\ker f=f^{-1}(0)$.
\end{definition}

\begin{proposition}[Relabelling]\label{prop:relabel}
For $F\in\Sym_N$, the unique code $R^{F}$ making $F$ an isomorphism $(\Ap{N},R)\to(\Ap{N},R^F)$ is
\begin{equation}\label{eq:relabel}
R^{F}=\sum_{a,b=0}^{N-1}\dg{N}{\dg{N}{aN+b}{R}}{F}\;N^{\,\dg{N}{a}{F}\,N+\dg{N}{b}{F}}.
\end{equation}
This is an action, $(R^{F})^{G}=R^{G\circ F}$, so $\Aut(\Ap{N},R)=\{F\in\Sym_N:R^{F}=R\}$ is a subgroup of $\Sym_N$.
\end{proposition}
\begin{proof}
\eqref{eq:hom} prescribes the digit of $R^F$ at position $F(a)N+F(b)$; apply~\eqref{eq:code}. $G\circ F$ is an isomorphism
$R\to(R^F)^G$, and the code making it so is unique. A stabiliser of an action is a subgroup.
\end{proof}

\section{Counting magmas up to symmetry}\label{sec:magma:count}
The relabelling action of Proposition~\ref{prop:relabel} moves each digit of a code twice: position $aN+b$ is relocated to $F(a)N+F(b)$, and the digit value is renamed by $F$. Counting its orbits is therefore a problem of the type treated in Section~\ref{ssec:polya}. For $F\in\Sym_N$ and $k\ge1$ let $\operatorname{fix}_k(F)$ be the number of $i\in\Ap{N}$ with $F^{k}(i)=i$, and for $\alpha,\beta\in\Sym_N$ let $e_k(\alpha,\beta)$ be the number of cycles of length $k$ of the permutation $(a,b)\mapsto(\alpha(a),\beta(b))$ of $\Ap{N}^{2}$. A cycle of length $i$ of $\alpha$ and a cycle of length $j$ of $\beta$ give $\gcd(i,j)$ cycles of length $\operatorname{lcm}(i,j)$ on their product, and this determines $e_k$ from the cycle types of $\alpha$ and $\beta$.

\begin{theorem}[Magmas up to isomorphism]\label{thm:magma:count}
The number of isomorphism classes of magmas on $\Ap{N}$ is
\[
\frac1{N!}\sum_{F\in\Sym_N}\ \prod_{k\ge1}\operatorname{fix}_k(F)^{\,e_k(F,F)} .
\]
\end{theorem}
\begin{proof}
By Proposition~\ref{prop:relabel} and~\eqref{eq:relabel}, $F$ fixes the code $R$ if and only if $\dg{N}{F(a)N+F(b)}{R}=F\big(\dg{N}{aN+b}{R}\big)$ for all $a,b$. Along a cycle of length $k$ of $(a,b)\mapsto(F(a),F(b))$ the digit at one position determines all digits on the cycle, and the condition closes up if and only if the starting digit $y$ satisfies $F^{k}(y)=y$. So $F$ fixes $\prod_k\operatorname{fix}_k(F)^{e_k(F,F)}$ codes, and the claim is Burnside's lemma (Theorem~\ref{thm:dtburnside}).
\end{proof}

Two magmas $R,S$ on $\Ap{N}$ are \emph{isotopic} if $S(\alpha(a),\beta(b))=\gamma(R(a,b))$ for some $\alpha,\beta,\gamma\in\Sym_N$; isomorphism is the case $\alpha=\beta=\gamma$. The group $\Sym_N^{3}$ acts on the codes, and the argument of the proof, with $(\alpha,\beta)$ acting on positions and $\gamma$ on values, is Theorem~\ref{thm:debruijn}.

\begin{theorem}[Magmas up to isotopy]\label{thm:magma:isotopy}
The number of isotopy classes of magmas on $\Ap{N}$ is
\[
\frac1{(N!)^{3}}\sum_{\alpha,\beta,\gamma\in\Sym_N}\ \prod_{k\ge1}\operatorname{fix}_k(\gamma)^{\,e_k(\alpha,\beta)} .
\]
\end{theorem}

\begin{proposition}[Magmas up to opposite]\label{prop:magma:opp}
The number of operators on $\Ap{N}$ up to the exchange $R(a,b)\leftrightarrow R(b,a)$, that is, of magmas up to passing to the opposite magma, is $\tfrac12\big(N^{N^{2}}+N^{N(N+1)/2}\big)$.
\end{proposition}
\begin{proof}
The exchange of the two coordinates is an involution of $\Ap{N}^{2}$ with $N$ fixed points and $N(N-1)/2$ two-cycles, hence $N(N+1)/2$ cycles. Apply Theorem~\ref{thm:polya} to $\mathbb Z_2$ acting on the $N^2$ positions with $N$ colours.
\end{proof}

\begin{table}[ht]
\centering
\small
\begin{tabular}{rrrrr}
$N$ & all $N^{N^{2}}$ & up to isomorphism & up to isotopy & up to opposite\\\hline
$1$ & $1$ & $1$ & $1$ & $1$\\
$2$ & $16$ & $10$ & $5$ & $12$\\
$3$ & $19683$ & $3330$ & $139$ & $10206$\\
$4$ & $4294967296$ & $178981952$ & $332034$ & $2147516416$
\end{tabular}
\caption{Magmas of small order up to symmetry. Verified by \texttt{scripts/polya\_verify.py} by formula and, for $N\le3$, by direct orbit enumeration.}\label{tab:magma:count}
\end{table}

\begin{example}
For $N=3$ the identity fixes all $3^{9}$ codes. A transposition has one fixed point and acts on $\Ap3^{2}$ with one fixed pair and four $2$-cycles, so it fixes $\operatorname{fix}_1^{1}\operatorname{fix}_2^{4}=1\cdot3^{4}=81$ codes; a $3$-cycle has three $3$-cycles on pairs and fixes $\operatorname{fix}_3^{3}=27$ codes. Hence $\frac16\big(19683+3\cdot81+2\cdot27\big)=3330$.
\end{example}

\section{Direct products}
\begin{definition}
Given magmas $*_h$ on $\Ap{p_h}$, $h<\ell$, their \emph{direct product} is the magma on $\Ap{\Pi_\ell}$
\begin{equation}\label{eq:dprod}
m*n=\sum_{h=0}^{\ell-1}\big(\md{h}{m}*_h\md{h}{n}\big)\Pi_h .
\end{equation}
If all factors equal $*$ on $\Ap{N}$, this is $m*n=\sum_h\big(\dg{N}{h}{m}*\dg{N}{h}{n}\big)N^{h}$.
\end{definition}

\begin{proposition}\label{prop:dprod}
Each projection $D_h$ is a surjective homomorphism onto $*_h$. The product is associative, commutative, a quasigroup,
has identity $0$, or is a group if and only if every factor has the same property.
\end{proposition}
\begin{proof}
By Lemma~\ref{lem:mixed}, $\md{h}{m*n}=\md{h}{m}*_h\md{h}{n}$, so every identity holds in the product exactly when it holds digit by digit.
\end{proof}

\section{\texorpdfstring{$r$}{r}-ary operators}\label{sec:rary}
\begin{definition}
An \emph{$r$-ary operator} on $\Ap{N}$ is a map $\Ap{N}^r\to\Ap{N}$, written by Theorem~\ref{thm:coding} as
$y=\dg{N}{\sum_{k<r}x_kN^k}{R}$ with code $R\in[0,N^{N^r})$. The argument $x_k$ is the $k$-th digit of the index. Magmas are the case $r=2$ with $x_1=a$ and $x_0=b$, as in~\eqref{eq:magma}, so that the left operand is the most significant digit; permutations and other unary maps are the case $r=1$.
\end{definition}

\begin{proposition}\label{prop:rarycomp}
The $j$-th projection $(x_0,\dots,x_{r-1})\mapsto x_j$ has code $\sum_{n<N^r}\dg{N}{j}{n}\,N^{n}$. If $R$ is $r$-ary and
$S_0,\dots,S_{r-1}$ are $s$-ary, then the composite
$x\mapsto R(S_0(x),\dots,S_{r-1}(x))$ has code
$\sum_{n<N^{s}}\dg{N}{\sum_{k<r}\dg{N}{n}{S_k}N^{k}}{R}\,N^{n}$.
\end{proposition}
\begin{proof}Both are~\eqref{eq:code} applied to the map in question, with $n=\sum_k x_kN^k$ and $x_j=\dg{N}{j}{n}$.\end{proof}

\subsection{Clones and identities}
\begin{definition}\label{def:clone}
A \emph{clone} on $\Ap{N}$ is a family of finitary operators on $\Ap{N}$ that contains all projections and is closed under the
composition of Proposition~\ref{prop:rarycomp}. The \emph{clone} $\mathrm{Clo}(R)$ of a magma $(\Ap{N},R)$ is the smallest clone
containing $*$; its members are the \emph{term operations} of the magma, and $\mathrm{Clo}_s(R)$ denotes its $s$-ary part.
\end{definition}
Since $a*b$ has $x_1=a$ and $x_0=b$, the term operation $S*T$ of two $s$-ary term operations has code $R\circ(T,S)$ in the notation of
Proposition~\ref{prop:rarycomp}. Thus $\mathrm{Clo}_s(R)$ is the closure of the $s$ projection codes under one explicit integer map.

\begin{proposition}\label{prop:identities}
Let $(\Ap{N},R)$ be a magma. An identity $u=v$ between two terms in $s$ variables holds if and only if the codes of the term
operations of $u$ and $v$ coincide. The set $\mathrm{Clo}_s(R)\subseteq[0,N^{N^s})$ is finite and is computed by finitely many
compositions.
\end{proposition}
\begin{proof}
By Theorem~\ref{thm:coding} two $s$-ary operators are equal if and only if their codes are equal. There are at most $N^{N^s}$ codes,
so the closure stabilises after finitely many steps.
\end{proof}
Associativity is the case $s=3$, $u=(x*y)*z$, $v=x*(y*z)$; the two codes are the ternary codes $L(R)$ and $R'(R)$
of~\eqref{dg:TLTR}, and~\eqref{dg:defect} counts the triples where they differ.

\begin{remark}[Varieties and pseudovarieties]
Abstract clones are equivalent to Lawvere theories~\cite{Lawvere1963}. By Birkhoff's theorem, a class of algebras is defined by identities
if and only if it is closed under homomorphic images, subalgebras and arbitrary direct products~\cite{Birkhoff1935}; groups form such a
variety in the signature $(*,{}^{-1},e)$, but not with $*$ alone. A book on finite groups needs the finite version. A
\emph{pseudovariety} is a class of finite algebras closed under homomorphic images, subalgebras and finite direct
products~\cite{EilenbergSchutzenberger1976}, and Reiterman's theorem characterises pseudovarieties by pseudoidentities~\cite{Reiterman1982}.
Finite groups form a pseudovariety but not a variety. In the finite case the binary law suffices for subalgebras: a nonempty subset of a
finite group closed under $*$ is a subgroup, because $x^{-1}$ is a positive power of $x$. All three closure operations are already
digital: homomorphisms are the solutions of~\eqref{eq:hom}, subalgebras are restrictions of the table, and finite direct products
are Proposition~\ref{prop:dprod}.
\end{remark}

\subsection{The ternary law of a group}
In a group, put
\begin{equation}\label{eq:malcev}
m(x,y,z)=x*y^{-1}*z .
\end{equation}
Then $m(x,y,y)=x$ and $m(x,x,y)=y$, so $m$ is a \emph{Mal'cev term}; by Mal'cev's theorem a variety has such a term if and only if any
two congruences of each of its algebras permute~\cite{Malcev1954}. The term $m$ also satisfies
\begin{equation}\label{eq:heap}
m(x,y,y)=x=m(y,y,x),\qquad m\big(m(x,y,z),u,v\big)=m\big(x,y,m(z,u,v)\big),
\end{equation}
and a nonempty set with a ternary law satisfying~\eqref{eq:heap} is a \emph{heap}~\cite{Pruefer1924,Baer1929}.

\begin{proposition}\label{prop:heap}
Let $m$ be a heap law on $\Ap{N}$ and $e\in\Ap{N}$. Then $x*y=m(x,e,y)$ is a group law with identity $e$ and inverse
$y^{-1}=m(e,y,e)$, and $m(x,y,z)=x*y^{-1}*z$. Conversely, every group law gives the heap law~\eqref{eq:malcev}.
\end{proposition}
\begin{proof}
Associativity of $*$ is the second identity of~\eqref{eq:heap} with $y=u=e$. The first identities give $m(e,e,y)=y$ and $m(x,e,e)=x$.
For the inverse, $y*m(e,y,e)=m(y,e,m(e,y,e))=m(m(y,e,e),y,e)=m(y,y,e)=e$, and similarly on the left. Finally
$x*y^{-1}*z=m(x,e,m(m(e,y,e),e,z))=m(x,e,m(e,y,z))=m(m(x,e,e),y,z)=m(x,y,z)$. The converse is a direct check.
\end{proof}
The heap forgets only the choice of identity: the $N$ choices of $e$ give isomorphic groups, by the translation $x\mapsto m(x,e,e')$.
In digital form the code of $m$ is obtained from the code $R$ and the inversion code by two compositions in
Proposition~\ref{prop:rarycomp}.

\subsection{Cochains}
Let $Q$ be a finite group acting on a finite abelian group $A$, with $|Q|=q$ and $|A|=n$, labelled by $\Ap{q}$ and $\Ap{n}$ with identities
at $0$. An \emph{$r$-cochain} is a map $f:Q^r\to A$; the set of them is $C^r(Q,A)=A^{Q^r}$. The proof of Theorem~\ref{thm:coding}, with input
base $q$ and output base $n$, gives
\begin{equation}\label{eq:cochaincode}
f(x_0,\dots,x_{r-1})=\dg{n}{\textstyle\sum_{k<r}x_kq^k}{F},\qquad 0\le F<n^{q^r},
\end{equation}
so $|C^r(Q,A)|=n^{q^r}$, and the \emph{normalised} cochains, which vanish when some argument is $0$, number $n^{(q-1)^r}$. The
carry code $F$ of Chapter~\ref{ch:diagrammatic} is~\eqref{eq:cochaincode} for $r=2$. The low degrees have the following meaning
\cite{EilenbergMacLane1947}: $1$-cocycles are crossed homomorphisms $f(xy)=f(x)+x\cdot f(y)$, and $H^1(Q,A)$ classifies the complements of
$A$ in $A\rtimes Q$ up to conjugacy; $2$-cocycles are the factor sets of extension theory, and $H^2(Q,A)$ classifies the extensions with the
given action. The action itself is not a cochain: it is the separate datum carried by the action code $\Phi$. A first example is
Proposition~\ref{prop:Bcocycle}: the carry $c_N:\Ap{N}\times\Ap{N}\to\Z$ is a normalised $2$-cocycle of $\Z_N$ with values in $\Z$ (trivial
action), and the extension it defines is $0\to\Z\xrightarrow{\,\cdot N\,}\Z\to\Z_N\to0$, that is, positional notation itself.

\subsection{Polyadic groups}
For $r\ge3$, an \emph{$r$-ary group} is a set $G$ with an $r$-ary law $f$ such that the value
\[
f\big(x_1,\dots,x_i,f(x_{i+1},\dots,x_{i+r}),x_{i+r+1},\dots,x_{2r-1}\big)
\]
does not depend on $i\in\{0,\dots,r-1\}$, and every equation
$f(a_1,\dots,a_{i-1},x,a_{i+1},\dots,a_r)=b$ has a unique solution $x$~\cite{Dornte1929}. Post's coset theorem embeds $G$ in a group $G^*$
with a normal subgroup $G_0$ such that $G^*/G_0$ is cyclic of order $r-1$, $G$ is a generating coset of $G_0$ and $f$ is the product
of $r$ factors in $G^*$~\cite{Post1940}.

\begin{theorem}[Hossz\'u--Gluskin~\cite{Hosszu1963,Gluskin1965}]\label{thm:hosszu}
Let $(G,f)$ be an $r$-ary group, $r\ge3$. There are a group law on $G$, an automorphism $\theta$ and an element $b$ with $\theta(b)=b$
and $\theta^{r-1}(x)=bxb^{-1}$ for all $x$, such that
\[
f(x_1,\dots,x_r)=x_1\,\theta(x_2)\,\theta^2(x_3)\cdots\theta^{r-1}(x_r)\,b .
\]
Conversely, any such data define an $r$-ary group.
\end{theorem}
For example, in an abelian group with $r=3$, $\theta(x)=-x$ and $b=0$ one gets $f(x,y,z)=x-y+z$, which is the heap law~\eqref{eq:malcev}.
Digitally, Theorem~\ref{thm:hosszu} replaces a search through all $N^{N^r}$ codes by a search through group codes on $\Ap{N}$,
automorphisms and elements.

\section{Exercises}
\begin{enumerate}
\item Show that the identity permutation code of degree $N$ is $\iota_N=\sum_{i<N}i\,N^i=N\big(1+(N-2)N^N\big)/(N-1)^2$, and check $\iota_3=21$.
\item Show that $\Sym_3=\{5,7,11,15,19,21\}$, identify the three transpositions, and write the composition table of these codes
using~\eqref{eq:comp}.
\item Show that the sign of a permutation code $R$ is $(-1)^{\iota(R)}$, where $\iota(R)$ is the number of pairs $i<j$ with
$\dg{N}{i}{R}>\dg{N}{j}{R}$, and compute the sign of each element of $\Sym_3$ (compare Theorem~\ref{thm:sign}).
\item Show that there are $N^{N^2}$ magmas on $\Ap N$, of which $N^{N(N+1)/2}$ are commutative and $N^{(N-1)^2}$ have $0$ as a two-sided identity.
\item Using the relabelling action of Proposition~\ref{prop:relabel}, show that there are exactly $10$ magmas of order $2$ up to isomorphism.
\item Show that addition modulo $3$ is the only group law on $\Ap3$ with identity $0$.
\item Write the code of the direct product $C_2\oplus C_2$ given by~\eqref{eq:dprod} and show that it differs from the code of $C_4$
although both laws live on $\Ap4$.
\item Show that there are exactly $5$ magmas on $\Ap{2}$ up to isotopy, and give the code of one representative of each class. Explain why the projections $a*b=a$ and $a*b=b$ are not isotopic, although they are exchanged by passing to the opposite magma.
\item Use Theorem~\ref{thm:magma:count} to show that there are $3330$ magmas of order $3$ up to isomorphism.
\item Show that every group law gives a heap by~\eqref{eq:malcev}, that is, verify~\eqref{eq:heap}. For a heap law $m$ on $\Ap{N}$ and
$e,e'\in\Ap{N}$, show that $x\mapsto m(x,e,e')$ is an isomorphism from the group $x*_e y=m(x,e,y)$ of Proposition~\ref{prop:heap} onto
$x*_{e'}y=m(x,e',y)$.
\item Let $(\Ap{N},+)$ be an abelian group. Show that $f(x,y,z)=x-y+z$ and $g(x,y,z)=x+y+z+1$ are ternary groups, and find the
automorphism $\theta$ and the element $b$ of Theorem~\ref{thm:hosszu} for each.
\item Let $Q=A=\Z_2$ with trivial action. Show that there are exactly two normalised $2$-cochains, give their codes $F$
by~\eqref{eq:cochaincode}, show that both are cocycles and that the only normalised coboundary is $0$, so that $H^2(\Z_2,\Z_2)\cong\Z_2$.
Identify the two extensions.
\end{enumerate}

\makeatletter
\@ifundefined{definition}{\theoremstyle{definition}\newtheorem{definition}{Definition}[chapter]}{}
\@ifundefined{example}{\theoremstyle{definition}\newtheorem{example}{Example}[chapter]}{}
\@ifundefined{exercise}{\theoremstyle{definition}\newtheorem{exercise}{Exercise}[chapter]}{}
\@ifundefined{proposition}{\theoremstyle{plain}\newtheorem{proposition}{Proposition}[chapter]}{}
\@ifundefined{remark}{\theoremstyle{remark}\newtheorem{remark}{Remark}[chapter]}{}
\makeatother
\tikzset{dgsolid/.style={circle,fill=black,inner sep=2.6pt},
         dghollow/.style={circle,draw,thick,fill=white,inner sep=2.6pt},
         dgarrow/.style={->,thick}}

\chapter{Diagrammatic construction of finite groups}\label{ch:diagrammatic}

The diagrammatic method for discrete operators~\cite{GM2018diag} exposes the layers into which an
operator ${}_r R_N$ on $N=p_0p_1\cdots p_{M-1}$ symbols decomposes, and represents the
couplings between layers by arrows. In this chapter we specialise the method to
groups. The aim is constructive: we want a single, uniform procedure that
produces finite groups of large order from a small amount of data, together with
an explicit translation of every diagram into digital mathematics, i.e.\ into a
formula built only from the digit function and integer codes. Questions of
classification and uniqueness are deliberately left aside.

The main message is the following. For a group, every node of a diagram is a
\emph{code} of a smaller group, and every arrow carries exactly \emph{two} codes: an
\emph{action code} $\Phi$ and a \emph{carry code} $F$. Once these codes are given, the code
$R$ of the whole group follows mechanically, even when $R$ itself is far too large
to be written down.

\section{Digital labelling and layers}

Throughout, a law of composition on $A_N=\{0,1,\dots,N-1\}$ is written, as in the
chapter on magmas,
\begin{equation}
a\ast b \;=\; d_N\!\left(aN+b,\;R\right), \qquad
d_N(k,x)=\left\lfloor \frac{x}{N^k}\right\rfloor - N\left\lfloor \frac{x}{N^{k+1}}\right\rfloor ,
\label{dg:law}
\end{equation}
so that the left operand $a$ is the most significant digit of the index $aN+b$ and
$R=\sum_{m=0}^{N^2-1}N^m\,a_m$ with $a_m=d_N(1,m)\ast d_N(0,m)$.

Let $N=p_0p_1\cdots p_{M-1}$. Every $x\in A_N$ has the mixed-radix (digital) expansion
\begin{equation}
x=\sum_{h=0}^{M-1}x^{(h)}\Pi_h ,\qquad
x^{(h)} = d_{p_h}\!\left(0,\;\frac{x}{\Pi_h}\right)\in A_{p_h},\qquad \Pi_h=\prod_{m=0}^{h-1}p_m ,
\label{dg:layers}
\end{equation}
where $x^{(h)}$ is the \emph{layer projector} of $x$ on layer $h$. Layer $h=0$ is the
ground layer, layers $h>0$ are lifted layers. We shall also use the \emph{upper block}
of $x$ above layer $h$,
\begin{equation}
U_h(x)=\left\lfloor \frac{x}{\Pi_{h+1}}\right\rfloor \in A_{Q_h},\qquad
Q_h=\prod_{m=h+1}^{M-1}p_m ,
\label{dg:upper}
\end{equation}
which collects all digits of $x$ on the layers $h+1,\dots,M-1$.

\begin{remark}
For groups, a layer (or a set of layers) is drawn as a solid node precisely when
its projector is a homomorphism onto the group of that layer. The kernel of such
a projector is a normal subgroup, and the solid node is the corresponding
quotient. Hollow nodes are the layers whose product is twisted by the layers
from which they receive arrows.
\end{remark}

\subsection{Graded operators}

The simplest diagrams have no arrows at all. They are characterised by the following notion, which applies to
operators of any arity. Recall $\Pi_h=\prod_{m<h}p_m$ from Eq.~(\ref{dg:layers}).

\begin{definition}\label{def:dggraded}
An $r$-ary operator $y=R(x_0,\dots,x_{r-1})$ on $A_N$, $N=p_0\cdots p_{M-1}$, is \emph{graded} (for this
ordering of the layers) if for every $h$ the output digit $y^{(h)}$ depends only on the input digits
$x_0^{(h)},\dots,x_{r-1}^{(h)}$ of the same layer, that is, if there are $r$-ary operators $A^{(h)}$ on $A_{p_h}$,
the \emph{layer operators}, with
\[
y^{(h)}=A^{(h)}\big(x_0^{(h)},\dots,x_{r-1}^{(h)}\big),\qquad h=0,\dots,M-1 .
\]
Otherwise the layers are \emph{coupled}.
\end{definition}

\begin{proposition}\label{prop:dggraded}
An $r$-ary operator $R$ on $A_N$ is graded if and only if it is the direct product of $r$-ary operators
$A^{(h)}$ on $A_{p_h}$,
\[
R(x_0,\dots,x_{r-1})=\sum_{h=0}^{M-1}A^{(h)}\big(x_0^{(h)},\dots,x_{r-1}^{(h)}\big)\,\Pi_h .
\]
The layer operators are then determined by $R$: their codes are
\[
A^{(h)}=\sum_{n=0}^{p_h^r-1}p_h^{\,n}\;d_{p_h}\!\left(0,\ \frac{R\big(d_{p_h}(0,n)\,\Pi_h,\dots,d_{p_h}(r-1,n)\,\Pi_h\big)}{\Pi_h}\right).
\]
For $r=2$ and groups this is the direct product of Proposition~\ref{prop:dprod}.
\end{proposition}
\begin{proof}
If $R$ is graded, expanding the output as $y=\sum_h y^{(h)}\Pi_h$ gives the displayed formula. Conversely, the
formula says precisely that the $h$-th output digit is $A^{(h)}$ applied to the $h$-th input digits. To read
off $A^{(h)}$, evaluate $R$ on inputs whose only nonzero digits lie on layer $h$, namely
$x_k=d_{p_h}(k,n)\Pi_h$, and take the layer projector of the output; by the coding theorem this gives the
$n$-th digit of the code of $A^{(h)}$.
\end{proof}

In the language of diagrams, a graded law is a diagram without arrows: carries and actions are exactly what
couple the layers. The carry-free sum $\oplus_p$ of Section~\ref{sec:bfunction} is graded with $M$ layers equal to
$p$, while ordinary addition modulo $p^M$ is coupled through its carries (Proposition~\ref{prop:Bcarry}).

\begin{remark}[Layer orderings]
The $M$ factors of $N$ can be listed in $M!$ orders (counting equal factors as distinct). Each ordering gives a
different mixed-radix reading of the labels in $A_N$, and the $M!$ readings are relabellings of the same
magma. Whether a given code is graded depends on the ordering, so a direct product may be visible as such
only in some of them. For example, the direct sum $C_2\oplus C_3$ with layers $p_0=2$, $p_1=3$ has the code of
$m\ast n=d_2(0,m+n)+2\,d_3\big(0,\lfloor m/2\rfloor+\lfloor n/2\rfloor\big)$, which is graded for this ordering;
the same code is not graded when $A_6$ is read with the layers in the order $p_0=3$, $p_1=2$ (checked
exhaustively).
\end{remark}

\section{The building block}

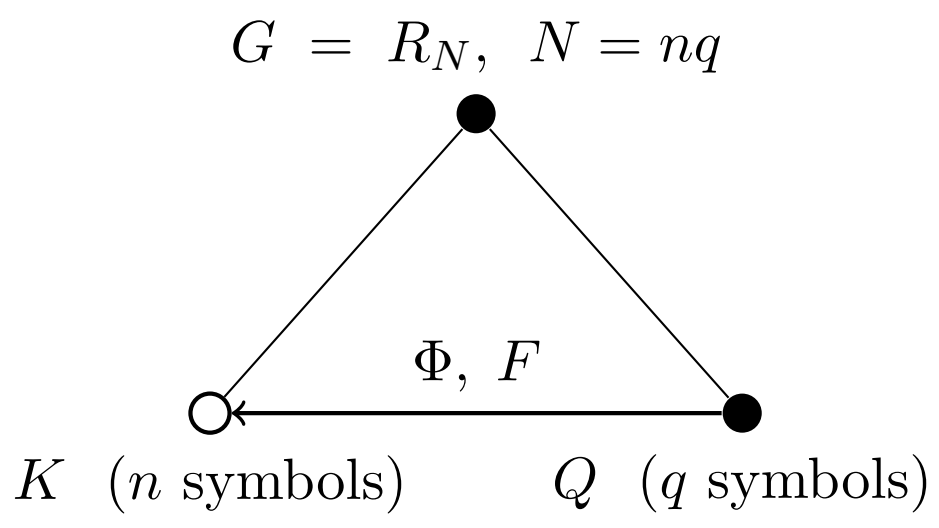
\begin{figure}[ht]
\centering
\begin{tikzpicture}[scale=1.1]
\node[dgsolid,label=above:{$G\;=\;R_{N}$, \ $N=nq$}] (G) at (0,1.8) {};
\node[dghollow,label=below:{$K$ \ ($n$ symbols)}] (K) at (-1.6,0) {};
\node[dgsolid,label=below:{$Q$ \ ($q$ symbols)}] (Q) at (1.6,0) {};
\draw (G)--(K); \draw (G)--(Q);
\draw[dgarrow] (Q) -- node[above]{$\Phi,\;F$} (K);
\end{tikzpicture}
\caption{The building block. The lifted layer $Q$ is a layer operator (solid node);
the ground layer $K$ (hollow node) receives one arrow carrying the action code
$\Phi$ and the carry code $F$.}
\label{fig:dgblock}
\end{figure}

The whole construction rests on one elementary diagram (Fig.~\ref{fig:dgblock}): a
model of order $N=nq$ with a hollow ground layer $K$ on $n$ symbols, a solid lifted
layer $Q$ on $q$ symbols, and one arrow from $Q$ to $K$.

\begin{definition}[Data of a diagram]\label{def:dgcodes}
The building block is specified by four integers:
\begin{center}
\begin{tabular}{|l|c|c|l|}
\hline
Element & Code & Range & Meaning \\ \hline
node $Q$ (solid) & $Q$ & $[0,q^{q^2})$ & $A\ast_Q B=d_q(B+qA,\,Q)$ \\
node $K$ (hollow) & $K$ & $[0,n^{n^2})$ & $u\ast_K v=d_n(v+nu,\,K)$ \\
action arrow & $\Phi$ & $[0,n^{nq})$ & $\varphi_A(v)=d_n(v+nA,\,\Phi)$ \\
carry arrow & $F$ & $[0,n^{q^2})$ & $f(A,B)=d_n(B+qA,\,F)$ \\ \hline
\end{tabular}
\end{center}
Thus $\Phi$ is the code of a map $A_q\times A_n\to A_n$ (for each $A\in A_q$ a unary
operator ${}_1R_n$ acting on the ground layer), and $F$ is the code of a map
$A_q\times A_q\to A_n$.
\end{definition}

For $a,b\in A_N$ we write $a^{(0)}=d_n(0,a)$, $a^{(1)}=d_q(0,a/n)$ and similarly for $b$.

\paragraph{What the lifted layer does.} Being a layer operator, it multiplies its own
digits with its own table:
\begin{equation}
y^{(1)} = d_q\!\left(b^{(1)}+q\,a^{(1)},\;Q\right).
\label{dg:y1}
\end{equation}

\paragraph{What the ground layer does.} It multiplies its own digits, twisted by the
information arriving through the arrow:
\begin{equation}
\boxed{\;
y^{(0)} = d_n\!\Big( d_n\!\big(b^{(1)}+q\,a^{(1)},\,F\big)
\;+\; n\, d_n\!\big( d_n\!\big(b^{(0)}+n\,a^{(1)},\,\Phi\big) + n\,a^{(0)},\; K\big),\; K\Big)\;}
\label{dg:y0}
\end{equation}
Read from the inside out: (i) the action arrow relabels the right operand,
$b^{(0)}\mapsto\varphi_{a^{(1)}}(b^{(0)})$; (ii) the ground layer multiplies,
$a^{(0)}\ast_K\varphi_{a^{(1)}}(b^{(0)})$; (iii) the carry arrow inserts
$f(a^{(1)},b^{(1)})$.

\paragraph{Reassembly.} By the decomposition property of the digit function,
\begin{equation}
a\ast b = y^{(0)}+n\,y^{(1)} .
\label{dg:reassembly}
\end{equation}
In ordinary notation, writing $a=(a^{(0)},a^{(1)})$, Eqs.~(\ref{dg:y1})--(\ref{dg:reassembly}) state
\begin{equation}
\big(a^{(0)},a^{(1)}\big)\ast\big(b^{(0)},b^{(1)}\big)
=\Big(a^{(0)}\ast_K\varphi_{a^{(1)}}\big(b^{(0)}\big)\ast_K f\big(a^{(1)},b^{(1)}\big),\;
a^{(1)}\ast_Q b^{(1)}\Big),
\label{dg:schreier}
\end{equation}
which is the general extension of $K$ by $Q$.

\paragraph{Link with the general arrow formula.} In the diagrammatic method a
hollow node is described by one code $A_{n'}$ for each configuration $n'$ of the
source layer. Here $n'=b^{(1)}+q\,a^{(1)}\in A_{q^2}$ and
\begin{equation}
y^{(0)}=\sum_{n'=0}^{q^2-1} d_n\!\big(b^{(0)}+n\,a^{(0)},\,A_{n'}\big)\;
d_q\!\big(n'-b^{(1)}-q\,a^{(1)},\,1\big),
\label{dg:arrowsum}
\end{equation}
where $d_q(k,1)=\delta_{k,0}$ for integer $k$, and
\begin{equation}
A_{n'}=\sum_{m=0}^{n^2-1} n^m\,
d_n\!\Big( d_n(n',F) + n\,d_n\!\big( d_n\big(d_n(0,m)+n\,d_q(1,n'),\Phi\big) + n\,d_n(1,m),\,K\big),\,K\Big).
\label{dg:Anprime}
\end{equation}
An arbitrary arrow may carry any family of codes $A_{n'}$; the arrows that produce
groups are those for which the whole family is generated by the three codes $K$,
$\Phi$, $F$ as in Eq.~(\ref{dg:Anprime}).

\paragraph{Code of the group.} Finally, the code of the model is
\begin{equation}
R=\sum_{m=0}^{N^2-1}N^m\Big[\,y^{(0)}(a,b)+n\,y^{(1)}(a,b)\Big]_{a=d_N(1,m),\;b=d_N(0,m)} .
\label{dg:R}
\end{equation}

\section{Admissible codes}

The four codes define a group when the following conditions hold.
\begin{enumerate}
\item \emph{Normalisation of the action:} $d_n(v,\Phi)=v$ for $0\le v<n$ (i.e.\ $\varphi_0=\mathrm{id}$).
\item \emph{Normalisation of the carry:} $d_n(m,F)=0$ whenever $m<q$ or $d_q(0,m)=0$
(i.e.\ $f(0,B)=f(A,0)=0$). Then $0$ is the neutral element.
\item \emph{Each $\varphi_A$ is an automorphism of $K$:}
\[
d_n\!\big(d_n(v+nu,K)+nA,\;\Phi\big)=d_n\!\big(d_n(v+nA,\Phi)+n\,d_n(u+nA,\Phi),\;K\big).
\]
\item \emph{Compatibility of action and carry:}
$\varphi_A\big(\varphi_B(v)\big)\ast_K f(A,B)=f(A,B)\ast_K\varphi_{A\ast_Q B}(v)$.
\item \emph{Cocycle condition:}
$f(A,B)\ast_K f(A\ast_Q B,C)=\varphi_A\big(f(B,C)\big)\ast_K f(A,B\ast_Q C)$.
\end{enumerate}
Conditions 4 and 5 become digital identities on substituting
$\varphi_A(v)=d_n(v+nA,\Phi)$, $f(A,B)=d_n(B+qA,F)$, $u\ast_K v=d_n(v+nu,K)$ and
$A\ast_QB=d_q(B+qA,Q)$. They are exactly what is obtained by expanding
$(a\ast b)\ast c=a\ast(b\ast c)$ with Eq.~(\ref{dg:schreier}), first with $b^{(0)}=c^{(0)}=0$ and then
with $c^{(1)}=0$.

\begin{remark}
For construction purposes it is often simplest to impose 1 and 2 and to test
associativity of the resulting table directly. Different admissible codes may
produce the same group: changing the labels of $K$ by a map depending on the
$Q$-digit changes $F$ without changing the group.
\end{remark}

\section{Kinds of arrows}

Let $\Phi_{\mathrm{id}}=\sum_{A=0}^{q-1}\sum_{v=0}^{n-1}v\,n^{v+nA}$ denote the trivial action code.
\begin{center}
\begin{tabular}{|l|l|l|}
\hline
Arrow & Codes & Resulting construction \\ \hline
no arrow & $\Phi=\Phi_{\mathrm{id}}$, $F=0$ & direct product $K\times Q$ (graded operator) \\
action only & $\Phi\ne\Phi_{\mathrm{id}}$, $F=0$ & semidirect product $K\rtimes Q$ \\
carry only & $\Phi=\Phi_{\mathrm{id}}$, $F\ne0$, $K$ abelian & central extension \\
action and carry & $\Phi\ne\Phi_{\mathrm{id}}$, $F\ne0$ & general extension \\ \hline
\end{tabular}
\end{center}
When $K=\mathbb{Z}_n$ and $\varphi_A(v)=\lambda_A v$ is multiplication by a unit, Eq.~(\ref{dg:y0})
simplifies to
\begin{equation}
y^{(0)}=d_n\!\Big(0,\;a^{(0)}+\lambda_{a^{(1)}}\,b^{(0)}+d_n\big(b^{(1)}+q\,a^{(1)},F\big)\Big),
\label{dg:metacyclic}
\end{equation}
since $d_n(0,x)$ is the residue of $x$ modulo $n$.

\section{Stacking layers: rank trees}

Larger groups are obtained by using an already built block as the solid node $Q$
of a new block. For $N=p_0p_1\cdots p_{M-1}$ the top layer $h=M-1$ is solid, and each
layer $h<M-1$ is a hollow node with code $K_h$ receiving arrows from the upper block
$U_h$, with an action code $\Phi_h\in[0,p_h^{\,p_hQ_h})$ and a carry code
$F_h\in[0,p_h^{\,Q_h^2})$. The layer outputs are
\begin{equation}
\boxed{\;
y^{(h)}=d_{p_h}\!\Big( d_{p_h}\!\big(U_h(b)+Q_hU_h(a),F_h\big)+p_h\,
d_{p_h}\!\big( d_{p_h}\big(b^{(h)}+p_hU_h(a),\Phi_h\big)+p_h\,a^{(h)},\,K_h\big),\,K_h\Big)\;}
\label{dg:stack}
\end{equation}
and
\begin{equation}
a\ast b=\sum_{h=0}^{M-1}y^{(h)}\Pi_h .
\label{dg:stacksum}
\end{equation}
For the top layer $Q_{M-1}=1$, $F_{M-1}=0$ and $\Phi_{M-1}$ is trivial, so that
Eq.~(\ref{dg:stack}) reduces to the layer operator $y^{(M-1)}=d_{p_{M-1}}(b^{(M-1)}+p_{M-1}a^{(M-1)},K_{M-1})$.
The evaluation runs from the top down: $y^{(h)}$ only needs the digits of the layers
$h'\ge h$. An arrow is drawn from layer $h'>h$ into layer $h$ if and only if $\Phi_h$ or
$F_h$ actually depends on the digit
$d_{p_{h'}}\big(0,\,U_h/\prod_{m=h+1}^{h'-1}p_m\big)$ of the upper block.
Any node may itself be a composite block (for instance a non-abelian simple
group), which is then treated as a single node with its own code.

\begin{figure}[ht]
\centering
\begin{tikzpicture}[scale=1.1]
\node[dgsolid,label=above:{$S_4=R_{24}$}] (G) at (0,1.9) {};
\node[dghollow,label=below:{$\mathbb{Z}_2^2$}] (V) at (-2.2,0) {};
\node[dghollow,label=below:{$\mathbb{Z}_3$}] (T) at (0,0) {};
\node[dgsolid,label=below:{$\mathbb{Z}_2$}] (S) at (2.2,0) {};
\draw (G)--(V); \draw (G)--(T); \draw (G)--(S);
\draw[dgarrow] (S) -- node[above]{$\Phi_1$} (T);
\draw[dgarrow] (T) -- node[above]{$\Phi_0$} (V);
\draw[dgarrow] (S) to[bend left=35] node[below]{$\Phi_0$} (V);
\end{tikzpicture}
\caption{A three-layer diagram for $S_4$ (Example~\ref{ex:S4}). Only action arrows occur.}
\label{fig:dgS4}
\end{figure}
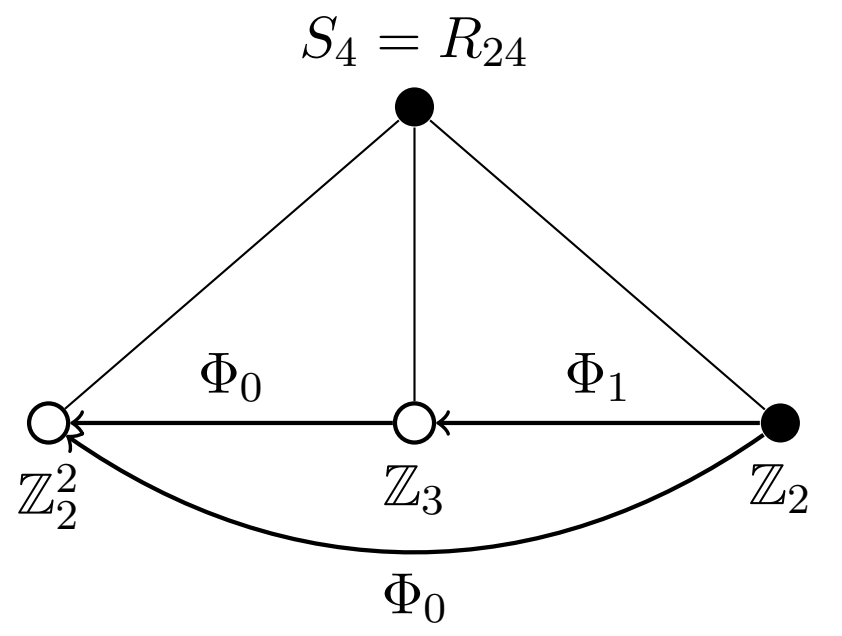

\section{Worked examples}

\begin{example}[Cyclic groups $\mathbb{Z}_{p^2}$: carry arrows]
Take $n=q=p$, $K=Q=$ the code of $\mathbb{Z}_p$, $\Phi=\Phi_{\mathrm{id}}$ and
$d_p(B+pA,F)=\lfloor (A+B)/p\rfloor$. Then
\[
y^{(1)}=d_p\big(0,a^{(1)}+b^{(1)}\big),\qquad
y^{(0)}=d_p\!\Big(0,\;a^{(0)}+b^{(0)}+\Big\lfloor\frac{a^{(1)}+b^{(1)}}{p}\Big\rfloor\Big),
\]
and $x\mapsto x^{(1)}+p\,x^{(0)}$ is an isomorphism onto $\mathbb{Z}_{p^2}$: the carry arrow is
nothing but the carry of base-$p$ addition. Iterating with three layers gives
$\mathbb{Z}_{27}$, with a code $R$ of $1044$ decimal digits.
\end{example}

\begin{example}[Dihedral and metacyclic groups: action arrows]
For $D_n$ take $K=\mathbb{Z}_n$, $Q=\mathbb{Z}_2$, $F=0$ and
$\Phi=\sum_{v=0}^{n-1}v\,n^v+\sum_{v=0}^{n-1}d_n(0,n-v)\,n^{n+v}$, i.e.\ $\varphi_1(v)=-v$. By
Eq.~(\ref{dg:metacyclic}), $y^{(0)}=d_n\big(0,a^{(0)}+(-1)^{a^{(1)}}b^{(0)}\big)$. With $K=\mathbb{Z}_7$,
$Q=\mathbb{Z}_3$ and $d_7(v+7A,\Phi)=d_7(0,2^Av)$ one obtains the Frobenius group
$\mathbb{Z}_7\rtimes\mathbb{Z}_3$ of order $21$.
\end{example}

\begin{example}[The quaternion group: action and carry]
Take $n=4$, $q=2$, $K=$ the code of $\mathbb{Z}_4$, $Q=6$ (the code of $\mathbb{Z}_2$), and
\[
\Phi=27876,\qquad F=128 .
\]
The base-$4$ digits of $\Phi$ are $(0,1,2,3,\,0,3,2,1)$ (identity, then $v\mapsto-v$) and the
only non-zero digit of $F$ is $d_4(3,F)=f(1,1)=2$. Equations (\ref{dg:y1}) and (\ref{dg:metacyclic}) give
\[
y^{(1)}=d_2\big(0,a^{(1)}+b^{(1)}\big),\qquad
y^{(0)}=d_4\big(0,\;a^{(0)}+(-1)^{a^{(1)}}b^{(0)}+2\,a^{(1)}b^{(1)}\big),
\]
and $a\ast b=y^{(0)}+4y^{(1)}$ is $Q_8$ (one element of order $2$, six of order $4$).
\end{example}

\begin{example}[$A_4$: a composite node]
Label the Klein group $\mathbb{Z}_2^2$ by $v=v_0+2v_1$, so that its code $K$ is that of
$v\ast w=d_2(0,v_0+w_0)+2\,d_2(0,v_1+w_1)$. The matrix $M(v_0,v_1)=(v_1,v_0+v_1)$ acts as the
permutation $(0,1,2,3)\mapsto(0,2,3,1)$. With $Q=\mathbb{Z}_3$, $F=0$ and $\Phi$ the code with
base-$4$ digits $(0,1,2,3,\;0,2,3,1,\;0,3,1,2)$ (i.e.\ $\varphi_A=M^A$) the block is $A_4$.
\end{example}

\begin{example}[$S_4$: three layers]\label{ex:S4}
First build $S_3$ from $K=\mathbb{Z}_3$, $Q=\mathbb{Z}_2$ with $\varphi_1(v)=-v$; an element
$u=u_0+3u_1$ of this block is $r^{u_0}s^{u_1}$. Then use it as the upper block of the Klein
layer, with $\varphi_u=M^{u_0}S^{u_1}$, where $S(v_0,v_1)=(v_1,v_0)$ acts as
$(0,1,2,3)\mapsto(0,2,1,3)$, and all carries zero. Eq.~(\ref{dg:stack}) with
$p_0=4$, $p_1=3$, $p_2=2$ gives $S_4$ (Fig.~\ref{fig:dgS4}); its code $R$ has $794$ decimal
digits, while the diagram is specified by a handful of small codes.
\end{example}

\begin{example}[$S_5$: a simple group as a node]
Let $K$ be the code of $A_5$ on $60$ symbols (for instance, from a labelling of the even
permutations of five letters), $Q=\mathbb{Z}_2$, $F=0$, and let $\varphi_1$ be conjugation by a
transposition. The block is $S_5$, of order $120$; its code $R$ has $29940$ decimal digits.
\end{example}

All the codes above were evaluated with Eqs.~(\ref{dg:y1})--(\ref{dg:R}) and (\ref{dg:stack})--(\ref{dg:stacksum}); in each case
the resulting $R$ reproduces, through Eq.~(\ref{dg:law}), an associative table with the expected
element orders:
\begin{center}
\begin{tabular}{|l|l|l|}
\hline
Group & Layers (top $\to$ ground) & Number of elements of order $1,2,3,\dots$ \\ \hline
$Q_8$ & $\mathbb{Z}_2\to\mathbb{Z}_4$ & ord.\ 1:1, 2:1, 4:6 \\
$\mathbb{Z}_{27}$ & $\mathbb{Z}_3\to\mathbb{Z}_3\to\mathbb{Z}_3$ & ord.\ 1:1, 3:2, 9:6, 27:18 \\
$S_4$ & $\mathbb{Z}_2\to\mathbb{Z}_3\to\mathbb{Z}_2^2$ & ord.\ 1:1, 2:9, 3:8, 4:6 \\
$S_5$ & $\mathbb{Z}_2\to A_5$ & ord.\ 1:1, 2:25, 3:20, 4:30, 5:24, 6:20 \\ \hline
\end{tabular}
\end{center}

\begin{example}[All diagrams of three $\mathbb{Z}_2$ layers]
Take $p_0=p_1=p_2=2$ with only carry arrows (the only automorphism of $\mathbb{Z}_2$ is the
identity). After normalisation, the carry into layer $1$ has one free digit and the
carry into layer $0$ has nine, giving $1024$ diagrams. Exactly $24$ of them are associative,
and they produce all groups of order $8$: $\mathbb{Z}_2^3$ (2 diagrams),
$\mathbb{Z}_4\times\mathbb{Z}_2$ (10), $\mathbb{Z}_8$ (4), $D_4$ (6) and $Q_8$ (2).
\end{example}

\section{Mixed layers and simple groups}\label{sec:dgsimple}

In the diagrammatic method a layer is nothing but a projector on a digit of the
mixed-radix representation, with a radix that divides the order of the group, and
an arrow means that layers are mixed: an incoming arrow indicates that a layer
contains material of another layer. Every group of order $N=p_0p_1\cdots p_{M-1}$ can be
labelled in this way for every factorisation of $N$ and for every bijection between
its elements and $A_N$. The arrows depend on the labelling, and the task of the
method is to find labellings with as few arrows as possible. This point of view
yields a characterisation of simple groups by their diagrams.

\subsection{Arrows as mixing}

\begin{definition}[Arrow]\label{def:arrowmix}
Let $a\ast b=d_N(aN+b,R)$, let $N=p_0p_1\cdots p_{M-1}$, and let $x^{(h)}$ denote the layer
projectors (\ref{dg:layers}). For $j\neq h$ there is an \emph{arrow} $j\to h$ if the output layer
$y^{(h)}=(a\ast b)^{(h)}$, regarded as a function of the $2M$ digits $a^{(m)},b^{(m)}$, is not
constant in $a^{(j)}$ or not constant in $b^{(j)}$. A layer receiving no arrow is a solid node,
a layer receiving arrows is a hollow node. A set $S$ of layers is \emph{closed} if there is
no arrow $j\to h$ with $j\notin S$ and $h\in S$.
\end{definition}

For $S\subseteq\{0,\dots,M-1\}$ let $\pi_S$ be the projector that keeps the digits of the layers in
$S$, and let $p_S=\prod_{h\in S}p_h$; $\pi_S$ maps $A_N$ onto a set of $p_S$ digit strings, which we
label in mixed radix by $A_{p_S}$.

\begin{proposition}[Closed sets of layers are quotients]\label{prop:closedquotient}
If $S$ is closed, then $\pi_S(a\ast b)$ depends only on $\pi_S(a)$ and $\pi_S(b)$, and
\[
R_S=\sum_{m=0}^{p_S^2-1}p_S^{\,m}\;\pi_S\big(a\ast b\big),\qquad \pi_S(a)=d_{p_S}(1,m),\ \pi_S(b)=d_{p_S}(0,m),
\]
is the code of a group on $p_S$ symbols. The projector $\pi_S$ is a surjective homomorphism,
and its kernel $K_S$ is a normal subgroup of order $N/p_S$.
\end{proposition}
\begin{proof}
Changing the digits of $a$ and $b$ outside $S$ one at a time does not change any layer
of $S$, because no arrow enters $S$ from outside; hence $\pi_S(a\ast b)$ depends only on
$\pi_S(a),\pi_S(b)$ and $R_S$ is well defined. By construction $\pi_S(a\ast b)=\pi_S(a)\ast_S\pi_S(b)$,
so $\ast_S$ inherits associativity, neutral element and inverses, and $\pi_S$ is a surjective
homomorphism. Its kernel is normal and all its fibres have $N/p_S$ elements.
\end{proof}

Conversely, let $K$ be a normal subgroup and label $x=x^{(0)}+|K|\,x^{(1)}$, where $x^{(1)}$ is the
coset $xK$ (in any numbering) and $x^{(0)}$ the position of $x$ inside its coset. Since
$(ab)K=(aK)(bK)$, the top layer receives no arrow: it is a solid node with the code of
$G/K$.

\begin{lemma}\label{lem:strong}
A directed graph on $M\ge2$ vertices has a closed set of vertices other than $\emptyset$
and the whole set if and only if it is not strongly connected.
\end{lemma}
\begin{proof}
If the graph is not strongly connected there are vertices $u,v$ such that $v$ cannot
be reached from $u$. The set of vertices from which $v$ can be reached contains $v$,
not $u$, and no arrow enters it from outside, since otherwise its source could reach $v$. Conversely,
no arrow leads from the complement of a closed proper set into it, so its vertices
cannot be reached from the complement.
\end{proof}

\begin{theorem}[Simple groups]\label{thm:simple}
A finite group $G$ of composite order is simple if and only if, for every factorisation
$|G|=p_0\cdots p_{M-1}$ with $M\ge2$ and every labelling of its elements, the arrow graph is
strongly connected.
\end{theorem}
\begin{proof}
If some arrow graph is not strongly connected, Lemma~\ref{lem:strong} gives a closed set $S$
with $\emptyset\neq S\neq\{0,\dots,M-1\}$. Since all $p_h\ge2$, we have $1<p_S<|G|$, and
Proposition~\ref{prop:closedquotient} gives a normal subgroup $K_S$ of order $|G|/p_S$, neither
trivial nor equal to $G$. Conversely, a normal subgroup $K\neq1,G$ gives, with the
labelling described above and $p_0=|K|$, $p_1=|G/K|$, a top layer without incoming
arrows.
\end{proof}

Thus the layers of a simple group can never be separated: whatever the factorisation and
the labelling, the mixing between layers can be reduced but not removed. For groups of
prime order there is only one layer and the statement is empty.

\subsection{Subgroup nodes and tower labellings}

Subgroups are a weaker notion than closed sets of layers. If the elements whose digits
vanish on all layers outside a set $L$ are closed under the product, they form a
subgroup, but this only constrains the product when both inputs have those digits
equal to zero, and mixing is still allowed. Simple groups therefore have plenty of
subgroup nodes, while they never have closed sets of layers.

A chain of subgroups $1=H_0<H_1<\cdots<H_M=G$ gives a \emph{tower labelling}. Choose for each
$k$ a transversal of the left cosets of $H_k$ in $H_{k+1}$ whose first element is the
identity, and write every element uniquely as $g=t_{M-1}t_{M-2}\cdots t_0$ with $t_k$ in the $k$-th
transversal. The digit $x^{(k)}$ is the index of $t_k$, with radix $p_k=[H_{k+1}:H_k]$. In this
labelling $H_{h+1}$ is exactly the set of elements with $U_h(x)=0$, and $U_h(x)$ labels the left
coset $xH_{h+1}$.

\begin{proposition}\label{prop:towerright}
In a tower labelling, $U_h(a\ast b)$ depends on $b$ only through $U_h(b)$. Hence the right
operand sends arrows only downwards, and every upward arrow comes from the left operand.
The upper block above layer $h$ is closed if and only if $H_{h+1}$ is normal in $G$.
\end{proposition}
\begin{proof}
$U_h(a\ast b)$ labels $abH_{h+1}=a\,(bH_{h+1})$, which depends on $a$ and on the coset $bH_{h+1}$, i.e.\ on
$U_h(b)$. If $H_{h+1}$ is normal, $abH_{h+1}=(aH_{h+1})(bH_{h+1})$ and no arrow enters the upper block. If
the upper block is closed, Proposition~\ref{prop:closedquotient} shows that $H_{h+1}$, the kernel of
$x\mapsto U_h(x)$, is normal.
\end{proof}

In a tower labelling, the upward arrows are given by the action of $G$ on the cosets of
the subgroups of the chain; changing the transversals only relabels the cosets and
cannot remove them unless the subgroups are normal.

\subsection{Examples}

The arrows of the following labellings were computed exactly from the full
multiplication tables, for $S_4$ and $A_5$ (radices listed from the ground layer up):
\begin{center}
\begin{tabular}{|p{0.30\textwidth}|p{0.56\textwidth}|}
\hline
Labelling & Arrow graph \\ \hline
$S_4$: tower $1<\mathbb{Z}_2<V<A_4<S_4$, radices $2,2,3,2$ &
$0\leftrightarrow1$, $2\to0,1$, $3\to0,1,2$; not strongly connected \\
$A_5$: tower $1<\mathbb{Z}_2<V<A_4<A_5$, radices $2,2,3,5$ &
all $j\to h$ with $j\neq h$; strongly connected \\
$S_4$: 285 random labellings over all 19 factorisations & all strongly connected \\
$A_5$: 645 random labellings over all 43 factorisations & all strongly connected \\ \hline
\end{tabular}
\end{center}
In the tower labelling of $S_4$, layer $3$ is closed (the sign homomorphism onto $\mathbb{Z}_2$) and
layers $\{2,3\}$ are closed (the quotient $S_4/V\cong S_3$), in agreement with
Proposition~\ref{prop:closedquotient}. In $A_5$ every layer contains material of every other layer.
The random labellings show that strong connectivity is the generic situation, also
for groups that are not simple: Theorem~\ref{thm:simple} states that for simple groups no
labelling escapes it.

\begin{example}[The tower of $M_{11}$]
The chain $1<\mathbb{Z}_2<\mathbb{Z}_4<Q_8<M_9<M_{10}<M_{11}$ of point stabilisers (with $Q_8$ refined by its
cyclic subgroups) gives the radices $2,2,2,9,10,11$. Inside the $Q_8$ node the two upper
layers receive no arrows and the ground layer receives a carry, which is the diagram of
the quaternion group obtained in the worked examples of this chapter. In a sample of 3000 random pairs, every layer of
$M_9$, $M_{10}$ and $M_{11}$ above the $Q_8$ node receives arrows from all lower layers, through the
left operand, as required by Proposition~\ref{prop:towerright} and Theorem~\ref{thm:simple}.
\end{example}

\section{A gallery of recipes}

\begin{itemize}
\item \emph{Generalised quaternion groups} $Q_{2^k}$: $K=\mathbb{Z}_{2^{k-1}}$, $Q=\mathbb{Z}_2$,
$\varphi_1(v)=-v$, $f(1,1)=2^{k-2}$.
\item \emph{Heisenberg groups} of order $p^3$: $Q=\mathbb{Z}_p\times\mathbb{Z}_p$ (a graded block, no
internal arrows), $K=\mathbb{Z}_p$, trivial action, $f\big((x_1,y_1),(x_2,y_2)\big)=d_p(0,x_1y_2)$.
\item \emph{Central extensions of a simple node}: $K=\mathbb{Z}_2$, $Q=A_5$, trivial action and a
non-trivial normalised cocycle $F$ give the double cover $SL(2,5)$ of order $120$.
\item \emph{Wreath products}: $K=A_5\times A_5$ (graded block on $3600$ symbols, $v=v_0+60v_1$),
$Q=\mathbb{Z}_2$, $F=0$ and $\varphi_1(v_0+60v_1)=v_1+60v_0$ give $A_5\wr\mathbb{Z}_2$, of order $7200$.
\item \emph{Simple groups}: by Theorem~\ref{thm:simple} their arrow graph is strongly connected
for every labelling, so no set of layers is ever closed. They are generated from permutation
codes (Proposition~\ref{prop:cayley}) or from matrix groups acting on layered vectors, e.g.\
$GL(3,2)\cong PSL(2,7)$ acting on $A_8\cong\mathbb{F}_2^3$, and their subgroup chains give tower
labellings with subgroup nodes and mixed layers (Section~\ref{sec:dgsimple}).
\end{itemize}

\section{Associativity by construction}\label{sec:dgassoc}

Associativity is the only axiom that involves three elements. Since ternary
operators can be built from binary operators acting twice, associativity can be
formulated as an identity between two ternary codes. This formulation suggests
two ways of \emph{constructing} associative codes instead of searching for them:
through the rows of the table, which are unary codes, and through the carry
arrows of a diagram, which are solutions of linear equations.

\subsection{Associativity as an identity between ternary codes}

Let $a\ast b=d_N(aN+b,R)$. Composing the binary operator with itself in the two
possible ways gives the ternary operators
\begin{equation}
T_L(a,b,c)=d_N\!\big(c+N\,d_N(b+Na,R),\,R\big),\qquad
T_R(a,b,c)=d_N\!\big(d_N(c+Nb,R)+Na,\,R\big),
\label{dg:TLTR}
\end{equation}
i.e.\ $T_L(a,b,c)=(a\ast b)\ast c$ and $T_R(a,b,c)=a\ast(b\ast c)$. Their codes, in the universal
normal form of ternary operators, are
\begin{align}
L(R)&=\sum_{m=0}^{N^3-1}N^m\,T_L\big(d_N(2,m),d_N(1,m),d_N(0,m)\big),\notag\\
R'(R)&=\sum_{m=0}^{N^3-1}N^m\,T_R\big(d_N(2,m),d_N(1,m),d_N(0,m)\big).
\label{dg:ternarycodes}
\end{align}
The operator with code $R$ is associative if and only if $L(R)=R'(R)$. The number of
positions in which the two codes differ,
\begin{equation}
D(R)=\sum_{m=0}^{N^3-1}\Big(1-d_N\big(d_N(m,L(R))-d_N(m,R'(R)),1\big)\Big),
\label{dg:defect}
\end{equation}
is the \emph{associativity defect}: $D(R)=0$ characterises associative operators, and
$D(R)$ can be used as an objective function when one explores magmas, quasigroups
or loops. As a search over all $N^{N^2}$ codes, however, this criterion is hopeless
except for very small $N$. The rest of the section shows how to avoid the search.

\subsection{Rows as unary codes}

The base-$N^N$ digits of $R$ are unary codes: since the index $b+Na$ has $a$ as its
most significant digit,
\begin{equation}
R=\sum_{a=0}^{N-1}\Lambda_a\,(N^N)^a,\qquad \Lambda_a=d_{N^N}(a,R),\qquad
a\ast b=d_N(b,\Lambda_a).
\label{dg:rows}
\end{equation}
$\Lambda_a$ is the code of left multiplication by $a$, i.e.\ of the row $a$ of the table.
The composition of two unary codes $\Lambda$, $\Lambda'$ is the unary code
\begin{equation}
\Lambda\circ\Lambda'=\sum_{v=0}^{N-1}N^v\,d_N\big(d_N(v,\Lambda'),\Lambda\big).
\label{dg:compose}
\end{equation}

\begin{proposition}\label{prop:dgrows}
The operator with code $R$ is associative if and only if
\[
\Lambda_{a\ast b}=\Lambda_a\circ\Lambda_b\qquad\text{for all } a,b\in A_N .
\]
\end{proposition}
\begin{proof}
By Eqs.~(\ref{dg:TLTR}) and (\ref{dg:rows}),
\[
T_L(a,b,c)=d_N(c,\Lambda_{a\ast b}),\qquad
T_R(a,b,c)=d_N\big(d_N(c,\Lambda_b),\Lambda_a\big)=d_N(c,\Lambda_a\circ\Lambda_b)
\]
by Eq.~(\ref{dg:compose}). The two ternary codes coincide for all $c$ exactly when the
two unary codes coincide.
\end{proof}

Read backwards, Proposition~\ref{prop:dgrows} is a construction, and associativity is
inherited from the composition of maps, which is always associative.

\begin{proposition}[Construction from unary codes]\label{prop:cayley}
Let $H$ be a set of $N$ unary codes on some alphabet $A_s$, containing the identity and
closed under composition (\ref{dg:compose}), and let $\phi:H\to A_N$ be a bijection with
$\phi(\mathrm{id})=0$. Then
\[
a\ast b=\phi\big(\phi^{-1}(a)\circ\phi^{-1}(b)\big)
\]
is associative with neutral element $0$, and its code is obtained from Eq.~(\ref{dg:rows}).
If the elements of $H$ are permutations, $\ast$ is a group law.
\end{proposition}

In practice $H$ is generated by a few permutation codes (for instance the generators
of a simple group), all products are formed by Eq.~(\ref{dg:compose}), and the elements are
labelled in the order in which they appear. This is how simple groups enter the
diagrams as single nodes.

\begin{proposition}[Generator test]\label{prop:light}
Let $X\subset A_N$ generate $A_N$ under $\ast$. If
$(x\ast g)\ast y=x\ast(g\ast y)$ for all $x,y\in A_N$ and all $g\in X$, then $\ast$ is associative.
Equivalently, it suffices to verify $\Lambda_{x\ast g}=\Lambda_x\circ\Lambda_g$ for $g\in X$.
\end{proposition}
\begin{proof}
If $g$ and $h$ satisfy the identity, then for all $x,y$
\[
x\ast\big((g\ast h)\ast y\big)=x\ast\big(g\ast(h\ast y)\big)=(x\ast g)\ast(h\ast y)
=\big((x\ast g)\ast h\big)\ast y=\big(x\ast(g\ast h)\big)\ast y ,
\]
so $g\ast h$ satisfies it as well. The set of elements satisfying the identity is thus
closed under $\ast$ and contains $X$, hence it is $A_N$.
\end{proof}
With $k$ generators the test needs about $kN^2$ evaluations of the digit function instead
of $N^3$.

\subsection{Carry arrows from linear algebra}

In a stacked diagram the associativity conditions of this chapter
separate layer by layer: the conditions for layer $h$ involve only its own codes and
the upper block, which has already been built and is associative. Consider the
building block with an abelian ground layer $K=\mathbb{Z}_n$, an upper block $Q$ and a fixed
action code $\Phi$ (each $\varphi_A$ an automorphism, i.e.\ multiplication by a unit). Then
the cocycle condition is \emph{linear} in the digits of the carry code $F$:
\begin{equation}
\begin{aligned}
&d_n(B+qA,F)+d_n\!\big(C+q\,d_q(B+qA,Q),F\big)\\
&\quad-d_n\!\big(d_n(C+qB,F)+nA,\Phi\big)-d_n\!\big(d_q(C+qB,Q)+qA,F\big)\equiv0 \pmod n ,
\end{aligned}
\label{dg:linearcocycle}
\end{equation}
for all $A,B,C\in A_q$. Together with the normalisation $d_n(m,F)=0$ for $m<q$ or
$d_q(0,m)=0$, this is a homogeneous linear system of $q^3$ equations modulo $n$ in the
$(q-1)^2$ unknown digits $d_n(B+qA,F)$, $A,B\ne0$. (Since $\varphi_A$ is additive,
$d_n(d_n(C+qB,F)+nA,\Phi)=\lambda_A\,d_n(C+qB,F)$ modulo $n$.) For abelian $K$, condition 4 of the
admissibility list holds automatically, so Eq.~(\ref{dg:linearcocycle}) is the only condition on $F$.

Define the digit-wise sum of two carry codes,
\begin{equation}
F\oplus F'=\sum_{m=0}^{q^2-1}n^m\,d_n\big(0,\;d_n(m,F)+d_n(m,F')\big).
\label{dg:Fsum}
\end{equation}

\begin{proposition}\label{prop:linearF}
With $Q$, $K=\mathbb{Z}_n$ and $\Phi$ fixed:
\begin{enumerate}
\item The admissible carry codes form a group $Z^2$ under $\oplus$, the solution space of
Eq.~(\ref{dg:linearcocycle}). Every admissible code is $F=c_1F_1\oplus\cdots\oplus c_rF_r$ for a basis
$F_1,\dots,F_r$ of solutions and coefficients $c_i\in A_n$.
\item For every unary code $g$ on $A_q\to A_n$ with $g(0)=0$, the carry code $\delta g$ with digits
\[
d_n(B+qA,\delta g)=d_n\!\Big(0,\;d_n\big(d_n(B,g)+nA,\Phi\big)-d_n\big(d_q(B+qA,Q),g\big)+d_n(A,g)\Big)
\]
is admissible. The set of these codes is a subgroup $B^2\subseteq Z^2$.
\item $F$ and $F\oplus\delta g$ define isomorphic groups: the relabelling
$x^{(0)}\mapsto x^{(0)}+g(x^{(1)})$ of the ground layer transforms one into the other.
\end{enumerate}
Hence one carry code per class of $H^2=Z^2/B^2$ suffices to produce all groups of the
diagram.
\end{proposition}

The action codes $\Phi$ are constructed in the same spirit, by choosing images of the
generators of $Q$ in $\mathrm{Aut}(K)$ that satisfy the relations of $Q$. When $K$ is
non-abelian Eq.~(\ref{dg:linearcocycle}) is no longer linear; such layers are usually
introduced as single nodes built with Proposition~\ref{prop:cayley}.

\subsection{Ordinary carries}\label{sec:dgordinarycarry}

The carry arrow takes its name from positional addition. Take a ground layer $K=\mathbb{Z}_p$ under an upper
layer $Q=\mathbb{Z}_p$ with trivial action, and the carry code with digits
\[
d_p(B+pA,F)=t\,c_p(A,B),\qquad t\in A_p,
\]
where $c_p$ is the carry of Proposition~\ref{prop:Bcarry}, $c_p(A,B)=\Bf\big(A+B-\tfrac{3p-2}{2},\tfrac{p-1}{2}\big)$.
By Proposition~\ref{prop:Bcocycle} it satisfies Eq.~(\ref{dg:linearcocycle}) for every $t$, so all these
diagrams are groups. For $t=0$ the diagram has no arrow and gives $C_p\oplus C_p$. For $t\ne0$ the element with
upper digit $1$ and ground digit $0$ has $p$-th power with upper digit $0$ and ground digit $t\neq0$, so it has
order $p^2$ and the group is $C_{p^2}$; for $t=1$ the law is ordinary addition modulo $p^2$, with the ground
layer playing the role of the higher digit. The codes $t\,c_p$ represent all classes of $H^2$, which has
order $p$: none of them with $t\neq0$ lies in $B^2$, since a coboundary gives a split, hence non-cyclic, group.
An exhaustive check confirms the counts: for $p=2$ there are $2$ admissible normalised carry codes and
$|H^2|=2$; for $p=3$ there are $9$, with $|B^2|=3$ and $|H^2|=3$, six of them giving $C_9$ and three giving
$C_3\oplus C_3$.

\subsection{Examples: orders 8 and 16}

\begin{example}[Order 8 without search]
For three $\mathbb{Z}_2$ layers the upper block is $\mathbb{Z}_2^2$ or $\mathbb{Z}_4$, and there are $9$ unknown
carry digits into the ground layer. Solving Eq.~(\ref{dg:linearcocycle}) modulo $2$ gives
$\dim Z^2=4$ over $\mathbb{Z}_2^2$ and $\dim Z^2=3$ over $\mathbb{Z}_4$, i.e.\ $16+8=24$ admissible
diagrams. This is exactly the number of associative diagrams found by testing all
$1024$ of them. Moreover $\dim B^2=1$ and $2$, so that $|H^2|=8$ and $2$ respectively.
\end{example}

\begin{example}[Boolean carry codes of order 8]\label{ex:dgboolean8}
Over the upper block $\mathbb{Z}_2^2$ write $A=a_0+2a_1$ and $B=b_0+2b_1$ with $a_i=d_2(i,A)$, $b_i=d_2(i,B)$, and take
the \emph{bilinear} carry codes
\[
d_2(B+4A,F)=d_2\big(0,\ c_{00}a_0b_0+c_{01}a_0b_1+c_{10}a_1b_0+c_{11}a_1b_1\big),\qquad c_{ij}\in A_2 .
\]
All $16$ are admissible and they represent all $8$ classes of $H^2$ (two codes per class). Listing
$c_{00}c_{01}c_{10}c_{11}$:
\begin{center}
\begin{tabular}{|l|l|}
\hline
Group & bilinear carry codes \\ \hline
$C_2^3$ & $0000$, $0110$ \\
$C_4\times C_2$ & $0001$, $0111$, $1000$, $1001$, $1110$, $1111$ \\
$D_4$ & $0010$, $0011$, $0100$, $0101$, $1010$, $1100$ \\
$Q_8$ & $1011$, $1101$ \\ \hline
\end{tabular}
\end{center}
The list is symmetric under exchanging the two operands and under exchanging the two bits. A single
product of bits is a carry: $a_0b_0$ is the carry of $\mathbb{Z}_4$ in the first bit, giving $C_4\times C_2$, while
the mixed product $a_1b_0$ gives $D_4$. The quaternion group needs three products,
$a_0b_0+a_1b_0+a_1b_1$. The code $a_0b_1+a_1b_0$ is the coboundary $\delta g$ of $g(A)=a_0a_1$ and gives the
elementary Abelian group. Each product $a_ib_j$ is a Boolean polynomial in the sense of
Corollary~\ref{cor:Bboolean}, so these carries are exactly the Boolean forms of degree two. An exhaustive
search over all $2^{16}$ carry tables confirms the classification: $8$ classes, of types $C_2^3$ (one class),
$C_4\times C_2$ (three), $D_4$ (three) and $Q_8$ (one).
\end{example}

\begin{example}[Order 16]
With four $\mathbb{Z}_2$ layers there are $24\times2^{49}\approx1.35\times10^{16}$ normalised diagrams with an
associative upper block, far too many to test. Solving Eq.~(\ref{dg:linearcocycle}) over each of the
$24$ associative diagrams of order $8$ gives $6912$ associative diagrams of order $16$:
\begin{center}
\begin{tabular}{|l|c|c|c|c|}
\hline
Upper block & diagrams of order 8 & $\dim Z^2$ & $\dim B^2$ & $|H^2|$ \\ \hline
$\mathbb{Z}_2^3$ & 2 & 10 & 4 & 64 \\
$\mathbb{Z}_4\times\mathbb{Z}_2$ & 10 & 8 & 5 & 8 \\
$D_4$ & 6 & 8 & 5 & 8 \\
$Q_8$ & 2 & 7 & 5 & 4 \\
$\mathbb{Z}_8$ & 4 & 7 & 6 & 2 \\ \hline
\end{tabular}
\end{center}
Indeed $2\cdot2^{10}+10\cdot2^8+6\cdot2^8+2\cdot2^7+4\cdot2^7=6912$. Taking one carry code per class of $H^2$
for each of the five groups of order $8$ yields $64+8+8+4+2=86$ tables, all associative,
which realise all $14$ groups of order $16$. This is to be expected: every group of
order $16$ has a central subgroup of order $2$, so one final $\mathbb{Z}_2$ layer with a carry
arrow reaches all of them.
\end{example}

The ternary codes $L(R)$ and $R'(R)$ remain useful as a final check of a construction,
or, through the defect (\ref{dg:defect}), when exploring structures beyond groups.

\section{Algorithm}

Given a diagram, the code of the group is obtained as follows.
\begin{enumerate}
\item Choose the layers $p_{M-1},\dots,p_0$ (top to ground) and the node codes $K_h$.
\item For each hollow layer choose the action code $\Phi_h$ and the carry code $F_h$;
draw an arrow from every upper layer on whose digit they depend.
\item Impose the normalisations $\varphi_0=\mathrm{id}$ and $f(0,\cdot)=f(\cdot,0)=0$.
\item For every pair $(a,b)$ compute $y^{(h)}$ from the top layer down with Eq.~(\ref{dg:stack})
and assemble $a\ast b$ with Eq.~(\ref{dg:stacksum}).
\item Check associativity (or conditions 3--5 layer by layer); for abelian layers construct
the carry codes directly as solutions of Eq.~(\ref{dg:linearcocycle}), one per class of $H^2$
(Section~\ref{sec:dgassoc}).
\item Compute $R$ from Eq.~(\ref{dg:R}) if the explicit code is required, and identify the
group through invariants such as the centre and the numbers of elements of each order.
\end{enumerate}

\section{Exercises}\label{sec:ch3-exercises}

\begin{exercise}
Construct $\mathbb{Z}_8$ with three $\mathbb{Z}_2$ layers and carry arrows only, and write the
carry codes $F_1$ and $F_0$ explicitly.
\end{exercise}

\begin{exercise}
Write the codes $\Phi$ of $D_n$ for $n=3,4,5$ and verify the formula for $y^{(0)}$.
\end{exercise}

\begin{exercise}
Show that the recipe for the Heisenberg group satisfies the cocycle condition and
compute its centre from Eq.~(\ref{dg:y0}).
\end{exercise}

\begin{exercise}
Verify that $Q_{16}$ is obtained from the generalised quaternion recipe and compare
its diagram with that of $D_8$.
\end{exercise}

\chapter{The main theorems of group theory in digital form}\label{ch:digitaltheorems}

In this chapter the classical theorems of finite group theory are formulated in
terms of the digit function, of codes and of the layers and arrows of the
diagrammatic method (Chapter~\ref{ch:diagrammatic}). In each case the statement is digital: subgroups
are ground blocks of mixed-radix labellings, cosets are fibres of layer projectors,
normal subgroups are closed layers, quotients are layer codes, actions are codes of
binary operators, and counting arguments are statements about base-$p$ digits. The
proofs follow the classical arguments, which become statements about digits.

\section{The digital setting}\label{sec:dtsetting}

A group of order $N$ is given by its code $R\in[0,N^{N^2})$,
\begin{equation}
a\ast b=d_N(aN+b,\;R),\qquad a,b\in \Ap{N} ,
\label{dt:law}
\end{equation}
and we always label the neutral element by $0$. The inverse of $a$ is denoted $\bar a$.
In this chapter the order of the group is written $N$ rather than $p$, because $p$ is reserved for the prime of the Sylow
theorems and $p_h$ for the radices of a labelling; otherwise the notation is that of the preceding chapters, with rows
$L_a$ and alphabets $\Ap{N}$.

\paragraph{Digital Kronecker delta.} For integers $k$ and $s\ge2$ one has
$d_s(k,1)=\lfloor s^{-k}\rfloor-s\lfloor s^{-k-1}\rfloor=\delta_{k,0}$. Hence, for $u,v\in \Ap{s}$,
\begin{equation}
\delta_{u,v}=d_s(u-v,\,1).
\label{dt:delta}
\end{equation}
With (\ref{dt:delta}) every counting statement becomes a sum of digits; for instance the
number of elements $x$ satisfying an equation $u(x)=v(x)$ is $\sum_x d_s\big(u(x)-v(x),1\big)$.

\paragraph{Labellings and isomorphism.} A \emph{labelling} is a permutation code $\lambda$ of $\Ap{N}$
with $\lambda(0)=0$. The relabelled code is
\begin{equation}
R^{\lambda}=\sum_{m=0}^{N^2-1}N^m\,\lambda\Big(d_N\big(\lambda^{-1}(d_N(1,m))\,N+\lambda^{-1}(d_N(0,m)),\;R\big)\Big),
\label{dt:relabel}
\end{equation}
and two codes $R,R'$ on $N$ symbols define isomorphic groups if and only if $R'=R^{\lambda}$ for
some labelling $\lambda$. Given a factorisation $N=p_0p_1\cdots p_{M-1}$, the labelled elements are
expanded in mixed radix, with layer projectors $x^{(h)}$ and upper blocks $U_h(x)$ as in
Eqs.~(\ref{dg:layers}) and (\ref{dg:upper}).

\paragraph{Homomorphisms.} A map $\psi:\Ap{N}\to \Ap{N'}$ is given by its code $\Psi\in[0,N'^{\,N})$,
$\psi(x)=d_{N'}(x,\Psi)$. It is a homomorphism from $R$ to $R'$ if and only if
\begin{equation}
d_{N'}\big(d_N(aN+b,R),\,\Psi\big)=d_{N'}\big(d_{N'}(a,\Psi)\,N'+d_{N'}(b,\Psi),\;R'\big)\qquad\forall a,b.
\label{dt:hom}
\end{equation}

\paragraph{Actions.} An action of $G$ on $\Ap{s}$ is given by a code $\Theta\in[0,s^{sN})$,
$\theta_a(x)=d_s(x+sa,\Theta)$, such that $\theta_0=\mathrm{id}$ and
\begin{equation}
d_s\big(d_s(x+sb,\Theta)+sa,\;\Theta\big)=d_s\big(x+s\,d_N(aN+b,R),\;\Theta\big)\qquad\forall a,b,x .
\label{dt:action}
\end{equation}
Thus an action is a binary operator $\Ap{N}\times \Ap{s}\to \Ap{s}$ that is compatible with $R$, exactly as
the action codes $\Phi$ of the arrows of Chapter~\ref{ch:diagrammatic}.

\paragraph{Adapted labellings.} Let $H$ be a subgroup of order $n$ and let $q=N/n$. A labelling $\lambda$
is \emph{adapted to $H$} (with radices $(n,q)$) if $\lambda(H)=\Ap{n}$ and the left cosets of $H$ are
mapped onto the fibres of the upper layer:
\[
\lambda(gH)=\{x\in \Ap{N}:\ \lfloor x/n\rfloor=\lfloor \lambda(g)/n\rfloor\}.
\]
In other words $H$ is the ground block, and two elements lie in the same left coset if
and only if their upper digits coincide.

\begin{center}
\begin{tabular}{|l|l|}
\hline
Group theory & Digital form \\ \hline
group of order $N$ & code $R\in[0,N^{N^2})$, $a\ast b=d_N(aN+b,R)$ \\
isomorphism & relabelling $R'=R^{\lambda}$ \\
left multiplication & rows $L_a=d_{N^N}(a,R)$ \\
subgroup $H$ & ground block of an adapted labelling \\
left cosets & fibres of the upper layer projector \\
index $[G:H]$ & radix of the upper layer \\
normal subgroup & ground block below a closed set of layers \\
quotient group & code $R_S$ of a closed set of layers \\
action on cosets & upward arrows of a tower labelling \\
subnormal series & tower labelling with locally closed layers \\
composition factors & layer codes of a composition labelling \\ \hline
\end{tabular}
\end{center}

\section{Cayley's theorem}

\begin{theorem}[Cayley, digital form]\label{thm:dtcayley}
Let $R$ be the code of a group of order $N$ and let $L_a=d_{N^N}(a,R)$ be its rows, so that
$R=\sum_{a}L_a\,N^{Na}$. Then each $L_a$ is a permutation code of $\Ap{N}$,
\[
a=d_N(0,L_a),\qquad L_{a\ast b}=L_a\circ L_b ,
\]
and $a\mapsto L_a$ is an injective homomorphism into the group of permutation codes of
$\Ap{N}$ under composition (\ref{dg:compose}). Every group of order $N$ is therefore isomorphic to a group
of $N$ unary codes ${}_1R_N$, and its code is recovered from the set of its rows.
\end{theorem}
\begin{proof}
By (\ref{dg:rows}), $d_N(b,L_a)=a\ast b$. Since $a\ast b=a\ast b'$ implies $b=b'$, the digits of $L_a$
are a permutation of $\Ap{N}$. Taking $b=0$ gives $d_N(0,L_a)=a$, so different elements have
different rows. Associativity gives $L_{a\ast b}=L_a\circ L_b$ by Proposition~\ref{prop:dgrows}.
\end{proof}

The converse, the construction of associative codes from closed sets of permutation
codes, is Proposition~\ref{prop:cayley}.

\section{Lagrange's theorem}

\begin{theorem}[Lagrange, digital form]\label{thm:dtlagrange}
For every subgroup $H$ of order $n$ of a group of order $N$ there is a labelling adapted to $H$.
Consequently $N=n\,q$: every subgroup is the ground block of a mixed-radix labelling
with radices $(n,q)$, its left cosets are the fibres of the upper layer, and the index
$[G:H]=q$ is the radix of the upper layer.
\end{theorem}
\begin{proof}
Choose a transversal $t_0=0,t_1,\dots,t_{q-1}$ of the left cosets and a bijection $\lambda_H:H\to \Ap{n}$
with $\lambda_H(0)=0$, and put $\lambda(t_u\ast h)=\lambda_H(h)+n\,u$. The left cosets are disjoint and
cover $G$, and $h\mapsto t_u\ast h$ is injective because the row $L_{t_u}$ is a permutation code
(Theorem~\ref{thm:dtcayley}). Hence every coset has $n$ elements, $\lambda$ is a bijection onto $A_{nq}$, and
$\lambda(t_uH)$ is the fibre $\{x:\lfloor x/n\rfloor=u\}$.
\end{proof}

\begin{corollary}\label{cor:dtorder}
The order of every element divides $N$. Equivalently, $L_a^{\circ N}=\sum_{v}v\,N^v$ (the identity
code) for every row of $R$.
\end{corollary}
\begin{proof}
Apply Theorem~\ref{thm:dtlagrange} to the cyclic subgroup generated by $a$, and Theorem~\ref{thm:dtcayley}.
\end{proof}

\begin{corollary}
If $N$ is prime, no labelling has more than one layer, and $R$ is a relabelling of the
code of $\mathbb{Z}_N$, $d_N(0,a+b)$.
\end{corollary}

Theorem~\ref{thm:dtlagrange} is the reason why the radices of every diagram of a group divide its
order, and why the tower labellings of Chapter~\ref{ch:diagrammatic} exist for every chain of subgroups.

\section{Group actions}

\begin{theorem}[Orbit--stabiliser, digital form]\label{thm:dtorbit}
Let $\Theta$ be the code of an action of $G$ on $\Ap{s}$, let $x\in \Ap{s}$, and let $G_x$ and $Gx$ be its
stabiliser and its orbit. There is a labelling of $G$ adapted to $G_x$, with radices
$(|G_x|,|Gx|)$, whose upper layer projector is the orbit map:
\[
\theta_g(x)=\iota\big(U_0(\lambda(g))\big)
\]
for a bijection $\iota:A_{|Gx|}\to Gx$. In particular $N=|G_x|\cdot|Gx|$, and the upward arrows of this
labelling (from the left operand into the upper layer) are the action itself:
$U_0\big(\lambda(a\ast g)\big)=\iota^{-1}\big(\theta_a(\iota(U_0(\lambda(g))))\big)$.
\end{theorem}
\begin{proof}
$\theta_g(x)=\theta_{g'}(x)$ if and only if $\theta_{\bar g\ast g'}(x)=x$, i.e.\ $g'\in gG_x$. Thus the fibres of the orbit map
are the left cosets of $G_x$, and Theorem~\ref{thm:dtlagrange} applied to $G_x$ gives the labelling, with
$\iota(u)=\theta_{t_u}(x)$. The last identity is (\ref{dt:action}).
\end{proof}

The number of points fixed by $a$ is, by (\ref{dt:delta}),
\begin{equation}
|\mathrm{Fix}(\theta_a)|=\sum_{x=0}^{s-1}d_s\big(d_s(x+sa,\Theta)-x,\;1\big).
\label{dt:fix}
\end{equation}

\begin{theorem}[Burnside's lemma, digital form]\label{thm:dtburnside}
The number of orbits of the action with code $\Theta$ is
\[
\frac1N\sum_{a=0}^{N-1}\sum_{x=0}^{s-1}d_s\big(d_s(x+sa,\Theta)-x,\;1\big).
\]
\end{theorem}
\begin{proof}
The double sum counts the pairs $(a,x)$ with $\theta_a(x)=x$, i.e.\ $\sum_x|G_x|=\sum_x N/|Gx|$ by
Theorem~\ref{thm:dtorbit}, and each orbit contributes $N$ to the last sum.
\end{proof}

\begin{lemma}[Ground-digit congruence]\label{lem:dtground}
Let a group $P$ of order $p^a$, $p$ prime, act on $\Ap{s}$, and let $\mathrm{Fix}$ be the set of fixed points. Then
\[
d_p\big(0,\,s\big)=d_p\big(0,\,|\mathrm{Fix}|\big).
\]
\end{lemma}
\begin{proof}
By Theorem~\ref{thm:dtorbit} every orbit has $p^k$ points, with $k\ge1$ for non-fixed points. Hence
$s=|\mathrm{Fix}|+\sum_ip^{k_i}$ with $k_i\ge1$, and the sum does not affect the digit $d_p(0,\cdot)$.
\end{proof}

\paragraph{Conjugation.} The conjugation action has the code $\Theta_c$ with $d_N(x+Na,\Theta_c)=a\ast x\ast\bar a$.
An element $x$ is central if and only if its row and its column coincide, so that
\begin{equation}
|Z(G)|=\sum_{x=0}^{N-1}\prod_{b=0}^{N-1}d_N\big(d_N(xN+b,R)-d_N(bN+x,R),\;1\big).
\label{dt:centre}
\end{equation}

\begin{theorem}[Class equation, digital form]\label{thm:dtclass}
Let $x_1,\dots,x_r$ represent the conjugacy classes with more than one element. Then
\[
N=|Z(G)|+\sum_{i=1}^r\frac{N}{|C(x_i)|},\qquad
|C(x)|=\sum_b d_N\big(d_N(xN+b,R)-d_N(bN+x,R),1\big).
\]
If $N=p^a$ with $a\ge1$, then $d_p\big(0,|Z(G)|\big)=0$, so the centre is non-trivial.
\end{theorem}
\begin{proof}
Apply Theorem~\ref{thm:dtorbit} to $\Theta_c$, whose fixed points form the centre and whose stabilisers
are the centralisers. For $N=p^a$, Lemma~\ref{lem:dtground} gives $d_p(0,|Z(G)|)=d_p(0,N)=0$, and $0\in Z(G)$.
\end{proof}

\begin{theorem}[Cauchy, digital form]\label{thm:dtcauchy}
If a prime $p$ divides $N$, then $G$ has an element of order $p$.
\end{theorem}
\begin{proof}
Read the integers $x\in A_{N^p}$ as strings of $p$ base-$N$ digits $g_i=d_N(i,x)$, and let $X$ be the set of
those $x$ with $g_{p-1}\ast g_{p-2}\ast\cdots\ast g_0=0$. Since $g_0$ is determined by the other digits,
$|X|=N^{p-1}$. The cyclic shift of digits
\[
\sigma(x)=\lfloor x/N\rfloor+N^{p-1}d_N(0,x)
\]
maps $X$ into itself, because $u\ast v=0$ implies $v\ast u=0$, and $\sigma^p=\mathrm{id}$, so $\mathbb{Z}_p$ acts on $X$. Its
fixed points are the repunits $x=g\,(N^p-1)/(N-1)$ with $g^{\ast p}=0$. By Lemma~\ref{lem:dtground},
$d_p(0,|\mathrm{Fix}|)=d_p(0,N^{p-1})=0$. Since $0\in\mathrm{Fix}$, $|\mathrm{Fix}|\ge p$, and there is $g\neq0$ with $g^{\ast p}=0$.
\end{proof}

\subsection{Transitive and sharply transitive actions}\label{sec:transitivity}
Orbits, stabilisers and the orbit--stabiliser theorem were introduced above (Theorem~\ref{thm:dtorbit}). Here they are used to
measure how strongly a group action mixes points, and then to explain the first digits in the Mathieu laws of Chapter~\ref{ch:sporadic}.

\begin{definition}\label{def:ktrans}
Let $1\le k\le n$. The action is \emph{$k$-transitive} if, for any two $k$-tuples of distinct points $(x_1,\dots,x_k)$ and
$(y_1,\dots,y_k)$, there is $g\in G$ with $g(x_i)=y_i$ for $i=1,\dots,k$. It is \emph{sharply $k$-transitive} if this $g$ is always
unique. A $1$-transitive action is called \emph{transitive}, and a sharply $1$-transitive action is called \emph{regular}.
\end{definition}

In words: a transitive group can move any point to any other point. A $k$-transitive group can move any $k$ points, in a given
order, to any other $k$ points, and a sharply $k$-transitive group can do so in exactly one way. A $k$-transitive action is
$l$-transitive for every $l\le k$, since an $l$-tuple of distinct points can be extended to a $k$-tuple.

\begin{lemma}\label{lem:sharp}
Let $G$ be $k$-transitive on $\Ap n$, let $x_1,\dots,x_k$ be distinct points, and let $G_{x_1,\dots,x_k}$ be their pointwise
stabiliser. Then
\[
|G|=n(n-1)\cdots(n-k+1)\,|G_{x_1,\dots,x_k}|,
\]
and the following are equivalent:
\begin{enumerate}
\item $G$ is sharply $k$-transitive;
\item $G_{x_1,\dots,x_k}=1$;
\item $g\mapsto(g(x_1),\dots,g(x_k))$ is a bijection from $G$ onto the $n(n-1)\cdots(n-k+1)$ $k$-tuples of distinct points.
\end{enumerate}
\end{lemma}
\begin{proof}
For $l\le k$, the stabiliser $G_{x_1,\dots,x_{l-1}}$ is transitive on the $n-l+1$ points other than $x_1,\dots,x_{l-1}$: given
such a point $y$, $k$-transitivity gives $g$ with $g(x_i)=x_i$ for $i<l$ and $g(x_l)=y$. Applying the orbit--stabiliser theorem
at each of the $k$ steps gives the order formula. Two elements $g,g'$ have the same images of $x_1,\dots,x_k$ exactly when
$g^{-1}g'\in G_{x_1,\dots,x_k}$. The map in (3) is onto by $k$-transitivity, so (2) and (3) are equivalent, and (1) implies (2).
Conversely, assume (2), and let $g(x_i')=g'(x_i')$ for some distinct $x_1',\dots,x_k'$. Choose $h$ with $h(x_i)=x_i'$. Then
$h^{-1}g^{-1}g'h$ fixes every $x_i$, so it is $1$ and $g=g'$.
\end{proof}

\begin{example}\label{ex:transitive}
\begin{itemize}
\item The cyclic group generated by $x\mapsto\dg{n}{0}{x+1}$ is regular on $\Ap n$.
\item $\Sym_n$ is sharply $n$-transitive, and also sharply $(n-1)$-transitive, since the image of the last point is forced.
      $\Alt_n$, $n\ge3$, is sharply $(n-2)$-transitive: of the two permutations with prescribed images of $n-2$ points,
      exactly one is even.
\item For a prime power $q$, the maps $x\mapsto ax+b$ with $a\ne0$ act sharply $2$-transitively on $\mathbb F_q$. The maps
      $x\mapsto(ax+b)/(cx+d)$ with $ad-bc\ne0$ form $PGL(2,q)$, which is sharply $3$-transitive on the $q+1$ points of the
      projective line.
\item $M_{11}$ and $M_{12}$ are sharply $4$- and $5$-transitive. $M_{22}$, $M_{23}$ and $M_{24}$ are $3$-, $4$- and
      $5$-transitive but not sharply so: the pointwise stabiliser of $3$, $4$ or $5$ points, respectively, has order $48$
      (Proposition~\ref{prop:mathieuchain}).
\end{itemize}
\end{example}

By a theorem of Jordan, for $k\ge4$ the only sharply $k$-transitive groups are $\Sym_k$, $\Sym_{k+1}$, $\Alt_{k+2}$ and, for
$k=4$ and $k=5$, the groups $M_{11}$ and $M_{12}$ \cite{Wilson}. The two small Mathieu groups are thus the only highly sharply
transitive groups that are not symmetric or alternating. By the classification of finite simple groups, every $6$-transitive
group is symmetric or alternating.

Transitivity is visible in the labels of Theorem~\ref{thm:chainlaw} in Chapter~\ref{ch:sporadic}. If $G$ is $k$-transitive and $b_1,\dots,b_k$ are
the first base points, Lemma~\ref{lem:sharp} gives $r_l=n-l+1$ for $l\le k$, so the first $k$ digits of a label have the
falling-factorial radices $n,n-1,\dots,n-k+1$. The group is sharply $k$-transitive exactly when the chain stops there, as for
$M_{11}$ and $M_{12}$ in Proposition~\ref{prop:mathieuchain}. For $M_{22}$, $M_{23}$ and $M_{24}$ two more digits, with radices
$16$ and $3$, code the pointwise stabiliser of order $48$.

\begin{remark}[Canonical labels]\label{rem:canonical}
For a sharply $k$-transitive group on $\Ap n$, Lemma~\ref{lem:sharp}(3) gives a label that needs no choice of representatives.
Put $y_l=g(l-1)$ for $l=1,\dots,k$, let $c_l$ be the number of points $y<y_l$ different from $y_1,\dots,y_{l-1}$, and set
\[
\kappa(g)=c_1+n\big(c_2+(n-1)\big(c_3+\dots+(n-k+2)c_k\big)\big).
\]
This is the falling-factorial rank of the word $y_1\cdots y_k$, and $\kappa$ is a bijection $G\to\Ap{|G|}$; for $M_{11}$ this
was checked on all $7920$ elements. For the elements of Examples~\ref{ex:m11dec} and~\ref{ex:m11mul}:
\begin{center}
\begin{tabular}{ccc}
chain label $\lambda$ & $(g(0),g(1),g(2),g(3))$ & canonical label $\kappa$\\ \hline
$1234$ & $(2,10,5,8)$ & $6481$\\
$5678$ & $(2,10,9,6)$ & $5931$\\
$2524=1234*5678$ & $(5,1,4,0)$ & $346$
\end{tabular}
\end{center}
For instance $\kappa=2+11\cdot9+110\cdot4+990\cdot6=6481$ for the first row. In the transported law, $6481*5931=346$. The
images of a product are not determined by the images of its factors: $(g\circ h)(i)=g(h(i))$ needs $g$ at the points
$h(0),\dots,h(3)=2,10,9,6$, and three of them are not among $0,1,2,3$. Finding these values means reconstructing $g$ from
its four images, which is what the chain of Theorem~\ref{thm:chainlaw} does. So the labels $\kappa$ are canonical, while the
labels $\lambda$ are the ones that are cheap to multiply.
\end{remark}

\subsection{Conjugacy classes}\label{ssec:dtclasses}
Theorem~\ref{thm:dtclass} gives the class equation for a group presented by a digital law. This subsection collects the general facts about conjugacy classes that the rest of the book uses, independently of the encoding: orbit--stabiliser for conjugation, the class equation, the number of classes, and the splitting of classes in a subgroup of index two. Power maps, rationality and algorithms follow in Chapter~\ref{ch:cc}, and the classes of $\Sym_n$ and $\Alt_n$ in Section~\ref{sec:sn:classes}.

\begin{definition}
Let $G$ be a group. For $g,x\in G$ write $g^{x}=x^{-1}gx$. The \emph{conjugacy class} of $g$ is $g^{G}=\{g^{x}:x\in G\}$, its \emph{centraliser} is $C_G(g)=\{x\in G:gx=xg\}$, and the \emph{centre} is $Z(G)=\{g:g^{G}=\{g\}\}$. A function $G\to\mathbb{C}$ constant on classes is a \emph{class function}.
\end{definition}

For permutation codes conjugation is a pure digit relocation: by Proposition~\ref{prop:conj}, $S\circ R\circ S^{-1}$ has digit $S(R(n))$ in position $S(n)$. Conjugating a code therefore costs $p$ digit look-ups, exactly as one composition~\eqref{eq:comp}; no multiplication table is needed.

\begin{theorem}[Orbit--stabiliser for conjugation]\label{thm:cc:os}
The map $C_G(g)x\mapsto g^{x}$ is a bijection from the right cosets of $C_G(g)$ onto $g^{G}$. Hence $|g^{G}|=|G:C_G(g)|$, and $|g^{G}|$ divides $|G|$.
\end{theorem}
\begin{proof}
$g^{x}=g^{y}$ iff $x^{-1}gx=y^{-1}gy$ iff $yx^{-1}$ commutes with $g$ iff $C_G(g)x=C_G(g)y$. So the map is well defined and injective; it is surjective by definition of $g^{G}$. Lagrange's theorem gives the divisibility.
\end{proof}

\begin{corollary}[Class equation]\label{cor:cc:classeq}
If $g_1,\dots,g_k$ represent the noncentral classes, then $|G|=|Z(G)|+\sum_{i=1}^{k}|G:C_G(g_i)|$, each summand dividing $|G|$ and exceeding $1$.
\end{corollary}
\begin{proof}
The classes partition $G$; the singleton classes are exactly the central elements. Apply Theorem~\ref{thm:cc:os}.
\end{proof}

Theorem~\ref{thm:dtclass} is the same equation, read off the digits of a Cayley table.

\begin{corollary}\label{cor:cc:pcentre}
A nontrivial finite $p$-group has nontrivial centre. A group of order $p^{2}$ is abelian.
\end{corollary}
\begin{proof}
In the class equation of $|G|=p^{a}$, $a\ge1$, every noncentral summand is a positive power of $p$, so $p$ divides $|Z(G)|\ge1$. If $|G|=p^{2}$ and $Z(G)\ne G$, then $G/Z(G)$ has order $p$, hence is cyclic, generated by $gZ(G)$; then every element is $g^{i}z$ with $z$ central, and any two such commute, a contradiction.
\end{proof}

\begin{theorem}[Counting classes]\label{thm:cc:count}
The number $k(G)$ of conjugacy classes satisfies
\[
k(G)=\frac{1}{|G|}\sum_{g\in G}|C_G(g)|=\frac{|\{(x,y)\in G^{2}:xy=yx\}|}{|G|}.
\]
\end{theorem}
\begin{proof}
The second equality counts commuting pairs by their first coordinate. For the first, group the sum by classes: on a class $C$ with representative $g$, $\sum_{x\in C}|C_G(x)|=|C|\,|C_G(g)|=|G|$ by Theorem~\ref{thm:cc:os}, because conjugate elements have conjugate, hence equipotent, centralisers. Summing $|G|$ over the $k(G)$ classes gives the claim.
\end{proof}

The ratio $k(G)/|G|$ is the probability that two uniformly random elements commute. It is the first hint of a recurring theme: centraliser orders, not class sizes, are the natural currency of computation.

\begin{theorem}[Splitting]\label{thm:cc:split}
Let $H\le G$ have index $2$ and let $h\in H$. Then $h^{G}\subseteq H$, and
\begin{enumerate}
\item if $C_G(h)\not\le H$, then $h^{G}=h^{H}$ is a single $H$-class;
\item if $C_G(h)\le H$, then $h^{G}$ is the disjoint union of two $H$-classes $h^{H}$ and $(h^{t})^{H}$, for any $t\in G\setminus H$, each of size $\tfrac12|h^{G}|$.
\end{enumerate}
\end{theorem}
\begin{proof}
$H$ is normal, so $h^{G}\subseteq H$. By Theorem~\ref{thm:cc:os}, $|h^{G}|=|G|/|C_G(h)|$ and $|h^{H}|=|H|/|C_H(h)|$ with $C_H(h)=C_G(h)\cap H$. If $C_G(h)\not\le H$, then $C_H(h)$ has index $2$ in $C_G(h)$, so $|h^{H}|=\frac{|G|/2}{|C_G(h)|/2}=|h^{G}|$. If $C_G(h)\le H$, then $|h^{H}|=\tfrac12|h^{G}|$. Since $G=H\cup tH$, $h^{G}=h^{H}\cup (h^{t})^{H}$, and the two pieces have equal size because conjugation by $t$ is a bijection carrying one onto the other; they are disjoint because their union has size $|h^{G}|=2|h^{H}|$.
\end{proof}

\begin{example}
In $\Alt_5\le\Sym_5$ a $5$-cycle has centraliser $\langle g\rangle$ of order $5$ in $\Sym_5$, which lies in $\Alt_5$, so its class of $24$ elements splits into two $\Alt_5$-classes of $12$. A $3$-cycle has centraliser of order $6$ in $\Sym_5$, which contains a transposition, so the $20$ three-cycles remain one class. Hence the class equation of $\Alt_5$ is $60=1+15+20+12+12$; see Theorem~\ref{thm:sn:altclasses} for the general rule.
\end{example}

\subsection{Counting colourings: \texorpdfstring{P\'olya--Redfield}{Polya--Redfield} and de Bruijn}\label{ssec:polya}
Burnside's lemma (Theorem~\ref{thm:dtburnside}) counts the orbits of a group on a set of points. Many counting problems ask instead for the orbits on the \emph{colourings} of a set of positions, and then the fixed points can be counted cycle by cycle. The two theorems below~\cite{redfield1927,polya1937,debruijn1959,debruijn1964} are used for magmas (Section~\ref{sec:magma:count}), necklaces and bracelets (Propositions~\ref{prop:necklaces} and~\ref{prop:bracelets}), the group $M_{12}$ (Section~\ref{sec:cc:m12}) and the symmetric groups (Section~\ref{sec:sn:classes}).

Let $G$ act on a finite set $X$ with $|X|=s$; for a digital action, $X=\Ap{s}$ and the action has the code $\Theta$. Each $g\in G$ induces a permutation of $X$. Write $c_k(g)$ for its number of cycles of length $k$ (fixed points are cycles of length $1$) and $c(g)=\sum_kc_k(g)$. The \emph{cycle index} of the action is
\[
Z(G;x_1,x_2,\dots)=\frac1{|G|}\sum_{g\in G}\ \prod_{k\ge1}x_k^{c_k(g)}.
\]
If $Y$ is a set of $m$ colours, $G$ acts on the set $Y^{X}$ of colourings by $(g\cdot f)(x)=f(g^{-1}x)$. In digits, a colouring is a string of $s$ base-$m$ digits and $g$ relocates the digits, just as conjugation relocates the digits of a permutation code (Proposition~\ref{prop:conj}).

\begin{theorem}[P\'olya--Redfield]\label{thm:polya}
Let $w:Y\to\mathcal R$ take values in a commutative ring $\mathcal R$ and give the colouring $f$ the weight $w(f)=\prod_{x\in X}w(f(x))$, which is constant on orbits. Then
\[
\sum_{O\in Y^{X}/G}w(O)=Z(G;p_1,p_2,\dots),\qquad p_k=\sum_{y\in Y}w(y)^{k}.
\]
In particular, the number of orbits of $G$ on $Y^{X}$ is $\dfrac1{|G|}\sum_{g\in G}m^{c(g)}=Z(G;m,m,\dots)$.
\end{theorem}
\begin{proof}
A colouring is fixed by $g$ if and only if it is constant on every cycle of $g$. Choosing one colour $y_\Gamma$ for each cycle $\Gamma$ gives such a colouring, of weight $\prod_\Gamma w(y_\Gamma)^{|\Gamma|}$, so the colourings fixed by $g$ have total weight $\prod_kp_k^{c_k(g)}$. On the other hand, Theorem~\ref{thm:dtorbit} for the action on $Y^{X}$ gives
\[
\sum_{g\in G}\ \sum_{f=g\cdot f}w(f)=\sum_{f}|G_f|\,w(f)=\sum_{O}\ \sum_{f\in O}\frac{|G|}{|O|}\,w(f)=|G|\sum_{O}w(O).
\]
Divide by $|G|$; the choice $w\equiv1$ gives the count.
\end{proof}

\begin{theorem}[de Bruijn]\label{thm:debruijn}
Let in addition a group $H$ act on the set $Y$ of $m$ colours, and let $G\times H$ act on $Y^{X}$ by $((g,h)\cdot f)(x)=h\big(f(g^{-1}x)\big)$. For $h\in H$ and $k\ge1$ let $\operatorname{fix}_k(h)=|\{y\in Y:h^{k}y=y\}|$, the number of points of $Y$ on cycles of $h$ whose length divides $k$. Then the number of orbits of $G\times H$ on $Y^{X}$ is
\[
\frac1{|G|\,|H|}\sum_{g\in G}\sum_{h\in H}\ \prod_{k\ge1}\operatorname{fix}_k(h)^{c_k(g)}
=\frac1{|H|}\sum_{h\in H}Z\big(G;\operatorname{fix}_1(h),\operatorname{fix}_2(h),\dots\big).
\]
\end{theorem}
\begin{proof}
A colouring $f$ is fixed by $(g,h)$ if and only if $f(gx)=h\,f(x)$ for all $x$. On a cycle $x_0,gx_0,\dots,g^{k-1}x_0$ of $g$ of length $k$ the value $y=f(x_0)$ determines $f(g^{i}x_0)=h^{i}y$, and the equation closes up exactly when $h^{k}y=y$; conversely every such $y$ gives a solution on the cycle. Hence $(g,h)$ fixes $\prod_{k}\operatorname{fix}_k(h)^{c_k(g)}$ colourings, and the claim is Burnside's lemma (Theorem~\ref{thm:dtburnside}) for $G\times H$.
\end{proof}

For trivial $H$ one has $\operatorname{fix}_k(h)=m$ and Theorem~\ref{thm:debruijn} is the unweighted case of Theorem~\ref{thm:polya}. For $H=\Sym(Y)$ it counts colourings up to renaming of the colours as well. The same argument, with $g$ acting on a set of pairs and $h$ on the values, counts operator tables up to relabelling (Section~\ref{sec:magma:count}).

\begin{example}[The square]
$D_8$ acts on the vertices $\Ap{4}$ of a square by $x\mapsto x+1$ and $x\mapsto-x$ modulo $4$. The identity has cycle type $1^{4}$, the two quarter turns $4$, the half turn and the two reflections through the midpoints of opposite edges $2^{2}$, and the two reflections through opposite vertices $2\,1^{2}$. Hence
\[
Z(D_8)=\tfrac18\big(x_1^{4}+2x_4+3x_2^{2}+2x_1^{2}x_2\big),
\]
and the number of orbits on $m$-colourings is $(m^{4}+2m^{3}+3m^{2}+2m)/8$, which is $6$ for $m=2$ and $21$ for $m=3$. With the weights $w(\text{white})=1$, $w(\text{black})=t$ one has $p_k=1+t^{k}$ and Theorem~\ref{thm:polya} gives $1+t+2t^{2}+t^{3}+t^{4}$: there are two inequivalent ways to colour two vertices black.
\end{example}

\section{Homomorphisms and the isomorphism theorems}

\begin{theorem}[First isomorphism theorem, digital form]\label{thm:dtfirst}
Let $\Psi$ be the code of a homomorphism $\psi$ from $R$ to $R'$ with kernel $K$ of order $k$, and let
$q=N/k$. There is a labelling of $G$ adapted to $K$ with radices $(k,q)$ in which the upper
layer is closed (it receives no arrows) and
\[
\psi(g)=\iota\big(U_0(\lambda(g))\big)
\]
for a bijection $\iota$ of $\Ap{q}$ onto the image. The code $R_{\{1\}}$ of the closed upper layer
(Proposition~\ref{prop:closedquotient}) is a relabelling of the code of $\psi(G)$. Conversely, every closed set of
layers defines a homomorphism whose kernel is the ground block below it.
\end{theorem}
\begin{proof}
$\psi(g)=\psi(g')$ if and only if $g'\in gK$, so by Theorem~\ref{thm:dtlagrange} the fibres of $\psi$ are the fibres of
the upper layer of a labelling adapted to $K$. Then $U_0(\lambda(a\ast b))$ is determined by
$\psi(a\ast b)=\psi(a)\psi(b)$, i.e.\ by $U_0(\lambda(a))$ and $U_0(\lambda(b))$: no arrow enters the upper layer.
The converse is Proposition~\ref{prop:closedquotient}.
\end{proof}

Hence normal subgroups, homomorphic images and closed sets of layers are the same
objects, and $G/K$ is the code of the closed upper layer.

\begin{theorem}[Second isomorphism theorem, digital form]\label{thm:dtsecond}
Let $H\le G$ and $K\trianglelefteq G$, and take a labelling of $G$ adapted to $K$ with closed upper layer. The
restriction of the upper projector $U_0$ to $H$ is a homomorphism with kernel $H\cap K$ whose image is
the set of upper digits of $HK$. Hence the subcode of $R_{\{1\}}$ on these digits is a relabelling of
the code of the closed upper layer of a labelling of $H$ adapted to $H\cap K$:
$HK/K\cong H/(H\cap K)$.
\end{theorem}
\begin{proof}
The restriction of a homomorphism is a homomorphism; its kernel is $H\cap K$ and its image is
$\{U_0(\lambda(h))\}=\{U_0(\lambda(hk))\}$. Apply Theorem~\ref{thm:dtfirst} to it.
\end{proof}

\begin{theorem}[Third isomorphism theorem, digital form]\label{thm:dtthird}
Let $K\le L$ be normal subgroups of $G$ and take the tower labelling of the chain $1<K<L<G$, with
radices $(|K|,[L:K],[G:L])$. Then $\{2\}$ and $\{1,2\}$ are closed sets of layers, $R_{\{1,2\}}$ is the code
of $G/K$, the layer $\{2\}$ is closed in $R_{\{1,2\}}$, and the code of $\{2\}$ computed from $R_{\{1,2\}}$
coincides with $R_{\{2\}}$ computed from $R$: closed projections compose,
$\pi_{\{2\}}=\pi_{\{2\}}\circ\pi_{\{1,2\}}$, i.e.\ $(G/K)/(L/K)\cong G/L$.
\end{theorem}
\begin{proof}
By Proposition~\ref{prop:towerright}, the upper block above layer $0$ (resp.\ $1$) is closed because $K$
(resp.\ $L$) is normal. Since $\pi_{\{2\}}(x)$ is a digit of $\pi_{\{1,2\}}(x)$ and both are homomorphisms, the
layer $2$ receives no arrow from layer $1$ in $R_{\{1,2\}}$, and the two codes coincide.
\end{proof}

\begin{theorem}[Correspondence theorem, digital form]\label{thm:dtcorr}
Let $K\trianglelefteq G$ and take a labelling adapted to $K$ with closed upper layer. The subgroups of $G$
containing $K$ are exactly the unions of fibres $\{x:\ U_0(x)\in S\}$ with $S$ a subgroup of $R_{\{1\}}$;
this correspondence preserves inclusions, indices (the radices) and normality.
\end{theorem}
\begin{proof}
A subgroup containing $K$ is a union of cosets of $K$, i.e.\ of fibres of $U_0$, and since $U_0$ is a
homomorphism the set of their upper digits is a subgroup of $R_{\{1\}}$; conversely the preimage of a
subgroup is a subgroup. Normality is preserved because $U_0$ is surjective.
\end{proof}

\section{Sylow theorems in digital form}\label{sec:dtsylow}

Let $p$ be a prime and let $a$ be the number of trailing zero digits of $N$ in base $p$:
\begin{equation}
d_p(k,N)=0\ \ (k<a),\qquad d_p(a,N)\neq0,\qquad m=N/p^a,\qquad d_p(0,m)=d_p(a,N).
\label{dt:sylowdigits}
\end{equation}

\begin{lemma}[Lucas]\label{lem:dtlucas}
For all integers $n,k\ge0$,
$\displaystyle\binom{n}{k}\equiv\prod_{i\ge0}\binom{d_p(i,n)}{d_p(i,k)}\pmod p$.
In particular $\displaystyle\binom{N}{p^a}\equiv d_p(a,N)\not\equiv0\pmod p$.
\end{lemma}
\begin{proof}
Modulo $p$, $(1+x)^p\equiv1+x^p$, hence $(1+x)^n\equiv\prod_i(1+x^{p^i})^{d_p(i,n)}$. By uniqueness of base-$p$
expansions, the coefficient of $x^k$ only receives the contribution of the digits $d_p(i,k)$. For
$k=p^a$ the only non-zero digit of $k$ is $d_p(a,k)=1$, which leaves the factor $d_p(a,N)$.
\end{proof}

\begin{remark}
By Kummer's theorem the exponent of $p$ in $\binom{n}{k}$ is the number of carries in the base-$p$
addition of $k$ and $n-k$, i.e.\ the number of times the carry arrows of $\mathbb{Z}_{p^2}$, with
$f(A,B)=\lfloor(A+B)/p\rfloor$, fire. For $n=N$ and $k=p^a$ there is no carry. This was checked for
$S_4$, $A_5$, $M_{11}$, $M_{12}$ and $J_1$ and every prime dividing their orders; for example
$\binom{7920}{5}\equiv d_5(1,7920)=4\pmod5$ for $M_{11}$.
\end{remark}

\begin{theorem}[Sylow I, digital form]\label{thm:dtsylow1}
For every prime $p$, $G$ has a labelling with radices $(p,\dots,p,m)$, with $a$ radices equal to $p$,
whose $a$ lowest layers form a subgroup node of order $p^a$ in which every upper block is
closed, so that inside the node all arrows point downwards.
\end{theorem}
\begin{proof}
Let $\Omega$ be the set of subsets of $G$ with $p^a$ elements, acted on by left multiplication. By
Lemma~\ref{lem:dtlucas}, $d_p(0,|\Omega|)\neq0$, so some orbit has a size with non-zero ground digit. If $X$ lies in
this orbit and $H$ is its stabiliser, Theorem~\ref{thm:dtorbit} gives $d_p(0,N/|H|)\neq0$, hence $p^a$ divides $|H|$.
On the other hand $Hx\subseteq X$ for $x\in X$, so $|H|\le p^a$. Thus $|H|=p^a$, and Theorem~\ref{thm:dtlagrange} gives a
labelling with radices $(p^a,m)$ whose ground block is $H$. By Theorem~\ref{thm:dtclass} a group of order $p^b$,
$b\ge1$, has a central subgroup of order $p$; by induction $H$ has a chain
$1=H_0<H_1<\cdots<H_a=H$ of subgroups that are normal in $H$, with all indices equal to $p$. Its tower
labelling refines the ground block, and its upper blocks inside $H$ are closed by
Proposition~\ref{prop:towerright}.
\end{proof}

\begin{theorem}[Sylow II, digital form]\label{thm:dtsylow2}
Let $Q$ be a Sylow $p$-subgroup and take a labelling adapted to $Q$, with radices $(p^a,m)$. For
every $p$-subgroup $P$ there is an upper digit $u\in A_m$ that is not moved by the upward arrow
from $P$, i.e.\ $U_0(\lambda(a\ast g))=u$ for all $a\in P$ whenever $U_0(\lambda(g))=u$. Consequently $P$ is contained in a
conjugate of $Q$, and all Sylow $p$-subgroups are conjugate.
\end{theorem}
\begin{proof}
By Theorem~\ref{thm:dtorbit} the upward arrow from the left operand is the action of $G$ on the cosets
of $Q$. Restricting it to $P$, Lemma~\ref{lem:dtground} gives $d_p(0,|\mathrm{Fix}|)=d_p(0,m)\neq0$, so some coset $gQ$ is
fixed: $PgQ=gQ$, i.e.\ $\bar g\ast P\ast g\subseteq Q$.
\end{proof}

\begin{theorem}[Sylow III, digital form]\label{thm:dtsylow3}
If $n_p$ is the number of Sylow $p$-subgroups, then
\[
d_p(0,n_p)=1,\qquad n_p\mid m.
\]
\end{theorem}
\begin{proof}
Let $Q$ act by conjugation on the set of Sylow $p$-subgroups. If $Q$ fixes $Q'$, then $Q$ and $Q'$ are
Sylow subgroups of the normaliser of $Q'$, in which $Q'$ is normal; by Theorem~\ref{thm:dtsylow2} they are
conjugate there, so $Q=Q'$. Lemma~\ref{lem:dtground} gives $d_p(0,n_p)=d_p(0,1)=1$. By Theorem~\ref{thm:dtorbit},
$n_p=N/|N_G(Q)|$, which divides $m$ because $Q\le N_G(Q)$.
\end{proof}

For $\Alt_5$, $n_2=5$, $n_3=10$ and $n_5=6$, whose ground digits in bases $2$, $3$ and $5$ are all equal to $1$.

\begin{proposition}[Labels of the Sylow subgroups]\label{prop:dtsylowlabels}
Let $N=p^am$ with $p\nmid m$, let $Q$ be a Sylow $p$-subgroup, and take the labelling of Theorem~\ref{thm:dtsylow1}, coarsened
to radices $(p^a,m)$, so that $Q=\Ap{p^a}$ and the upper digit $u\in\Ap{m}$ of a label is its left coset of $Q$. Let $g_u$ be
the element of label $p^au$ and put
\[
Q_u=\{\,g_u\ast h\ast\bar g_u : h\in\Ap{p^a}\,\},\qquad u\in\Ap{m}.
\]
Then every $Q_u$ is a Sylow $p$-subgroup, every Sylow $p$-subgroup is one of the $Q_u$, and $Q_u=Q_v$ if and only if
$\bar g_v\ast g_u\in N_G(Q)$. Each Sylow $p$-subgroup occurs $|N_G(Q):Q|=m/n_p$ times in the list $(Q_u)_{u\in\Ap m}$.
\end{proposition}
\begin{proof}
Conjugation is an automorphism, so $Q_u=g_uQ\bar g_u$ is a subgroup of order $p^a$. If $Q'$ is a Sylow $p$-subgroup,
Theorem~\ref{thm:dtsylow2} gives $Q'=gQ\bar g$ for some $g$; if $u$ is the upper digit of $\lambda(g)$, then $g=g_u\ast h'$ with
$h'\in Q$, and $h'Q\bar h'=Q$ gives $Q'=Q_u$. Next, $Q_u=Q_v$ if and only if $\bar g_v\ast g_u$ normalises $Q$. The
elements $g_u$ represent the $m$ left cosets of $Q$, and $g_u$, $g_v$ give the same conjugate exactly when they lie in the same
left coset of $N_G(Q)$, which is a union of $|N_G(Q):Q|$ cosets of $Q$. By Theorem~\ref{thm:dtsylow3},
$|N_G(Q):Q|=(N/n_p)/p^a=m/n_p$.
\end{proof}

The proposition gives the Sylow subgroups in a labelling adapted to $p$. In any other labelling they are still finite sets of
labels, computable from the law: for $p$ dividing $N$ only once they are the sets $\{0,x,x\ast x,\dots\}$ generated by the labels
$x$ of order $p$.

\begin{example}[Sylow labels of $\Alt_5$]\label{ex:dtsylowa5}
In the Lehmer law of $\Alt_5$ (Tables~\ref{tab:an:law5a} and~\ref{tab:an:law5b}) the Sylow $3$- and $5$-subgroups are the
cyclic subgroups generated by the labels of order $3$ and $5$, and the Sylow $2$-subgroups are the Klein groups
$\{0,a,b,a\ast b\}$ formed by two commuting labels of order $2$. Computing them from the law gives:
\begin{center}\small
\begin{tabular}{rrp{0.72\textwidth}}
\hline
$p$ & $n_p$ & labels of the Sylow $p$-subgroups\tabularnewline
\hline
$2$ & $5$ & \raggedright $\{0,1,58,59\}$, $\{0,5,50,55\}$, $\{0,7,48,49\}$, $\{0,12,21,33\}$, $\{0,17,34,41\}$\tabularnewline
$3$ & $10$ & \raggedright $\{0,2,6\}$, $\{0,8,27\}$, $\{0,9,47\}$, $\{0,10,25\}$, $\{0,15,45\}$, $\{0,20,40\}$, $\{0,23,31\}$, $\{0,24,36\}$, $\{0,28,29\}$, $\{0,30,35\}$\tabularnewline
$5$ & $6$ & \raggedright $\{0,3,19,46,57\}$, $\{0,4,18,26,52\}$, $\{0,11,38,42,54\}$, $\{0,13,37,44,51\}$, $\{0,14,32,43,56\}$, $\{0,16,22,39,53\}$\tabularnewline
\hline
\end{tabular}
\end{center}
The five Sylow $2$-subgroups contain the $15$ labels of order $2$, the ten Sylow $3$-subgroups the $20$ labels of order $3$,
and the six Sylow $5$-subgroups the $24$ labels of order $5$; together with $0$ these are the $60$ labels, in accordance with
the class equation $60=1+15+20+12+12$. Distinct Sylow subgroups for the same prime meet only in $0$.
\end{example}

\section{Series: Zassenhaus, Schreier and Jordan--H\"older}

\begin{definition}[Subnormal and composition labellings]
A \emph{subnormal labelling} is the tower labelling of a chain
$1=H_0\trianglelefteq H_1\trianglelefteq\cdots\trianglelefteq H_M=G$ in which each subgroup is normal in the next. Inside the
block $H_{k+1}$ (the labels with $U_k=0$), layer $k$ is closed, by Proposition~\ref{prop:towerright} applied to $H_{k+1}$;
its code $R^{(k)}$ on $p_k=[H_{k+1}:H_k]$ symbols is the \emph{layer code}, the code of $H_{k+1}/H_k$. A
\emph{refinement} factorises radices by inserting layers. A \emph{composition labelling} is a
subnormal labelling that admits no proper refinement, i.e.\ all of whose layer codes are
simple: by Theorem~\ref{thm:simple}, each layer code has prime radix or a strongly connected arrow
graph for every labelling.
\end{definition}

\begin{lemma}[Zassenhaus butterfly, digital form]\label{lem:dtzassenhaus}
Let $A'\trianglelefteq A\le G$ and $B'\trianglelefteq B\le G$, and let $D=(A'\cap B)(A\cap B')$. Then $D\trianglelefteq A\cap B$, and the
following three codes are relabellings of each other:
\begin{enumerate}
\item the closed upper layer of a labelling of $A'(A\cap B)$ adapted to $A'(A\cap B')$;
\item the closed upper layer of a labelling of $B'(A\cap B)$ adapted to $B'(A'\cap B)$;
\item the closed upper layer of a labelling of $A\cap B$ adapted to $D$.
\end{enumerate}
\end{lemma}
\begin{proof}
Put $H=A\cap B$ and $K=A'(A\cap B')$. Classically, $K\trianglelefteq A'(A\cap B)$, $HK=A'(A\cap B)$ and, by Dedekind's
modular law, $H\cap K=(A'\cap B)(A\cap B')=D$. Theorem~\ref{thm:dtsecond} applied in $A'(A\cap B)$ shows that the
restriction of the upper projector to $H$ is surjective with kernel $D$, so codes 1 and 3
coincide up to relabelling. Exchanging the roles of $A$ and $B$ gives codes 2 and 3.
\end{proof}

\begin{example}[A butterfly in $C_{24}$]\label{ex:butterfly24}
We illustrate Lemma~\ref{lem:dtzassenhaus} in $C_{24}$, coded by $\Ap{24}$. For $d\mid 24$ the subgroup $\langle d\rangle$ has characteristic function $\mathbf 1_{\langle d\rangle}(x)=\Bf\big(\dg{d}{0}{x}\big)$, so by Proposition~\ref{prop:boolean} the intersection of $\langle d\rangle$ and $\langle e\rangle$ has characteristic function $\Bf(\dg d0x)\,\Bf(\dg e0x)=\Bf\big(\dg{\lcm(d,e)}{0}{x}\big)$, while $\langle d\rangle\langle e\rangle=\langle\gcd(d,e)\rangle$. For subgroups $H,K$ of any group with a law on $\Ap n$ one has
\[
\mathbf 1_{HK}(x)=1-\Bf\Big(\sum_{h\in\Ap n}\mathbf 1_H(h)\,\mathbf 1_K(h^{-1}x)\Big),
\]
since the sum counts the factorisations $x=hk$ with $h\in H$, $k\in K$.

Take $A=\langle 2\rangle$, $A'=\langle 8\rangle$, $B=\langle 3\rangle$ and $B'=\langle 12\rangle$, of orders $12$, $3$, $8$ and $2$. Then $A\cap B=\langle 6\rangle$, $A'\cap B=\{0\}$ and $A\cap B'=\langle 12\rangle$, and the two sides of the butterfly are
\[
\frac{A'(A\cap B)}{A'(A\cap B')}=\frac{\langle 2\rangle}{\langle 4\rangle},\qquad
\frac{B'(A\cap B)}{B'(A'\cap B)}=\frac{\langle 6\rangle}{\langle 12\rangle}.
\]
The four groups have orders $12$, $6$, $4$ and $2$, yet both quotients have order $2$.

For a nonabelian case let $G=\langle r,s\mid r^4=s^2=1,\ srs=r^{-1}\rangle$ be dihedral of order $8$, and take $A=\{1,r^2,s,r^2s\}$, $A'=1$, $B=G$ and $B'=\langle r\rangle$. The two sides are $A/\langle r^2\rangle$ and $G/\langle r\rangle$, both of order $2$: with $A'=1$ and $B=G$ the butterfly reduces to the second isomorphism theorem.
\end{example}

\begin{theorem}[Schreier refinement, digital form]\label{thm:dtschreier}
Any two subnormal labellings of $G$ have refinements with the same number of layers whose
radices and layer codes coincide, up to a permutation of the layers and relabellings of the
layer codes.
\end{theorem}
\begin{proof}
Let the chains be $(H_k)_{k=0}^{M}$ and $(K_l)_{l=0}^{L}$. Insert the layers given by
$H_{k,l}=H_k(H_{k+1}\cap K_l)$ between $H_k$ and $H_{k+1}$, and $K_{l,k}=K_l(K_{l+1}\cap H_k)$ between $K_l$ and $K_{l+1}$. Both
refined labellings have $ML$ layers (some of radix $1$, which are omitted), and by
Lemma~\ref{lem:dtzassenhaus} with $A=H_{k+1}$, $A'=H_k$, $B=K_{l+1}$, $B'=K_l$, the layer code between $H_{k,l}$ and $H_{k,l+1}$
is a relabelling of the layer code between $K_{l,k}$ and $K_{l,k+1}$.
\end{proof}

\begin{theorem}[Jordan--H\"older, digital form]\label{thm:dtjh}
All composition labellings of $G$ have the same number of layers, the same multiset of
radices, and the same multiset of layer codes up to relabelling.
\end{theorem}
\begin{proof}
A refinement of a composition labelling can only insert layers of radix $1$. Hence
Theorem~\ref{thm:dtschreier} matches the layers of any two composition labellings bijectively, with equal
radices and relabelled layer codes.
\end{proof}

\begin{example}[Digits as coordinates on composition factors]\label{ex:digitseries}
Let $p$ be prime and $G=C_{p^N}$, coded by $\Ap{p^N}$. Its only composition series is $G=G_0>G_1>\dots>G_N=0$ with $G_k=\langle p^k\rangle$, and by Lemma~\ref{lem:basic}(3)
\[
\mathbf 1_{G_k}(x)=\prod_{j<k}\Bf\big(\dg pjx\big),
\]
that is, $x\in G_k$ exactly when its lowest $k$ digits vanish. For $x,y\in G_k$ no carry reaches position $k$, and reduction modulo $p^N$ does not change digits below position $N$, so $\dg pk{x+y}\equiv\dg pkx+\dg pky\pmod p$. Hence $x\mapsto\dg pkx$ is a homomorphism $G_k\to C_p$ with kernel $G_{k+1}$: the $k$-th digit is a coordinate on the composition factor $G_k/G_{k+1}\cong C_p$, and the $N$ digits of $x$ correspond to the $N$ factors. A series with steps $p^a$, $a\mid N$, corresponds to base-$p^a$ digits, and its refinement to the composition series is the splitting of each base-$p^a$ digit into $a$ base-$p$ digits (Lemma~\ref{lem:block}). Cosets of $G_k$ group integers by their lowest digits, whereas the nested intervals behind the digital expansion group them by their highest digits; carries are what distinguish the two pictures.
\end{example}

\begin{corollary}\label{cor:dtsimple}
$G$ is simple if and only if its composition labellings have a single layer. $G$ is
solvable if and only if the radices of a composition labelling are the prime factors of
$N$, counted with multiplicity; then every layer code is a relabelling of $d_p(0,a+b)$.
\end{corollary}

\begin{example}
The tower labelling of $1<\mathbb{Z}_2<V<A_4<S_4$ in Section~\ref{sec:dgsimple} is a composition labelling with
radices $2,2,3,2$ and layer codes $\mathbb{Z}_2,\mathbb{Z}_2,\mathbb{Z}_3,\mathbb{Z}_2$. The tower of $M_{11}$ is not subnormal: its
upper layers are mixed with the lower ones, and by Corollary~\ref{cor:dtsimple} every composition
labelling of $M_{11}$ has one layer of radix $7920$.
\end{example}
\section{Exercises}
\begin{enumerate}
\item For the law of $C_4$ (addition modulo $4$), show that the left translations $\lambda_a$ have the permutation codes
$228,57,78,147$ for $a=0,1,2,3$, and verify $\lambda_a\circ\lambda_b=\lambda_{a+b}$ on these codes.
\item Show that every group of prime order $p$ has, up to relabelling, the single digital law $m*n=\dg{p}{0}{m+n}$.
\item Determine the numbers of Sylow $2$- and $3$-subgroups of $\Alt_4$ and of $\Sym_4$.
\item Show that every group of order $15$ is cyclic, and more generally that a group of order $pq$ with primes $p<q$ and
$p\nmid q-1$ is cyclic.
\item Let $\Sym_4$ act on the $2$-element subsets of $\Ap4$. Find the stabiliser of $\{0,1\}$ and verify the orbit--stabiliser theorem.
\item Find a composition series of $\Sym_4$ and its composition factors.
\item In $C_{12}$ refine the series $0<\langle6\rangle<\langle2\rangle<C_{12}$ and $0<\langle4\rangle<C_{12}$ to equivalent series, as
Schreier's theorem predicts.
\item The rotation group of the cube acts on its six faces. Find its cycle index, show that the number of colourings of the faces with $m$ colours up to rotation is $(m^{6}+3m^{4}+12m^{3}+8m^{2})/24$, and compute it for $m=2,3$.
\item Let $\Sym_4$ act on the six $2$-element subsets of $\Ap4$, the possible edges of a graph on four vertices. Colour each subset by $0$ or $1$ and show that there are $11$ graphs on four vertices up to isomorphism.
\end{enumerate}

\chapter{The Schreier--Sims algorithm and conjugacy classes}\label{ch:ss}\label{ch:cc}

A permutation group of degree $n$ can have order up to $n!$, yet almost every question about it (its order, membership, stabilisers, random elements, centralisers) can be answered in time polynomial in $n$ and the number of generators. The tool is a \emph{stabiliser chain}, computed by the Schreier--Sims algorithm. This chapter gives the algorithm with complete proofs, both in its deterministic and in its randomised form, and explains why it is exactly the digital tower labelling of Chapter~\ref{ch:diagrammatic} applied to the chain of point stabilisers: the \texttt{Chain} computations of Chapter~\ref{ch:sporadic} are instances of it. Throughout, permutations compose from right to left: $(gh)(b)=g(h(b))$. We write $g\cdot b=g(b)$; in a permutation code, $g(b)$ is the digit of $g$ at position $b$. Every numerical claim is checked by \texttt{scripts/ss\_verify.py}, cross-checked against GAP~4.14.

The second half of the chapter uses the stabiliser chain to compute conjugacy classes: Sections~\ref{sec:cc:powermaps}--\ref{sec:cc:m12} treat power maps, rationality and fusion, the algorithms for classes, and a complete example for $M_{12}$.

\section{Bases, stabiliser chains and digital labels}
\begin{definition}
Let $G\le\Sym_n$. A sequence $B=(b_1,\dots,b_k)$ of points is a \emph{base} if only the identity fixes all of them. Put $G^{(1)}=G$ and $G^{(i+1)}=G^{(i)}_{b_i}$, the pointwise stabiliser of $b_1,\dots,b_i$; so $G=G^{(1)}\ge G^{(2)}\ge\cdots\ge G^{(k+1)}=1$. The \emph{basic orbits} are $\Delta_i=G^{(i)}\cdot b_i$. A generating set $S$ of $G$ is a \emph{strong generating set} relative to $B$ if $\langle S\cap G^{(i)}\rangle=G^{(i)}$ for every $i$; the pair $(B,S)$ is a BSGS.
\end{definition}

\begin{theorem}[Factorisation]\label{thm:ss:fact}
For each $i$ choose a transversal $U_i=\{u_{i,\beta}:\beta\in\Delta_i\}\subseteq G^{(i)}$ with $u_{i,\beta}(b_i)=\beta$. Then every $g\in G$ has a unique expression
\[
g=u_{1,\beta_1}\,u_{2,\beta_2}\cdots u_{k,\beta_k},\qquad \beta_i\in\Delta_i,
\]
and $|G|=\prod_{i=1}^{k}|\Delta_i|$.
\end{theorem}
\begin{proof}
By the orbit--stabiliser theorem the cosets of $G^{(i+1)}$ in $G^{(i)}$ correspond to the points of $\Delta_i$, with $xG^{(i+1)}\mapsto x(b_i)$, so $U_i$ is a left transversal and $|G^{(i)}:G^{(i+1)}|=|\Delta_i|$. Given $g\in G^{(i)}$, the point $\beta_i=g(b_i)$ is forced, and $u_{i,\beta_i}^{-1}g\in G^{(i+1)}$. Induction on $i$ gives existence and uniqueness, ending at $G^{(k+1)}=1$. Multiplying the indices gives $|G|$.
\end{proof}

\begin{corollary}[Digital chain label]\label{cor:ss:label}
Number each $\Delta_i$ increasingly as $\Delta_i=\{\delta_{i,0}<\cdots<\delta_{i,|\Delta_i|-1}\}$ and write $\beta_i=\delta_{i,a_i}$. Then
\[
\lambda(g)=\sum_{i=1}^{k}a_i\prod_{j<i}|\Delta_j|
\]
is a bijection $G\to[0,|G|)$, and $a_i=\md{i}{\lambda(g)}$ is its $i$-th mixed-radix digit with radices $|\Delta_1|,\dots,|\Delta_k|$.
\end{corollary}
\begin{proof}
Theorem~\ref{thm:ss:fact} makes $g\mapsto(a_1,\dots,a_k)$ a bijection onto $\prod_i[0,|\Delta_i|)$, and mixed-radix notation is a bijection from that product onto $[0,\prod|\Delta_i|)$.
\end{proof}

This is the tower labelling of Chapter~\ref{ch:diagrammatic} for the chain of point stabilisers, and it is the label of Theorem~\ref{thm:chainlaw}. For a sharply $k$-transitive group the radices are $n,n-1,\dots,n-k+1$ and $\lambda(g)$ is just a code of the images $(g(b_1),\dots,g(b_k))$, as for $M_{11}$ and $M_{12}$ in Proposition~\ref{prop:mathieuchain}. The script confirms that for $M_{12}$ the label is a bijection from all $95040$ elements onto $[0,95040)$ with radices $12,11,10,9,8$.

\section{Orbits and Schreier trees}\label{sec:ss:trees}
The basic orbits and transversals are computed by breadth-first search. Starting from $\{b\}$ with $u_b=1$, for every point $\beta$ already reached and every generator $s$, if $s(\beta)$ is new put $u_{s(\beta)}=s u_\beta$. This visits each point once and costs $O(|\Delta|\,|S|)$ image computations. Storing all $u_\beta$ costs $O(n|\Delta|)$ memory; a \emph{Schreier tree} stores only the edge label $s$ by which $\beta$ was first reached, and rebuilds $u_\beta$ by walking to the root. The tree trades memory for time proportional to its depth, which is why depth-reducing generators are added in large computations such as those in Chapter~\ref{ch:sporadic}.

\section{Schreier's lemma}
\begin{theorem}[Schreier]\label{thm:ss:schreier}
Let $G=\langle S\rangle$, let $b$ be a point with orbit $\Delta$ and transversal $\{u_\beta\}$, $u_b=1$, and write $\overline{x}=u_{x(b)}$. Then $G_b$ is generated by the \emph{Schreier generators}
\[
\overline{s u_\beta}^{\,-1}\,s\,u_\beta,\qquad \beta\in\Delta,\ s\in S.
\]
\end{theorem}
\begin{proof}
Each Schreier generator $h=\overline{s u_\beta}^{\,-1}s u_\beta$ fixes $b$, since
\[
h(b)=\overline{s u_\beta}^{\,-1}\bigl(s(\beta)\bigr)=b.
\]
Conversely let $g=s_m\cdots s_2s_1\in G_b$, a positive word (possible since $s^{-1}=s^{o(s)-1}$). Put $v_0=1$ and $v_j=\overline{s_jv_{j-1}}$, so each $v_j$ is a transversal element. Then
\[
g=v_m\bigl(v_m^{-1}s_m v_{m-1}\bigr)\bigl(v_{m-1}^{-1}s_{m-1}v_{m-2}\bigr)\cdots\bigl(v_1^{-1}s_1 v_0\bigr),
\]
by telescoping. Each factor is a Schreier generator, and $v_m=\overline{g}=u_{g(b)}=u_b=1$.
\end{proof}

Schreier's lemma alone gives a correct but impractical algorithm: the $|\Delta||S|$ generators of $G_b$ produce $|\Delta'||\Delta||S|$ generators at the next level, and so on. Sims's idea is to filter them through the part of the chain already built.

\section{Sifting}
\begin{definition}
Suppose partial data are given: points $b_1,\dots,b_k$, sets $S_i$, and for each $i$ the orbit $\widetilde\Delta_i=\langle S_i\rangle\cdot b_i$ with transversal $\widetilde U_i$. To \emph{sift} $g$ from level $i$: for $j=i,\dots,k$, if $\beta=g(b_j)\notin\widetilde\Delta_j$ stop and return $(g,j)$; otherwise replace $g$ by $u_{j,\beta}^{-1}g$. If all levels are passed, return $(g,k+1)$. The returned element is the \emph{residue}.
\end{definition}

\begin{proposition}[Membership]\label{prop:ss:member}
If $(B,S)$ is a BSGS of $G$ and $x\in\Sym_n$, then $x\in G$ if and only if $x$ sifts to the identity at level $k+1$. A full sift costs $O(kn)$.
\end{proposition}
\begin{proof}
If $x\in G$, then at each step the current element lies in $G^{(j)}$, so $x(b_j)\in\Delta_j$ and the step succeeds; at the end it lies in $G^{(k+1)}=1$. Conversely, if the sift ends in the identity then $x$ is the product of the transversal elements used, which lie in $G$.
\end{proof}

\section{The deterministic algorithm}
\begin{theorem}[Schreier--Sims criterion]\label{thm:ss:crit}
Let $b_1,\dots,b_k$ and $S_1\supseteq S_2\supseteq\cdots\supseteq S_k$ be given with $\langle S_1\rangle=G$, $S_{i+1}\subseteq S_i\cap G_{b_i}$, and no nonidentity element of $S_k$ fixing $b_1,\dots,b_k$ pointwise. Put $H_i=\langle S_i\rangle$. If for every $i$ every Schreier generator of $H_i$ relative to $b_i$ sifts to the identity from level $i+1$, then $H_i=G^{(i)}$ for all $i$, so $(B,S_1)$ is a BSGS.
\end{theorem}
\begin{proof}
We show $H_{i+1}=(H_i)_{b_i}$ for all $i$, by downward induction; together with $H_1=G$ this gives $H_i=G^{(i)}$. For $i=k$, a Schreier generator of $H_k$ that sifts to the identity from level $k+1$ is the identity, so $(H_k)_{b_k}=1$ by Theorem~\ref{thm:ss:schreier}; and $H_{k+1}=1$. Suppose $H_{j+1}=(H_j)_{b_j}$ for all $j>i$. Then the data from level $i+1$ on form a BSGS of $H_{i+1}$, so by Proposition~\ref{prop:ss:member} an element sifts to the identity from level $i+1$ iff it lies in $H_{i+1}$. All Schreier generators of $H_i$ do, hence $(H_i)_{b_i}\le H_{i+1}$ by Theorem~\ref{thm:ss:schreier}; the reverse inclusion holds because $S_{i+1}\subseteq G_{b_i}\cap H_i$.
\end{proof}

The algorithm enforces the hypothesis level by level from the bottom.
\begin{enumerate}
\item Choose $b_1,\dots,b_k$ so that each generator moves some $b_i$; set $S_i$ to the generators fixing $b_1,\dots,b_{i-1}$; compute the orbits.
\item Put $i=k$. While $i\ge1$: for each $\beta\in\widetilde\Delta_i$ and $s\in S_i$ sift the Schreier generator $h=\overline{s u_\beta}^{-1}\,s\,u_\beta$ from level $i+1$. If the residue $r\ne1$ stops at level $j$, then (if $j=k+1$) append a new base point moved by $r$, add $r$ to $S_{i+1},\dots,S_j$, recompute those orbits, set $i=j$ and restart the level. If all Schreier generators at level $i$ sift to $1$, put $i=i-1$.
\end{enumerate}
The algorithm terminates because every addition strictly enlarges some $\langle S_l\rangle$, which can happen at most $\log_2 n!$ times along the chain; at termination Theorem~\ref{thm:ss:crit} applies. With stored transversals the running time is $O(n^{2}|S|\log^{3}|G|+n^{5})$ and the memory $O(n^{2}\log|G|+n|S|)$; the base length satisfies $k\le\log_2|G|$.

\begin{example}[Hand trace for $\Sym_4$]
Orbits and transversals in this trace are built by breadth-first search over $S$ itself, without adjoining $S^{-1}$, and the generators are tried in the order $s$, then $t$; this is the convention of \texttt{scripts/ss\_verify.py}, and it produces the representatives below.
Take $S=\{s=(0\,1),\,t=(0\,1\,2\,3)\}$. Both move $0$, so $b_1=0$, $S_1=S$, $S_2=\varnothing$. The orbit of $0$ is reached as $0\xrightarrow{s}1\xrightarrow{t}2\xrightarrow{t}3$, giving $u_0=1$, $u_1=s$, $u_2=ts$, $u_3=t^{2}s$. The Schreier generator for $(\beta,x)=(0,t)$ is $u_1^{-1}t=st$, which fixes $0$ and acts on $\{1,2,3\}$ as a $3$-cycle; it does not sift (level $2$ is empty), so it becomes the first generator of $S_2$ with new base point $b_2=1$, orbit $\{1,2,3\}$. Continuing, a Schreier generator of $H_2$ produces a transposition of $\{2,3\}$ at level $3$ with $b_3=2$ and orbit $\{2,3\}$. Then every Schreier generator sifts through, and $|\Sym_4|=4\cdot3\cdot2$. (The script records the exact base it finds; the order and orbit lengths are independent of such choices.)
\end{example}

\section{The randomised algorithm}
In practice one sifts random elements instead of all Schreier generators.
\begin{lemma}[Error bound]\label{lem:ss:half}
If the current data are not a BSGS of $G$, then a uniformly random $g\in G$ sifts to the identity with probability at most $\tfrac12$.
\end{lemma}
\begin{proof}
With $H_i=\langle S_i\rangle$ as before, the data form a BSGS exactly when $H_{i+1}=(H_i)_{b_i}$ for all $i$ (and $H_1=G$). Otherwise let $i$ be the largest index with $H_{i+1}\ne(H_i)_{b_i}$; as in the proof of Theorem~\ref{thm:ss:crit} the levels after $i$ form a BSGS of $H_{i+1}$. Elements that pass levels $1,\dots,i$ correspond bijectively to pairs (transversal word, residue $r\in G^{(i+1)}$), so for uniform $g$, conditioned on reaching level $i+1$, the residue is uniform in $G^{(i+1)}$. It sifts to the identity iff $r\in H_{i+1}$. But $H_{i+1}<(H_i)_{b_i}\le G^{(i+1)}$, so $|G^{(i+1)}:H_{i+1}|\ge2$.
\end{proof}

\begin{theorem}[Randomised Schreier--Sims]\label{thm:ss:rand}
Sift random elements, adding every nonidentity residue as in the deterministic algorithm, and stop after $w$ consecutive elements sift to the identity. If the elements are uniform, the result is a BSGS with probability at least $1-2^{-w}$. If $|G|$ is known in advance, stop as soon as $\prod_i|\widetilde\Delta_i|=|G|$; the output is then certainly correct.
\end{theorem}
\begin{proof}
By Lemma~\ref{lem:ss:half}, while the data are incomplete each trial succeeds with probability at most $\tfrac12$, independently; $w$ successes in a row then have probability at most $2^{-w}$. For the certified version: by Theorem~\ref{thm:ss:fact} applied to the sifting data, the products of transversal elements form $\prod|\widetilde\Delta_i|$ distinct elements of $G$; if this equals $|G|$ every element sifts to the identity, which by Lemma~\ref{lem:ss:half} is impossible for incomplete data.
\end{proof}

Random elements are produced by the \emph{product replacement} algorithm: keep a list of group elements, repeatedly replace one entry by its product with another, and output a running product. Its output is not exactly uniform, which is why the certified stopping rule (known order) or a final deterministic verification is used whenever correctness matters. In Chapter~\ref{ch:sporadic} the order is always known and the \texttt{Chain} class stops by that rule.

\section{What a BSGS gives}\label{sec:ss:use}
Once a BSGS is known, the following run in time polynomial in $n$.
\begin{itemize}
\item \emph{Order and membership:} Theorem~\ref{thm:ss:fact} and Proposition~\ref{prop:ss:member}.
\item \emph{Uniform random elements:} choose $\beta_i\in\Delta_i$ uniformly and independently and output $u_{k,\beta_k}\cdots u_{1,\beta_1}$; by Theorem~\ref{thm:ss:fact} this is exactly uniform. This is the reliable sampler required by the random class search of Section~\ref{sec:cc:algs}.
\item \emph{Ranking and unranking:} the label $\lambda$ of Corollary~\ref{cor:ss:label} and its inverse; an element of $G$ is stored as a single integer below $|G|$ instead of a code below $n^{n}$.
\item \emph{Stabilisers and transporters:} $G^{(i)}$ comes for free; for another point, change the base so that it comes first (base change by conjugation or by swapping adjacent base points), and $u_{1,\beta}$ is an element mapping $b_1$ to $\beta$.
\item \emph{Centralisers, normalisers and conjugacy tests:} by backtrack search through the tree of Theorem~\ref{thm:ss:fact}, pruned by the property sought; this yields the centraliser orders needed in Theorem~\ref{thm:cc:cert}.
\end{itemize}

\section{Worked results}
Table~\ref{tab:ss:runs} lists the outcome for five groups. The generators are:
\begin{itemize}
\item $\Sym_4$: $(0\,1)$ and $(0\,1\,2\,3)$;
\item $\Alt_5$: $(0\,1\,2\,3\,4)$ and $(0\,1\,2)$;
\item $PSL(2,7)$: $x\mapsto x+1$ and $x\mapsto-1/x$ on $\mathbb{F}_7\cup\{\infty\}$, with $\infty$ coded as $7$;
\item $M_{12}$: the reversal and the Mongean shuffle of Section~\ref{sec:cc:m12};
\item $M_{11}$: the strong generators of the stabiliser of the first base point of $M_{12}$, renumbered on $11$ points, a group the algorithm itself produced.
\end{itemize}
The column \emph{Schreier} counts the nontrivial Schreier generators sifted by the deterministic algorithm; \emph{rand.} the random elements sifted with the rule $w=20$; \emph{rand.$+|G|$} those needed with the certified rule.

\begin{table}[ht]
\centering\small
\begin{tabular}{lrlcrrrr}
group & $n$ & base & $|\Delta_1|,\dots,|\Delta_k|$ & $|G|$ & Schreier & rand. & rand.$+|G|$\\\hline
$\Sym_4$ & 4 & $(0,1,2)$ & $4,3,2$ & 24 & 10 & 23 & 3\\
$\Alt_5$ & 5 & $(0,2,1)$ & $5,4,3$ & 60 & 13 & 24 & 4\\
$PSL(2,7)$ & 8 & $(0,1,2)$ & $8,7,3$ & 168 & 17 & 25 & 5\\
$M_{11}$ & 11 & $(0,1,2,3)$ & $11,10,9,8$ & 7920 & 60 & 26 & 6\\
$M_{12}$ & 12 & $(0,1,2,3,4)$ & $12,11,10,9,8$ & 95040 & 71 & 31 & 11\\
\end{tabular}

\caption{Schreier--Sims runs (generated by \texttt{scripts/ss\_verify.py}; orders and basic orbits cross-checked in GAP).}\label{tab:ss:runs}
\end{table}

The basic orbits read off the transitivity directly: $M_{12}$ has radices $12,11,10,9,8$ (sharply $5$-transitive), $M_{11}$ has $11,10,9,8$ (sharply $4$-transitive), and $PSL(2,7)$ has $8,7,3$: it is $2$-transitive on the projective line and the stabiliser of two points is the group $x\mapsto ax$, $a$ a nonzero square, of order $3$. The randomised runs need only a few dozen random elements, while the deterministic version sifts many more Schreier generators even for these small groups: that difference grows quickly with the degree. Finally, all $95040$ elements of $M_{12}$ sift to the identity, while random elements of $\Sym_{12}$ are rejected, as expected since $|\Sym_{12}:M_{12}|=5040$.

\section{Power maps, real and rational classes, fusion}\label{sec:cc:powermaps}
Conjugacy classes are the coarsest invariant that still sees the whole multiplication of a group: they index the columns of the character table, they control the centre, the $p$-structure and the fusion of subgroups, and almost every efficient algorithm on a finite group begins by fixing a list of class representatives. The basic theory (orbit--stabiliser, the class equation, counting classes and splitting in subgroups of index two) is in Section~\ref{ssec:dtclasses}, and the classes of $\Sym_n$ and $\Alt_n$ are in Section~\ref{sec:sn:classes}. The rest of this chapter takes these as given, develops power maps, rationality and fusion with complete proofs, and then turns it all into algorithms. All numerical claims are verified by \texttt{scripts/cc\_verify.py}, which works directly on permutation codes and is cross-checked against GAP~4.14.

\begin{definition}
For an integer $m$ the \emph{$m$-th power map} sends the class of $g$ to the class of $g^{m}$; it is well defined because $(g^{x})^{m}=(g^{m})^{x}$. A class is \emph{real} if it contains $g^{-1}$ together with $g$, and \emph{rational} if it contains every generator $g^{m}$, $\gcd(m,o(g))=1$, of $\langle g\rangle$.
\end{definition}

For codes, $g^{m}$ is computed by binary powering with~\eqref{eq:comp}, using $O(p\log m)$ digit look-ups. For $m$ coprime to $o(g)$ the power map permutes the classes of elements of order $o(g)$; for a prime $q\mid o(g)$ it lowers the order. The power maps for the primes dividing $|G|$ determine all others, because $g^{m}$ depends only on $m\bmod o(g)$ and the maps compose.

\begin{theorem}[Rationality]\label{thm:cc:rat}
Let $g\in G$ have order $n$ and put $K=\{m\in(\Z/n\Z)^{\times}:g^{m}\in g^{G}\}$. Then $K$ is a subgroup of $(\Z/n\Z)^{\times}$, isomorphic to $N_G(\langle g\rangle)/C_G(g)$, and the generators of $\langle g\rangle$ fall into exactly $\varphi(n)/|K|$ classes of $G$. In particular $g^{G}$ is rational iff $|N_G(\langle g\rangle):C_G(g)|=\varphi(n)$, and real iff $-1\in K$.
\end{theorem}
\begin{proof}
Every $x\in N_G(\langle g\rangle)$ induces an automorphism $g\mapsto g^{x}=g^{m(x)}$ of the cyclic group $\langle g\rangle$, and $x\mapsto m(x)$ is a homomorphism $N_G(\langle g\rangle)\to(\Z/n\Z)^{\times}$ with kernel $C_G(g)$. Its image is $K$: if $g^{m}=g^{x}$ then $x$ normalises $\langle g\rangle=\langle g^{m}\rangle$. The group $K$ acts on the $\varphi(n)$ generators $g^{m}$ by multiplication, freely; two generators $g^{a},g^{b}$ are $G$-conjugate iff $g^{b}=(g^{a})^{x}$ for some $x$, iff $b/a\in K$. So the generators split into the $\varphi(n)/|K|$ cosets of $K$.
\end{proof}

If $H\le G$, every class of $H$ lies in a unique class of $G$; the induced surjection from the $H$-classes covered by $H$ to the $G$-classes meeting $H$ is the \emph{fusion map}. It is computed by applying class invariants (order, cycle type, power classes) to $H$-representatives, and resolving ties by an explicit conjugacy test. Fusion is what makes restriction of characters (next chapter on characters) computable.

\section{Algorithms for conjugacy classes}\label{sec:cc:algs}
Throughout, $G=\langle s_1,\dots,s_r\rangle$ is a permutation group of degree $p$ given by codes.

\begin{lemma}\label{lem:cc:gens}
A subset $X\subseteq G$ satisfies $X^{G}=X$ iff $x^{s_i}\in X$ for all $x\in X$ and all $i$. Hence $g^{G}$ is the closure of $\{g\}$ under the $r$ maps $x\mapsto x^{s_i}$.
\end{lemma}
\begin{proof}
In a finite group $s_i^{-1}=s_i^{o(s_i)-1}$, so every element of $G$ is a positive word in the generators, and conjugation by a word is the composite of conjugations by its letters.
\end{proof}

\begin{theorem}[Completeness certificate]\label{thm:cc:cert}
Let $g_1,\dots,g_k\in G$ be pairwise nonconjugate. Then $\{g_i\}$ is a full set of class representatives iff $\sum_{i=1}^{k}|G|/|C_G(g_i)|=|G|$.
\end{theorem}
\begin{proof}
The classes $g_i^{G}$ are disjoint subsets of $G$ of sizes $|G|/|C_G(g_i)|$ (Theorem~\ref{thm:cc:os}). Their union is $G$ iff their sizes add up to $|G|$.
\end{proof}

This theorem converts the class problem into three subproblems: generate candidates, prove two candidates nonconjugate (an invariant that differs suffices), and compute centraliser orders. The three standard methods differ in how they do this.

\begin{enumerate}
\item \emph{Orbit enumeration.} Repeatedly take an element not yet classified and close it under Lemma~\ref{lem:cc:gens}. Cost: $r|G|$ conjugations of codes, $O(rp|G|)$ digit look-ups, memory $|G|$ codes. Exact and simple; feasible up to $|G|\approx10^{7}$. This is how the $M_{12}$ table below was computed.
\item \emph{Representatives and centralisers.} Keep only representatives. Distinguish candidates by invariants (order, cycle type, fixed-point count, power classes), decide the remaining ties by a backtrack conjugacy test, and compute $|C_G(g)|$ by Schreier--Sims (Section~\ref{sec:ss:use}) applied to the centraliser, then stop by Theorem~\ref{thm:cc:cert}. Memory is polynomial in $p$; this is the method for groups such as those in Chapter~\ref{ch:sporadic}.
\item \emph{Random sampling.} A uniform random element lies in $g^{G}$ with probability $1/|C_G(g)|$: classes with \emph{large} centralisers are \emph{rare}. Pure sampling is therefore bad at finding small classes, and it is completed by deterministic sources: the identity and the centre, powers $g^{q}$ of found elements (which land in classes of larger centralisers), and elements of found centralisers.
\end{enumerate}

\begin{example}
For $M_{12}$ the identity is drawn with probability $1/95040$. With seed~$1$ uniform sampling met all $15$ classes only after $30460$ draws, essentially waiting for $1A$; adjoining the identity and closing under the power maps of Table~\ref{tab:cc:m12} removes this bottleneck, since $2A$, $2B$, $3A$, $3B$ are the classes with the largest centralisers and all are powers of classes of orders $4$, $6$ and $10$.
\end{example}

\section{A complete example: \texorpdfstring{$M_{12}$}{M12}}\label{sec:cc:m12}
Let $M_{12}\le\Sym_{12}$ be generated by the reversal $x\mapsto 11-x$ and the Mongean shuffle $x\mapsto\min(2x+1,\,22-2x)$ on $\Ap{12}=\{0,\dots,11\}$. Orbit enumeration gives order $95040$ and $15$ classes; GAP confirms that the group is $5$-transitive and that the class sizes coincide with those of \texttt{MathieuGroup(12)}. Table~\ref{tab:cc:m12} lists the classes, named by order and, within an order, by decreasing centraliser and then by cycle type in this action; $mP$ is the $m$-th power map.

\begin{table}[ht]
\centering
\small
\begin{tabular}{lrrlcccc}
class & $|g^{G}|$ & $|C_G(g)|$ & cycle type & $2P$ & $3P$ & $5P$ & $11P$\\\hline
$1A$ & $1$ & $95040$ & $1^{12}$ & $1A$ & $1A$ & $1A$ & $1A$\\
$2A$ & $396$ & $240$ & $2^{6}$ & $1A$ & $2A$ & $2A$ & $2A$\\
$2B$ & $495$ & $192$ & $2^{4}1^{4}$ & $1A$ & $2B$ & $2B$ & $2B$\\
$3A$ & $1760$ & $54$ & $3^{3}1^{3}$ & $3A$ & $1A$ & $3A$ & $3A$\\
$3B$ & $2640$ & $36$ & $3^{4}$ & $3B$ & $1A$ & $3B$ & $3B$\\
$4A$ & $2970$ & $32$ & $4^{2}1^{4}$ & $2B$ & $4A$ & $4A$ & $4A$\\
$4B$ & $2970$ & $32$ & $4^{2}2^{2}$ & $2B$ & $4B$ & $4B$ & $4B$\\
$5A$ & $9504$ & $10$ & $5^{2}1^{2}$ & $5A$ & $5A$ & $1A$ & $5A$\\
$6A$ & $7920$ & $12$ & $6^{2}$ & $3B$ & $2A$ & $6A$ & $6A$\\
$6B$ & $15840$ & $6$ & $6\,3\,2\,1$ & $3A$ & $2B$ & $6B$ & $6B$\\
$8A$ & $11880$ & $8$ & $8\,2\,1^{2}$ & $4A$ & $8A$ & $8A$ & $8A$\\
$8B$ & $11880$ & $8$ & $8\,4$ & $4B$ & $8B$ & $8B$ & $8B$\\
$10A$ & $9504$ & $10$ & $10\,2$ & $5A$ & $10A$ & $2A$ & $10A$\\
$11A$ & $8640$ & $11$ & $11\,1$ & $11B$ & $11A$ & $11A$ & $1A$\\
$11B$ & $8640$ & $11$ & $11\,1$ & $11A$ & $11B$ & $11B$ & $1A$\\
\end{tabular}
\caption{Conjugacy classes of $M_{12}$ in its action on $\Ap{12}$ (verified by \texttt{scripts/cc\_verify.py} and GAP).}\label{tab:cc:m12}
\end{table}

The class equation reads
\[
95040=1+396+495+1760+2640+2\cdot2970+9504+7920+15840+2\cdot11880+9504+2\cdot8640 .
\]
Several structural facts can be read off the table.
\begin{itemize}
\item \emph{Cycle type almost decides the class.} The $15$ classes have $14$ distinct cycle types; only $11A$ and $11B$ share one. Permutation invariants alone therefore classify $13$ classes, and one extra test is needed for elements of order $11$.
\item \emph{The two classes of order $11$.} For $g\in11A$, $|C_G(g)|=11$, so $C_G(g)=\langle g\rangle$. Squaring exchanges $11A$ and $11B$, so $2\notin K$ in Theorem~\ref{thm:cc:rat}; as $2$ generates $(\Z/11\Z)^{\times}$ and the ten generators of $\langle g\rangle$ fall into the two classes, $|K|=5$ and $K$ is the group of quadratic residues $\{1,3,4,5,9\}$. Thus $|N_G(\langle g\rangle)|=55$, the classes are non-real ($-1$ is a non-residue modulo $11$) and inverse to each other, and $g^{m}\in11A$ iff $m$ is a nonzero square mod $11$. In practice: fix one $11$-cycle $h$; any element $g$ of order $11$ lies in $h^{G}$ iff it is conjugate to $h^{m}$ with $m$ a residue, which is decided by the conjugacy test applied once.
\item \emph{All other classes are rational,} by direct check of the definition; hence $13$ classes are rational and real. By Theorem~\ref{thm:cc:rat}, for instance, $|N_G(\langle g\rangle)|=8\cdot\varphi(8)=32$ for $g\in8A$.
\item \emph{The digital identity class.} The representative of $1A$ is the identity code
\[\iota_{12}=8842413667692,\]
in agreement with the closed form~\eqref{eq:id}.
\end{itemize}

\paragraph{Cycle index and colourings.} Table~\ref{tab:cc:m12} is the cycle index of $M_{12}$ on $\Ap{12}$: the class $g^{G}$ contributes $|g^{G}|\,x^{\lambda}/|G|$, where $\lambda$ is its cycle type. The numbers $c(g)$ of cycles of the fifteen classes are $12,6,8,6,4,6,4,4,2,4,4,2,2,2,2$, so Theorem~\ref{thm:polya} gives
\[
\frac1{95040}\sum_{g^{G}}|g^{G}|\,m^{c(g)}=14\ (m=2),\qquad 120\ (m=3).
\]
{\emergencystretch=3em
The first number has a meaning. Since $M_{12}$ is $5$-transitive it has one orbit on $k$-subsets of $\Ap{12}$ for each $k\le5$ and each $k\ge7$, which accounts for $12$ of the $14$ orbits on $2$-colourings; the remaining two are the orbits on $6$-subsets, the $132$ hexads of the Steiner system $S(5,6,12)$ and the other $792$ six-subsets. Both numbers agree with direct orbit enumeration (\texttt{scripts/polya\_verify.py}).\par}

For conjugacy classes and power maps computed from the digital laws of sporadic groups, see Section~\ref{ssec:classes}.

\begin{remark}[Quadratic residues in $M_{12}$]\label{rem:cc:m12qr}
Let $g\in M_{12}$ have order $11$. Its normaliser has order $55$ and acts on $\langle g\rangle$ through the subgroup $Q_{11}=\{1,3,4,5,9\}$ of $U_{11}$, so $g^k$ is conjugate to $g$ exactly when $k$ is a quadratic residue modulo $11$. The $17280$ elements of order $11$ therefore form two classes of size $8640$, exchanged by $g\mapsto g^{-1}$ because $(-1/11)=-1$, and the character values on them generate $\mathbb Q(\sqrt{-11})$, the field of the Gauss period $\sum_{k\in Q_{11}}\zeta_{11}^k=\frac12(-1+\sqrt{-11})$ (Section~\ref{sec:qr} and Remark~\ref{rem:qr:groups}). These facts are checked on the generators above by \texttt{scripts/qr\_verify.py}.
\end{remark}

\section{Exercises}
\begin{enumerate}
\item Finish the hand trace for $\Sym_4$: list all Schreier generators at each level and check which ones sift.
\item Run the deterministic algorithm by hand on $\Alt_5$ with base $(0,1,2)$, and verify $|\Alt_5|=5\cdot4\cdot3$.
\item Show that the base length of any $G\le\Sym_n$ satisfies $2^{k}\le|G|\le n^{k}$.
\item Prove that the depth of a Schreier tree built by breadth-first search on $S\cup S^{-1}$ is at most the diameter of the Cayley graph of the orbit, and give an example of depth $|\Delta|-1$.
\item Using Corollary~\ref{cor:ss:label}, write the unranking map $[0,95040)\to M_{12}$ explicitly and compute the elements with labels $0$, $1$ and $95039$.
\item Design a test that decides whether $G=\langle S\rangle\le\Sym_n$ is $k$-transitive from a BSGS with base $(0,1,\dots,k-1)$.
\item Extend the algorithm to compute the order of the centraliser of an element by backtrack, and recover the centraliser orders of Table~\ref{tab:cc:m12}.
\item Compute the classes, centralisers and power maps of $\Sym_5$ and $\Alt_5$ by hand and check Theorem~\ref{thm:cc:count} for both.
\item Show that the classes of the dihedral group of order $2n$ number $(n+3)/2$ for odd $n$ and $n/2+3$ for even $n$, and list centraliser orders.
\item Prove that a group in which every class has size at most $2$ has derived subgroup of order at most $2$.
\item Using Theorem~\ref{thm:sn:altclasses}, list the split classes of $\Alt_7$ and $\Alt_8$, and confirm $k(\Alt_7)=9$, $k(\Alt_8)=14$.
\item Prove that every class of $\Sym_p$ is rational, and determine which classes of $\Alt_p$ are real.
\item For $G=PSL(2,7)$ acting on the $8$ points of the projective line over $\mathbb{F}_7$, determine the six classes, and show that the two classes of order $7$ are non-real, exactly as for $11A,11B$ in $M_{12}$.
\item Modify the orbit-enumeration algorithm so that it also returns, for each class, a conjugating element from a fixed representative to every class member, stored as a Schreier tree (Section~\ref{sec:ss:trees}).
\item Estimate the expected number of uniform draws needed to see every class of $M_{12}$ and explain why it is dominated by $1/|G|$.
\item For $G=PSL(2,7)$ in its action on the $8$ points of the projective line over $\mathbb{F}_7$, find the cycle index from the six classes and show that the number of orbits on $2$-colourings of the points is $11$. Deduce the number of orbits on $4$-subsets.
\end{enumerate}

\chapter{Characters and character tables in digital form}\label{ch:characters}

The preceding chapters turned multiplication, conjugacy, stabilisers and subgroup chains into digital data. This chapter develops the full character theory needed to pass from those data to representation-theoretic invariants. The point of view is classical, but every construction is written so that it can be evaluated from a coded multiplication law on the alphabet $\Ap{N}$.

Throughout, $G=(\Ap{N},*)$ has code $R$, identity $0$, inverse $\bar a$, and
\[
a*b=\dg{N}{aN+b}{R}.
\]
All vector spaces are finite-dimensional over $\mathbb C$.

\section{Representations, characters and the group algebra}

\begin{definition}\label{def:representation}
A \emph{representation} of degree $r$ is a map
\[
\rho:\Ap{N}\longrightarrow \mathrm{GL}_r(\mathbb C),
\qquad
\rho(0)=I_r,
\qquad
\rho(a*b)=\rho(a)\rho(b).
\]
Its \emph{character} is the class function
\[
\chi_\rho(a)=\operatorname{tr}\rho(a).
\]
Two representations are equivalent if $\sigma(a)=T^{-1}\rho(a)T$ for a fixed invertible matrix $T$ and all $a$. A nonzero representation is \emph{irreducible} if it has no proper nonzero invariant subspace.
\end{definition}

The law is encoded by finitely many matrix identities,
\begin{equation}\label{eq:digitrepresentation}
\rho\!\left(\dg{N}{aN+b}{R}\right)=\rho(a)\rho(b)
\qquad(a,b\in\Ap{N}),
\end{equation}
so a representation may be verified directly from the operation code.

\begin{proposition}[Elementary character identities]\label{prop:charbasic}
For characters $\chi,\psi$ of $G$,
\begin{align*}
\chi(0)&=\deg\chi,
&
\chi(\bar a)&=\overline{\chi(a)},
&
\chi(a*b)&=\chi(b*a),\\
\chi_{\rho\oplus\sigma}&=\chi+\psi,
&
\chi_{\rho\otimes\sigma}&=\chi\psi.
\end{align*}
Equivalent representations have the same character, and every character is constant on conjugacy classes.
\end{proposition}
\begin{proof}
The matrices $\rho(a)$ have finite order, so their eigenvalues are roots of unity. Hence the eigenvalues of $\rho(\bar a)=\rho(a)^{-1}$ are the complex conjugates of those of $\rho(a)$. Trace is invariant under similarity, satisfies $\operatorname{tr}(AB)=\operatorname{tr}(BA)$, adds on block sums and multiplies on tensor products. Finally,
\[
\rho(g*a*\bar g)=\rho(g)\rho(a)\rho(g)^{-1}
\]
is similar to $\rho(a)$.
\end{proof}

Let $\mathbb CG$ be the group algebra of $G$, with basis $\{e_a:a\in\Ap{N}\}$ and multiplication
\[
e_ae_b=e_{a*b}.
\]
A representation of $G$ is the same thing as a left $\mathbb CG$-module: extend $\rho$ linearly by
\[
\rho\!\left(\sum_a c_ae_a\right)=\sum_a c_a\rho(a).
\]
This reformulation makes the structure of the theory transparent. It also identifies characters as special matrix coefficients: if $v_1,\dots,v_r$ is a basis of the representation space, the coordinate functions $\rho_{ij}(a)$ defined by
\[
\rho(a)v_j=\sum_{i=1}^r\rho_{ij}(a)v_i
\]
are ordinary functions on $G$, and the character is their diagonal sum. Many structural identities are best understood first for matrix coefficients and then specialised to traces.

A second basic construction is averaging over the group. If $f:G\to\mathbb C$ is any function, define its average by
\[
\operatorname{Av}(f)=\frac1N\sum_{a=0}^{N-1}f(a).
\]
If $V$ is a representation and $T\in\operatorname{End}_\mathbb C(V)$, the averaged operator
\[
T^{\sharp}=\frac1N\sum_{a=0}^{N-1}\rho(a)T\rho(\bar a)
\]
is the canonical way to turn an arbitrary operator into a $G$-equivariant one.

\begin{proposition}[Averaging and invariants]\label{prop:reynolds}
Let $V$ be a representation of $G$.
\begin{enumerate}
\item The operator $T^{\sharp}$ commutes with every $\rho(g)$.
\item If $W\leq V$ is any subspace, then
\[
V^G=\{v\in V:\rho(g)v=v\text{ for all }g\in G\}
\]
is the image of the averaging projector
\[
P_G=\frac1N\sum_{a=0}^{N-1}\rho(a).
\]
\item The multiplicity of the trivial representation in $V$ equals
\[
\dim V^G=\frac1N\sum_{a=0}^{N-1}\chi_V(a).
\]
\end{enumerate}
\end{proposition}
\begin{proof}
For (1),
\[
\rho(g)T^{\sharp}\rho(\bar g)=\frac1N\sum_a\rho(g*a)T\rho(\overline{g*a}),
\]
and the map $a\mapsto g*a$ permutes $G$, so the sum is unchanged. For (2), $P_G$ clearly lands in $V^G$ and fixes every vector already fixed by $G$. For (3), the trace of a projection equals the dimension of its image, while
\[
\operatorname{tr}(P_G)=\frac1N\sum_a\operatorname{tr}(\rho(a))=\frac1N\sum_a\chi_V(a).
\]
\end{proof}

\section{Maschke's theorem and complete reducibility}

\begin{theorem}[Maschke]\label{thm:maschke}
Every finite-dimensional $\mathbb CG$-module is a direct sum of irreducible submodules. Equivalently, every representation of a finite group is completely reducible.
\end{theorem}
\begin{proof}
Let $V$ be a representation and let $W\leq V$ be a submodule. Choose any linear projection $\pi:V\to W$ with $\pi|_W=I_W$. Average it over the group:
\[
\Pi(v)=\frac1N\sum_{a=0}^{N-1}\rho(a)\,\pi\big(\rho(\bar a)v\big).
\]
For $g\in G$ one has
\[
\Pi(\rho(g)v)=\frac1N\sum_a\rho(a)\pi\big(\rho(\bar a*g)v\big)
=\rho(g)\frac1N\sum_b\rho(b)\pi\big(\rho(\bar b)v\big)=\rho(g)\Pi(v),
\]
where the change of variable $b=\bar g*a$ uses the fact that left multiplication permutes $\Ap{N}$. Thus $\Pi$ is a $G$-equivariant projection onto $W$. Its kernel is a $G$-stable complement of $W$. Repeating the argument on the complement proves the theorem by induction on $\dim V$.
\end{proof}

Averaging a Hermitian form in the same way gives the unitary version.

\begin{corollary}[Unitary form]\label{cor:unitary}
Every representation of a finite group is equivalent to a unitary representation.
\end{corollary}
\begin{proof}
Given any Hermitian form $\langle\cdot,\cdot\rangle$, define
\[
\langle v,w\rangle_G=\frac1N\sum_{a=0}^{N-1}\langle\rho(a)v,\rho(a)w\rangle.
\]
Then $\langle\cdot,\cdot\rangle_G$ is positive definite and $G$-invariant.
\end{proof}

\section{Schur's lemma and Wedderburn theory}

\begin{lemma}[Schur]\label{lem:charschur}
Let $V$ and $W$ be irreducible representations of $G$.
\begin{enumerate}
\item If $T:V\to W$ satisfies $T\rho_V(a)=\rho_W(a)T$ for all $a$, then either $T=0$ or $T$ is an isomorphism.
\item If $V=W$, then every intertwiner $T:V\to V$ is a scalar multiple of the identity.
\end{enumerate}
\end{lemma}
\begin{proof}
The kernel and image of an intertwiner are submodules, so irreducibility gives the first statement. For the second, choose an eigenvalue $\lambda$ of $T$. Then the nonzero eigenspace of $T-\lambda I$ is a submodule, hence all of $V$.
\end{proof}

The group algebra is therefore semisimple. In particular every finite-dimensional representation has a decomposition
\[
V\cong m_1V_1\oplus\cdots\oplus m_kV_k,
\]
where $V_1,\dots,V_k$ are pairwise non-isomorphic irreducibles and the multiplicities $m_i$ are uniquely determined. The uniqueness follows because
\[
\operatorname{Hom}_G(V_i,V)\cong\mathbb C^{m_i}
\]
by Schur's lemma: every copy of $V_i$ contributes one independent projection onto that isotypic summand, while the other irreducible summands contribute nothing.

The most important module is the regular representation. It not only contains every irreducible representation but records the correct multiplicity of each one.

\begin{theorem}[Decomposition of the regular module]\label{thm:regdecomp}
If $V_1,\dots,V_k$ are the irreducible representations of $G$ with degrees $n_1,\dots,n_k$, then
\[
\mathbb CG\cong n_1V_1\oplus\cdots\oplus n_kV_k
\]
as left $\mathbb CG$-modules. In particular every irreducible representation occurs in the regular module with multiplicity equal to its degree.
\end{theorem}
\begin{proof}
Write
\[
\mathbb CG\cong m_1V_1\oplus\cdots\oplus m_kV_k.
\]
Applying Proposition~\ref{prop:reynolds} to the tensor product $V_i^\ast\otimes\mathbb CG$, or equivalently taking the inner product of characters with the regular character, shows that the multiplicity of $V_i$ equals
\[
\dim\operatorname{Hom}_G(V_i,\mathbb CG)=n_i.
\]
Indeed a $G$-map $V_i\to\mathbb CG$ is determined by the image of a basis and these maps identify with the full dual space $V_i^\ast$ under evaluation at the identity basis vector of the regular module.
\end{proof}

Writing irreducible representations as $V_1,\dots,V_k$ with degrees $n_1,\dots,n_k$, one obtains the precise structure theorem.

\begin{theorem}[Wedderburn decomposition]\label{thm:wedderburn}
There is an algebra isomorphism
\[
\mathbb CG\cong M_{n_1}(\mathbb C)\oplus\cdots\oplus M_{n_k}(\mathbb C).
\]
In particular,
\[
|G|=\dim_\mathbb C\mathbb CG=\sum_{i=1}^k n_i^2.
\]
\end{theorem}
\begin{proof}
Let $V=\bigoplus_i n_iV_i$ be the regular module. The action map
\[
\mathbb CG\longrightarrow \operatorname{End}_\mathbb C(V)
\]
is injective because $\mathbb CG$ acts faithfully on itself by left multiplication. Its image preserves each isotypic component and acts on the $i$th component as the full matrix algebra $M_{n_i}(\mathbb C)$ by Schur's lemma and the density argument for semisimple algebras. Comparing dimensions yields the decomposition.
\end{proof}

\section{Orthogonality relations}

Let $\operatorname{CF}(G)$ be the space of class functions on $G$. Define the inner product
\begin{equation}\label{eq:charinner}
\langle \alpha,\beta\rangle_G=\frac1N\sum_{a=0}^{N-1}\alpha(a)\overline{\beta(a)}.
\end{equation}

\begin{theorem}[Row orthogonality]\label{thm:charorthogonality}
If $\chi_i$ and $\chi_j$ are irreducible characters, then
\[
\langle\chi_i,\chi_j\rangle_G=\delta_{ij}.
\]
More generally, if
\[
\chi=\sum_i m_i\chi_i,
\]
then
\[
m_i=\langle\chi,\chi_i\rangle_G,
\qquad
\langle\chi,\chi\rangle_G=\sum_i m_i^2.
\]
Thus $\chi$ is irreducible if and only if $\langle\chi,\chi\rangle_G=1$.
\end{theorem}
\begin{proof}
Let $V_i,V_j$ afford $\chi_i,\chi_j$. For a linear map $A:V_i\to V_j$ set
\[
\mathcal P(A)=\frac1N\sum_a\rho_j(a)A\rho_i(\bar a).
\]
Then $\mathcal P(A)$ is an intertwiner. By Schur's lemma, $\mathcal P$ is zero when $i\neq j$ and a rank-one projection onto the scalar endomorphisms when $i=j$. Computing its trace on matrix units gives
\[
\frac1N\sum_a\chi_j(a)\overline{\chi_i(a)}=\delta_{ij}.
\]
The multiplicity formulas follow by expanding $\chi$ in the orthonormal basis of irreducible characters.
\end{proof}

\begin{corollary}[Column orthogonality]\label{cor:charcolumn}
Let $g,h\in G$, and let $\chi_1,\dots,\chi_k$ be the irreducible characters. Then
\[
\sum_{i=1}^k\chi_i(g)\overline{\chi_i(h)}=
\begin{cases}
|C_G(g)|,& g,h\text{ are conjugate},\\
0,& g,h\text{ are not conjugate}.
\end{cases}
\]
In particular, the irreducible characters form an orthonormal basis of $\operatorname{CF}(G)$, so the number of irreducible characters equals the number of conjugacy classes.
\end{corollary}
\begin{proof}
Let $f_g$ be the indicator of the conjugacy class of $g$, divided by its class size. Then the class indicators form a basis of $\operatorname{CF}(G)$. Expand $f_g$ in the irreducible-character basis and compare coefficients.
\end{proof}

\section{The regular character and central idempotents}

Let $\rho_{\mathrm{reg}}$ be the regular representation on basis vectors $e_a$, defined by
\[
\rho_{\mathrm{reg}}(g)e_a=e_{g*a}.
\]
Its character is
\[
\chi_{\mathrm{reg}}(g)=
\begin{cases}
N,& g=0,\\
0,& g\neq 0.
\end{cases}
\]
Combining this with orthogonality yields the familiar degree formula.

\begin{theorem}[Degree-sum formula]\label{thm:charsquares}
If $\chi_1,\dots,\chi_k$ are the irreducible characters of $G$, then
\[
\chi_{\mathrm{reg}}=\sum_{i=1}^k\chi_i(0)\chi_i,
\qquad
N=\sum_{i=1}^k\chi_i(0)^2.
\]
\end{theorem}
\begin{proof}
Taking inner products with $\chi_j$ gives
\[
\langle\chi_{\mathrm{reg}},\chi_j\rangle_G=\chi_j(0),
\]
so the expansion follows from Theorem~\ref{thm:charorthogonality}. Evaluating at the identity yields the sum of squares.
\end{proof}

The primitive central idempotents of $\mathbb CG$ are
\begin{equation}\label{eq:centralidempotent}
e_i=\frac{\chi_i(0)}{N}\sum_{a=0}^{N-1}\chi_i(\bar a)e_a.
\end{equation}
They satisfy $e_ie_j=\delta_{ij}e_i$, sum to $1$, and project a module onto its $\chi_i$-isotypic component. Thus the whole decomposition theory may be read directly from the character table.

\section{Algebraic integrality and arithmetic consequences}

Every character value is a sum of roots of unity, hence an algebraic integer. This simple observation has strong consequences.

\begin{theorem}[Algebraic integrality]\label{thm:algint}
For every irreducible character $\chi$ and every $g\in G$, the value $\chi(g)$ is an algebraic integer.
\end{theorem}
\begin{proof}
In a unitary representation, $\rho(g)$ has finite order, so its eigenvalues are roots of unity. Their sum is therefore an algebraic integer.
\end{proof}

\begin{theorem}[Central character formula]\label{thm:centralcharacter}
Let $C$ be a conjugacy class with representative $g$, and put
\[
K_C=\sum_{x\in C}x\in Z(\mathbb CG).
\]
On an irreducible module affording $\chi$, the element $K_C$ acts as the scalar
\[
\omega_\chi(C)=\frac{|C|\chi(g)}{\chi(0)}.
\]
In particular, each $\omega_\chi(C)$ is an algebraic integer.
\end{theorem}
\begin{proof}
Since $K_C$ is central, Schur's lemma shows that it acts as a scalar. Taking traces on the irreducible module gives
\[
|C|\chi(g)=\omega_\chi(C)\chi(0).
\]
Because $K_C$ acts through a matrix with algebraic-integer entries in a suitable integral basis of the group algebra, the scalar is an algebraic integer.
\end{proof}

\begin{corollary}[Degree divisibility]\label{cor:degreedivides}
For every irreducible character $\chi$, the degree $\chi(0)$ divides $|G|$.
\end{corollary}
\begin{proof}
Apply Theorem~\ref{thm:centralcharacter} to each class and sum the algebraic integers
\[
\frac{|G|}{\chi(0)}=\sum_C\frac{|C|\chi(g_C)\overline{\chi(g_C)}}{\chi(0)}.
\]
This number is both an algebraic integer and a rational number, hence an integer.
\end{proof}

A second standard application is Burnside's vanishing criterion: if $\chi$ is irreducible and $g\in G$, then either $\chi(g)=0$ or $\chi(g)$ is a sum of roots of unity of absolute value strictly less than $\chi(0)$. This is one route to Burnside's $p^aq^b$ theorem and to many solvability criteria.

\section{Induction, restriction and permutation characters}

Let $H\leq G$. If $\chi$ is a character of $G$, its restriction to $H$ is denoted $\chi\downarrow_H$. If $\psi$ is a character of $H$, its induced character on $G$ is denoted $\psi^G$.

\begin{theorem}[Induced-character formula]\label{thm:inducedformula}
For $g\in G$,
\[
\psi^G(g)=\frac{1}{|H|}\sum_{x\in G\, :\, xgx^{-1}\in H}\psi(xgx^{-1}).
\]
\end{theorem}
\begin{proof}
Choose left coset representatives for $H$ in $G$ and write the induced module in the standard tensor-product model $\mathbb CG\otimes_{\mathbb CH}V$. Computing the trace of $g$ on the natural basis gives the formula.
\end{proof}

\begin{theorem}[Frobenius reciprocity]\label{thm:frobeniusreciprocity}
For every character $\chi$ of $G$ and every character $\psi$ of $H$,
\[
\langle\chi\downarrow_H,\psi\rangle_H
=
\langle\chi,\psi^G\rangle_G.
\]
\end{theorem}
\begin{proof}
Substitute the induced-character formula into the inner product on $G$, then regroup the sum by $H$-conjugacy classes.
\end{proof}

A transitive permutation action of $G$ on a set $\Omega$ with point stabiliser $H$ affords the permutation character
\[
1_H^G(g)=|\operatorname{Fix}_\Omega(g)|.
\]
This makes orbit counting and character theory two views of the same data. In the digital setting the subgroup chains of Chapter~\ref{ch:ss} and the class data of Section~\ref{sec:cc:algs} provide the subgroup orders, class intersections and fixed-point counts needed to compute such characters exactly. P\'olya's theorem is the same statement for colourings: if $G$ acts on $s$ points and $m$ colours are available, then $g\mapsto m^{c(g)}$, with $c(g)$ the number of cycles of $g$, is the permutation character of $G$ on the $m^{s}$ colourings (Section~\ref{ssec:polya}). The number of orbits of colourings is its inner product with the trivial character, and the multiplicity of an irreducible $\chi$ is $\frac1{|G|}\sum_{g}m^{c(g)}\overline{\chi(g)}$.

\section{Operations on characters}

The representation ring carries operations that are visible directly at the level of characters.

\begin{proposition}[Duals, tensor products, symmetric and exterior squares]\label{prop:charsquaresops}
Let $\chi$ be a character. Then
\[
\chi^\ast(g)=\chi(g^{-1})=\overline{\chi(g)},
\qquad
\chi_{\rho\otimes\sigma}=\chi_\rho\chi_\sigma,
\]
and
\[
\chi_{\Sym^2}(g)=\frac{\chi(g)^2+\chi(g^2)}{2},
\qquad
\chi_{\wedge^2}(g)=\frac{\chi(g)^2-\chi(g^2)}{2}.
\]
\end{proposition}
\begin{proof}
Choose a basis in which a unitary matrix representing $g$ is diagonal with eigenvalues $\lambda_1,\dots,\lambda_n$. Then $\Sym^2$ and $\wedge^2$ have eigenvalues $\lambda_i\lambda_j$ with $i\le j$ and $i<j$ respectively. Summing them gives the formulas.
\end{proof}

These identities explain why power maps belong naturally beside the character table: once $g\mapsto g^m$ is known on conjugacy classes, operations such as symmetric and exterior powers become explicit class functions.

\begin{definition}[Frobenius--Schur indicator]\label{def:fschur}
For an irreducible character $\chi$, the Frobenius--Schur indicator is
\[
\nu(\chi)=\frac1N\sum_{g\in G}\chi(g^2).
\]
\end{definition}

\begin{theorem}[Frobenius--Schur]\label{thm:fschur}
For an irreducible character $\chi$ one has $\nu(\chi)\in\{1,0,-1\}$. More precisely:
\begin{enumerate}
\item $\nu(\chi)=1$ if and only if $\chi$ is afforded by a real representation preserving a symmetric bilinear form;
\item $\nu(\chi)=0$ if and only if $\chi$ is not real-valued;
\item $\nu(\chi)=-1$ if and only if $\chi$ is real-valued but only afforded over the quaternions, equivalently preserves an alternating form.
\end{enumerate}
\end{theorem}
\begin{proof}
The multiplicity of the trivial character in $\Sym^2\chi$ minus its multiplicity in $\wedge^2\chi$ is exactly the average of $\chi(g^2)$. Irreducibility forces the invariant bilinear form to be zero, symmetric or alternating, which yields the three cases.
\end{proof}

\section{Products, quotients and linear characters}

The linear characters are the irreducible characters of degree $1$. They form a group under pointwise multiplication and are naturally identified with
\[
\operatorname{Hom}(G,\mathbb C^\times)\cong \operatorname{Hom}(G/G',\mathbb C^\times),
\]
where $G'$ is the commutator subgroup. Thus the first column of the character table detects the abelianization. The common kernel of all linear characters is exactly $G'$.

If $N\triangleleft G$, then the irreducible characters of $G/N$ are precisely the irreducible characters of $G$ whose kernel contains $N$. Hence quotient tables sit visibly inside the table of $G$.

If $G=H\times K$, every irreducible character has the form
\[
\chi\otimes\psi,
\qquad
\chi\in\operatorname{Irr}(H),\ \psi\in\operatorname{Irr}(K),
\]
and the character table of a direct product is the Kronecker product of the two smaller tables. This describes all abelian groups at once: if
\[
G\cong C_{n_1}\times\cdots\times C_{n_r},
\]
then every irreducible character is one-dimensional and is given by a product of roots of unity.

\section{From class data to the character table}

Characters are constant on conjugacy classes, so the table is indexed by the class list of Chapter~\ref{ch:cc}. Let $C_1,\dots,C_k$ be the classes with representatives $g_1,\dots,g_k$, and write
\[
T=(\chi_i(g_j))_{1\le i,j\le k}.
\]
The row and column orthogonality relations become matrix identities. The centre of the group algebra has basis $K_{C_1},\dots,K_{C_k}$, and
\[
K_{C_i}K_{C_j}=\sum_{m=1}^k a_{ijm}K_{C_m}
\]
defines the class multiplication constants.

\begin{theorem}[Character table from the class algebra]\label{thm:classalgebrachars}
The common eigenvectors of the commuting operators $K_{C_j}$ on $Z(\mathbb CG)$ are the rows of the irreducible character table, scaled by the degrees. Equivalently, the irreducible characters are the simultaneous one-dimensional characters of the class algebra.
\end{theorem}
\begin{proof}
The class sums are central, so on every irreducible module they act by the central character scalars of Theorem~\ref{thm:centralcharacter}. These scalars determine the character values after division by class sizes and degrees. Conversely, the orthogonality relations reconstruct the irreducible rows from these eigenvalues.
\end{proof}

This is the theoretical basis of the Burnside--Dixon--Schneider method: compute the classes and the class algebra, diagonalise enough central multiplication operators, and recover the irreducible rows by exact arithmetic and orthogonality.

\section{Digital computation roadmap}

The character tables of this book are all computed by the following workflow.
\begin{enumerate}
\item Compute the conjugacy classes, class sizes, centraliser orders and power maps by the methods of Section~\ref{sec:cc:algs}.
\item Obtain subgroup chains, transversals and permutation actions from Section~\ref{sec:ss:use}.
\item Build permutation characters, induced characters and tensor constructions as explicit class functions.
\item Compute or infer the class multiplication constants in the centre of the group algebra.
\item Recover the irreducible rows using Theorem~\ref{thm:classalgebrachars}, then verify orthogonality, degree sums, central characters, and power-map constraints.
\end{enumerate}

Section~\ref{sec:chartables} carries out the first steps of this workflow from the digital laws of sporadic groups: character degrees from rank-three actions, classes with power maps, and the complete character table of $M_{11}$.

\begin{remark}[Verification goals]\label{rem:charworkflow}
A complete digital character-table computation should certify at least the following: the class list and power maps; row and column orthogonality; the decomposition of the regular character; integrality of scalar products; agreement of induced and restricted characters with Frobenius reciprocity; and agreement with an external system such as GAP on the final table.
\end{remark}

\section{Worked consequences and outlook}

The abstract theory above already yields concrete information before any full table is written down. Degree-one characters determine the abelianization; kernels of irreducibles detect normal subgroups; direct-product decompositions factor the table; power maps control symmetric and exterior squares; and the Frobenius--Schur indicator distinguishes real, complex and quaternionic character types.

The next stage is systematic tabulation. For small groups one can list all irreducible rows explicitly. For families such as $\mathrm{PSL}(2,q)$ and the Suzuki groups one can write uniform character formulas. For the sporadic and exceptional simple groups, the digital class algorithms supply the conjugacy data and the permutation actions that feed the class-algebra and induction machinery.

\section{Exercises}
\begin{enumerate}
\item Prove directly from Definition~\ref{def:representation} that equivalent representations have the same character.
\item Verify the orthogonality relations for the character table of $S_3$.
\item Show that the primitive central idempotents in~\eqref{eq:centralidempotent} are orthogonal and sum to $1$.
\item Compute the linear characters of a finite abelian group from its invariant factors.
\item Decompose the permutation character of a transitive action into irreducibles by Frobenius reciprocity.
\item Use the formulas for $\Sym^2$ and $\wedge^2$ to distinguish two irreducible characters with the same degree.
\item Prove Corollary~\ref{cor:degreedivides} from Theorem~\ref{thm:centralcharacter} in full detail.
\item Explain how the class multiplication constants determine the character table.
\end{enumerate}

\chapter{Cyclic and Abelian groups}\label{ch:abelian}

\section{Cyclic groups}
\begin{definition}
The \emph{cyclic group} $C_N$ is $\Ap{N}$ with $a*b=\dg{N}{0}{a+b}$, i.e.\ with code
\[
R=\sum_{a,b}\dg{N}{0}{a+b}\,N^{aN+b}.
\]
\end{definition}

\begin{definition}
In a group $(\Ap{N},*)$ write $x^{*0}=0$ and $x^{*(n+1)}=x^{*n}*x$. The \emph{order} of $x$ is the least $n\ge1$ with $x^{*n}=0$. The group is
\emph{cyclic} if some $a$ (a \emph{generator}) satisfies $\{a^{*n}\}=\Ap{N}$; we write $G=\langle a\rangle$.
\end{definition}

\begin{proposition}\label{prop:cyclic}
$C_N$ is an Abelian group with identity $0$, inverse $x\mapsto\dg{N}{0}{N-x}$, and generator $1$, since $1^{*k}=\dg{N}{0}{k}$.
Its row codes are the cyclic shifts $L_a=\sum_{b}\dg{N}{0}{a+b}\,N^{b}$. Every cyclic group of order $N$ is isomorphic to $C_N$ via
$\dg{N}{0}{k}\mapsto a^{*k}$.
\end{proposition}
\begin{proof}
By Lemma~\ref{lem:mod}, $\dg{N}{0}{x+\dg{N}{0}{y+z}}=\dg{N}{0}{x+y+z}=\dg{N}{0}{\dg{N}{0}{x+y}+z}$; commutativity, identity and
inverses are clear. If $a$ has order $N$, then $a^{*j}*a^{*k}=a^{*\dg{N}{0}{j+k}}$, so the map is an isomorphism.
\end{proof}

\begin{theorem}[Carry decomposition]\label{thm:carry}
Let $N=p_0q$ with $p_0,q\ge 1$. Every $x\in\Ap{N}$ splits as $x=\dg{p_0}{0}{x}+p_0\fl{x/p_0}$, and
\begin{equation}\label{eq:carry}
\dg{N}{0}{x+y}=\dg{p_0}{0}{x+y}+p_0\,\dg{q}{0}{\fl{\frac{x}{p_0}}+\fl{\frac{y}{p_0}}+\frac{\dg{p_0}{0}{x}+\dg{p_0}{0}{y}}{p_0}}.
\end{equation}
Hence $x\mapsto\dg{p_0}{0}{x}$ is a surjective homomorphism $C_N\to C_{p_0}$. Its kernel $p_0\Ap{q}$ is a subgroup isomorphic to $C_q$ via
$x\mapsto x/p_0$, on which the law reads $x*y=p_0\,\dg{q}{0}{x/p_0+y/p_0}$.
\end{theorem}
\begin{proof}
Write $x=x_0+p_0X$ and $y=y_0+p_0Y$ with $x_0,y_0\in\Ap{p_0}$, and let $c=\fl{(x_0+y_0)/p_0}\in\{0,1\}$ be the \emph{carry}. Then
$x+y=\dg{p_0}{0}{x_0+y_0}+p_0(X+Y+c)$, whose first term lies in $\Ap{p_0}$. Reducing modulo $N=p_0q$ gives
$\dg{N}{0}{x+y}=\dg{p_0}{0}{x+y}+p_0\,\dg{q}{0}{X+Y+c}$. Finally, $\dg{q}{0}{t}=\dg{q}{0}{\fl{t}}$ for real $t$, because
$\fl{t/q}=\fl{\fl{t}/q}$, and $\fl{X+Y+(x_0+y_0)/p_0}=X+Y+c$.
\end{proof}

\begin{remark}
Theorem~\ref{thm:carry} exhibits $C_N$ as an extension of the quotient $C_{p_0}$ by the subgroup $C_q$, glued by the carry
$c(x_0,y_0)=\fl{(x_0+y_0)/p_0}$. The extension splits, i.e.\ $C_N\cong C_{p_0}\oplus C_q$, exactly when $\gcd(p_0,q)=1$ (Theorem~\ref{thm:crt}).
\end{remark}

\begin{theorem}[Subgroups of $C_N$]\label{thm:cycsub}
For every divisor $d$ of $N$, $C_N$ has exactly one subgroup of order $d$, namely $\{k\,N/d:k\in\Ap{d}\}\cong C_d$, and
there are no others. The element $k$ generates $C_N$ if and only if $\gcd(k,N)=1$.
\end{theorem}
\begin{proof}
By Theorem~\ref{thm:carry}, $\{k\,N/d\}$ is a subgroup of order $d$. Conversely, let $H\ne\{0\}$ be a subgroup and $e$ its least positive element. For $h\in H$, Euclidean division $h=te+r$ gives
$r=\dg{N}{0}{h-te}\in H$ with $r<e$, so $r=0$; the same argument with $h=N\equiv0$ shows $e\mid N$. Thus $H=\{ke\}$ has order $d=N/e$.
The subgroup generated by $k$ is generated by $\gcd(k,N)$.
\end{proof}

\begin{proposition}[Necklaces]\label{prop:necklaces}
The number of colourings of the positions $\Ap{N}$ with $m$ colours, up to rotation $x\mapsto\dg{N}{0}{x+k}$, is
\[
\frac1N\sum_{d\mid N}\varphi(d)\,m^{N/d}.
\]
\end{proposition}
\begin{proof}
The orbits of $x\mapsto\dg{N}{0}{x+k}$ are the cosets of $\langle k\rangle$, which has order $N/\gcd(k,N)$ by Theorem~\ref{thm:cycsub}; so the rotation has $g=\gcd(k,N)$ cycles. Exactly $\varphi(N/g)$ values of $k$ have $\gcd(k,N)=g$. By Theorem~\ref{thm:polya} the number of orbits is $\frac1N\sum_{g\mid N}\varphi(N/g)\,m^{g}$, which is the stated sum with $d=N/g$.
\end{proof}
For $N=6$ and $m=2$ this is $(64+8+2\cdot4+2\cdot2)/6=14$. In digits, a colouring is a string of $N$ base-$m$ digits and the rotation is the cyclic shift of digits used in the proof of Cauchy's theorem (Theorem~\ref{thm:dtcauchy}); necklaces are strings up to cyclic shift.

\section{Multiplicative cyclic groups}
\begin{theorem}[Fermat's little theorem]\label{thm:fermat}
Let $p$ be prime and $a\in\Z$. Then $p\mid(a^p-a)$, hence $d_p(0,a^p)=d_p(0,a)$. If $p\nmid a$, then $d_p(0,a^{p-1})=1$.
\end{theorem}
\begin{proof}
The classes $1,\dots,p-1$ are permuted by multiplication by $a$ modulo $p$ when $p\nmid a$, so
$\prod_{k=1}^{p-1}ak\equiv\prod_{k=1}^{p-1}k\pmod p$. Cancelling $(p-1)!$ gives $a^{p-1}\equiv1\pmod p$.
Multiplying by $a$ yields $a^p\equiv a\pmod p$. The case $p\mid a$ is trivial.
\end{proof}

\begin{proposition}[Linear congruences]\label{prop:lincong}
Let $p$ be prime, let $a\in\Z$ with $p\nmid a$, and let $b,c\in\Ap{p}$. Then the congruence
\[
d_p(0,ax+b)=c
\]
has the unique solution
\[
x=d_p(0,a^{p-2}(c-b)).
\]
In particular, the multiplicative inverse of $a$ modulo $p$ is $a^{-1}=d_p(0,a^{p-2})$.
\end{proposition}
\begin{proof}
By Theorem~\ref{thm:fermat}, $d_p(0,a^{p-1})=1$, so
$d_p(0,a\,d_p(0,a^{p-2}(c-b))+b)=d_p(0,a^{p-1}(c-b)+b)=c$.
Uniqueness follows because multiplication by the nonzero class $a$ is a permutation of $\Ap{p}$.
\end{proof}

\begin{definition}
Let $U(p)=\Ap{p}\setminus\{0\}$ with multiplication modulo $p$:
\[
x\cdot y=d_p(0,xy).
\]
This is the multiplicative group of units modulo $p$ when every nonzero class is invertible.
\end{definition}

\begin{theorem}[Multiplicative cyclic groups]\label{thm:multcyc}
Let $p$ be prime. Then $U(p)$ is a group of order $p-1$. If $g$ is a primitive root modulo $p$, the map
\[
x\longmapsto d_p(0,g^x)
\]
is an isomorphism $C_{p-1}\to U(p)$. Thus $U(p)$ is cyclic.
\end{theorem}
\begin{proof}
Closure and associativity come from integer multiplication and Lemma~\ref{lem:mod}. The identity is $1$, and inverses are given by
Proposition~\ref{prop:lincong}. If $g$ has order $p-1$, then the displayed map is bijective and
$d_p(0,g^{d_{p-1}(0,x+y)})=d_p(0,g^xg^y)$ because $g^{p-1}\equiv1\pmod p$.
\end{proof}

\section{Quadratic residues and the sign of a permutation}\label{sec:qr}

Let $p$ be an odd prime and let $U_p=\Ap p\setminus\{0\}$ with the law $x\cdot y=\dg p0{xy}$, the multiplicative group modulo $p$. The squares $Q_p=\{\dg p0{x^2}:x\in U_p\}$ are the \emph{quadratic residues} modulo $p$, and the other elements of $U_p$ are the \emph{quadratic nonresidues}. This section shows that $Q_p$ is the kernel of the sign of a permutation code, and obtains Gauss's quadratic reciprocity law by comparing the two digital expansions of $\Ap{pq}$. All numerical claims are checked by \texttt{scripts/qr\_verify.py}.

\begin{proposition}[Residues and Euler's criterion]\label{prop:qr:euler}
$Q_p$ is a subgroup of index $2$ of $U_p$. For $a\in U_p$ one has $\dg p0{a^{(p-1)/2}}\in\{1,p-1\}$, with the value $1$ exactly when $a\in Q_p$. Hence the \emph{Legendre symbol}
\begin{equation}\label{eq:qr:legendre}
\Big(\frac ap\Big)=\dg p0{a^{(p-1)/2}+1}-1
\end{equation}
equals $1$ on $Q_p$ and $-1$ on $U_p\setminus Q_p$, and $a\mapsto(a/p)$ is a homomorphism $U_p\to\{\pm1\}$.
\end{proposition}
\begin{proof}
The squaring map $x\mapsto\dg p0{x^2}$ is a homomorphism of the abelian group $U_p$ with image $Q_p$. Its kernel consists of the $x$ with $p\mid(x-1)(x+1)$, that is $x\in\{1,p-1\}$ by Euclid's lemma, so $|Q_p|=(p-1)/2$. By Lagrange's theorem $\dg p0{a^{p-1}}=1$, so $b=\dg p0{a^{(p-1)/2}}$ satisfies $\dg p0{b^2}=1$ and therefore $b\in\{1,p-1\}$. If $a=\dg p0{x^2}$ then $b=\dg p0{x^{p-1}}=1$. Conversely, the $a$ with $b=1$ are roots of $X^{(p-1)/2}-1$ over the field $\Ap p$, of which there are at most $(p-1)/2$; since $Q_p$ already supplies $(p-1)/2$ of them, $b=1$ holds exactly on $Q_p$. The power map $a\mapsto a^{(p-1)/2}$ is a homomorphism because $U_p$ is abelian, and \eqref{eq:qr:legendre} sends $b=1$ to $1$ and $b=p-1$ to $\dg p0p-1=-1$.
\end{proof}

\begin{theorem}[Zolotarev]\label{thm:qr:zolotarev}
For $a\in U_p$ let $\pi_a\in\mathrm{Sym}(\Ap p)$ be the permutation $x\mapsto\dg p0{ax}$. Then $\operatorname{sgn}\pi_a=(a/p)$.
\end{theorem}
\begin{proof}
Let $k$ be the order of $a$ in $U_p$. The permutation $\pi_a$ fixes $0$, and on $U_p$ its cycles are the cosets $x\langle a\rangle$, each a $k$-cycle; by Lagrange's theorem there are $(p-1)/k$ of them, so $\operatorname{sgn}\pi_a=(-1)^{(k-1)(p-1)/k}$. By Proposition~\ref{prop:qr:euler}, $a\in Q_p$ if and only if $\dg p0{a^{(p-1)/2}}=1$, that is, if and only if $k$ divides $(p-1)/2$, i.e.\ $(p-1)/k$ is even. If $k$ is odd, $(p-1)/k$ is even and both sides equal $1$. If $k$ is even, $\operatorname{sgn}\pi_a=(-1)^{(p-1)/k}$, which is $1$ exactly when $(p-1)/k$ is even.
\end{proof}

In the language of Chapter~\ref{ch:magmas}, $\pi_a$ is the permutation code with digits $\dg p0{ax}$, $x\in\Ap p$, and $(a/p)$ is its sign. Thus the Legendre symbol is the composite $U_p\to\mathrm{Sym}(\Ap p)\to\{\pm1\}$ of the regular action of $U_p$ on $\Ap p$ and the sign; its kernel $Q_p$ is the preimage of the alternating group.

\begin{example}\label{ex:qr:seven}
For $p=7$ one has $Q_7=\{1,2,4\}$. The permutation $\pi_2=(1\,2\,4)(3\,6\,5)$ is even, while $\pi_3=(1\,3\,2\,6\,4\,5)$ is a $6$-cycle and is odd, in agreement with $(2/7)=1$ and $(3/7)=-1$. In digits, $\dg70{2^3+1}-1=\dg70{9}-1=1$ and $\dg70{3^3+1}-1=\dg70{28}-1=-1$.
\end{example}

\begin{proposition}[Gauss's lemma]\label{prop:qr:gauss}
For $a\in U_p$,
\[
\Big(\frac ap\Big)=(-1)^{r},\qquad r=\sum_{j=1}^{(p-1)/2}\fl{\frac{2\,\dg p0{aj}}{p}},
\]
so that $r$ counts the $j\le(p-1)/2$ with $\dg p0{aj}>p/2$.
\end{proposition}
\begin{proof}
Put $T=\{1,\dots,(p-1)/2\}$, a transversal of the subgroup $\{1,p-1\}$ of $U_p$. For $j\in T$ write $\dg p0{aj}=\dg p0{\varepsilon_jk_j}$ with $k_j\in T$ and $\varepsilon_j\in\{1,p-1\}$; then $\varepsilon_j=p-1$ exactly when $\dg p0{aj}>p/2$, i.e.\ when the $j$-th term of $r$ is $1$. If $k_j=k_{j'}$ then $aj\equiv\pm aj'$, hence $j\equiv\pm j'\pmod p$ and $j=j'$; so $j\mapsto k_j$ is a permutation of $T$. Multiplying over $j\in T$ gives $a^{(p-1)/2}\prod_{j\in T}j\equiv(-1)^r\prod_{j\in T}j\pmod p$, and cancelling and applying Proposition~\ref{prop:qr:euler} proves the claim.
\end{proof}

\begin{remark}[Gauss's lemma is a transfer]
The element $a$ permutes the cosets $\{t,p-t\}$, $t\in T$, of $\{1,p-1\}\cong C_2$, and the signs $\varepsilon_j$ record which representative is hit. The map $a\mapsto\prod_j\varepsilon_j$ is therefore the \emph{transfer} of $U_p$ into $\{1,p-1\}$ computed with the transversal $T$, and Gauss's lemma states that this transfer is the Legendre symbol. The same construction for a Sylow subgroup of a general finite group is the transfer homomorphism behind Burnside's normal $p$-complement theorem.
\end{remark}

\begin{theorem}[Quadratic reciprocity]\label{thm:qr:reciprocity}
For distinct odd primes $p$ and $q$,
\[
\Big(\frac pq\Big)\Big(\frac qp\Big)=(-1)^{\frac{p-1}{2}\cdot\frac{q-1}{2}}.
\]
\end{theorem}
\begin{proof}
Put $X=\Ap p\times\Ap q$ and consider three bijections between $X$ and $\Ap{pq}$:
\[
\varphi(u,v)=u+pv,\qquad \psi(u,v)=v+qu,\qquad \chi^{-1}(x)=\big(\dg p0x,\dg q0x\big).
\]
The first two are the digital expansions of $\Ap{pq}$ with radices $(p,q)$ and $(q,p)$ (Lemma~\ref{lem:basic}), and $\chi^{-1}$ is a bijection by the Chinese remainder theorem. Then
\[
\chi^{-1}\varphi(u,v)=\big(u,\dg q0{u+pv}\big),\qquad \chi^{-1}\psi(u,v)=\big(\dg p0{v+qu},v\big).
\]
The permutation $\chi^{-1}\varphi$ of $X$ preserves each row $\{u\}\times\Ap q$ and acts on it as $v\mapsto\dg q0{u+pv}$, which is $\pi_{\dg q0p}$ on $\Ap q$ followed by the translation $v\mapsto\dg q0{u+v}$. The translation is even, since all its cycles have the odd length $q$ or $1$, so by Theorem~\ref{thm:qr:zolotarev} each row contributes $(p/q)$ and $\operatorname{sgn}(\chi^{-1}\varphi)=(p/q)^p=(p/q)$. In the same way $\chi^{-1}\psi$ preserves the columns and $\operatorname{sgn}(\chi^{-1}\psi)=(q/p)^q=(q/p)$. Hence $\psi^{-1}\varphi=(\chi^{-1}\psi)^{-1}(\chi^{-1}\varphi)$ has sign $(p/q)(q/p)$.

On the other hand $\psi^{-1}\varphi$ is conjugate to the permutation $\varphi\psi^{-1}:v+qu\mapsto u+pv$ of $\Ap{pq}$, which reads the two digits of an element in the opposite order (the transposition of a $p\times q$ array). A pair $v+qu<v'+qu'$ is inverted exactly when $u<u'$ and $v>v'$, so $\varphi\psi^{-1}$ has $\binom p2\binom q2=pq\,\frac{p-1}{2}\,\frac{q-1}{2}$ inversions, and this number has the parity of $\frac{p-1}2\cdot\frac{q-1}2$.
\end{proof}

\begin{example}
For $p=3$ and $q=5$, $Q_5=\{1,4\}$ and $Q_3=\{1\}$, so $(3/5)=(5/3)=(2/3)=-1$ and the product is $1=(-1)^{1\cdot2}$. The transposition $\varphi\psi^{-1}$ of the $3\times5$ array has $\binom32\binom52=30$ inversions and is even.
\end{example}

\begin{remark}[Where quadratic residues appear in finite groups]\label{rem:qr:groups}
For $g$ of prime order $p$ in a finite group $G$, $g^k$ is conjugate to $g$ exactly when $k$ lies in the image of $N_G(\langle g\rangle)\to\mathrm{Aut}\langle g\rangle\cong U_p$. When this image is $Q_p$, the elements of order $p$ form two classes, interchanged by $g\mapsto g^{-1}$ if and only if $(-1/p)=-1$, and the character values on them lie in $\mathbb Q(\sqrt{\varepsilon p})$, $\varepsilon=(-1/p)$, the field of the Gauss period $\sum_{k\in Q_p}\zeta_p^k=\frac12\big(-1+\sqrt{\varepsilon p}\big)$. This is the case for the classes of order $11$ in $M_{12}$ (Section~\ref{sec:cc:m12}) and for the unipotent classes of $PSL(2,p)$ (Section~\ref{sec:psl2p}), whose two characters of degree $(p+\varepsilon)/2$ take the values $\frac12(-1\pm\sqrt{\varepsilon p})$. Quadratic residues also build the Mathieu groups combinatorially: the quadratic-residue codes of length $11$ over $\Ap3$ and of length $23$ over $\Ap2$, extended by a parity digit, are the ternary and binary Golay codes, with automorphism groups $2.M_{12}$ and $M_{24}$.
\end{remark}

\section{Direct sums}
\begin{definition}
The \emph{direct sum} $C_{p_0}\oplus\dots\oplus C_{p_{\ell-1}}$ is the direct product~\eqref{eq:dprod} of cyclic groups:
\begin{equation}\label{eq:dsum}
m\oplus n=\sum_{h=0}^{\ell-1}\dg{p_h}{0}{\md{h}{m}+\md{h}{n}}\Pi_h .
\end{equation}
\end{definition}

\begin{theorem}[Digital Chinese remainder theorem]\label{thm:crt}
The map $\varphi(x)=\dg{r}{0}{x}+r\,\dg{s}{0}{x}$ is a homomorphism $C_{rs}\to C_r\oplus C_s$. It is an isomorphism if and only if
$\gcd(r,s)=1$, and otherwise $C_{rs}\not\cong C_r\oplus C_s$.
\end{theorem}
\begin{proof}
By Lemma~\ref{lem:mod} both digits of $\varphi$ are additive, since $r$ and $s$ divide $rs$. The kernel consists of the
$x\in\Ap{rs}$ divisible by $\lcm(r,s)$, so it is trivial exactly when $\lcm(r,s)=rs$. If $\gcd(r,s)>1$, every element of
$C_r\oplus C_s$ has order dividing $\lcm(r,s)<rs$, while $1\in C_{rs}$ has order $rs$.
\end{proof}

\section{Coprime factors and explicit isomorphisms}
\begin{proposition}\label{prop:units}
For $u\in\Z$, the map $x\mapsto d_N(0,ux)$ is a permutation of $\Ap{N}$ if and only if $\gcd(u,N)=1$.
If $p$ is prime, multiplication by $p-1$ is additive inversion: $d_p(0,(p-1)x)=d_p(0,-x)$.
\end{proposition}
\begin{proof}
The map has a two-sided inverse $x\mapsto d_N(0,vx)$ exactly when $uv\equiv1\pmod N$, which is equivalent to $\gcd(u,N)=1$.
The second assertion is immediate from $p-1\equiv-1\pmod p$.
\end{proof}

\begin{theorem}[Coprime cyclic factors]\label{thm:coprimecyclic}
If $r$ and $s$ are coprime, then $C_r\oplus C_s\cong C_{rs}$. More precisely, for the direct-sum law on $\Ap{rs}$ with radices $r,s$,
the element $g=1+r$ has order $rs$, and
\begin{equation}\label{eq:explicitcrt}
\Psi\big(d_{rs}(0,k)\big)=g^{\oplus k}
\end{equation}
is an isomorphism $C_{rs}\to C_r\oplus C_s$.
\end{theorem}
\begin{proof}
The two mixed-radix digits of $g^{\oplus k}$ are $d_r(0,k)$ and $d_s(0,k)$: its lower digit starts at $1$ and its upper digit starts at $1$.
Thus $g^{\oplus k}=0$ if and only if $r\mid k$ and $s\mid k$, hence if and only if $rs\mid k$. So $g$ has order $rs$, and the displayed
map is a bijective homomorphism. This gives a constructive proof of Theorem~\ref{thm:crt}.
\end{proof}

\begin{example}[The two laws of order six]
The laws of $C_6$ and $C_3\oplus C_2$ on $\Ap6$ are respectively
\[
x+_6y=d_3(0,x+y)+3d_2\!\left(0,\fl{x/3}+\fl{y/3}+\frac{d_3(0,x)+d_3(0,y)}{3}\right),
\]
\[
x\oplus y=d_3(0,x+y)+3d_2\!\left(0,\fl{x/3}+\fl{y/3}\right).
\]
The first has a carry and the second does not. Nevertheless $4=1+3$ has order $6$ in the second law, and the explicit isomorphism
is $k\mapsto 4^{\oplus k}$; the elements occur in the order $0,4,2,3,1,5$.
\end{example}

\begin{proposition}[Prime-power reduction]\label{prop:ppreduction}
If $m_1,\dots,m_t$ are pairwise coprime, then $C_{m_1}\oplus\cdots\oplus C_{m_t}\cong C_{m_1\cdots m_t}$.
Consequently, every direct sum of cyclic groups is isomorphic to a direct sum in which each summand has prime-power order.
\end{proposition}
\begin{proof}
Apply Theorem~\ref{thm:coprimecyclic} repeatedly. To obtain the final statement, factor each cyclic order into its pairwise-coprime
prime-power factors and split it in the reverse direction.
\end{proof}

\begin{theorem}[Fundamental theorem, digital form]\label{thm:fundamental}
Every finite Abelian group is isomorphic to a unique direct sum~\eqref{eq:dsum} with radices $p_0\mid p_1\mid\dots\mid p_{N-1}$ and
$p_0\ge2$ (invariant factors); equivalently, by Theorem~\ref{thm:crt}, to a direct sum of cyclic groups of prime-power order. The number of
Abelian groups of order $\prod_i\ell_i^{e_i}$ is $\prod_i\pi(e_i)$, where $\pi$ is the partition function.
\end{theorem}
\begin{proof}
Let $G=(\Ap{n},*)$ be Abelian of order $n=\prod_i\ell_i^{e_i}$; powers $x^{*k}$ with $k<0$ denote powers of the inverse.

\emph{Step 0: $x^{*n}=0$.} The row code $L_x$ is a permutation of $\Ap{n}$, so $\{x*g:g\in\Ap{n}\}=\Ap{n}$. Multiplying all elements
in both forms and using commutativity gives $x^{*n}*\prod_g g=\prod_g g$, hence $x^{*n}=0$.

\emph{Step 1: primary decomposition.} For each prime $\ell=\ell_i$ let $G_\ell=\{x:x^{*\ell^{e_i}}=0\}$, a subgroup because $G$ is
Abelian, and let $m_\ell=n/\ell^{e_i}$. The $m_\ell$ have no common prime factor, so there are integers $u_\ell$ with
$\sum_\ell u_\ell m_\ell=1$. By Step~0, $x^{*u_\ell m_\ell}\in G_\ell$ and $x=\prod_\ell x^{*u_\ell m_\ell}$, so the map
$(x_\ell)_\ell\mapsto\prod_\ell x_\ell$ from $\bigoplus_\ell G_\ell$ to $G$ is a surjective homomorphism. It is injective: if
$\prod_\ell x_\ell=0$, raising to the power $m_\ell$ kills every factor except $x_\ell$, so $x_\ell^{*m_\ell}=0$; since also
$x_\ell^{*\ell^{e_i}}=0$ and $\gcd(m_\ell,\ell^{e_i})=1$, we get $x_\ell=0$. Thus $G\cong\bigoplus_\ell G_\ell$, and every element of
$G_\ell$ has order a power of $\ell$.

\emph{Step 2: an Abelian group $A$ in which every element has $\ell$-power order is a direct sum of cyclic groups.} Induct on
$|A|$. Let $a$ have maximal order $\ell^k$, so $y^{*\ell^k}=0$ for all $y\in A$, and let $K$ be a subgroup maximal with
$K\cap\langle a\rangle=\{0\}$. Suppose $\langle a\rangle K\ne A$. Some $y\notin\langle a\rangle K$ has a power $y^{*\ell^j}\in\langle a\rangle K$; taking $j$
minimal and $x=y^{*\ell^{j-1}}$ gives $x\notin\langle a\rangle K$ and $x^{*\ell}=a^{*t}*c$ with $c\in K$. Then
$0=x^{*\ell^k}=a^{*t\ell^{k-1}}*c^{*\ell^{k-1}}$, so $a^{*t\ell^{k-1}}\in K\cap\langle a\rangle=\{0\}$, hence $\ell\mid t$, say $t=\ell s$. The element
$z=x*a^{*(-s)}$ satisfies $z\notin\langle a\rangle K$ and $z^{*\ell}=c\in K$. Let $K'$ be generated by $K$ and $z$. If
$a^{*r}=c'*z^{*j}$ with $c'\in K$, then either $\ell\mid j$, so $z^{*j}\in K$ and $a^{*r}\in K\cap\langle a\rangle=\{0\}$, or $\gcd(j,\ell)=1$,
and writing $1=\alpha j+\beta\ell$ gives $z=(z^{*j})^{*\alpha}*(z^{*\ell})^{*\beta}\in\langle a\rangle K$, a contradiction. So
$K'\cap\langle a\rangle=\{0\}$ with $K'\supsetneq K$, against maximality. Hence $A=\langle a\rangle K$ with trivial intersection, so
$A\cong\langle a\rangle\oplus K\cong C_{\ell^k}\oplus K$ by Proposition~\ref{prop:cyclic}, and the induction hypothesis applies to $K$.

\emph{Digital form of Steps 1--2.} Steps 1 and 2 give elements $g_h$ of prime-power orders $q_h$, $h<N$, with $\prod_hq_h=n$, and the labelling
\[
\Psi(m)=g_0^{*\md{0}{m}}*g_1^{*\md{1}{m}}*\cdots*g_{N-1}^{*\md{N-1}{m}},\qquad m\in\Ap{n},
\]
is an isomorphism from the direct-sum law~\eqref{eq:dsum} with radices $q_h$ onto $G$: mixed-radix digits add modulo $q_h$
under $\oplus$, and $g_h^{*q_h}=0$. Thus $G$ is a direct sum of cyclic groups of prime-power order.

\emph{Step 3: invariant factors.} For each prime $\ell$ list the exponents of the $\ell$-power radices in decreasing order,
padding with zeros. Let $p_{N'-1}$ be the product over $\ell$ of the largest $\ell$-powers, $p_{N'-2}$ the product of the
second largest, and so on, discarding factors equal to $1$. By Proposition~\ref{prop:ppreduction} each $p_h$ is the order of
a cyclic summand obtained by combining coprime prime-power summands, and by construction $p_0\mid p_1\mid\dots\mid p_{N'-1}$ with
$p_0\ge2$.

\emph{Step 4: uniqueness.} For a prime $\ell$ and $j\ge0$ let
\[
c_\ell(j)=\big|\{x:x^{*\ell^j}=0\}\big|=\sum_{x\in\Ap{n}}d_n\big(x^{*\ell^j},1\big),
\]
an isomorphism invariant, computed with the digital Kronecker delta $d_n(u,1)=\delta_{u,0}$ for $u\in\Ap{n}$ (equivalently $\Bf(u)$, Proposition~\ref{prop:Bvalues}). In $C_{\ell^k}$ the
equation $\ell^jx=0$ has $\ell^{\min(j,k)}$ solutions, by Theorem~\ref{thm:cycsub}; in $C_q$ with $\ell\nmid q$ it has one solution; and
solution counts multiply over direct summands. Hence, if the $\ell$-power summands have exponents $k_1,\dots,k_r$,
\[
c_\ell(j)=\ell^{\sum_i\min(j,k_i)},\qquad \frac{c_\ell(j+1)}{c_\ell(j)}=\ell^{\,|\{i:\ k_i>j\}|},
\]
so the multiset $\{k_i\}$ is determined by $G$. This proves uniqueness of the prime-power decomposition. In an
invariant-factor decomposition the exponent of $\ell$ in $p_h$ is nondecreasing in $h$, so splitting each $p_h$ into prime
powers recovers the sorted exponent lists, and conversely those lists determine the $p_h$ as in Step~3. Hence the invariant
factors are unique too.

\emph{Step 5: counting.} By Steps 1 and 4, the isomorphism class of $G$ is determined by, and determines, for each $i$ a
multiset of positive exponents with sum $e_i$, since $|G_{\ell_i}|=\ell_i^{e_i}$. Every such choice occurs. These multisets are the
partitions of $e_i$, which gives $\prod_i\pi(e_i)$ classes.
\end{proof}

\begin{example}[Order $1176$]
Since
\[
1176=2^3\cdot3\cdot7^2,
\qquad
\pi(3)\pi(1)\pi(2)=3\cdot1\cdot2=6,
\]
the fundamental theorem gives six Abelian isomorphism classes. A prime-power list is
\[
\begin{array}{lll}
C_8\oplus C_3\oplus C_{49},& C_4\oplus C_2\oplus C_3\oplus C_{49},& C_2^3\oplus C_3\oplus C_{49},\\
C_8\oplus C_3\oplus C_7^2,& C_4\oplus C_2\oplus C_3\oplus C_7^2,& C_2^3\oplus C_3\oplus C_7^2.
\end{array}
\]
By Proposition~\ref{prop:ppreduction}, factors of coprime order may be combined; for instance
$C_{1176}\cong C_{49}\oplus C_{24}\cong C_8\oplus C_3\oplus C_{49}$.
\end{example}

\begin{example}[Abelian groups of order at most $32$]
Every order not listed below has only the cyclic group.
\[
\begin{array}{c|l}
p & \text{invariant factor decompositions}\\\hline
4 & C_4,\ C_2^2\\
8 & C_8,\ C_4\oplus C_2,\ C_2^3\\
9 & C_9,\ C_3^2\\
12 & C_{12},\ C_6\oplus C_2\\
16 & C_{16},\ C_8\oplus C_2,\ C_4^2,\ C_4\oplus C_2^2,\ C_2^4\\
18 & C_{18},\ C_6\oplus C_3\\
20 & C_{20},\ C_{10}\oplus C_2\\
24 & C_{24},\ C_{12}\oplus C_2,\ C_6\oplus C_2^2\\
25 & C_{25},\ C_5^2\\
27 & C_{27},\ C_9\oplus C_3,\ C_3^3\\
28 & C_{28},\ C_{14}\oplus C_2\\
32 & C_{32},\ C_{16}\oplus C_2,\ C_8\oplus C_4,\ C_8\oplus C_2^2,\ C_4^2\oplus C_2,\ C_4\oplus C_2^3,\ C_2^5
\end{array}
\]
(Here the summands are written with decreasing radices.)
\end{example}

\section{Exercises}
\begin{enumerate}
\item Show that the subgroups of $C_n$ are the sets $\{0,d,2d,\dots\}$ with $d\mid n$ (compare Theorem~\ref{thm:cycsub}), and that
$C_n$ has exactly $\varphi(d)$ elements of order $d$.
\item Find a generator of the multiplicative group $(\Z/p)^\times$ for $p=7,11,13$, and write the corresponding discrete-logarithm
relabelling of $\Ap{p-1}$.
\item Show that $m=a+3b\mapsto\dg{12}{0}{4a+9b}$ $(a\in\Ap3,\ b\in\Ap4)$ is an isomorphism $C_3\oplus C_4\to C_{12}$, and compare
with Theorem~\ref{thm:crt}.
\item Show that $C_m\oplus C_n\cong C_{mn}$ if and only if $\gcd(m,n)=1$.
\item List the Abelian groups of order $72$ up to isomorphism.
\item Show that $(\Z/2^k)^\times\cong C_2\oplus C_{2^{k-2}}$ for $k\ge3$, generated by $-1$ and $5$, and check the case $k=4$ explicitly.
\item Use Theorem~\ref{thm:qr:zolotarev} to show that $(-1/p)=(-1)^{(p-1)/2}$: the permutation $\pi_{p-1}:x\mapsto\dg p0{-x}$ is a product of $(p-1)/2$ transpositions.
\item Apply Proposition~\ref{prop:qr:gauss} with $a=2$ to show that $r=\frac{p-1}2-\fl{p/4}$ and deduce $(2/p)=(-1)^{(p^2-1)/8}$.
\item Let $C_p$ act by rotation on the words of length $p$ over $\Ap a$, coded by $\Ap{a^p}$. Use Burnside's lemma to show that there are $(a^p+(p-1)a)/p$ orbits, and deduce Fermat's theorem $\dg p0{a^p}=\dg p0a$.
\item Show that the product of all elements of a finite abelian group is the identity unless the group has exactly one element of order $2$, in which case it is that element. Applying this to $U_p$, prove Wilson's theorem $\dg p0{(p-1)!}=p-1$.
\item For distinct odd primes $p,q$, count the lattice points $(j,k)$ with $1\le j\le\frac{p-1}2$, $1\le k\le\frac{q-1}2$ on either side of the line $qx=py$ to show
$\sum_{j=1}^{(p-1)/2}\fl{qj/p}+\sum_{k=1}^{(q-1)/2}\fl{pk/q}=\frac{p-1}2\cdot\frac{q-1}2$. Show that $(q/p)=(-1)^{\sum_j\fl{qj/p}}$ (Eisenstein), and obtain a second proof of Theorem~\ref{thm:qr:reciprocity}.
\item Prove Kummer's theorem: the exponent of $p$ in $\binom{a+b}{a}$ is the number of carries in the base-$p$ addition of $a$ and $b$. Compare the carries with the composition factors in Example~\ref{ex:digitseries}.
\end{enumerate}

\chapter{Nonabelian groups}\label{ch:nonabelian}

\section{Semidirect products}
\begin{theorem}\label{thm:semidirect}
Let $G_0=(\Ap{p_0},*_0)$ and $G_1=(\Ap{p_1},*_1)$ be groups with identity $0$, and let $c\mapsto\phi_c$ be a homomorphism
$G_1\to\Aut(G_0)$. With $m_h=\md{h}{m}$ (radices $p_0,p_1$), the law
\begin{equation}\label{eq:semidirect}
m*n=\big(m_0*_0\phi_{m_1}(n_0)\big)+p_0\big(m_1*_1n_1\big)
\end{equation}
makes $\Ap{p_0p_1}$ a group $G_0\rtimes_\phi G_1$ with identity $0$. It is the direct product if and only if $\phi$ is trivial.
\end{theorem}
\begin{proof}
Both $(m*n)*k$ and $m*(n*k)$ have upper digit $m_1*_1n_1*_1k_1$ and lower digits
$m_0*_0\phi_{m_1}(n_0)*_0\phi_{m_1*_1n_1}(k_0)$ and $m_0*_0\phi_{m_1}\big(n_0*_0\phi_{n_1}(k_0)\big)$; these agree because
$\phi_{m_1}$ is an automorphism and $\phi_{m_1*_1n_1}=\phi_{m_1}\circ\phi_{n_1}$. Since $\phi_0=\mathrm{id}$ and $\phi_c(0)=0$, $0$ is an
identity. The element with digits $\big(\phi_{m_1^{-1}}(m_0^{-1}),m_1^{-1}\big)$ is a right inverse of $m$, and a finite monoid with right
inverses is a group.
\end{proof}

\section{Dicyclic groups}
\begin{theorem}[Dicyclic groups]\label{thm:dicyclic}
Let $q\ge 1$. On $\Ap{4q}$ define, with $\epsilon(m)=d_2\!\left(0,\frac{m}{2q}\right)$,
\begin{equation}\label{eq:dicyclic}
m*n=d_{2q}\!\left(0,m+(-1)^{\epsilon(m)}n+q\,\epsilon(m)\epsilon(n)\right)+2q\,d_2\!\left(0,\epsilon(m)+\epsilon(n)\right).
\end{equation}
Then $\Ap{4q}$ is a group, denoted $Q_{4q}$. It is Abelian only for $q=1$, and for $q=2$ it is the quaternion group $Q_8$.
\end{theorem}
\begin{proof}
Write $m=a+2q\alpha$ and $n=b+2q\beta$ with $a,b\in\Ap{2q}$ and $\alpha,\beta\in\Ap2$.
Then~\eqref{eq:dicyclic} becomes
\[
(a,\alpha)*(b,\beta)=\big(a+(-1)^\alpha b+q\alpha\beta\bmod 2q,\ \alpha+\beta\bmod2\big).
\]
This is the standard semidirect-type model of the dicyclic group with presentation generated by $x=(1,0)$ and $y=(0,1)$,
satisfying $x^{2q}=0$, $y^2=x^q$, and $yxy^{-1}=x^{-1}$. Associativity follows by direct substitution from the law,
or from that presentation; the identity is $(0,0)$ and
$(a,\alpha)^{-1}=((-1)^\alpha(q\alpha-a)\bmod 2q,\alpha)$. If $q=1$, the law reduces to $C_2\oplus C_2$; if $q\ge2$, then
$(0,1)*(1,0)=(2q-1,1)$ while $(1,0)*(0,1)=(1,1)$.
\end{proof}

\begin{example}[Quaternion group]
For $q=2$, Theorem~\ref{thm:dicyclic} gives $Q_8$. The element $x=(1,0)$ has order $4$, the element $y=(0,1)$ satisfies $y^2=x^2$,
and every noncentral element outside $\langle x\rangle$ has order $4$.
\end{example}

\begin{example}[The quaternion group $Q_8$ as a Cayley table]\label{ex:q8}
For $q=2$ the law~\eqref{eq:dicyclic} on $\Ap{8}$ reads
\[
m*n=\dg{4}{0}{m+(-1)^{\epsilon(m)}n+2\,\epsilon(m)\epsilon(n)}+4\,\dg{2}{0}{\epsilon(m)+\epsilon(n)},\qquad
\epsilon(m)=\dg{2}{0}{\tfrac m4}.
\]
With $x=1$ and $y=4$ every element is uniquely $m=x^{m_0}y^{\epsilon(m)}$, where $m_0=\dg{4}{0}{m}$, and $y^2=x^2$, $yxy^{-1}=x^{-1}$.
Identifying $i=1$, $j=4$, $k=5$ and $-1=2$ gives the labels
\[
1\mapsto 0,\quad i\mapsto1,\quad -1\mapsto2,\quad -i\mapsto3,\quad j\mapsto4,\quad k\mapsto5,\quad -j\mapsto6,\quad -k\mapsto7,
\]
and Table~\ref{tab:q8} is the multiplication table of the quaternion units. The element $2$ is the only involution and generates
the centre; the six elements $1,3,4,5,6,7$ all have order $4$.
\end{example}

\begin{table}[htbp]
\centering
\footnotesize
\setlength{\tabcolsep}{3pt}
\renewcommand{\arraystretch}{1.1}
\begin{tabular}{r|*{4}{r}|*{4}{r}}
$*$ & 0 & 1 & 2 & 3 & 4 & 5 & 6 & 7 \\
\hline
0 & 0 & 1 & 2 & 3 & 4 & 5 & 6 & 7 \\
1 & 1 & 2 & 3 & 0 & 5 & 6 & 7 & 4 \\
2 & 2 & 3 & 0 & 1 & 6 & 7 & 4 & 5 \\
3 & 3 & 0 & 1 & 2 & 7 & 4 & 5 & 6 \\
\hline
4 & 4 & 7 & 6 & 5 & 2 & 1 & 0 & 3 \\
5 & 5 & 4 & 7 & 6 & 3 & 2 & 1 & 0 \\
6 & 6 & 5 & 4 & 7 & 0 & 3 & 2 & 1 \\
7 & 7 & 6 & 5 & 4 & 1 & 0 & 3 & 2 \\
\end{tabular}
\caption{Cayley table of $Q_8$ on $\Ap{8}$ given by the digital law~\eqref{eq:dicyclic} with $q=2$; the entry in row $m$, column $n$ is $m*n$. The rules separate the cosets $\{0,\dots,3\}$ and $\{4,\dots,7\}$ of $\langle x\rangle$, i.e.\ the value of $\epsilon(m)$.}
\label{tab:q8}
\end{table}

\begin{example}[Computing an entry of Table~\ref{tab:q8}]\label{ex:q8-entry}
We compute $5*6$ directly from the law, using only the digit function~\eqref{eq:digit}.
\begin{enumerate}
\item \emph{Signs.}
\[
\epsilon(5)=\dg{2}{0}{\tfrac54}=\fl{\tfrac54}-2\fl{\tfrac58}=1-0=1,\qquad
\epsilon(6)=\dg{2}{0}{\tfrac64}=\fl{\tfrac64}-2\fl{\tfrac68}=1-0=1.
\]
\item \emph{Lower digit of the product.} Since $\epsilon(5)=1$ the sign is $(-1)^1=-1$ and the correction is $2\cdot1\cdot1=2$, so
\[
\dg{4}{0}{5-6+2}=\dg{4}{0}{1}=1.
\]
\item \emph{Upper digit of the product.}
\[
\dg{2}{0}{\epsilon(5)+\epsilon(6)}=\dg{2}{0}{2}=2-2\fl{1}=0.
\]
\item \emph{Assembly.} $5*6=1+4\cdot0=1$, which is the entry in row $5$, column $6$ of the table.
\end{enumerate}
\emph{Check via the normal form.} Since $5=xy$, $6=x^2y$ and $yx^2y^{-1}=x^{-2}=x^2$,
\[
5*6=xy\,x^2y=x\,(yx^2y^{-1})\,y^2=x\cdot x^{2}\cdot x^2=x^{5}=x=1.
\]
In quaternion notation this is $k\cdot(-j)=-(kj)=i$.
\emph{Non-commutativity.} In the reverse order both signs are again $1$ and
\[
6*5=\dg{4}{0}{6-5+2}+4\,\dg{2}{0}{2}=\dg{4}{0}{3}+0=3\neq1,
\]
that is, $(-j)\cdot k=-(jk)=-i$.
\end{example}

\begin{example}[The dicyclic group $Q_{12}$]\label{ex:dic12}
Take $q=3$ in Theorem~\ref{thm:dicyclic}. The law on $\Ap{12}$ is
\begin{equation}\label{eq:dic12}
m*n=\dg{6}{0}{m+(-1)^{\epsilon(m)}n+3\,\epsilon(m)\epsilon(n)}+6\,\dg{2}{0}{\epsilon(m)+\epsilon(n)},\qquad
\epsilon(m)=\dg{2}{0}{\tfrac m6}.
\end{equation}
With $x=1$ and $y=6$ every element is uniquely $m=x^{m_0}y^{\epsilon(m)}$, where $m_0=\dg{6}{0}{m}$, and $x^6=0$, $y^2=x^3$,
$yxy^{-1}=x^{-1}$. The resulting Cayley table is Table~\ref{tab:dic12}. The elements $1,5$ have order $6$, the elements $2,4$ have
order $3$, the element $3$ is the unique involution and generates the centre, and the six elements $6,\dots,11$ outside
$\langle x\rangle$ all have order $4$. The group is not a semidirect product $C_6\rtimes C_2$, because its only involution is
central and lies in $\langle x\rangle$.
\end{example}

\begin{table}[htbp]
\centering
\footnotesize
\setlength{\tabcolsep}{3pt}
\renewcommand{\arraystretch}{1.1}
\begin{tabular}{r|*{6}{r}|*{6}{r}}
$*$ & 0 & 1 & 2 & 3 & 4 & 5 & 6 & 7 & 8 & 9 & 10 & 11 \\
\hline
0 & 0 & 1 & 2 & 3 & 4 & 5 & 6 & 7 & 8 & 9 & 10 & 11 \\
1 & 1 & 2 & 3 & 4 & 5 & 0 & 7 & 8 & 9 & 10 & 11 & 6 \\
2 & 2 & 3 & 4 & 5 & 0 & 1 & 8 & 9 & 10 & 11 & 6 & 7 \\
3 & 3 & 4 & 5 & 0 & 1 & 2 & 9 & 10 & 11 & 6 & 7 & 8 \\
4 & 4 & 5 & 0 & 1 & 2 & 3 & 10 & 11 & 6 & 7 & 8 & 9 \\
5 & 5 & 0 & 1 & 2 & 3 & 4 & 11 & 6 & 7 & 8 & 9 & 10 \\
\hline
6 & 6 & 11 & 10 & 9 & 8 & 7 & 3 & 2 & 1 & 0 & 5 & 4 \\
7 & 7 & 6 & 11 & 10 & 9 & 8 & 4 & 3 & 2 & 1 & 0 & 5 \\
8 & 8 & 7 & 6 & 11 & 10 & 9 & 5 & 4 & 3 & 2 & 1 & 0 \\
9 & 9 & 8 & 7 & 6 & 11 & 10 & 0 & 5 & 4 & 3 & 2 & 1 \\
10 & 10 & 9 & 8 & 7 & 6 & 11 & 1 & 0 & 5 & 4 & 3 & 2 \\
11 & 11 & 10 & 9 & 8 & 7 & 6 & 2 & 1 & 0 & 5 & 4 & 3 \\
\end{tabular}
\caption{Cayley table of $Q_{12}$ on $\Ap{12}$ given by the digital law~\eqref{eq:dic12}; the entry in row $m$, column $n$ is $m*n$. The rules separate the cosets $\{0,\dots,5\}$ and $\{6,\dots,11\}$ of $\langle x\rangle$, i.e.\ the value of $\epsilon(m)$.}
\label{tab:dic12}
\end{table}

\begin{example}[Computing an entry of Table~\ref{tab:dic12}]\label{ex:dic12-entry}
We compute $8*10$ directly from~\eqref{eq:dic12}, using only the digit function~\eqref{eq:digit}.
\begin{enumerate}
\item \emph{Signs.}
\begin{gather*}
\epsilon(8)=\dg{2}{0}{\tfrac86}=\fl{\tfrac86}-2\fl{\tfrac8{12}}=1-0=1,\\
\epsilon(10)=\dg{2}{0}{\tfrac{10}6}=\fl{\tfrac{10}6}-2\fl{\tfrac{10}{12}}=1-0=1.
\end{gather*}
\item \emph{Lower digit of the product.} The sign is $(-1)^1=-1$ and the correction is $3\cdot1\cdot1=3$, so
\[
\dg{6}{0}{8-10+3}=\dg{6}{0}{1}=\fl{1}-6\fl{\tfrac16}=1.
\]
By Lemma~\ref{lem:mod} one may reduce first: $m_0=2$, $n_0=4$, and $\dg{6}{0}{2-4+3}=1$.
\item \emph{Upper digit of the product.}
\[
\dg{2}{0}{\epsilon(8)+\epsilon(10)}=\dg{2}{0}{2}=2-2\fl{1}=0.
\]
\item \emph{Assembly.} $8*10=1+6\cdot0=1$, which is the entry in row $8$, column $10$ of the table.
\end{enumerate}
\emph{Check via the normal form.} Since $8=x^2y$, $10=x^4y$, $yx^4y^{-1}=x^{-4}$ and $y^2=x^3$,
\[
8*10=x^2y\,x^4y=x^2\,(yx^4y^{-1})\,y^2=x^{2-4+3}=x=1.
\]
\emph{Non-commutativity.} In the reverse order both signs are again $1$ and
\[
10*8=\dg{6}{0}{10-8+3}+6\,\dg{2}{0}{2}=\dg{6}{0}{5}+0=5\neq1 .
\]
\end{example}

\section{Dihedral and metacyclic groups}
\begin{theorem}[Metacyclic groups]\label{thm:metacyclic}
Let $q,s\ge2$ and let $r$ satisfy $\gcd(r,q)=1$ and $r^{s}\equiv1\pmod q$. With $m_1=\dg{s}{0}{m/q}$, the law
\begin{equation}\label{eq:metacyclic}
m*n=\dg{q}{0}{m+r^{m_1}n}+q\,\dg{s}{0}{m_1+n_1}
\end{equation}
makes $\Ap{qs}$ a group $C_q\rtimes_r C_s$. It is nonabelian if and only if $r\not\equiv1\pmod q$.
\end{theorem}
\begin{proof}
Take $\phi_c(x)=\dg{q}{0}{r^{c}x}$ in Theorem~\ref{thm:semidirect}: it is an automorphism of $C_q$ because $r$ is a unit, and
$c\mapsto\phi_c$ is a homomorphism from $C_s$ because $r^s\equiv1$. By Lemma~\ref{lem:mod}, $m\equiv m_0$ and $n\equiv n_0 \pmod q$, so
\eqref{eq:semidirect} becomes~\eqref{eq:metacyclic}. Finally $1*q=1+q$ while $q*1=\dg{q}{0}{r}+q$.
\end{proof}

\begin{corollary}[Dihedral groups]\label{cor:dihedral}
For $q\ge 2$, $r=q-1$, $s=2$, the law
\begin{equation}\label{eq:dihedral}
m*n=\dg{q}{0}{m+n\,(q-1)^{\dg{2}{0}{m/q}}}+q\,\dg{2}{0}{\dg{2}{0}{\tfrac{m}{q}}+\dg{2}{0}{\tfrac{n}{q}}}
\end{equation}
defines the dihedral group $D_{2q}$ of order $2q$. It is nonabelian if and only if $q\ge3$; $D_2\cong C_2$ and $D_4\cong C_2\oplus C_2$.
\end{corollary}

\begin{remark}
The dihedral law is exactly the semidirect product $C_q\rtimes C_2$ in which the nontrivial element of $C_2$ acts on $C_q$ by inversion.
Equivalently, writing $m=a+q\alpha$ and $n=b+q\beta$ with $a,b\in\Ap q$ and $\alpha,\beta\in\Ap2$, one has
$(a,\alpha)*(b,\beta)=(a+(-1)^\alpha b\bmod q,\ \alpha+\beta\bmod2)$. This makes transparent why $D_6\cong \Sym_3$ and why the law is
noncommutative as soon as $q\ge3$.
\end{remark}

\begin{proposition}[Presentations and generators]\label{prop:presentations}
In $D_{2q}$ let $r=(1,0)$ and $s=(0,1)$ in the two-digit model of the preceding remark. Then
\[
r^q=s^2=1,\qquad srs^{-1}=r^{-1},
\]
and every element is uniquely $r^a$ or $r^as$ with $a\in\Ap{q}$. Hence
\[
D_{2q}=\langle r,s\mid r^q=s^2=1,\ srs^{-1}=r^{-1}\rangle
=\langle r,s\mid r^q=s^2=(rs)^2=1\rangle.
\]
In $Q_{4q}$, with the generators $x=(1,0)$ and $y=(0,1)$ of Theorem~\ref{thm:dicyclic},
\[
Q_{4q}=\langle x,y\mid x^{2q}=1,\ y^2=x^q,\ yxy^{-1}=x^{-1}\rangle.
\]
\end{proposition}
\begin{proof}
The displayed digit laws give $r^a=(a,0)$, $r^as=(a,1)$, $x^a=(a,0)$ and $x^ay=(a,1)$ directly.
Substitution in the respective laws verifies the relations; normal forms give the stated orders and uniqueness.
\end{proof}

\begin{remark}[Geometric form]
$D_{2q}$ is the symmetry group of a regular $q$-gon: $r$ is rotation by one vertex and $s$ is a reflection. Thus $D_6$ is the
symmetry group of an equilateral triangle and $D_8$ that of a square. The first non-Abelian dicyclic group is $Q_8$, not $D_8$:
the two groups have the same order but different element-order structure.
\end{remark}

\begin{example}[The metacyclic group $C_7\rtimes_2 C_3$]\label{ex:metacyclic21}
Take $q=7$, $s=3$ and $r=2$ in Theorem~\ref{thm:metacyclic}; the hypotheses hold because $\gcd(2,7)=1$ and
$2^{3}=8\equiv1\pmod 7$. Writing out $m_1=\dg{3}{0}{m/7}$ and $n_1=\dg{3}{0}{n/7}$, the law~\eqref{eq:metacyclic} on $\Ap{21}$ is
\begin{equation}\label{eq:metacyclic21}
m*n=\dg{7}{0}{m+2^{\dg{3}{0}{m/7}}\,n}+7\,\dg{3}{0}{\dg{3}{0}{\tfrac{m}{7}}+\dg{3}{0}{\tfrac{n}{7}}}.
\end{equation}
With $x=1$ and $y=7$ every element is uniquely $m=x^{m_0}y^{m_1}$, where $m_0=\dg{7}{0}{m}$, and $yxy^{-1}=x^{2}$.
The resulting Cayley table is Table~\ref{tab:metacyclic21}. The six elements $1,\dots,6$ have order $7$, the fourteen elements
$7,\dots,20$ have order $3$, and the centre is trivial: this is the Frobenius group of order $21$, the smallest
non-Abelian group of odd order.
\end{example}

\begin{table}[htbp]
\centering
\footnotesize
\setlength{\tabcolsep}{2.2pt}
\renewcommand{\arraystretch}{1.1}
\begin{tabular}{r|*{7}{r}|*{7}{r}|*{7}{r}}
$*$ & 0 & 1 & 2 & 3 & 4 & 5 & 6 & 7 & 8 & 9 & 10 & 11 & 12 & 13 & 14 & 15 & 16 & 17 & 18 & 19 & 20 \\
\hline
0 & 0 & 1 & 2 & 3 & 4 & 5 & 6 & 7 & 8 & 9 & 10 & 11 & 12 & 13 & 14 & 15 & 16 & 17 & 18 & 19 & 20 \\
1 & 1 & 2 & 3 & 4 & 5 & 6 & 0 & 8 & 9 & 10 & 11 & 12 & 13 & 7 & 15 & 16 & 17 & 18 & 19 & 20 & 14 \\
2 & 2 & 3 & 4 & 5 & 6 & 0 & 1 & 9 & 10 & 11 & 12 & 13 & 7 & 8 & 16 & 17 & 18 & 19 & 20 & 14 & 15 \\
3 & 3 & 4 & 5 & 6 & 0 & 1 & 2 & 10 & 11 & 12 & 13 & 7 & 8 & 9 & 17 & 18 & 19 & 20 & 14 & 15 & 16 \\
4 & 4 & 5 & 6 & 0 & 1 & 2 & 3 & 11 & 12 & 13 & 7 & 8 & 9 & 10 & 18 & 19 & 20 & 14 & 15 & 16 & 17 \\
5 & 5 & 6 & 0 & 1 & 2 & 3 & 4 & 12 & 13 & 7 & 8 & 9 & 10 & 11 & 19 & 20 & 14 & 15 & 16 & 17 & 18 \\
6 & 6 & 0 & 1 & 2 & 3 & 4 & 5 & 13 & 7 & 8 & 9 & 10 & 11 & 12 & 20 & 14 & 15 & 16 & 17 & 18 & 19 \\
\hline
7 & 7 & 9 & 11 & 13 & 8 & 10 & 12 & 14 & 16 & 18 & 20 & 15 & 17 & 19 & 0 & 2 & 4 & 6 & 1 & 3 & 5 \\
8 & 8 & 10 & 12 & 7 & 9 & 11 & 13 & 15 & 17 & 19 & 14 & 16 & 18 & 20 & 1 & 3 & 5 & 0 & 2 & 4 & 6 \\
9 & 9 & 11 & 13 & 8 & 10 & 12 & 7 & 16 & 18 & 20 & 15 & 17 & 19 & 14 & 2 & 4 & 6 & 1 & 3 & 5 & 0 \\
10 & 10 & 12 & 7 & 9 & 11 & 13 & 8 & 17 & 19 & 14 & 16 & 18 & 20 & 15 & 3 & 5 & 0 & 2 & 4 & 6 & 1 \\
11 & 11 & 13 & 8 & 10 & 12 & 7 & 9 & 18 & 20 & 15 & 17 & 19 & 14 & 16 & 4 & 6 & 1 & 3 & 5 & 0 & 2 \\
12 & 12 & 7 & 9 & 11 & 13 & 8 & 10 & 19 & 14 & 16 & 18 & 20 & 15 & 17 & 5 & 0 & 2 & 4 & 6 & 1 & 3 \\
13 & 13 & 8 & 10 & 12 & 7 & 9 & 11 & 20 & 15 & 17 & 19 & 14 & 16 & 18 & 6 & 1 & 3 & 5 & 0 & 2 & 4 \\
\hline
14 & 14 & 18 & 15 & 19 & 16 & 20 & 17 & 0 & 4 & 1 & 5 & 2 & 6 & 3 & 7 & 11 & 8 & 12 & 9 & 13 & 10 \\
15 & 15 & 19 & 16 & 20 & 17 & 14 & 18 & 1 & 5 & 2 & 6 & 3 & 0 & 4 & 8 & 12 & 9 & 13 & 10 & 7 & 11 \\
16 & 16 & 20 & 17 & 14 & 18 & 15 & 19 & 2 & 6 & 3 & 0 & 4 & 1 & 5 & 9 & 13 & 10 & 7 & 11 & 8 & 12 \\
17 & 17 & 14 & 18 & 15 & 19 & 16 & 20 & 3 & 0 & 4 & 1 & 5 & 2 & 6 & 10 & 7 & 11 & 8 & 12 & 9 & 13 \\
18 & 18 & 15 & 19 & 16 & 20 & 17 & 14 & 4 & 1 & 5 & 2 & 6 & 3 & 0 & 11 & 8 & 12 & 9 & 13 & 10 & 7 \\
19 & 19 & 16 & 20 & 17 & 14 & 18 & 15 & 5 & 2 & 6 & 3 & 0 & 4 & 1 & 12 & 9 & 13 & 10 & 7 & 11 & 8 \\
20 & 20 & 17 & 14 & 18 & 15 & 19 & 16 & 6 & 3 & 0 & 4 & 1 & 5 & 2 & 13 & 10 & 7 & 11 & 8 & 12 & 9 \\
\end{tabular}
\caption{Cayley table of $C_7\rtimes_2 C_3$ on $\Ap{21}$ given by the digital law~\eqref{eq:metacyclic21}; the entry in row $m$,
column $n$ is $m*n$. The rules separate the cosets $\{0,\dots,6\}$, $\{7,\dots,13\}$, $\{14,\dots,20\}$ of
$\langle x\rangle$, i.e.\ the value of the upper digit $m_1$.}
\label{tab:metacyclic21}
\end{table}

\begin{example}[Computing an entry of Table~\ref{tab:metacyclic21}]\label{ex:metacyclic21-entry}
We compute $12*17$ directly from~\eqref{eq:metacyclic21}, using only the digit function~\eqref{eq:digit}.
\begin{enumerate}
\item \emph{Upper digits.}
\[
m_1=\dg{3}{0}{\tfrac{12}{7}}=\fl{\tfrac{12}{7}}-3\fl{\tfrac{12}{21}}=1-0=1,\qquad
n_1=\dg{3}{0}{\tfrac{17}{7}}=\fl{\tfrac{17}{7}}-3\fl{\tfrac{17}{21}}=2-0=2.
\]
\item \emph{Twisting factor.} $r^{m_1}=2^{1}=2$.
\item \emph{Lower digit of the product.}
\[
\dg{7}{0}{12+2\cdot17}=\dg{7}{0}{46}=\fl{46}-7\fl{\tfrac{46}{7}}=46-42=4.
\]
By Lemma~\ref{lem:mod} one may reduce first: $m_0=5$, $n_0=3$, and $\dg{7}{0}{5+2\cdot3}=\dg{7}{0}{11}=4$.
\item \emph{Upper digit of the product.}
\[
\dg{3}{0}{m_1+n_1}=\dg{3}{0}{3}=3-3\fl{1}=0.
\]
\item \emph{Assembly.} $12*17=4+7\cdot0=4$, which is the entry in row $12$, column $17$ of the table.
\end{enumerate}
\emph{Check via the normal form.} Since $12=x^{5}y$, $17=x^{3}y^{2}$ and $yx^{3}y^{-1}=x^{6}$,
\[
12*17=x^{5}y\,x^{3}y^{2}=x^{5}\,(yx^{3}y^{-1})\,y^{3}=x^{11}=x^{4}=4.
\]
\emph{Non-commutativity.} In the reverse order $m_1=2$, so $r^{m_1}=4$ and
\[
17*12=\dg{7}{0}{17+4\cdot12}+7\,\dg{3}{0}{2+1}=\dg{7}{0}{65}+0=2\neq4 .
\]
\end{example}

\begin{proposition}[Bracelets]\label{prop:bracelets}
Let $D_{2q}$ act on the vertices $\Ap{q}$ of a regular $q$-gon, generated by $x\mapsto x+1$ and $x\mapsto-x$ modulo $q$, and let $N_q(m)=\frac1q\sum_{d\mid q}\varphi(d)m^{q/d}$ be the number of necklaces (Proposition~\ref{prop:necklaces}). The number of $m$-colourings of the vertices up to $D_{2q}$ is
\[
\frac12N_q(m)+\begin{cases}\tfrac12\,m^{(q+1)/2},&q\text{ odd},\\[2pt] \tfrac14\,m^{q/2}(m+1),&q\text{ even}.\end{cases}
\]
\end{proposition}
\begin{proof}
The $q$ rotations contribute $\frac1{2q}\sum_{\text{rot}}m^{c}=\frac12N_q(m)$. For odd $q$ each of the $q$ reflections fixes one vertex and swaps the others in pairs, so it has $(q+1)/2$ cycles, contributing $q\,m^{(q+1)/2}/(2q)$. For even $q$, half of the reflections fix two vertices ($q/2+1$ cycles) and half fix none ($q/2$ cycles), contributing $\frac q2\big(m^{q/2+1}+m^{q/2}\big)/(2q)$. Apply Theorem~\ref{thm:polya}.
\end{proof}
For $q=4,6$ and $m=2$ one gets $6$ and $13$ bracelets, and for $q=4$, $m=3$ one gets $21$, as in the example of the square after Theorem~\ref{thm:debruijn}.

\section{Characters of metacyclic and dicyclic groups}\label{sec:metacyclic-chars}
Throughout this section $\zeta_n=e^{2\pi i/n}$, $\omega=\zeta_3$, and characters are ordinary complex characters.
Let $G=C_q\rtimes_r C_s$ be the group of Theorem~\ref{thm:metacyclic} and put $x=1$, $y=q$. Since $r^0=1$, the law
\eqref{eq:metacyclic} gives $x^a=a$ and $x^ay^\alpha=a+q\alpha$ for $a\in\Ap q$ and $\alpha\in\Ap s$, so the two digits
\[
m_0=\dg{q}{0}{m},\qquad m_1=\dg{s}{0}{m/q}
\]
of $m$ are exactly the exponents in the normal form $m=x^{m_0}y^{m_1}$. Moreover $yxy^{-1}=x^r$, and $A=\langle x\rangle=\Ap q$ is a
cyclic normal subgroup with cyclic quotient $G/A\cong C_s$. The characters of $A$ are $\lambda_j(x^a)=\zeta_q^{ja}$, $j\in\Z/q$,
and conjugation by $y$ acts on them by $\lambda_j^{\,y}(x^a)=\lambda_j(yx^ay^{-1})=\lambda_{rj}(x^a)$. Hence the cyclic group
$\langle r\rangle\le(\Z/q)^\times$ acts on $\Z/q$ by multiplication; for an orbit $O=j\langle r\rangle$ write $e_O=|O|$, which divides
$s$, and $t_O=s/e_O$.

\begin{theorem}[Characters of metacyclic groups]\label{thm:metacyclic-chars}
For every orbit $O$ of $\langle r\rangle$ on $\Z/q$ and every $\mu\in\mathbb{C}$ with $\mu^{t_O}=1$, the function
\begin{equation}\label{eq:metacyclic-char}
\chi_{O,\mu}(m)=
\begin{cases}
\mu^{m_1/e_O}\displaystyle\sum_{k\in O}\zeta_q^{\,k m_0}, & e_O\mid m_1,\\[1ex]
0, & e_O\nmid m_1,
\end{cases}
\end{equation}
is an irreducible character of $C_q\rtimes_r C_s$ of degree $e_O$. These characters are pairwise distinct and exhaust
$\operatorname{Irr}(C_q\rtimes_r C_s)$. In particular
\[
k(C_q\rtimes_r C_s)=\sum_{O}\frac{s}{e_O},
\]
and the linear characters are the $\chi_{\{j\},\mu}$ with $(r-1)j\equiv0\pmod q$.
\end{theorem}
\begin{proof}
Fix $j\in O$, $e=e_O$, $t=t_O$, and let $S=\{m: e\mid m_1\}=\langle x,y^e\rangle$, a normal subgroup of index $e$ containing $A$.
By~\eqref{eq:metacyclic}, $x^ay^{e\beta}*x^by^{e\beta'}=x^{a+r^{e\beta}b}y^{e(\beta+\beta')}$, and $\zeta_q^{jr^{e\beta}b}=\zeta_q^{jb}$ because
$r^ej\equiv j$. As $\mu^t=1$, the exponent $\beta$ may be read modulo $t$, so $\lambda(x^ay^{e\beta})=\zeta_q^{ja}\mu^\beta$ is a linear
character of $S$ extending $\lambda_j$. Induce it over the transversal $\{y^c:0\le c<e\}$. Since $y^cx^ay^{e\beta}y^{-c}=x^{r^ca}y^{e\beta}$,
\[
\operatorname{Ind}_S^G\lambda\,(x^ay^{e\beta})=\mu^\beta\sum_{c=0}^{e-1}\zeta_q^{\,jr^ca}=\mu^\beta\sum_{k\in O}\zeta_q^{\,ka},
\]
while the induced character vanishes off the normal subgroup $S$; this is~\eqref{eq:metacyclic-char}. Its restriction to $S$ is the sum
of the $e$ conjugates $\lambda^{y^c}$, which are distinct because their restrictions $\lambda_{r^cj}$ to $A$ are distinct. Frobenius
reciprocity gives $\langle\chi_{O,\mu},\chi_{O,\mu}\rangle=\langle\lambda,\chi_{O,\mu}|_S\rangle_S=1$, so $\chi_{O,\mu}$ is irreducible.
The restriction $\sum_{k\in O}\lambda_k$ of $\chi_{O,\mu}$ to $A$ determines $O$, and $\chi_{O,\mu}(y^{e})=e\mu$ determines $\mu$.
Finally $\sum_O t_Oe_O^2=\sum_O s\,e_O=sq=|G|$, so no irreducible character is missing.
\end{proof}

The formula is digital in the strict sense: $\chi_{O,\mu}(m)$ is computed from the two digits $m_0$, $m_1$ of $m$ alone, without
multiplying in the group.

\begin{corollary}\label{cor:metacyclic-prime-chars}
Let $q$ be prime and let $r$ have multiplicative order exactly $s$ modulo $q$. Then $C_q\rtimes_r C_s$ has the $s$ linear characters
$m\mapsto\mu^{m_1}$ $(\mu^s=1)$ and $(q-1)/s$ characters of degree $s$, which vanish off $A$ and take the values
$\sum_{k\in j\langle r\rangle}\zeta_q^{km_0}$ on $A$. Hence $k(C_q\rtimes_r C_s)=s+(q-1)/s$.
\end{corollary}
\begin{proof}
The orbits are $\{0\}$ and the $(q-1)/s$ cosets of $\langle r\rangle$ in $(\Z/q)^\times$, each of size $s$.
\end{proof}

\begin{example}
For $C_5\rtimes_2C_4$ (order $20$) the nonzero residues form one orbit, so there is a single character of degree $4$; it equals $4$ at
$0$, $\zeta_5+\zeta_5^2+\zeta_5^3+\zeta_5^4=-1$ on $A\setminus\{0\}$ and $0$ elsewhere. For $C_7\rtimes_2C_3$ (order $21$) the orbits
$\{1,2,4\}$ and $\{3,5,6\}$ give two characters of degree $3$ whose values on $x$ are the quadratic Gauss periods
$\zeta_7+\zeta_7^2+\zeta_7^4=\tfrac{-1+i\sqrt7}{2}$ and its complex conjugate. See Tables~\ref{tab:nonab-chars-C5x2C4}
and~\ref{tab:nonab-chars-C7x2C3}.
\end{example}

\begin{corollary}[Characters of dihedral groups]\label{cor:dihedral-chars}
Let $q\ge3$ and $D_{2q}=C_q\rtimes_{q-1}C_2$ as in Corollary~\ref{cor:dihedral}. Its irreducible characters are the linear
characters $m\mapsto(\pm1)^{m_1}$, together with $m\mapsto(-1)^{m_0}(\pm1)^{m_1}$ when $q$ is even, and the characters of degree $2$
\[
\psi_j(m)=\begin{cases}2\cos(2\pi jm_0/q),& m_1=0,\\ 0,& m_1=1,\end{cases}\qquad 1\le j<q/2.
\]
Thus $k(D_{2q})=(q+3)/2$ for odd $q$ and $(q+6)/2$ for even $q$, and every irreducible character of $D_{2q}$ has
Frobenius--Schur indicator $+1$.
\end{corollary}
\begin{proof}
Here $r=-1$, whose orbits on $\Z/q$ are $\{0\}$, $\{q/2\}$ for even $q$, and $\{j,-j\}$ for $1\le j<q/2$; apply
Theorem~\ref{thm:metacyclic-chars}. Every element $x^ay$ is an involution, so the indicator
$\nu(\psi_j)=\frac1{2q}\sum_{g}\psi_j(g^2)$ equals $\frac1{2q}\big(\sum_{a}2\cos(4\pi ja/q)+q\,\psi_j(0)\big)=1$, the first sum
vanishing because $q\nmid 2j$. The linear characters are real, hence have indicator $1$.
\end{proof}

The dicyclic group is not a split extension of $\langle x\rangle$ for $q\ge2$, so Theorem~\ref{thm:metacyclic-chars} does not apply
verbatim, but the same induction argument does.

\begin{proposition}[Characters of dicyclic groups]\label{prop:dicyclic-chars}
Let $q\ge2$ and write $m\in Q_{4q}$ as $m=x^{m_0}y^{m_1}$ with $m_0=\dg{2q}{0}{m}$ and $m_1=\dg{2}{0}{m/2q}$, where $x=1$ and
$y=2q$ as in Theorem~\ref{thm:dicyclic}. Then $\operatorname{Irr}(Q_{4q})$ consists of
\begin{enumerate}
\item the four linear characters $\lambda(m)=\varepsilon^{m_0}\delta^{m_1}$ with $\varepsilon=\pm1$ and $\delta^2=\varepsilon^q$;
\item the $q-1$ characters of degree $2$
\[
\psi_j(m)=\begin{cases}2\cos(\pi jm_0/q),& m_1=0,\\ 0,& m_1=1,\end{cases}\qquad 1\le j\le q-1.
\]
\end{enumerate}
Hence $k(Q_{4q})=q+3$. The indicator of $\psi_j$ is $(-1)^j$; a linear character has indicator $1$ if it is real and $0$ if
$\delta=\pm i$, which happens exactly when $q$ is odd and $\varepsilon=-1$.
\end{proposition}
\begin{proof}
The subgroup $A=\langle x\rangle=\Ap{2q}$ is cyclic of index $2$ and $yxy^{-1}=x^{-1}$ by Proposition~\ref{prop:presentations}.
For $j\in\Z/2q\setminus\{0,q\}$ the characters $\lambda_j$ and $\lambda_j^{\,y}=\lambda_{-j}$ differ, so
$\psi_j=\operatorname{Ind}_A^G\lambda_j$ is irreducible by Frobenius reciprocity, with the stated values. Since $\psi_j=\psi_{2q-j}$ and the
restrictions $\lambda_j+\lambda_{-j}$ are distinct for $1\le j\le q-1$, this gives $q-1$ characters of degree $2$. A linear character
is determined by $\varepsilon=\lambda(x)$ and $\delta=\lambda(y)$; the relations of Proposition~\ref{prop:presentations} hold if and
only if $\varepsilon^2=1$ and $\delta^2=\varepsilon^q$, which gives four characters. As $4\cdot1+(q-1)\cdot4=4q$, the list is complete.
For the indicators note that $(x^ay)^2=x^q$ by~\eqref{eq:dicyclic}, so
\[
\nu(\psi_j)=\frac1{4q}\Big(\sum_{a\in\Ap{2q}}\psi_j(x^{2a})+2q\,\psi_j(x^q)\Big)=\frac{0+2q\cdot2\cos(\pi j)}{4q}=(-1)^j,
\]
the sum vanishing because $q\nmid j$. A linear character has indicator $1$ if $\lambda^2=1$ and $0$ otherwise.
\end{proof}

\begin{example}[$D_8$ versus $Q_8$]
By Corollary~\ref{cor:dihedral-chars} and Proposition~\ref{prop:dicyclic-chars}, $D_8$ and $Q_8$ both have four linear characters
and one character of degree $2$, and their character tables coincide. The indicators separate them: the degree-$2$ character of $D_8$
has $\nu=1$, while $\psi_1$ of $Q_8$ has $\nu=-1$, so it is afforded by no real representation. More generally $Q_{4q}$ has exactly
$\lfloor q/2\rfloor$ irreducible characters with indicator $-1$.
\end{example}

Table~\ref{tab:nonab-char-summary} lists the character data of nine examples and Tables~\ref{tab:nonab-chars-Q16},
\ref{tab:nonab-chars-C5x2C4} and~\ref{tab:nonab-chars-C7x2C3} display three complete character tables. All of them are generated
from the digit laws by \texttt{scripts/\allowbreak nonab\_char\_verify.py}, which builds each group, computes its classes and the characters of
Theorem~\ref{thm:metacyclic-chars} and Proposition~\ref{prop:dicyclic-chars}, verifies the orthogonality relations exactly in
$\Z[\zeta_N]$, and checks every character against GAP's \texttt{Irr}.


\begin{table}[htbp]
\centering\small
\caption{Character data of the metacyclic, dihedral and dicyclic examples: law parameters, class number $k(G)$, number $|G:G'|$ of linear characters, degree multiset and numbers of characters with Frobenius--Schur indicator $+1,0,-1$.}
\label{tab:nonab-char-summary}
\begin{tabular}{llcccc}
\hline
$G$ & law & $k(G)$ & $|G:G'|$ & degrees & $\nu=+1,0,-1$\\
\hline
$D_{6}$ & $q=3$, $r=2$, $s=2$ & 3 & 2 & $1^{2}\,2^{1}$ & 3, 0, 0\\
$D_{8}$ & $q=4$, $r=3$, $s=2$ & 5 & 4 & $1^{4}\,2^{1}$ & 5, 0, 0\\
$D_{10}$ & $q=5$, $r=4$, $s=2$ & 4 & 2 & $1^{2}\,2^{2}$ & 4, 0, 0\\
$D_{12}$ & $q=6$, $r=5$, $s=2$ & 6 & 4 & $1^{4}\,2^{2}$ & 6, 0, 0\\
$Q_{8}$ & $q=2$ & 5 & 4 & $1^{4}\,2^{1}$ & 4, 0, 1\\
$Q_{12}$ & $q=3$ & 6 & 4 & $1^{4}\,2^{2}$ & 3, 2, 1\\
$Q_{16}$ & $q=4$ & 7 & 4 & $1^{4}\,2^{3}$ & 5, 0, 2\\
$C_{5}\rtimes_{2}C_{4}$ & $q=5$, $r=2$, $s=4$ & 5 & 4 & $1^{4}\,4^{1}$ & 3, 2, 0\\
$C_{7}\rtimes_{2}C_{3}$ & $q=7$, $r=2$, $s=3$ & 5 & 3 & $1^{3}\,3^{2}$ & 1, 4, 0\\
\hline
\end{tabular}
\end{table}


\begin{table}[htbp]
\centering\small
\caption{Character table of $Q_{16}$ ($q=4$). Columns are indexed by class representatives $x^ay^\alpha=a+|A|\alpha$; $\nu$ is the Frobenius--Schur indicator.}
\label{tab:nonab-chars-Q16}
\begin{tabular}{l|ccccccc|c}
\hline
 & $1$ & $x^{4}$ & $x^{2}$ & $y$ & $xy$ & $x$ & $x^{3}$ & \\
$|K|$ & 1 & 1 & 2 & 4 & 4 & 2 & 2 & \\
$o(g)$ & 1 & 2 & 4 & 4 & 4 & 8 & 8 & $\nu$\\
\hline
$\chi_{1}$ & $1$ & $1$ & $1$ & $1$ & $1$ & $1$ & $1$ & $1$\\
$\chi_{2}$ & $1$ & $1$ & $1$ & $-1$ & $-1$ & $1$ & $1$ & $1$\\
$\chi_{3}$ & $1$ & $1$ & $1$ & $-1$ & $1$ & $-1$ & $-1$ & $1$\\
$\chi_{4}$ & $1$ & $1$ & $1$ & $1$ & $-1$ & $-1$ & $-1$ & $1$\\
$\chi_{5}$ & $2$ & $-2$ & $0$ & $0$ & $0$ & $-\sqrt{2}$ & $\sqrt{2}$ & $-1$\\
$\chi_{6}$ & $2$ & $-2$ & $0$ & $0$ & $0$ & $\sqrt{2}$ & $-\sqrt{2}$ & $-1$\\
$\chi_{7}$ & $2$ & $2$ & $-2$ & $0$ & $0$ & $0$ & $0$ & $1$\\
\hline
\end{tabular}
\end{table}

\begin{table}[htbp]
\centering\small
\caption{Character table of $C_{5}\rtimes_{2}C_{4}$ ($q=5$, $r=2$, $s=4$). Columns are indexed by class representatives $x^ay^\alpha=a+|A|\alpha$; $\nu$ is the Frobenius--Schur indicator.}
\label{tab:nonab-chars-C5x2C4}
\begin{tabular}{l|ccccc|c}
\hline
 & $1$ & $y^{2}$ & $y$ & $y^{3}$ & $x$ & \\
$|K|$ & 1 & 5 & 5 & 5 & 4 & \\
$o(g)$ & 1 & 2 & 4 & 4 & 5 & $\nu$\\
\hline
$\chi_{1}$ & $1$ & $1$ & $1$ & $1$ & $1$ & $1$\\
$\chi_{2}$ & $1$ & $-1$ & $-i$ & $i$ & $1$ & $0$\\
$\chi_{3}$ & $1$ & $-1$ & $i$ & $-i$ & $1$ & $0$\\
$\chi_{4}$ & $1$ & $1$ & $-1$ & $-1$ & $1$ & $1$\\
$\chi_{5}$ & $4$ & $0$ & $0$ & $0$ & $-1$ & $1$\\
\hline
\end{tabular}
\end{table}

\begin{table}[htbp]
\centering\small
\caption{Character table of $C_{7}\rtimes_{2}C_{3}$ ($q=7$, $r=2$, $s=3$). Columns are indexed by class representatives $x^ay^\alpha=a+|A|\alpha$; $\nu$ is the Frobenius--Schur indicator.}
\label{tab:nonab-chars-C7x2C3}
\begin{tabular}{l|ccccc|c}
\hline
 & $1$ & $y$ & $y^{2}$ & $x$ & $x^{3}$ & \\
$|K|$ & 1 & 7 & 7 & 3 & 3 & \\
$o(g)$ & 1 & 3 & 3 & 7 & 7 & $\nu$\\
\hline
$\chi_{1}$ & $1$ & $1$ & $1$ & $1$ & $1$ & $1$\\
$\chi_{2}$ & $1$ & $\bar\omega$ & $\omega$ & $1$ & $1$ & $0$\\
$\chi_{3}$ & $1$ & $\omega$ & $\bar\omega$ & $1$ & $1$ & $0$\\
$\chi_{4}$ & $3$ & $0$ & $0$ & $\tfrac{-1+i\sqrt7}{2}$ & $\tfrac{-1-i\sqrt7}{2}$ & $0$\\
$\chi_{5}$ & $3$ & $0$ & $0$ & $\tfrac{-1-i\sqrt7}{2}$ & $\tfrac{-1+i\sqrt7}{2}$ & $0$\\
\hline
\end{tabular}
\end{table}

\section{Wreath products}
\begin{theorem}\label{thm:wreath}
Let $G=(\Ap{N},*)$ and $H=(\Ap{s},*_H)$ be groups and $h\mapsto S_h$ a homomorphism $H\to\Sym_\ell$. For $m\in\Ap{N^{\ell}s}$ put
$m^{\natural}=\dg{s}{0}{m/N^{\ell}}$ (control digit) and $\pi_m=S_{m^\natural}$. Then
\begin{equation}\label{eq:wreath}
m*n=\sum_{i=0}^{\ell-1}\Big(\dg{N}{i}{m}*\dg{N}{\pi_m^{-1}(i)}{n}\Big)N^{i}+N^{\ell}\big(m^{\natural}*_Hn^{\natural}\big)
\end{equation}
is a group, the \emph{wreath product} $G\wr H$ with \emph{fibre} $G$ and \emph{control} $H$, of order $N^{\ell}s$.
\end{theorem}
\begin{proof}
This is Theorem~\ref{thm:semidirect} with $G_0=G^\ell$, the direct product of
Eq.~\eqref{eq:dprod}. Define
\[
\phi_h(f)=\sum_i\dg{N}{S_h^{-1}(i)}{f}N^{i}.
\]
Each $\phi_h$ permutes coordinates, so it is an automorphism, and
$\phi_h\circ\phi_k=\phi_{h*_Hk}$ because $S_k^{-1}\circ S_h^{-1}=(S_h\circ S_k)^{-1}$.
\end{proof}

\begin{remark}
The transfer acts only on the right factor. If $\pi_m$ also permuted the digits of $m$, the law would fail to be associative.
The control element is recovered from $m*n$ through
\[
d_N(i,m*n)=d_N(i,m)*d_N(\pi_m^{-1}(i),n).
\]
Thus the control state determines exactly which coordinate of the right fibre enters each coordinate of the product. This is the
correct digital form of the transfer mechanism in the original notes.
\end{remark}

\section{Wreath products and \texorpdfstring{P\'olya}{Polya} counting}\label{sec:wreath:polya}
The wreath product is the group in which P\'olya's theory lives: $G\wr H$ is the symmetry group of $\ell$ blocks, each carrying a pattern that may be moved by $G$, where the blocks themselves may be permuted by $H$. This section proves that its cycle index is the composition of the cycle indices of $G$ and $H$ (Theorem~\ref{thm:wreath:ci}), derives the corresponding orbit counts (Corollary~\ref{cor:wreath:colour} and Theorem~\ref{thm:wreath:product}), and shows how the groups of de Bruijn's theorem sit inside a wreath product (Proposition~\ref{prop:wreath:debruijn}). The notation of Section~\ref{ssec:polya} is used throughout.

\subsection*{The action on blocks}
Let $G$ act on a finite set $X$ with $|X|=r$, let $H$ be a group with a homomorphism $h\mapsto S_h\in\Sym_\ell$, and write the elements of $G\wr H$ as pairs $(f,h)$ with $f=(f_0,\dots,f_{\ell-1})\in G^{\ell}$. In the digital form of Theorem~\ref{thm:wreath}, $G=(\Ap N,*)$, the pair $(f,h)$ is the number $m$ with digits $f_i=\dg{N}{i}{m}$ and control digit $m^{\natural}=h$, and~\eqref{eq:wreath} is the product
\begin{equation}\label{eq:wr:pair}
(f,h)(f',h')=\big(f\cdot(S_hf'),\,hh'\big),\qquad (S_hf')_i=f'_{S_h^{-1}(i)} .
\end{equation}
The \emph{imprimitive action} of $G\wr H$ on the $r\ell$ points $X\times\Ap{\ell}$, the set $X$ repeated in $\ell$ blocks, is
\begin{equation}\label{eq:wr:act}
(f,h)\cdot(x,i)=\big(f_{S_h(i)}\cdot x,\;S_h(i)\big).
\end{equation}
That this is an action follows from~\eqref{eq:wr:pair}: applying $(f',h')$ and then $(f,h)$ sends $(x,i)$ to $\big(f_{S_hS_{h'}(i)}\,f'_{S_{h'}(i)}\cdot x,\;S_hS_{h'}(i)\big)$, and the coordinate of $f\cdot(S_hf')$ at $S_hS_{h'}(i)$ is $f_{S_hS_{h'}(i)}\,f'_{S_{h'}(i)}$. When $G=(\Ap N,*)$ acts on $X=\Ap N$ by left multiplication, the action of the number $m$ is
\[
(x,i)\longmapsto\big(\dg{N}{\pi_m(i)}{m}*x,\;\pi_m(i)\big),\qquad \pi_m=S_{m^{\natural}} ,
\]
so the control digit moves the block and the digit of $m$ at the new block acts on the point inside it.

\begin{lemma}[Cycles of an element]\label{lem:wreath:cycles}
Let $(f,h)\in G\wr H$, and let $C=(i_0,i_1,\dots,i_{k-1})$ be a cycle of $S_h$ of length $k$, $i_t=S_h^{t}(i_0)$. Put $g_C=f_{i_0}\,f_{i_{k-1}}\cdots f_{i_1}\in G$, the \emph{cycle product}.
\begin{enumerate}
\item The cycles of $(f,h)$ on $X\times C$ correspond bijectively to the cycles of $g_C$ on $X$, and a cycle of length $j$ of $g_C$ corresponds to a cycle of length $jk$.
\item For fixed $h$, the map $f\mapsto(g_C)_C$ from $G^{\ell}$ to $G^{c(S_h)}$ is surjective and each of its fibres has $|G|^{\ell-c(S_h)}$ elements.
\end{enumerate}
\end{lemma}
\begin{proof}
(1) By~\eqref{eq:wr:act}, $(f,h)^{t}\cdot(x,i_0)=\big(f_{i_t}\cdots f_{i_1}x,\;i_t\big)$, and for $t=k$ this is $(g_Cx,i_0)$. Hence the cycle of $(x,i_0)$ meets the fibre over $i_0$ exactly in the $g_C$-cycle of $x$; if that has length $j$, the cycle has length $jk$. Every cycle on $X\times C$ meets this fibre.
(2) Distinct cycles of $S_h$ involve disjoint sets of coordinates, and every coordinate lies on one cycle. Within a cycle, each coordinate occurs once in $g_C$, so when all other coordinates are fixed the map $f_{i}\mapsto g_C$ is a bijection $G\to G$. Hence each $g_C$ is attained by $|G|^{k-1}$ choices of the coordinates on $C$, and the fibre over $(g_C)_C$ has $\prod_C|G|^{k_C-1}=|G|^{\ell-c(S_h)}$ elements.
\end{proof}

\subsection*{The cycle index}
\begin{theorem}[Cycle index of a wreath product]\label{thm:wreath:ci}
Let $Z(H)=Z(H;y_1,y_2,\dots)$ be the cycle index of $H$ acting on $\Ap{\ell}$ through $S$, and put
\[
z_k=Z(G;x_k,x_{2k},x_{3k},\dots)=\frac1{|G|}\sum_{g\in G}\ \prod_{j\ge1}x_{jk}^{c_j(g)} .
\]
Then the cycle index of $G\wr H$ in its action on $X\times\Ap{\ell}$ is the composition
\[
Z(G\wr H;x_1,x_2,\dots)=Z(H;z_1,z_2,z_3,\dots)=\frac1{|H|}\sum_{h\in H}\ \prod_{k\ge1}z_k^{\,c_k(S_h)} .
\]
\end{theorem}
\begin{proof}
By Lemma~\ref{lem:wreath:cycles}(1), the monomial of $(f,h)$ is $\prod_C\prod_jx_{jk_C}^{c_j(g_C)}$, where $C$ runs over the cycles of $S_h$ and $k_C$ is the length of $C$. Since $|G\wr H|=|G|^{\ell}|H|$, Lemma~\ref{lem:wreath:cycles}(2) gives
\[
Z(G\wr H)=\frac1{|G|^{\ell}|H|}\sum_{h}\ |G|^{\ell-c(S_h)}\prod_C\ \sum_{g\in G}\prod_jx_{jk_C}^{c_j(g)}
=\frac1{|H|}\sum_h\prod_C z_{k_C},
\]
and $\prod_Cz_{k_C}=\prod_kz_k^{c_k(S_h)}$.
\end{proof}

The substitution $x_j\mapsto x_{jk}$ in the $k$-th slot is the reason for the name \emph{composition}: a $k$-cycle of blocks lengthens every cycle inside a block by the factor $k$.

\begin{corollary}[Colourings of the blocks]\label{cor:wreath:colour}
Let $Y$ be a set of $m$ colours with weights $w$ and $p_k=\sum_yw(y)^{k}$. The weighted number of orbits of $G\wr H$ on the colourings of $X\times\Ap{\ell}$ is $Z(H;\omega_1,\omega_2,\dots)$, where $\omega_k=Z(G;p_k,p_{2k},p_{3k},\dots)$. In particular, with $N_K(t)=Z(K;t,t,\dots)$ the number of orbits of a group $K$ on $t$-colourings of its set, the number of orbits is $N_H\big(N_G(m)\big)$.
\end{corollary}
\begin{proof}
Theorem~\ref{thm:polya} says that the weighted count is the cycle index at $x_j=p_j$, which turns $z_k$ into $\omega_k$. For $w\equiv1$ all $p_j$ equal $m$, so every $z_k$ equals $N_G(m)$ and $Z(H;N_G(m),N_G(m),\dots)=N_H(N_G(m))$.
\end{proof}

The unweighted formula says that a colouring of the blocks is the same as an $\ell$-tuple of patterns, each pattern being a $G$-orbit of colourings of $X$, and $H$ then permutes the tuple. The next theorem proves this reading directly and for arbitrary $G$-sets.

\subsection*{Product action}
If $G$ acts on a set $\Omega$, then $G\wr H$ acts on $\Omega^{\ell}$ by the \emph{product action}
\begin{equation}\label{eq:wr:prod}
\big((f,h)\cdot u\big)_j=f_j\cdot u_{S_h^{-1}(j)} .
\end{equation}
For $\Omega=Y^{X}$ this is the action on colourings of $X\times\Ap{\ell}$: a colouring $c$ corresponds to $u_i=c(\,\cdot\,,i)$, and for $\gamma=(f,h)$ one has $\gamma^{-1}(y,j)=(f_j^{-1}y,S_h^{-1}(j))$ by~\eqref{eq:wr:act}, so $(\gamma c)(y,j)=u_{S_h^{-1}(j)}(f_j^{-1}y)=\big(f_j\cdot u_{S_h^{-1}(j)}\big)(y)$.

\begin{theorem}[P\'olya counting for the product action]\label{thm:wreath:product}
Let $\omega$ be the number of orbits of $G$ on $\Omega$. Then the number of orbits of $G\wr H$ on $\Omega^{\ell}$ is $N_H(\omega)=\frac1{|H|}\sum_{h\in H}\omega^{c(S_h)}$.
\end{theorem}
\begin{proof}
A tuple $u$ is fixed by $(f,h)$ if and only if $u_j=f_j\cdot u_{S_h^{-1}(j)}$ for all $j$. Along a cycle $C=(j_0,\dots,j_{k-1})$ of $S_h$ this gives $u_{j_t}=f_{j_t}u_{j_{t-1}}$, hence $u_{j_0}=g_C\cdot u_{j_0}$ with the cycle product $g_C$ of Lemma~\ref{lem:wreath:cycles}, and the value $u_{j_0}\in\operatorname{Fix}_\Omega(g_C)$ determines $u$ on $C$. Conversely every such value gives a fixed tuple. So $(f,h)$ fixes $\prod_C|\operatorname{Fix}_\Omega(g_C)|$ tuples. By Lemma~\ref{lem:wreath:cycles}(2) and Burnside's lemma (Theorem~\ref{thm:dtburnside}),
\[
\frac1{|G|^{\ell}}\sum_{f\in G^{\ell}}\ \prod_C|\operatorname{Fix}_\Omega(g_C)|=\prod_C\ \frac1{|G|}\sum_{g\in G}|\operatorname{Fix}_\Omega(g)|=\omega^{c(S_h)} .
\]
Burnside's lemma for $G\wr H$ now gives $\frac1{|H|}\sum_h\omega^{c(S_h)}$.
\end{proof}

For $\Omega=Y^{X}$ one has $\omega=N_G(m)$, which is the second proof of Corollary~\ref{cor:wreath:colour}. The two proofs are two readings of one computation: the cycle product $g_C$ acts on a block, and the fixed colourings of that block are counted either by the cycles of $g_C$ or by its fixed points.

\subsection*{De Bruijn's group inside a wreath product}
\begin{proposition}\label{prop:wreath:debruijn}
Let $G$ act on $X$ and $H$ on $Y$, with $\ell=|X|$. The map
\[
(g,h)\mapsto\big((h,h,\dots,h),g\big)
\]
is an injective homomorphism $G\times H\to H\wr G$, where $G$ permutes the $\ell$ blocks through its action on $X$, and it carries the action of Theorem~\ref{thm:debruijn} on $Y^{X}$ to the restriction of the product action of $H\wr G$ on $Y^{X}=Y^{\ell}$.
\end{proposition}
\begin{proof}
A constant tuple is unchanged by any permutation of its coordinates, so by~\eqref{eq:wr:pair} $\big((h,\dots,h),g\big)\big((h',\dots,h'),g'\big)=\big((hh',\dots,hh'),gg'\big)$. In the product action~\eqref{eq:wr:prod} the element acts by $u\mapsto\big(h\cdot u_{g^{-1}(x)}\big)_x$, which is $((g,h)\cdot u)(x)=h\big(u(g^{-1}x)\big)$.
\end{proof}

Thus de Bruijn's group is the subgroup of $H\wr G$ whose fibre element is the same at every position. The full wreath product allows an independent element of $H$ at each position and therefore has fewer orbits: by Theorem~\ref{thm:wreath:product} it has $N_G(\omega_H)$, where $\omega_H$ is the number of $H$-orbits on $Y$. For $G=H=C_2$ acting on $X=Y=\Ap2$ the subgroup has $2$ orbits on the four functions $X\to Y$ (Exercise~\ref{ex:wr:debruijn}) and the wreath product has one.

\subsection*{Examples}
\begin{example}[$C_2\wr C_2=D_8$]
With $Z(C_2)=\frac12(x_1^{2}+x_2)$ on two points one has $z_1=\frac12(x_1^{2}+x_2)$ and $z_2=\frac12(x_2^{2}+x_4)$, so
\[
Z(C_2\wr C_2)=\tfrac12\big(z_1^{2}+z_2\big)=\tfrac18\big(x_1^{4}+2x_1^{2}x_2+3x_2^{2}+2x_4\big),
\]
the cycle index of the square in the example following Theorem~\ref{thm:debruijn}: the blocks are the two diagonals $\{0,2\}$ and $\{1,3\}$, the fibre group swaps the two vertices of a diagonal, and the control swaps the diagonals. The number of orbits on $m$-colourings is $N_{C_2}\big(N_{C_2}(m)\big)=q(q+1)/2$ with $q=m(m+1)/2$, which is $6$ for $m=2$ and $21$ for $m=3$.
\end{example}

\begin{example}[The hyperoctahedral group]
The group $C_2\wr\Sym_n$ of order $2^{n}n!$ acts on $2n$ points, $n$ blocks of two. By Corollary~\ref{cor:wreath:colour} and Corollary~\ref{cor:sn:cycleindex}, the number of orbits on $m$-colourings is $N_{\Sym_n}(q)=\binom{q+n-1}{n}$ with $q=N_{C_2}(m)=m(m+1)/2$: the orbits are the multisets of $n$ unordered pairs of colours. For $n=3$ and $m=2$ there are $\binom53=10$. More generally, for every $G$ the generating function of the wreath products with symmetric groups is
\[
\sum_{n\ge0}Z(G\wr\Sym_n)\,t^{n}=\exp\Big(\sum_{k\ge1}\frac{t^{k}}{k}\,Z(G;x_k,x_{2k},\dots)\Big),
\]
by Theorem~\ref{thm:wreath:ci} and the exponential formula of Section~\ref{sec:sn:classes}.
\end{example}

\begin{example}[Sylow subgroups of symmetric groups]
For a prime $p$, $C_p\wr C_p$ is a Sylow $p$-subgroup of $\Sym_{p^{2}}$. With $Z(C_p;x)=\frac1p\big(x_1^{p}+(p-1)x_p\big)$ one gets $z_1=\frac1p\big(x_1^{p}+(p-1)x_p\big)$ and $z_p=\frac1p\big(x_p^{p}+(p-1)x_{p^{2}}\big)$, and the number of orbits on $m$-colourings of the $p^{2}$ points is $\frac1p\big(q^{p}+(p-1)q\big)$ with $q=\frac1p\big(m^{p}+(p-1)m\big)$. For $p=3$ and $m=2$ this is $24$ orbits, and $451$ for $m=3$. Iterating, the number of orbits of the $n$-fold wreath power $C_p\wr\cdots\wr C_p$ on $p^{n}$ points is obtained by applying $q\mapsto\frac1p\big(q^{p}+(p-1)q\big)$ to $m$ a total of $n$ times.
\end{example}

\begin{table}[ht]
\centering
\small
\begin{tabular}{lrrrr}
Group & order & points & $m=2$ & $m=3$\\\hline
$C_2\wr C_2$ & $8$ & $4$ & $6$ & $21$\\
$C_3\wr C_2$ & $18$ & $6$ & $10$ & $66$\\
$C_2\wr\Sym_3$ & $48$ & $6$ & $10$ & $56$\\
$\Sym_3\wr C_2$ & $72$ & $6$ & $10$ & $55$\\
$C_3\wr C_3$ & $81$ & $9$ & $24$ & $451$
\end{tabular}
\caption{Orbits of wreath products on colourings. Each group was generated explicitly in \texttt{scripts/polya\_verify.py}; the order, the identity of Theorem~\ref{thm:wreath:ci} at random points, and the counts $N_H(N_G(m))$ were checked.}\label{tab:wreath:count}
\end{table}

\section{Transitive actions and stabiliser chains}\label{sec:transitive}
A transitive action labels each group element by a point and an element of the point stabiliser. Iterating this along a
chain of point stabilisers gives mixed-radix labels whose radices are orbit sizes. These laws are used for the symmetric and
alternating groups in Chapter~\ref{ch:simple} and for the sporadic groups in Chapter~\ref{ch:sporadic}.

\begin{theorem}[Transitive law]\label{thm:transitive}
Let $G$ act transitively on $\Ap{n}$, and let $K$ be the stabiliser of $0$, coded by a group law $(\Ap{k},*_K)$ with identity $0$. Choose
$s_t\in G$ with $s_t(0)=t$ and $s_0=1$. Code $g=s_t\,h$, with $h\in K$ of code $\kappa$, by $m=t+n\kappa$. Then the law
\begin{equation}\label{eq:transitive}
(t,h)*(u,k)=\big(t',\;s_{t'}^{-1}s_t\,h\,s_u\,k\big),\qquad t'=s_t\big(h(u)\big),
\end{equation}
makes $\Ap{nk}$ a group isomorphic to $G$, with identity $0$.
\end{theorem}
\begin{proof}
Every $g\in G$ can be written uniquely as $g=s_{g(0)}\cdot\big(s_{g(0)}^{-1}g\big)$, with the second factor in $K$. Hence
$g\mapsto\big(g(0),s_{g(0)}^{-1}g\big)$ is a bijection $G\to\Ap n\times K$, and $|G|=nk$. For $g=s_th$ and $g'=s_uk$ we have
$gg'(0)=s_th(s_u(0))=s_t(h(u))=t'$, so $gg'=s_{t'}\cdot\big(s_{t'}^{-1}s_ths_uk\big)$. Thus~\eqref{eq:transitive} is the multiplication
of $G$ transported along this bijection. Since $s_0=1$ and $0$ codes $1\in K$, the code $0$ is the identity.
\end{proof}

\begin{remark}
If $K$ is normal in $G$ and the representatives $s_t$ form a subgroup $T$, then $h(u)=u$ for $h\in K$, so $s_{t'}=s_ts_u$ and
\eqref{eq:transitive} reduces to $(t,h)*(u,k)=\big(t',\;s_u^{-1}hs_u\,k\big)$, the multiplication of the semidirect product
$G=K\rtimes T$ of Theorem~\ref{thm:semidirect}.
\end{remark}

\begin{theorem}[Chain law]\label{thm:chainlaw}
Let $G$ act faithfully on $\Ap n$, and let $(b_1,\dots,b_k)$ be a base: only the identity fixes every $b_l$. Put $G_0=G$, let $G_l$
be the pointwise stabiliser of $b_1,\dots,b_l$, let $\Delta_l=G_{l-1}(b_l)$ and $r_l=|\Delta_l|$. List $\Delta_l$ as
$\delta_{l,0}=b_l,\delta_{l,1},\dots,\delta_{l,r_l-1}$, and choose $s_{l,t}\in G_{l-1}$ with $s_{l,t}(b_l)=\delta_{l,t}$ and
$s_{l,0}=1$. Then every $g\in G$ has a unique expression
\[
g=s_{1,t_1}\circ s_{2,t_2}\circ\cdots\circ s_{k,t_k},\qquad t_l\in\Ap{r_l},
\]
so $|G|=r_1\cdots r_k$. The digits are read off from the top: with $g_1=g$, $t_l$ is the position of $g_l(b_l)$ in $\Delta_l$ and
$g_{l+1}=s_{l,t_l}^{-1}\circ g_l$. The code
\[
\lambda(g)=t_1+r_1\big(t_2+r_2\big(t_3+\dots+r_{k-1}t_k\big)\big)
\]
is a bijection $G\to\Ap{|G|}$ with $\lambda(1)=0$, and $m*n=\lambda\big(\lambda^{-1}(m)\circ\lambda^{-1}(n)\big)$ is a digital law
of $G$.
\end{theorem}

\begin{proof}
$G_{l-1}$ acts transitively on $\Delta_l$, and the stabiliser of $b_l$ in $G_{l-1}$ is $G_l$. Relabel $\Delta_l$ by positions,
so that $b_l$ becomes $0$, and apply Theorem~\ref{thm:transitive} to $G_{l-1}$ with $K=G_l$. This writes each element of
$G_{l-1}$ uniquely as $s_{l,t}h$ with $h\in G_l$, coded by $t+r_l\kappa$. Induction on $l$, ending with $G_k=1$ because the
$b_l$ form a base, gives the unique expression and the mixed-radix code. Since $s_{l,0}=1$ for all $l$, the identity has code $0$.
\end{proof}

\section{Simple groups and the classification}\label{sec:simpleclass}
The groups of this chapter are assembled from smaller ones: a semidirect product adds a layer of digits on top of a normal
subgroup, and a wreath product stacks a whole alphabet of copies. A group that cannot be assembled in this way is
\emph{simple}. In digital terms, Corollary~\ref{cor:dtsimple} says that $G$ is simple if and only if its composition
labellings have a single layer, and that $G$ is solvable if and only if the radices of a composition labelling are the
prime factors of $N=|G|$. By Theorem~\ref{thm:dtjh} the radices and layer codes of a composition labelling do not depend on
the labelling, so the simple groups are the building blocks of all finite groups. This section uses
Chapters~\ref{ch:digitaltheorems} and~\ref{ch:characters}, and it anticipates the later chapters, where a digital law is
given for every finite simple group.

\subsection{The classification}
\begin{theorem}[Classification of the finite simple groups]\label{thm:cfsg}
Every finite simple group is isomorphic to one of the following:
\begin{enumerate}
\item a cyclic group $C_p$ of prime order (Chapter~\ref{ch:abelian});
\item an alternating group $\Alt_n$, $n\ge5$ (Chapter~\ref{ch:simple});
\item a group of Lie type, the Tits group ${}^2F_4(2)'$ being counted among them (Chapter~\ref{ch:lie});
\item one of the $26$ sporadic groups (Chapter~\ref{ch:sporadic}).
\end{enumerate}
\end{theorem}
The proof is spread over several hundred articles and is not reproduced here. What the following chapters add is a
digital law $(\Ap{N},*)$ for every entry of the list. The smallest non-Abelian simple group is $\Alt_5$, of order $60$,
and every non-Abelian group in the list has even order. That this is no accident is the content of the next theorem.

\subsection{Groups of odd order}
\begin{theorem}[Feit--Thompson]\label{thm:oddorder}
Every group of odd order is solvable.
\end{theorem}
By Corollary~\ref{cor:dtsimple}, Theorem~\ref{thm:oddorder} says that every law $(\Ap{N},*)$ with $N$ odd can be
relabelled into a mixed radix $(p_1,\dots,p_r)$ of primes, $N=p_1\cdots p_r$, in which every truncation is a homomorphism and
every layer code is a relabelling of $\dg{p_i}{0}{a+b}$. The original proof occupies a whole issue of the Pacific Journal of
Mathematics~\cite{FeitThompson}, and it has been checked in full by computer in the Coq proof assistant~\cite{Gonthier2013}. Its proof is not reproduced here either;
instead we record what the digital language sees of the theorem, and where it stops.

Write $N=2q+1$. The number $2q$ enters first through inversion.

\begin{proposition}[Inversion as negation]\label{prop:oddinversion}
If $N$ is odd, $G$ has a labelling $\lambda\colon G\to\Ap{N}$ with $\lambda(1)=0$ and $\lambda(g^{-1})=N-\lambda(g)$ for $g\ne1$.
In such a labelling the inverse of $m\ne0$ is $N-m$, as in $C_N$, and $m\mapsto m*m$ is a bijection of $\Ap{N}$.
\end{proposition}
\begin{proof}
A group of odd order has no element of order $2$, so $g\ne g^{-1}$ for $g\ne1$, and inversion splits $G\setminus\{1\}$ into $q$
pairs $\{g,g^{-1}\}$. Label the $j$-th pair by $j$ and $N-j$, $1\le j\le q$. Squaring is injective, because $g^2=h^2$ gives
$g=(g^2)^{q+1}=(h^2)^{q+1}=h$, using $g^{N}=1$.
\end{proof}

\begin{lemma}[Burnside]\label{lem:oddreal}
If $N$ is odd, the trivial character is the only real-valued irreducible character of $G$.
\end{lemma}
\begin{proof}
The number of real-valued irreducible characters equals the number of conjugacy classes $C$ with $C=C^{-1}$, since complex
conjugation permutes the rows and the columns of the character table with the same number of fixed points
\cite[Problem~3.16]{Isaacs}. Suppose $xgx^{-1}=g^{-1}$. Then $x^2$ commutes with $g$, and since $x$ has odd order, $x$ is a power
of $x^2$; hence $x$ commutes with $g$ and $g=g^{-1}$, so $g=1$. The only real class is $\{1\}$.
\end{proof}

\begin{proposition}[The decomposition of $q$]\label{prop:qdecomp}
Let $N=2q+1$ be odd. The non-trivial irreducible characters of $G$ fall into $s$ pairs $\{\chi_i,\overline{\chi_i}\}$, and with
$d_i=\chi_i(0)$,
\[
q=\sum_{i=1}^{s}d_i^2 ,
\]
where each $d_i$ is an odd divisor of $N$. The number of conjugacy classes is $k=2s+1$, and $k\equiv N \pmod{16}$.
\end{proposition}
\begin{proof}
By Lemma~\ref{lem:oddreal}, $\chi\ne\overline{\chi}$ for $\chi\ne1$, which gives the pairs. Theorem~\ref{thm:charsquares} gives
$N=1+2\sum_i d_i^2$, and $d_i$ divides $N$ by Corollary~\ref{cor:degreedivides}, so it is odd. The number of classes equals
the number of irreducible characters, $1+2s$. Odd squares are $\equiv1\pmod 8$, so $q\equiv s\pmod 8$ and
$N=2q+1\equiv 2s+1=k\pmod{16}$.
\end{proof}

\begin{example}[Order $21$]
For $N=21$, $q=10$. The cyclic group $C_{21}$ gives $q=1^2+\dots+1^2$ ($s=10$, $k=21$). The metacyclic group
$C_7\rtimes_2C_3$ of Section~\ref{sec:metacyclic-chars} has degrees $1,1,1,3,3$, so $q=1^2+3^2$, $s=2$ and $k=5$. In both
cases $k\equiv21\pmod{16}$.
\end{example}

The commutator subgroup $G'$ is the natural normal subgroup to construct. In the labelling of
Proposition~\ref{prop:oddinversion} it is generated by the labels $(N-m)*(N-n)*m*n$, with $N-0$ read as $0$, so it can be
computed from the law alone, and $|G:G'|$ is the number of linear characters (Chapter~\ref{ch:characters}).

\begin{theorem}[Equivalent forms of the Odd Order Theorem]\label{thm:oddforms}
The following statements are equivalent.
\begin{enumerate}
\item Every group of odd order is solvable.
\item Every non-trivial group of odd order satisfies $G'\ne G$.
\item For every non-trivial group of odd order $N=2q+1$, some $d_i$ in Proposition~\ref{prop:qdecomp} equals $1$.
\end{enumerate}
\end{theorem}
\begin{proof}
(1)$\Rightarrow$(2) because a solvable non-trivial group has a non-trivial Abelian quotient. (2)$\Rightarrow$(1): every
subgroup of a group of odd order has odd order, so the derived series $G\ge G'\ge G''\ge\cdots$ decreases strictly until it
reaches $1$. (2)$\Leftrightarrow$(3): the non-trivial linear characters are non-real by Lemma~\ref{lem:oddreal}, so
$|G:G'|=1+2\,\#\{i: d_i=1\}$, which exceeds $1$ exactly when some $d_i=1$.
\end{proof}

\subsection{What the order alone cannot see}
It is tempting to look for the normal subgroup in the digits of $N=2q+1$ in radix $2q$. This cannot work.

\begin{proposition}[Casting out $2q$]\label{prop:castout}
For every $x\ge0$, $x\equiv\sum_k\dg{N}{k}{x}\pmod{N-1}$. If $R$ is the code of any law on $\Ap{N}$ whose Cayley table is a
Latin square (in particular of any group), and $N=2q+1$, then $\dg{2q}{0}{R}=q$, whatever the order in which the entries are
placed in $R$.
\end{proposition}
\begin{proof}
Since $N\equiv1\pmod{N-1}$, $x=\sum_k\dg{N}{k}{x}N^k\equiv\sum_k\dg{N}{k}{x}$. The digits of $R$ are the $N^2$ entries of the table;
each row is a permutation of $\Ap{N}$ with sum $N(N-1)/2=Nq$, so the digit sum is $N^2q=(4q^2+4q+1)q\equiv q\pmod{2q}$.
\end{proof}
The only digit of $R$ in radix $2q$ that is determined by digits in radix $N$ is therefore the same for all groups of order
$N$. More simply, $\gcd(N,2q)=1$, so no non-trivial subgroup has an order dividing $2q$: normal subgroups correspond to
factorisations $N=ab$ and truncations in radix $(a,b)$, not to additive decompositions of $N$.

The arithmetic of Proposition~\ref{prop:qdecomp} does not decide the question either. A perfect group of odd order would need
a decomposition $q=\sum d_i^2$ with every $d_i$ an odd divisor of $N$, $3\le d_i\le\sqrt N$, and such decompositions exist:
\begin{center}
\begin{tabular}{rrl}
\hline
$N$ & $q$ & degrees $d_i$ of the pairs\\
\hline
$105$ & $52$ & $3,3,3,5$\\
$189$ & $94$ & $3,3,3,3,3,7$\\
$255$ & $127$ & $3,3,3,5,5,5,5$\\
$297$ & $148$ & $3,3,3,11$\\
\hline
\end{tabular}
\end{center}
For $N=105$ the decomposition gives $k=9\equiv105\pmod{16}$, so even Burnside's congruence is satisfied. Yet no group of order
$105=3\cdot5\cdot7$ is perfect: if neither the Sylow $5$-subgroup nor the Sylow $7$-subgroup were normal there would be $21$ of
the first and $15$ of the second (Section~\ref{sec:dtsylow}), hence $84+90>105$ elements of order $5$ or $7$. Excluding the
decompositions of $q$ requires information about subgroups, not only about $N$ and the degrees; this is the local analysis of
the Feit--Thompson proof.

\subsection{Cases within reach, and an arithmetic question}
Two classical special cases have short proofs in the present language. If $N$ is square-free, all Sylow subgroups are
cyclic, $G$ is metacyclic by a theorem of H\"older and Burnside \cite{Hall}, and its law is the two-digit rule of
Section~\ref{sec:metacyclic-chars}; it is solvable by inspection. If $N=p^aq^b$, Burnside's vanishing criterion of
Chapter~\ref{ch:characters} gives solvability for all $N$ of this form, odd or even.

The last step of the Feit--Thompson proof is arithmetic, and it is digital in nature. For distinct primes $p$ and $q$, the
number $(p^q-1)/(p-1)=\sum_{k=0}^{q-1}p^k$ is the \emph{repunit} whose $q$ digits in base $p$ are all equal to $1$:
$\dg{p}{k}{(p^q-1)/(p-1)}=1$ for $0\le k<q$. Feit and Thompson~\cite{FeitThompson1962} conjectured that this repunit never divides the repunit with $p$
digits in base $q$; a proof would shorten the final chapter of their argument considerably. The conjecture is known for
$q=2$ and $q=3$. The stronger statement that the two repunits are coprime is false: for $p=17$ and $q=3313$ they share the
factor $2pq+1=112643$~\cite{Stephens1971}.

\section{Exercises}
\begin{enumerate}
\item Compute the Cayley table of $D_6$ from~\eqref{eq:dihedral} and give an explicit isomorphism onto $\Sym_3=\{5,7,11,15,19,21\}$.
\item Show that $C_q\rtimes_rC_s\cong C_q\rtimes_{r'}C_s$ whenever $r'\equiv r^k\pmod q$ with $\gcd(k,s)=1$.
\item Show that the centre of $D_{2q}$ is trivial for odd $q$ and is $\{0,q/2\}$ for even $q$.
\item Show that $Q_8$ has a unique element of order $2$ and that all its subgroups are normal; deduce $Q_8\not\cong D_8$.
\item Let $q$ be prime and let $r$ have order $s>1$ modulo $q$. Show that $C_q\rtimes_rC_s$ has trivial centre.
\item Using Theorem~\ref{thm:metacyclic-chars}, show that $C_{13}\rtimes_3C_3$ has $7$ conjugacy classes and character degrees $1,1,1,3,3,3,3$.
\item Show that the group of Theorem~\ref{thm:wreath} has order $N^{\ell}s$, and identify $C_2\wr C_2$ (with $C_2$ permuting two
coordinates) with $D_8$.
\item Compute the cycle index of $C_3\wr C_2$ acting on six points and show that it has $10$ orbits on $2$-colourings.
\item Show that the number of orbits of $C_2\wr\Sym_n$ on the $m$-colourings of its $2n$ points is $\binom{q+n-1}{n}$ with $q=m(m+1)/2$.
\item For a prime $p$, show that $C_p\wr C_p$ has $\frac1p\big(q^{p}+(p-1)q\big)$ orbits on the $m$-colourings of $p^{2}$ points, where $q=\frac1p\big(m^{p}+(p-1)m\big)$, and evaluate it for $p=3$, $m=2$.
\item\label{ex:wr:debruijn} Let $C_2\times C_2$ act on the four functions $\Ap2\to\Ap2$ as in Theorem~\ref{thm:debruijn}, the first factor swapping the two positions and the second complementing the values. Count the orbits with Theorem~\ref{thm:debruijn}, and compare with the orbits of $C_2\wr C_2$ on the same set (Proposition~\ref{prop:wreath:debruijn}).
\item Compute $\dg{20}{0}{R}$ for the codes of $C_{21}$ and $C_7\rtimes_2C_3$, and check Proposition~\ref{prop:castout}.
\item List all odd $N<400$ for which $q=(N-1)/2$ is a sum of squares of odd divisors $3\le d\le\sqrt N$ of $N$.
\item Construct the labelling of Proposition~\ref{prop:oddinversion} for $C_7\rtimes_2C_3$ and compute $G'$ from the law.
\item Show that every group of order $105$ has a normal subgroup of order $35$.
\end{enumerate}

\chapter{Small finite groups}\label{ch:small}

\begin{theorem}[Lagrange]
If $H$ is a subgroup of a group $(\Ap{N},R)$, then $|H|$ divides $N$.
\end{theorem}
\begin{proof}
The left cosets $\lambda_a(H)$ are images of $H$ under the bijections $\lambda_a$ (Theorem~\ref{thm:cayley}), so all have $|H|$
elements, and two cosets are equal or disjoint.
\end{proof}

\begin{corollary}\label{cor:primeorder}
If $p$ is prime, every group of order $p$ is isomorphic to $C_p$.
\end{corollary}
\begin{proof}
Any $a\neq0$ generates a subgroup of order greater than $1$ dividing $p$, hence the whole group, and $k\mapsto a^{k}$ is an isomorphism from $C_p$.
\end{proof}

\begin{example}[Order $4$]
There are two groups: $C_4$, with $a*b=\dg{4}{0}{a+b}$, and $V_4=C_2\oplus C_2$, with
$a*b=\dg{2}{0}{a+b}+2\,\dg{2}{0}{\dg{2}{1}{a}+\dg{2}{1}{b}}$ (digitwise addition mod $2$).
\end{example}

\begin{example}[Order $6$]
There are two groups: $C_6\cong C_3\oplus C_2$ (Theorem~\ref{thm:crt}), and the nonabelian $D_6\cong\Sym_3$
(Corollary~\ref{cor:dihedral} with $q=3$).
\end{example}

\section{Classification of the groups of order at most \texorpdfstring{$15$}{15}}\label{sec:cat15}
This section lists every isomorphism type of group of order $n\le15$, each realised on the alphabet $\Ap{n}$ with identity $0$
by one of the digital laws already constructed: the cyclic law, the digitwise mixed-radix law of direct sums, the metacyclic law
\eqref{eq:metacyclic}, the dihedral law~\eqref{eq:dihedral}, the dicyclic law~\eqref{eq:dicyclic}, and the semidirect law
\eqref{eq:semidirect}. The only group needing a new automorphism is $\Alt_4$.

\subsection{The alternating group \texorpdfstring{$\Alt_4\cong V_4\rtimes C_3$}{Alt4 = V4:C3}}
\begin{proposition}\label{prop:A4digital}
Write $m,n\in\Ap{12}$ in the mixed radix $(2,2,3)$,
\[
\begin{aligned}
a_0&=\dg{2}{0}{m},& a_1&=\dg{2}{1}{m},& c&=\dg{3}{0}{\tfrac m4},\\
b_0&=\dg{2}{0}{n},& b_1&=\dg{2}{1}{n},& e&=\dg{3}{0}{\tfrac n4},
\end{aligned}
\]
and let $f_k=\dg{3}{0}{k^2}$, so that $f_k=0$ if $3\mid k$ and $f_k=1$ otherwise. Then the law
\begin{equation}\label{eq:A4}
m*n=\dg{2}{0}{a_0+f_{c+2}\,b_0+f_c\,b_1}+2\,\dg{2}{0}{a_1+f_c\,b_0+f_{c+1}\,b_1}+4\,\dg{3}{0}{c+e}
\end{equation}
makes $\Ap{12}$ a group isomorphic to $\Alt_4\cong V_4\rtimes C_3$. The two lower layers carry the Klein group
$V_4=C_2\oplus C_2$ (digitwise addition modulo $2$), the top layer carries $C_3$, and the coefficients $f$ are the
action of $C_3$ on $V_4$: Eq.~\eqref{eq:A4} is $m*n=\big(\dg{4}{0}{m}\oplus\phi^{c}(\dg{4}{0}{n})\big)+4\,\dg{3}{0}{c+e}$, where
$\oplus$ is the bitwise exclusive-or on $\Ap4$ and $\phi=\left(\begin{smallmatrix}0&1&2&3\\0&2&3&1\end{smallmatrix}\right)$ is the
linear map $(v_0,v_1)\mapsto(v_1,v_0+v_1)$ of $\mathbb{F}_2^2$ written on the code $v_0+2v_1$.
\end{proposition}
\begin{proof}
On the coordinates $(v_0,v_1)$ the map $\phi$ has matrix $M=\left(\begin{smallmatrix}0&1\\1&1\end{smallmatrix}\right)$ over
$\mathbb{F}_2$. Since $M^2=M+I$, its powers are $M^c=\left(\begin{smallmatrix}F_{c-1}&F_c\\F_c&F_{c+1}\end{smallmatrix}\right)$
with $F_k$ the Fibonacci numbers, and modulo $2$ these have period $3$ with values $0,1,1$, that is, $F_k\equiv f_k$
and $F_{c-1}\equiv f_{c+2}$. Hence the two lower digits of~\eqref{eq:A4} are $\dg{4}{0}{m}\oplus\phi^{c}(\dg{4}{0}{n})$
and the two forms agree. As $M^3=I$ and $M$ is invertible, $c\mapsto\phi^c$ is a homomorphism $C_3\to\Aut(V_4)$, and
Theorem~\ref{thm:semidirect} gives a group of order $12$. Its normal subgroup $\Ap4$ (the elements with $c=0$) is a Klein
group. The eight elements with $c\neq0$ have order $3$, because $M^2+M+I=0$ over $\mathbb{F}_2$. A group of order $12$ with
eight elements of order $3$ is $\Alt_4$ by Theorem~\ref{thm:order12} below.
\end{proof}

\begin{remark}[The same law as a code]
For a group this small the action can also be stored in a single code, by the coding principle. Writing the three
permutations $\phi^0,\phi^1,\phi^2$ one after the other in base $4$ (digits $0123$, $0231$, $0312$ from position $0$),
\[
\Phi=\sum_{c=0}^{2}\sum_{v=0}^{3}\phi^c(v)\,4^{4c+v}=10254564,\qquad \phi^c(v)=\dg{4}{4c+v}{\Phi}.
\]
With the exclusive-or written digitally, $u\oplus w=\sum_{j=0}^{1}2^j\,\dg{2}{0}{\dg{2}{j}{u}+\dg{2}{j}{w}}$, the law becomes
\[
m*n=\sum_{j=0}^{1}2^j\,\dg{2}{0}{\dg{2}{j}{m}+\dg{2}{j}{w}}+4\,\dg{3}{0}{c+e},\qquad w=\dg{4}{4c+\dg{4}{0}{n}}{\Phi}.
\]
This form looks the action up in a table instead of computing it; Eq.~\eqref{eq:A4} shows its arithmetic, which is what
generalises.
\end{remark}

\subsection{Classification theorems}
\begin{theorem}[Order $p^2$]\label{thm:orderp2}
If $|G|=p^2$, then $G\cong C_{p^2}$ or $G\cong C_p\oplus C_p$.
\end{theorem}
\begin{proof}
The class equation gives $Z(G)\neq\{0\}$. If $|Z(G)|=p$, then $G/Z(G)$ is cyclic, which forces $G$ Abelian, a contradiction.
Hence $G$ is Abelian, and the classification of finite Abelian groups gives the two types.
\end{proof}

\begin{theorem}[Order $pq$]\label{thm:orderpq}
Let $p<q$ be primes. If $p\nmid q-1$, then $G\cong C_{pq}$. If $p\mid q-1$, then $G$ is either $C_{pq}$ or the nonabelian group
$C_q\rtimes_rC_p$ of Theorem~\ref{thm:metacyclic}, where $r$ has multiplicative order $p$ modulo $q$; this group does not depend on the
choice of $r$.
\end{theorem}
\begin{proof}
By Sylow's theorems $n_q\equiv1\pmod q$ and $n_q\mid p<q$, so $n_q=1$, the Sylow subgroup $C_q$ is normal, and
$G\cong C_q\rtimes_\phi C_p$ with $\phi:C_p\to\Aut(C_q)\cong C_{q-1}$. If $p\nmid q-1$, $\phi$ is trivial. Otherwise the nontrivial
homomorphisms send a generator to the elements $r,r^2,\dots,r^{p-1}$ of order $p$; replacing the generator of $C_p$ by its $k$th power
turns $r$ into $r^k$, so all of them give isomorphic groups.
\end{proof}
\noindent For $p=2$ the nonabelian case is the dihedral group $D_{2q}$ of Corollary~\ref{cor:dihedral}. This covers $n=6,10,14$.
Since $3\nmid4$, every group of order $15$ is cyclic.

\begin{theorem}[Order $8$]\label{thm:order8}
There are exactly five groups of order $8$: $C_8$, $C_4\oplus C_2$, $C_2\oplus C_2\oplus C_2$, $D_8=C_4\rtimes_3C_2$ and $Q_8$.
\end{theorem}
\begin{proof}
The Abelian types come from the partitions of $3$. Let $G$ be nonabelian. Then no element has order $8$, and not every element has
order $\le2$ (otherwise $G$ is Abelian). Take $a$ of order $4$ and $b\notin\langle a\rangle$. The subgroup $\langle a\rangle$ has index $2$,
so it is normal, and conjugation by $b$ is an automorphism of $C_4$. It is not the identity, since then $G$ would be Abelian, so
$bab^{-1}=a^{-1}$. Next $b^2\in\langle a\rangle$, and $b^2\neq a^{\pm1}$ because $b$ would then have order $8$, so $b^2\in\{0,a^2\}$.
If some $b\notin\langle a\rangle$ has $b^2=0$, the relations are those of $D_8$ (Proposition~\ref{prop:presentations}). Otherwise
$b^2=a^2$, which gives $Q_8$. The two groups are not isomorphic: $D_8$ has five involutions and $Q_8$ has one.
\end{proof}

\begin{theorem}[Order $12$]\label{thm:order12}
There are exactly five groups of order $12$: $C_{12}$, $C_6\oplus C_2$, $D_{12}=C_6\rtimes_5C_2$, $\Alt_4$ (Proposition~\ref{prop:A4digital})
and $Q_{12}=C_3\rtimes_2C_4$.
\end{theorem}
\begin{proof}
By Sylow's theorems $n_3\in\{1,4\}$. If $n_3=4$, the four Sylow $3$-subgroups contain eight elements of order $3$. The remaining four
elements form the unique, hence normal, Sylow $2$-subgroup $S$, so $G\cong S\rtimes_\phi C_3$. Here $\phi$ is nontrivial, since otherwise
the Sylow $3$-subgroup would be normal. Because $\Aut(C_4)\cong C_2$ has no element of order $3$, $S=V_4$. The two nontrivial choices
$\phi^{\pm1}$ differ by inverting the generator of $C_3$, so $G\cong\Alt_4$.

If $n_3=1$, then $G\cong C_3\rtimes_\psi S$ with $S\in\{C_4,V_4\}$ and $\psi:S\to\Aut(C_3)\cong C_2$. For $S=C_4$, the trivial $\psi$ gives
$C_{12}$ and the surjective one gives $C_3\rtimes_2C_4=Q_{12}$. For $S=V_4$, the trivial $\psi$ gives $C_6\oplus C_2$. The three surjective
$\psi$ are permuted transitively by $\Aut(V_4)\cong\Sym_3$, so they give one group, $D_6\times C_2\cong D_{12}$.
The element-order profiles of Table~\ref{tab:profiles15} show that the five groups are pairwise nonisomorphic.
\end{proof}

\begin{remark}
The group $C_3\rtimes_2C_4$ is the dicyclic group $Q_{4\cdot3}$ of Theorem~\ref{thm:dicyclic}: both are nonabelian, of order $12$, with a
normal $C_3$ and a cyclic Sylow $2$-subgroup. The dicyclic law~\eqref{eq:dicyclic} with $q=3$ and the metacyclic law
\eqref{eq:metacyclic} with $(q,s,r)=(3,4,2)$ are two digital codings of the same group.
\end{remark}

\subsection{Catalogue}
\begin{table}[ht]
\centering
\begin{tabular}{c|c|l}
$n$ & number & groups and digital law\\ \hline
$1$ & $1$ & $C_1$\\
$2,3,5,7,11,13$ & $1$ & $C_n$ (Corollary~\ref{cor:primeorder})\\
$4$ & $2$ & $C_4,\ C_2\oplus C_2$\\
$6$ & $2$ & $C_6,\ D_6=C_3\rtimes_2C_2\cong\Sym_3$\\
$8$ & $5$ & $C_8,\ C_4\oplus C_2,\ C_2^{3},\ D_8=C_4\rtimes_3C_2,\ Q_8$\\
$9$ & $2$ & $C_9,\ C_3\oplus C_3$\\
$10$ & $2$ & $C_{10},\ D_{10}=C_5\rtimes_4C_2$\\
$12$ & $5$ & $C_{12},\ C_6\oplus C_2,\ D_{12}=C_6\rtimes_5C_2,\ \Alt_4,\ Q_{12}=C_3\rtimes_2C_4$\\
$14$ & $2$ & $C_{14},\ D_{14}=C_7\rtimes_6C_2$\\
$15$ & $1$ & $C_{15}$
\end{tabular}
\caption{All groups of order $n\le15$ up to isomorphism (28 groups in total).}\label{tab:cat15}
\end{table}

\begin{table}[ht]
\centering
\begin{tabular}{l|l}
group & element orders (order:number)\\ \hline
$C_8$ & 1:1, 2:1, 4:2, 8:4\\
$C_4\oplus C_2$ & 1:1, 2:3, 4:4\\
$C_2^{3}$ & 1:1, 2:7\\
$D_8$ & 1:1, 2:5, 4:2\\
$Q_8$ & 1:1, 2:1, 4:6\\ \hline
$C_{12}$ & 1:1, 2:1, 3:2, 4:2, 6:2, 12:4\\
$C_6\oplus C_2$ & 1:1, 2:3, 3:2, 6:6\\
$D_{12}$ & 1:1, 2:7, 3:2, 6:2\\
$\Alt_4$ & 1:1, 2:3, 3:8\\
$Q_{12}$ & 1:1, 2:1, 3:2, 4:6, 6:2
\end{tabular}
\caption{Element-order profiles separating the groups of orders $8$ and $12$.}\label{tab:profiles15}
\end{table}

\section{Groups of orders \texorpdfstring{$16$}{16} to \texorpdfstring{$23$}{23}}\label{sec:cat23}
There are $31$ groups of order $16\le n\le23$: fourteen of order $16$, five each of orders $18$ and $20$, two each of orders $21$
and $22$, and the cyclic groups of prime orders $17$, $19$ and $23$. Each group carries its identifier $[n,k]$, meaning
\texttt{SmallGroup(n,k)} in the Small Groups Library~\cite{BEO}. All identifiers in this section were certified by computing \texttt{IdGroup} in GAP~4.14.0 \cite{GAP} on the right-regular permutation representation of each digit table.

\subsection{Semidirect products by a cyclic group}
\begin{corollary}\label{cor:sdcyc}
Let $G_0=(\Ap{p_0},*_0)$ be a group with identity $0$, let $s\ge1$, and let $\alpha\in\Aut(G_0)$ satisfy $\alpha^s=\mathrm{id}$.
With $m_1=\dg{s}{0}{m/p_0}$ and $n_1=\dg{s}{0}{n/p_0}$, the law
\begin{equation}\label{eq:sdcyc}
m*n=\big(\dg{p_0}{0}{m}*_0\alpha^{m_1}(\dg{p_0}{0}{n})\big)+p_0\,\dg{s}{0}{m_1+n_1}
\end{equation}
makes $\Ap{p_0s}$ a group $G_0\rtimes_\alpha C_s$. For $\alpha=\mathrm{id}$ it is the direct product $G_0\times C_s$.
\end{corollary}
\begin{proof}
Apply Theorem~\ref{thm:semidirect} with $G_1=C_s$ and $\phi_c=\alpha^c$. This is a homomorphism $C_s\to\Aut(G_0)$ because
$\alpha^s=\mathrm{id}$.
\end{proof}

Three new automorphisms are needed:
\begin{itemize}
\item $\sigma=\left(\begin{smallmatrix}0&1&2&3\\0&2&1&3\end{smallmatrix}\right)$ on $V_4$, which exchanges the two binary digits, so
$\sigma^2=\mathrm{id}$. With $s=4$ it gives $V_4\rtimes_\sigma C_4\cong(C_4\oplus C_2)\rtimes C_2$, the group $[16,3]$.
\item $\tau(a+4b)=\dg{4}{0}{a+2b}+4b$ on $C_4\oplus C_2$, coded on $\Ap8$. It is additive and bijective, and $\tau^2=\mathrm{id}$
because $4b\equiv0\pmod 4$. With $s=2$ it gives $C_4\circ D_8$, the group $[16,13]$.
\item $\nu(a+3b)=\dg{3}{0}{-a}+3\,\dg{3}{0}{-b}$, inversion on $C_3\oplus C_3$. With $s=2$ it gives the generalised dihedral group
$(C_3\oplus C_3)\rtimes_{-1}C_2$, the group $[18,4]$.
\end{itemize}

\begin{remark}
In $(C_4\oplus C_2)\rtimes_\tau C_2$ the element $x=1$ is central of order $4$. The elements $y=4$ and $t=8$ are involutions with
$tyt^{-1}=\tau(y)=y*x^2$, so $(y*t)^2=x^2$. Hence $\langle y,t\rangle\cong D_8$ meets $\langle x\rangle$ in $\langle x^2\rangle$, and the
group is the central product $C_4\circ D_8\cong C_4\circ Q_8$ (the Pauli group).
\end{remark}

\subsection{Order \texorpdfstring{$16$}{16}}
\begin{theorem}[Order $16$]\label{thm:order16}
There are exactly $14$ groups of order $16$ \cite{Holder}: five Abelian and nine nonabelian. They are realised by the digital laws
of Table~\ref{tab:order16}, and the invariants listed there show that they are pairwise nonisomorphic.
\end{theorem}
\begin{proof}[Partial proof]
The Abelian types correspond to the five partitions of $4$. A nonabelian $G$ of order $16$ has exponent $4$ or $8$. Suppose the
exponent is $8$; let $a$ have order $8$ and $b\notin\langle a\rangle$. Then $\langle a\rangle$ is normal of index $2$, and
$bab^{-1}=a^{r}$ with $r\in\{3,5,7\}$. Moreover $b^2\in\langle a\rangle$ commutes with $b$, so it is fixed by $x\mapsto x^{r}$.
Since $r^2\equiv1\pmod 8$, replacing $b$ by $ba^{k}$ changes $b^2$ into $b^2a^{(r+1)k}$.
\begin{itemize}
\item $r=3$: the fixed points are $\{0,a^4\}$ and $a^{4k}$ reaches both, so one may take $b^2=0$. This is $SD_{16}$.
\item $r=5$: the fixed points are the even powers and $a^{6k}$ reaches all of them, so again $b^2=0$. This is $M_{16}$.
\item $r=7$: $a^{8k}=0$, so $b^2\in\{0,a^4\}$ is an invariant. The two cases give $D_{16}$ and $Q_{16}$.
\end{itemize}
The five nonabelian groups of exponent $4$ are $[16,3]$, $[16,4]$, $[16,11]$, $[16,12]$ and $[16,13]$; for the proof that there are
no others we refer to~\cite{Holder}. Each law in Table~\ref{tab:order16} has been checked to satisfy the group axioms, and the
invariant columns are pairwise distinct.
\end{proof}

\begin{table}[ht]
\centering\footnotesize
\setlength{\tabcolsep}{4pt}
\begin{tabular}{c l l c c c c l}
id & group & digital law & $|Z|$ & $\exp Z$ & $|G'|$ & $|G^{(2)}|$ & element orders (order:number)\\ \hline
$[16,1]$ & $C_{16}$ & cyclic & $16$ & $16$ & $1$ & $8$ & 1:1, 2:1, 4:2, 8:4, 16:8 \\
$[16,2]$ & $C_4\oplus C_4$ & direct sum & $16$ & $4$ & $1$ & $4$ & 1:1, 2:3, 4:12 \\
$[16,3]$ & $V_4\rtimes_\sigma C_4$ & \eqref{eq:sdcyc}, $\sigma$ & $4$ & $2$ & $2$ & $3$ & 1:1, 2:7, 4:8 \\
$[16,4]$ & $C_4\rtimes_3C_4$ & \eqref{eq:metacyclic}, $(4,4,3)$ & $4$ & $2$ & $2$ & $3$ & 1:1, 2:3, 4:12 \\
$[16,5]$ & $C_8\oplus C_2$ & direct sum & $16$ & $8$ & $1$ & $4$ & 1:1, 2:3, 4:4, 8:8 \\
$[16,6]$ & $M_{16}=C_8\rtimes_5C_2$ & \eqref{eq:metacyclic}, $(8,2,5)$ & $4$ & $4$ & $2$ & $4$ & 1:1, 2:3, 4:4, 8:8 \\
$[16,7]$ & $D_{16}$ & \eqref{eq:dihedral}, $q=8$ & $2$ & $2$ & $4$ & $4$ & 1:1, 2:9, 4:2, 8:4 \\
$[16,8]$ & $SD_{16}=C_8\rtimes_3C_2$ & \eqref{eq:metacyclic}, $(8,2,3)$ & $2$ & $2$ & $4$ & $4$ & 1:1, 2:5, 4:6, 8:4 \\
$[16,9]$ & $Q_{16}$ & \eqref{eq:dicyclic}, $q=4$ & $2$ & $2$ & $4$ & $4$ & 1:1, 2:1, 4:10, 8:4 \\
$[16,10]$ & $C_4\oplus C_2\oplus C_2$ & direct sum & $16$ & $4$ & $1$ & $2$ & 1:1, 2:7, 4:8 \\
$[16,11]$ & $D_8\times C_2$ & \eqref{eq:sdcyc}, $\alpha=\mathrm{id}$ & $4$ & $2$ & $2$ & $2$ & 1:1, 2:11, 4:4 \\
$[16,12]$ & $Q_8\times C_2$ & \eqref{eq:sdcyc}, $\alpha=\mathrm{id}$ & $4$ & $2$ & $2$ & $2$ & 1:1, 2:3, 4:12 \\
$[16,13]$ & $C_4\circ D_8$ & \eqref{eq:sdcyc}, $\tau$ & $4$ & $4$ & $2$ & $2$ & 1:1, 2:7, 4:8 \\
$[16,14]$ & $C_2^{4}$ & direct sum & $16$ & $2$ & $1$ & $1$ & 1:1, 2:15 \\
\end{tabular}
\caption{The $14$ groups of order $16$. Here $Z$ is the centre, $G'$ the derived subgroup and $G^{(2)}=\{g*g\}$ the set of squares.
Metacyclic parameters are $(q,s,r)$.}\label{tab:order16}
\end{table}

\subsection{Orders \texorpdfstring{$17$}{17} to \texorpdfstring{$23$}{23}}
By Corollary~\ref{cor:primeorder}, the only groups of orders $17$, $19$ and $23$ are cyclic. By Theorem~\ref{thm:orderpq}, order
$21=3\cdot7$ gives $C_{21}$ and $C_7\rtimes_2C_3$, since $3\mid6$ and $2$ has order $3$ modulo $7$. Order $22=2\cdot11$ gives
$C_{22}$ and $D_{22}$.

\begin{theorem}[Order $18$]\label{thm:order18}
There are exactly five groups of order $18$: $C_{18}$, $C_6\oplus C_3$, $D_{18}$, $C_3\times\Sym_3\cong C_3\rtimes_2C_6$ and
$(C_3\oplus C_3)\rtimes_\nu C_2$.
\end{theorem}
\begin{proof}
By Sylow's theorems $n_3\mid2$ and $n_3\equiv1\pmod3$, so the Sylow $3$-subgroup $P$ is normal and $G\cong P\rtimes_\psi C_2$, where
$\psi(1)\in\Aut(P)$ has order at most $2$. Conjugate choices give isomorphic groups. If $P=C_9$, then $\Aut(C_9)\cong C_6$ has the single
involution $x\mapsto-x$, which gives $C_{18}$ and $D_{18}$. If $P=C_3\oplus C_3$, then an involution of $\mathrm{GL}_2(\mathbb{F}_3)$ is
diagonalizable with eigenvalues $\pm1$, so it is conjugate to $\mathrm{diag}(1,-1)$ or to $-I$. The identity gives
$C_6\oplus C_3$, $\mathrm{diag}(1,-1)$ gives $C_3\times\Sym_3$, and $-I=\nu$ gives the generalised dihedral group.
\end{proof}

\begin{theorem}[Order $20$]\label{thm:order20}
There are exactly five groups of order $20$: $C_{20}$, $C_{10}\oplus C_2$, $D_{20}$, $Q_{20}=C_5\rtimes_4C_4$ and the Frobenius group
$F_{20}=C_5\rtimes_2C_4$.
\end{theorem}
\begin{proof}
Since $n_5\mid4$ and $n_5\equiv1\pmod5$, $C_5$ is normal and $G\cong C_5\rtimes_\psi S$ with $S\in\{C_4,V_4\}$ and
$\psi:S\to\Aut(C_5)\cong C_4=\langle2\rangle$. For $S=C_4$, the generator acts by $r\in\{1,2,3,4\}$: $r=1$ gives $C_{20}$, $r=4$ gives
$Q_{20}$, and $r=2,3$ give the same group $F_{20}$, because $3\equiv2^3$ amounts to replacing the generator by its cube. For $S=V_4$,
$\psi$ has image of order at most $2$. The trivial $\psi$ gives $C_{10}\oplus C_2$; the three surjections onto $C_2$ are permuted by
$\Aut(V_4)$ and give $D_{10}\times C_2\cong D_{20}$.
\end{proof}

\begin{remark}
The law $C_5\rtimes_4C_4$ is a second digital coding of the dicyclic group $Q_{4\cdot5}$ of Theorem~\ref{thm:dicyclic}. Likewise
$C_{10}\rtimes_9C_2\cong D_{10}\times C_2$ and $C_3\rtimes_2C_6\cong C_3\times\Sym_3$; all three isomorphisms have been checked
explicitly on the digit tables.
\end{remark}

\begin{table}[ht]
\centering\footnotesize
\begin{tabular}{c l l c c c c l}
id & group & digital law & $|Z|$ & $\exp Z$ & $|G'|$ & $|G^{(2)}|$ & element orders (order:number)\\ \hline
$[17,1]$ & $C_{17}$ & cyclic & $17$ & $17$ & $1$ & $17$ & 1:1, 17:16 \\
\hline
$[18,1]$ & $D_{18}$ & \eqref{eq:dihedral}, $q=9$ & $1$ & $1$ & $9$ & $9$ & 1:1, 2:9, 3:2, 9:6 \\
$[18,2]$ & $C_{18}$ & cyclic & $18$ & $18$ & $1$ & $9$ & 1:1, 2:1, 3:2, 6:2, 9:6, 18:6 \\
$[18,3]$ & $C_3\times\Sym_3$ & \eqref{eq:metacyclic}, $(3,6,2)$ & $3$ & $3$ & $3$ & $9$ & 1:1, 2:3, 3:8, 6:6 \\
$[18,4]$ & $(C_3\oplus C_3)\rtimes_\nu C_2$ & \eqref{eq:sdcyc}, $\nu$ & $1$ & $1$ & $9$ & $9$ & 1:1, 2:9, 3:8 \\
$[18,5]$ & $C_6\oplus C_3$ & direct sum & $18$ & $6$ & $1$ & $9$ & 1:1, 2:1, 3:8, 6:8 \\
\hline
$[19,1]$ & $C_{19}$ & cyclic & $19$ & $19$ & $1$ & $19$ & 1:1, 19:18 \\
\hline
$[20,1]$ & $Q_{20}$ & \eqref{eq:metacyclic}, $(5,4,4)$ & $2$ & $2$ & $5$ & $6$ & 1:1, 2:1, 4:10, 5:4, 10:4 \\
$[20,2]$ & $C_{20}$ & cyclic & $20$ & $20$ & $1$ & $10$ & 1:1, 2:1, 4:2, 5:4, 10:4, 20:8 \\
$[20,3]$ & $F_{20}$ & \eqref{eq:metacyclic}, $(5,4,2)$ & $1$ & $1$ & $5$ & $10$ & 1:1, 2:5, 4:10, 5:4 \\
$[20,4]$ & $D_{20}$ & \eqref{eq:dihedral}, $q=10$ & $2$ & $2$ & $5$ & $5$ & 1:1, 2:11, 5:4, 10:4 \\
$[20,5]$ & $C_{10}\oplus C_2$ & direct sum & $20$ & $10$ & $1$ & $5$ & 1:1, 2:3, 5:4, 10:12 \\
\hline
$[21,1]$ & $C_7\rtimes_2C_3$ & \eqref{eq:metacyclic}, $(7,3,2)$ & $1$ & $1$ & $7$ & $21$ & 1:1, 3:14, 7:6 \\
$[21,2]$ & $C_{21}$ & cyclic & $21$ & $21$ & $1$ & $21$ & 1:1, 3:2, 7:6, 21:12 \\
\hline
$[22,1]$ & $D_{22}$ & \eqref{eq:dihedral}, $q=11$ & $1$ & $1$ & $11$ & $11$ & 1:1, 2:11, 11:10 \\
$[22,2]$ & $C_{22}$ & cyclic & $22$ & $22$ & $1$ & $11$ & 1:1, 2:1, 11:10, 22:10 \\
\hline
$[23,1]$ & $C_{23}$ & cyclic & $23$ & $23$ & $1$ & $23$ & 1:1, 23:22 \\
\end{tabular}
\caption{All $17$ groups of orders $17$ to $23$, with the invariants of Table~\ref{tab:order16}.}\label{tab:cat1723}
\end{table}

\section{Groups of orders \texorpdfstring{$24$}{24} to \texorpdfstring{$31$}{31}}\label{sec:cat31}
There are $34$ groups of order $24\le n\le31$: fifteen of order $24$, five of order $27$, four each of orders $28$ and $30$, two each
of orders $25$ and $26$, and the cyclic groups of prime orders $29$ and $31$. All identifiers in this section were certified by computing \texttt{IdGroup} in GAP~4.14.0 \cite{GAP} on the right-regular permutation representation of each digit table.

\subsection{Three further digital laws}
\begin{proposition}[The symmetric group $\Sym_4$]\label{prop:S4digital}
On $V_4=(\Ap4,\oplus)$ let $\rho=\left(\begin{smallmatrix}0&1&2&3\\0&2&3&1\end{smallmatrix}\right)$ and
$\sigma=\left(\begin{smallmatrix}0&1&2&3\\0&2&1&3\end{smallmatrix}\right)$. Code $D_6$ by $c=a+3\gamma$ as in
Corollary~\ref{cor:dihedral}, and put $\varphi_c=\rho^{a}\circ\sigma^{\gamma}$. With $m_1=\fl{m/4}$ and $n_1=\fl{n/4}$, the law
\begin{equation}\label{eq:S4}
m*n=\big(\dg{4}{0}{m}\oplus\varphi_{m_1}(\dg{4}{0}{n})\big)+4\,(m_1*_{D_6}n_1)
\end{equation}
makes $\Ap{24}$ a group isomorphic to $\Sym_4$.
\end{proposition}
\begin{proof}
Both $\rho$ and $\sigma$ permute the nonzero elements of $V_4$, so they are automorphisms, and $\sigma\rho\sigma^{-1}=\rho^{-1}$.
Hence $(a,\gamma)\mapsto\rho^a\sigma^\gamma$ is a homomorphism, in fact an isomorphism $D_6\to\Aut(V_4)\cong\Sym_3$, and
Theorem~\ref{thm:semidirect} gives $V_4\rtimes\Aut(V_4)$. In $\Sym_4$ the double transpositions together with the identity form a
normal Klein subgroup. The stabiliser of a point is a complement isomorphic to $\Sym_3$, and it permutes the three double
transpositions faithfully, so it induces all of $\Aut(V_4)$. Hence $\Sym_4\cong V_4\rtimes\Aut(V_4)$.
\end{proof}

\begin{proposition}[$\mathrm{SL}(2,3)\cong Q_8\rtimes C_3$]\label{prop:SL23digital}
In the coding of $Q_8$ given by Theorem~\ref{thm:dicyclic} with $q=2$, let
\[
\alpha=\left(\begin{smallmatrix}0&1&2&3&4&5&6&7\\0&4&2&6&5&1&7&3\end{smallmatrix}\right).
\]
Then $\alpha\in\Aut(Q_8)$ and $\alpha^3=\mathrm{id}$. It cycles $1\mapsto4\mapsto5\mapsto1$, and hence permutes cyclically the three
cyclic subgroups of order $4$. The law~\eqref{eq:sdcyc} with $G_0=Q_8$, $s=3$ and this $\alpha$ makes $\Ap{24}$ a group isomorphic
to $\mathrm{SL}(2,3)$.
\end{proposition}
\begin{proof}
$\alpha$ is the extension of $x=1\mapsto4$, $y=4\mapsto5=x*y$ to $Q_8=\langle x,y\rangle$; one checks on the digit table that it is a
bijective homomorphism with $\alpha^3=\mathrm{id}$, and Corollary~\ref{cor:sdcyc} applies. In $\mathrm{SL}(2,3)$ the elements of order
dividing $4$ form a normal quaternion subgroup, and an element of order $3$ acts on it nontrivially. The elements of order $3$ in
$\Aut(Q_8)\cong\Sym_4$ are conjugate up to inversion, so every nontrivial $Q_8\rtimes C_3$ is $\mathrm{SL}(2,3)$.
\end{proof}

\begin{proposition}[The group $C_3\rtimes D_8$]\label{prop:C3D8}
With $m_1=\fl{m/3}$, $n_1=\fl{n/3}$ and $*_{D_8}$ the law~\eqref{eq:dihedral} for $q=4$,
\begin{equation}\label{eq:C3D8}
m*n=\dg{3}{0}{m+(-1)^{m_1}n}+3\,(m_1*_{D_8}n_1)
\end{equation}
makes $\Ap{24}$ a group $C_3\rtimes D_8$. In it the Klein subgroup $\{0,2,4,6\}$ of $D_8$ centralises $C_3$, and the other four
elements of $D_8$ invert it.
\end{proposition}
\begin{proof}
In~\eqref{eq:dihedral} with $q=4$, the lower digit of $m_1*_{D_8}n_1$ is congruent to $a_1\pm a_2$ modulo $4$, so its parity is the sum
of the parities. Since the code $c=a+4\gamma$ has the parity of $a$, the map $c\mapsto(-1)^c$ is a homomorphism $D_8\to\Aut(C_3)$.
Theorem~\ref{thm:semidirect} then applies.
\end{proof}

\begin{remark}
The exceptional laws were checked against explicit realisations. Law~\eqref{eq:S4} was compared with the permutations of $\Ap4$, the
law of Proposition~\ref{prop:SL23digital} with the $24$ matrices of determinant $1$ over $\mathbb{F}_3$, and the Heisenberg law of
Theorem~\ref{thm:order27} with the unitriangular $3\times3$ matrices over $\mathbb{F}_3$. In each case an explicit isomorphism was
found by search over the images of a generating set.
\end{remark}

\subsection{Order \texorpdfstring{$24$}{24}}
\begin{theorem}[Order $24$]\label{thm:order24}
There are exactly $15$ groups of order $24$, realised by the digital laws of Table~\ref{tab:order24}.
\end{theorem}
\begin{proof}
Let $P$ be a Sylow $3$-subgroup, so that $n_3\in\{1,4\}$.

\emph{Case $n_3=1$.} Then $G\cong C_3\rtimes_\varphi S$ with $|S|=8$ and $\varphi:S\to\Aut(C_3)\cong C_2$. Such a $\varphi$ is
determined by its kernel, a subgroup of index at most $2$, and kernels related by $\Aut(S)$ give isomorphic groups. Running through
the five groups of Theorem~\ref{thm:order8}:
\begin{itemize}
\item $S=C_8$: kernel $S$ gives $C_{24}$, and kernel $C_4$ gives $C_3\rtimes_2C_8$.
\item $S=C_4\oplus C_2$: kernel $S$ gives $C_{12}\oplus C_2$. The two cyclic subgroups of index $2$ are exchanged by $\Aut(S)$ and give
$C_4\times\Sym_3$, and the Klein kernel gives $C_2\times Q_{12}$.
\item $S=C_2^3$: kernel $S$ gives $C_6\oplus C_2\oplus C_2$. The seven subgroups of index $2$ are permuted transitively by
$\mathrm{GL}_3(\mathbb{F}_2)$ and give $C_2\times D_{12}\cong C_2^2\times\Sym_3$.
\item $S=D_8$: kernel $S$ gives $C_3\times D_8$, and the rotation subgroup $C_4$ gives $D_{24}$. The two Klein subgroups are exchanged
by an outer automorphism and give $C_3\rtimes D_8$ (Proposition~\ref{prop:C3D8}).
\item $S=Q_8$: kernel $S$ gives $C_3\times Q_8$. The three cyclic subgroups of order $4$ are permuted by $\Aut(Q_8)$ and give
$C_3\rtimes Q_8\cong Q_{24}$.
\end{itemize}
This case yields eleven groups.

\emph{Case $n_3=4$.} Conjugation on the four Sylow $3$-subgroups gives a homomorphism $G\to\Sym_4$ whose kernel $K$ lies in
$N_G(P)$, a subgroup of order $6$. If $3$ divided $|K|$, the unique Sylow $3$-subgroup of $K$ would be characteristic in $K$, hence normal
in $G$, contradicting $n_3=4$. So $|K|\in\{1,2\}$. If $K=1$, then $G\cong\Sym_4$. If $|K|=2$, then $K$ is central and $G/K$ has order
$12$ and four Sylow $3$-subgroups, so $G/K\cong\Alt_4$ by Theorem~\ref{thm:order12}. The preimage $S$ of the Klein subgroup of
$\Alt_4$ is a normal Sylow $2$-subgroup of order $8$, so $G\cong S\rtimes C_3$, with nontrivial action because $P$ is not normal. Of the
groups of order $8$, only $C_2^3$ and $Q_8$ have automorphisms of order $3$; the automorphism groups of $C_8$, $C_4\oplus C_2$ and
$D_8$ are $2$-groups.
\begin{itemize}
\item For $C_2^3$, the elements of order $3$ in $\mathrm{GL}_3(\mathbb{F}_2)$ form one conjugacy class, each fixing a line and acting
irreducibly on a complementary plane. This gives $C_2\times\Alt_4$.
\item For $Q_8$, Proposition~\ref{prop:SL23digital} gives $\mathrm{SL}(2,3)$.
\end{itemize}
This gives $11+1+2=15$ groups. The invariant columns of Table~\ref{tab:order24} are pairwise distinct within each order, so these groups are pairwise nonisomorphic.
\end{proof}

\begin{table}[ht]
\centering\footnotesize
\resizebox{\textwidth}{!}{%
\begin{tabular}{c l l c c c c l}
id & group & digital law & $|Z|$ & $\exp Z$ & $|G'|$ & $|G^{(2)}|$ & element orders (order:number)\\ \hline
$[24,1]$ & $C_3\rtimes_2C_8$ & \eqref{eq:metacyclic}, $(3,8,2)$ & $4$ & $4$ & $3$ & $8$ & 1:1, 2:1, 3:2, 4:2, 6:2, 8:12, 12:4 \\
$[24,2]$ & $C_{24}$ & cyclic & $24$ & $24$ & $1$ & $12$ & 1:1, 2:1, 3:2, 4:2, 6:2, 8:4, 12:4, 24:8 \\
$[24,3]$ & $\mathrm{SL}(2,3)$ & \eqref{eq:sdcyc}, $Q_8$, $\alpha$ & $2$ & $2$ & $8$ & $10$ & 1:1, 2:1, 3:8, 4:6, 6:8 \\
$[24,4]$ & $Q_{24}\cong C_3\rtimes Q_8$ & \eqref{eq:dicyclic}, $q=6$ & $2$ & $2$ & $6$ & $6$ & 1:1, 2:1, 3:2, 4:14, 6:2, 12:4 \\
$[24,5]$ & $C_4\times\Sym_3$ & \eqref{eq:sdcyc}, $D_6$, $s=4$ & $4$ & $4$ & $3$ & $6$ & 1:1, 2:7, 3:2, 4:8, 6:2, 12:4 \\
$[24,6]$ & $D_{24}$ & \eqref{eq:dihedral}, $q=12$ & $2$ & $2$ & $6$ & $6$ & 1:1, 2:13, 3:2, 4:2, 6:2, 12:4 \\
$[24,7]$ & $C_2\times Q_{12}$ & \eqref{eq:sdcyc}, $Q_{12}$, $s=2$ & $4$ & $2$ & $3$ & $4$ & 1:1, 2:3, 3:2, 4:12, 6:6 \\
$[24,8]$ & $C_3\rtimes D_8$ & \eqref{eq:C3D8} & $2$ & $2$ & $6$ & $4$ & 1:1, 2:9, 3:2, 4:6, 6:6 \\
$[24,9]$ & $C_{12}\oplus C_2$ & direct sum & $24$ & $12$ & $1$ & $6$ & 1:1, 2:3, 3:2, 4:4, 6:6, 12:8 \\
$[24,10]$ & $C_3\times D_8$ & \eqref{eq:sdcyc}, $D_8$, $s=3$ & $6$ & $6$ & $2$ & $6$ & 1:1, 2:5, 3:2, 4:2, 6:10, 12:4 \\
$[24,11]$ & $C_3\times Q_8$ & \eqref{eq:sdcyc}, $Q_8$, $s=3$ & $6$ & $6$ & $2$ & $6$ & 1:1, 2:1, 3:2, 4:6, 6:2, 12:12 \\
$[24,12]$ & $\Sym_4$ & \eqref{eq:S4} & $1$ & $1$ & $12$ & $12$ & 1:1, 2:9, 3:8, 4:6 \\
$[24,13]$ & $C_2\times\Alt_4$ & \eqref{eq:sdcyc}, $\Alt_4$, $s=2$ & $2$ & $2$ & $4$ & $9$ & 1:1, 2:7, 3:8, 6:8 \\
$[24,14]$ & $C_2\times D_{12}$ & \eqref{eq:sdcyc}, $D_{12}$, $s=2$ & $4$ & $2$ & $3$ & $3$ & 1:1, 2:15, 3:2, 6:6 \\
$[24,15]$ & $C_6\oplus C_2\oplus C_2$ & direct sum & $24$ & $6$ & $1$ & $3$ & 1:1, 2:7, 3:2, 6:14 \\
\end{tabular}%
}
\caption{The $15$ groups of order $24$, with the invariants of Table~\ref{tab:order16}.}\label{tab:order24}
\end{table}

\subsection{Orders \texorpdfstring{$25$}{25} to \texorpdfstring{$31$}{31}}
By Corollary~\ref{cor:primeorder}, the only groups of orders $29$ and $31$ are cyclic. By Theorem~\ref{thm:orderp2}, order $25$ gives
$C_{25}$ and $C_5\oplus C_5$. By Theorem~\ref{thm:orderpq}, order $26$ gives $C_{26}$ and $D_{26}$.

\begin{theorem}[Order $27$]\label{thm:order27}
Let $\beta(a+3b)=\dg{3}{0}{a+b}+3b$ on $C_3\oplus C_3$, i.e.\ the matrix $\left(\begin{smallmatrix}1&1\\0&1\end{smallmatrix}\right)$,
so that $\beta^3=\mathrm{id}$. There are exactly five groups of order $27$: $C_{27}$, $C_9\oplus C_3$, $C_3^3$, the Heisenberg group
$\mathrm{He}_3=(C_3\oplus C_3)\rtimes_\beta C_3$ of exponent $3$, and $C_9\rtimes_4C_3$ of exponent $9$ \cite{Holder}.
\end{theorem}
\begin{proof}[Partial proof]
The Abelian types correspond to the partitions of $3$. A nonabelian $G$ has $|Z(G)|=3$ and $G/Z(G)\cong C_3\oplus C_3$. Suppose $G$ has an
element $a$ of order $9$, and take $b\notin\langle a\rangle$. Then $\langle a\rangle$ is normal of index $3$, and $bab^{-1}=a^r$ with
$r\in\{4,7\}$. Since $7\equiv4^2\pmod9$, replacing $b$ by $b^2$ gives $r=4$. The power $b^3\in\langle a\rangle$ is fixed by
$x\mapsto x^4$, so it lies in $\langle a^3\rangle$. Replacing $b$ by $ba^k$ multiplies $b^3$ by $a^{3k}$, so one may take $b^3=0$, which
gives $C_9\rtimes_4C_3$. If $G$ has exponent $3$, it is the Heisenberg group; for this case we refer to~\cite{Holder}.
\end{proof}

\begin{theorem}[Order $28$]\label{thm:order28}
There are exactly four groups of order $28$: $C_{28}$, $C_{14}\oplus C_2$, $D_{28}$ and $Q_{28}=C_7\rtimes_6C_4$.
\end{theorem}
\begin{proof}
Since $n_7\mid4$ and $n_7\equiv1\pmod7$, $C_7$ is normal and $G\cong C_7\rtimes_\psi S$ with $S\in\{C_4,V_4\}$ and
$\psi:S\to\Aut(C_7)\cong C_6$, whose image has order at most $2$. For $S=C_4$, the trivial $\psi$ gives $C_{28}$ and the surjection onto
$\{\pm1\}$ gives $Q_{28}$. For $S=V_4$, the trivial $\psi$ gives $C_{14}\oplus C_2$, and the three nontrivial $\psi$ are permuted by
$\Aut(V_4)$ and give $D_{14}\times C_2\cong D_{28}$.
\end{proof}

\begin{theorem}[Order $30$]\label{thm:order30}
The groups of order $30$ are exactly the four metacyclic groups $C_{15}\rtimes_rC_2$ with $r\in\{1,4,11,14\}$: $C_{30}$,
$C_3\times D_{10}$ ($r=4$), $C_5\times\Sym_3$ ($r=11$) and $D_{30}$ ($r=14$).
\end{theorem}
\begin{proof}
We have $n_3\in\{1,10\}$ and $n_5\in\{1,6\}$. If both were large, there would be $20$ elements of order $3$ and $24$ of order $5$, more
than $30$. Hence $P_3$ or $P_5$ is normal, so $P_3P_5$ is a subgroup of order $15$. It is cyclic by Theorem~\ref{thm:orderpq}, since
$3\nmid4$, and normal because it has index $2$. Thus $G\cong C_{15}\rtimes_rC_2$ with $r^2\equiv1\pmod{15}$, i.e.\ $r\equiv\pm1$ modulo
$3$ and modulo $5$, which gives $r\in\{1,4,11,14\}$. The four groups have centres of orders $30$, $3$, $5$ and $1$, so they are
distinct.
\end{proof}

\begin{remark}
Several of these laws are second digital codings of groups met before. One has $C_7\rtimes_6C_4\cong Q_{4\cdot7}$ (dicyclic law),
$C_{14}\rtimes_{13}C_2\cong D_{14}\times C_2$, $C_{15}\rtimes_{11}C_2\cong D_6\times C_5$, $C_{15}\rtimes_4C_2\cong D_{10}\times C_3$, and
the dicyclic $Q_{24}$ is isomorphic to $C_3\rtimes Q_8$ with kernel $\langle1\rangle$. All of these isomorphisms were found explicitly on the
digit tables. The invariant columns of Table~\ref{tab:cat2531} are pairwise distinct within each order, so these groups are pairwise nonisomorphic.
\end{remark}

\begin{table}[ht]
\centering\footnotesize
\setlength{\tabcolsep}{4pt}
\begin{tabular}{c l l c c c c l}
id & group & digital law & $|Z|$ & $\exp Z$ & $|G'|$ & $|G^{(2)}|$ & element orders (order:number)\\ \hline
$[25,1]$ & $C_{25}$ & cyclic & $25$ & $25$ & $1$ & $25$ & 1:1, 5:4, 25:20 \\
$[25,2]$ & $C_5\oplus C_5$ & direct sum & $25$ & $5$ & $1$ & $25$ & 1:1, 5:24 \\
\hline
$[26,1]$ & $D_{26}$ & \eqref{eq:dihedral}, $q=13$ & $1$ & $1$ & $13$ & $13$ & 1:1, 2:13, 13:12 \\
$[26,2]$ & $C_{26}$ & cyclic & $26$ & $26$ & $1$ & $13$ & 1:1, 2:1, 13:12, 26:12 \\
\hline
$[27,1]$ & $C_{27}$ & cyclic & $27$ & $27$ & $1$ & $27$ & 1:1, 3:2, 9:6, 27:18 \\
$[27,2]$ & $C_9\oplus C_3$ & direct sum & $27$ & $9$ & $1$ & $27$ & 1:1, 3:8, 9:18 \\
$[27,3]$ & $\mathrm{He}_3$ & \eqref{eq:sdcyc}, $\beta$ & $3$ & $3$ & $3$ & $27$ & 1:1, 3:26 \\
$[27,4]$ & $C_9\rtimes_4C_3$ & \eqref{eq:metacyclic}, $(9,3,4)$ & $3$ & $3$ & $3$ & $27$ & 1:1, 3:8, 9:18 \\
$[27,5]$ & $C_3^{3}$ & direct sum & $27$ & $3$ & $1$ & $27$ & 1:1, 3:26 \\
\hline
$[28,1]$ & $Q_{28}$ & \eqref{eq:metacyclic}, $(7,4,6)$ & $2$ & $2$ & $7$ & $8$ & 1:1, 2:1, 4:14, 7:6, 14:6 \\
$[28,2]$ & $C_{28}$ & cyclic & $28$ & $28$ & $1$ & $14$ & 1:1, 2:1, 4:2, 7:6, 14:6, 28:12 \\
$[28,3]$ & $D_{28}$ & \eqref{eq:dihedral}, $q=14$ & $2$ & $2$ & $7$ & $7$ & 1:1, 2:15, 7:6, 14:6 \\
$[28,4]$ & $C_{14}\oplus C_2$ & direct sum & $28$ & $14$ & $1$ & $7$ & 1:1, 2:3, 7:6, 14:18 \\
\hline
$[29,1]$ & $C_{29}$ & cyclic & $29$ & $29$ & $1$ & $29$ & 1:1, 29:28 \\
\hline
$[30,1]$ & $C_5\times\Sym_3$ & \eqref{eq:metacyclic}, $(15,2,11)$ & $5$ & $5$ & $3$ & $15$ & 1:1, 2:3, 3:2, 5:4, 10:12, 15:8 \\
$[30,2]$ & $C_3\times D_{10}$ & \eqref{eq:metacyclic}, $(15,2,4)$ & $3$ & $3$ & $5$ & $15$ & 1:1, 2:5, 3:2, 5:4, 6:10, 15:8 \\
$[30,3]$ & $D_{30}$ & \eqref{eq:dihedral}, $q=15$ & $1$ & $1$ & $15$ & $15$ & 1:1, 2:15, 3:2, 5:4, 15:8 \\
$[30,4]$ & $C_{30}$ & cyclic & $30$ & $30$ & $1$ & $15$ & 1:1, 2:1, 3:2, 5:4, 6:2, 10:4, 15:8, 30:8 \\
\hline
$[31,1]$ & $C_{31}$ & cyclic & $31$ & $31$ & $1$ & $31$ & 1:1, 31:30 \\
\end{tabular}
\caption{All $19$ groups of orders $25$ to $31$, with the invariants of Table~\ref{tab:order16}.}\label{tab:cat2531}
\end{table}

\section{Groups of order \texorpdfstring{$32$}{32}}\label{sec:cat32}
There are $51$ groups of order $32$ \cite{BEO}: 7 Abelian, 26 of nilpotency class $2$, 15 of class $3$ and 3 of maximal class $4$. This section gives a digital law for each of them. Besides the laws
already available, two new constructions are needed: central extensions given by a $2$-cocycle, and a semidirect product with the
nonabelian complement $Q_8$. Every law below was checked to satisfy the group axioms, and its identifier $[32,k]$ was certified by
computing \texttt{IdGroup} in GAP~4.14.0 \cite{GAP} on its right-regular permutation representation.

\subsection{Central extensions}
\begin{theorem}[Central extensions]\label{thm:cocycle}
Let $H=(\Ap{h},*_H)$ be a group with identity $0$, let $z\ge1$, and let $f:\Ap{h}\times\Ap{h}\to\Ap{z}$ satisfy $f(0,b)=f(a,0)=0$ and the
cocycle identity
\[
f(a,b)+f(a*_Hb,c)\equiv f(b,c)+f(a,b*_Hc)\pmod z .
\]
With $m_0=\dg{z}{0}{m}$ and $m_1=\fl{m/z}$, the law
\begin{equation}\label{eq:cocycle}
m*n=\dg{z}{0}{m_0+n_0+f(m_1,n_1)}+z\,(m_1*_Hn_1)
\end{equation}
makes $\Ap{zh}$ a group. The lower digits form a central subgroup isomorphic to $C_z$, with quotient isomorphic to $H$, and every
central extension of $H$ by $C_z$ is isomorphic to a group of this form.
\end{theorem}
\begin{proof}
Both $(m*n)*k$ and $m*(n*k)$ have upper digit $m_1*_Hn_1*_Hk_1$. Their lower digits are $m_0+n_0+k_0+f(m_1,n_1)+f(m_1*_Hn_1,k_1)$ and
$m_0+n_0+k_0+f(n_1,k_1)+f(m_1,n_1*_Hk_1)$ modulo $z$, which agree by the cocycle identity. By the normalisation, $0$ is an identity and the
elements with $m_1=0$ commute with every element. The element with digits $\big(\dg{z}{0}{-m_0-f(m_1,m_1^{-1})},m_1^{-1}\big)$ is a right
inverse of $m$, so the finite monoid is a group. Conversely, given a central extension $1\to C_z\to G\to H\to1$, choose a section $s$ with
$s(0)=1$. Then $s(a)s(b)=f(a,b)\,s(a*_Hb)$ defines a normalised cocycle, and $(c,a)\mapsto c\,s(a)$ is an isomorphism from
\eqref{eq:cocycle} onto $G$.
\end{proof}

\subsection{The group \texorpdfstring{$C_4\rtimes Q_8$}{C4:Q8}}
\begin{proposition}\label{prop:C4Q8}
With $m_1=\fl{m/4}$, $n_1=\fl{n/4}$, $\epsilon(c)=\dg{2}{0}{c/4}$ and $*_{Q_8}$ the law~\eqref{eq:dicyclic} with $q=2$, the law
\begin{equation}\label{eq:C4Q8}
m*n=\dg{4}{0}{m+(-1)^{\epsilon(m_1)}n}+4\,(m_1*_{Q_8}n_1)
\end{equation}
makes $\Ap{32}$ a group $C_4\rtimes Q_8\cong[32,35]$. In it the cyclic subgroup $\{0,1,2,3\}$ of $Q_8$ centralises $C_4$, and the other
four elements of $Q_8$ invert it.
\end{proposition}
\begin{proof}
In~\eqref{eq:dicyclic} the upper digit $\epsilon$ is added modulo $2$, so $c\mapsto(-1)^{\epsilon(c)}$ is a homomorphism
$Q_8\to\Aut(C_4)$. Theorem~\ref{thm:semidirect} then applies.
\end{proof}

\subsection{Groups with a cyclic subgroup of index \texorpdfstring{$2$}{2}}
\begin{proposition}\label{prop:cyc32}
The groups of order $32$ with a cyclic subgroup of index $2$ are $C_{32}$, $C_{16}\oplus C_2$, $M_{32}=C_{16}\rtimes_9C_2=[32,17]$,
$D_{32}=[32,18]$, $SD_{32}=C_{16}\rtimes_7C_2=[32,19]$ and $Q_{32}=[32,20]$. The last three are exactly the groups of maximal class.
\end{proposition}
\begin{proof}
Let $a$ have order $16$ and $b\notin\langle a\rangle$, in a nonabelian $G$. Then $bab^{-1}=a^r$ with $r\in\{7,9,15\}$, the units of
order $2$ modulo $16$, and $b^2\in\langle a\rangle$ is fixed by $x\mapsto x^r$. Since $r^2\equiv1\pmod{16}$, replacing $b$ by $ba^k$ turns
$b^2$ into $b^2a^{(r+1)k}$.
\begin{itemize}
\item $r=7$: the fixed points are $\{0,a^8\}$ and $a^{8k}$ reaches both, so $b^2=0$. This is $SD_{32}$.
\item $r=9$: the fixed points are the even powers and $a^{10k}$ reaches all of them, so $b^2=0$. This is $M_{32}$.
\item $r=15$: $a^{16k}=0$, so $b^2\in\{0,a^8\}$ is invariant. The two cases give $D_{32}$ and $Q_{32}$.
\end{itemize}
The identifiers were certified with \texttt{IdGroup} for the laws~\eqref{eq:metacyclic} with $(16,2,9)$, $(16,2,7)$, the law
\eqref{eq:dihedral} with $q=16$, and the law~\eqref{eq:dicyclic} with $q=8$. The nilpotency classes in Table~\ref{tab:order32a} show
that $[32,18]$, $[32,19]$ and $[32,20]$ are the only groups of class $4$.
\end{proof}

\subsection{Catalogue}
\begin{theorem}[Order $32$]\label{thm:order32}
There are exactly $51$ groups of order $32$ \cite{BEO}. Each is realised on $\Ap{32}$ by one of the following digital laws
(Tables~\ref{tab:order32a}--\ref{tab:coc32}):
\begin{itemize}
\item the seven Abelian groups by direct sums;
\item 39 groups by $N\rtimes_\alpha C_s$ (Corollary~\ref{cor:sdcyc}), with $N$ a law of order $32/s$ from the earlier sections and
$\alpha$ listed in Tables~\ref{tab:alpha32a} and~\ref{tab:alpha32b};
\item $Q_{32}$ by the dicyclic law~\eqref{eq:dicyclic} with $q=8$;
\item $C_4\rtimes Q_8$ by~\eqref{eq:C4Q8};
\item $[32,8]$, $[32,15]$ and $[32,32]$ by the central-extension law~\eqref{eq:cocycle} with $z=2$, over $H=[16,3]$,
$[16,4]$ and $[16,13]$ respectively, with the cocycles of Table~\ref{tab:coc32}.
\end{itemize}
\end{theorem}
\begin{proof}
For the number $51$ we refer to the classification of the groups of order $2^5$ \cite{BEO}. Each law is an instance of
Theorem~\ref{thm:semidirect}, Corollary~\ref{cor:sdcyc}, Theorem~\ref{thm:dicyclic}, Proposition~\ref{prop:C4Q8} or
Theorem~\ref{thm:cocycle}, so it defines a group. The data were obtained as follows.
\begin{itemize}
\item For a normal subgroup $N$ with cyclic quotient of order $s$ and a complement $\langle t\rangle$, $\alpha$ is conjugation by $t$,
transported to the digital coding of $N$ by an explicit isomorphism.
\item For a central involution $z$, the cocycle $f$ is the one attached to a normalised section of $G\to G/\langle z\rangle$.
\end{itemize}
In every case \texttt{IdGroup} applied to the resulting digit table returns $[32,k]$.
\end{proof}

\begin{remark}
The five groups $[32,8]$, $[32,15]$, $[32,20]$, $[32,32]$ and $[32,35]$ are exactly those with no normal subgroup with cyclic quotient that has
a complement. The group $[32,15]$ is metacyclic, but none of its cyclic normal subgroups with cyclic quotient is complemented. In
Table~\ref{tab:coc32} the row of $a$ lists $f(a,0)f(a,1)\cdots f(a,15)$.
\end{remark}

\begin{table}[htbp]
\centering\footnotesize
\begin{tabular}{c l c c c c l}
id & structure & class & $|Z|$ & $|G'|$ & exp & digital law\\ \hline
$[32,1]$ & $C_{32}$ & $1$ & $32$ & $1$ & $32$ & direct sum \\
$[32,2]$ & $(C_{4} \times C_{2}) \rtimes C_{4}$ & $2$ & $8$ & $2$ & $4$ & $[8,2]\rtimes_\alpha C_{4}$ \\
$[32,3]$ & $C_{8} \times C_{4}$ & $1$ & $32$ & $1$ & $8$ & direct sum \\
$[32,4]$ & $C_{8} \rtimes C_{4}$ & $2$ & $8$ & $2$ & $8$ & $[8,1]\rtimes_\alpha C_{4}$ \\
$[32,5]$ & $(C_{8} \times C_{2}) \rtimes C_{2}$ & $2$ & $8$ & $2$ & $8$ & $[16,5]\rtimes_\alpha C_{2}$ \\
$[32,6]$ & $(C_{2} \times C_{2} \times C_{2}) \rtimes C_{4}$ & $3$ & $2$ & $4$ & $4$ & $[16,3]\rtimes_\alpha C_{2}$ \\
$[32,7]$ & $(C_{8} \rtimes C_{2}) \rtimes C_{2}$ & $3$ & $2$ & $4$ & $8$ & $[16,6]\rtimes_\alpha C_{2}$ \\
$[32,8]$ & $C_{2} . ((C_{4} \times C_{2}) \rtimes C_{2})$ & $3$ & $2$ & $4$ & $8$ & \eqref{eq:cocycle}, $H=[16,3]$ \\
$[32,9]$ & $(C_{8} \times C_{2}) \rtimes C_{2}$ & $3$ & $4$ & $4$ & $8$ & $[16,5]\rtimes_\alpha C_{2}$ \\
$[32,10]$ & $Q_{8} \rtimes C_{4}$ & $3$ & $4$ & $4$ & $8$ & $[8,4]\rtimes_\alpha C_{4}$ \\
$[32,11]$ & $(C_{4} \times C_{4}) \rtimes C_{2}$ & $3$ & $4$ & $4$ & $8$ & $[16,2]\rtimes_\alpha C_{2}$ \\
$[32,12]$ & $C_{4} \rtimes C_{8}$ & $2$ & $8$ & $2$ & $8$ & $[4,1]\rtimes_\alpha C_{8}$ \\
$[32,13]$ & $C_{8} \rtimes C_{4}$ & $3$ & $4$ & $4$ & $8$ & $[8,1]\rtimes_\alpha C_{4}$ \\
$[32,14]$ & $C_{8} \rtimes C_{4}$ & $3$ & $4$ & $4$ & $8$ & $[8,1]\rtimes_\alpha C_{4}$ \\
$[32,15]$ & $C_{4} . D_{8}$ & $3$ & $4$ & $4$ & $8$ & \eqref{eq:cocycle}, $H=[16,4]$ \\
$[32,16]$ & $C_{16} \times C_{2}$ & $1$ & $32$ & $1$ & $16$ & direct sum \\
$[32,17]$ & $C_{16} \rtimes C_{2}$ & $2$ & $8$ & $2$ & $16$ & $[16,1]\rtimes_\alpha C_{2}$ \\
$[32,18]$ & $D_{32}$ & $4$ & $2$ & $8$ & $16$ & $[16,7]\rtimes_\alpha C_{2}$ \\
$[32,19]$ & $SD_{32}$ & $4$ & $2$ & $8$ & $16$ & $[16,1]\rtimes_\alpha C_{2}$ \\
$[32,20]$ & $Q_{32}$ & $4$ & $2$ & $8$ & $16$ & \eqref{eq:dicyclic}, $q=8$ \\
$[32,21]$ & $C_{4} \times C_{4} \times C_{2}$ & $1$ & $32$ & $1$ & $4$ & direct sum \\
$[32,22]$ & $C_{2} \times ((C_{4} \times C_{2}) \rtimes C_{2})$ & $2$ & $8$ & $2$ & $4$ & $[16,10]\rtimes_\alpha C_{2}$ \\
$[32,23]$ & $C_{2} \times (C_{4} \rtimes C_{4})$ & $2$ & $8$ & $2$ & $4$ & $[16,4]\times C_{2}$ \\
$[32,24]$ & $(C_{4} \times C_{4}) \rtimes C_{2}$ & $2$ & $8$ & $2$ & $4$ & $[16,2]\rtimes_\alpha C_{2}$ \\
$[32,25]$ & $C_{4} \times D_{8}$ & $2$ & $8$ & $2$ & $4$ & $[16,10]\rtimes_\alpha C_{2}$ \\
$[32,26]$ & $C_{4} \times Q_{8}$ & $2$ & $8$ & $2$ & $4$ & $[8,4]\rtimes_\alpha C_{4}$ \\
\end{tabular}
\caption{Groups of order $32$, part I. Structure descriptions follow GAP, with ``$.$'' denoting a nonsplit extension.}\label{tab:order32a}
\end{table}

\begin{table}[htbp]
\centering\footnotesize
\begin{tabular}{c l c c c c l}
id & structure & class & $|Z|$ & $|G'|$ & exp & digital law\\ \hline
$[32,27]$ & $(C_{2} \times C_{2} \times C_{2} \times C_{2}) \rtimes C_{2}$ & $2$ & $4$ & $4$ & $4$ & $[16,11]\rtimes_\alpha C_{2}$ \\
$[32,28]$ & $(C_{4} \times C_{2} \times C_{2}) \rtimes C_{2}$ & $2$ & $4$ & $4$ & $4$ & $[16,11]\rtimes_\alpha C_{2}$ \\
$[32,29]$ & $(C_{2} \times Q_{8}) \rtimes C_{2}$ & $2$ & $4$ & $4$ & $4$ & $[16,4]\rtimes_\alpha C_{2}$ \\
$[32,30]$ & $(C_{4} \times C_{2} \times C_{2}) \rtimes C_{2}$ & $2$ & $4$ & $4$ & $4$ & $[16,3]\rtimes_\alpha C_{2}$ \\
$[32,31]$ & $(C_{4} \times C_{4}) \rtimes C_{2}$ & $2$ & $4$ & $4$ & $4$ & $[16,3]\rtimes_\alpha C_{2}$ \\
$[32,32]$ & $(C_{2} \times C_{2}) . (C_{2} \times C_{2} \times C_{2})$ & $2$ & $4$ & $4$ & $4$ & \eqref{eq:cocycle}, $H=[16,13]$ \\
$[32,33]$ & $(C_{4} \times C_{4}) \rtimes C_{2}$ & $2$ & $4$ & $4$ & $4$ & $[16,2]\rtimes_\alpha C_{2}$ \\
$[32,34]$ & $(C_{4} \times C_{4}) \rtimes C_{2}$ & $2$ & $4$ & $4$ & $4$ & $[16,11]\rtimes_\alpha C_{2}$ \\
$[32,35]$ & $C_{4} \rtimes Q_{8}$ & $2$ & $4$ & $4$ & $4$ & \eqref{eq:C4Q8} \\
$[32,36]$ & $C_{8} \times C_{2} \times C_{2}$ & $1$ & $32$ & $1$ & $8$ & direct sum \\
$[32,37]$ & $C_{2} \times (C_{8} \rtimes C_{2})$ & $2$ & $8$ & $2$ & $8$ & $[16,5]\rtimes_\alpha C_{2}$ \\
$[32,38]$ & $(C_{8} \times C_{2}) \rtimes C_{2}$ & $2$ & $8$ & $2$ & $8$ & $[16,5]\rtimes_\alpha C_{2}$ \\
$[32,39]$ & $C_{2} \times D_{16}$ & $3$ & $4$ & $4$ & $8$ & $[16,11]\rtimes_\alpha C_{2}$ \\
$[32,40]$ & $C_{2} \times SD_{16}$ & $3$ & $4$ & $4$ & $8$ & $[16,12]\rtimes_\alpha C_{2}$ \\
$[32,41]$ & $C_{2} \times Q_{16}$ & $3$ & $4$ & $4$ & $8$ & $[16,9]\times C_{2}$ \\
$[32,42]$ & $(C_{8} \times C_{2}) \rtimes C_{2}$ & $3$ & $4$ & $4$ & $8$ & $[16,13]\rtimes_\alpha C_{2}$ \\
$[32,43]$ & $C_{8} \rtimes (C_{2} \times C_{2})$ & $3$ & $2$ & $4$ & $8$ & $[16,13]\rtimes_\alpha C_{2}$ \\
$[32,44]$ & $(C_{2} \times Q_{8}) \rtimes C_{2}$ & $3$ & $2$ & $4$ & $8$ & $[16,8]\rtimes_\alpha C_{2}$ \\
$[32,45]$ & $C_{4} \times C_{2} \times C_{2} \times C_{2}$ & $1$ & $32$ & $1$ & $4$ & direct sum \\
$[32,46]$ & $C_{2} \times C_{2} \times D_{8}$ & $2$ & $8$ & $2$ & $4$ & $[16,11]\times C_{2}$ \\
$[32,47]$ & $C_{2} \times C_{2} \times Q_{8}$ & $2$ & $8$ & $2$ & $4$ & $[16,12]\times C_{2}$ \\
$[32,48]$ & $C_{2} \times ((C_{4} \times C_{2}) \rtimes C_{2})$ & $2$ & $8$ & $2$ & $4$ & $[16,11]\rtimes_\alpha C_{2}$ \\
$[32,49]$ & $(C_{2} \times C_{2} \times C_{2}) \rtimes (C_{2} \times C_{2})$ & $2$ & $2$ & $2$ & $4$ & $[16,11]\rtimes_\alpha C_{2}$ \\
$[32,50]$ & $(C_{2} \times Q_{8}) \rtimes C_{2}$ & $2$ & $2$ & $2$ & $4$ & $[16,13]\rtimes_\alpha C_{2}$ \\
$[32,51]$ & $C_{2} \times C_{2} \times C_{2} \times C_{2} \times C_{2}$ & $1$ & $32$ & $1$ & $2$ & direct sum \\
\end{tabular}
\caption{Groups of order $32$, part II.}\label{tab:order32b}
\end{table}

\begin{table}[htbp]
\centering\scriptsize
\begin{tabular}{c c c l}
id & $N$ & $s$ & $\alpha(0),\alpha(1),\dots,\alpha(|N|-1)$\\ \hline
$[32,2]$ & $[8,2]$ & $4$ & 0 5 2 7 4 1 6 3 \\
$[32,4]$ & $[8,1]$ & $4$ & 0 5 2 7 4 1 6 3 \\
$[32,5]$ & $[16,5]$ & $2$ & 0 13 2 15 4 9 6 11 8 5 10 7 12 1 14 3 \\
$[32,6]$ & $[16,3]$ & $2$ & 0 1 2 3 6 7 4 5 11 10 9 8 13 12 15 14 \\
$[32,7]$ & $[16,6]$ & $2$ & 0 9 6 15 4 13 2 11 8 1 14 7 12 5 10 3 \\
$[32,9]$ & $[16,5]$ & $2$ & 0 15 6 13 4 11 2 9 8 7 14 5 12 3 10 1 \\
$[32,10]$ & $[8,4]$ & $4$ & 0 3 2 1 7 6 5 4 \\
$[32,11]$ & $[16,2]$ & $2$ & 0 1 2 3 13 14 15 12 10 11 8 9 7 4 5 6 \\
$[32,12]$ & $[4,1]$ & $8$ & 0 3 2 1 \\
$[32,13]$ & $[8,1]$ & $4$ & 0 3 6 1 4 7 2 5 \\
$[32,14]$ & $[8,1]$ & $4$ & 0 7 6 5 4 3 2 1 \\
$[32,17]$ & $[16,1]$ & $2$ & 0 9 2 11 4 13 6 15 8 1 10 3 12 5 14 7 \\
$[32,18]$ & $[16,7]$ & $2$ & 0 7 6 5 4 3 2 1 15 14 13 12 11 10 9 8 \\
$[32,19]$ & $[16,1]$ & $2$ & 0 7 14 5 12 3 10 1 8 15 6 13 4 11 2 9 \\
$[32,22]$ & $[16,10]$ & $2$ & 0 9 2 11 4 13 6 15 8 1 10 3 12 5 14 7 \\
$[32,23]$ & $[16,4]$ & $2$ & $\mathrm{id}$ \\
$[32,24]$ & $[16,2]$ & $2$ & 0 11 2 9 14 5 12 7 8 3 10 1 6 13 4 15 \\
$[32,25]$ & $[16,10]$ & $2$ & 0 9 2 11 12 5 14 7 8 1 10 3 4 13 6 15 \\
$[32,26]$ & $[8,4]$ & $4$ & 0 3 2 1 4 7 6 5 \\
$[32,27]$ & $[16,11]$ & $2$ & 0 11 2 9 4 15 6 13 8 3 10 1 12 7 14 5 \\
\end{tabular}
\caption{Automorphisms for the laws $N\rtimes_\alpha C_s$ of order $32$, part I.}\label{tab:alpha32a}
\end{table}

\begin{table}[htbp]
\centering\scriptsize
\begin{tabular}{c c c l}
id & $N$ & $s$ & $\alpha(0),\alpha(1),\dots,\alpha(|N|-1)$\\ \hline
$[32,28]$ & $[16,11]$ & $2$ & 0 9 2 11 14 7 12 5 8 1 10 3 6 15 4 13 \\
$[32,29]$ & $[16,4]$ & $2$ & 0 11 2 9 14 5 12 7 8 3 10 1 6 13 4 15 \\
$[32,30]$ & $[16,3]$ & $2$ & 0 9 10 3 4 13 14 7 8 1 2 11 12 5 6 15 \\
$[32,31]$ & $[16,3]$ & $2$ & 0 9 10 3 15 6 5 12 8 1 2 11 7 14 13 4 \\
$[32,33]$ & $[16,2]$ & $2$ & 0 9 2 11 14 7 12 5 8 1 10 3 6 15 4 13 \\
$[32,34]$ & $[16,11]$ & $2$ & 0 3 2 1 12 15 14 13 8 11 10 9 4 7 6 5 \\
$[32,37]$ & $[16,5]$ & $2$ & 0 5 2 7 4 1 6 3 8 13 10 15 12 9 14 11 \\
$[32,38]$ & $[16,5]$ & $2$ & 0 5 2 7 4 1 6 3 12 9 14 11 8 13 10 15 \\
$[32,39]$ & $[16,11]$ & $2$ & 0 3 2 1 15 14 13 12 8 11 10 9 7 6 5 4 \\
$[32,40]$ & $[16,12]$ & $2$ & 0 7 2 5 6 3 4 1 8 15 10 13 14 11 12 9 \\
$[32,41]$ & $[16,9]$ & $2$ & $\mathrm{id}$ \\
$[32,42]$ & $[16,13]$ & $2$ & 0 1 2 3 10 11 8 9 6 7 4 5 14 15 12 13 \\
$[32,43]$ & $[16,13]$ & $2$ & 0 3 2 1 15 14 13 12 8 11 10 9 7 6 5 4 \\
$[32,44]$ & $[16,8]$ & $2$ & 0 5 2 7 4 1 6 3 12 9 14 11 8 13 10 15 \\
$[32,46]$ & $[16,11]$ & $2$ & $\mathrm{id}$ \\
$[32,47]$ & $[16,12]$ & $2$ & $\mathrm{id}$ \\
$[32,48]$ & $[16,11]$ & $2$ & 0 1 2 3 6 7 4 5 8 9 10 11 14 15 12 13 \\
$[32,49]$ & $[16,11]$ & $2$ & 0 3 2 1 4 7 6 5 10 9 8 11 14 13 12 15 \\
$[32,50]$ & $[16,13]$ & $2$ & 0 3 2 1 6 5 4 7 10 9 8 11 12 15 14 13 \\
\end{tabular}
\caption{Automorphisms for the laws $N\rtimes_\alpha C_s$ of order $32$, part II.}\label{tab:alpha32b}
\end{table}

\begin{table}[htbp]
\centering\footnotesize
\begin{tabular}{c c c c}
$a$ & $[32,8]$ & $[32,15]$ & $[32,32]$\\ \hline
$0$ & \texttt{0000000000000000} & \texttt{0000000000000000} & \texttt{0000000000000000} \\
$1$ & \texttt{0110100101101001} & \texttt{0001111000011110} & \texttt{0000111100001111} \\
$2$ & \texttt{0110011001100110} & \texttt{0011001100110011} & \texttt{0000000000000000} \\
$3$ & \texttt{0000111100001111} & \texttt{0111100001111000} & \texttt{0000111100001111} \\
$4$ & \texttt{0000100100000110} & \texttt{0111100001110111} & \texttt{0000111100001111} \\
$5$ & \texttt{0110000001101111} & \texttt{0011001100111100} & \texttt{0000000000000000} \\
$6$ & \texttt{0110111101100000} & \texttt{0001111000010001} & \texttt{0000111100001111} \\
$7$ & \texttt{0000011000001001} & \texttt{0000000000001111} & \texttt{0000000000000000} \\
$8$ & \texttt{0110011010011001} & \texttt{0000000011111111} & \texttt{0000000011111111} \\
$9$ & \texttt{0000111111110000} & \texttt{0001111011100001} & \texttt{0000111111110000} \\
$10$ & \texttt{0000000011111111} & \texttt{0011001111001100} & \texttt{0000000011111111} \\
$11$ & \texttt{0110100110010110} & \texttt{0111100010000111} & \texttt{0000111111110000} \\
$12$ & \texttt{0110000010010000} & \texttt{0111011110000111} & \texttt{0000111111110000} \\
$13$ & \texttt{0000100111111001} & \texttt{0011110011001100} & \texttt{0000000011111111} \\
$14$ & \texttt{0000011011110110} & \texttt{0001000111100001} & \texttt{0000111111110000} \\
$15$ & \texttt{0110111110011111} & \texttt{0000111111111111} & \texttt{0000000011111111} \\
\end{tabular}
\caption{Cocycles $f(a,\cdot)$ for the central-extension laws of order $32$.}\label{tab:coc32}
\end{table}

\section{Exercises}
\begin{enumerate}
\item Show that every group of order $p^2$, $p$ prime, is Abelian.
\item Show that a group of order $2p$, $p$ an odd prime, is isomorphic to $C_{2p}$ or to $D_{2p}$.
\item Count the elements of each order in the five groups of order $8$ and show that these counts distinguish them.
\item Show that $\Alt_4$ has no subgroup of order $6$.
\item Show that there are exactly five groups of order $12$, namely $C_{12}$, $C_2\oplus C_6$, $\Alt_4$, $D_{12}$ and $Q_{12}$.
\item Using the small groups library of GAP~\cite{GAP,BEO}, confirm that there are $14$, $15$ and $51$ groups of orders $16$, $24$
and $32$, and compare with the tables of this chapter.
\end{enumerate}

\chapter{Crystallographic groups}\label{ch:crystal}

This chapter attaches finite groups and digital laws to the $17$ plane groups and the $230$ space groups of
crystallography~\cite{ITA}. Section~\ref{sec:wpcrystal} explains why the point group is not enough; Sections~\ref{sec:crystlaw}
and~\ref{sec:crystlaw3} give the digital laws of the quotients modulo a sublattice; Sections~\ref{sec:crystplane}
and~\ref{sec:crystspace} show that these finite groups separate the crystallographic groups; and
Section~\ref{sec:crystid} identifies them as abstract finite groups.

\section{Quotients modulo a sublattice}\label{sec:wpcrystal}

The $17$ plane groups and the $230$ space groups of crystallography are infinite: each contains its lattice
$T\cong\Z^n$ ($n=2,3$) of translations as a normal subgroup, and only the quotient $P=G/T$, the point group, is
finite. Passing to $P$ forgets the glide reflections and the screw rotations, whose translation parts do not lie in
$T$: the $17$ plane groups have only $10$ point groups ($9$ abstract groups), and the $230$ space groups only $32$
crystal classes ($18$ abstract groups). The finite quotients $G/mT$ keep this information, and they carry digital
laws of the kind used throughout the book.

\section{The digital law of a plane group}\label{sec:crystlaw}

In a lattice basis every element of $G$ is an affine map $\{A\mid v_A+t\}\colon x\mapsto Ax+v_A+t$
with $A\in P\subset\GL(n,\Z)$, $t\in\Z^n$, and a fixed fractional translation $v_A$ for each $A$. For $m\ge1$ the
subgroup $mT$ is normal in $G$ and $|G/mT|=m^n|P|$. Numbering $P=\{A_0=1,A_1,\dots\}$, in the plane we label
$\{A_i\mid v_i+t\}$, $t\in(\Z/m)^2$, by
\begin{equation}\label{eq:crystlabel}
 k=i\,m^2+t_1m+t_2 ,
\end{equation}
and the product is
\begin{equation}\label{eq:crystlaw}
 \{A_i\mid v_i+t\}\{A_j\mid v_j+s\}=\{A_iA_j\mid v_{ij}+A_is+t+c(i,j)\},\qquad c(i,j)=A_iv_j+v_i-v_{ij}\in\Z^n,
\end{equation}
where $v_{ij}$ is the fractional translation of $A_iA_j$ and translations are reduced modulo $m$. The carry $c$ is a
$2$-cocycle representing the class of the extension $1\to T\to G\to P\to1$ in $H^2(P,T)$; it vanishes for the
symmorphic groups. If $|P|$ divides $m$, the exact sequence $H^2(P,T)\xrightarrow{\;m\;}H^2(P,T)\to H^2(P,T/mT)$
together with $|P|\,H^2(P,T)=0$ shows that the class survives in $H^2(P,T/mT)$, so the extension
$1\to T/mT\to G/mT\to P\to1$ determines $G$.

\section{A worked product: \texorpdfstring{$pm$ and $pg$}{pm and pg}}\label{sec:crystpmpg}

Both groups have the point group $\{1,A_1\}$, $A_1=\mathrm{diag}(-1,1)$,
the mirror $x\mapsto(-x_1,x_2)$; in $pg$ it carries $v_1=(0,\tfrac12)$, a glide. In $G/6T$ the label of
$\{A_1\mid v_1\}$ is $k=1\cdot 6^2=36$ by~\eqref{eq:crystlabel}. In $pm$ we have $c(1,1)=0$ and $36*36=0$, the
identity. In $pg$, $c(1,1)=A_1v_1+v_1=(0,1)$ and $36*36=1$, the label of the translation $\{1\mid(0,1)\}$: the
glide squares to a lattice translation, and the carry records exactly this. Already modulo $2T$ the two groups
differ: $pm/2T\cong C_2^3$, whereas $pg/2T\cong C_4\times C_2$.

\section{Separating the plane groups}\label{sec:crystplane}

Table~\ref{tab:wallpaperm} counts, for small $m$, the isomorphism types among the $17$
quotients $G/mT$, as detected by the distribution of element orders, the order of the centre, the number of
conjugacy classes and the order of the derived subgroup. For $m=6$ all $17$ are distinct: the factor $2$ detects
the glides, and the factor $3$ separates $p3m1$ from $p31m$, which differ only in the position of the mirrors
relative to the hexagonal lattice. For $m=1$ the quotients are the point groups.

\begin{table}[ht]\centering
\caption{Plane groups modulo $mT$: isomorphism types among the $17$ quotients $G/mT$.}\label{tab:wallpaperm}
\begin{tabular}{rrr}\hline
$m$ & $|G/mT|$ & types\\\hline
$1$ & $1$--$12$ & $9$\\
$2$ & $4$--$48$ & $13$\\
$3$ & $9$--$108$ & $10$\\
$4$ & $16$--$192$ & $16$\\
$6$ & $36$--$432$ & $17$\\
\hline
\end{tabular}
\end{table}

\section{Digital laws of the 230 space groups}\label{sec:crystlaw3}

The law~\eqref{eq:crystlaw} carries over to three dimensions without change. Let $\{R_i\mid t_i\}$, $0\le i<r$, be the
symmetry operations of the conventional cell, centring translations included, with translations written in units
of $\tfrac1{24}$ and $0\le t_i<24$, and let $L=\Z^3$ be the conventional lattice. The element $\{R_i\mid t_i+24s\}$ of
$G/mL$, $s\in(\Z/m)^3$, gets the label
\begin{equation}\label{eq:spacelabel}
 k=i\,m^3+s_1m^2+s_2m+s_3 ,
\end{equation}
and the product is
\begin{equation}\label{eq:spacelaw}
\begin{gathered}
 \{R_i\mid t_i+24s\}\{R_j\mid t_j+24s'\}=\{R_k\mid t_k+24(R_is'+s+c(i,j))\},\\
 c(i,j)=\tfrac1{24}\bigl(R_it_j+t_i-t_k\bigr)\in\Z^3,
\end{gathered}
\end{equation}
where $k$ is the operation with rotation $R_iR_j$ and translation congruent to $R_it_j+t_i$ modulo $24$. Thus every
one of the $230$ space groups, and every modulus $m$, gives a digital law on $\{0,\dots,m^3r-1\}$: the digit $i$
names a symmetry operation, the base-$m$ digits $s_1s_2s_3$ a translation, and the carry table $c$, of size
$r\times r$, is the whole crystallographic content of the law. In \texttt{multiply.py}, \texttt{inverse.py} and
\texttt{conjugate.py} these laws are called \texttt{SG$N$/$m$}; for example \texttt{SG62/2} is $Pnma$ modulo $2L$,
a group of order $64$, and \texttt{SG225/2} is $Fm\bar3m$ modulo $2L$, of order $1536$.

\paragraph{By hand: $P2_1/c$ (No.~14).} The operation lists are the triplets of the International Tables, kept in
\texttt{data/space\_\allowbreak{}group\_\allowbreak{}ops.txt}; nothing else is needed. For $P2_1/c$ (setting $P\,1\,2_1/c\,1$)
they are
\[
\begin{aligned}
 &0:\ x,y,z, &\qquad &1:\ -x,y+\tfrac12,-z+\tfrac12,\\
 &2:\ -x,-y,-z, &\qquad &3:\ x,-y+\tfrac12,z+\tfrac12 ,
\end{aligned}
\]
that is, $R_0=I$, $R_1=\operatorname{diag}(-1,1,-1)$, $R_2=-I$, $R_3=\operatorname{diag}(1,-1,1)$, $t_0=t_2=0$ and
$t_1=t_3=(0,12,12)$ in units of $\tfrac1{24}$. Composing the sixteen pairs gives the table of $(k,c(i,j))$:
\[
\setlength{\arraycolsep}{4pt}
\begin{array}{c|cccc}
 i\backslash j & 0 & 1 & 2 & 3\\\hline
 0 & 0,(0,0,0) & 1,(0,0,0) & 2,(0,0,0) & 3,(0,0,0)\\
 1 & 1,(0,0,0) & 0,(0,1,0) & 3,(0,0,0) & 2,(0,1,0)\\
 2 & 2,(0,0,0) & 3,(0,-1,-1) & 0,(0,0,0) & 1,(0,-1,-1)\\
 3 & 3,(0,0,0) & 2,(0,0,1) & 1,(0,0,0) & 0,(0,0,1)
\end{array}
\]
The entry $c(1,1)=(0,1,0)$ says that the screw $2_1$ squared is the translation by $b$; with all carries zero one
would have the symmorphic group $P2/m$. Modulo $2L$ there are $32$ labels $k=8i+4s_1+2s_2+s_3$. For instance $13$ is
$i=1$, $s=(1,0,1)$, the map $-x+1,\,y+\tfrac12,\,-z+\tfrac32$, and $27$ is $i=3$, $s=(0,1,1)$; then
$k=2$, $c(1,3)=(0,1,0)$, and $R_1(0,1,1)+(1,0,1)+(0,1,0)=(1,2,0)\equiv(1,0,0)$, so $13\cdot27=8\cdot2+4=20$
(\texttt{multiply.py SG14/2 13 27 --check}). Likewise $13\cdot13=2$, a translation by $b$, so $13$ has order $4$.

\section{Separating the space groups}\label{sec:crystspace}

For the $230$ space-group types, in the standard settings of the International Tables,
we use the lattice $L\subseteq T$ of the conventional cell, which is normal in $G$, and the quotients $G/mL$ of order
$m^3|G/L|$. For $g=\{R\mid t\}$ with $R$ of order $k$ we have $g^k=\{1\mid(1+R+\dots+R^{k-1})t\}$, so the order of
every element of $G/mL$ is obtained in closed form. So is its centralizer: $\{R\mid t\}$ and $\{S\mid u\}$ commute if
and only if $RS=SR$ and $(R-1)u\equiv(S-1)t$ modulo $mL$, and for fixed $g$ and $S$ the admissible $u$ form either
the empty set or a coset of the kernel of $R-1$ on $L/mL$. No multiplication table is needed, even for
$|G/12L|=331\,776$. Table~\ref{tab:spacem} counts the types separated by the distribution of element orders, by the
joint distribution of $(\operatorname{ord}g,|C(g)|)$ (both are isomorphism invariants of $G/mL$; the second contains
the order of the centre and the number of conjugacy classes), and by the latter refined by $(\det R,\operatorname{tr}R)$,
which also uses the action of the point group. Two space groups are isomorphic if and only if they are affinely
conjugate (Bieberbach), so the $11$ enantiomorphic pairs, such as $P3_121$ and $P3_221$, have isomorphic quotients
for every $m$: no invariant of the abstract group separates more than $219$ types, and the remaining distinction
requires an orientation of $L/mL$.

\begin{table}[ht]\centering
\caption{Space groups modulo $mL$: types separated among the $230$ (at most $219$ up to isomorphism).}\label{tab:spacem}
\begin{tabular}{rrrrr}\hline
$m$ & $\max|G/mL|$ & $\operatorname{ord}g$ & $(\operatorname{ord}g,|C(g)|)$ & with $(\det R,\operatorname{tr}R)$\\\hline
$1$ & $192$ & $37$ & $38$ & $73$\\
$2$ & $1536$ & $88$ & $110$ & $189$\\
$3$ & $5184$ & $71$ & $78$ & $82$\\
$4$ & $12288$ & $156$ & $198$ & $198$\\
$6$ & $41472$ & $171$ & $204$ & $204$\\
$12$ & $331776$ & $183$ & $213$ & $213$\\
\hline
\end{tabular}
\end{table}

For $m=12$ the joint distribution separates $213$ types. Besides the enantiomorphic pairs it leaves five sets,
which account for the remaining $6$: $\{Pca2_1,Pna2_1\}$ (Nos.~29, 33), $\{Pnc2,Pba2,Pnn2\}$ (30, 32, 34),
$\{Pnnn,Pban\}$ (48, 50), $\{P\bar4b2,P\bar4n2\}$ (117, 118) and $\{P6_3cm,P6_3mc\}$ (185, 186); the same sets are
unseparated for $m=4$ and $m=6$.
The full multiplication tables (\texttt{space\_\allowbreak{}residual.py}) separate most of them by finer
invariants: $Pca2_1$ and $Pna2_1$ already modulo $2L$, where the subgroup generated by the squares has order $4$,
resp.\ $8$; $Pnnn$ and $Pban$ modulo $4L$, where the third term of the lower central series has order $8$, resp.\ $4$;
$P\bar4b2$ and $P\bar4n2$ modulo $4L$ ($|\gamma_3|=8$, resp.\ $16$); $Pnc2$ from $Pba2$ and $Pnn2$ modulo $4L$
($|\gamma_3|=2$, resp.\ $4$); and $P6_3cm$ and $P6_3mc$ modulo $6L$, by the fibres of the squaring map. For $Pba2$ and $Pnn2$ (Nos.~32, 34) none of these invariants differs for $m\le8$, but an exhaustive
backtracking search over images of generators shows that their quotients modulo $4L$ are not isomorphic (for $m=2$ the search does find isomorphisms). Hence
the finite quotients $G/mL$ with $m\in\{2,4,6,12\}$ (modulo $12L$ for the bulk of the types, modulo $2L$, $4L$ and $6L$
for the residual sets) separate all $219$ isomorphism types of space groups, and only the orientation of the
enantiomorphic pairs is invisible to them. The last step, an exhaustive search in groups of order $256$, is confirmed independently by GAP
(\texttt{space\_\allowbreak{}gap.py}): modulo $2L$ the quotients of $Pnc2$, $Pba2$ and $Pnn2$ are all the small group
$[32,22]$, whereas modulo $4L$ they are $[256,2662]$, $[256,1216]$ and $[256,557]$; modulo $8L$, of order $2048$,
\texttt{IsomorphismGroups} also finds no isomorphism between the quotients of $Pba2$ and $Pnn2$.

The computations are made by \texttt{wallpaper\_\allowbreak{}law.py} and \texttt{space\_\allowbreak{}quotients.py}
(Chapter~\ref{ch:software}); the $17$ laws~\eqref{eq:crystlaw} are available in \texttt{multiply.py}, \texttt{inverse.py}
and \texttt{conjugate.py} under the names \texttt{pg/6}, \texttt{p4g/6}, \dots, with any modulus (\texttt{pg/4}, \dots).

\section{The quotients as abstract finite groups}\label{sec:crystid}

\paragraph{Extensions.} For every $m$ the quotient is an extension
\begin{equation}\label{eq:crystext}
 1\to T/mT\to G/mT\to P\to 1
\end{equation}
of the abelian group $T/mT\cong(\Z/m)^n$ by the point group $P$ (for a centred space group and the conventional
lattice $L$, of the group of all translations modulo $mL$ by $P$). If $G$ is symmorphic, \eqref{eq:crystext} splits and
$G/mT\cong(\Z/m)^n\rtimes P$, with $P$ acting through its integer matrices reduced modulo $m$. If $G$ is not
symmorphic and $|P|$ divides $m$, the extension does not split, since the class of $G$ is non-zero and survives in
$H^2(P,T/mT)$ (Section~\ref{sec:crystlaw}). For smaller $m$ the class may die: modulo $2T$ the glide group $pg$
gives $C_4\times C_2$, a non-split extension of $C_2^2$ by $C_2$, but $p3m1$ and $p31m$ both give $S_4$, and they
are told apart only modulo $6T$.

\paragraph{Dicyclic and metacyclic quotients.} None of these groups is dicyclic: a dicyclic group has a unique
involution, whereas every quotient in Table~\ref{tab:planeid} has at least $3$ involutions, and every quotient
$G/mL$ of a space group with $m$ even contains $L/2L$-translations forming $C_2^3$, hence at least $7$ involutions.
Metacyclic groups are rare too. Every subgroup of a metacyclic group is generated by two elements, so no quotient
$G/mL$ of a space group with $m\ge2$ is metacyclic: it contains $(\Z/m)^3$. In the plane, GAP finds exactly
$p1$, $pg$ and $cm$ metacyclic modulo $2T$, and $p1$ and $pg$ modulo $6T$; for example $pg/6T\cong C_6\times(C_3\rtimes C_4)$.

\paragraph{The tables.} Tables~\ref{tab:planeid} and~\ref{tab:spaceid} identify the quotients with GAP
(\texttt{space\_\allowbreak{}table.py}): $[N,i]$ is the group \texttt{SmallGroup($N$,$i$)} of the SmallGroups library, and
the structure is GAP's \texttt{StructureDescription}, in which $D_{2k}$ is the dihedral group of order $2k$, $A\rtimes B$ a
split extension and $A\cdot B$ a non-split one. A structure description names a group only up to such products, so
the library number is the precise identification; it is not available for the order $1536$ ($4$ face-centred cubic
types), marked ---. Each row is also a digital law: \texttt{multiply.py pg/2}, \texttt{multiply.py SG62/2}, and so on.

\begin{table}[ht]\centering\small
\caption{The plane groups $G$, their point groups $P$ and their quotients modulo $2T$ and $6T$.}\label{tab:planeid}
\resizebox{\textwidth}{!}{%
\begin{tabular}{llllll}\hline
$G$ & $P$ & $G/2T$ & structure & $G/6T$ & structure\\\hline
$p1$ & $1$ & $[4,2]$ & $C_{2}\times C_{2}$ & $[36,14]$ & $C_{6}\times C_{6}$\\
$p2$ & $C_{2}$ & $[8,5]$ & $C_{2}\times C_{2}\times C_{2}$ & $[72,49]$ & $C_{2}\times C_{2}\times ((C_{3}\times C_{3})\rtimes C_{2})$\\
$pm$ & $C_{2}$ & $[8,5]$ & $C_{2}\times C_{2}\times C_{2}$ & $[72,48]$ & $C_{2}\times C_{6}\times S_{3}$\\
$pg$ & $C_{2}$ & $[8,2]$ & $C_{4}\times C_{2}$ & $[72,29]$ & $C_{6}\times (C_{3}\rtimes C_{4})$\\
$cm$ & $C_{2}$ & $[8,3]$ & $D_{8}$ & $[72,30]$ & $C_{3}\times ((C_{6}\times C_{2})\rtimes C_{2})$\\
$pmm$ & $C_{2}\times C_{2}$ & $[16,14]$ & $C_{2}\times C_{2}\times C_{2}\times C_{2}$ & $[144,192]$ & $C_{2}\times C_{2}\times S_{3}\times S_{3}$\\
$pmg$ & $C_{2}\times C_{2}$ & $[16,11]$ & $C_{2}\times D_{8}$ & $[144,151]$ & $C_{2}\times ((C_{6}\times S_{3})\rtimes C_{2})$\\
$pgg$ & $C_{2}\times C_{2}$ & $[16,3]$ & $(C_{4}\times C_{2})\rtimes C_{2}$ & $[144,65]$ & $(C_{3}\times C_{3})\rtimes ((C_{4}\times C_{2})\rtimes C_{2})$\\
$cmm$ & $C_{2}\times C_{2}$ & $[16,11]$ & $C_{2}\times D_{8}$ & $[144,154]$ & $(C_{3}\times C_{3})\rtimes (C_{2}\times D_{8})$\\
$p4$ & $C_{4}$ & $[16,3]$ & $(C_{4}\times C_{2})\rtimes C_{2}$ & $[144,136]$ & $(C_{3}\times C_{3})\rtimes ((C_{4}\times C_{2})\rtimes C_{2})$\\
$p4m$ & $D_{8}$ & $[32,27]$ & $(C_{2}\times C_{2}\times C_{2}\times C_{2})\rtimes C_{2}$ & $[288,889]$ & $(C_{3}\times C_{3})\rtimes ((C_{2}\times C_{2}\times C_{2}\times C_{2})\rtimes C_{2})$\\
$p4g$ & $D_{8}$ & $[32,6]$ & $(C_{2}\times C_{2}\times C_{2})\rtimes C_{4}$ & $[288,386]$ & $(C_{3}\times C_{3})\rtimes ((C_{2}\times C_{2}\times C_{2})\rtimes C_{4})$\\
$p3$ & $C_{3}$ & $[12,3]$ & $A_{4}$ & $[108,22]$ & $(C_{2}\times C_{2})\rtimes ((C_{3}\times C_{3})\rtimes C_{3})$\\
$p3m1$ & $S_{3}$ & $[24,12]$ & $S_{4}$ & $[216,95]$ & $((C_{2}\times C_{2})\rtimes ((C_{3}\times C_{3})\rtimes C_{3}))\rtimes C_{2}$\\
$p31m$ & $S_{3}$ & $[24,12]$ & $S_{4}$ & $[216,92]$ & $((C_{3}\times A_{4})\rtimes C_{2})\rtimes C_{3}$\\
$p6$ & $C_{6}$ & $[24,13]$ & $C_{2}\times A_{4}$ & $[216,99]$ & $(C_{2}\times C_{2}\times ((C_{3}\times C_{3})\rtimes C_{2}))\rtimes C_{3}$\\
$p6m$ & $D_{12}$ & $[48,48]$ & $C_{2}\times S_{4}$ & $[432,523]$ & $((C_{2}\times C_{2}\times ((C_{3}\times C_{3})\rtimes C_{2}))\rtimes C_{3})\rtimes C_{2}$\\
\hline
\end{tabular}}
\end{table}

{\footnotesize
\setlength{\tabcolsep}{3pt}
\begin{longtable}{@{}rlllp{0.34\textwidth}l@{}}
\caption{The $230$ space groups $G$ (International Tables number and symbol), their point groups $P$ and their
quotients modulo $2L$, $L$ the conventional lattice.}\label{tab:spaceid}\\
\hline No. & $G$ & $P$ & $|G/2L|$ & structure of $G/2L$ & SmallGroup\\\hline\endfirsthead
\hline No. & $G$ & $P$ & $|G/2L|$ & structure of $G/2L$ & SmallGroup\\\hline\endhead
\hline\endfoot
1 & $P1$ & $1$ & $8$ & $C_{2}\times C_{2}\times C_{2}$ & $[8,5]$\\
2 & $P\bar{1}$ & $C_{2}$ & $16$ & $C_{2}\times C_{2}\times C_{2}\times C_{2}$ & $[16,14]$\\
3 & $P121$ & $C_{2}$ & $16$ & $C_{2}\times C_{2}\times C_{2}\times C_{2}$ & $[16,14]$\\
4 & $P12_{1}1$ & $C_{2}$ & $16$ & $C_{4}\times C_{2}\times C_{2}$ & $[16,10]$\\
5 & $C121$ & $C_{2}$ & $32$ & $C_{2}\times ((C_{4}\times C_{2})\rtimes C_{2})$ & $[32,22]$\\
6 & $P1m1$ & $C_{2}$ & $16$ & $C_{2}\times C_{2}\times C_{2}\times C_{2}$ & $[16,14]$\\
7 & $P1c1$ & $C_{2}$ & $16$ & $C_{4}\times C_{2}\times C_{2}$ & $[16,10]$\\
8 & $C1m1$ & $C_{2}$ & $32$ & $C_{2}\times ((C_{4}\times C_{2})\rtimes C_{2})$ & $[32,22]$\\
9 & $C1c1$ & $C_{2}$ & $32$ & $(C_{4}\times C_{2})\rtimes C_{4}$ & $[32,2]$\\
10 & $P12/m1$ & $C_{2}\times C_{2}$ & $32$ & $C_{2}\times C_{2}\times C_{2}\times C_{2}\times C_{2}$ & $[32,51]$\\
11 & $P12_{1}/m1$ & $C_{2}\times C_{2}$ & $32$ & $C_{2}\times C_{2}\times D_{8}$ & $[32,46]$\\
12 & $C12/m1$ & $C_{2}\times C_{2}$ & $64$ & $C_{2}\times ((C_{2}\times C_{2}\times C_{2}\times C_{2})\rtimes C_{2})$ & $[64,202]$\\
13 & $P12/c1$ & $C_{2}\times C_{2}$ & $32$ & $C_{2}\times C_{2}\times D_{8}$ & $[32,46]$\\
14 & $P12_{1}/c1$ & $C_{2}\times C_{2}$ & $32$ & $C_{2}\times ((C_{4}\times C_{2})\rtimes C_{2})$ & $[32,22]$\\
15 & $C12/c1$ & $C_{2}\times C_{2}$ & $64$ & $(C_{2}\times C_{2}\times D_{8})\rtimes C_{2}$ & $[64,73]$\\
16 & $P222$ & $C_{2}\times C_{2}$ & $32$ & $C_{2}\times C_{2}\times C_{2}\times C_{2}\times C_{2}$ & $[32,51]$\\
17 & $P222_{1}$ & $C_{2}\times C_{2}$ & $32$ & $C_{2}\times C_{2}\times D_{8}$ & $[32,46]$\\
18 & $P2_{1}2_{1}2$ & $C_{2}\times C_{2}$ & $32$ & $C_{2}\times ((C_{4}\times C_{2})\rtimes C_{2})$ & $[32,22]$\\
19 & $P2_{1}2_{1}2_{1}$ & $C_{2}\times C_{2}$ & $32$ & $(C_{4}\times C_{2})\rtimes C_{4}$ & $[32,2]$\\
20 & $C222_{1}$ & $C_{2}\times C_{2}$ & $64$ & $(C_{2}\times ((C_{4}\times C_{2})\rtimes C_{2}))\rtimes C_{2}$ & $[64,75]$\\
21 & $C222$ & $C_{2}\times C_{2}$ & $64$ & $C_{2}\times ((C_{2}\times C_{2}\times C_{2}\times C_{2})\rtimes C_{2})$ & $[64,202]$\\
22 & $F222$ & $C_{2}\times C_{2}$ & $128$ & $(C_{2}\times ((C_{2}\times C_{2}\times C_{2}\times C_{2})\rtimes C_{2}))\rtimes C_{2}$ & $[128,1411]$\\
23 & $I222$ & $C_{2}\times C_{2}$ & $64$ & $(C_{2}\times ((C_{4}\times C_{2})\rtimes C_{2}))\rtimes C_{2}$ & $[64,60]$\\
24 & $I2_{1}2_{1}2_{1}$ & $C_{2}\times C_{2}$ & $64$ & $(C_{2}\times ((C_{4}\times C_{2})\rtimes C_{2}))\rtimes C_{2}$ & $[64,73]$\\
25 & $Pmm2$ & $C_{2}\times C_{2}$ & $32$ & $C_{2}\times C_{2}\times C_{2}\times C_{2}\times C_{2}$ & $[32,51]$\\
26 & $Pmc2_{1}$ & $C_{2}\times C_{2}$ & $32$ & $C_{4}\times C_{2}\times C_{2}\times C_{2}$ & $[32,45]$\\
27 & $Pcc2$ & $C_{2}\times C_{2}$ & $32$ & $C_{4}\times C_{2}\times C_{2}\times C_{2}$ & $[32,45]$\\
28 & $Pma2$ & $C_{2}\times C_{2}$ & $32$ & $C_{2}\times C_{2}\times D_{8}$ & $[32,46]$\\
29 & $Pca2_{1}$ & $C_{2}\times C_{2}$ & $32$ & $C_{2}\times (C_{4}\rtimes C_{4})$ & $[32,23]$\\
30 & $Pnc2$ & $C_{2}\times C_{2}$ & $32$ & $C_{2}\times ((C_{4}\times C_{2})\rtimes C_{2})$ & $[32,22]$\\
31 & $Pmn2_{1}$ & $C_{2}\times C_{2}$ & $32$ & $C_{2}\times ((C_{4}\times C_{2})\rtimes C_{2})$ & $[32,22]$\\
32 & $Pba2$ & $C_{2}\times C_{2}$ & $32$ & $C_{2}\times ((C_{4}\times C_{2})\rtimes C_{2})$ & $[32,22]$\\
33 & $Pna2_{1}$ & $C_{2}\times C_{2}$ & $32$ & $(C_{4}\times C_{2})\rtimes C_{4}$ & $[32,2]$\\
34 & $Pnn2$ & $C_{2}\times C_{2}$ & $32$ & $C_{2}\times ((C_{4}\times C_{2})\rtimes C_{2})$ & $[32,22]$\\
35 & $Cmm2$ & $C_{2}\times C_{2}$ & $64$ & $C_{2}\times ((C_{2}\times C_{2}\times C_{2}\times C_{2})\rtimes C_{2})$ & $[64,202]$\\
36 & $Cmc2_{1}$ & $C_{2}\times C_{2}$ & $64$ & $(C_{2}\times (C_{4}\rtimes C_{4}))\rtimes C_{2}$ & $[64,66]$\\
37 & $Ccc2$ & $C_{2}\times C_{2}$ & $64$ & $((C_{4}\times C_{2})\rtimes C_{4})\rtimes C_{2}$ & $[64,67]$\\
38 & $Amm2$ & $C_{2}\times C_{2}$ & $64$ & $C_{2}\times C_{2}\times ((C_{4}\times C_{2})\rtimes C_{2})$ & $[64,193]$\\
39 & $Abm2$ & $C_{2}\times C_{2}$ & $64$ & $C_{2}\times C_{4}\times D_{8}$ & $[64,196]$\\
40 & $Ama2$ & $C_{2}\times C_{2}$ & $64$ & $(C_{4}\times C_{2}\times C_{2}\times C_{2})\rtimes C_{2}$ & $[64,67]$\\
41 & $Aba2$ & $C_{2}\times C_{2}$ & $64$ & $(C_{4}\times C_{4}\times C_{2})\rtimes C_{2}$ & $[64,69]$\\
42 & $Fmm2$ & $C_{2}\times C_{2}$ & $128$ & $(C_{2}\times C_{2}\times ((C_{4}\times C_{2})\rtimes C_{2}))\rtimes C_{2}$ & $[128,1090]$\\
43 & $Fdd2$ & $C_{2}\times C_{2}$ & $128$ & $((C_{4}\times C_{2})\rtimes C_{8})\rtimes C_{2}$ & $[128,2]$\\
44 & $Imm2$ & $C_{2}\times C_{2}$ & $64$ & $(C_{2}\times ((C_{4}\times C_{2})\rtimes C_{2}))\rtimes C_{2}$ & $[64,60]$\\
45 & $Iba2$ & $C_{2}\times C_{2}$ & $64$ & $((C_{4}\times C_{2})\rtimes C_{4})\rtimes C_{2}$ & $[64,66]$\\
46 & $Ima2$ & $C_{2}\times C_{2}$ & $64$ & $(C_{2}\times ((C_{4}\times C_{2})\rtimes C_{2}))\rtimes C_{2}$ & $[64,67]$\\
47 & $Pmmm$ & $C_{2}\times C_{2}\times C_{2}$ & $64$ & $C_{2}\times C_{2}\times C_{2}\times C_{2}\times C_{2}\times C_{2}$ & $[64,267]$\\
48 & $Pnnn$ & $C_{2}\times C_{2}\times C_{2}$ & $64$ & $C_{2}\times ((C_{2}\times C_{2}\times C_{2}\times C_{2})\rtimes C_{2})$ & $[64,202]$\\
49 & $Pccm$ & $C_{2}\times C_{2}\times C_{2}$ & $64$ & $C_{2}\times C_{2}\times C_{2}\times D_{8}$ & $[64,261]$\\
50 & $Pban$ & $C_{2}\times C_{2}\times C_{2}$ & $64$ & $C_{2}\times ((C_{2}\times C_{2}\times C_{2}\times C_{2})\rtimes C_{2})$ & $[64,202]$\\
51 & $Pmma$ & $C_{2}\times C_{2}\times C_{2}$ & $64$ & $C_{2}\times C_{2}\times C_{2}\times D_{8}$ & $[64,261]$\\
52 & $Pnna$ & $C_{2}\times C_{2}\times C_{2}$ & $64$ & $(C_{2}\times C_{2}\times D_{8})\rtimes C_{2}$ & $[64,73]$\\
53 & $Pmna$ & $C_{2}\times C_{2}\times C_{2}$ & $64$ & $C_{2}\times ((C_{2}\times C_{2}\times C_{2}\times C_{2})\rtimes C_{2})$ & $[64,202]$\\
54 & $Pcca$ & $C_{2}\times C_{2}\times C_{2}$ & $64$ & $C_{2}\times ((C_{4}\times C_{2}\times C_{2})\rtimes C_{2})$ & $[64,203]$\\
55 & $Pbam$ & $C_{2}\times C_{2}\times C_{2}$ & $64$ & $C_{2}\times C_{2}\times ((C_{4}\times C_{2})\rtimes C_{2})$ & $[64,193]$\\
56 & $Pccn$ & $C_{2}\times C_{2}\times C_{2}$ & $64$ & $(C_{2}\times ((C_{4}\times C_{2})\rtimes C_{2}))\rtimes C_{2}$ & $[64,67]$\\
57 & $Pbcm$ & $C_{2}\times C_{2}\times C_{2}$ & $64$ & $C_{2}\times ((C_{4}\times C_{2}\times C_{2})\rtimes C_{2})$ & $[64,203]$\\
58 & $Pnnm$ & $C_{2}\times C_{2}\times C_{2}$ & $64$ & $(C_{2}\times ((C_{4}\times C_{2})\rtimes C_{2}))\rtimes C_{2}$ & $[64,60]$\\
59 & $Pmmn$ & $C_{2}\times C_{2}\times C_{2}$ & $64$ & $C_{2}\times ((C_{2}\times C_{2}\times C_{2}\times C_{2})\rtimes C_{2})$ & $[64,202]$\\
60 & $Pbcn$ & $C_{2}\times C_{2}\times C_{2}$ & $64$ & $(C_{2}\times ((C_{4}\times C_{2})\rtimes C_{2}))\rtimes C_{2}$ & $[64,75]$\\
61 & $Pbca$ & $C_{2}\times C_{2}\times C_{2}$ & $64$ & $((C_{4}\times C_{2})\rtimes C_{4})\rtimes C_{2}$ & $[64,77]$\\
62 & $Pnma$ & $C_{2}\times C_{2}\times C_{2}$ & $64$ & $(C_{2}\times ((C_{4}\times C_{2})\rtimes C_{2}))\rtimes C_{2}$ & $[64,67]$\\
63 & $Cmcm$ & $C_{2}\times C_{2}\times C_{2}$ & $128$ & $(C_{2}\times ((C_{2}\times C_{2}\times C_{2}\times C_{2})\rtimes C_{2}))\rtimes C_{2}$ & $[128,1142]$\\
64 & $Cmca$ & $C_{2}\times C_{2}\times C_{2}$ & $128$ & $(C_{2}\times ((C_{2}\times C_{2}\times C_{2}\times C_{2})\rtimes C_{2}))\rtimes C_{2}$ & $[128,1156]$\\
65 & $Cmmm$ & $C_{2}\times C_{2}\times C_{2}$ & $128$ & $C_{2}\times C_{2}\times ((C_{2}\times C_{2}\times C_{2}\times C_{2})\rtimes C_{2})$ & $[128,2163]$\\
66 & $Cccm$ & $C_{2}\times C_{2}\times C_{2}$ & $128$ & $((C_{2}\times C_{2}\times D_{8})\rtimes C_{2})\rtimes C_{2}$ & $[128,1135]$\\
67 & $Cmma$ & $C_{2}\times C_{2}\times C_{2}$ & $128$ & $C_{2}\times D_{8}\times D_{8}$ & $[128,2194]$\\
68 & $Ccca$ & $C_{2}\times C_{2}\times C_{2}$ & $128$ & $((C_{2}\times C_{2}\times D_{8})\rtimes C_{2})\rtimes C_{2}$ & $[128,1401]$\\
69 & $Fmmm$ & $C_{2}\times C_{2}\times C_{2}$ & $256$ & $(C_{2}\times C_{2}\times ((C_{2}\times C_{2}\times C_{2}\times C_{2})\rtimes C_{2}))\rtimes C_{2}$ & $[256,51978]$\\
70 & $Fddd$ & $C_{2}\times C_{2}\times C_{2}$ & $256$ & $((C_{2}\times ((C_{4}\times C_{4})\rtimes C_{2}))\rtimes C_{2})\rtimes C_{2}$ & $[256,722]$\\
71 & $Immm$ & $C_{2}\times C_{2}\times C_{2}$ & $128$ & $(C_{2}\times ((C_{2}\times C_{2}\times C_{2}\times C_{2})\rtimes C_{2}))\rtimes C_{2}$ & $[128,1578]$\\
72 & $Ibam$ & $C_{2}\times C_{2}\times C_{2}$ & $128$ & $((C_{2}\times C_{2}\times D_{8})\rtimes C_{2})\rtimes C_{2}$ & $[128,1142]$\\
73 & $Ibca$ & $C_{2}\times C_{2}\times C_{2}$ & $128$ & $((C_{2}\times C_{2}\times D_{8})\rtimes C_{2})\rtimes C_{2}$ & $[128,1407]$\\
74 & $Imma$ & $C_{2}\times C_{2}\times C_{2}$ & $128$ & $(C_{2}\times ((C_{2}\times C_{2}\times C_{2}\times C_{2})\rtimes C_{2}))\rtimes C_{2}$ & $[128,1135]$\\
75 & $P4$ & $C_{4}$ & $32$ & $C_{2}\times ((C_{4}\times C_{2})\rtimes C_{2})$ & $[32,22]$\\
76 & $P4_{1}$ & $C_{4}$ & $32$ & $(C_{8}\times C_{2})\rtimes C_{2}$ & $[32,5]$\\
77 & $P4_{2}$ & $C_{4}$ & $32$ & $C_{2}\times ((C_{4}\times C_{2})\rtimes C_{2})$ & $[32,22]$\\
78 & $P4_{3}$ & $C_{4}$ & $32$ & $(C_{8}\times C_{2})\rtimes C_{2}$ & $[32,5]$\\
79 & $I4$ & $C_{4}$ & $64$ & $(C_{4}\times C_{2}\times C_{2})\rtimes C_{4}$ & $[64,23]$\\
80 & $I4_{1}$ & $C_{4}$ & $64$ & $((C_{8}\times C_{2})\rtimes C_{2})\rtimes C_{2}$ & $[64,4]$\\
81 & $P\bar{4}$ & $C_{4}$ & $32$ & $C_{2}\times ((C_{4}\times C_{2})\rtimes C_{2})$ & $[32,22]$\\
82 & $I\bar{4}$ & $C_{4}$ & $64$ & $(C_{4}\times C_{2}\times C_{2})\rtimes C_{4}$ & $[64,23]$\\
83 & $P4/m$ & $C_{4}\times C_{2}$ & $64$ & $C_{2}\times C_{2}\times ((C_{4}\times C_{2})\rtimes C_{2})$ & $[64,193]$\\
84 & $P4_{2}/m$ & $C_{4}\times C_{2}$ & $64$ & $(C_{2}\times ((C_{4}\times C_{2})\rtimes C_{2}))\rtimes C_{2}$ & $[64,60]$\\
85 & $P4/n$ & $C_{4}\times C_{2}$ & $64$ & $C_{2}\times ((C_{2}\times C_{2}\times C_{2})\rtimes C_{4})$ & $[64,90]$\\
86 & $P4_{2}/n$ & $C_{4}\times C_{2}$ & $64$ & $C_{2}\times ((C_{2}\times C_{2}\times C_{2})\rtimes C_{4})$ & $[64,90]$\\
87 & $I4/m$ & $C_{4}\times C_{2}$ & $128$ & $((C_{4}\times C_{2}\times C_{2})\rtimes C_{4})\rtimes C_{2}$ & $[128,630]$\\
88 & $I4_{1}/a$ & $C_{4}\times C_{2}$ & $128$ & $(((C_{8}\times C_{2})\rtimes C_{2})\rtimes C_{2})\rtimes C_{2}$ & $[128,75]$\\
89 & $P422$ & $D_{8}$ & $64$ & $C_{2}\times ((C_{2}\times C_{2}\times C_{2}\times C_{2})\rtimes C_{2})$ & $[64,202]$\\
90 & $P42_{1}2$ & $D_{8}$ & $64$ & $C_{2}\times ((C_{2}\times C_{2}\times C_{2})\rtimes C_{4})$ & $[64,90]$\\
91 & $P4_{1}22$ & $D_{8}$ & $64$ & $(C_{2}\times C_{2}\times D_{8})\rtimes C_{2}$ & $[64,128]$\\
92 & $P4_{1}2_{1}2$ & $D_{8}$ & $64$ & $((C_{8}\times C_{2})\rtimes C_{2})\rtimes C_{2}$ & $[64,8]$\\
93 & $P4_{2}22$ & $D_{8}$ & $64$ & $C_{2}\times ((C_{2}\times C_{2}\times C_{2}\times C_{2})\rtimes C_{2})$ & $[64,202]$\\
94 & $P4_{2}2_{1}2$ & $D_{8}$ & $64$ & $C_{2}\times ((C_{2}\times C_{2}\times C_{2})\rtimes C_{4})$ & $[64,90]$\\
95 & $P4_{3}22$ & $D_{8}$ & $64$ & $(C_{2}\times C_{2}\times D_{8})\rtimes C_{2}$ & $[64,128]$\\
96 & $P4_{3}2_{1}2$ & $D_{8}$ & $64$ & $((C_{8}\times C_{2})\rtimes C_{2})\rtimes C_{2}$ & $[64,8]$\\
97 & $I422$ & $D_{8}$ & $128$ & $((C_{2}\times ((C_{4}\times C_{2})\rtimes C_{2}))\rtimes C_{2})\rtimes C_{2}$ & $[128,753]$\\
98 & $I4_{1}22$ & $D_{8}$ & $128$ & $((C_{2}\times ((C_{4}\times C_{2})\rtimes C_{2}))\rtimes C_{2})\rtimes C_{2}$ & $[128,327]$\\
99 & $P4mm$ & $D_{8}$ & $64$ & $C_{2}\times ((C_{2}\times C_{2}\times C_{2}\times C_{2})\rtimes C_{2})$ & $[64,202]$\\
100 & $P4bm$ & $D_{8}$ & $64$ & $C_{2}\times ((C_{2}\times C_{2}\times C_{2})\rtimes C_{4})$ & $[64,90]$\\
101 & $P4_{2}cm$ & $D_{8}$ & $64$ & $(C_{4}\times C_{2}\times C_{2}\times C_{2})\rtimes C_{2}$ & $[64,67]$\\
102 & $P4_{2}nm$ & $D_{8}$ & $64$ & $C_{2}\times ((C_{2}\times C_{2}\times C_{2})\rtimes C_{4})$ & $[64,90]$\\
103 & $P4cc$ & $D_{8}$ & $64$ & $((C_{4}\times C_{2})\rtimes C_{4})\rtimes C_{2}$ & $[64,66]$\\
104 & $P4nc$ & $D_{8}$ & $64$ & $((C_{4}\times C_{2})\rtimes C_{2})\rtimes C_{4}$ & $[64,23]$\\
105 & $P4_{2}mc$ & $D_{8}$ & $64$ & $(C_{2}\times ((C_{4}\times C_{2})\rtimes C_{2}))\rtimes C_{2}$ & $[64,60]$\\
106 & $P4_{2}bc$ & $D_{8}$ & $64$ & $((C_{4}\times C_{2})\rtimes C_{2})\rtimes C_{4}$ & $[64,23]$\\
107 & $I4mm$ & $D_{8}$ & $128$ & $((C_{2}\times ((C_{4}\times C_{2})\rtimes C_{2}))\rtimes C_{2})\rtimes C_{2}$ & $[128,621]$\\
108 & $I4cm$ & $D_{8}$ & $128$ & $(((C_{4}\times C_{2})\rtimes C_{4})\rtimes C_{2})\rtimes C_{2}$ & $[128,599]$\\
109 & $I4_{1}md$ & $D_{8}$ & $128$ & $(((C_{8}\times C_{2})\rtimes C_{2})\rtimes C_{2})\rtimes C_{2}$ & $[128,48]$\\
110 & $I4_{1}cd$ & $D_{8}$ & $128$ & $((C_{4}\times C_{2})\rtimes C_{2})\rtimes C_{8}$ & $[128,58]$\\
111 & $P\bar{4}2m$ & $D_{8}$ & $64$ & $C_{2}\times ((C_{2}\times C_{2}\times C_{2}\times C_{2})\rtimes C_{2})$ & $[64,202]$\\
112 & $P\bar{4}2c$ & $D_{8}$ & $64$ & $(C_{2}\times ((C_{4}\times C_{2})\rtimes C_{2}))\rtimes C_{2}$ & $[64,60]$\\
113 & $P\bar{4}2_{1}m$ & $D_{8}$ & $64$ & $C_{2}\times ((C_{2}\times C_{2}\times C_{2})\rtimes C_{4})$ & $[64,90]$\\
114 & $P\bar{4}2_{1}c$ & $D_{8}$ & $64$ & $((C_{4}\times C_{2})\rtimes C_{2})\rtimes C_{4}$ & $[64,23]$\\
115 & $P\bar{4}m2$ & $D_{8}$ & $64$ & $C_{2}\times ((C_{2}\times C_{2}\times C_{2}\times C_{2})\rtimes C_{2})$ & $[64,202]$\\
116 & $P\bar{4}c2$ & $D_{8}$ & $64$ & $(C_{4}\times C_{2}\times C_{2}\times C_{2})\rtimes C_{2}$ & $[64,67]$\\
117 & $P\bar{4}b2$ & $D_{8}$ & $64$ & $C_{2}\times ((C_{2}\times C_{2}\times C_{2})\rtimes C_{4})$ & $[64,90]$\\
118 & $P\bar{4}n2$ & $D_{8}$ & $64$ & $C_{2}\times ((C_{2}\times C_{2}\times C_{2})\rtimes C_{4})$ & $[64,90]$\\
119 & $I\bar{4}m2$ & $D_{8}$ & $128$ & $((C_{2}\times ((C_{4}\times C_{2})\rtimes C_{2}))\rtimes C_{2})\rtimes C_{2}$ & $[128,753]$\\
120 & $I\bar{4}c2$ & $D_{8}$ & $128$ & $(((C_{4}\times C_{2})\rtimes C_{4})\rtimes C_{2})\rtimes C_{2}$ & $[128,776]$\\
121 & $I\bar{4}2m$ & $D_{8}$ & $128$ & $((C_{2}\times ((C_{4}\times C_{2})\rtimes C_{2}))\rtimes C_{2})\rtimes C_{2}$ & $[128,621]$\\
122 & $I\bar{4}2d$ & $D_{8}$ & $128$ & $((C_{4}\times C_{2}\times C_{2})\rtimes C_{4})\rtimes C_{2}$ & $[128,75]$\\
123 & $P4/mmm$ & $C_{2}\times D_{8}$ & $128$ & $C_{2}\times C_{2}\times ((C_{2}\times C_{2}\times C_{2}\times C_{2})\rtimes C_{2})$ & $[128,2163]$\\
124 & $P4/mcc$ & $C_{2}\times D_{8}$ & $128$ & $(C_{2}\times C_{2}\times C_{2}\times D_{8})\rtimes C_{2}$ & $[128,1142]$\\
125 & $P4/nbm$ & $C_{2}\times D_{8}$ & $128$ & $C_{2}\times (((C_{2}\times C_{2}\times C_{2})\rtimes C_{4})\rtimes C_{2})$ & $[128,1755]$\\
126 & $P4/nnc$ & $C_{2}\times D_{8}$ & $128$ & $(C_{2}\times ((C_{2}\times C_{2}\times C_{2})\rtimes C_{4}))\rtimes C_{2}$ & $[128,753]$\\
127 & $P4/mbm$ & $C_{2}\times D_{8}$ & $128$ & $C_{2}\times C_{2}\times ((C_{2}\times C_{2}\times C_{2})\rtimes C_{4})$ & $[128,1613]$\\
128 & $P4/mnc$ & $C_{2}\times D_{8}$ & $128$ & $(C_{2}\times C_{2}\times ((C_{4}\times C_{2})\rtimes C_{2}))\rtimes C_{2}$ & $[128,630]$\\
129 & $P4/nmm$ & $C_{2}\times D_{8}$ & $128$ & $C_{2}\times (((C_{2}\times C_{2}\times C_{2})\rtimes C_{4})\rtimes C_{2})$ & $[128,1755]$\\
130 & $P4/ncc$ & $C_{2}\times D_{8}$ & $128$ & $(C_{2}\times ((C_{2}\times C_{2}\times C_{2})\rtimes C_{4}))\rtimes C_{2}$ & $[128,776]$\\
131 & $P4_{2}/mmc$ & $C_{2}\times D_{8}$ & $128$ & $(C_{2}\times C_{2}\times C_{2}\times C_{2}\times C_{2}\times C_{2})\rtimes C_{2}$ & $[128,1578]$\\
132 & $P4_{2}/mcm$ & $C_{2}\times D_{8}$ & $128$ & $(C_{2}\times C_{2}\times C_{2}\times D_{8})\rtimes C_{2}$ & $[128,1135]$\\
133 & $P4_{2}/nbc$ & $C_{2}\times D_{8}$ & $128$ & $(C_{2}\times ((C_{2}\times C_{2}\times C_{2})\rtimes C_{4}))\rtimes C_{2}$ & $[128,753]$\\
134 & $P4_{2}/nnm$ & $C_{2}\times D_{8}$ & $128$ & $C_{2}\times (((C_{2}\times C_{2}\times C_{2})\rtimes C_{4})\rtimes C_{2})$ & $[128,1755]$\\
135 & $P4_{2}/mbc$ & $C_{2}\times D_{8}$ & $128$ & $((C_{2}\times ((C_{4}\times C_{2})\rtimes C_{2}))\rtimes C_{2})\rtimes C_{2}$ & $[128,630]$\\
136 & $P4_{2}/mnm$ & $C_{2}\times D_{8}$ & $128$ & $((C_{2}\times ((C_{4}\times C_{2})\rtimes C_{2}))\rtimes C_{2})\rtimes C_{2}$ & $[128,513]$\\
137 & $P4_{2}/nmc$ & $C_{2}\times D_{8}$ & $128$ & $(C_{2}\times ((C_{2}\times C_{2}\times C_{2})\rtimes C_{4}))\rtimes C_{2}$ & $[128,753]$\\
138 & $P4_{2}/ncm$ & $C_{2}\times D_{8}$ & $128$ & $(C_{2}\times ((C_{2}\times C_{2}\times C_{2})\rtimes C_{4}))\rtimes C_{2}$ & $[128,522]$\\
139 & $I4/mmm$ & $C_{2}\times D_{8}$ & $256$ & $((C_{2}\times ((C_{2}\times C_{2}\times C_{2}\times C_{2})\rtimes C_{2}))\rtimes C_{2})\rtimes C_{2}$ & $[256,16888]$\\
140 & $I4/mcm$ & $C_{2}\times D_{8}$ & $256$ & $(((C_{2}\times C_{2}\times D_{8})\rtimes C_{2})\rtimes C_{2})\rtimes C_{2}$ & $[256,17392]$\\
141 & $I4_{1}/amd$ & $C_{2}\times D_{8}$ & $256$ & $((((C_{8}\times C_{2})\rtimes C_{2})\rtimes C_{2})\rtimes C_{2})\rtimes C_{2}$ & $[256,5078]$\\
142 & $I4_{1}/acd$ & $C_{2}\times D_{8}$ & $256$ & $((((C_{8}\times C_{2})\rtimes C_{2})\rtimes C_{2})\rtimes C_{2})\rtimes C_{2}$ & $[256,5136]$\\
143 & $P3$ & $C_{3}$ & $24$ & $C_{2}\times A_{4}$ & $[24,13]$\\
144 & $P3_{1}$ & $C_{3}$ & $24$ & $C_{2}\times A_{4}$ & $[24,13]$\\
145 & $P3_{2}$ & $C_{3}$ & $24$ & $C_{2}\times A_{4}$ & $[24,13]$\\
146 & $R3$ & $C_{3}$ & $72$ & $C_{6}\times A_{4}$ & $[72,47]$\\
147 & $P\bar{3}$ & $C_{6}$ & $48$ & $C_{2}\times C_{2}\times A_{4}$ & $[48,49]$\\
148 & $R\bar{3}$ & $C_{6}$ & $144$ & $C_{2}\times S_{3}\times A_{4}$ & $[144,190]$\\
149 & $P312$ & $S_{3}$ & $48$ & $C_{2}\times S_{4}$ & $[48,48]$\\
150 & $P321$ & $S_{3}$ & $48$ & $C_{2}\times S_{4}$ & $[48,48]$\\
151 & $P3_{1}12$ & $S_{3}$ & $48$ & $C_{2}\times S_{4}$ & $[48,48]$\\
152 & $P3_{1}21$ & $S_{3}$ & $48$ & $C_{2}\times S_{4}$ & $[48,48]$\\
153 & $P3_{2}12$ & $S_{3}$ & $48$ & $C_{2}\times S_{4}$ & $[48,48]$\\
154 & $P3_{2}21$ & $S_{3}$ & $48$ & $C_{2}\times S_{4}$ & $[48,48]$\\
155 & $R32$ & $S_{3}$ & $144$ & $C_{2}\times ((C_{3}\times A_{4})\rtimes C_{2})$ & $[144,189]$\\
156 & $P3m1$ & $S_{3}$ & $48$ & $C_{2}\times S_{4}$ & $[48,48]$\\
157 & $P31m$ & $S_{3}$ & $48$ & $C_{2}\times S_{4}$ & $[48,48]$\\
158 & $P3c1$ & $S_{3}$ & $48$ & $A_{4}\rtimes C_{4}$ & $[48,30]$\\
159 & $P31c$ & $S_{3}$ & $48$ & $A_{4}\rtimes C_{4}$ & $[48,30]$\\
160 & $R3m$ & $S_{3}$ & $144$ & $C_{6}\times S_{4}$ & $[144,188]$\\
161 & $R3c$ & $S_{3}$ & $144$ & $C_{3}\times (A_{4}\rtimes C_{4})$ & $[144,123]$\\
162 & $P\bar{3}1m$ & $D_{12}$ & $96$ & $C_{2}\times C_{2}\times S_{4}$ & $[96,226]$\\
163 & $P\bar{3}1c$ & $D_{12}$ & $96$ & $(C_{2}\times S_{4})\rtimes C_{2}$ & $[96,195]$\\
164 & $P\bar{3}m1$ & $D_{12}$ & $96$ & $C_{2}\times C_{2}\times S_{4}$ & $[96,226]$\\
165 & $P\bar{3}c1$ & $D_{12}$ & $96$ & $(C_{2}\times S_{4})\rtimes C_{2}$ & $[96,195]$\\
166 & $R\bar{3}m$ & $D_{12}$ & $288$ & $C_{2}\times S_{3}\times S_{4}$ & $[288,1028]$\\
167 & $R\bar{3}c$ & $D_{12}$ & $288$ & $(C_{2}\times ((C_{3}\times A_{4})\rtimes C_{2}))\rtimes C_{2}$ & $[288,858]$\\
168 & $P6$ & $C_{6}$ & $48$ & $C_{2}\times C_{2}\times A_{4}$ & $[48,49]$\\
169 & $P6_{1}$ & $C_{6}$ & $48$ & $C_{4}\times A_{4}$ & $[48,31]$\\
170 & $P6_{5}$ & $C_{6}$ & $48$ & $C_{4}\times A_{4}$ & $[48,31]$\\
171 & $P6_{2}$ & $C_{6}$ & $48$ & $C_{2}\times C_{2}\times A_{4}$ & $[48,49]$\\
172 & $P6_{4}$ & $C_{6}$ & $48$ & $C_{2}\times C_{2}\times A_{4}$ & $[48,49]$\\
173 & $P6_{3}$ & $C_{6}$ & $48$ & $C_{4}\times A_{4}$ & $[48,31]$\\
174 & $P\bar{6}$ & $C_{6}$ & $48$ & $C_{2}\times C_{2}\times A_{4}$ & $[48,49]$\\
175 & $P6/m$ & $C_{6}\times C_{2}$ & $96$ & $C_{2}\times C_{2}\times C_{2}\times A_{4}$ & $[96,228]$\\
176 & $P6_{3}/m$ & $C_{6}\times C_{2}$ & $96$ & $A_{4}\times D_{8}$ & $[96,197]$\\
177 & $P622$ & $D_{12}$ & $96$ & $C_{2}\times C_{2}\times S_{4}$ & $[96,226]$\\
178 & $P6_{1}22$ & $D_{12}$ & $96$ & $(C_{4}\times A_{4})\rtimes C_{2}$ & $[96,187]$\\
179 & $P6_{5}22$ & $D_{12}$ & $96$ & $(C_{4}\times A_{4})\rtimes C_{2}$ & $[96,187]$\\
180 & $P6_{2}22$ & $D_{12}$ & $96$ & $C_{2}\times C_{2}\times S_{4}$ & $[96,226]$\\
181 & $P6_{4}22$ & $D_{12}$ & $96$ & $C_{2}\times C_{2}\times S_{4}$ & $[96,226]$\\
182 & $P6_{3}22$ & $D_{12}$ & $96$ & $(C_{4}\times A_{4})\rtimes C_{2}$ & $[96,187]$\\
183 & $P6mm$ & $D_{12}$ & $96$ & $C_{2}\times C_{2}\times S_{4}$ & $[96,226]$\\
184 & $P6cc$ & $D_{12}$ & $96$ & $C_{2}\times (A_{4}\rtimes C_{4})$ & $[96,194]$\\
185 & $P6_{3}cm$ & $D_{12}$ & $96$ & $C_{4}\times S_{4}$ & $[96,186]$\\
186 & $P6_{3}mc$ & $D_{12}$ & $96$ & $C_{4}\times S_{4}$ & $[96,186]$\\
187 & $P\bar{6}m2$ & $D_{12}$ & $96$ & $C_{2}\times C_{2}\times S_{4}$ & $[96,226]$\\
188 & $P\bar{6}c2$ & $D_{12}$ & $96$ & $(C_{2}\times C_{2}\times A_{4})\rtimes C_{2}$ & $[96,195]$\\
189 & $P\bar{6}2m$ & $D_{12}$ & $96$ & $C_{2}\times C_{2}\times S_{4}$ & $[96,226]$\\
190 & $P\bar{6}2c$ & $D_{12}$ & $96$ & $(C_{2}\times C_{2}\times A_{4})\rtimes C_{2}$ & $[96,195]$\\
191 & $P6/mmm$ & $C_{2}\times C_{2}\times S_{3}$ & $192$ & $C_{2}\times C_{2}\times C_{2}\times S_{4}$ & $[192,1537]$\\
192 & $P6/mcc$ & $C_{2}\times C_{2}\times S_{3}$ & $192$ & $C_{2}\times ((C_{2}\times S_{4})\rtimes C_{2})$ & $[192,1488]$\\
193 & $P6_{3}/mcm$ & $C_{2}\times C_{2}\times S_{3}$ & $192$ & $S_{4}\times D_{8}$ & $[192,1472]$\\
194 & $P6_{3}/mmc$ & $C_{2}\times C_{2}\times S_{3}$ & $192$ & $S_{4}\times D_{8}$ & $[192,1472]$\\
195 & $P23$ & $A_{4}$ & $96$ & $C_{2}\times ((C_{2}\times C_{2}\times C_{2}\times C_{2})\rtimes C_{3})$ & $[96,229]$\\
196 & $F23$ & $A_{4}$ & $384$ & $((C_{2}\times ((C_{2}\times C_{2}\times C_{2}\times C_{2})\rtimes C_{2}))\rtimes C_{2})\rtimes C_{3}$ & $[384,5863]$\\
197 & $I23$ & $A_{4}$ & $192$ & $((C_{2}\times ((C_{4}\times C_{2})\rtimes C_{2}))\rtimes C_{2})\rtimes C_{3}$ & $[192,191]$\\
198 & $P2_{1}3$ & $A_{4}$ & $96$ & $((C_{4}\times C_{2})\rtimes C_{4})\rtimes C_{3}$ & $[96,3]$\\
199 & $I2_{1}3$ & $A_{4}$ & $192$ & $((C_{2}\times ((C_{4}\times C_{2})\rtimes C_{2}))\rtimes C_{2})\rtimes C_{3}$ & $[192,194]$\\
200 & $Pm\bar{3}$ & $C_{2}\times A_{4}$ & $192$ & $C_{2}\times C_{2}\times ((C_{2}\times C_{2}\times C_{2}\times C_{2})\rtimes C_{3})$ & $[192,1540]$\\
201 & $Pn\bar{3}$ & $C_{2}\times A_{4}$ & $192$ & $C_{2}\times (((C_{2}\times C_{2}\times C_{2}\times C_{2})\rtimes C_{2})\rtimes C_{3})$ & $[192,1000]$\\
202 & $Fm\bar{3}$ & $C_{2}\times A_{4}$ & $768$ & $((C_{2}\times C_{2}\times ((C_{2}\times C_{2}\times C_{2}\times C_{2})\rtimes C_{2}))\rtimes C_{2})\rtimes C_{3}$ & $[768,1084555]$\\
203 & $Fd\bar{3}$ & $C_{2}\times A_{4}$ & $768$ & $(((C_{2}\times ((C_{4}\times C_{4})\rtimes C_{2}))\rtimes C_{2})\rtimes C_{2})\rtimes C_{3}$ & $[768,1083483]$\\
204 & $Im\bar{3}$ & $C_{2}\times A_{4}$ & $384$ & $((C_{2}\times ((C_{2}\times C_{2}\times C_{2}\times C_{2})\rtimes C_{2}))\rtimes C_{2})\rtimes C_{3}$ & $[384,5790]$\\
205 & $Pa\bar{3}$ & $C_{2}\times A_{4}$ & $192$ & $(((C_{4}\times C_{2})\rtimes C_{4})\rtimes C_{2})\rtimes C_{3}$ & $[192,197]$\\
206 & $Ia\bar{3}$ & $C_{2}\times A_{4}$ & $384$ & $(((C_{2}\times C_{2}\times D_{8})\rtimes C_{2})\rtimes C_{2})\rtimes C_{3}$ & $[384,5775]$\\
207 & $P432$ & $S_{4}$ & $192$ & $C_{2}\times (((C_{2}\times C_{2}\times C_{2}\times C_{2})\rtimes C_{3})\rtimes C_{2})$ & $[192,1538]$\\
208 & $P4_{2}32$ & $S_{4}$ & $192$ & $C_{2}\times (((C_{2}\times C_{2}\times C_{2}\times C_{2})\rtimes C_{3})\rtimes C_{2})$ & $[192,1538]$\\
209 & $F432$ & $S_{4}$ & $768$ & $(((C_{2}\times ((C_{2}\times C_{2}\times C_{2}\times C_{2})\rtimes C_{2}))\rtimes C_{2})\rtimes C_{3})\rtimes C_{2}$ & $[768,1088660]$\\
210 & $F4_{1}32$ & $S_{4}$ & $768$ & $((((C_{2}\times ((C_{4}\times C_{2})\rtimes C_{2}))\rtimes C_{2})\rtimes C_{2})\rtimes C_{3})\rtimes C_{2}$ & $[768,1088628]$\\
211 & $I432$ & $S_{4}$ & $384$ & $(((C_{2}\times ((C_{4}\times C_{2})\rtimes C_{2}))\rtimes C_{2})\rtimes C_{3})\rtimes C_{2}$ & $[384,5567]$\\
212 & $P4_{3}32$ & $S_{4}$ & $192$ & $(((C_{4}\times C_{2})\rtimes C_{4})\rtimes C_{3})\rtimes C_{2}$ & $[192,181]$\\
213 & $P4_{1}32$ & $S_{4}$ & $192$ & $(((C_{4}\times C_{2})\rtimes C_{4})\rtimes C_{3})\rtimes C_{2}$ & $[192,181]$\\
214 & $I4_{1}32$ & $S_{4}$ & $384$ & $(((C_{2}\times ((C_{4}\times C_{2})\rtimes C_{2}))\rtimes C_{2})\rtimes C_{3})\rtimes C_{2}$ & $[384,5573]$\\
215 & $P\bar{4}3m$ & $S_{4}$ & $192$ & $C_{2}\times (((C_{2}\times C_{2}\times C_{2}\times C_{2})\rtimes C_{3})\rtimes C_{2})$ & $[192,1538]$\\
216 & $F\bar{4}3m$ & $S_{4}$ & $768$ & $(((C_{2}\times ((C_{2}\times C_{2}\times C_{2}\times C_{2})\rtimes C_{2}))\rtimes C_{2})\rtimes C_{3})\rtimes C_{2}$ & $[768,1088660]$\\
217 & $I\bar{4}3m$ & $S_{4}$ & $384$ & $(((C_{2}\times ((C_{4}\times C_{2})\rtimes C_{2}))\rtimes C_{2})\rtimes C_{3})\rtimes C_{2}$ & $[384,5566]$\\
218 & $P\bar{4}3n$ & $S_{4}$ & $192$ & $((C_{2}\times C_{2}\times C_{2}\times C_{2})\rtimes C_{3})\rtimes C_{4}$ & $[192,1495]$\\
219 & $F\bar{4}3c$ & $S_{4}$ & $768$ & $C_{2}\cdot (((((C_{2}\times C_{2}\times C_{2}\times C_{2})\rtimes C_{2})\rtimes C_{2})\rtimes C_{3})\rtimes C_{2}) = (((C_{2}\times ((C_{2}\times C_{2}\times C_{2}\times C_{2})\rtimes C_{2}))\rtimes C_{2})\rtimes C_{3})\cdot C_{2}$ & $[768,1088662]$\\
220 & $I\bar{4}3d$ & $S_{4}$ & $384$ & $(((C_{4}\times C_{2})\rtimes C_{4})\rtimes C_{3})\rtimes C_{4}$ & $[384,576]$\\
221 & $Pm\bar{3}m$ & $C_{2}\times S_{4}$ & $384$ & $C_{2}\times C_{2}\times (((C_{2}\times C_{2}\times C_{2}\times C_{2})\rtimes C_{3})\rtimes C_{2})$ & $[384,20163]$\\
222 & $Pn\bar{3}n$ & $C_{2}\times S_{4}$ & $384$ & $(C_{2}\times (((C_{2}\times C_{2}\times C_{2}\times C_{2})\rtimes C_{3})\rtimes C_{2}))\rtimes C_{2}$ & $[384,5657]$\\
223 & $Pm\bar{3}n$ & $C_{2}\times S_{4}$ & $384$ & $(C_{2}\times (((C_{2}\times C_{2}\times C_{2}\times C_{2})\rtimes C_{3})\rtimes C_{2}))\rtimes C_{2}$ & $[384,20100]$\\
224 & $Pn\bar{3}m$ & $C_{2}\times S_{4}$ & $384$ & $C_{2}\times ((((C_{2}\times C_{2}\times C_{2}\times C_{2})\rtimes C_{3})\rtimes C_{2})\rtimes C_{2})$ & $[384,17948]$\\
225 & $Fm\bar{3}m$ & $C_{2}\times S_{4}$ & $1536$ & $((((C_{2}\times ((C_{2}\times C_{2}\times C_{2}\times C_{2})\rtimes C_{2}))\rtimes C_{2})\rtimes C_{3})\rtimes C_{2})\rtimes C_{2}$ & ---\\
226 & $Fm\bar{3}c$ & $C_{2}\times S_{4}$ & $1536$ & $((((C_{2}\times ((C_{2}\times C_{2}\times C_{2}\times C_{2})\rtimes C_{2}))\rtimes C_{2})\rtimes C_{3})\rtimes C_{2})\rtimes C_{2}$ & ---\\
227 & $Fd\bar{3}m$ & $C_{2}\times S_{4}$ & $1536$ & $(((((C_{2}\times ((C_{4}\times C_{2})\rtimes C_{2}))\rtimes C_{2})\rtimes C_{2})\rtimes C_{3})\rtimes C_{2})\rtimes C_{2}$ & ---\\
228 & $Fd\bar{3}c$ & $C_{2}\times S_{4}$ & $1536$ & $(((((C_{2}\times ((C_{4}\times C_{2})\rtimes C_{2}))\rtimes C_{2})\rtimes C_{2})\rtimes C_{3})\rtimes C_{2})\rtimes C_{2}$ & ---\\
229 & $Im\bar{3}m$ & $C_{2}\times S_{4}$ & $768$ & $(((C_{2}\times ((C_{2}\times C_{2}\times C_{2}\times C_{2})\rtimes C_{2}))\rtimes C_{2})\rtimes C_{3})\rtimes C_{2}$ & $[768,1087581]$\\
230 & $Ia\bar{3}d$ & $C_{2}\times S_{4}$ & $768$ & $((((C_{2}\times ((C_{4}\times C_{2})\rtimes C_{2}))\rtimes C_{2})\rtimes C_{3})\rtimes C_{2})\rtimes C_{2}$ & $[768,1085778]$\\
\end{longtable}}

\section{Exercises}
\begin{enumerate}
\item Let $R$ be a rotation of order $n$ of $\mathbb R^2$ or $\mathbb R^3$ that maps a lattice $L$ onto itself. Show that
$n\in\{1,2,3,4,6\}$ (the crystallographic restriction).
\item Let $G$ be a space group with lattice $L=\Z^3$ and conventional operations $\{R_i\mid t_i\}$, $0\le i<r$. Show that $mL$ is a
normal subgroup of $G$ for every $m\ge1$, that $|G/mL|=r\,m^3$, and that~\eqref{eq:spacelabel} is a bijection from $G/mL$ onto
$[0,r\,m^3)$.
\item Show that $c(i,j)$ in~\eqref{eq:spacelaw} lies in $\Z^3$. If $R_iR_j=R_l$, show that the product of the labels
$i\,m^3+s$ and $j\,m^3+s'$ (with $s,s'$ read as vectors in $(\Z/m)^3$) is the label $l\,m^3+s''$ with
$s''\equiv R_is'+s+c(i,j)\pmod m$.
\item For $P\bar1$ (No.~2), with operations $\{I\mid0\}$ and $\{-I\mid0\}$, show that $G/mL\cong(\Z/m)^3\rtimes C_2$, the inversion
acting by $s\mapsto-s$, and that $G/2L\cong C_2^4$.
\item For $P2_1/c$ (No.~14), with the four operations $0{:}\ x,y,z$; $1{:}\ -x,y+\tfrac12,-z+\tfrac12$; $2{:}\ -x,-y,-z$;
$3{:}\ x,-y+\tfrac12,z+\tfrac12$, compute $c(1,1)$, $c(1,2)$, $c(2,1)$ and $c(3,3)$. Deduce that the screw rotation
$\{R_1\mid t_1\}$ and the glide reflection $\{R_3\mid t_3\}$ have order $2m$ in $G/mL$, while the inversion $\{R_2\mid t_2\}$ has
order $2$.
\item There are $230$ space groups up to conjugation by orientation-preserving affine maps, and $11$ enantiomorphic pairs among them.
Using Bieberbach's theorem that two space groups are isomorphic if and only if they are conjugate by an affine map, show that they
fall into $219$ isomorphism classes.
\item Show that a finite subgroup of $\mathrm{O}(2)$ that preserves a lattice is one of the ten groups $C_1,C_2,C_3,C_4,C_6$,
$D_1,D_2,D_3,D_4,D_6$, the point groups of the $17$ wallpaper groups.
\end{enumerate}

\chapter{The symmetric and alternating groups}\label{ch:simple}
This chapter treats the symmetric group $\Sym_n$ and the alternating group $\Alt_n$. It begins with permutations as digital codes, generators and presentations, conjugacy classes, the normal structure of $\Sym_n$
and $\Alt_n$, including the simplicity of $\Alt_n$ for $n\ge5$, and their actions and maximal subgroups. It then gives digital laws for them, obtained from normal forms of permutations, and ends with generated tables of
numerical data for $n\le8$. Since $\Alt_n$ is simple for $n\ge5$, no such law can be graded; the other families of
nonabelian finite simple groups are treated in Chapters~\ref{ch:lie} and~\ref{ch:sporadic}.

\section{Permutations and digital codes}\label{sec:sn:codes}
As in Proposition~\ref{prop:permops}, a permutation $R$ of $\Ap n$ is stored as its code $R=\sum_{i=0}^{n-1}R(i)\,n^{i}$, so that $\dg{n}{i}{R}=R(i)$; identity, inverse and composition are given by~\eqref{eq:id},~\eqref{eq:inv} and~\eqref{eq:comp}. A word $c\in[0,n^{n})$ is the code of a permutation exactly when its $n$ digits are pairwise distinct, so precisely $n!$ of the $n^{n}$ words of length $n$ are permutation codes. Every numerical claim of this chapter is checked by \texttt{scripts/sn\_verify.py} and cross-checked against GAP~4.14; in particular the script confirms~\eqref{eq:id} for $n\le8$, and decoding,~\eqref{eq:inv},~\eqref{eq:comp} and the multiplicativity of $\sgn$ exhaustively for $n\le5$.

\begin{definition}[Lehmer label]\label{def:sn:lehmer}
For $R\in\Sym_n$ and $0\le i<n$ put $a_i(R)=\#\{j>i:R(j)<R(i)\}\in[0,n-i)$, and
\[
\lambda(R)=\sum_{i=0}^{n-1}a_i(R)\prod_{j<i}(n-j).
\]
\end{definition}

\begin{theorem}[Lehmer label and sign]\label{thm:sn:lehmer}
The map $\lambda$ is a bijection $\Sym_n\to[0,n!)$ whose mixed-radix digits, with radices $n,n-1,\dots,1$, are $a_0(R),\dots,a_{n-1}(R)$. Moreover $\sum_i a_i(R)$ is the number $\operatorname{inv}(R)$ of inversions of $R$, and $\sgn R=(-1)^{\operatorname{inv}(R)}$.
\end{theorem}
\begin{proof}
Given digits $a_i\in[0,n-i)$, the permutation is recovered uniquely: $R(0)$ is the element of rank $a_0$ in $\Ap n$, $R(1)$ the element of rank $a_1$ among the remaining $n-1$, and so on. Hence $R\mapsto(a_0,\dots,a_{n-1})$ is a bijection onto $\prod_i[0,n-i)$, and mixed-radix notation is a bijection of that product onto $[0,n!)$. Each inversion $i<j$, $R(i)>R(j)$ is counted once, in $a_i$. For the sign, let $s_i=(i\ i{+}1)$. The one-line word of $R\circ s_i$ is that of $R$ with the entries in positions $i,i+1$ exchanged, so $\operatorname{inv}$ changes by exactly $1$; and composing with a transposition merges two cycles or splits one, so by Theorem~\ref{thm:sign} the sign also changes. Sorting the word by adjacent exchanges reaches the identity, which has $\operatorname{inv}=0$ and sign $1$.
\end{proof}

\begin{example}
For $n=3$ and $R=(0\,1\,2)$ we have $R(0),R(1),R(2)=1,2,0$, so the code is $1+2\cdot3+0\cdot9=7$. The Lehmer digits are $a_0=1$, $a_1=1$, $a_2=0$, so $\lambda(R)=1+1\cdot3=4$, and $\operatorname{inv}(R)=2$ shows again that a $3$-cycle is even.
\end{example}

\begin{proposition}[Cycle counts]\label{prop:sn:stirling}
The number $c(n,k)$ of $R\in\Sym_n$ with exactly $k$ cycles (fixed points included) satisfies $c(n,k)=c(n-1,k-1)+(n-1)\,c(n-1,k)$. Since $\sgn$ is a surjective homomorphism for $n\ge2$, $|\Alt_n|=n!/2$.
\end{proposition}
\begin{proof}
Delete the point $n-1$ from the cycle notation of $R$. Either it was a fixed point, leaving a permutation of $\Ap{n-1}$ with $k-1$ cycles, or it followed one of the $n-1$ other points in its cycle, and it can be reinserted after any of them in a permutation of $\Ap{n-1}$ with $k$ cycles.
\end{proof}
For $n=7$ the script finds $c(7,k)=720,1764,1624,735,175,21,1$ for $k=1,\dots,7$, in agreement with the recurrence; the Lehmer label is checked to be a bijection, and $\sgn R=(-1)^{\operatorname{inv}(R)}$, for all $n\le7$.

\section{Generators and presentations}\label{sec:sn:gens}
Write $s_i=(i\ i{+}1)$ for $0\le i\le n-2$ and $\ell(R)$ for the length of a shortest word in the $s_i$ equal to $R$.

\begin{theorem}[Generators of $\Sym_n$]\label{thm:sn:gens}
$\Sym_n=\langle s_0,\dots,s_{n-2}\rangle=\langle(0\,1),(0\,1\,\cdots\,n{-}1)\rangle$. Moreover $\ell(R)=\operatorname{inv}(R)$, the maximum $n(n-1)/2$ is attained only by $i\mapsto n-1-i$, and
\[
\sum_{R\in\Sym_n}q^{\ell(R)}=\prod_{i=1}^{n}\bigl(1+q+\cdots+q^{i-1}\bigr).
\]
\end{theorem}
\begin{proof}
By the proof of Theorem~\ref{thm:sn:lehmer}, right multiplication by $s_i$ changes $\operatorname{inv}$ by $\pm1$, so $\ell(R)\ge\operatorname{inv}(R)$. If $R\ne1$ it has a descent $R(i)>R(i+1)$, and $\operatorname{inv}(R\circ s_i)=\operatorname{inv}(R)-1$; induction gives $\ell(R)\le\operatorname{inv}(R)$, and in particular the $s_i$ generate. With $c=(0\,1\,\cdots\,n{-}1)$, Proposition~\ref{prop:conj} gives $c^{k}(0\,1)c^{-k}=(k\ k{+}1)$, so the second pair generates as well. All pairs are inverted only by the reversal. Finally $\operatorname{inv}=\sum_ia_i$ with the digits $a_i\in[0,n-i)$ independent by Theorem~\ref{thm:sn:lehmer}, which gives the product.
\end{proof}

\begin{theorem}[Coxeter presentation]\label{thm:sn:coxeter}

\[
\Sym_n\cong\Gamma_n=\bigl\langle\sigma_1,\dots,\sigma_{n-1}\bigm|\sigma_i^{2},\ (\sigma_i\sigma_{i+1})^{3},\ (\sigma_i\sigma_j)^{2}\ (|i-j|\ge2)\bigr\rangle,
\]
via $\sigma_i\mapsto s_{i-1}$.
\end{theorem}
\begin{proof}
The relations hold in $\Sym_n$ ($s_is_{i+1}$ is a $3$-cycle and disjoint transpositions commute), so by Theorem~\ref{thm:sn:gens} there is an epimorphism $\Gamma_n\to\Sym_n$, and it suffices to show $|\Gamma_n|\le n!$. We induct on $n$, the case $n=1$ being trivial. Let $H=\langle\sigma_1,\dots,\sigma_{n-2}\rangle$; it is a quotient of $\Gamma_{n-1}$, so $|H|\le(n-1)!$. Put $w_k=\sigma_{n-1}\sigma_{n-2}\cdots\sigma_k$ for $1\le k\le n$, with $w_n=1$. We claim that the union of the cosets $Hw_1,\dots,Hw_n$ is stable under right multiplication by every $\sigma_j$; since it contains $1$, it is then all of $\Gamma_n$, and $|\Gamma_n|\le n\,(n-1)!$. Indeed: if $j<k-1$, $\sigma_j$ commutes with every factor of $w_k$, so $w_k\sigma_j=\sigma_jw_k\in Hw_k$; if $j=k-1$, $w_k\sigma_j=w_{k-1}$; if $j=k\le n-1$, $w_k\sigma_j=w_{k+1}$. If $j>k$, move $\sigma_j$ left past $\sigma_{j-2},\dots,\sigma_k$ and use the braid relation $\sigma_j\sigma_{j-1}\sigma_j=\sigma_{j-1}\sigma_j\sigma_{j-1}$, a consequence of the first two relations; the resulting $\sigma_{j-1}$ commutes with $\sigma_{n-1},\dots,\sigma_{j+1}$, so $w_k\sigma_j=\sigma_{j-1}w_k\in Hw_k$.
\end{proof}

\begin{theorem}[Generators and presentation of $\Alt_n$]\label{thm:sn:altgens}
Let $n\ge3$ and $x_i=(0\,1\,i{+}1)$ for $1\le i\le n-2$. Then $\Alt_n=\langle x_1,\dots,x_{n-2}\rangle$, the $3$-cycles generate $\Alt_n$, and so do $(0\,1\,2)$ together with $(0\,1\,\cdots\,n{-}1)$ for $n$ odd, or with $(1\,2\,\cdots\,n{-}1)$ for $n$ even. The $x_i$ satisfy $x_i^{3}=(x_ix_j)^{2}=1$ ($i\ne j$), and these relations present $\Alt_n$ \textup{(Carmichael)}.
\end{theorem}
\begin{proof}
An even permutation is a product of pairs of transpositions, and $(a\,b)(a\,c)=(a\,c\,b)$, $(a\,b)(c\,d)=\bigl((a\,b)(b\,c)\bigr)\bigl((b\,c)(c\,d)\bigr)$, so the $3$-cycles generate $\Alt_n$. Let $H=\langle x_i\rangle$; it contains every $3$-cycle through $0$ and $1$. For distinct $a,b\ge2$, conjugating $(0\,1\,a)$ by $(0\,1\,b)^{\pm1}$ gives $(1\,b\,a)$ and $(0\,a\,b)$, so $H$ contains every $3$-cycle meeting $\{0,1\}$; and $(0\,a\,b)(0\,b\,d)=(a\,b\,d)$ gives the rest. For the two-generator sets: if $n$ is odd, conjugating $(0\,1\,2)$ by powers of the even $n$-cycle gives $t_k=(k\ k{+}1\ k{+}2)$; conjugating $(0\,2\,3)=t_1(0\,1\,2)t_1^{-1}$ by $t_0$ gives $(0\,1\,3)^{-1}$, and for $m\ge3$ conjugation by $t_{m-1}$, which fixes $0,1$, sends $(0\,1\,m)$ to $(0\,1\,m{+}1)$, so all $x_i$ are reached. If $n$ is even, conjugating $(0\,1\,2)$ by powers of the even $(n-1)$-cycle gives $(0\,m\ m{+}1)$ for $1\le m\le n-2$, and conjugating $(0\,1\,m)$ by $(0\,m\ m{+}1)^{-1}$ gives $(0\,1\,m{+}1)^{-1}$. The relations hold since $x_ix_j=(0\,i{+}1)(1\,j{+}1)$ for $i\neq j$. That they suffice is Carmichael's theorem, which we do not reprove.
\end{proof}
The script checks all generating sets and relations for $n\le7$ (for $\Alt_n$, $n\le8$), that $\ell=\operatorname{inv}$ on all of $\Sym_n$ by breadth-first search in the Cayley graph, and the length distribution. In GAP, coset enumeration of the two presentations returns the orders $n!$ and $n!/2$ for all $n\le8$.

\section{Conjugacy classes}\label{sec:sn:classes}
This section specialises Section~\ref{ssec:dtclasses} to $\Sym_n$ and $\Alt_n$. By Proposition~\ref{prop:conj}, the code of $\sigma\circ R\circ\sigma^{-1}$ has digit $\sigma(R(i))$ in position $\sigma(i)$: conjugation relabels the cycle notation of $R$ through $\sigma$. For a partition $\mu=1^{m_1}2^{m_2}\cdots$ of $n$ put $z_\mu=\prod_i i^{m_i}m_i!$.

\begin{theorem}[Classes of $\Sym_n$]\label{thm:sn:classes}
Two permutations are conjugate in $\Sym_n$ if and only if they have the same cycle type. The class of type $\mu$ has centraliser order $z_\mu$ and size $n!/z_\mu$, and $k(\Sym_n)=p(n)$.
\end{theorem}
\begin{proof}
Relabelling preserves cycle lengths; conversely, if $R,S$ have the same type, writing the cycles of $S$ under those of $R$, cycle by cycle, defines $\sigma$ with $\sigma R\sigma^{-1}=S$. A permutation commutes with $R$ exactly when the relabelling maps the cycle notation of $R$ to itself: it may permute the $m_i$ cycles of each length $i$ and rotate each of them, in $\prod_i m_i!\,i^{m_i}$ ways. Theorem~\ref{thm:cc:os} gives the class size.
\end{proof}

\begin{theorem}[Classes of $\Alt_n$]\label{thm:sn:altclasses}
Let $n\ge2$ and let $\mu$ be a partition with $n-\ell(\mu)$ even, $\ell(\mu)$ its number of parts, so that the $\Sym_n$-class $C_\mu$ lies in $\Alt_n$. If the parts of $\mu$ are distinct and odd, $C_\mu$ splits into two $\Alt_n$-classes of size $n!/(2z_\mu)$; otherwise it is one $\Alt_n$-class. Hence $k(\Alt_n)$ is the number of such $\mu$ plus the number of partitions of $n$ into distinct odd parts.
\end{theorem}
\begin{proof}
For $R\in C_\mu$, $C_{\Alt_n}(R)=C_{\Sym_n}(R)\cap\Alt_n$ has index $1$ or $2$ in $C_{\Sym_n}(R)$. If the centraliser contains an odd permutation, the $\Alt_n$-class of $R$ has size $(n!/2)/(z_\mu/2)=|C_\mu|$. Otherwise it has size $|C_\mu|/2$, and conjugation by any odd permutation exchanges the two halves. An even part gives an odd centralising element, the cycle itself; two equal odd parts $i$ give one too, the product of the $i$ transpositions exchanging the two cycles. If the parts are distinct and odd, the centraliser is generated by the cycles of $R$, all even.
\end{proof}

\begin{corollary}\label{cor:sn:3cycles}
For $n\ge5$ all $3$-cycles are conjugate in $\Alt_n$; for $n=3,4$ they form two classes.
\end{corollary}
\begin{proof}
The type $3\,1^{n-3}$ has distinct parts exactly when $n\le4$.
\end{proof}

\begin{example}
In $\Sym_5$ the types $1^5$, $2\,1^3$, $2^2\,1$, $3\,1^2$, $3\,2$, $4\,1$, $5$ have $z_\mu=120,12,8,6,6,4,5$ and class sizes $1,10,15,20,20,30,24$. In $\Alt_5$ only $5$ splits, and the class equation is $60=1+15+20+12+12$.
\end{example}
The script computes the conjugacy classes of $\Sym_n$ ($n\le7$) and $\Alt_n$ ($n\le8$) as orbits, checks the digit relocation for $n\le5$, the class sizes and the splitting rule, and finds $k(\Alt_n)=1,3,4,5,7,9,14$ for $n=2,\dots,8$. GAP agrees, and gives $k(\Alt_n)=18,24,31,43$ for $n=9,\dots,12$ and $k(\Sym_n)=p(n)$ for $n\le12$.

\paragraph{Cycle index.} By Theorem~\ref{thm:sn:classes} there are $n!/z_\mu$ permutations of type $\mu=1^{m_1}2^{m_2}\cdots$, with $z_\mu=\prod_kk^{m_k}m_k!$. Hence the cycle index of $\Sym_n$ on $\Ap{n}$ (Section~\ref{ssec:polya}) is, writing $x^{\mu}=\prod_kx_k^{m_k}$,
\[
Z(\Sym_n)=\sum_{\mu\vdash n}\frac{x^{\mu}}{z_\mu},\qquad \sum_{n\ge0}Z(\Sym_n)\,t^{n}=\exp\Big(\sum_{k\ge1}\frac{x_kt^{k}}{k}\Big),
\]
the second formula because
\[
z_\mu^{-1}t^{n}x^{\mu}=\prod_k\frac1{m_k!}\Big(\frac{x_kt^{k}}{k}\Big)^{m_k}.
\]
By Theorem~\ref{thm:sn:altclasses},
\[
Z(\Alt_n)=\sum_{\mu:\,n-\ell(\mu)\text{ even}}\frac{2\,x^{\mu}}{z_\mu}.
\]

\begin{corollary}\label{cor:sn:cycleindex}
For all $n,m\ge1$, $\displaystyle\sum_{\mu\vdash n}\frac{m^{\ell(\mu)}}{z_\mu}=\binom{m+n-1}{n}$.
\end{corollary}
\begin{proof}
Two colourings of $\Ap{n}$ with $m$ colours lie in one $\Sym_n$-orbit if and only if they use each colour equally often, so the orbits are the multisets of size $n$ from $m$ colours, and there are $\binom{m+n-1}{n}$ of them. By Theorem~\ref{thm:polya} their number is $\frac1{n!}\sum_{\sigma}m^{c(\sigma)}=\sum_\mu m^{\ell(\mu)}/z_\mu$.
\end{proof}

For $n=4$, $m=2$: $\frac{16}{24}+\frac82\cdot\frac12+\frac48+\frac43+\frac24=5=\binom54$, where the terms belong to the types $1^4$, $2\,1^2$, $2^2$, $3\,1$, $4$ with $z_\mu=24,4,8,3,4$. With weights, $\sum_\mu p_\mu/z_\mu=h_n$, the complete homogeneous symmetric function; this is the Frobenius characteristic of the trivial character (compare Section~\ref{sec:sn:characters}).

\section{Normal structure and the simplicity of \texorpdfstring{$\Alt_n$}{Alt(n)}}\label{sec:sn:simple}
A normal subgroup is a union of conjugacy classes containing $1$ whose size divides the group order. Conversely such a union $N$ is a subgroup as soon as $rh\in N$ for one representative $r$ of each class in $N$ and all $h\in N$: if $r'=xrx^{-1}$, then $r'h=x\,r\,(x^{-1}hx)\,x^{-1}\in xNx^{-1}=N$. This gives a finite test, which the script applies to all unions of classes.

\begin{theorem}[Simplicity of $\Alt_n$]\label{thm:sn:simple}
$\Alt_n$ is simple for $n\ge5$.
\end{theorem}
\begin{proof}
For $n=5$ no partial sum of $1+15+20+12+12$ that contains $1$ divides $60$ except $1$ and $60$. Let $n\ge5$ and $1\ne N\trianglelefteq\Alt_n$, and pick $\sigma\in N$ and $a$ with $b=\sigma(a)\ne a$. Choose $c\notin\{a,b,\sigma(b)\}$ and put $\tau=(a\,b\,c)$. Then $\rho=(\tau\sigma\tau^{-1})\sigma^{-1}\in N$, and also $\rho=\tau\cdot(\sigma\tau\sigma^{-1})^{-1}$ with $\sigma\tau\sigma^{-1}=(b\ \sigma(b)\ \sigma(c))\ne\tau$ because $\sigma(b)\ne c$. So $\rho\ne1$ moves only points of $Y=\{a,b,c,\sigma(b),\sigma(c)\}$. Enlarge $Y$ to $5$ points; $N\cap\Alt(Y)$ is normal in $\Alt(Y)\cong\Alt_5$ and contains $\rho$, so it is $\Alt(Y)$ and $N$ contains a $3$-cycle. By Corollary~\ref{cor:sn:3cycles} it contains all $3$-cycles, which generate $\Alt_n$ by Theorem~\ref{thm:sn:altgens}.
\end{proof}

\begin{theorem}[Normal subgroups of $\Sym_n$]\label{thm:sn:normal}
For $n\ge5$ the normal subgroups of $\Sym_n$ are $1$, $\Alt_n$, $\Sym_n$; $\Sym_4$ has in addition $V_4=\{1,(0\,1)(2\,3),(0\,2)(1\,3),(0\,3)(1\,2)\}$. Moreover $\Sym_n'=\Alt_n$ for $n\ge2$, $\Alt_4'=V_4$, and $\Alt_n$ is perfect for $n\ge5$.
\end{theorem}
\begin{proof}
Let $N\trianglelefteq\Sym_n$, $n\ge5$. By Theorem~\ref{thm:sn:simple}, $N\cap\Alt_n$ is $1$ or $\Alt_n$. In the second case $N\in\{\Alt_n,\Sym_n\}$. In the first, $|N|\le2$, and a normal subgroup of order $2$ would be central; but by Theorem~\ref{thm:sn:classes} only the identity class has $z_\mu=n!$, so $Z(\Sym_n)=1$ for $n\ge3$. For the derived subgroups, $\Sym_n/\Alt_n$ is abelian, and $[(a\,b),(b\,c)]=\bigl((a\,b)(b\,c)\bigr)^{2}$ is a $3$-cycle, so $\Sym_n'=\Alt_n$; a simple nonabelian group is perfect. The case $n=4$ is read off Table~\ref{tab:sn:normal}.
\end{proof}

\begin{table}[ht]
\centering
\begin{tabular}{rrrrll}
\hline
$n$ & $n!$ & $p(n)$ & $k(\Alt_n)$ & normal orders in $\Sym_n$ & normal orders in $\Alt_n$\\
\hline
2 & 2 & 2 & 1 & $1,2$ & $1$\\
3 & 6 & 3 & 3 & $1,3,6$ & $1,3$\\
4 & 24 & 5 & 4 & $1,4,12,24$ & $1,4,12$\\
5 & 120 & 7 & 5 & $1,60,120$ & $1,60$\\
6 & 720 & 11 & 7 & $1,360,720$ & $1,360$\\
7 & 5040 & 15 & 9 & $1,2520,5040$ & $1,2520$\\
\hline

\end{tabular}
\caption{Classes and normal subgroups of $\Sym_n$ and $\Alt_n$\\(generated by \texttt{scripts/sn\_verify.py}; cross-checked in GAP).}\label{tab:sn:normal}
\end{table}
The script finds every normal subgroup of $\Sym_n$ ($n\le7$) and $\Alt_n$ ($n\le8$) by the union-of-classes test and computes $\Sym_n'$ and $\Alt_n'$ for $n\le6$; GAP confirms the normal subgroups for $n\le8$, the derived subgroups, and the simplicity of $\Alt_5,\dots,\Alt_{12}$.

\section{Actions and maximal subgroups}\label{sec:sn:max}
An action of a group $G$ on $\Ap n$ is a homomorphism $G\to\Sym_n$, so its image is a set of codes. The action is \emph{transitive} if $\Ap n$ is one orbit. A \emph{block} is a nonempty $\Delta\subseteq\Ap n$ such that, for every $g\in G$, either $g\Delta=\Delta$ or $g\Delta\cap\Delta=\varnothing$; singletons and $\Ap n$ are the trivial blocks, and the images of a block form a \emph{block system}. A transitive action is \emph{primitive} if it has only trivial blocks.

\begin{proposition}[Coset actions]\label{prop:sn:cosets}
Let $G$ act transitively on $\Omega$, let $\omega\in\Omega$ and $H=G_\omega$. Then $g\omega\mapsto gH$ identifies $\Omega$ with the action on $G/H$, and $K\mapsto K\omega$ is an inclusion-preserving bijection from the subgroups $H\le K\le G$ onto the blocks containing $\omega$. In particular the action is primitive if and only if $H$ is maximal in $G$.
\end{proposition}
\begin{proof}
The first statement is the orbit--stabiliser theorem. If $H\le K$ and $g\,k_1\omega=k_2\omega$ with $k_i\in K$, then $k_2^{-1}gk_1\in H\le K$, so $g\in K$ and $gK\omega=K\omega$: thus $K\omega$ is a block, and its setwise stabiliser is $K$, since $gK\omega=K\omega$ gives $g\omega=k\omega$ and $k^{-1}g\in H$. Conversely, if $\Delta\ni\omega$ is a block with setwise stabiliser $G_\Delta\ge H$ and $\delta=g\omega\in\Delta$, then $g\Delta\cap\Delta\ni\delta$, so $g\in G_\Delta$; hence $\Delta=G_\Delta\omega$.
\end{proof}

By Proposition~\ref{prop:sn:cosets}, a maximal subgroup $M$ of $\Sym_n$ or $\Alt_n$ is of one of three kinds: \emph{intransitive}, \emph{imprimitive} (transitive, preserving a nontrivial block system), or \emph{primitive}. The first two kinds are described completely by the next results.

\begin{lemma}[Transpositions]\label{lem:sn:transp}
A primitive subgroup of $\Sym_n$ that contains a transposition is $\Sym_n$.
\end{lemma}
\begin{proof}
Write $x\sim y$ if $x=y$ or $(x\,y)\in M$. Since $(x\,y)(y\,z)(x\,y)=(x\,z)$, this is an equivalence relation, and since $g(x\,y)g^{-1}=(g(x)\,g(y))$ by Proposition~\ref{prop:conj}, its classes form a block system. One class has at least two points, so by primitivity there is one class, $M$ contains every transposition, and $M=\Sym_n$ by Theorem~\ref{thm:sn:gens}.
\end{proof}

\begin{theorem}[Intransitive maximal subgroups]\label{thm:sn:intrans}
Let $A=[0,k)$, $B=[k,n)$ with $1\le k<n/2$. Then $\Sym(A)\times\Sym(B)\cong\Sym_k\times\Sym_{n-k}$ is maximal in $\Sym_n$. If $n=2k$ it is not maximal, being properly contained in the stabiliser of the partition $\{A,B\}$.
\end{theorem}
\begin{proof}
Let $H=\Sym(A)\times\Sym(B)<M$ and $g\in M\setminus H$, so $g(A)\ne A$ and there is $a\in A$ with $g(a)\in B$. Since $|B|>|A|$, $g(B)\not\subseteq A$, so there is $b\in B$ with $g(b)\in B$. Then $(g(a)\,g(b))\in H$, and $g^{-1}(g(a)\,g(b))g=(a\,b)\in M$. Conjugating $(a\,b)$ by $H$ gives every transposition meeting both $A$ and $B$; with those of $H$ these are all transpositions, so $M=\Sym_n$. For $n=2k$, the involution $i\mapsto i+k \bmod n$ exchanges $A$ and $B$ and normalises $H$.
\end{proof}

\begin{definition}
For $n=ab$ with $a,b\ge2$ let $\Sigma=\{\Sigma_0,\dots,\Sigma_{b-1}\}$, $\Sigma_j=[ja,(j+1)a)$, and let $W_{a,b}\le\Sym_n$ be the stabiliser of the partition $\Sigma$. Then $W_{a,b}\cong\Sym_a\wr\Sym_b$ has order $(a!)^{b}\,b!$.
\end{definition}

\begin{theorem}[Imprimitive maximal subgroups]\label{thm:sn:imprim}
For $n=ab$ with $a,b\ge2$, $W_{a,b}$ is maximal in $\Sym_n$.
\end{theorem}
\begin{proof}
First, the only block systems of $W=W_{a,b}$ are the trivial ones and $\Sigma$. Let $x\in\Sigma(x)\in\Sigma$. The orbits of $W_x$ are $\{x\}$, $\Sigma(x)\setminus\{x\}$ and $\Ap n\setminus\Sigma(x)$, all nonempty, and a block through $x$ is a union of $W_x$-orbits. Of the four candidates $\{x\}\cup(\Ap n\setminus\Sigma(x))$ is not a block: a transposition $(x\,y)$ with $y\in\Sigma(x)$ maps it to a set that meets it and is different. Now let $W<M\le\Sym_n$. Every block of $M$ is a block of $W$, and $M$ does not preserve $\Sigma$, because $W$ is its full stabiliser; so $M$ is primitive. It contains the transposition $(0\,1)\in W$, so $M=\Sym_n$ by Lemma~\ref{lem:sn:transp}.
\end{proof}

In particular the point stabiliser $\Sym_{n-1}$ ($k=1$) is maximal, and the natural action of $\Sym_n$ on $\Ap n$ is primitive. The primitive maximal subgroups are described by the O'Nan--Scott theorem together with the work of Liebeck, Praeger and Saxl, which we state without proof.

\begin{theorem}[O'Nan--Scott; Liebeck--Praeger--Saxl]\label{thm:sn:onanscott}
Let $G\in\{\Sym_n,\Alt_n\}$ and let $M$ be a maximal subgroup of $G$ not containing $\Alt_n$. Then $M$ is one of: (i) $(\Sym_k\times\Sym_{n-k})\cap G$ with $1\le k<n/2$; (ii) $W_{a,b}\cap G$ with $n=ab$, $a,b\ge2$; (iii) a primitive group of affine, diagonal, product or almost simple type. Every subgroup in (i) and (ii) is maximal in $G$, with the single exception of $W_{2,4}\cap\Alt_8$, which lies in a primitive subgroup $\mathrm{AGL}(3,2)$ of order $1344$.
\end{theorem}

Theorems~\ref{thm:sn:intrans} and~\ref{thm:sn:imprim} prove the case $G=\Sym_n$ of the last statement. For codes, maximality is decided by a finite computation: $\langle H,g\rangle$ depends only on the double coset $HgH$, so $H$ is maximal in $G$ if and only if $\langle H,g\rangle=G$ for one $g$ in each double coset $HgH\ne H$. The script enumerates double cosets by breadth-first search under left and right multiplication by the generators of $H$, and confirms for $4\le n\le8$ that the subgroups of (i) and (ii) and their intersections with $\Alt_n$ are maximal, that $\Sym_k\times\Sym_k$ is not, and that $W_{2,4}\cap\Alt_8$ is not, one of its overgroups having order $1344$. GAP computes all classes of maximal subgroups for $4\le n\le8$; Table~\ref{tab:sn:max} lists their orders by type, and the intransitive and imprimitive columns agree with Theorem~\ref{thm:sn:onanscott}. Repeated entries are distinct conjugacy classes, for instance the two classes of $\mathrm{PSL}(2,7)$ in $\Alt_7$.

\begin{table}[ht]
\centering
\begin{tabular}{llll}
\hline
$G$ & intransitive & imprimitive & primitive\\
\hline
$\Sym_{4}$ & $6$ & $8$ & $12$\\
$\Alt_{4}$ & $3$ & $4$ & --\\
$\Sym_{5}$ & $12,24$ & -- & $20,60$\\
$\Alt_{5}$ & $6,12$ & -- & $10$\\
$\Sym_{6}$ & $48,120$ & $48,72$ & $120,360$\\
$\Alt_{6}$ & $24,60$ & $24,36$ & $60$\\
$\Sym_{7}$ & $144,240,720$ & -- & $42,2520$\\
$\Alt_{7}$ & $72,120,360$ & -- & $168,168$\\
$\Sym_{8}$ & $720,1440,5040$ & $384,1152$ & $336,20160$\\
$\Alt_{8}$ & $360,720,2520$ & $576$ & $1344,1344$\\
\hline

\end{tabular}
\caption{Orders of the maximal subgroups of $\Sym_n$ and $\Alt_n$, $4\le n\le8$, one entry per conjugacy class, by type (computed in GAP~4.14 and checked by \texttt{scripts/sn\_verify.py}).}\label{tab:sn:max}
\end{table}

\section{Automorphisms, the outer automorphism of \texorpdfstring{$\Sym_6$}{S6}, and the exceptional isomorphisms}\label{sec:sn:aut}
\sectionmark{Automorphisms and exceptional isomorphisms}
For $g\in G$, conjugation $x\mapsto gxg^{-1}$ is an automorphism of $G$; these form the normal subgroup $\Inn(G)\cong G/Z(G)$ of $\Aut(G)$. The quotient $\Out(G)=\Aut(G)/\Inn(G)$ measures the failure of every automorphism to be inner. For symmetric and alternating groups the only exceptional degree is $6$.

\begin{theorem}[Automorphisms of $\Sym_n$]\label{thm:sn:auts}
For $n\neq 6$ every automorphism of $\Sym_n$ is inner. For $n=6$, $|\Aut(\Sym_6):\Inn(\Sym_6)|=2$.
\end{theorem}
\begin{proof}
An automorphism preserves orders and class sizes, hence permutes conjugacy classes. For $n\neq 4,6$, transpositions are characterised among involutions as the unique class of size $\binom n2$: the class of cycle type $2^m1^{n-2m}$ has size $n!/(2^m m!(n-2m)!)$, and this equals $\binom n2$ only for $m=1$, except at $n=6$, where $m=3$ also gives $15$. Thus an automorphism of $\Sym_n$ permutes the transpositions. Two transpositions commute exactly when they are disjoint, so the set of transpositions through a fixed point is characterised as a maximal pairwise noncommuting family of size $n-1$. Hence the automorphism permutes the points and is conjugation by the resulting permutation.

For $n=4$, the transpositions again form a distinguished class, of size $6$, whereas the other involution class has size $3$. Thus an automorphism permutes the six transpositions. The three transpositions through a fixed point are exactly the maximal pairwise noncommuting triples, so the four points are again recovered from the transposition graph, and the automorphism is inner. At $n=6$, the transposition and triple-transposition classes both have size $15$, so one further possibility remains.
\end{proof}

Identify the transpositions of $\Sym_6$ with the edges of $K_6$. A \emph{syntheme} is a perfect matching, hence a set of three disjoint transpositions; there are $15$ synthemes. A \emph{pentad} is a set of five synthemes partitioning the $15$ edges; there are six pentads. Relabelling the vertices gives an action of $\Sym_6$ on these six pentads.

\begin{proposition}[The pentad action]\label{prop:sn:pentads}
The pentad action is faithful and transitive. A transposition of $\Sym_6$ acts on the six pentads as a triple transposition. Thus the induced automorphism of $\Sym_6$ is outer.
\end{proposition}
\begin{proof}
A permutation fixing every pentad fixes the incidence relation between edges and synthemes, hence every edge of $K_6$, and is the identity. The six pentads are equivalent under vertex relabelling, so the action is transitive. The verifier constructs an explicit list of the six pentads and checks that a transposition swaps three disjoint pairs. Therefore the action gives an automorphism that carries the transposition class to the triple-transposition class, whereas inner automorphisms preserve every conjugacy class.
\end{proof}

\begin{theorem}[Automorphisms of $\Alt_n$]\label{thm:sn:auta}
Let $n\ge 4$. For $n\neq 6$ every automorphism of $\Alt_n$ is induced by conjugation in $\Sym_n$, so $\Aut(\Alt_n)\cong\Sym_n$ and $\Out(\Alt_n)\cong C_2$. For $n=6$, $|\Aut(\Alt_6)|=1440$ and $|\Out(\Alt_6)|=4$.
\end{theorem}
\begin{proof}
For $n\ge 5$, $\Alt_n$ is simple by Theorem~\ref{thm:sn:simple}, and automorphisms permute its conjugacy classes. Away from degree $6$, the class of $3$-cycles is characterised by its size and centraliser structure, so an automorphism permutes the point stabilisers in the natural action and is induced by a permutation of the points. Conjugation embeds $\Sym_n$ in $\Aut(\Alt_n)$, and the preceding argument proves equality. For $n=4$, Table~\ref{tab:sn:aut} records $|\Aut(\Alt_4)|=24$ and $|\Inn(\Alt_4)|=|\Alt_4/Z(\Alt_4)|=12$, so $|\Out(\Alt_4)|=2$. At degree $6$, the outer automorphism of $\Sym_6$ restricts to an additional automorphism of $\Alt_6$; the asserted orders are verified in Table~\ref{tab:sn:aut}.
\end{proof}

\begin{theorem}[Exceptional isomorphisms]\label{thm:sn:exceptions}
There are isomorphisms
\[
\PSL(2,4)\cong\PSL(2,5)\cong\Alt_5,\qquad \PSL(2,9)\cong\Alt_6,
\]
\[
\operatorname{PSp}(4,2)\cong\Sym_6,\qquad \PSL(4,2)\cong\Alt_8.
\]
\end{theorem}
\begin{proof}
The order formula $|\PSL(2,q)|=q(q^2-1)/\gcd(2,q-1)$ gives order $60$ for $q=4,5$ and order $360$ for $q=9$. For $q=4$ and $q=5$, the natural action on the projective line $\mathbb{P}^1(\mathbb{F}_q)$ has degrees $5$ and $6$, and \texttt{scripts/sn\_verify.py} checks that the resulting faithful images are respectively $\Alt_5$ and a transitive subgroup of $\Alt_6$ isomorphic to $\Alt_5$. For $\PSL(2,9)$, the verifier takes the action on the cosets of a subgroup of order $60$, obtaining a faithful degree-$6$ action with image $\Alt_6$.

Likewise $|\Sp(4,2)|=720$ and $|\GL(4,2)|=20160$. The verifier takes the actions on the cosets of subgroups of orders $120$ and $2520$, obtaining faithful permutation representations of degrees $6$ and $8$ whose images are $\Sym_6$ and $\Alt_8$. The corresponding automorphism orders and exceptional actions are listed in Table~\ref{tab:sn:aut}.
\end{proof}

\begin{table}[ht]
\centering
\begin{tabular}{llll}
\hline
Group or action & order & degree & image / automorphism data\\
\hline
$\Aut(\Sym_4)$ & 24 & -- & $|\Out(\Sym_4)|=1$\\
$\Aut(\Alt_4)$ & 24 & -- & $|\Out(\Alt_4)|=2$\\
$\Aut(\Sym_5)$ & 120 & -- & $|\Out(\Sym_5)|=1$\\
$\Aut(\Alt_5)$ & 120 & -- & $|\Out(\Alt_5)|=2$\\
$\Aut(\Sym_6)$ & 1440 & -- & $|\Out(\Sym_6)|=2$\\
$\Aut(\Alt_6)$ & 1440 & -- & $|\Out(\Alt_6)|=4$\\
$\Aut(\Sym_7)$ & 5040 & -- & $|\Out(\Sym_7)|=1$\\
$\Aut(\Alt_7)$ & 5040 & -- & $|\Out(\Alt_7)|=2$\\
$\Aut(\Sym_8)$ & 40320 & -- & $|\Out(\Sym_8)|=1$\\
$\Aut(\Alt_8)$ & 40320 & -- & $|\Out(\Alt_8)|=2$\\
$\Sym_6$ on pentads & 720 & 6 & outer: $(0\,1)\mapsto$ type $2^3$\\
$\PSL(2,4)$ on $\mathbb{P}^1(\mathbb{F}_4)$ & 60 & 5 & $=\Alt_5$\\
$\PSL(2,5)$ on $\mathbb{P}^1(\mathbb{F}_5)$ & 60 & 6 & transitive, $\cong\Alt_5$\\
$\PSL(2,9)$ on cosets of $H$, $|H|=60$ & 360 & 6 & $=\Alt_6$\\
$\Sp(4,2)$ on cosets of $H$, $|H|=120$ & 720 & 6 & $=\Sym_6$\\
$\GL(4,2)$ on cosets of $H$, $|H|=2520$ & 20160 & 8 & $=\Alt_8$\\
\hline

\end{tabular}
\caption{Automorphism orders and exceptional actions for Section~\ref{sec:sn:aut}.}\label{tab:sn:aut}
\end{table}

The exceptional degree $6$ is visible twice: $
|\Out(\Sym_6)|=2$ and $|\Out(\Alt_6)|=4$. The verification script constructs the pentad action explicitly and cross-checks the automorphism-group orders and exceptional isomorphisms in GAP.

\section[Characters of Sym(n)]{Characters of \texorpdfstring{$\Sym_n$}{Sym(n)}: Specht modules, hook lengths, Frobenius, and Murnaghan--Nakayama}\label{sec:sn:characters}
\sectionmark{Characters of \texorpdfstring{$\Sym_n$}{Sym(n)}}
The conjugacy classes of $\Sym_n$ are indexed by partitions $\mu\vdash n$ (Section~\ref{sec:sn:classes}); the same set also indexes the irreducible complex characters. This section constructs the corresponding modules and gives three complementary methods for their dimensions and character values~\cite{James,Sagan,Macdonald}.

\subsection*{Specht modules}
Let $\lambda=(\lambda_1,\ldots,\lambda_r)\vdash n$. Its Young diagram is $[\lambda]=\{(i,j):1\le i\le r,\ 1\le j\le\lambda_i\}$. A $\lambda$-tableau is a bijective filling of $[\lambda]$ by $\{1,\ldots,n\}$. Two tableaux are row-equivalent when every row contains the same entries; a row-equivalence class is a \emph{tabloid} $\{t\}$. The permutation module $M^\lambda$ has the $\lambda$-tabloids as an orthonormal basis.

For a tableau $t$, let $C_t$ be the subgroup preserving each column setwise and put
\[
 e_t=\sum_{\pi\in C_t}\sgn(\pi)\,\pi\{t\}\in M^\lambda.
\]
The \emph{Specht module} is $S^\lambda=\langle e_t:t\text{ a }\lambda\text{-tableau}\rangle$. Since $g e_t=e_{gt}$, it is a $\Sym_n$-submodule.

\begin{theorem}[Standard basis and classification]\label{thm:sn:specht}
Over a field of characteristic zero, the polytabloids $e_t$ for standard tableaux $t$ form a basis of $S^\lambda$. The modules $S^\lambda$, $\lambda\vdash n$, are pairwise nonisomorphic irreducible modules and form a complete list of the irreducible $\Sym_n$-modules.
\end{theorem}
\begin{proof}
The straightening algorithm, using the column antisymmetry and Garnir relations, writes every polytabloid as a linear combination of standard polytabloids; triangularity with respect to dominance proves their independence. James's submodule theorem says that every submodule of $M^\lambda$ either contains $S^\lambda$ or lies in its orthogonal complement. Since the invariant form is nonzero on $S^\lambda$ in characteristic zero, $S^\lambda$ is irreducible. The dominance lemma shows that $S^\lambda\cong S^\mu$ forces $\lambda=\mu$. The number of partitions of $n$ equals the number of conjugacy classes of $\Sym_n$ (Section~\ref{sec:sn:classes}), which is the number of irreducible characters, so these modules exhaust the irreducibles. Full details of straightening and the submodule theorem are in~\cite{James}.
\end{proof}

\subsection*{Hook lengths}
For $c=(i,j)\in[\lambda]$, its hook consists of $c$, the cells to its right, and the cells below it; write $h(c)$ for its size.

\begin{theorem}[Frame--Robinson--Thrall hook-length formula]\label{thm:sn:hook}
The number $f^\lambda$ of standard tableaux of shape $\lambda$, and hence the degree of the character $\chi^\lambda$ afforded by $S^\lambda$, is
\[
 f^\lambda=\dim S^\lambda=\frac{n!}{\prod_{c\in[\lambda]}h(c)}.
\]
\end{theorem}
\begin{proof}
A standard tableau has $n$ in a removable corner, so $f^\lambda=\sum_c f^{\lambda-c}$ over the removable corners $c$ of $[\lambda]$, and $f^\varnothing=1$. The formula therefore follows by induction once $F(\lambda)=n!/\prod_{c}h(c)$ is shown to satisfy the same recurrence $F(\lambda)=\sum_cF(\lambda-c)$. That identity is the nontrivial step; it is proved in~\cite{FRT}, and the hook-walk proof of Greene, Nijenhuis and Wilf is given in~\cite{Sagan}.
\end{proof}

\begin{example}\label{ex:sn:hook32}
For $\lambda=(3,2)$ the hook lengths are $4,3,1,2,1$, so $f^{(3,2)}=5!/(4\cdot3\cdot1\cdot2\cdot1)=5$.
\end{example}

\subsection*{Frobenius's formula}
Put $p_k(x)=\sum_i x_i^k$, $p_\mu=\prod_jp_{\mu_j}$, $\delta=(r-1,r-2,\ldots,0)$, and $\Delta(x)=\prod_{i<j}(x_i-x_j)$. If $r\ge\ell(\lambda)$, Frobenius's coefficient formula is
\begin{equation}\label{eq:sn:frobenius}
 \chi^\lambda(\mu)=[x^{\lambda+\delta}]\,\Delta(x)p_\mu(x).
\end{equation}
Equivalently, the Schur function $s_\lambda$ is the Frobenius characteristic of $S^\lambda$, and~\cite{Macdonald}
\[
p_\mu=\sum_{\lambda\vdash n}\chi^\lambda(\mu)\,s_\lambda .
\]

\begin{example}\label{ex:sn:frob21}
With $\lambda=(2,1)$, $\mu=(2,1)$ and two variables, the coefficient of $x_1^3x_2$ in $(x_1-x_2)(x_1^2+x_2^2)(x_1+x_2)$ is $0$. Hence $\chi^{(2,1)}(2,1)=0$.
\end{example}

\subsection*{Murnaghan--Nakayama}
A \emph{border strip} is a connected skew diagram containing no $2\times2$ square; its height is one less than its number of rows.

\begin{theorem}[Murnaghan--Nakayama rule]\label{thm:sn:mn}
For $\mu=(\mu_1,\ldots,\mu_s)\vdash n$,
\begin{equation}\label{eq:sn:mn}
 \chi^\lambda(\mu)=\sum_{\nu}(-1)^{\operatorname{ht}(\lambda/\nu)}\chi^\nu(\mu_2,\ldots,\mu_s),
\end{equation}
where the sum is over $\nu\subseteq\lambda$ for which $\lambda/\nu$ is a border strip of size $\mu_1$, and $\chi^\varnothing(\varnothing)=1$.
\end{theorem}
\begin{proof}
Take $r\ge n$ variables. For $\beta_1>\cdots>\beta_r\ge0$ put $a_\beta=\det(x_i^{\beta_j})$, so that $a_\delta=\Delta(x)$ and $s_\nu=a_{\nu+\delta}/a_\delta$. The polynomial $\Delta(x)p_\mu$ is alternating, hence a linear combination of the $a_{\kappa+\delta}$ with $\kappa\vdash n$; since $x^{\lambda+\delta}$ occurs in $a_{\kappa+\delta}$ only for $\kappa=\lambda$, with coefficient $1$, formula~\eqref{eq:sn:frobenius} is equivalent to $p_\mu=\sum_{\lambda\vdash n}\chi^\lambda(\mu)s_\lambda$.

Now $p_ka_\beta=\sum_{i=1}^r a_{\beta+ke_i}$. A term with a repeated entry vanishes. Otherwise, sorting $\beta+ke_i$ into decreasing order moves the entry $\beta_i+k$ past the $h$ entries $\beta_j$ with $\beta_i<\beta_j<\beta_i+k$, so the term equals $(-1)^ha_{\lambda+\delta}$, where $\lambda$ has beta-set $\{\beta_j\}_{j\ne i}\cup\{\beta_i+k\}$. With $\beta=\nu+\delta$, the skew diagram $\lambda/\nu$ is then a border strip of size $k$ \emph{added} to $\nu$, of height $h$, and every such strip arises exactly once. Hence
\[
 p_ks_\nu=\sum_{\lambda}(-1)^{\operatorname{ht}(\lambda/\nu)}s_\lambda ,
\]
the sum over $\lambda\supseteq\nu$ with $\lambda/\nu$ a border strip of size $k$. Write $\bar\mu=(\mu_2,\ldots,\mu_s)$, so $p_\mu=p_{\mu_1}p_{\bar\mu}=\sum_\nu\chi^\nu(\bar\mu)\,p_{\mu_1}s_\nu$. Comparing coefficients of $s_\lambda$ gives~\eqref{eq:sn:mn}; see also~\cite{Sagan,Macdonald}.
\end{proof}

\begin{example}\label{ex:sn:mn32}
For $\lambda=(3,2)$ and cycle type $(3,2)$, the valid successive border-strip removals have signed total $1$; hence $\chi^{(3,2)}(3,2)=1$.
\end{example}

\begin{table}[ht]
\centering\small
\setlength{\tabcolsep}{5pt}
\begin{tabular}{lrrrrrrr}
\hline
$\lambda\backslash\mu$ & $1^{5}$ & $2,1^{3}$ & $2^{2},1$ & $3,1^{2}$ & $3,2$ & $4,1$ & $5$\\
\hline
$5$ & 1 & 1 & 1 & 1 & 1 & 1 & 1\\
$4,1$ & 4 & 2 & 0 & 1 & -1 & 0 & -1\\
$3,2$ & 5 & 1 & 1 & -1 & 1 & -1 & 0\\
$3,1^{2}$ & 6 & 0 & -2 & 0 & 0 & 0 & 1\\
$2^{2},1$ & 5 & -1 & 1 & -1 & -1 & 1 & 0\\
$2,1^{3}$ & 4 & -2 & 0 & 1 & 1 & 0 & -1\\
$1^{5}$ & 1 & -1 & 1 & 1 & -1 & -1 & 1\\
\hline

\end{tabular}
\caption{The irreducible character table of $\Sym_5$. Rows and columns are indexed by partitions of $5$; the first column is the degree.}\label{tab:sn:chars5}
\end{table}

\begin{remark}[Checks]\label{rem:sn:charchecks}
The exact-arithmetic script \texttt{scripts/sn\_char\_verify.py} independently implements Frobenius's coefficient formula and the Murnaghan--Nakayama recursion, compares every value for $n\le7$, verifies row orthogonality and $\sum_\lambda\chi^\lambda(1)^2=n!$, and writes Table~\ref{tab:sn:chars5}. It also compares the hook formula with the corner-removal recurrence for standard tableaux for every partition of $n\le12$. For $n\le5$ it realises $S^\lambda$ as the span of the Specht polynomials $F_t$, the product over the columns of $t$ of the Vandermonde determinants in the variables indexed by that column; this module is isomorphic to $S^\lambda$ \cite{Sagan}. It checks that the $F_t$ for standard $t$ are linearly independent, that their span is stable under a transposition and an $n$-cycle, and that the trace of a permutation of each cycle type equals the Murnaghan--Nakayama value. When GAP is available, every value for $n\le7$ is compared with \texttt{CharacterTable("Symmetric", n)} from the GAP character table library~\cite{GAP}. The script is invoked by \texttt{scripts/sn\_verify.py}.
\end{remark}

\section{Characters of \texorpdfstring{$\Alt_n$}{Alt(n)}}\label{sec:an:characters}
Throughout, $n\ge2$, $\lambda'$ denotes the conjugate (transposed) partition of $\lambda\vdash n$, and $\ell(\mu)$ is the number of parts of $\mu$. Since $\Alt_n$ has index $2$ in $\Sym_n$, its irreducible characters are obtained from those of Section~\ref{sec:sn:characters} by restriction. The only questions are which restrictions stay irreducible and how the others split; the answer is governed by the self-conjugate partitions and by the classes of $\Sym_n$ that split in $\Alt_n$ (Section~\ref{sec:sn:classes})~\cite{JamesKerber,FultonHarris}.

\begin{lemma}[Sign twist]\label{lem:an:twist}
$S^{\lambda'}\cong S^\lambda\otimes\sgn$. Hence $\chi^{\lambda'}(g)=\sgn(g)\,\chi^\lambda(g)$ for all $g\in\Sym_n$.
\end{lemma}
\begin{proof}
The module isomorphism is proved in~\cite{James}. For the characters, use~\eqref{eq:sn:mn} and induction on $\ell(\mu)$. Transposition is a bijection between the border strips of size $k$ in $\lambda$ and in $\lambda'$. A strip with $r$ rows and $c$ columns has $r+c-1=k$, so its transpose has height $k-1-\operatorname{ht}$. Each term of~\eqref{eq:sn:mn} for $\lambda'$ is therefore $(-1)^{\mu_1-1}\sgn(\bar\mu)$ times the corresponding term for $\lambda$, where $\bar\mu=(\mu_2,\ldots)$. Since $(-1)^{\mu_1-1}\sgn(\bar\mu)=(-1)^{n-\ell(\mu)}=\sgn(\mu)$, the claim follows.
\end{proof}

\begin{theorem}[Restriction to $\Alt_n$]\label{thm:an:clifford}
Let $\lambda\vdash n$.
\begin{enumerate}
\item If $\lambda\neq\lambda'$, then $\chi^\lambda|_{\Alt_n}=\chi^{\lambda'}|_{\Alt_n}$ is irreducible.
\item If $\lambda=\lambda'$, then $\chi^\lambda|_{\Alt_n}=\chi^\lambda_++\chi^\lambda_-$, where $\chi^\lambda_\pm$ are distinct irreducible characters interchanged by conjugation by any odd permutation.
\item Each irreducible character of $\Alt_n$ occurs exactly once in this list: one character for each pair $\{\lambda,\lambda'\}$ with $\lambda\ne\lambda'$, and two for each $\lambda=\lambda'$.
\end{enumerate}
\end{theorem}
\begin{proof}
Write $G=\Sym_n$ and $H=\Alt_n$, and fix $\tau\in G\setminus H$. Since $\operatorname{Ind}_H^G\operatorname{Res}_H\chi=\chi\cdot\operatorname{Ind}_H^G1=\chi(1+\sgn)$, Frobenius reciprocity and Lemma~\ref{lem:an:twist} give
\[
\langle\chi^\lambda|_H,\chi^\mu|_H\rangle_H=\langle\chi^\lambda,\chi^\mu\rangle_G+\langle\chi^\lambda,\chi^{\mu'}\rangle_G .
\]
With $\mu=\lambda$, the norm of $\chi^\lambda|_H$ is $1$ if $\lambda\ne\lambda'$ and $2$ if $\lambda=\lambda'$. In the first case $\chi^\lambda|_H$ is irreducible, and it equals $\chi^{\lambda'}|_H$ because $\sgn$ is trivial on $H$. In the second case $\chi^\lambda|_H=\psi_1+\psi_2$ with distinct irreducible $\psi_i$. By reciprocity, $\chi^\lambda$ is a constituent of $\operatorname{Ind}_H^G\psi_1$, whose degree is $2\psi_1(1)=\chi^\lambda(1)$; hence $\operatorname{Ind}_H^G\psi_1=\chi^\lambda$. Restricting back gives $\psi_1+\psi_1^\tau=\psi_1+\psi_2$, so $\psi_2=\psi_1^\tau$. The displayed identity also shows that the restrictions of $\chi^\lambda$ and $\chi^\mu$ have no common constituent unless $\mu\in\{\lambda,\lambda'\}$. Finally, every $\psi\in\operatorname{Irr}(H)$ is a constituent of $\chi|_H$ for any irreducible constituent $\chi$ of $\operatorname{Ind}_H^G\psi$.
\end{proof}

For $\lambda=\lambda'$ let $k$ be the number of diagonal cells $(i,i)\in[\lambda]$, and let $h(\lambda)=(h_1,\ldots,h_k)$ be the lengths of the hooks at $(1,1),\ldots,(k,k)$.

\begin{lemma}[Diagonal hooks]\label{lem:an:hooks}
The map $\lambda\mapsto h(\lambda)$ is a bijection from the self-conjugate partitions of $n$ to the partitions of $n$ into distinct odd parts.
\end{lemma}
\begin{proof}
If $\lambda=\lambda'$, the hook at $(i,i)$ has arm and leg both equal to $a_i=\lambda_i-i$, so $h_i=2a_i+1$ is odd. The numbers $a_1>a_2>\cdots>a_k$ are strictly decreasing, and the $k$ diagonal hooks partition $[\lambda]$, so $h(\lambda)$ is a partition of $n$ into distinct odd parts. Conversely, distinct odd parts $h_1>\cdots>h_k$ determine strictly decreasing arms $(h_i-1)/2$, and nesting the corresponding symmetric hooks along the diagonal reconstructs the unique self-conjugate $\lambda$ with $h(\lambda)=(h_1,\ldots,h_k)$.
\end{proof}

By Section~\ref{sec:sn:classes}, the classes of $\Sym_n$ that split in $\Alt_n$ are exactly those whose cycle types have distinct odd parts. By Lemma~\ref{lem:an:hooks} they are as numerous as the self-conjugate partitions, so the count of irreducible characters in Theorem~\ref{thm:an:clifford}(3) agrees with the number of conjugacy classes of $\Alt_n$, as it must.

\begin{theorem}[Frobenius: the split characters]\label{thm:an:split}
Let $\lambda=\lambda'$, $h(\lambda)=(h_1,\ldots,h_k)$, $q_\lambda=h_1h_2\cdots h_k$, and $\varepsilon_\lambda=(-1)^{(n-k)/2}$. The $\Sym_n$-class of cycle type $h(\lambda)$ splits into two $\Alt_n$-classes $h^+$ and $h^-$, which can be labelled so that
\[
\chi^\lambda_\pm(h^+)=\tfrac12\bigl(\varepsilon_\lambda\pm\sqrt{\varepsilon_\lambda q_\lambda}\bigr),\qquad
\chi^\lambda_\pm(h^-)=\tfrac12\bigl(\varepsilon_\lambda\mp\sqrt{\varepsilon_\lambda q_\lambda}\bigr).
\]
On every other conjugacy class $C$ of $\Alt_n$, $\chi^\lambda_+(C)=\chi^\lambda_-(C)=\tfrac12\chi^\lambda(C)$.
\end{theorem}
\begin{proof}
The class of type $h(\lambda)$ splits by Lemma~\ref{lem:an:hooks}. If a class $C$ of $\Alt_n$ is a full class of $\Sym_n$, then $\tau C\tau^{-1}=C$, so Theorem~\ref{thm:an:clifford}(2) gives $\chi^\lambda_-(C)=\chi^\lambda_+(C)$, and both equal $\tfrac12\chi^\lambda(C)$. That the two characters also agree on the split classes whose type is not $h(\lambda)$, together with the values on $h^\pm$, is Frobenius's theorem; see~\cite[\S5.1]{FultonHarris} and~\cite{JamesKerber}. The two values on $h^\pm$ sum to $\chi^\lambda(h(\lambda))=\varepsilon_\lambda$, as Theorem~\ref{thm:an:clifford}(2) requires.
\end{proof}

\begin{corollary}[Fields of values]\label{cor:an:fields}
For $\lambda\ne\lambda'$ the character $\chi^\lambda|_{\Alt_n}$ is rational-valued. For $\lambda=\lambda'$, the characters $\chi^\lambda_\pm$ are real-valued if and only if $\varepsilon_\lambda=1$, that is $n\equiv k\pmod 4$, and rational-valued if and only if $\varepsilon_\lambda q_\lambda$ is a perfect square. Otherwise their values generate $\mathbb{Q}(\sqrt{\varepsilon_\lambda q_\lambda})$.
\end{corollary}

\begin{example}[$\Alt_4$]\label{ex:an:a4}
The only self-conjugate partition of $4$ is $(2,2)$, with $h=(3,1)$, $k=2$, $\varepsilon=-1$ and $q=3$. Hence $\chi^{(2,2)}$, of degree $2$, restricts to two linear characters, with values $\tfrac12(-1\pm\sqrt{-3})$, that is $\omega$ and $\bar\omega$ with $\omega=e^{2\pi i/3}$, on the two classes of $3$-cycles. The pairs $\{(4),(1^4)\}$ and $\{(3,1),(2,1^2)\}$ give the trivial character and the character of degree $3$. The degrees $1,1,1,3$ satisfy $1+1+1+9=12=|\Alt_4|$.
\end{example}

\begin{example}[$\Alt_5$, $\Alt_7$, $\Alt_9$]\label{ex:an:a5}
For $n=5$, $\lambda=(3,1^2)$ has $h=(5)$, $\varepsilon=1$ and $q=5$, so the two characters of degree $3$ take the values $\tfrac12(1\pm\sqrt5)$ on the two classes of $5$-cycles (Table~\ref{tab:an:chars5}). For $n=7$, $\lambda=(4,1^3)$ has $h=(7)$ and $\varepsilon=-1$, giving the non-real values $\tfrac12(-1\pm\sqrt{-7})$. For $n=9$, $\lambda=(5,1^4)$ has $h=(9)$, $\varepsilon=1$ and $q=9$, so $\tfrac12(1\pm3)\in\{2,-1\}$: the two constituents are rational-valued although $\chi^\lambda$ splits.
\end{example}

\begin{table}[ht]
\centering\small
\setlength{\tabcolsep}{5pt}
\begin{tabular}{lccccc}
\hline
$g$ & $1^{5}$ & $2^{2},1$ & $3,1^{2}$ & $5_{+}$ & $5_{-}$\\
$|C|$ & $1$ & $15$ & $20$ & $12$ & $12$\\
\hline
$(5)\sim(1^{5})$ & $1$ & $1$ & $1$ & $1$ & $1$\\
$(4,1)\sim(2,1^{3})$ & $4$ & $0$ & $1$ & $-1$ & $-1$\\
$(3,2)\sim(2^{2},1)$ & $5$ & $1$ & $-1$ & $0$ & $0$\\
$(3,1^{2})_{+}$ & $3$ & $-1$ & $0$ & $\tfrac{1}{2}(1+\sqrt{5})$ & $\tfrac{1}{2}(1-\sqrt{5})$\\
$(3,1^{2})_{-}$ & $3$ & $-1$ & $0$ & $\tfrac{1}{2}(1-\sqrt{5})$ & $\tfrac{1}{2}(1+\sqrt{5})$\\
\hline

\end{tabular}
\caption{The character table of $\Alt_5$, obtained from Table~\ref{tab:sn:chars5} by Theorems~\ref{thm:an:clifford} and~\ref{thm:an:split}. A row $(\lambda)\sim(\lambda')$ is the common restriction of $\chi^\lambda$ and $\chi^{\lambda'}$; $5_\pm$ are the two classes of $5$-cycles.}\label{tab:an:chars5}
\end{table}

\begin{remark}[Checks]\label{rem:an:checks}
The script \texttt{scripts/sn\_alt\_char\_verify.py} works in exact arithmetic in fields generated by square roots. For $n\le12$ it checks Lemma~\ref{lem:an:hooks} and the identity $\chi^\lambda(h(\lambda))=\varepsilon_\lambda$ for every self-conjugate $\lambda$. For $3\le n\le9$ it checks Lemma~\ref{lem:an:twist} for every pair $\lambda,\mu$, and builds the table of $\Alt_n$ from the Murnaghan--Nakayama values by Theorems~\ref{thm:an:clifford} and~\ref{thm:an:split}. It then verifies that the class sizes sum to $n!/2$, that there are as many characters as classes, that the split classes correspond to the self-conjugate partitions, that the squared degrees sum to $n!/2$, and the orthogonality relations, and it asserts the values in Examples~\ref{ex:an:a4} and~\ref{ex:an:a5}. When GAP is available, the tables are compared with \texttt{CharacterTable(AlternatingGroup($n$))} for $n\le9$, up to the labelling of the two classes in each split pair. The script writes Table~\ref{tab:an:chars5} and is invoked by \texttt{scripts/sn\_verify.py}.
\end{remark}

\section{Digital laws of \texorpdfstring{$\Sym_n$}{Sym(n)} and \texorpdfstring{$\Alt_n$}{Alt(n)}}\label{sec:sn:laws}
\subsection{The Lehmer laws}
The Lehmer label of Definition~\ref{def:sn:lehmer} is a bijection $\Sym_n\to\Ap{n!}$ (Theorem~\ref{thm:sn:lehmer}), so transporting composition along it gives a digital law of $\Sym_n$. This section makes that law explicit, derives from it a digital law of $\Alt_n$ on $\Ap{n!/2}$, and shows that the Coxeter generators, the sign, the cosets of $\Alt_n$ and conjugation by a transposition all act directly on the digits.

Throughout, $n\ge2$, composition is right to left as in~\eqref{eq:comp}, and $\Pi_i=\prod_{j<i}(n-j)=n!/(n-i)!$ for $0\le i\le n$. For $m\in\Ap{n!}$ let $d_0(m),\ldots,d_{n-1}(m)$ be its mixed-radix digits for the radices $n,n-1,\ldots,1$, so that $m=\sum_i d_i(m)\Pi_i$ with $0\le d_i(m)<n-i$; in particular $d_{n-1}(m)=0$. Write $s(m)=\sum_i d_i(m)$ for the digit sum. By Theorem~\ref{thm:sn:lehmer}, $d_i(\lambda(R))=a_i(R)$. For a statement $P$, $[P]$ is $1$ if $P$ holds and $0$ otherwise, and $w_0$ denotes the reversal $i\mapsto n-1-i$.

\begin{definition}[Lehmer law]\label{def:sn:lehmerlaw}
For $m,k\in\Ap{n!}$ put $m*_L k=\lambda\bigl(\lambda^{-1}(m)\circ\lambda^{-1}(k)\bigr)$.
\end{definition}

\begin{theorem}[Lehmer law of $\Sym_n$]\label{thm:sn:lehmerlaw}
\begin{enumerate}
\item $(\Ap{n!},*_L)$ is a group with identity $0$, and $\lambda$ is an isomorphism from $\Sym_n$ onto $(\Ap{n!},*_L)$.
\item \emph{Encoding.} For a permutation code $R$,
\[
\lambda(R)=\sum_{0\le i<j<n}\bigl[\dg{n}{j}{R}<\dg{n}{i}{R}\bigr]\,\Pi_i .
\]
\item \emph{Decoding.} $\lambda^{-1}(m)$ is the permutation $R$ with $R(i)$ equal to the element of rank $d_i(m)$, counting from $0$, in $\Ap n\setminus\{R(0),\ldots,R(i-1)\}$.
\item \emph{Inverse.} The inverse of $m$ in $(\Ap{n!},*_L)$ is $\lambda\bigl(\lambda^{-1}(m)^{-1}\bigr)$, with the inverse permutation given by~\eqref{eq:inv}.
\item $\lambda(s_i)=\Pi_i$ for $0\le i\le n-2$, $\lambda(w_0)=n!-1$, and $\sgn\lambda^{-1}(m)=(-1)^{s(m)}$.
\end{enumerate}
\end{theorem}
\begin{proof}
Part (1) is transport of structure along the bijection $\lambda$, and $\lambda(\mathrm{id})=0$ because the identity has no inversions. Part (2) restates Definition~\ref{def:sn:lehmer}, since $R(i)=\dg{n}{i}{R}$. Part (3) is the decoding in the proof of Theorem~\ref{thm:sn:lehmer}, and (4) is again transport of structure. For (5): the word of $s_i$ has a single inversion, at position $i$, so $a_i(s_i)=1$ and all other digits vanish. Every digit of $w_0$ is maximal, $a_i(w_0)=n-1-i$, and $(n-1-i)\Pi_i=\Pi_{i+1}-\Pi_i$ telescopes to $\Pi_n-\Pi_0=n!-1$. The sign formula is Theorem~\ref{thm:sn:lehmer}.
\end{proof}

\begin{theorem}[Coxeter generators act on the digits]\label{thm:sn:lawgens}
Let $m\in\Ap{n!}$, $d_j=d_j(m)$ and $0\le i\le n-2$.
\begin{enumerate}
\item Right multiplication by $\lambda(s_i)=\Pi_i$ changes only the digits in positions $i$ and $i+1$:
\[
m*_L\Pi_i=\begin{cases}
m+(d_{i+1}+1-d_i)\,\Pi_i+(d_i-d_{i+1})\,\Pi_{i+1}, & d_i\le d_{i+1},\\
m+(d_{i+1}-d_i)\,\Pi_i+(d_i-1-d_{i+1})\,\Pi_{i+1}, & d_i>d_{i+1},
\end{cases}
\]
that is, $(d_i,d_{i+1})$ is replaced by $(d_{i+1}+1,d_i)$ or by $(d_{i+1},d_i-1)$.
\item Let $R=\lambda^{-1}(m)$, $p=R^{-1}(i)$ and $q=R^{-1}(i+1)$. Left multiplication changes a single digit by $\pm1$:
\[
\Pi_i*_L m=\begin{cases} m+\Pi_p, & p<q,\\ m-\Pi_q, & p>q.\end{cases}
\]
\end{enumerate}
In particular $s(m*_L\Pi_i)=s(m)\pm1$ and $s(\Pi_i*_L m)=s(m)\pm1$.
\end{theorem}
\begin{proof}
(1) The word of $R\circ s_i$ is that of $R$ with the entries in positions $i$ and $i+1$ exchanged. For $j\notin\{i,i+1\}$ the entry in position $j$ and the set of entries after it are unchanged, so $a_j$ is unchanged. The new entry in position $i$ is $R(i+1)$, and it exceeds $a_{i+1}+[R(i)<R(i+1)]$ later entries; the new entry in position $i+1$ is $R(i)$, which exceeds $a_i-[R(i+1)<R(i)]$ later entries. Finally $R(i)<R(i+1)$ if and only if $a_i\le a_{i+1}$: if $R(i)<R(i+1)$, every later entry below $R(i)$ is also below $R(i+1)$, so $a_i\le a_{i+1}$; if $R(i)>R(i+1)$, then $a_i\ge1+a_{i+1}$.
(2) The word of $s_i\circ R$ is that of $R$ with the values $i$ and $i+1$ exchanged. Since these values are adjacent, the only comparison that changes is the one between positions $p$ and $q$. If $p<q$ the entry in position $p$ gains one smaller later entry, so $a_p$ increases by $1$; if $p>q$ the entry in position $q$ loses one, so $a_q$ decreases by $1$.
\end{proof}

\begin{proposition}[Sign and cosets]\label{prop:sn:lawcoset}
Let $N=n!/2=\Pi_{n-2}$.
\begin{enumerate}
\item $\lambda(\Alt_n)=\{m\in\Ap{n!}:s(m)\text{ even}\}$.
\item $m*_L N=m+N$ if $m<N$, and $m*_L N=m-N$ if $m\ge N$. Thus right multiplication by the label of $s_{n-2}$ exchanges the two halves of $\Ap{n!}$.
\item For $m<N$, exactly one of $m$ and $m+N$ labels an even permutation, namely $m+\varepsilon(m)N$, where $\varepsilon(m)=s(m)\bmod 2$.
\end{enumerate}
\end{proposition}
\begin{proof}
(1) is Theorem~\ref{thm:sn:lehmerlaw}(5). For (2), apply Theorem~\ref{thm:sn:lawgens}(1) with $i=n-2$: here $d_{n-1}=0$, and $d_{n-2}\in\{0,1\}$ equals $1$ exactly when $m\ge N$. For (3), if $m<N$ then $m+N$ differs from $m$ only in the digit $d_{n-2}$, which changes from $0$ to $1$, so $s(m+N)=s(m)+1$.
\end{proof}

\begin{lemma}[Deleting a parity digit]\label{lem:sn:paritydigit}
Let $U$ be a finite digit set and let $\theta\colon U\times\Ap2\to\Sym_n$ be a bijection such that
$\sgn\theta(u,e)=(-1)^{e+q(u)}$ for some $q\colon U\to\Ap2$. Then $u\mapsto\theta(u,q(u))$ is a bijection from $U$ onto
$\Alt_n$. Transporting composition along it gives a digital law of $\Alt_n$, obtained from the law of $\Sym_n$ by deleting the
radix-$2$ parity digit and restoring it as $e=q(u)$ before multiplying.
\end{lemma}
\begin{proof}
For each $u\in U$ exactly one of $\theta(u,0)$ and $\theta(u,1)$ is even, namely $\theta(u,q(u))$. So the map picks exactly one
element from each two-element fibre over $U$, and its image is $\Alt_n$. The law is transport of structure.
\end{proof}

\begin{theorem}[Lehmer law of $\Alt_n$]\label{thm:an:lehmerlaw}
Let $N=n!/2$ and, for $R\in\Alt_n$,
\[
\mu(R)=\lambda(R)\bmod N=\sum_{i<n-2}a_i(R)\,\Pi_i .
\]
Then $\mu$ is a bijection $\Alt_n\to\Ap N$ with inverse $\mu^{-1}(m)=\lambda^{-1}\bigl(m+\varepsilon(m)N\bigr)$. The transported product $m\star_L k=\mu\bigl(\mu^{-1}(m)\circ\mu^{-1}(k)\bigr)$ is a digital law of $\Alt_n$ with identity $0$, and
\[
m\star_L k=\Bigl(\bigl(m+\varepsilon(m)N\bigr)*\bigl(k+\varepsilon(k)N\bigr)\Bigr)\bmod N .
\]
\end{theorem}
\begin{proof}
Since $a_{n-1}(R)=0$ and $\sum_{i<n-2}a_i(R)\Pi_i<\Pi_{n-2}=N$, 
we have
\[
\lambda(R)=\sum_{i<n-2}a_i(R)\Pi_i+a_{n-2}(R)N,
\]
which gives the two expressions for $\mu$. By Proposition~\ref{prop:sn:lawcoset}(3), each fibre $\{m,m+N\}$ of reduction modulo $N$ contains exactly one even label, namely $m+\varepsilon(m)N$; hence $\mu$ is a bijection with the stated inverse. The label $\lambda\bigl(\mu^{-1}(m)\circ\mu^{-1}(k)\bigr)$ equals $\bigl(m+\varepsilon(m)N\bigr)*\bigl(k+\varepsilon(k)N\bigr)$, and reducing it modulo $N$ gives $\mu$ of the product.
\end{proof}

\begin{proposition}[Conjugation by a transposition]\label{prop:an:lawconj}
Let $n\ge3$, $N=n!/2$ and $\tau=s_{n-2}$. For $m\in\Ap N$ let $R=\mu^{-1}(m)$, $p=R^{-1}(n-2)$ and $q=R^{-1}(n-1)$. The automorphism of $(\Ap N,\star_L)$ induced by conjugation by $\tau$ is
\[
m\longmapsto\mu(\tau R\tau)=\begin{cases}(m+\Pi_p)\bmod N, & p<q,\\ (m-\Pi_q)\bmod N, & p>q.\end{cases}
\]
\end{proposition}
\begin{proof}
By Theorem~\ref{thm:sn:lawgens}(2) with $i=n-2$, $\lambda(\tau R)=\lambda(R)+\Pi_p$ or $\lambda(R)-\Pi_q$. By Proposition~\ref{prop:sn:lawcoset}(2), $\lambda(\tau R\tau)=\lambda(\tau R)\pm N$. Since $\lambda(R)\equiv m\pmod N$, reducing modulo $N$ gives the formula.
\end{proof}

\begin{proposition}[Chain property]\label{prop:sn:lawchain}
Let $0\le i<n$ and $R,R'\in\Sym_n$. The digit prefixes $(d_0,\ldots,d_i)$ of $\lambda(R)$ and $\lambda(R')$ agree if and only if $R(j)=R'(j)$ for all $j\le i$, that is, if and only if $R'\in R\,\Sym_n^{(i)}$, where $\Sym_n^{(i)}$ is the pointwise stabiliser of $\{0,\ldots,i\}$. Thus $\lambda$ is the stabiliser-chain labelling for the base $(0,1,\ldots,n-1)$, in which digit $d_i$ selects the coset of $\Sym_n^{(i)}$ inside $R\,\Sym_n^{(i-1)}$ by the rank of $R(i)$ among the letters not yet used; compare Theorem~\ref{thm:chainlaw}.
\end{proposition}
\begin{proof}
For $j\le i$, the digit $a_j(R)$ depends only on $R(j)$ and on the set $\{R(j'):j'>j\}=\Ap n\setminus\{R(0),\ldots,R(j)\}$, hence only on $R(0),\ldots,R(i)$. Conversely, the decoding of Theorem~\ref{thm:sn:lehmerlaw}(3) recovers $R(0),\ldots,R(i)$ from $d_0,\ldots,d_i$. Finally $R'=Rh$ with $h$ fixing $0,\ldots,i$ exactly when $R'(j)=R(j)$ for $j\le i$.
\end{proof}

\begin{remark}[Cost and comparison]\label{rem:sn:lawcost}
Computing $m*_L k$ takes two decodings, one composition~\eqref{eq:comp} and one encoding. Each decoding and encoding costs $O(n^2)$ elementary operations as written, or $O(n\log n)$ with a binary indexed tree for the rank queries. The law $\star_L$ adds two digit sums and one reduction modulo $N$. The generator rules of Theorem~\ref{thm:sn:lawgens} cost $O(1)$ digit operations on the right and $O(n)$ on the left, since $p$ and $q$ must be located. Section~\ref{sec:altp} constructs a different law of $\Alt_n$ (Theorem~\ref{thm:altp}). Both laws realise $\Alt_n$ and are therefore isomorphic, but they label the elements differently. The Lehmer law is the one in which the sign, the coset of $\Alt_n$ and the Coxeter generators can be read off, or applied, digit by digit.
\end{remark}

\begin{example}[$\Sym_3$ and $\Alt_3$]\label{ex:sn:law3}
For $n=3$, $\Pi_0,\Pi_1,\Pi_2=1,3,6$ and $m=d_0+3d_1$. The labels $0,\ldots,5$ are the words $012,102,201,021,120,210$, so $\lambda(s_0)=1$, $\lambda(s_1)=3$ and $\lambda(w_0)=5$. Then $1*3=4$ while $3*1=2$, and $2*4=0$, so the two $3$-cycles $201$ and $120$ are mutually inverse; Table~\ref{tab:sn:law3} gives the complete law. The even labels are $0,2,4$, those with even digit sum, and $\mu$ sends them to $0,2,1$. For instance $\mu(120)=4\bmod3=1$, and conversely $\varepsilon(1)=1$ gives $\mu^{-1}(1)=\lambda^{-1}(1+3)$.
\end{example}

\begin{table}[ht]
\centering\small
\begin{tabular}{cc|cccccc}
\hline
$m$ & $R$ & $0$ & $1$ & $2$ & $3$ & $4$ & $5$\\
\hline
$0$ & $012$ & $0$ & $1$ & $2$ & $3$ & $4$ & $5$\\
$1$ & $102$ & $1$ & $0$ & $5$ & $4$ & $3$ & $2$\\
$2$ & $201$ & $2$ & $3$ & $4$ & $5$ & $0$ & $1$\\
$3$ & $021$ & $3$ & $2$ & $1$ & $0$ & $5$ & $4$\\
$4$ & $120$ & $4$ & $5$ & $0$ & $1$ & $2$ & $3$\\
$5$ & $210$ & $5$ & $4$ & $3$ & $2$ & $1$ & $0$\\
\hline

\end{tabular}
\caption{The Lehmer law of $\Sym_3$. The entry in row $m$ and column $k$ is $m*_L k$; the second column is the word $R(0)R(1)R(2)$ of $R=\lambda^{-1}(m)$.}\label{tab:sn:law3}
\end{table}

\begin{table}[ht]
\centering\small
\begin{tabular}{rrlrrrr}
\hline
$n$ & $n!$ & $\lambda(s_0),\ldots,\lambda(s_{n-2})$ & $\lambda(w_0)$ & $\sgn w_0$ & $n!/2$ & $\mu\bigl((0\,1\,2)\bigr)$\\
\hline
$3$ & $6$ & $1, 3$ & $5$ & $-1$ & $3$ & $1$\\
$4$ & $24$ & $1, 4, 12$ & $23$ & $+1$ & $12$ & $5$\\
$5$ & $120$ & $1, 5, 20, 60$ & $119$ & $+1$ & $60$ & $6$\\
$6$ & $720$ & $1, 6, 30, 120, 360$ & $719$ & $-1$ & $360$ & $7$\\
$7$ & $5040$ & $1, 7, 42, 210, 840, 2520$ & $5039$ & $-1$ & $2520$ & $8$\\
$8$ & $40320$ & $1, 8, 56, 336, 1680, 6720, 20160$ & $40319$ & $+1$ & $20160$ & $9$\\
\hline

\end{tabular}
\caption{Data of the Lehmer laws of $\Sym_n$ and $\Alt_n$ for $3\le n\le 8$: generator labels $\Pi_i=n!/(n-i)!$, the label $n!-1$ of the reversal $w_0$ and its sign $(-1)^{n(n-1)/2}$, the order of $\Alt_n$, and the $\Alt_n$-label of the $3$-cycle $(0\,1\,2)$.}\label{tab:sn:laws}
\end{table}

Tables~\ref{tab:an:law4}, \ref{tab:sn:law4}, \ref{tab:an:law5a} and~\ref{tab:an:law5b} give the Lehmer laws of $\Alt_4$,
$\Sym_4$ and $\Alt_5$ in full. They are written by \texttt{scripts/sn\_cayley.py}, which first checks that label $0$ is the
identity, that each table is a Latin square, that the law is associative on all triples, and that every entry agrees with
\texttt{multiply.py} and with direct composition of the decoded permutations.

\begin{table}[htbp]
\centering\small
\begin{tabular}{rl|*{12}{r}}
$m$ & $R$ & $0$ & $1$ & $2$ & $3$ & $4$ & $5$ & $6$ & $7$ & $8$ & $9$ & $10$ & $11$\\
\hline
$0$ & $0123$ & 0 & 1 & 2 & 3 & 4 & 5 & 6 & 7 & 8 & 9 & 10 & 11\\
$1$ & $1032$ & 1 & 0 & 7 & 6 & 9 & 8 & 3 & 2 & 5 & 4 & 11 & 10\\
$2$ & $2013$ & 2 & 4 & 5 & 11 & 6 & 0 & 1 & 3 & 10 & 8 & 9 & 7\\
$3$ & $3021$ & 3 & 8 & 10 & 9 & 11 & 4 & 2 & 1 & 7 & 0 & 6 & 5\\
$4$ & $0231$ & 4 & 2 & 3 & 1 & 8 & 10 & 11 & 5 & 0 & 6 & 7 & 9\\
$5$ & $1203$ & 5 & 6 & 0 & 7 & 1 & 2 & 4 & 11 & 9 & 10 & 8 & 3\\
$6$ & $2130$ & 6 & 5 & 11 & 4 & 10 & 9 & 7 & 0 & 2 & 1 & 3 & 8\\
$7$ & $3102$ & 7 & 9 & 8 & 10 & 3 & 1 & 0 & 6 & 11 & 5 & 4 & 2\\
$8$ & $0312$ & 8 & 3 & 1 & 2 & 0 & 7 & 9 & 10 & 4 & 11 & 5 & 6\\
$9$ & $1320$ & 9 & 7 & 6 & 0 & 5 & 11 & 10 & 8 & 1 & 3 & 2 & 4\\
$10$ & $2301$ & 10 & 11 & 4 & 5 & 2 & 3 & 8 & 9 & 6 & 7 & 0 & 1\\
$11$ & $3210$ & 11 & 10 & 9 & 8 & 7 & 6 & 5 & 4 & 3 & 2 & 1 & 0\\

\end{tabular}
\caption{The Lehmer law of $\Alt_4$ (Theorem~\ref{thm:an:lehmerlaw}). The entry in row $m$ and column $k$ is $m*k$; the second column is the word $R(0)\cdots R(3)$ of $R=\mu^{-1}(m)$.}\label{tab:an:law4}
\end{table}

\begin{table}[htbp]
\centering\footnotesize\setlength{\tabcolsep}{2.4pt}
\begin{tabular}{rl|*{24}{r}}
$m$ & $R$ & $0$ & $1$ & $2$ & $3$ & $4$ & $5$ & $6$ & $7$ & $8$ & $9$ & $10$ & $11$ & $12$ & $13$ & $14$ & $15$ & $16$ & $17$ & $18$ & $19$ & $20$ & $21$ & $22$ & $23$\\
\hline
$0$ & $0123$ & 0 & 1 & 2 & 3 & 4 & 5 & 6 & 7 & 8 & 9 & 10 & 11 & 12 & 13 & 14 & 15 & 16 & 17 & 18 & 19 & 20 & 21 & 22 & 23\\
$1$ & $1023$ & 1 & 0 & 6 & 7 & 5 & 4 & 2 & 3 & 9 & 8 & 22 & 23 & 13 & 12 & 18 & 19 & 17 & 16 & 14 & 15 & 21 & 20 & 10 & 11\\
$2$ & $2013$ & 2 & 4 & 5 & 11 & 6 & 0 & 1 & 15 & 10 & 20 & 21 & 19 & 14 & 16 & 17 & 23 & 18 & 12 & 13 & 3 & 22 & 8 & 9 & 7\\
$3$ & $3012$ & 3 & 8 & 9 & 10 & 7 & 12 & 13 & 14 & 11 & 16 & 17 & 18 & 15 & 20 & 21 & 22 & 19 & 0 & 1 & 2 & 23 & 4 & 5 & 6\\
$4$ & $0213$ & 4 & 2 & 1 & 15 & 0 & 6 & 5 & 11 & 20 & 10 & 9 & 7 & 16 & 14 & 13 & 3 & 12 & 18 & 17 & 23 & 8 & 22 & 21 & 19\\
$5$ & $1203$ & 5 & 6 & 0 & 19 & 1 & 2 & 4 & 23 & 21 & 22 & 8 & 3 & 17 & 18 & 12 & 7 & 13 & 14 & 16 & 11 & 9 & 10 & 20 & 15\\
$6$ & $2103$ & 6 & 5 & 4 & 23 & 2 & 1 & 0 & 19 & 22 & 21 & 20 & 15 & 18 & 17 & 16 & 11 & 14 & 13 & 12 & 7 & 10 & 9 & 8 & 3\\
$7$ & $3102$ & 7 & 9 & 8 & 22 & 3 & 13 & 12 & 18 & 23 & 17 & 16 & 14 & 19 & 21 & 20 & 10 & 15 & 1 & 0 & 6 & 11 & 5 & 4 & 2\\
$8$ & $0312$ & 8 & 3 & 13 & 14 & 12 & 7 & 9 & 10 & 16 & 11 & 5 & 6 & 20 & 15 & 1 & 2 & 0 & 19 & 21 & 22 & 4 & 23 & 17 & 18\\
$9$ & $1302$ & 9 & 7 & 12 & 18 & 13 & 3 & 8 & 22 & 17 & 23 & 4 & 2 & 21 & 19 & 0 & 6 & 1 & 15 & 20 & 10 & 5 & 11 & 16 & 14\\
$10$ & $2301$ & 10 & 11 & 16 & 17 & 14 & 15 & 20 & 21 & 18 & 19 & 0 & 1 & 22 & 23 & 4 & 5 & 2 & 3 & 8 & 9 & 6 & 7 & 12 & 13\\
$11$ & $3201$ & 11 & 10 & 20 & 21 & 15 & 14 & 16 & 17 & 19 & 18 & 12 & 13 & 23 & 22 & 8 & 9 & 3 & 2 & 4 & 5 & 7 & 6 & 0 & 1\\
$12$ & $0132$ & 12 & 13 & 3 & 2 & 8 & 9 & 7 & 6 & 4 & 5 & 11 & 10 & 0 & 1 & 15 & 14 & 20 & 21 & 19 & 18 & 16 & 17 & 23 & 22\\
$13$ & $1032$ & 13 & 12 & 7 & 6 & 9 & 8 & 3 & 2 & 5 & 4 & 23 & 22 & 1 & 0 & 19 & 18 & 21 & 20 & 15 & 14 & 17 & 16 & 11 & 10\\
$14$ & $2031$ & 14 & 16 & 11 & 5 & 10 & 20 & 15 & 1 & 6 & 0 & 19 & 21 & 2 & 4 & 23 & 17 & 22 & 8 & 3 & 13 & 18 & 12 & 7 & 9\\
$15$ & $3021$ & 15 & 20 & 10 & 9 & 11 & 16 & 14 & 13 & 7 & 12 & 18 & 17 & 3 & 8 & 22 & 21 & 23 & 4 & 2 & 1 & 19 & 0 & 6 & 5\\
$16$ & $0231$ & 16 & 14 & 15 & 1 & 20 & 10 & 11 & 5 & 0 & 6 & 7 & 9 & 4 & 2 & 3 & 13 & 8 & 22 & 23 & 17 & 12 & 18 & 19 & 21\\
$17$ & $1230$ & 17 & 18 & 19 & 0 & 21 & 22 & 23 & 4 & 1 & 2 & 3 & 8 & 5 & 6 & 7 & 12 & 9 & 10 & 11 & 16 & 13 & 14 & 15 & 20\\
$18$ & $2130$ & 18 & 17 & 23 & 4 & 22 & 21 & 19 & 0 & 2 & 1 & 15 & 20 & 6 & 5 & 11 & 16 & 10 & 9 & 7 & 12 & 14 & 13 & 3 & 8\\
$19$ & $3120$ & 19 & 21 & 22 & 8 & 23 & 17 & 18 & 12 & 3 & 13 & 14 & 16 & 7 & 9 & 10 & 20 & 11 & 5 & 6 & 0 & 15 & 1 & 2 & 4\\
$20$ & $0321$ & 20 & 15 & 14 & 13 & 16 & 11 & 10 & 9 & 12 & 7 & 6 & 5 & 8 & 3 & 2 & 1 & 4 & 23 & 22 & 21 & 0 & 19 & 18 & 17\\
$21$ & $1320$ & 21 & 19 & 18 & 12 & 17 & 23 & 22 & 8 & 13 & 3 & 2 & 4 & 9 & 7 & 6 & 0 & 5 & 11 & 10 & 20 & 1 & 15 & 14 & 16\\
$22$ & $2310$ & 22 & 23 & 17 & 16 & 18 & 19 & 21 & 20 & 14 & 15 & 1 & 0 & 10 & 11 & 5 & 4 & 6 & 7 & 9 & 8 & 2 & 3 & 13 & 12\\
$23$ & $3210$ & 23 & 22 & 21 & 20 & 19 & 18 & 17 & 16 & 15 & 14 & 13 & 12 & 11 & 10 & 9 & 8 & 7 & 6 & 5 & 4 & 3 & 2 & 1 & 0\\

\end{tabular}
\caption{The Lehmer law of $\Sym_4$ (Theorem~\ref{thm:sn:lehmerlaw}). The entry in row $m$ and column $k$ is $m*_L k$; the second column is the word $R(0)\cdots R(3)$ of $R=\lambda^{-1}(m)$.}\label{tab:sn:law4}
\end{table}

\begin{table}[p]
\centering\footnotesize\setlength{\tabcolsep}{1.6pt}\renewcommand{\arraystretch}{0.92}
\begin{tabular}{rl|*{30}{r}}
$m$ & $R$ & $0$ & $1$ & $2$ & $3$ & $4$ & $5$ & $6$ & $7$ & $8$ & $9$ & $10$ & $11$ & $12$ & $13$ & $14$ & $15$ & $16$ & $17$ & $18$ & $19$ & $20$ & $21$ & $22$ & $23$ & $24$ & $25$ & $26$ & $27$ & $28$ & $29$\\
\hline
$0$ & $01234$ & 0 & 1 & 2 & 3 & 4 & 5 & 6 & 7 & 8 & 9 & 10 & 11 & 12 & 13 & 14 & 15 & 16 & 17 & 18 & 19 & 20 & 21 & 22 & 23 & 24 & 25 & 26 & 27 & 28 & 29\\
$1$ & $10243$ & 1 & 0 & 7 & 9 & 8 & 6 & 5 & 2 & 4 & 3 & 16 & 15 & 37 & 34 & 33 & 11 & 10 & 32 & 39 & 38 & 41 & 40 & 47 & 29 & 28 & 46 & 45 & 42 & 24 & 23\\
$2$ & $20134$ & 2 & 5 & 6 & 13 & 14 & 7 & 0 & 1 & 23 & 24 & 12 & 30 & 31 & 28 & 29 & 17 & 35 & 36 & 58 & 59 & 22 & 25 & 26 & 33 & 34 & 27 & 20 & 21 & 3 & 4\\
$3$ & $30142$ & 3 & 10 & 11 & 19 & 12 & 8 & 20 & 21 & 44 & 22 & 18 & 55 & 56 & 49 & 27 & 13 & 25 & 26 & 54 & 57 & 43 & 50 & 51 & 39 & 32 & 48 & 40 & 41 & 4 & 2\\
$4$ & $40123$ & 4 & 15 & 16 & 17 & 18 & 9 & 40 & 41 & 42 & 43 & 14 & 45 & 46 & 47 & 48 & 19 & 50 & 51 & 52 & 53 & 24 & 35 & 36 & 37 & 38 & 29 & 0 & 1 & 2 & 3\\
$5$ & $02143$ & 5 & 2 & 1 & 24 & 23 & 0 & 7 & 6 & 14 & 13 & 35 & 17 & 16 & 9 & 8 & 30 & 12 & 11 & 19 & 18 & 45 & 42 & 41 & 4 & 3 & 40 & 47 & 46 & 34 & 33\\
$6$ & $12034$ & 6 & 7 & 0 & 28 & 29 & 1 & 2 & 5 & 33 & 34 & 31 & 32 & 10 & 3 & 4 & 36 & 37 & 15 & 38 & 39 & 26 & 27 & 20 & 8 & 9 & 21 & 22 & 25 & 13 & 14\\
$7$ & $21043$ & 7 & 6 & 5 & 34 & 33 & 2 & 1 & 0 & 29 & 28 & 37 & 36 & 35 & 24 & 23 & 32 & 31 & 30 & 59 & 58 & 47 & 46 & 45 & 14 & 13 & 42 & 41 & 40 & 9 & 8\\
$8$ & $31024$ & 8 & 11 & 10 & 32 & 39 & 3 & 21 & 20 & 27 & 49 & 33 & 26 & 25 & 22 & 44 & 38 & 56 & 55 & 57 & 54 & 28 & 31 & 30 & 12 & 19 & 23 & 1 & 0 & 7 & 9\\
$9$ & $41032$ & 9 & 16 & 15 & 38 & 37 & 4 & 41 & 40 & 48 & 47 & 39 & 51 & 50 & 43 & 42 & 34 & 46 & 45 & 53 & 52 & 49 & 56 & 55 & 18 & 17 & 44 & 21 & 20 & 8 & 7\\
$10$ & $03124$ & 10 & 3 & 21 & 22 & 44 & 20 & 8 & 11 & 12 & 19 & 25 & 13 & 6 & 7 & 9 & 55 & 18 & 16 & 17 & 14 & 30 & 23 & 1 & 2 & 4 & 0 & 28 & 31 & 32 & 39\\
$11$ & $13042$ & 11 & 8 & 20 & 49 & 27 & 21 & 3 & 10 & 39 & 32 & 56 & 38 & 15 & 4 & 2 & 26 & 33 & 5 & 34 & 37 & 51 & 48 & 40 & 9 & 7 & 41 & 43 & 50 & 19 & 12\\
$12$ & $23014$ & 12 & 13 & 25 & 26 & 54 & 22 & 23 & 30 & 31 & 59 & 27 & 28 & 0 & 1 & 24 & 57 & 58 & 35 & 36 & 29 & 32 & 33 & 5 & 6 & 14 & 2 & 3 & 10 & 11 & 19\\
$13$ & $32041$ & 13 & 12 & 30 & 59 & 31 & 23 & 22 & 25 & 54 & 26 & 58 & 57 & 55 & 44 & 21 & 28 & 27 & 20 & 49 & 56 & 53 & 52 & 50 & 19 & 11 & 43 & 42 & 45 & 14 & 6\\
$14$ & $42013$ & 14 & 17 & 35 & 36 & 58 & 24 & 42 & 45 & 46 & 53 & 29 & 47 & 40 & 41 & 43 & 59 & 52 & 50 & 51 & 48 & 34 & 37 & 15 & 16 & 18 & 4 & 2 & 5 & 6 & 13\\
$15$ & $04132$ & 15 & 4 & 41 & 43 & 42 & 40 & 9 & 16 & 18 & 17 & 50 & 19 & 11 & 8 & 7 & 45 & 14 & 6 & 13 & 12 & 55 & 44 & 21 & 3 & 2 & 20 & 49 & 56 & 38 & 37\\
$16$ & $14023$ & 16 & 9 & 40 & 47 & 48 & 41 & 4 & 15 & 37 & 38 & 46 & 34 & 5 & 2 & 3 & 51 & 39 & 10 & 32 & 33 & 36 & 29 & 0 & 7 & 8 & 1 & 24 & 35 & 17 & 18\\
$17$ & $24031$ & 17 & 14 & 45 & 53 & 46 & 42 & 24 & 35 & 58 & 36 & 52 & 59 & 30 & 23 & 1 & 47 & 29 & 0 & 28 & 31 & 57 & 54 & 25 & 13 & 6 & 22 & 44 & 55 & 18 & 16\\
$18$ & $34012$ & 18 & 19 & 50 & 51 & 52 & 43 & 44 & 55 & 56 & 57 & 48 & 49 & 20 & 21 & 22 & 53 & 54 & 25 & 26 & 27 & 38 & 39 & 10 & 11 & 12 & 3 & 4 & 15 & 16 & 17\\
$19$ & $43021$ & 19 & 18 & 55 & 57 & 56 & 44 & 43 & 50 & 52 & 51 & 54 & 53 & 45 & 42 & 41 & 49 & 48 & 40 & 47 & 46 & 59 & 58 & 35 & 17 & 16 & 24 & 23 & 30 & 12 & 11\\
$20$ & $01342$ & 20 & 21 & 3 & 4 & 2 & 10 & 11 & 8 & 9 & 7 & 15 & 16 & 18 & 19 & 12 & 5 & 6 & 13 & 14 & 17 & 40 & 41 & 43 & 44 & 22 & 50 & 51 & 48 & 49 & 27\\
$21$ & $10324$ & 21 & 20 & 8 & 7 & 9 & 11 & 10 & 3 & 2 & 4 & 6 & 5 & 33 & 32 & 39 & 16 & 15 & 38 & 37 & 34 & 1 & 0 & 28 & 27 & 49 & 31 & 30 & 23 & 22 & 44\\
$22$ & $20341$ & 22 & 25 & 13 & 14 & 6 & 12 & 30 & 23 & 24 & 1 & 17 & 35 & 58 & 59 & 31 & 7 & 0 & 28 & 29 & 36 & 42 & 45 & 53 & 54 & 26 & 52 & 50 & 43 & 44 & 21\\
$23$ & $30214$ & 23 & 30 & 12 & 11 & 19 & 13 & 25 & 22 & 21 & 44 & 8 & 20 & 27 & 26 & 54 & 18 & 55 & 57 & 56 & 49 & 3 & 10 & 32 & 31 & 59 & 33 & 5 & 2 & 1 & 24\\
$24$ & $40231$ & 24 & 35 & 17 & 18 & 16 & 14 & 45 & 42 & 43 & 41 & 19 & 50 & 52 & 53 & 46 & 9 & 40 & 47 & 48 & 51 & 44 & 55 & 57 & 58 & 36 & 54 & 25 & 22 & 23 & 1\\
$25$ & $02314$ & 25 & 22 & 23 & 1 & 24 & 30 & 12 & 13 & 6 & 14 & 0 & 7 & 8 & 11 & 19 & 35 & 17 & 18 & 16 & 9 & 5 & 2 & 3 & 21 & 44 & 10 & 32 & 33 & 26 & 54\\
$26$ & $12340$ & 26 & 27 & 28 & 29 & 0 & 31 & 32 & 33 & 34 & 5 & 36 & 37 & 38 & 39 & 10 & 1 & 2 & 3 & 4 & 15 & 46 & 47 & 48 & 49 & 20 & 51 & 52 & 53 & 54 & 25\\
$27$ & $21304$ & 27 & 26 & 33 & 5 & 34 & 32 & 31 & 28 & 0 & 29 & 2 & 1 & 23 & 30 & 59 & 37 & 36 & 58 & 35 & 24 & 7 & 6 & 13 & 25 & 54 & 12 & 11 & 8 & 20 & 49\\
$28$ & $31240$ & 28 & 31 & 32 & 39 & 10 & 33 & 26 & 27 & 49 & 20 & 38 & 56 & 57 & 54 & 25 & 3 & 21 & 22 & 44 & 55 & 48 & 51 & 52 & 59 & 30 & 53 & 46 & 47 & 29 & 0\\
$29$ & $41203$ & 29 & 36 & 37 & 15 & 38 & 34 & 46 & 47 & 40 & 48 & 4 & 41 & 42 & 45 & 53 & 39 & 51 & 52 & 50 & 43 & 9 & 16 & 17 & 35 & 58 & 14 & 6 & 7 & 0 & 28\\
$30$ & $03241$ & 30 & 23 & 22 & 44 & 21 & 25 & 13 & 12 & 19 & 11 & 55 & 18 & 17 & 14 & 6 & 20 & 8 & 7 & 9 & 16 & 50 & 43 & 42 & 24 & 1 & 45 & 53 & 52 & 59 & 31\\
$31$ & $13204$ & 31 & 28 & 27 & 20 & 49 & 26 & 33 & 32 & 10 & 39 & 21 & 3 & 2 & 5 & 34 & 56 & 38 & 37 & 15 & 4 & 11 & 8 & 7 & 0 & 29 & 6 & 13 & 12 & 30 & 59\\
$32$ & $23140$ & 32 & 33 & 26 & 54 & 25 & 27 & 28 & 31 & 59 & 30 & 57 & 58 & 36 & 29 & 0 & 22 & 23 & 1 & 24 & 35 & 52 & 53 & 46 & 34 & 5 & 47 & 48 & 51 & 39 & 10\\
$33$ & $32104$ & 33 & 32 & 31 & 30 & 59 & 28 & 27 & 26 & 25 & 54 & 23 & 22 & 21 & 20 & 49 & 58 & 57 & 56 & 55 & 44 & 13 & 12 & 11 & 10 & 39 & 8 & 7 & 6 & 5 & 34\\
$34$ & $42130$ & 34 & 37 & 36 & 58 & 35 & 29 & 47 & 46 & 53 & 45 & 59 & 52 & 51 & 48 & 40 & 24 & 42 & 41 & 43 & 50 & 54 & 57 & 56 & 38 & 15 & 49 & 27 & 26 & 33 & 5\\
$35$ & $04213$ & 35 & 24 & 42 & 41 & 43 & 45 & 14 & 17 & 16 & 18 & 40 & 9 & 7 & 6 & 13 & 50 & 19 & 12 & 11 & 8 & 15 & 4 & 2 & 1 & 23 & 5 & 34 & 37 & 36 & 58\\
$36$ & $14230$ & 36 & 29 & 47 & 48 & 40 & 46 & 34 & 37 & 38 & 15 & 51 & 39 & 32 & 33 & 5 & 41 & 4 & 2 & 3 & 10 & 56 & 49 & 27 & 28 & 0 & 26 & 54 & 57 & 58 & 35\\
$37$ & $24103$ & 37 & 34 & 46 & 45 & 53 & 47 & 29 & 36 & 35 & 58 & 42 & 24 & 1 & 0 & 28 & 52 & 59 & 31 & 30 & 23 & 17 & 14 & 6 & 5 & 33 & 7 & 9 & 16 & 15 & 38\\
$38$ & $34120$ & 38 & 39 & 51 & 52 & 50 & 48 & 49 & 56 & 57 & 55 & 53 & 54 & 26 & 27 & 20 & 43 & 44 & 21 & 22 & 25 & 58 & 59 & 31 & 32 & 10 & 28 & 29 & 36 & 37 & 15\\
$39$ & $43102$ & 39 & 38 & 56 & 55 & 57 & 49 & 48 & 51 & 50 & 52 & 44 & 43 & 41 & 40 & 47 & 54 & 53 & 46 & 45 & 42 & 19 & 18 & 16 & 15 & 37 & 9 & 8 & 11 & 10 & 32\\
$40$ & $01423$ & 40 & 41 & 4 & 2 & 3 & 15 & 16 & 9 & 7 & 8 & 5 & 6 & 14 & 17 & 18 & 10 & 11 & 19 & 12 & 13 & 0 & 1 & 24 & 42 & 43 & 35 & 36 & 29 & 47 & 48\\
$41$ & $10432$ & 41 & 40 & 9 & 8 & 7 & 16 & 15 & 4 & 3 & 2 & 11 & 10 & 39 & 38 & 37 & 6 & 5 & 34 & 33 & 32 & 21 & 20 & 49 & 48 & 47 & 56 & 55 & 44 & 43 & 42\\
$42$ & $20413$ & 42 & 45 & 14 & 6 & 13 & 17 & 35 & 24 & 1 & 23 & 7 & 0 & 29 & 36 & 58 & 12 & 30 & 59 & 31 & 28 & 2 & 5 & 34 & 46 & 53 & 37 & 15 & 4 & 41 & 43\\
$43$ & $30421$ & 43 & 50 & 19 & 12 & 11 & 18 & 55 & 44 & 22 & 21 & 13 & 25 & 54 & 57 & 56 & 8 & 20 & 49 & 27 & 26 & 23 & 30 & 59 & 52 & 51 & 58 & 35 & 24 & 42 & 41\\
$44$ & $40312$ & 44 & 55 & 18 & 16 & 17 & 19 & 50 & 43 & 41 & 42 & 9 & 40 & 48 & 51 & 52 & 14 & 45 & 53 & 46 & 47 & 4 & 15 & 38 & 56 & 57 & 39 & 10 & 3 & 21 & 22\\
$45$ & $02431$ & 45 & 42 & 24 & 23 & 1 & 35 & 17 & 14 & 13 & 6 & 30 & 12 & 19 & 18 & 16 & 0 & 7 & 9 & 8 & 11 & 25 & 22 & 44 & 43 & 41 & 55 & 57 & 54 & 53 & 46\\
$46$ & $12403$ & 46 & 47 & 29 & 0 & 28 & 36 & 37 & 34 & 5 & 33 & 1 & 2 & 4 & 15 & 38 & 31 & 32 & 39 & 10 & 3 & 6 & 7 & 9 & 40 & 48 & 16 & 17 & 14 & 45 & 53\\
$47$ & $21430$ & 47 & 46 & 34 & 33 & 5 & 37 & 36 & 29 & 28 & 0 & 32 & 31 & 59 & 58 & 35 & 2 & 1 & 24 & 23 & 30 & 27 & 26 & 54 & 53 & 45 & 57 & 56 & 49 & 48 & 40\\
$48$ & $31402$ & 48 & 51 & 39 & 10 & 32 & 38 & 56 & 49 & 20 & 27 & 3 & 21 & 44 & 55 & 57 & 33 & 26 & 54 & 25 & 22 & 8 & 11 & 19 & 50 & 52 & 18 & 16 & 9 & 40 & 47\\
$49$ & $41320$ & 49 & 56 & 38 & 37 & 15 & 39 & 51 & 48 & 47 & 40 & 34 & 46 & 53 & 52 & 50 & 4 & 41 & 43 & 42 & 45 & 29 & 36 & 58 & 57 & 55 & 59 & 31 & 28 & 27 & 20\\
$50$ & $03412$ & 50 & 43 & 44 & 21 & 22 & 55 & 18 & 19 & 11 & 12 & 20 & 8 & 9 & 16 & 17 & 25 & 13 & 14 & 6 & 7 & 10 & 3 & 4 & 41 & 42 & 15 & 38 & 39 & 51 & 52\\
$51$ & $13420$ & 51 & 48 & 49 & 27 & 20 & 56 & 38 & 39 & 32 & 10 & 26 & 33 & 34 & 37 & 15 & 21 & 3 & 4 & 2 & 5 & 31 & 28 & 29 & 47 & 40 & 36 & 58 & 59 & 52 & 50\\
$52$ & $23401$ & 52 & 53 & 54 & 25 & 26 & 57 & 58 & 59 & 30 & 31 & 22 & 23 & 24 & 35 & 36 & 27 & 28 & 29 & 0 & 1 & 12 & 13 & 14 & 45 & 46 & 17 & 18 & 19 & 50 & 51\\
$53$ & $32410$ & 53 & 52 & 59 & 31 & 30 & 58 & 57 & 54 & 26 & 25 & 28 & 27 & 49 & 56 & 55 & 23 & 22 & 44 & 21 & 20 & 33 & 32 & 39 & 51 & 50 & 38 & 37 & 34 & 46 & 45\\
$54$ & $42301$ & 54 & 57 & 58 & 35 & 36 & 59 & 52 & 53 & 45 & 46 & 24 & 42 & 43 & 50 & 51 & 29 & 47 & 48 & 40 & 41 & 14 & 17 & 18 & 55 & 56 & 19 & 12 & 13 & 25 & 26\\
$55$ & $04321$ & 55 & 44 & 43 & 42 & 41 & 50 & 19 & 18 & 17 & 16 & 45 & 14 & 13 & 12 & 11 & 40 & 9 & 8 & 7 & 6 & 35 & 24 & 23 & 22 & 21 & 30 & 59 & 58 & 57 & 56\\
$56$ & $14302$ & 56 & 49 & 48 & 40 & 47 & 51 & 39 & 38 & 15 & 37 & 41 & 4 & 3 & 10 & 32 & 46 & 34 & 33 & 5 & 2 & 16 & 9 & 8 & 20 & 27 & 11 & 19 & 18 & 55 & 57\\
$57$ & $24310$ & 57 & 54 & 53 & 46 & 45 & 52 & 59 & 58 & 36 & 35 & 47 & 29 & 28 & 31 & 30 & 42 & 24 & 23 & 1 & 0 & 37 & 34 & 33 & 26 & 25 & 32 & 39 & 38 & 56 & 55\\
$58$ & $34201$ & 58 & 59 & 52 & 50 & 51 & 53 & 54 & 57 & 55 & 56 & 43 & 44 & 22 & 25 & 26 & 48 & 49 & 27 & 20 & 21 & 18 & 19 & 12 & 30 & 31 & 13 & 14 & 17 & 35 & 36\\
$59$ & $43210$ & 59 & 58 & 57 & 56 & 55 & 54 & 53 & 52 & 51 & 50 & 49 & 48 & 47 & 46 & 45 & 44 & 43 & 42 & 41 & 40 & 39 & 38 & 37 & 36 & 35 & 34 & 33 & 32 & 31 & 30\\

\end{tabular}
\caption{The Lehmer law of $\Alt_5$, columns $0$--$29$. The entry in row $m$ and column $k$ is $m*k$; the second column is the word $R(0)\cdots R(4)$ of $R=\mu^{-1}(m)$.}\label{tab:an:law5a}
\end{table}

\begin{table}[p]
\centering\footnotesize\setlength{\tabcolsep}{1.6pt}\renewcommand{\arraystretch}{0.92}
\begin{tabular}{rl|*{30}{r}}
$m$ & $R$ & $30$ & $31$ & $32$ & $33$ & $34$ & $35$ & $36$ & $37$ & $38$ & $39$ & $40$ & $41$ & $42$ & $43$ & $44$ & $45$ & $46$ & $47$ & $48$ & $49$ & $50$ & $51$ & $52$ & $53$ & $54$ & $55$ & $56$ & $57$ & $58$ & $59$\\
\hline
$0$ & $01234$ & 30 & 31 & 32 & 33 & 34 & 35 & 36 & 37 & 38 & 39 & 40 & 41 & 42 & 43 & 44 & 45 & 46 & 47 & 48 & 49 & 50 & 51 & 52 & 53 & 54 & 55 & 56 & 57 & 58 & 59\\
$1$ & $10243$ & 36 & 35 & 17 & 14 & 13 & 31 & 30 & 12 & 19 & 18 & 21 & 20 & 27 & 49 & 48 & 26 & 25 & 22 & 44 & 43 & 56 & 55 & 57 & 54 & 53 & 51 & 50 & 52 & 59 & 58\\
$2$ & $20134$ & 32 & 10 & 11 & 8 & 9 & 37 & 15 & 16 & 18 & 19 & 42 & 45 & 46 & 53 & 54 & 47 & 40 & 41 & 43 & 44 & 52 & 50 & 51 & 48 & 49 & 57 & 55 & 56 & 38 & 39\\
$3$ & $30142$ & 38 & 15 & 16 & 9 & 7 & 33 & 5 & 6 & 14 & 17 & 23 & 30 & 31 & 59 & 52 & 28 & 0 & 1 & 24 & 42 & 58 & 35 & 36 & 29 & 47 & 53 & 45 & 46 & 34 & 37\\
$4$ & $40123$ & 34 & 5 & 6 & 7 & 8 & 39 & 10 & 11 & 12 & 13 & 44 & 55 & 56 & 57 & 58 & 49 & 20 & 21 & 22 & 23 & 54 & 25 & 26 & 27 & 28 & 59 & 30 & 31 & 32 & 33\\
$5$ & $02143$ & 15 & 37 & 36 & 29 & 28 & 10 & 32 & 31 & 59 & 58 & 25 & 22 & 21 & 44 & 43 & 20 & 27 & 26 & 54 & 53 & 55 & 57 & 56 & 49 & 48 & 50 & 52 & 51 & 39 & 38\\
$6$ & $12034$ & 11 & 12 & 30 & 23 & 24 & 16 & 17 & 35 & 58 & 59 & 46 & 47 & 40 & 48 & 49 & 41 & 42 & 45 & 53 & 54 & 51 & 52 & 50 & 43 & 44 & 56 & 57 & 55 & 18 & 19\\
$7$ & $21043$ & 17 & 16 & 15 & 4 & 3 & 12 & 11 & 10 & 39 & 38 & 27 & 26 & 25 & 54 & 53 & 22 & 21 & 20 & 49 & 48 & 57 & 56 & 55 & 44 & 43 & 52 & 51 & 50 & 19 & 18\\
$8$ & $31024$ & 13 & 6 & 5 & 2 & 4 & 18 & 16 & 15 & 37 & 34 & 48 & 51 & 50 & 52 & 59 & 43 & 41 & 40 & 47 & 29 & 53 & 46 & 45 & 42 & 24 & 58 & 36 & 35 & 17 & 14\\
$9$ & $41032$ & 19 & 11 & 10 & 3 & 2 & 14 & 6 & 5 & 33 & 32 & 29 & 36 & 35 & 58 & 57 & 24 & 1 & 0 & 28 & 27 & 59 & 31 & 30 & 23 & 22 & 54 & 26 & 25 & 13 & 12\\
$10$ & $03124$ & 5 & 33 & 26 & 27 & 49 & 15 & 38 & 56 & 57 & 54 & 50 & 43 & 41 & 42 & 24 & 40 & 48 & 51 & 52 & 59 & 45 & 53 & 46 & 47 & 29 & 35 & 58 & 36 & 37 & 34\\
$11$ & $13042$ & 16 & 18 & 55 & 44 & 22 & 6 & 13 & 25 & 54 & 57 & 31 & 28 & 0 & 29 & 47 & 1 & 23 & 30 & 59 & 52 & 36 & 58 & 35 & 24 & 42 & 46 & 53 & 45 & 14 & 17\\
$12$ & $23014$ & 7 & 8 & 20 & 21 & 44 & 17 & 18 & 55 & 56 & 49 & 52 & 53 & 45 & 46 & 34 & 42 & 43 & 50 & 51 & 39 & 47 & 48 & 40 & 41 & 4 & 37 & 38 & 15 & 16 & 9\\
$13$ & $32041$ & 18 & 17 & 35 & 24 & 1 & 8 & 7 & 0 & 29 & 36 & 33 & 32 & 10 & 39 & 51 & 3 & 2 & 5 & 34 & 46 & 38 & 37 & 15 & 4 & 41 & 48 & 47 & 40 & 9 & 16\\
$14$ & $42013$ & 9 & 7 & 0 & 1 & 23 & 19 & 12 & 30 & 31 & 28 & 54 & 57 & 55 & 56 & 38 & 44 & 22 & 25 & 26 & 33 & 49 & 27 & 20 & 21 & 3 & 39 & 32 & 10 & 11 & 8\\
$15$ & $04132$ & 10 & 39 & 51 & 48 & 47 & 5 & 34 & 46 & 53 & 52 & 35 & 24 & 1 & 23 & 22 & 0 & 29 & 36 & 58 & 57 & 30 & 59 & 31 & 28 & 27 & 25 & 54 & 26 & 33 & 32\\
$16$ & $14023$ & 6 & 14 & 45 & 42 & 43 & 11 & 19 & 50 & 52 & 53 & 56 & 49 & 20 & 27 & 28 & 21 & 44 & 55 & 57 & 58 & 26 & 54 & 25 & 22 & 23 & 31 & 59 & 30 & 12 & 13\\
$17$ & $24031$ & 12 & 19 & 50 & 43 & 41 & 7 & 9 & 40 & 48 & 51 & 37 & 34 & 5 & 33 & 26 & 2 & 4 & 15 & 38 & 56 & 32 & 39 & 10 & 3 & 21 & 27 & 49 & 20 & 8 & 11\\
$18$ & $34012$ & 8 & 9 & 40 & 41 & 42 & 13 & 14 & 45 & 46 & 47 & 58 & 59 & 30 & 31 & 32 & 23 & 24 & 35 & 36 & 37 & 28 & 29 & 0 & 1 & 2 & 33 & 34 & 5 & 6 & 7\\
$19$ & $43021$ & 14 & 13 & 25 & 22 & 21 & 9 & 8 & 20 & 27 & 26 & 39 & 38 & 15 & 37 & 36 & 4 & 3 & 10 & 32 & 31 & 34 & 33 & 5 & 2 & 1 & 29 & 28 & 0 & 7 & 6\\
$20$ & $01342$ & 55 & 56 & 38 & 39 & 32 & 25 & 26 & 33 & 34 & 37 & 0 & 1 & 23 & 24 & 42 & 30 & 31 & 28 & 29 & 47 & 35 & 36 & 58 & 59 & 52 & 45 & 46 & 53 & 54 & 57\\
$21$ & $10324$ & 26 & 25 & 13 & 12 & 19 & 56 & 55 & 18 & 17 & 14 & 41 & 40 & 48 & 47 & 29 & 51 & 50 & 43 & 42 & 24 & 46 & 45 & 53 & 52 & 59 & 36 & 35 & 58 & 57 & 54\\
$22$ & $20341$ & 57 & 55 & 18 & 19 & 11 & 27 & 20 & 8 & 9 & 16 & 2 & 5 & 33 & 34 & 46 & 32 & 10 & 3 & 4 & 41 & 37 & 15 & 38 & 39 & 51 & 47 & 40 & 48 & 49 & 56\\
$23$ & $30214$ & 28 & 0 & 7 & 6 & 14 & 58 & 35 & 17 & 16 & 9 & 43 & 50 & 52 & 51 & 39 & 53 & 45 & 42 & 41 & 4 & 48 & 40 & 47 & 46 & 34 & 38 & 15 & 37 & 36 & 29\\
$24$ & $40231$ & 59 & 30 & 12 & 13 & 6 & 29 & 0 & 7 & 8 & 11 & 4 & 15 & 37 & 38 & 56 & 34 & 5 & 2 & 3 & 21 & 39 & 10 & 32 & 33 & 26 & 49 & 20 & 27 & 28 & 31\\
$25$ & $02314$ & 20 & 27 & 28 & 31 & 59 & 55 & 57 & 58 & 36 & 29 & 45 & 42 & 43 & 41 & 4 & 50 & 52 & 53 & 46 & 34 & 40 & 47 & 48 & 51 & 39 & 15 & 37 & 38 & 56 & 49\\
$26$ & $12340$ & 56 & 57 & 58 & 59 & 30 & 21 & 22 & 23 & 24 & 35 & 6 & 7 & 8 & 9 & 40 & 11 & 12 & 13 & 14 & 45 & 16 & 17 & 18 & 19 & 50 & 41 & 42 & 43 & 44 & 55\\
$27$ & $21304$ & 22 & 21 & 3 & 10 & 39 & 57 & 56 & 38 & 15 & 4 & 47 & 46 & 53 & 45 & 14 & 52 & 51 & 48 & 40 & 9 & 42 & 41 & 43 & 50 & 19 & 17 & 16 & 18 & 55 & 44\\
$28$ & $31240$ & 58 & 36 & 37 & 34 & 5 & 23 & 1 & 2 & 4 & 15 & 8 & 11 & 12 & 19 & 50 & 13 & 6 & 7 & 9 & 40 & 18 & 16 & 17 & 14 & 45 & 43 & 41 & 42 & 24 & 35\\
$29$ & $41203$ & 24 & 1 & 2 & 5 & 33 & 59 & 31 & 32 & 10 & 3 & 49 & 56 & 57 & 55 & 18 & 54 & 26 & 27 & 20 & 8 & 44 & 21 & 22 & 25 & 13 & 19 & 11 & 12 & 30 & 23\\
$30$ & $03241$ & 35 & 58 & 57 & 54 & 26 & 0 & 28 & 27 & 49 & 56 & 10 & 3 & 2 & 4 & 41 & 5 & 33 & 32 & 39 & 51 & 15 & 38 & 37 & 34 & 46 & 40 & 48 & 47 & 29 & 36\\
$31$ & $13204$ & 1 & 23 & 22 & 25 & 54 & 36 & 58 & 57 & 55 & 44 & 51 & 48 & 47 & 40 & 9 & 46 & 53 & 52 & 50 & 19 & 41 & 43 & 42 & 45 & 14 & 16 & 18 & 17 & 35 & 24\\
$32$ & $23140$ & 37 & 38 & 56 & 49 & 20 & 2 & 3 & 21 & 44 & 55 & 12 & 13 & 6 & 14 & 45 & 7 & 8 & 11 & 19 & 50 & 17 & 18 & 16 & 9 & 40 & 42 & 43 & 41 & 4 & 15\\
$33$ & $32104$ & 3 & 2 & 1 & 0 & 29 & 38 & 37 & 36 & 35 & 24 & 53 & 52 & 51 & 50 & 19 & 48 & 47 & 46 & 45 & 14 & 43 & 42 & 41 & 40 & 9 & 18 & 17 & 16 & 15 & 4\\
$34$ & $42130$ & 39 & 32 & 31 & 28 & 0 & 4 & 2 & 1 & 23 & 30 & 14 & 17 & 16 & 18 & 55 & 9 & 7 & 6 & 13 & 25 & 19 & 12 & 11 & 8 & 20 & 44 & 22 & 21 & 3 & 10\\
$35$ & $04213$ & 0 & 29 & 47 & 46 & 53 & 30 & 59 & 52 & 51 & 48 & 55 & 44 & 22 & 21 & 3 & 25 & 54 & 57 & 56 & 38 & 20 & 49 & 27 & 26 & 33 & 10 & 39 & 32 & 31 & 28\\
$36$ & $14230$ & 31 & 59 & 52 & 53 & 45 & 1 & 24 & 42 & 43 & 50 & 16 & 9 & 7 & 8 & 20 & 6 & 14 & 17 & 18 & 55 & 11 & 19 & 12 & 13 & 25 & 21 & 44 & 22 & 23 & 30\\
$37$ & $24103$ & 2 & 4 & 41 & 40 & 48 & 32 & 39 & 51 & 50 & 43 & 57 & 54 & 26 & 25 & 13 & 27 & 49 & 56 & 55 & 18 & 22 & 44 & 21 & 20 & 8 & 12 & 19 & 11 & 10 & 3\\
$38$ & $34120$ & 33 & 34 & 46 & 47 & 40 & 3 & 4 & 41 & 42 & 45 & 18 & 19 & 11 & 12 & 30 & 8 & 9 & 16 & 17 & 35 & 13 & 14 & 6 & 7 & 0 & 23 & 24 & 1 & 2 & 5\\
$39$ & $43102$ & 4 & 3 & 21 & 20 & 27 & 34 & 33 & 26 & 25 & 22 & 59 & 58 & 36 & 35 & 17 & 29 & 28 & 31 & 30 & 12 & 24 & 23 & 1 & 0 & 7 & 14 & 13 & 6 & 5 & 2\\
$40$ & $01423$ & 45 & 46 & 34 & 37 & 38 & 50 & 51 & 39 & 32 & 33 & 20 & 21 & 44 & 22 & 23 & 55 & 56 & 49 & 27 & 28 & 25 & 26 & 54 & 57 & 58 & 30 & 31 & 59 & 52 & 53\\
$41$ & $10432$ & 51 & 50 & 19 & 18 & 17 & 46 & 45 & 14 & 13 & 12 & 1 & 0 & 29 & 28 & 27 & 36 & 35 & 24 & 23 & 22 & 31 & 30 & 59 & 58 & 57 & 26 & 25 & 54 & 53 & 52\\
$42$ & $20413$ & 47 & 40 & 9 & 16 & 18 & 52 & 50 & 19 & 11 & 8 & 22 & 25 & 54 & 26 & 33 & 57 & 55 & 44 & 21 & 3 & 27 & 20 & 49 & 56 & 38 & 32 & 10 & 39 & 51 & 48\\
$43$ & $30421$ & 53 & 45 & 14 & 17 & 16 & 48 & 40 & 9 & 7 & 6 & 3 & 10 & 39 & 32 & 31 & 38 & 15 & 4 & 2 & 1 & 33 & 5 & 34 & 37 & 36 & 28 & 0 & 29 & 47 & 46\\
$44$ & $40312$ & 49 & 20 & 8 & 11 & 12 & 54 & 25 & 13 & 6 & 7 & 24 & 35 & 58 & 36 & 37 & 59 & 30 & 23 & 1 & 2 & 29 & 0 & 28 & 31 & 32 & 34 & 5 & 33 & 26 & 27\\
$45$ & $02431$ & 50 & 52 & 59 & 58 & 36 & 40 & 47 & 29 & 28 & 31 & 5 & 2 & 4 & 3 & 21 & 15 & 37 & 34 & 33 & 26 & 10 & 32 & 39 & 38 & 56 & 20 & 27 & 49 & 48 & 51\\
$46$ & $12403$ & 41 & 42 & 24 & 35 & 58 & 51 & 52 & 59 & 30 & 23 & 26 & 27 & 49 & 20 & 8 & 56 & 57 & 54 & 25 & 13 & 21 & 22 & 44 & 55 & 18 & 11 & 12 & 19 & 50 & 43\\
$47$ & $21430$ & 52 & 51 & 39 & 38 & 15 & 42 & 41 & 4 & 3 & 10 & 7 & 6 & 14 & 13 & 25 & 17 & 16 & 9 & 8 & 20 & 12 & 11 & 19 & 18 & 55 & 22 & 21 & 44 & 43 & 50\\
$48$ & $31402$ & 43 & 41 & 4 & 15 & 37 & 53 & 46 & 34 & 5 & 2 & 28 & 31 & 59 & 30 & 12 & 58 & 36 & 29 & 0 & 7 & 23 & 1 & 24 & 35 & 17 & 13 & 6 & 14 & 45 & 42\\
$49$ & $41320$ & 54 & 26 & 33 & 32 & 10 & 44 & 21 & 3 & 2 & 5 & 9 & 16 & 18 & 17 & 35 & 19 & 11 & 8 & 7 & 0 & 14 & 6 & 13 & 12 & 30 & 24 & 1 & 23 & 22 & 25\\
$50$ & $03412$ & 40 & 48 & 49 & 56 & 57 & 45 & 53 & 54 & 26 & 27 & 30 & 23 & 24 & 1 & 2 & 35 & 58 & 59 & 31 & 32 & 0 & 28 & 29 & 36 & 37 & 5 & 33 & 34 & 46 & 47\\
$51$ & $13420$ & 46 & 53 & 54 & 57 & 55 & 41 & 43 & 44 & 22 & 25 & 11 & 8 & 9 & 7 & 0 & 16 & 18 & 19 & 12 & 30 & 6 & 13 & 14 & 17 & 35 & 1 & 23 & 24 & 42 & 45\\
$52$ & $23401$ & 42 & 43 & 44 & 55 & 56 & 47 & 48 & 49 & 20 & 21 & 32 & 33 & 34 & 5 & 6 & 37 & 38 & 39 & 10 & 11 & 2 & 3 & 4 & 15 & 16 & 7 & 8 & 9 & 40 & 41\\
$53$ & $32410$ & 48 & 47 & 29 & 36 & 35 & 43 & 42 & 24 & 1 & 0 & 13 & 12 & 19 & 11 & 10 & 18 & 17 & 14 & 6 & 5 & 8 & 7 & 9 & 16 & 15 & 3 & 2 & 4 & 41 & 40\\
$54$ & $42301$ & 44 & 22 & 23 & 30 & 31 & 49 & 27 & 28 & 0 & 1 & 34 & 37 & 38 & 15 & 16 & 39 & 32 & 33 & 5 & 6 & 4 & 2 & 3 & 10 & 11 & 9 & 7 & 8 & 20 & 21\\
$55$ & $04321$ & 25 & 54 & 53 & 52 & 51 & 20 & 49 & 48 & 47 & 46 & 15 & 4 & 3 & 2 & 1 & 10 & 39 & 38 & 37 & 36 & 5 & 34 & 33 & 32 & 31 & 0 & 29 & 28 & 27 & 26\\
$56$ & $14302$ & 21 & 44 & 43 & 50 & 52 & 26 & 54 & 53 & 45 & 42 & 36 & 29 & 28 & 0 & 7 & 31 & 59 & 58 & 35 & 17 & 1 & 24 & 23 & 30 & 12 & 6 & 14 & 13 & 25 & 22\\
$57$ & $24310$ & 27 & 49 & 48 & 51 & 50 & 22 & 44 & 43 & 41 & 40 & 17 & 14 & 13 & 6 & 5 & 12 & 19 & 18 & 16 & 15 & 7 & 9 & 8 & 11 & 10 & 2 & 4 & 3 & 21 & 20\\
$58$ & $34201$ & 23 & 24 & 42 & 45 & 46 & 28 & 29 & 47 & 40 & 41 & 38 & 39 & 32 & 10 & 11 & 33 & 34 & 37 & 15 & 16 & 3 & 4 & 2 & 5 & 6 & 8 & 9 & 7 & 0 & 1\\
$59$ & $43210$ & 29 & 28 & 27 & 26 & 25 & 24 & 23 & 22 & 21 & 20 & 19 & 18 & 17 & 16 & 15 & 14 & 13 & 12 & 11 & 10 & 9 & 8 & 7 & 6 & 5 & 4 & 3 & 2 & 1 & 0\\

\end{tabular}
\caption{The Lehmer law of $\Alt_5$, columns $30$--$59$, continuing Table~\ref{tab:an:law5a}.}\label{tab:an:law5b}
\end{table}

\subsection{The rotation laws}\label{sec:altp}
Throughout, $n\ge3$, permutations act on $\Ap n$, and $(\sigma\circ\tau)(i)=\sigma(\tau(i))$, the right-to-left composition of~\eqref{eq:comp}.

\begin{definition}\label{def:prefixrot}
For $2\le k\le n$ the \emph{prefix rotation} $\rho_k$ is the $k$-cycle $(0\,1\,\cdots\,k-1)$ of $\Ap n$. Its powers are
\begin{equation}\label{eq:prefixrot}
\rho_k^{\ell}(i)=i+\Bf\big(\fl{i/k}\big)\big(\dg{k}{0}{i+\ell}-i\big),\qquad i\in\Ap n .
\end{equation}
\end{definition}
\noindent Since $\Bf(\fl{i/k})=1$ exactly when $i<k$, the letters $i\ge k$ are fixed and the first $k$ letters are shifted
cyclically by $\ell$.

\begin{lemma}[Rotation normal form]\label{lem:rotnf}
Every $\sigma\in\Sym_n$ has a unique expression
\[
\sigma=\rho_n^{\ell_{n-1}}\circ\cdots\circ\rho_3^{\ell_2}\circ\rho_2^{\ell_1},\qquad \ell_h\in\Ap{h+1}.
\]
The exponents are read off from the top down: with $\pi_n=\sigma$, for $k=n,n-1,\dots,2$,
\begin{equation}\label{eq:rotdecode}
\ell_{k-1}=\dg{k}{0}{\pi_k(k-1)+1},\qquad \pi_{k-1}(i)=\dg{k}{0}{\pi_k(i)-\ell_{k-1}}\quad(i<k-1).
\end{equation}
Moreover $\sgn\sigma=(-1)^{\ell_1+\ell_3+\ell_5+\cdots}$.
\end{lemma}
\begin{proof}
Since $\rho_k^{\ell}(k-1)=\dg{k}{0}{\ell-1}$, the $k$ powers of $\rho_k$ send $k-1$ to distinct letters of $\Ap k$. Let $\pi_k$
be a permutation of $\Ap k$ (fixing the letters $\ge k$). Exactly one $\ell\in\Ap k$ satisfies
$\rho_k^{-\ell}\big(\pi_k(k-1)\big)=k-1$, namely $\ell=\dg{k}{0}{\pi_k(k-1)+1}$, and then $\pi_{k-1}=\rho_k^{-\ell}\circ\pi_k$ is a
permutation of $\Ap{k-1}$. Its values are those in~\eqref{eq:rotdecode}, because $\rho_k^{-\ell}(j)=\dg{k}{0}{j-\ell}$ for $j<k$.
Induction from $k=n$ down to $k=2$, where $\pi_1$ is the identity, gives existence and uniqueness; as a check, there are
$2\cdot3\cdots n=n!$ exponent tuples. A $k$-cycle has sign $(-1)^{k-1}$, so $\sgn\sigma=(-1)^{\sum_h h\ell_h}$, and only the odd
$h$ contribute.
\end{proof}

\begin{theorem}[Rotation law of $\Sym_n$]\label{thm:sn:rotlaw}
For $\sigma\in\Sym_n$ let $\ell_h(\sigma)$ be the exponents of Lemma~\ref{lem:rotnf}, and put
\[
r(\sigma)=\sum_{h=1}^{n-1}\ell_h(\sigma)\,h!\in\Ap{n!}.
\]
Then $r$ is a bijection, and
\[
x*_\rho y=r\bigl(r^{-1}(x)\circ r^{-1}(y)\bigr)
\]
makes $\Ap{n!}$ a group isomorphic to $\Sym_n$, with identity $0$. The digit of place value $h!$ is $\ell_h$. In particular
$r(\rho_2)=1$ and
\[
x*_\rho1=\begin{cases}x+1,&x\text{ even},\\x-1,&x\text{ odd},\end{cases}
\]
so right multiplication by the transposition $\rho_2$ toggles the radix-$2$ digit.
\end{theorem}
\begin{proof}
Lemma~\ref{lem:rotnf} identifies $\Sym_n$ bijectively with $\Ap2\times\Ap3\times\cdots\times\Ap n$, and mixed-radix evaluation
with place values $1!,2!,\dots,(n-1)!$ identifies this product with $\Ap{n!}$ (Lemma~\ref{lem:mixed}). The law is transport of
structure. Since $\rho_2$ is the rightmost factor of the normal form, right multiplication by $\rho_2$ replaces $\ell_1$ by
$\ell_1+1\pmod 2$ and leaves the other exponents unchanged.
\end{proof}

\begin{theorem}[Rotation law of $\Alt_n$]\label{thm:altp}
Let $n\ge3$. Write $a\in\Ap{n!/2}$ in the mixed radix $(3,4,\dots,n)$ and put
\begin{equation}\label{eq:altlabel}
\ell_h=\md{h-2}{a}\quad(h=2,\dots,n-1),\qquad \ell_1=\dg{2}{0}{\ell_3+\ell_5+\cdots}.
\end{equation}
Let $\sigma_a=\rho_n^{\ell_{n-1}}\circ\cdots\circ\rho_2^{\ell_1}$, evaluated with~\eqref{eq:prefixrot}. Then $a\mapsto\sigma_a$ is a
bijection $\Ap{n!/2}\to\Alt_n$, and
\[
a=\sum_{h=2}^{n-1}\ell_h\,\frac{h!}{2}=\Big\lfloor\frac{r(\sigma_a)}{2}\Big\rfloor.
\]
For $a,b\in\Ap{n!/2}$ let $u_h=\ell_h(\sigma_a\circ\sigma_b)$ be the exponents given by~\eqref{eq:rotdecode}. Then
\begin{equation}\label{eq:altlaw}
a\star_\rho b=\sum_{h=2}^{n-1}u_h\,\frac{h!}{2}
\end{equation}
makes $\Ap{n!/2}$ a group isomorphic to $\Alt_n$, with identity $0$.
\end{theorem}
\begin{proof}
By Lemma~\ref{lem:rotnf}, $\sgn\sigma=(-1)^{\ell_1+\ell_3+\ell_5+\cdots}$. Apply Lemma~\ref{lem:sn:paritydigit} with
$U=\Ap3\times\cdots\times\Ap n$, parity digit $e=\ell_1$ and $q=\ell_3+\ell_5+\cdots\pmod2$. Mixed-radix evaluation of the remaining
digits has place values $h!/2$, which gives the bijection and the first expression for $a$. Since $r(\sigma_a)=\ell_1+2a$ with
$\ell_1\in\Ap2$, also $a=\lfloor r(\sigma_a)/2\rfloor$. The product $\sigma_a\circ\sigma_b$ is even, so $u_1$ is its restored parity
digit, and~\eqref{eq:altlaw} is its label. For $a=0$ all exponents vanish and $\sigma_0$ is the identity.
\end{proof}

\begin{example}
For $n=5$ the labels $a\in\Ap{60}$ have digits $\ell_2\in\Ap3$, $\ell_3\in\Ap4$, $\ell_4\in\Ap5$, with $a=\ell_2+3\ell_3+12\ell_4$
and $\ell_1=\dg{2}{0}{\ell_3}$. The first labels, with $\sigma_a$ written as $\sigma_a(0)\cdots\sigma_a(4)$, are
\begin{center}
\begin{tabular}{c|cccccccc}
$a$ & $0$ & $1$ & $2$ & $3$ & $4$ & $5$ & $6$ & $7$\\ \hline
$\sigma_a$ & $01234$ & $12034$ & $20134$ & $21304$ & $32104$ & $13204$ & $23014$ & $30214$
\end{tabular}
\end{center}
The labels $0,1,2$ form the subgroup $\Alt_3$, and the labels $0,\dots,11$ the subgroup $\Alt_4$ fixing the letter $4$. A
computer check of~\eqref{eq:altlaw} confirms associativity for $n\le5$ and gives the conjugacy class sizes $1,12,12,15,20$ for
$n=5$ and $1,40,40,45,72,72,90$ for $n=6$, those of the simple groups $\Alt_5$ and $\Alt_6$.
\end{example}

A third chain labelling comes directly from Theorem~\ref{thm:transitive}, with cyclic shifts as coset representatives.

\begin{corollary}[Symmetric and alternating groups]\label{cor:alternating}
Let $r(x)=\dg{n}{0}{x+1}$ on $\Ap n$. For $\Sym_n$, and for $\Alt_n$ with $n$ odd, take $s_t=r^t$. Then~\eqref{eq:transitive} becomes
\begin{equation}\label{eq:shiftdiff}
(t,h)*(u,k)=\big(\dg{n}{0}{t+h(u)},\;h^{(u)}\circ k\big),\qquad h^{(u)}(x)=\dg{n}{0}{h(x+u)-h(u)} .
\end{equation}
The stabiliser of $0$ is $\Sym_{n-1}$ (respectively $\Alt_{n-1}$) acting on $\{1,\dots,n-1\}$. Relabelling $x\mapsto x-1$ and repeating gives a
mixed-radix code
\[
m=t_n+n\big(t_{n-1}+(n-1)(t_{n-2}+\cdots)\big),\qquad t_j\in\Ap j ,
\]
with radices $n,n-1,\dots,2$ for $\Sym_n$ and $n,n-1,\dots,3$ for $\Alt_n$. For $\Alt_n$ with $n$ even, take $s_t=r^t$ for $t$ even and
$s_t=r^t\circ(1\,2)$ for $t$ odd.
\end{corollary}
\begin{proof}
For $s_t=r^t$ we get $t'=t+h(u)$ and $s_{t'}^{-1}s_ths_u=r^{-h(u)}hr^u$, which sends $x$ to $h(x+u)-h(u)$. For odd $n$ the $n$-cycle $r$
is even, so $s_t\in \Alt_n$. For even $n$ the given $s_t$ are even and satisfy $s_t(0)=t$, because $(1\,2)$ fixes $0$. Now apply
Theorem~\ref{thm:transitive}.
\end{proof}

\begin{remark}
Corollary~\ref{cor:alternating} also uses rotations as coset representatives, but it peels off the stabiliser of a point and
needs an extra transposition when $n$ is even. The prefix rotations together with the parity digit $\ell_1$ give one
formula for all $n$.
\end{remark}

\begin{remark}[Why the decoding is a recursion]
For $n\ge5$ the group $\Alt_n$ is simple, so no layer projector of a labelling of $\Alt_n$ onto a layer of size
strictly between $1$ and $n!/2$ is a homomorphism: every layer is coupled to the layers above it. Correspondingly the decoding
\eqref{eq:rotdecode} runs from the top layer down, each exponent depending on those above it, much as carries propagate in
ordinary addition.
\end{remark}

\subsection{Comparison and checks}
\begin{remark}[Two chain labellings]\label{rem:sn:twolaws}
The Lehmer and rotation labels are two stabiliser-chain coordinate systems on the same groups.
\begin{center}
\begin{tabular}{lll}
 & Lehmer ($*_L$, $\star_L$) & rotation ($*_\rho$, $\star_\rho$)\\ \hline
radices & $n,n-1,\dots,2$ & $2,3,\dots,n$\\
place values & $n!/(n-i)!$ & $h!$\\
radix-$2$ digit & last, $d_{n-2}$ & first, $\ell_1$\\
sign & parity of all digits & parity of $\ell_1+\ell_3+\cdots$\\
parity toggle & right mult.\ by $s_{n-2}$: $m\mapsto m\pm n!/2$ & right mult.\ by $\rho_2$: $x\mapsto x\pm1$\\
$\Alt_n$ label & $\lambda\bmod n!/2$ & $\lfloor r/2\rfloor$
\end{tabular}
\end{center}
In both cases Lemma~\ref{lem:sn:paritydigit} produces the $\Alt_n$ law by deleting the radix-$2$ digit. The relabellings
\[
\Phi=r\circ\lambda^{-1}\colon(\Ap{n!},*_L)\to(\Ap{n!},*_\rho),\qquad
\Phi_A=\lfloor r\circ\mu^{-1}/2\rfloor\colon(\Ap{n!/2},\star_L)\to(\Ap{n!/2},\star_\rho)
\]
are explicit group isomorphisms. The Lehmer coordinates make the Coxeter generators, the sign and the cosets of $\Alt_n$
transparent. The rotation coordinates make the nested subgroups $\Sym_2<\Sym_3<\cdots<\Sym_n$ and $\Alt_3<\cdots<\Alt_n$ transparent:
they are initial segments of labels.
\end{remark}

\begin{remark}[Checks]\label{rem:sn:lawchecks}
The script \texttt{scripts/sn\_law\_verify.py} checks both labellings exhaustively for $n\le6$. For the Lehmer laws it verifies
encoding and decoding, identities and inverses, both generator rules, the sign and coset rule, the $\Alt_n$ product formula, the
conjugation rule and the chain property. For the rotation laws it verifies Lemma~\ref{lem:rotnf}, Theorems~\ref{thm:sn:rotlaw}
and~\ref{thm:altp}, the sign formula, the identity $a=\lfloor r(\sigma_a)/2\rfloor$, the parity toggle and the relabelling
between the two $\Alt_n$ laws. For $n\le5$ it forms the full Cayley tables and checks associativity. When GAP is available the
right regular representations of all four laws are checked with \texttt{IsomorphismGroups} against
\texttt{SymmetricGroup($n$)} and \texttt{AlternatingGroup($n$)} for $n\le6$. The script writes Tables~\ref{tab:sn:law3}
and~\ref{tab:sn:laws} and is called by \texttt{scripts/sn\_verify.py}. GAP~4.14.0 \cite{GAP} confirms the shift labelling of Corollary~\ref{cor:alternating} on small cases: the codes of
$\Alt_5$ and $\Alt_6$ have \texttt{IdGroup} $[60,5]$ and $[360,118]$, and the code of $\Alt_7$ gives a simple group isomorphic to $\Alt_7$.
\end{remark}

\section{Worked products, inverses and conjugates}\label{sec:sn:worked}
This section computes by hand a product in the Lehmer laws of $\Sym_5$ and of $\Alt_5$, an inverse and a conjugate, in the manner
of Chapter~\ref{ch:products}. Labels are decoded with Proposition~\ref{prop:seqdigits} and Theorem~\ref{thm:sn:lehmer}, and
permutations are composed from right to left, as in~\eqref{eq:comp}.

\begin{example}[A product in $\Sym_5$]\label{ex:sn:wprod}
We compute $68*101$ in the Lehmer law of $\Sym_5$ (Theorem~\ref{thm:sn:lehmerlaw}), whose place values are
$\Pi_i=1,5,20,60,120$. By Proposition~\ref{prop:seqdigits},
\[
68=13\cdot5+3,\quad 13=3\cdot4+1,\quad 3=1\cdot3+0,\quad 1=0\cdot2+1,\quad 0=0\cdot1+0,
\]
so $68$ has digits $(a_0,\dots,a_4)=(3,1,0,1,0)$. Reading each digit as a rank among the values not yet used
(Theorem~\ref{thm:sn:lehmer}) gives $R=\lambda^{-1}(68)=(3,1,0,4,2)$ in one-line notation. In the same way $101$ has digits $(1,0,2,1,0)$
and $S=\lambda^{-1}(101)=(1,0,4,3,2)$. Composing from right to left, $(R\circ S)(i)=R(S(i))$, gives $R\circ S=(1,3,2,4,0)$, whose inversion
counts are $(1,2,1,1,0)$. Hence
\[
68*101=\lambda(R\circ S)=1\cdot1+2\cdot5+1\cdot20+1\cdot60=91.
\]
\end{example}

\begin{example}[A product in $\Alt_5$]\label{ex:an:wprod}
We compute $23*41$ in the Lehmer law of $\Alt_5$ (Theorem~\ref{thm:an:lehmerlaw}), with $N=60$. A label $m\in\Ap{60}$ is
decoded through its Lehmer digits in $\Sym_5$: if $\varepsilon(m)$ is the parity of their sum, then $m+60\,\varepsilon(m)$ is the
Lehmer label of an even permutation, and a product is encoded by reducing its Lehmer label modulo $60$. Here $23$ has digits
$(3,0,1,0,0)$ with sum $4$, so $\varepsilon(23)=0$ and $R=\lambda^{-1}(23)=(3,0,2,1,4)$. Similarly $41$ has digits $(1,0,2,0,0)$ with sum
$3$, so $\varepsilon(41)=1$ and $S=\lambda^{-1}(101)=(1,0,4,3,2)$. Then $R\circ S=(0,3,4,1,2)$ has digits $(0,2,2,0,0)$, so
\[
\lambda(R\circ S)=2\cdot5+2\cdot20=50,\qquad 23*41=50\bmod 60=50.
\]
\end{example}

\begin{example}[An inverse in $\Sym_5$]\label{ex:sn:winv}
Let $R=\lambda^{-1}(68)=(3,1,0,4,2)$ be as in Example~\ref{ex:sn:wprod}. Since $R^{-1}(R(i))=i$, the word of $R^{-1}$ lists the
positions of $0,1,\dots,4$ in the word of $R$, so $R^{-1}=(2,1,4,0,3)$. Its inversion counts are $(2,1,2,0,0)$, hence
$68^{-1}=\lambda(R^{-1})=2\cdot1+1\cdot5+2\cdot20=47$. The script \texttt{inverse.py} reaches the same label inside the law, without decoding
(Remark~\ref{rem:wpinverses}): the powers of $68$ return to $68$ after $4$ steps, so $68$ has order $4$ and
$68^{-1}=68^{3}=47$.
\end{example}

\begin{example}[A conjugate in $\Sym_5$]\label{ex:sn:wconj}
With $R$ and $S$ as in Example~\ref{ex:sn:wprod}, the conjugate $R\circ S\circ R^{-1}$ maps $R(i)$ to $R(S(i))$. It is therefore
obtained from the cycles of $S$ by applying $R$ to every entry: $S=(0\,1)(2\,4)$ gives $R\circ S\circ R^{-1}=(3\,1)(0\,2)$, with word $(2,3,0,1,4)$
and inversion counts $(2,2,0,0,0)$. Hence
\[
68*101*68^{-1}=2\cdot1+2\cdot5=12,
\]
which is the output of \texttt{conjugate.py} \texttt{Sym(5)} \texttt{68} \texttt{101}. The other convention gives
$68^{-1}*101*68=65$ (option \texttt{--right}).
\end{example}

For larger degrees the same computations are carried out by the scripts. The following outputs were obtained with
\texttt{--check} and agree with direct composition of the decoded permutations.
\begingroup\small
\begin{center}
\begin{tabular}{@{}ll@{}}
command&output\\ \hline
\texttt{multiply.py} \allowbreak{}\texttt{Sym(20)} \allowbreak{}\texttt{123456789012345} \allowbreak{}\texttt{98765432109876}&\texttt{34297188027031591}\\
\texttt{inverse.py} \allowbreak{}\texttt{Sym(20)} \allowbreak{}\texttt{123456789012345}&\texttt{508822672168494726}\\
\texttt{conjugate.py} \allowbreak{}\texttt{Sym(20)} \allowbreak{}\texttt{123456789012345} \allowbreak{}\texttt{98765432109876}&\texttt{558695374284871541}\\
\texttt{multiply.py} \allowbreak{}\texttt{Alt(12)} \allowbreak{}\texttt{123456789} \allowbreak{}\texttt{98765432}&\texttt{125411668}\\
\texttt{inverse.py} \allowbreak{}\texttt{Alt(12)} \allowbreak{}\texttt{123456789}&\texttt{67842571}\\
\texttt{conjugate.py} \allowbreak{}\texttt{Alt(12)} \allowbreak{}\texttt{123456789} \allowbreak{}\texttt{98765432}&\texttt{43983445}\\
\end{tabular}
\end{center}
\endgroup

\begin{remark}[Checks]\label{rem:sn:wchecks}
The products, inverses and conjugates of this section are reproduced by \texttt{multiply.py}, \texttt{inverse.py} and
\texttt{conjugate.py}, run from the directory \texttt{scripts} as in Remarks~\ref{rem:wpchecks} and~\ref{rem:wpinverses}. The
group names are \texttt{Sym(n)} and \texttt{Alt(n)}. In \texttt{zsh}, quote a name containing parentheses, for example
\texttt{'Sym(20)'}, or use the aliases \texttt{Sym20}, \texttt{S20}, \texttt{Alt12} and \texttt{A12}. With \texttt{--check},
\texttt{multiply.py} compares each product of
$\Sym_n$ or $\Alt_n$ with a separately coded ranking of permutations. Labels of any size are accepted; decimal labels with more
than $4300$ digits, which occur for $n$ larger than about $1560$, require raising the default limit of Python on the conversion of
long decimal strings.
\end{remark}

\begin{remark}[Conjugates]\label{rem:sn:wconj}
The conjugate of a label $n$ by a label $m$, in any group of the archive, is computed by
\begin{center}\ttfamily\small python3 conjugate.py GROUP m n [--right] [--check]\end{center}
which prints $m*n*m^{-1}$, the image of $n$ under the inner automorphism $x\mapsto mxm^{-1}$; with \texttt{--right} it prints
$m^{-1}*n*m$. The script uses only the law of \texttt{multiply.py} and the inverse of \texttt{inverse.py}, so it applies without
change to every group for which these scripts are available. With \texttt{--check} it verifies $m*m^{-1}=m^{-1}*m=e$, the
identity $(m*n*m^{-1})*m=m*n$, and that $n$ and its conjugate have the same order.
\end{remark}

\section{Generated tables}\label{sec:sn:tables}
The numbers in this section depend only on the partitions of $n$. A conjugacy class of $\Sym_n$ is given by a cycle type
$\lambda=(1^{m_1}2^{m_2}\cdots)$ with $\ell(\lambda)$ parts (Section~\ref{sec:sn:classes}). For $g$ in the class $C_\lambda$,
\[
|C_\lambda|=\frac{n!}{z_\lambda},\qquad z_\lambda=\prod_i i^{m_i}m_i!,\qquad o(g)=\operatorname{lcm}(\lambda_1,\lambda_2,\dots),\qquad
\sgn g=(-1)^{n-\ell(\lambda)},
\]
where $z_\lambda$ is the order of the centraliser of $g$. The class lies in $\Alt_n$ when $n-\ell(\lambda)$ is even, and then it
splits into two classes of $\Alt_n$ exactly when the parts of $\lambda$ are distinct and odd. The irreducible character
$\chi^\lambda$ of $\Sym_n$ has degree $f^\lambda=n!/\prod_x h(x)$, the product running over the hook lengths of the diagram of
$\lambda$ (Section~\ref{sec:sn:characters}). On $\Alt_n$, $\chi^\lambda$ and $\chi^{\lambda'}$ restrict to the same irreducible
character when $\lambda\ne\lambda'$, and $\chi^\lambda$ splits into two characters of degree $f^\lambda/2$ when $\lambda=\lambda'$
(Section~\ref{sec:an:characters}). We write $k(G)$ for the number of conjugacy classes of $G$.

In Table~\ref{tab:sn:orders}, $e(G)$ is the exponent of $G$, $g(n)$ is Landau's function, the largest order of an element of
$\Sym_n$, $g_A(n)$ is the largest order of an element of $\Alt_n$, and $t(G)$ is the number of $x\in G$ with $x^2=1$, the
identity included. Every irreducible representation of $\Sym_n$ is realised over $\mathbb Q$, so every Frobenius--Schur
indicator is $1$ and $t(\Sym_n)$ is the sum of the degrees in Table~\ref{tab:sn:degrees}. It also satisfies
$t(\Sym_n)=t(\Sym_{n-1})+(n-1)\,t(\Sym_{n-2})$, since a solution either fixes the letter $n-1$ or exchanges it with one of the
other $n-1$ letters.

\begin{table}[ht]
\centering
\begin{tabular}{rrrrrrr}
$n$ & $e(\Sym_n)$ & $g(n)$ & $e(\Alt_n)$ & $g_A(n)$ & $t(\Sym_n)$ & $t(\Alt_n)$\\ \hline
2 & 2 & 2 & 1 & 1 & 2 & 1 \\
3 & 6 & 3 & 3 & 3 & 4 & 1 \\
4 & 12 & 4 & 6 & 3 & 10 & 4 \\
5 & 60 & 6 & 30 & 5 & 26 & 16 \\
6 & 60 & 6 & 60 & 5 & 76 & 46 \\
7 & 420 & 12 & 420 & 7 & 232 & 106 \\
8 & 840 & 15 & 420 & 15 & 764 & 316 \\

\end{tabular}
\caption{Exponents, largest element orders and solutions of $x^2=1$.}\label{tab:sn:orders}
\end{table}

\begin{table}[ht]
\centering\small
\begin{tabular}{lrl}
$G$ & $k(G)$ & degrees of the irreducible characters\\ \hline
$\Sym_{2}$ & 2 & $1^{2}$ \\
$\Sym_{3}$ & 3 & $1^{2}, 2$ \\
$\Sym_{4}$ & 5 & $1^{2}, 2, 3^{2}$ \\
$\Sym_{5}$ & 7 & $1^{2}, 4^{2}, 5^{2}, 6$ \\
$\Sym_{6}$ & 11 & $1^{2}, 5^{4}, 9^{2}, 10^{2}, 16$ \\
$\Sym_{7}$ & 15 & $1^{2}, 6^{2}, 14^{4}, 15^{2}, 20, 21^{2}, 35^{2}$ \\
$\Sym_{8}$ & 22 & $1^{2}, 7^{2}, 14^{2}, 20^{2}, 21^{2}, 28^{2}, 35^{2}, 42, 56^{2}, 64^{2}, 70^{2}, 90$ \\
\hline
$\Alt_{3}$ & 3 & $1^{3}$ \\
$\Alt_{4}$ & 4 & $1^{3}, 3$ \\
$\Alt_{5}$ & 5 & $1, 3^{2}, 4, 5$ \\
$\Alt_{6}$ & 7 & $1, 5^{2}, 8^{2}, 9, 10$ \\
$\Alt_{7}$ & 9 & $1, 6, 10^{2}, 14^{2}, 15, 21, 35$ \\
$\Alt_{8}$ & 14 & $1, 7, 14, 20, 21^{3}, 28, 35, 45^{2}, 56, 64, 70$ \\

\end{tabular}
\caption{Character degrees of $\Sym_n$ and $\Alt_n$; exponents are multiplicities.}\label{tab:sn:degrees}
\end{table}

\begin{table}[ht]
\centering
\begin{tabular}{lrrrcc}
$\lambda$ & $|C_\lambda|$ & $z_\lambda$ & $o(g)$ & $\sgn g$ & classes in $\Alt_8$\\ \hline
$8$ & 5040 & 8 & 8 & $-$ & -- \\
$7,1$ & 5760 & 7 & 7 & $+$ & 2 \\
$6,2$ & 3360 & 12 & 6 & $+$ & 1 \\
$6,1^{2}$ & 3360 & 12 & 6 & $-$ & -- \\
$5,3$ & 2688 & 15 & 15 & $+$ & 2 \\
$5,2,1$ & 4032 & 10 & 10 & $-$ & -- \\
$5,1^{3}$ & 1344 & 30 & 5 & $+$ & 1 \\
$4^{2}$ & 1260 & 32 & 4 & $+$ & 1 \\
$4,3,1$ & 3360 & 12 & 12 & $-$ & -- \\
$4,2^{2}$ & 1260 & 32 & 4 & $-$ & -- \\
$4,2,1^{2}$ & 2520 & 16 & 4 & $+$ & 1 \\
$4,1^{4}$ & 420 & 96 & 4 & $-$ & -- \\
$3^{2},2$ & 1120 & 36 & 6 & $-$ & -- \\
$3^{2},1^{2}$ & 1120 & 36 & 3 & $+$ & 1 \\
$3,2^{2},1$ & 1680 & 24 & 6 & $+$ & 1 \\
$3,2,1^{3}$ & 1120 & 36 & 6 & $-$ & -- \\
$3,1^{5}$ & 112 & 360 & 3 & $+$ & 1 \\
$2^{4}$ & 105 & 384 & 2 & $+$ & 1 \\
$2^{3},1^{2}$ & 420 & 96 & 2 & $-$ & -- \\
$2^{2},1^{4}$ & 210 & 192 & 2 & $+$ & 1 \\
$2,1^{6}$ & 28 & 1440 & 2 & $-$ & -- \\
$1^{8}$ & 1 & 40320 & 1 & $+$ & 1 \\

\end{tabular}
\caption{The $22$ conjugacy classes of $\Sym_8$.}\label{tab:sn:classes8}
\end{table}

\begin{remark}[Checks]\label{rem:sn:tablechecks}
The script \texttt{scripts/sn\_tables.py} writes Tables~\ref{tab:sn:orders},~\ref{tab:sn:degrees} and~\ref{tab:sn:classes8}, and a
data file with the same data for $n\le10$. It writes them only after the following checks pass for every $n\le10$. The class
sizes add up to $n!$, and the number of partitions agrees with Euler's pentagonal recurrence. The hook formula agrees with the
branching rule $f^\lambda=\sum f^{\lambda-\square}$, and the squared degrees add up to $n!$ for $\Sym_n$ and to $n!/2$ for
$\Alt_n$. The numbers of classes and of irreducible characters of $\Alt_n$ agree, and there are as many self-conjugate partitions
as partitions into distinct odd parts. The three values of $t(\Sym_n)$, from the classes, from the recurrence and as the sum of the
degrees, agree, and $g(n)$ equals the largest product of powers of distinct primes with sum at most $n$. For $n\le8$ the script
also enumerates $\Sym_n$: it recomputes cycle types, orders, signs, exponents and $t$, and checks the splitting of each class of
$\Alt_n$ by computing its orbit under conjugation by the $3$-cycles $(0\,1\,i)$. It compares the class numbers with
Table~\ref{tab:sn:normal}. When GAP is available it compares the degrees with \texttt{CharacterTable("Symmetric",$n$)} and
\texttt{CharacterTable("Alternating",$n$)} for $n\le10$.
\end{remark}

\section{Exercises}
\begin{enumerate}
\item Show that $k(\Sym_6)=11$ and list the classes of $\Sym_6$ with their sizes.
\item Show that $\Sym_7$ has $8$ classes of even permutations, that only the class of $7$-cycles splits in $\Alt_7$, and deduce $k(\Alt_7)=9$.
\item Use the hook length formula to show that the irreducible characters of $\Sym_5$ have degrees $1,1,4,4,5,5,6$, and check that
the squares add up to $120$.
\item Show that $\Alt_n$ is generated by its $3$-cycles for $n\ge3$.
\item Show that the classes of cycle types $2$ and $2^3$, $3$ and $3^2$, and $6$ and $3\,2\,1$ of $\Sym_6$ have equal sizes in pairs,
as required by the outer automorphism.
\item Show that the number of permutations of $\Ap n$ with exactly $k$ inversions is the coefficient of $t^k$ in
$\prod_{i=1}^{n}(1+t+\dots+t^{i-1})$.
\end{enumerate}

\chapter{Finite simple groups of Lie type}\label{ch:lie}
This chapter gives digital laws for the finite simple groups of Lie type: the classical groups, the Suzuki and Ree groups,
and the exceptional groups up to $E_8(q)$. Simplicity has a direct meaning for digital laws. A layer is drawn as a solid node exactly
when its layer projector is a homomorphism, and then its kernel is a normal subgroup. A simple group has no proper nontrivial normal
subgroups, so no labelling of it can be graded, and it cannot be built by stacking hollow layers on solid ones. Its laws must come
from another principle. For the alternating groups of Chapter~\ref{ch:simple} it was a normal form of permutations; here it is a
normal form of matrices, mostly from the Bruhat decomposition, which is uniform in $q$. After the field arithmetic, the smallest
family $PSL(2,p)$ is treated in Section~\ref{sec:psl2p}.

\section{Dynkin diagrams, Weyl data and digital laws}\label{sec:dynkindigital}
A Dynkin diagram is the compact combinatorial record of the root data that organise this chapter. Its vertices are the simple roots, and its edges encode the Cartan integers. Thus it determines the Cartan matrix, the Coxeter presentation of the Weyl group, the positive roots and, in particular, the Bruhat cells used below. It does not replace the matrix calculations: it predicts their finite combinatorics and says which calculations belong to the same family.

For a split group of rank $r$, let $\Phi^+$ be its positive roots and let $W$ be its Weyl group. The Borel subgroup has the form $B=UT$, with one field coordinate for each positive-root subgroup in $U$ and $r$ multiplicative coordinates in $T$ (subject to the central quotient when one passes to a simple group). The Bruhat decomposition
\[
 G=\coprod_{w\in W}BwB
\]
therefore gives the architecture of a digital label: first a Weyl digit selecting $w$, then coordinates for the two unipotent pieces and the torus. The explicit normal forms of Sections~\ref{sec:psln}--\ref{sec:e8} are implementations of this architecture, with the elimination procedure determining the digits rather than merely counting them.

\begin{table}[ht]
\centering\small
\begin{tabular}{@{}l|p{.19\linewidth}|p{.19\linewidth}|c|p{.30\linewidth}@{}}
type & diagram & $|W|$ & $|\Phi^+|$ & treatment in this chapter\\ \hline
$A_r$ & $\circ-\cdots-\circ$ & $(r+1)!$ & $r(r+1)/2$ & $\PSL_{r+1}(q)$, Section~\ref{sec:psln}\\
$B_r,C_r$ & $\circ-\cdots-\circ\Rightarrow\circ$ & $2^r r!$ & $r^2$ & odd orthogonal and symplectic groups, Section~\ref{sec:classical}\\
$D_r$ & fork at the last vertex & $2^{r-1}r!$ & $r(r-1)$ & even orthogonal groups, Sections~\ref{sec:classical}, \ref{sec:twistedorth}\\
$G_2$ & $\circ\Leftarrow\!\!\!\!\circ$ & $12$ & $6$ & $G_2(q)$ and ${}^2G_2(q)$, Sections~\ref{sec:g2}, \ref{sec:ree}\\
$F_4$ & $\circ-\circ\Leftarrow\circ-\circ$ & $1152$ & $24$ & $F_4(q)$ and ${}^2F_4(q)$, Sections~\ref{sec:f4}, \ref{sec:2f4}\\
$E_6,E_7,E_8$ & exceptional simply laced trees & $51840$, $2903040$, $696729600$ & $36$, $63$, $120$ & Sections~\ref{sec:e6}--\ref{sec:e8}
\end{tabular}
\caption{Dynkin data behind the digital laws of Chapter~\ref{ch:lie}. A double or triple edge is directed toward the short root.}\label{tab:dynkindigital}
\end{table}

The classical diagrams explain the matrix restrictions. Type $A_r$ has Weyl group $\Sym_{r+1}$, so the pivot permutations in Gaussian elimination are ordinary permutations. Types $B_r$ and $C_r$ have signed permutations, and type $D_r$ has the even signed permutations; these are exactly the permutations preserving the antidiagonal forms in Section~\ref{sec:classical}. The explicit Coxeter generators of Chapter~\ref{ch:simple} are therefore the type-$A$ prototype of the Weyl digits used here.

Twisting has an equally simple diagrammatic meaning. A symmetry of the diagram may be composed with the Frobenius map $t\mapsto t^q$: this gives ${}^2A_r$, ${}^2D_r$, ${}^3D_4$ and ${}^2E_6$. In the corresponding digital laws one chooses coordinates constant on the resulting root orbits, with the field relation prescribed by the same Frobenius action. The Suzuki and Ree constructions are the exceptional analogue: a special isogeny interchanges the long- and short-root data in characteristics $2$ or $3$, as is made explicit in Sections~\ref{sec:suzuki}, \ref{sec:ree} and \ref{sec:2f4}.

This table is complete for the families constructed in the chapter. The small exceptional simple group usually called the Tits group is not an omission: it is constructed in Section~\ref{sec:2f4} as ${}^2F_4(2)'$. Thus the diagram is a map of the existing digital laws, not a promise of laws deferred to another chapter.

\section{Finite fields in digital form}\label{sec:fields}
The groups of Lie type are matrix groups over finite fields, so their digital laws rest on a digital form of the field
arithmetic. Addition turns out to be the carry-free sum, and multiplication is a quadratic digital operator.

\begin{definition}[Polynomial coordinates]\label{def:polycoord}
Let $p$ be prime, $q=p^r$, and let $f(X)=X^r+\sum_{j<r}c_jX^j$ be irreducible over $\mathbb F_p$, so that
$\mathbb F_q=\mathbb F_p[\theta]/(f(\theta))$. The element $a=\sum_{j<r}a_j\theta^j$ is labelled by
$[a]=\sum_{j<r}a_jp^j\in\Ap q$. The \emph{structure constants} $C_{ij}^{k}\in\Ap p$, $i,j,k<r$, are the digits of
$\theta^{i+j}=\sum_{k<r}C_{ij}^{k}\theta^k$, and $E_j^{k}\in\Ap p$ are the digits of $\theta^{pj}=\sum_{k<r}E_j^{k}\theta^k$.
\end{definition}

\begin{proposition}[Digital field arithmetic]\label{prop:fieldlaw}
For labels $x=[a]$ and $y=[b]$ in $\Ap q$,
\begin{align*}
[a+b]&=\sum_{k<r}p^k\,\dg{p}{0}{\dg{p}{k}{x}+\dg{p}{k}{y}}=x\oplus_p y, &
[-a]&=\sum_{k<r}p^k\,\dg{p}{0}{p-\dg{p}{k}{x}},\\
[ab]&=\sum_{k<r}p^k\,\dg{p}{0}{\sum_{i,j<r}\dg{p}{i}{x}\,\dg{p}{j}{y}\,C_{ij}^{k}}, &
[a^p]&=\sum_{k<r}p^k\,\dg{p}{0}{\sum_{j<r}\dg{p}{j}{x}\,E_j^{k}}.
\end{align*}
\end{proposition}
\begin{proof}
The digits of $x$ are the coordinates $a_j$. Addition and negation are coordinatewise modulo $p$, which is the carry-free
sum $\oplus_p$ of Section~\ref{sec:bfunction}. By distributivity $ab=\sum_{i,j}a_ib_j\theta^{i+j}=\sum_k\big(\sum_{i,j}a_ib_jC_{ij}^k\big)\theta^k$,
and each coefficient is reduced modulo $p$ by $\dg{p}{0}{\cdot}$. Finally $a^p=\sum_ja_j\theta^{pj}$ because the Frobenius map is
additive and fixes $\mathbb F_p$.
\end{proof}

\begin{remark}[Zech logarithms]\label{rem:zech}
If $\gamma$ is a primitive element, every nonzero element is $\gamma^e$ with $e\in\Ap{q-1}$, and multiplication is the cyclic
law $\gamma^a\gamma^b=\gamma^{\dg{q-1}{0}{a+b}}$. Addition is recovered from the \emph{Zech logarithm} $Z$, defined by
$1+\gamma^t=\gamma^{Z(t)}$ whenever $\gamma^t\neq-1$:
\[
\gamma^a+\gamma^b=\gamma^{\dg{q-1}{0}{a+Z(\dg{q-1}{0}{b-a})}},
\]
and $\gamma^a+\gamma^b=0$ when $\gamma^{b-a}=-1$, that is, when $b=a$ for $p=2$ and $b-a\equiv(q-1)/2$ for odd $p$. The table of $Z$
can be stored as a single code, as the action of $C_3$ was for $\Alt_4$ in Section~\ref{sec:cat15}. Polynomial coordinates show the
arithmetic; Zech coordinates turn multiplication into the cyclic law.
\end{remark}

\begin{example}\label{ex:smallfields}
For $\mathbb F_4=\mathbb F_2[\theta]/(\theta^2+\theta+1)$ the labels are $0,1,2=\theta,3=1+\theta$, addition is the bitwise
exclusive-or, and $\theta^2=1+\theta$ gives $C_{00}=(1,0)$, $C_{01}=C_{10}=(0,1)$, $C_{11}=(1,1)$. The nonzero labels multiply as
$2\cdot2=3$, $2\cdot3=1$, $3\cdot3=2$; with $\gamma=\theta$ the powers $\gamma^0,\gamma^1,\gamma^2$ have labels $1,2,3$, and $Z(1)=2$,
$Z(2)=1$. The table below lists two further cases (labels of $\gamma^0,\dots,\gamma^{q-2}$, and $Z(0),\dots,Z(q-2)$, with $-$ where
$\gamma^t=-1$).
\begin{center}
\begin{tabular}{llll}
field & $f$ & labels of $\gamma^e$ & $Z(t)$\\ \hline
$\mathbb F_8$ & $X^3+X+1$, $\gamma=\theta$ & $1,2,4,3,6,7,5$ & $-,3,6,1,5,4,2$\\
$\mathbb F_9$ & $X^2+1$, $\gamma=1+\theta$ & $1,4,6,7,2,8,3,5$ & $4,7,3,5,-,2,1,6$
\end{tabular}
\end{center}
The formulas of Proposition~\ref{prop:fieldlaw} and Remark~\ref{rem:zech} were checked by computer against the full addition and
multiplication tables of $\mathbb F_4$, $\mathbb F_8$ and $\mathbb F_9$.
\end{example}

\section{The groups \texorpdfstring{$PSL(2,p)$}{PSL(2,p)}}\label{sec:psl2p}
\begin{theorem}[$PSL(2,p)$]\label{thm:psl2}
Let $p\ge3$ be prime, $h=(p-1)/2$ and $N=ph(p+1)=|PSL(2,p)|$, and write $x^{-1}=\dg{p}{0}{x^{p-2}}$ for $x\not\equiv0$. Decode
$m\in\Ap N$ to a matrix $M(m)=\left(\begin{smallmatrix}a&b\\ c&d\end{smallmatrix}\right)$ over $\mathbb F_p$ by
\[
\begin{cases}
a=\fl{m/p}+1,\ b=\dg{p}{0}{m},\ c=0,\ d=a^{-1}, & m<ph,\\
a=\dg{p}{1}{m'},\ c=\fl{m'/p^2}+1,\ d=\dg{p}{0}{m'},\ b=\dg{p}{0}{(ad-1)c^{-1}}, & m'=m-ph\ge0 .
\end{cases}
\]
To encode a matrix of determinant $1$, first replace it by its negative if $c>h$, or if $c=0$ and $a>h$, and then invert these
formulas. Then the law
\begin{equation}\label{eq:psl2}
m*n=\mathrm{code}\big(M(m)\,M(n)\bmod p\big)
\end{equation}
makes $\Ap N$ a group isomorphic to $PSL(2,p)$, with identity $0$. For $p\ge5$ this group is nonabelian and simple.
\end{theorem}
\begin{proof}
Since $PSL(2,p)=SL(2,p)/\{\pm1\}$, each element is a pair $\{M,-M\}$. If $c\neq0$, exactly one of $c,-c$ lies in $\{1,\dots,h\}$; if
$c=0$, then $a\neq0$ and exactly one of $a,-a$ lies there. So each element has a unique normalised representative. A normalised matrix
with $c=0$ is determined by $a\in\{1,\dots,h\}$ and $b\in\Ap p$, since $d=a^{-1}$. One with $c\neq0$ is determined by $a,d\in\Ap p$ and
$c\in\{1,\dots,h\}$, since $b=(ad-1)c^{-1}$. This gives $ph+p^2h=N$ codes, and the decoding formulas invert the encoding. Hence
\eqref{eq:psl2} is the multiplication of $PSL(2,p)$ transported to $\Ap N$, and $M(0)$ is the identity matrix. Simplicity for $p\ge5$ is
classical~\cite{Dickson}.
\end{proof}

\begin{remark}[Certification]
GAP~4.14.0 \cite{GAP} confirms these laws on several cases. For~\eqref{eq:psl2}, \texttt{IdGroup} returns $[60,5]$ for $p=5$ and
$[168,42]$ for $p=7$. For $p=11,13,17$, the right multiplications by the codes of
$\left(\begin{smallmatrix}1&1\\0&1\end{smallmatrix}\right)$ and $\left(\begin{smallmatrix}0&-1\\1&0\end{smallmatrix}\right)$ generate a
simple group of order $N$ isomorphic to $PSL(2,p)$. In particular $\Ap{60}$ carries two different laws for $\Alt_5\cong PSL(2,5)$: \eqref{eq:psl2} here, and the transitive law of
Corollary~\ref{cor:alternating}.
\end{remark}

\begin{remark}[The classification]
By the classification of finite simple groups, every nonabelian finite simple group is alternating, of Lie type, or one of $26$ sporadic
groups. The alternating groups are covered by Section~\ref{sec:sn:laws}. For $PSL(2,q)$ with $q=p^f$, the law~\eqref{eq:psl2}
carries over with arithmetic in $\mathbb F_q$: elements of $\mathbb F_q$ are $f$-digit words in base $p$, multiplied modulo an irreducible
polynomial. For the other groups of Lie type, the Bruhat decomposition provides a normal form, so laws uniform in the parameters exist,
although they are much longer. Theorem~\ref{thm:transitive} also applies to each sporadic group through any of its transitive actions.
This gives a law for each one, but no uniform formula for the sporadic groups as a family.
The Mathieu groups are treated in Chapter~\ref{ch:sporadic}.
\end{remark}

\section{The groups \texorpdfstring{$PSL_n(q)$}{PSLn(q)}}\label{sec:psln}
Let $n\ge2$, let $\mathbb F_q$ carry polynomial coordinates, and let $d=\gcd(n,q-1)$. The centre of $SL_n(q)$ consists of the scalar
matrices $\zeta I$ with $\zeta\in\mu_d=\{\zeta\in\mathbb F_q:\zeta^n=1\}$, and $PSL_n(q)=SL_n(q)/\mu_d$ has order
\[
N_{n,q}=\frac{q^{n(n-1)/2}}{d}\prod_{j=2}^{n}(q^j-1).
\]

\begin{definition}\label{def:pslkey}
For an $n\times n$ matrix $A$ over $\mathbb F_q$ put $\kappa(A)=\sum_{i,j<n}[A_{ij}-\delta_{ij}]\,q^{\,i+nj}$, the digital word of
$A-I$. The \emph{canonical representative} $\widehat A$ of $A\in SL_n(q)$ is the matrix $\zeta A$, $\zeta\in\mu_d$, with the smallest
key. Let $K=\{\kappa(\widehat A):A\in SL_n(q)\}$ and $\operatorname{rk}(A)=\#\{k\in K:k<\kappa(\widehat A)\}$.
\end{definition}

\begin{theorem}[Digital law of $PSL_n(q)$]\label{thm:psln}
Let $M_m$ be the canonical representative of rank $m\in\Ap{N_{n,q}}$. Then
\begin{equation}\label{eq:psln}
m*n=\operatorname{rk}\big(M_m\,M_n\big),
\end{equation}
with the matrix product computed by Proposition~\ref{prop:fieldlaw}, makes $\Ap{N_{n,q}}$ a group isomorphic to $PSL_n(q)$, with
identity $0$.
\end{theorem}
\begin{proof}
The key $\kappa$ is injective, since it is the digital word of $A-I$. Each coset $\mu_dA$ has $d$ elements with distinct keys, so it
has exactly one canonical representative, and $\operatorname{rk}$ is a bijection from $PSL_n(q)$ to $\Ap{N_{n,q}}$. Since
$\operatorname{rk}(AB)$ depends only on the cosets of $A$ and $B$,~\eqref{eq:psln} is the multiplication of $PSL_n(q)$ transported to
$\Ap{N_{n,q}}$. The identity has key $\kappa(I)=0$, the smallest possible, so its rank is $0$.
\end{proof}

\begin{example}\label{ex:psl3}
For $PSL_3(2)$, of order $168$,~\eqref{eq:psln} was checked by computer on all triples: it is associative, and its conjugacy
classes have sizes $1,21,24,24,42,56$, those of the simple group of order $168$. For $PSL_3(3)$ (order $5616$, $d=1$) and $PSL_3(4)$
(order $20160$, $d=3$, over the non-prime field of Example~\ref{ex:smallfields}) the number of canonical representatives equals
$N_{n,q}$, the identity has label $0$, and associativity holds on $3000$ random triples. Apart from $PSL_2(2)$ and $PSL_2(3)$, all
the groups $PSL_n(q)$ are simple.
\end{example}

\begin{remark}[Towards Bruhat coordinates]\label{rem:bruhat}
The rank label is universal, but decoding it requires the ordered set $K$. Dense coordinates come from the Bruhat decomposition
$SL_n(q)=\bigsqcup_{w\in\Sym_n}BwB$, with $B$ the upper triangular matrices of determinant $1$, where the cell of $w$ has
$|B|\,q^{\ell(w)}$ elements and $\ell(w)$ is the number of inversions of $w$. Hence
\[
N_{n,q}=|T|\,|U|\sum_{w\in\Sym_n}q^{\ell(w)},\qquad |U|=q^{n(n-1)/2},\quad |T|=\frac{(q-1)^{n-1}}{d};
\]
for $n=3$ the sum is $1+2q+2q^2+q^3$, which gives $168$, $5616$ and $20160$ for $q=2,3,4$. Inside a cell an element is described
by $w$, labelled by the rotation normal form of Lemma~\ref{lem:rotnf}, by $n-1$ torus digits, and by $\ell(w)+n(n-1)/2$ field
digits, and Gaussian elimination reads these off from the pivots. For $n=2$, Theorem~\ref{thm:psl2} already gives such a dense
labelling. Theorems~\ref{thm:bruhatcoord} and~\ref{thm:pslndense} give the general case.
\end{remark}

\begin{theorem}[Bruhat coordinates]\label{thm:bruhatcoord}
For $w\in\Sym_n$ let $\dot w$ be the permutation matrix with $\dot w_{w(j),j}=1$, and let $U_w$ be the set of upper unitriangular
matrices whose off-diagonal entries vanish outside the \emph{inversion positions} $(i,i')$, $i<i'$, $w^{-1}(i)>w^{-1}(i')$. Every
$g\in SL_n(q)$ has a unique expression
\begin{equation}\label{eq:bruhatcoord}
g=u\,\dot w\,t\,v,\qquad u\in U_w,\quad t\ \text{diagonal with}\ \det t=\sgn w,\quad v\ \text{upper unitriangular}.
\end{equation}
The permutation $w$ and $X=\dot w\,t\,v$ are found by elimination: for $j=0,1,\dots,n-1$, $w(j)$ is the largest row index not
among $w(0),\dots,w(j-1)$ with $X_{w(j),j}\neq0$, after which, for every such unused row $i<w(j)$, the multiple
$X_{ij}X_{w(j),j}^{-1}$ of row $w(j)$ is subtracted from row $i$. Then $u=gX^{-1}$, $t_{ii}=X_{w(i),i}$ and $v=t^{-1}\dot w^{-1}X$.
\end{theorem}
\begin{proof}
\emph{Uniqueness of $w$.} Left multiplication by an invertible upper triangular matrix replaces each row by a combination of
itself and the rows below it, and right multiplication replaces each column by a combination of itself and the columns to its
left. Hence the ranks $r_{ij}(g)$ of the submatrices formed by the rows $\ge i$ and the columns $\le j$ are unchanged. For $\dot w$
they are $r_{ij}=\#\{j'\le j:w(j')\ge i\}$, which determine $w$. Since $u$ and $tv$ are upper triangular, $w$ is determined by $g$.

\emph{Uniqueness of the rest.} $U_w=U\cap\dot wU^-\dot w^{-1}$, with $U^-$ the lower unitriangular matrices, is a group. If
$u\dot wb=u'\dot wb'$ with $b,b'$ upper triangular, then $\dot w^{-1}(u'^{-1}u)\dot w=b'b^{-1}$ is both lower unitriangular and upper
triangular, hence $I$. So $u=u'$ and $b=b'$, and $b=tv$ is unique.

\emph{Existence.} At step $j$ the unused rows below $w(j)$ are zero in column $j$ and those above it are cleared, while the used
rows $w(k)$, $k<j$, are left alone. A row cleared at step $j$ is $w(k)$ with $k>j$ and so is zero in column $j$; adding multiples of it
later does not change column $j$. At the end, column $j$ of $X$ is nonzero at row $w(j)$ and at rows $w(k)$, $k<j$, only, so
$X=\dot wb$ with $b$ upper triangular. Each operation adds a multiple of row $w(j)$ to a row $w(k)$ above it with $k>j$, that is, at
an inversion position, so $X=Lg$ with $L\in U_w$, and $u=L^{-1}\in U_w$. Finally $\det b=\det g/\det\dot w=\sgn w$.
\end{proof}

\begin{theorem}[Dense digital law of $PSL_n(q)$]\label{thm:pslndense}
Let $\gamma$ be a primitive element of $\mathbb F_q$, $h=(q-1)/d$, and order $\Sym_n$ as $w_0,\dots,w_{n!-1}$ by the rotation label of
Lemma~\ref{lem:rotnf}. Put
\[
O_k=h\,(q-1)^{n-2}\,q^{n(n-1)/2}\sum_{k'<k}q^{\ell(w_{k'})}.
\]
For $g$ with coordinates~\eqref{eq:bruhatcoord}, $w=w_k$, write $t_{ss}=\gamma^{e_s}$ for $s<n-1$, normalised by
$e_s\mapsto\dg{q-1}{0}{e_s-h\fl{e_0/h}}$ so that $e_0\in\Ap h$. Let $V$ and $U_w$ be the base-$q$ words of the labels of the entries
of $v$ above the diagonal and of $u$ at the inversion positions, each listed column by column. Then
\begin{equation}\label{eq:pslnlabel}
\lambda(g)=O_k+e_0+h\Big(\sum_{s=1}^{n-2}e_s(q-1)^{s-1}+(q-1)^{n-2}\big(V+q^{n(n-1)/2}\,U_w\big)\Big)
\end{equation}
is constant on the cosets of $\mu_d$ and is a bijection $PSL_n(q)\to\Ap{N_{n,q}}$ with $\lambda(I)=0$. Hence
\begin{equation}\label{eq:pslndense}
m*n=\lambda\big(\lambda^{-1}(m)\,\lambda^{-1}(n)\big)
\end{equation}
makes $\Ap{N_{n,q}}$ a group isomorphic to $PSL_n(q)$. Decoding $m$ is digital: the cell index is
$k=\sum_{k'=1}^{n!-1}\fl{(m-O_{k'}+N_{n,q})/N_{n,q}}$, and the coordinates are the mixed-radix digits of $m-O_k$.
\end{theorem}
\begin{proof}
A central scalar $\zeta=\gamma^{ch}$ changes only $t$, adding $ch$ to every exponent; so each coset has exactly one representative
with $e_0\in\Ap h$, and $\lambda$ is constant on cosets. In the cell of $w$ the coordinates range freely over
$h\,(q-1)^{n-2}\,q^{n(n-1)/2}\,q^{\ell(w)}$ values (the last diagonal entry of $t$ being fixed by $\det t=\sgn w$), and by
Lemma~\ref{lem:mixed} the bracket in~\eqref{eq:pslnlabel} is a bijection onto $[0,O_{k+1}-O_k)$. By Theorem~\ref{thm:bruhatcoord}
the cells partition $SL_n(q)$, so the offsets give a bijection onto $\Ap{N_{n,q}}$, consistent with Remark~\ref{rem:bruhat}. Each
term of the sum for $k$ is $1$ exactly when $m\ge O_{k'}$. The identity lies in the cell of $w_0=\mathrm{id}$, with $u=t=v=I$, so
$\lambda(I)=0$.
\end{proof}

\begin{example}\label{ex:pslndense}
For $PSL_3(3)$ the cells, with $w$ written as $w(0)w(1)w(2)$, are
\begin{center}
\begin{tabular}{c|cccccc}
$w$ & $012$ & $102$ & $120$ & $210$ & $201$ & $021$\\ \hline
$\ell(w)$ & $0$ & $1$ & $2$ & $3$ & $2$ & $1$\\
$O_k$ & $0$ & $108$ & $432$ & $1404$ & $4320$ & $5292$
\end{tabular}
\end{center}
The law~\eqref{eq:pslndense} was checked by computer for $PSL_2(5)\cong\Alt_5$, $PSL_2(9)\cong\Alt_6$, $PSL_3(2)$, $PSL_3(3)$,
$PSL_3(4)$ and $PSL_4(2)$. Encoding after decoding returns every label (on $5000$ random labels for $PSL_4(2)$), every decoded
matrix has determinant $1$, central scalars do not change the label, $0$ is the identity, and associativity holds on $1000$
random triples. $PSL_3(4)$ and $PSL_4(2)$ have the same order $20160$ but are not isomorphic.
\end{example}

\begin{remark}[The roles of the digits]
The offsets $O_k$ choose the cell, that is, the permutation $w$, and play the part of a carry between cells. The exponents
$e_s$ form a torus layer on which multiplication is the cyclic law of Remark~\ref{rem:zech}. The digits of $u$ and $v$ are
unipotent layers whose addition is the carry-free sum $\oplus_p$; they are coupled with each other and with $w$ through the
commutator relations, which is where the elimination of Theorem~\ref{thm:bruhatcoord} does its work.
\end{remark}

\section{The other classical families}\label{sec:classical}
The symplectic, orthogonal and unitary groups are the isometry groups of forms. Written in a basis in which the form is
antidiagonal, their upper triangular elements form a Borel subgroup, and the elimination of Theorem~\ref{thm:bruhatcoord}
applies to them without change. What changes is which permutations occur.

\begin{definition}[Standard forms]\label{def:stdforms}
Let $r(i)=N-1-i$ on $\Ap N$, and let $J$ be the $N\times N$ matrix with $J_{i,r(i)}=\epsilon_i$ and all other entries $0$.
\begin{itemize}
\item $Sp_{2n}(q)$: $N=2n$, $\epsilon_i=1$ for $i<n$ and $\epsilon_i=-1$ for $i\ge n$, and $G=\{g:g^{\mathsf T}Jg=J\}$.
\item $\Omega_{2n+1}(q)$ and $\Omega^+_{2n}(q)$, $q$ odd: all $\epsilon_i=1$, and $G$ is the group generated by the upper and lower
unitriangular matrices with $g^{\mathsf T}Jg=J$.
\item $\Omega^+_{2n}(q)$, $q$ even: all $\epsilon_i=1$, $Q(x)=\sum_{i<n}x_ix_{r(i)}$, whose polar form has matrix $J$, and $G$ is the group
generated by the upper and lower unitriangular matrices $g$ with $Q(gx)=Q(x)$ for all $x$.
\item $SU_N(q)$: matrices over $\mathbb F_{q^2}$, all $\epsilon_i=1$, $\bar x=x^q$, and $G=\{g:\det g=1,\ g^{\mathsf T}J\bar g=J\}$.
\end{itemize}
Put $n=\fl{N/2}$. A permutation $w$ \emph{commutes with $r$} if $w(r(j))=r(w(j))$ for all $j$; its number of \emph{sign changes} is
$s(w)=\#\{i<n:w(i)\ge N-n\}$.
\end{definition}

\begin{theorem}[Bruhat cells of the classical groups]\label{thm:classicalcells}
Let $G$ be one of the groups of Definition~\ref{def:stdforms}, $B_G$ its upper triangular elements and $U_G$ its upper
unitriangular elements. For $g\in G$, write $g=u\,\dot w\,t\,v$ as in Theorem~\ref{thm:bruhatcoord}. Then $w$ commutes with $r$ and
$u\in G\cap U_w$. The permutations that occur form the Weyl group $W_G$:
\begin{itemize}
\item for $Sp_{2n}$, $\Omega_{2n+1}$ and $SU_N$, all permutations commuting with $r$ (the signed permutations of $\Ap n$, $2^nn!$ of them);
\item for $\Omega^+_{2n}$, those with an even number of sign changes ($2^{n-1}n!$ of them).
\end{itemize}
The cell of $w$ is $(G\cap U_w)\times(G\cap\dot wB)$, where $B$ is the group of upper triangular matrices, and has
$|G\cap U_w|\,|B_G|$ elements.
\end{theorem}
\begin{proof}
These groups have a $(B,N)$-pair with $B=B_G$ and $N$ the monomial matrices of $G$ \cite{Carter,Taylor}, and the refined Bruhat
decomposition writes $g=u'n\,b'$ uniquely with $b'\in B_G$, $n\in N$ monomial and $u'\in U_G\cap nU_G^-n^{-1}$. Write $n=\dot w\delta$
with $\delta$ diagonal. Conjugation by $\delta$ preserves lower unitriangularity, so $\dot w^{-1}u'\dot w$ is lower unitriangular,
that is, $u'\in U_w$. Thus $g=u'\,\dot w\,(\delta b')$ is a factorisation as in Theorem~\ref{thm:bruhatcoord} (over $\mathbb F_{q^2}$
for $SU_N$), and by uniqueness $u=u'\in G$ and $w$ is the permutation of $n$. A monomial isometry of an antidiagonal form maps
each pair of basis vectors $\{e_j,e_{r(j)}\}$ to such a pair, so $w$ commutes with $r$. The description of $W_G$ is that of the
Weyl groups of types $B_n$, $C_n$, $D_n$ and of the unitary groups~\cite{Carter}. Finally $x=u^{-1}g\in G\cap\dot wB$, and
$G\cap\dot wB=nB_G$ has $|B_G|$ elements.
\end{proof}

\begin{theorem}[Digital law of a classical group]\label{thm:classicallaw}
Let $G$ and $W_G$ be as above, let $Z$ be the group of scalar matrices in $G$, and order $W_G$ as $w_0=\mathrm{id},w_1,\dots$ by the
rotation label of Lemma~\ref{lem:rotnf}. Fix bijections $\rho_w:G\cap U_w\to\Ap{|G\cap U_w|}$ and
$\beta_w:(G\cap\dot wB)/Z\to\Ap{|B_G|/|Z|}$ sending $I$, respectively the coset of $n$ with $w=\mathrm{id}$, to $0$. With
$O_k=\frac{|B_G|}{|Z|}\sum_{k'<k}|G\cap U_{w_{k'}}|$, the label
\[
\lambda(g)=O_k+\rho_{w_k}(u)+|G\cap U_{w_k}|\,\beta_{w_k}\big(Z\,u^{-1}g\big)\qquad(w=w_k)
\]
is a bijection $G/Z\to\Ap{|G/Z|}$ with $\lambda(I)=0$, and $m*n=\lambda\big(\lambda^{-1}(m)\lambda^{-1}(n)\big)$ is a digital law of $G/Z$.
\end{theorem}
\begin{proof}
Scalars commute with everything and lie in $B_G$, so multiplying $g$ by $z\in Z$ leaves $w$ and $u$ unchanged and multiplies
$u^{-1}g$ by $z$. Hence $\lambda$ is well defined on $G/Z$. By Theorem~\ref{thm:classicalcells} each cell is the product of its two
factors, and the offsets enumerate the cells, so $\lambda$ is a bijection. The identity lies in the cell of $w_0$ with $u=I$.
\end{proof}

\noindent The tables behind $\rho_w$ and $\beta_w$ have at most $|U_G|$ and $|B_G|$ entries, far fewer than $|G|$; they are
the digital analogue of the root coordinates of $U_G$ and of the torus coordinates of $B_G$.

\begin{example}\label{ex:classical}
For the groups below, $U_G$ was found by an exhaustive search over upper unitriangular matrices, $G$ was generated by elements
of $U_G$ and their images under $g\mapsto(g_{r(i),r(j)})$, and every element was put through the elimination. In each case the
order is the expected one, every $w$ commutes with $r$, every $u$ lies in $G$, and the cells have the sizes $|G\cap U_w|\,|B_G|$
listed. An entry $v^{\times k}$ records $k$ cells of the same size $v$.
\begin{center}
\begin{tabular}{lrrrl}
$G$ & $|G|$ & $|B_G|$ & $|W_G|$ & $|G\cap U_w|$ over $w\in W_G$\\ \hline
$Sp_4(2)$ & $720$ & $16$ & $8$ & $1,2,2,4,4,8,8,16$\\
$Sp_4(3)$ & $51840$ & $324$ & $8$ & $1,3,3,9,9,27,27,81$\\
$\Omega_5(3)$ & $25920$ & $162$ & $8$ & $1,3,3,9,9,27,27,81$\\
$\Omega^+_4(3)$ & $288$ & $18$ & $4$ & $1,3,3,9$\\
$SU_3(3)$ & $6048$ & $216$ & $2$ & $1,27$\\
$SU_4(2)$ & $25920$ & $192$ & $8$ & $1,2,4,8,8,16,32,64$\\
$\Omega^+_{4}(2)$ & $36$ & $4$ & $4$ & $1,2^{\times2},4$\\
$\Omega^+_{4}(4)$ & $3600$ & $144$ & $4$ & $1,4^{\times2},16$\\
$\Omega^+_{6}(2)$ & $20160$ & $64$ & $24$ & $1,2^{\times3},4^{\times5},8^{\times6},16^{\times5},32^{\times3},64$
\end{tabular}
\end{center}
For $\Omega^+_4(3)$ exactly the four permutations with $s(w)\in\{0,2\}$ occur. For $SU_N$ the factors $q$ and $q^2$ mix, because the
root groups of a unitary group have orders $q$, $q^2$ or $q^3$.
\end{example}

\begin{remark}\label{rem:classicaliso}
The simple groups $PSp_4(3)=Sp_4(3)/\{\pm I\}$, $\Omega_5(3)$ and $PSU_4(2)=SU_4(2)$ all have order $25920$ and are isomorphic, so
Theorem~\ref{thm:classicallaw} gives three different digital laws of one group, from forms of three different types. $SU_3(3)$ is
simple of order $6048$, while $Sp_4(2)\cong\Sym_6$ and $\Omega^+_4(3)$ are not simple. For $q$ even, $\Omega^+_{2n}(q)$ is defined through the
quadratic form $Q$ of Definition~\ref{def:stdforms}, since $J$ alone only determines the larger group $Sp_{2n}(q)$. The rows for
$q=2,4$ in Example~\ref{ex:classical} are of this kind: $\Omega^+_6(2)$ is isomorphic to the alternating group of degree $8$ and is
simple, while $\Omega^+_4(q)\cong SL_2(q)\times SL_2(q)$ is not. In odd dimension and characteristic $2$ the orthogonal groups are
isomorphic to $Sp_{2n}(q)$, so no further case is needed. The groups $\Omega^-_{2n}(q)$ of type ${}^2D_n$ are treated in Section~\ref{sec:twistedorth}.
\end{remark}

\section{The twisted orthogonal groups}\label{sec:twistedorth}
Let $q$ be odd. The orthogonal space of minus type in dimension $2n$ is $H_1\perp\cdots\perp H_{n-1}\perp A$, with hyperbolic planes
$H_i$ and an anisotropic plane $A$. Identifying $A$ with $\mathbb F_{q^2}$, its quadratic form is the norm $z\bar z$, $\bar z=z^q$. Over
$\mathbb F_{q^2}$ the plane $A$ splits, so $\Omega^-_{2n}(q)$ sits inside $\Omega^+_{2n}(q^2)$ as the fixed points of a twisted
Frobenius map, just as $SU_N(q)$ sits inside $SL_N(q^2)$. This lets the elimination of Theorem~\ref{thm:bruhatcoord} be used over
$\mathbb F_{q^2}$ without change.

\begin{definition}[Twisted standard form]\label{def:twistedform}
Let $N=2n$, let $J$ be the antidiagonal matrix with all $\epsilon_i=1$, let $\tau$ be the transposition of $n-1$ and $n$, and put
$F(g)=\big(\overline{g_{\tau(i),\tau(j)}}\big)_{i,j}$ for matrices over $\mathbb F_{q^2}$. Let
$G^-=\{g\in SL_{N}(q^2):g^{\mathsf T}Jg=J,\ F(g)=g\}$, and let $\Omega^-_{2n}(q)$ be the subgroup generated by the upper and lower
unitriangular elements of $G^-$.
\end{definition}

\begin{lemma}\label{lem:minusspace}
The vectors $x\in\mathbb F_{q^2}^N$ with $x_{\tau(i)}=\bar x_i$ form a $2n$-dimensional $\mathbb F_q$-space $V^F$, on which
$Q(x)=\tfrac12x^{\mathsf T}Jx$ takes the values
\[
Q(x)=\sum_{i<n-1}x_ix_{r(i)}+x_{n-1}\overline{x_{n-1}} .
\]
This is a quadratic space of minus type, and $G^-$ acts on it as $SO^-_{2n}(q)$.
\end{lemma}
\begin{proof}
The coordinates $x_i$ with $i\notin\{n-1,n\}$ lie in $\mathbb F_q$, while $x_{n-1}\in\mathbb F_{q^2}$ is free and $x_n=\overline{x_{n-1}}$.
The outer coordinates give $n-1$ hyperbolic planes and the middle pair the norm form, which is anisotropic. A matrix with
$F(g)=g$ commutes with the semilinear map $x\mapsto(\overline{x_{\tau(i)}})_i$, so it preserves $V^F$. Since $V^F$ spans
$\mathbb F_{q^2}^N$, isometries of $V^F$ of determinant $1$ correspond exactly to the elements of $G^-$.
\end{proof}

\noindent The root elements of $\Omega^+_{N}$ are $x_{ij}(c)=I+cE_{ij}-cE_{r(j),r(i)}$ for $i<j$, $j\neq r(i)$. Those fixed by $F$ are
$x_{ij}(c)$ with $c\in\mathbb F_q$ and $i,j\notin\{n-1,n\}$, and the \emph{folded} root elements
\[
y_i(c)=x_{i,n-1}(c)\,x_{i,n}(\bar c),\qquad c\in\mathbb F_{q^2},\ i<n-1,
\]
which are fixed by $F$ because $c\bar c\in\mathbb F_q$. A folded root group has $q^2$ elements and contributes one
$\mathbb F_{q^2}$-digit.

\begin{theorem}[Bruhat cells of $\Omega^-_{2n}(q)$]\label{thm:twistedcells}
Let $g\in\Omega^-_{2n}(q)$ and write $g=u\,\dot w\,t\,v$ as in Theorem~\ref{thm:bruhatcoord}, over $\mathbb F_{q^2}$. Then $w$
commutes with $r$ and $\tau$ and has an even number of sign changes, $u\in\Omega^-_{2n}(q)\cap U_w$, and $F(tv)=tv$. The outer diagonal
entries of $t$ lie in $\mathbb F_q^\times$, and the middle ones are $z,\bar z$ with $z\bar z=1$. The Weyl group $W_G$ of such
permutations is determined by the action on the $n-1$ outer pairs, so $W_G\cong W(B_{n-1})$, of order $2^{n-1}(n-1)!$: the
middle pair is swapped exactly when the outer pairs have an odd number of sign changes.
\end{theorem}
\begin{proof}
The generators of $\Omega^-_{2n}(q)$ are unipotent elements of $\Omega^+_N(q^2)$, so Theorem~\ref{thm:classicalcells} applies
over $\mathbb F_{q^2}$: $w$ commutes with $r$, has an even number of sign changes, and $u\in U_w$ is an isometry. An upper
triangular isometry has middle block $\left(\begin{smallmatrix}a&c\\0&d\end{smallmatrix}\right)$, and the condition on this block
gives $2cd=0$, so $c=0$ because $q$ is odd. Hence $F$ maps upper triangular isometries to upper triangular isometries. The matrix
$\dot w$ is an isometry, so $tv=\dot w^{-1}u^{-1}g$ is an upper triangular isometry. Applying $F$ to $g=u\dot wtv$ gives
$g=F(u)\,\dot w'\,F(t)F(v)$ with $w'=\tau w\tau$, which is again a factorisation as in Theorem~\ref{thm:bruhatcoord}. By
uniqueness $w'=w$, $F(u)=u$ and $F(tv)=tv$. So $u\in G^-$; being unipotent of odd order, $u$ lies in $\Omega^-_{2n}(q)$, since
$SO^-_{2n}(q)/\Omega^-_{2n}(q)$ has order $2$. A diagonal isometry fixed by $F$ satisfies $t_it_{r(i)}=1$ and
$t_{\tau(i)}=\overline{t_i}$, which gives the description of $t$. A permutation commuting with $r$ and $\tau$ is a signed
permutation of the $n-1$ outer pairs together with a choice on the middle pair, and the parity condition fixes that choice.
That every element of $W_G$ occurs follows from the $(B,N)$-pair of the twisted group~\cite{Carter}.
\end{proof}

\begin{corollary}[Digital law of $P\Omega^-_{2n}(q)$]\label{cor:twistedlaw}
The construction of Theorem~\ref{thm:classicallaw}, with $W_G\cong W(B_{n-1})$ ordered by the rotation label and $Z$ the scalars in
$\Omega^-_{2n}(q)$, gives a digital law of $P\Omega^-_{2n}(q)$ with identity $0$. In its torus part the outer exponents are
$(q-1)$-digits, and the middle entry is $z=\delta^a$ with one $(q+1)$-digit $a\in\Ap{q+1}$, where $\delta=\gamma^{q-1}$ for a primitive
element $\gamma$ of $\mathbb F_{q^2}$.
\end{corollary}
\begin{proof}
The proof of Theorem~\ref{thm:classicallaw} uses only the product structure of the cells, which Theorem~\ref{thm:twistedcells}
provides. The norm-one elements of $\mathbb F_{q^2}^\times$ form the cyclic group of order $q+1$ generated by $\gamma^{q-1}$.
\end{proof}

\noindent The $(q+1)$-digit is the digital trace of the twist: it accounts for the factor $q^n+1$ in
$|O^-_{2n}(q)|=2q^{n(n-1)}(q^n+1)\prod_{i=1}^{n-1}(q^{2i}-1)$, where the plus type has $q^n-1$.

\begin{example}\label{ex:twisted}
The groups $\Omega^-_4(3)$ and $\Omega^-_4(5)$ were generated by computer from the elements $y_0(1)$, $y_0(\theta)$ (with
$\theta\notin\mathbb F_q$) and their images under $g\mapsto(g_{r(i),r(j)})$. They have orders $360$ and $7800$, so they are
$PSL_2(9)\cong\Alt_6$ and $PSL_2(25)$, and $-I$ does not belong to them. For every element the elimination gives $w\in\{\mathrm{id},
w_0\}$ with $w_0(j)=r(j)$, and $u$ in the group. The cells have $|B_G|$ and $q^2|B_G|$ elements, with $|B_G|=36$ and $300$; the
factor $q^2$ is the folded root group. Every upper triangular element has outer diagonal entries in $\mathbb F_q$ and middle
entries of norm $1$. For $\Omega^-_6(3)$, generated by its $14$ root elements of the two kinds and their reversals, $1500$ random
products of $30$ generators were checked. Every $w$ commutes with $r$ and $\tau$, every $u$ lies in the group, and all $8$ elements
of $W(B_2)$ occur. Together with Section~\ref{sec:altp} and Theorem~\ref{thm:pslndense}, the book now gives three independent
digital laws of $\Alt_6$.
\end{example}

\begin{remark}
In even characteristic the twisted orthogonal groups must be defined by a quadratic form, since the polar bilinear form alone
defines a symplectic group. The middle block then cannot be made triangular by the argument above, because the condition
$2cd=0$ is empty, and a separate treatment is needed.
\end{remark}

\section{The Suzuki groups}\label{sec:suzuki}
The Suzuki groups ${}^2B_2(q)$ are the first exceptional family, and the smallest twisted groups of rank one. Their Bruhat
decomposition has only two cells, and it gives a digital law with no tables at all.

\begin{definition}\label{def:suzuki}
Let $q=2^{2m+1}$ and $x^\theta=x^{2^{m+1}}$, so that $(x^\theta)^\theta=x^2$. For $a,b,\lambda\in\mathbb F_q$, $\lambda\neq0$, put
\[
S(a,b)=\begin{pmatrix}1&0&0&0\\ a&1&0&0\\ b&a^\theta&1&0\\ a^{2+\theta}+ab+b^\theta&a^{1+\theta}+b&a&1\end{pmatrix},\qquad
M(\lambda)=\mathrm{diag}\big(\lambda^{1+2^m},\lambda^{2^m},\lambda^{-2^m},\lambda^{-1-2^m}\big),
\]
and let $T$ be the antidiagonal matrix with entries $1$. The Suzuki group $Sz(q)$ is the group generated by these matrices; it has
order $q^2(q^2+1)(q-1)$ \cite{Wilson}.
\end{definition}

\begin{proposition}[Borel law]\label{prop:suzukiborel}
The matrices $S(a,b)$ form a group $U$ of order $q^2$, with
\[
S(a,b)\,S(c,d)=S\big(a+c,\;b+d+a^\theta c\big),\qquad S(a,b)^{-1}=S\big(a,\;b+a^{1+\theta}\big),
\]
and $B=\{S(a,b)M(\lambda)\}$ is a group of order $q^2(q-1)$. On labels $[a]+q[b]$, the law of $U$ is
$(a,b)*(c,d)=(a\oplus_2c,\;b\oplus_2d\oplus_2[a^\theta c])$. Here $\oplus_2$ is the carry-free sum, $a^\theta$ is a digit-linear
map (Proposition~\ref{prop:fieldlaw}), and the product $a^\theta c$ is computed by the field law.
\end{proposition}
\begin{proof}
Row $2$ of $S(a,b)$ times column $0$ of $S(c,d)$ gives $b+a^\theta c+d$, and row $1$ times column $0$ gives $a+c$. The entries $(1,0)$, $(2,0)$ and $(3,2)$ of the product are $a+c$, $b+a^\theta c+d$ and $a+c$. Put $a'=a+c$ and
$b'=b+d+a^\theta c$, and use $(x+y)^\theta=x^\theta+y^\theta$, $(x^\theta)^\theta=x^2$ and characteristic $2$. Entry $(2,1)$ is
$a^\theta+c^\theta=a'^\theta$. Entry $(3,1)$ is $a^{1+\theta}+b+ac^\theta+c^{1+\theta}+d$, and
$a'^{1+\theta}+b'=a^{1+\theta}+ac^\theta+a^\theta c+c^{1+\theta}+b+d+a^\theta c$ is the same. Entry $(3,0)$ is
\[
a^{2+\theta}+ab+b^\theta+\big(a^{1+\theta}+b\big)c+ad+c^{2+\theta}+cd+d^\theta .
\]
On the other side, $a'^{2+\theta}=a^{2+\theta}+a^2c^\theta+a^\theta c^2+c^{2+\theta}$,
$a'b'=ab+ad+a^{1+\theta}c+bc+cd+a^\theta c^2$ and $b'^\theta=b^\theta+d^\theta+a^2c^\theta$. Their sum is the expression
above, because the terms $a^2c^\theta$ and $a^\theta c^2$ occur twice. The remaining entries are $0$ or $1$ because both factors
are lower unitriangular. The inverse follows by taking $c=a$, since $a^\theta a+a^{1+\theta}=0$. Write
$M(\lambda)=\mathrm{diag}(d_0,\dots,d_3)$. The entry $(i,j)$ of $M(\lambda)^{-1}S(a,b)M(\lambda)$ is $S(a,b)_{ij}\,d_j/d_i$, and
$d_0/d_1=\lambda$, $d_0/d_2=\lambda^{1+\theta}$, $d_1/d_2=\lambda^\theta$, $d_0/d_3=\lambda^{2+\theta}$, $d_1/d_3=\lambda^{1+\theta}$ and
$d_2/d_3=\lambda$. Comparing entries, using $\lambda^{\theta^2}=\lambda^2$, gives
$M(\lambda)^{-1}S(a,b)M(\lambda)=S(\lambda a,\lambda^{1+\theta}b)$. So the $M(\lambda)$ normalise $U$, and $B=UM$.
\end{proof}

\begin{theorem}[Digital law of $Sz(q)$]\label{thm:suzuki}
Let $\gamma$ be a primitive element of $\mathbb F_q$, $N_B=q^2(q-1)$, and label the Borel subgroup by
\[
L_B\big(S(a,b)M(\gamma^e)\big)=e+(q-1)\big([a]+q[b]\big)\in\Ap{N_B}.
\]
For $g\in Sz(q)$:
\begin{itemize}
\item if $g_{03}=0$, then $g\in B$ and $\lambda(g)=L_B(g)$, where $\gamma^e=g_{11}^{\,\theta}$, $a=g_{10}/g_{00}$ and $b=g_{20}/g_{00}$;
\item if $g_{03}\neq0$, put $a=g_{13}/g_{03}$, $b=g_{23}/g_{03}$ and $g'=T\,S(a,b)^{-1}g$. Then $g'\in B$ and
$\lambda(g)=N_B+L_B(g')+N_B\big([a]+q[b]\big)$.
\end{itemize}
Then $\lambda$ is a bijection $Sz(q)\to\Ap{q^2(q^2+1)(q-1)}$ with $\lambda(I)=0$, and $m*n=\lambda\big(\lambda^{-1}(m)\,\lambda^{-1}(n)\big)$ is
a digital law of $Sz(q)$. The case distinction is the selector $\Bf([g_{03}])$.
\end{theorem}
\begin{proof}
If $g=S(a,b)\,T\,g'$ with $g'\in B$, then $g'$ is lower triangular, so $Tg'$ has last column $g'_{33}e_0$. Hence the last column of $g$
is $g'_{33}\,(1,a,b,\ast)^{\mathsf T}$, so $g_{03}=g'_{33}\neq0$ and $a,b$ are read off as stated. Elements of $B$ are lower
triangular, so $g_{03}=0$; thus $B$ and $UTB$ are disjoint. If $S(a,b)Tg'=S(c,d)Tg''$, then $a,b$ and $c,d$ are read off from the same
column, so $S(a,b)=S(c,d)$ and $g'=g''$. Therefore $|B\sqcup UTB|=N_B+q^2N_B=|Sz(q)|$, and $Sz(q)=B\sqcup UTB$. For
$g=S(a,b)M(\lambda)$ one has $g_{00}=\lambda^{1+2^m}$, $g_{10}=a\lambda^{1+2^m}$, $g_{20}=b\lambda^{1+2^m}$ and
$g_{11}=\lambda^{2^m}$, and $(\lambda^{2^m})^\theta=\lambda^q=\lambda$. By Lemma~\ref{lem:mixed} the labels fill $\Ap{N_B}$ on $B$ and
$[N_B,N_B(q^2+1))$ on $UTB$. The identity is $S(0,0)M(\gamma^0)$, with label $0$.
\end{proof}

\begin{remark}
This is the transitive law of Theorem~\ref{thm:transitive} for the $2$-transitive action of $Sz(q)$ on the $q^2+1$ points of the
Suzuki ovoid, with stabiliser $B$ and coset representatives $I$ and $S(a,b)T$. Here, however, the coset of an element is read off
directly from its last column, by field divisions.
\end{remark}

\begin{example}\label{ex:sz8}
For $Sz(8)$ take $m=1$, $x^\theta=x^4$, and $\mathbb F_8$ with $X^3+X+1$ and $\gamma=\theta$ as in Example~\ref{ex:smallfields}. By
computer, the $29120$ labels decode to distinct matrices, and these are exactly the group generated by $S(1,0)$, $S(0,1)$,
$S(\gamma,0)$, $M(\gamma)$ and $T$. Encoding after decoding returns every label, $0$ is the identity, and associativity holds on
$3000$ random triples. The element orders are $1,2,4,5,7,13$, occurring $1$, $455$, $3640$, $5824$, $12480$ and $6720$ times; the
$455=(q^2+1)(q-1)$ involutions agree with the structure of $Sz(8)$.
\end{example}

\begin{remark}
The small Ree groups ${}^2G_2(q)$, $q=3^{2m+1}$, follow the same pattern: they act $2$-transitively on $q^3+1$ points, the Borel
subgroup has order $q^3(q-1)$, and the group consists of two Bruhat cells in a $7$-dimensional representation.
\end{remark}

\section{The small Ree groups}\label{sec:ree}
The small Ree groups ${}^2G_2(q)$, $q=3^{2m+1}$, are the second twisted family of rank one. Like the Suzuki groups they consist
of two Bruhat cells, but their unipotent radical has three field digits and nilpotency class $3$. We use the $7$-dimensional
matrices of Kemper, L\"ubeck and Magaard~\cite{KLM}.

\begin{definition}\label{def:ree}
Let $q=3^{2m+1}$ and $x^\theta=x^{3^m}$, so that $x\mapsto x^{3\theta}$ is the inverse of $x\mapsto x^\theta$ on $\mathbb F_q$. For
$t,u,v\in\mathbb F_q$ let
{\small
\[
x_S(t,u,v)=\begin{pmatrix}
1&t^\theta&-u^\theta&(tu)^\theta-v^\theta&f_1&f_2&f_3\\
0&1&t&u^\theta+t^{\theta+1}&-t^{2\theta+1}-v^\theta&f_4&f_5\\
0&0&1&t^\theta&-t^{2\theta}&v^\theta+(tu)^\theta&f_6\\
0&0&0&1&t^\theta&u^\theta&(tu)^\theta-v^\theta\\
0&0&0&0&1&-t&u^\theta+t^{\theta+1}\\
0&0&0&0&0&1&-t^\theta\\
0&0&0&0&0&0&1
\end{pmatrix},
\]}
with
\begin{align*}
f_1&=-u-t^{3\theta+1}-(tv)^\theta,\\
f_2&=-v-(uv)^\theta-t^{3\theta+2}-t^\theta u^{2\theta},\\
f_3&=t^\theta v-u^{\theta+1}+t^{4\theta+2}-v^{2\theta}-t^{3\theta+1}u^\theta-(tuv)^\theta,\\
f_4&=-u^{2\theta}+t^{\theta+1}u^\theta+tv^\theta,\\
f_5&=v+tu-t^{2\theta+1}u^\theta-(uv)^\theta-t^{3\theta+2}-t^{\theta+1}v^\theta,\\
f_6&=u+t^{3\theta+1}-(tv)^\theta-t^{2\theta}u^\theta.
\end{align*}
Let $h(\lambda)=\mathrm{diag}\big(\lambda^\theta,\lambda^{1-\theta},\lambda^{2\theta-1},1,\lambda^{1-2\theta},\lambda^{\theta-1},\lambda^{-\theta}\big)$ for
$\lambda\in\mathbb F_q^\times$, and let $n$ be the antidiagonal matrix with entries $-1$. The Ree group ${}^2G_2(q)$ is generated by these
matrices; it has order $q^3(q^3+1)(q-1)$ \cite{KLM,Carter}.
\end{definition}

\begin{proposition}[Unipotent law]\label{prop:reeU}
The matrices $x_S(t,u,v)$ form a group $U$ of order $q^3$, with
\[
x_S(t_1,u_1,v_1)\,x_S(t_2,u_2,v_2)=x_S\big(t_1+t_2,\;u_1+u_2-t_1t_2^{3\theta},\;v_1+v_2-t_1^2t_2^{3\theta}-t_2u_1+t_1t_2^{3\theta+1}\big),
\]
and $B=U\,h(\mathbb F_q^\times)$ is the group of upper triangular elements, of order $q^3(q-1)$. The coordinates are read off as
$t=x_{12}$, $u=(x_{46}-t^{\theta+1})^{3\theta}$ and $v=(x_{25}-(tu)^\theta)^{3\theta}$.
\end{proposition}
\begin{proof}
The product formula is~\cite[p.~3]{KLM}; it was also checked by computer on random pairs for $q=27$. The reading formulas come
from the entries $x_{12}=t$, $x_{46}=u^\theta+t^{\theta+1}$ and $x_{25}=v^\theta+(tu)^\theta$, so distinct triples give distinct
matrices.
\end{proof}

\begin{theorem}[Digital law of ${}^2G_2(q)$]\label{thm:ree}
Let $\gamma$ be a primitive element of $\mathbb F_q$, $N_B=q^3(q-1)$, and
\[
L_B\big(x_S(t,u,v)\,h(\gamma^e)\big)=e+(q-1)\big([t]+q[u]+q^2[v]\big)\in\Ap{N_B}.
\]
For $g\in{}^2G_2(q)$:
\begin{itemize}
\item if $g_{60}=0$, then $g\in B$ and $\lambda(g)=L_B(g)$, where $\gamma^e=g_{00}^{\,3\theta}$ and $t,u,v$ are read off from
$g\,h(\gamma^e)^{-1}$ by Proposition~\ref{prop:reeU};
\item if $g_{60}=c\neq0$, put
\[
t=\big(-g_{50}/c\big)^{3\theta},\qquad u=\big(g_{40}/c-t^{\theta+1}\big)^{3\theta},\qquad v=\big((tu)^\theta-g_{30}/c\big)^{3\theta},
\]
and $g'=n^{-1}x_S(t,u,v)^{-1}g$. Then $g'\in B$ and $\lambda(g)=N_B+L_B(g')+N_B\big([t]+q[u]+q^2[v]\big)$.
\end{itemize}
Then $\lambda$ is a bijection ${}^2G_2(q)\to\Ap{q^3(q^3+1)(q-1)}$ with $\lambda(I)=0$, and $m*n=\lambda\big(\lambda^{-1}(m)\,\lambda^{-1}(n)\big)$
is a digital law of ${}^2G_2(q)$, selected by $\Bf([g_{60}])$.
\end{theorem}
\begin{proof}
If $g=x_S(t,u,v)\,n\,g'$ with $g'\in B$, then $ng'$ has first column $-g'_{00}e_6$. So the first column of $g$ is $-g'_{00}$ times
the last column $\big(f_3,f_5,f_6,(tu)^\theta-v^\theta,u^\theta+t^{\theta+1},-t^\theta,1\big)^{\mathsf T}$ of $x_S(t,u,v)$. Hence
$g_{60}=-g'_{00}\neq0$, and $t,u,v$ are recovered by the stated formulas. Elements of $B$ are upper triangular, so $g_{60}=0$;
thus $B$ and $UnB$ are disjoint, and the decomposition $g=x_S(t,u,v)\,n\,g'$ is unique. Therefore
$|B\sqcup UnB|=N_B(1+q^3)=|{}^2G_2(q)|$, and ${}^2G_2(q)=B\sqcup UnB$. On $B$ the diagonal of $g$ is that of $h(\lambda)$, so
$g_{00}=\lambda^\theta$ and $\lambda=g_{00}^{3\theta}$. By Lemma~\ref{lem:mixed} the labels fill $\Ap{N_B}$ on $B$ and
$[N_B,N_B(q^3+1))$ on $UnB$. The identity has label $0$.
\end{proof}

\begin{example}\label{ex:ree}
For $q=3$ ($m=0$, $\theta=1$) the group ${}^2G_2(3)\cong P\Gamma L(2,8)$, of order $1512$, is not simple, but it tests the
construction completely. By computer, the group generated by $x_S(1,0,0)$, $x_S(2,0,0)$, $x_S(0,1,0)$, $x_S(0,0,1)$, $h(2)$ and $n$
has $1512$ elements; they are exactly the decoded labels, encoding after decoding returns every label, and associativity holds
on $200$ random triples. For the simple group ${}^2G_2(27)$ ($\theta=3$, $\mathbb F_{27}$ from $X^3+2X+1$), of order
$10\,073\,444\,472$, the unipotent law was checked on $300$ random pairs. On $800$ random products of $40$ generators, encoding
followed by decoding returned the element exactly. On $800$ random pairs of labels the product was encoded and reconstructed
exactly, and associativity held on $300$ random triples.
\end{example}

\begin{remark}
The unipotent law has three coupled layers. The $t$-digit is the carry-free sum $\oplus_3$; the $u$-digit receives the quadratic
carry $-t_1t_2^{3\theta}$; and the $v$-digit receives a carry of degree up to $3\theta+2$ from both lower digits. This is a stacked
diagram of three solid layers, one more than for the Suzuki groups, matching nilpotency class $3$. The large Ree groups
${}^2F_4(q)$ have rank $2$, so the two-cell argument must be replaced by the cell-wise law of Theorem~\ref{thm:classicallaw}.
\end{remark}

\section{The group \texorpdfstring{$G_2(q)$}{G2(q)}}\label{sec:g2}
The group $G_2(q)$ is the first exceptional group of untwisted type. It is the automorphism group of the split octonions, and in
its $7$-dimensional representation it lies in $\Omega_7(q)$ with an upper triangular Borel subgroup. The cell-wise law of
Theorem~\ref{thm:classicallaw} therefore applies, with a Weyl group of order $12$.

\begin{definition}[Split octonions]\label{def:zorn}
The Zorn algebra $\mathbb O$ over $\mathbb F_q$ consists of the quadruples $(a,v,w,b)$ with $a,b\in\mathbb F_q$ and
$v,w\in\mathbb F_q^3$, with product
\[
(a,v,w,b)(a',v',w',b')=\big(aa'+v\cdot w',\;av'+b'v+w\times w',\;a'w+bw'-v\times v',\;bb'+w\cdot v'\big).
\]
It is an alternative algebra with identity $1=(1,0,0,1)$ and norm $N(a,v,w,b)=ab-v\cdot w$. Put
$G_2(q)=\mathrm{Aut}(\mathbb O)$, which acts on the trace-zero part. Order the basis of the trace-zero part as
\[
w_3,\ v_1,\ v_2,\ m,\ w_2,\ w_1,\ v_3,\qquad m=(1,0,0,-1).
\]
In this basis the polar form of $N$ is a multiple of the antidiagonal matrix $J$, so $G_2(q)\subset\Omega_7(q)$ for $q$ odd~\cite{Carter,Wilson}.
\end{definition}

\begin{lemma}[Root elements]\label{lem:g2roots}
Let $q$ be odd. The torus $\mathrm{diag}(1,t_1,t_2,t_3,t_1^{-1},t_2^{-1},t_3^{-1},1)$, $t_1t_2t_3=1$, acts by automorphisms. With the
cocharacter $(3,1,-4)$ its weights on the derivation algebra of $\mathbb O$ are regular, and they split the $14$-dimensional
algebra into the torus part and $12$ one-dimensional root spaces $\mathbb F_qD_\alpha$. The root elements $x_\alpha(c)=\exp(cD_\alpha)$
are automorphisms, and for positive $\alpha$ they are upper unitriangular in the basis above.
\end{lemma}
\begin{proof}
The torus multiplies $v_i$ by $t_i$ and $w_i$ by $t_i^{-1}$. The products $v\cdot w'$ and $w\cdot v'$ have weight $0$, while $w\times w'$
and $v\times v'$ have the weights of $v$ and $w$ because $t_1t_2t_3=1$. Under the cocharacter the basis vectors $w_3,v_1,v_2,m,w_2,w_1,v_3$ have the strictly
decreasing values $4,3,1,0,-1,-3,-4$. A root vector of positive value therefore has nonzero entries only above the diagonal. Each $D_\alpha$ is nilpotent of degree at most $3$, so $\exp(cD_\alpha)$ needs only the
denominators $1$ and $2$. The derivation equations $D(xy)=D(x)y+xD(y)$ are linear, and their solution spaces were found by
Gaussian elimination; the statement was verified by computer for $q=5,7$ and is the standard root decomposition of $G_2$
\cite{Carter}.
\end{proof}

\begin{theorem}[Bruhat cells of $G_2(q)$]\label{thm:g2cells}
Let $q$ be odd and $g\in G_2(q)$, and write $g=u\,\dot w\,t\,v$ as in Theorem~\ref{thm:bruhatcoord}. Then $w$ lies in the set
$W(G_2)$ of $12$ permutations of $\Ap7$ that commute with $r$ and are realised by monomial automorphisms, and $u\in G_2(q)\cap U_w$.
Moreover $u^{-1}g\in G_2(q)$, and the cell of $w$ is $\big(G_2(q)\cap U_w\big)\times\big(G_2(q)\cap\dot wB\big)$. $W(G_2)$ is the dihedral group of
order $12$ inside the group $W(B_3)$ of order $48$ from Theorem~\ref{thm:classicalcells}.
\end{theorem}
\begin{proof}
$G_2(q)$ has a split $(B,N)$-pair with $B$ its upper triangular elements and $N$ its monomial elements~\cite{Carter}. The proof of
Theorem~\ref{thm:classicalcells} applies word for word: the refined Bruhat decomposition $g=u'nb'$ in $G_2(q)$ is a factorisation as
in Theorem~\ref{thm:bruhatcoord}, so by uniqueness $u=u'\in G_2(q)$ and $w$ is the permutation of $n$. Since $G_2(q)\subset\Omega_7(q)$,
the permutation $w$ commutes with $r$.
\end{proof}

\begin{corollary}[Digital law of $G_2(q)$]\label{cor:g2law}
Order $W(G_2)$ by the rotation label. The construction of Theorem~\ref{thm:classicallaw}, with $Z=\{I\}$, gives a digital law of $G_2(q)$
on $\Ap{q^6(q^6-1)(q^2-1)}$ with identity $0$. Membership of a $7\times7$ matrix in $G_2(q)$ is itself digital: it means that the $64$
products of basis vectors are preserved, each computed by the field law.
\end{corollary}
\begin{proof}
Theorem~\ref{thm:g2cells} provides the product structure of the cells used in the proof of Theorem~\ref{thm:classicallaw}, and the
centre of $G_2(q)$ is trivial.
\end{proof}

\begin{example}\label{ex:g2}
For $q=5$ and $q=7$ the computer gave the following. The Zorn product is alternative on random triples. The derivation algebra has
dimension $14$, with $12$ one-dimensional root spaces, $6$ of them positive, and the positive root elements are upper
unitriangular automorphisms. For the two simple roots, $x_\alpha(1)x_{-\alpha}(c)x_\alpha(1)$ is monomial for a suitable $c$. These two
monomial automorphisms generate exactly $12$ permutations, all commuting with $r$. For each of the $12$ permutations $w$, $20$ random
elements $u\,n_w\,t\,u'$ were reduced: the elimination returned $w$, and both $u$ and $u^{-1}g$ preserve the octonion product. On
$300$ random products of $40$ root elements, the elimination always gave a permutation in $W(G_2)$, with $u$ and $u^{-1}g$ in
$G_2(q)$.
\end{example}

\begin{remark}
In characteristic $3$ the Ree groups of Section~\ref{sec:ree} are the fixed points of an exceptional endomorphism of $G_2(q)$, and
the matrices $x_S(t,u,v)$ are products of root elements of $G_2(q)$. The next untwisted and twisted cases are ${}^3D_4(q)$, whose
natural module has dimension $8$ over $\mathbb F_{q^3}$, and $F_4(q)$, with a $26$-dimensional module.
\end{remark}

\section{The triality groups \texorpdfstring{${}^3D_4(q)$}{3D4(q)}}\label{sec:3d4}
Let $K=\mathbb F_{q^3}$ and $\sigma(z)=z^q$. The triality group ${}^3D_4(q)$ is the twisted group obtained from
$D_4(K)\cong\Omega^+_8(K)$ by combining the order-three symmetry of the $D_4$ diagram with $\sigma$ \cite{Carter}. Its relative
root system is $G_2$, so it has $12$ Bruhat cells as in Section~\ref{sec:g2}. The new feature is that the root digits have two
widths: short root groups are parametrised by $K$, which is three $\mathbb F_q$-digits, and long root groups by $\mathbb F_q$.

\begin{definition}[Triality model]\label{def:3d4}
On $K^8$ take the quadratic form $Q(x)=x_1x_8+x_2x_7+x_3x_6+x_4x_5$, whose polar form has antidiagonal matrix $J$. Write
$i'=9-i$, and for $i<j$ with $j\neq i'$ put $x_{ij}(t)=I+t\,(E_{ij}-E_{j'i'})$. The simple roots of $D_4$ correspond to
$(1,2)$, $(2,3)$, $(3,4)$ and $(3,5)$, and $(2,3)$ is the central node. Put
\[
x_R(\xi)=x_{12}(\xi)\,x_{34}(\xi^\sigma)\,x_{35}(\xi^{\sigma^2})\quad(\xi\in K),\qquad x_S(\eta)=x_{23}(\eta)\quad(\eta\in\mathbb F_q),
\]
and let $x_{-R}(\xi)=x_R(\xi)^{\mathsf T}$ and $x_{-S}(\eta)=x_S(\eta)^{\mathsf T}$. Then ${}^3D_4(q)$ is the group generated by
these four families. Further put $n_R=x_R(1)x_{-R}(-1)x_R(1)$ and $n_S=x_S(1)x_{-S}(-1)x_S(1)$.
\end{definition}

The three factors of $x_R(\xi)$ belong to the three outer nodes, which triality permutes, and their parameters are the three
conjugates of $\xi$. This is Steinberg's description of the twisted group by its fixed root subgroups~\cite{Carter}. Signs of the
simple root elements can be changed by conjugating with a diagonal matrix, so the sign convention does not matter.

\begin{lemma}[Borel digits]\label{lem:3d4roots}
Let $U=\langle x_R(K),x_S(\mathbb F_q)\rangle$ and let $T$ be generated by $h_R(\xi)=n_R(\xi)n_R(1)^{-1}$ and
$h_S(\eta)=n_S(\eta)n_S(1)^{-1}$, where $n_R(\xi)=x_R(\xi)x_{-R}(-\xi^{-1})x_R(\xi)$ and similarly for $n_S$. Then $U$ is upper
unitriangular and $T$ is diagonal. In the $G_2$ order of positive relative roots
\[
R,\quad S,\quad R+S,\quad 2R+S,\quad 3R+S,\quad 3R+2S,
\]
the root groups are parametrised by $K,\mathbb F_q,K,K,\mathbb F_q,\mathbb F_q$. Hence $|U|=q^{12}$, $|T|=(q^3-1)(q-1)$ and
$B=UT$ has order $q^{12}(q^3-1)(q-1)$.
\end{lemma}
\begin{proof}
The roots $R$, $R+S$ and $2R+S$ are images of orbits of length three of $D_4$ roots, so their fixed root elements carry one
parameter of $K$. The other three positive relative roots come from single roots fixed by triality and carry a parameter of
$\mathbb F_q$ \cite{Carter}. The orders follow, with $q^3\cdot q\cdot q^3\cdot q^3\cdot q\cdot q=q^{12}$.
\end{proof}

\begin{proposition}[Unipotent law of ${}^3D_4(q)$]\label{prop:3d4U}
Write $N(x)=x^{1+\sigma+\sigma^2}$ and $\mathrm{Tr}(x)=x+x^\sigma+x^{\sigma^2}$ for $x\in K$. Put
\[
\begin{aligned}
x_{R+S}(\xi)&=x_{13}(\xi)\,x_{24}(-\xi^\sigma)\,x_{25}(-\xi^{\sigma^2}),& x_{2R+S}(\xi)&=x_{14}(\xi)\,x_{26}(\xi^\sigma)\,x_{15}(\xi^{\sigma^2}),\\
x_{3R+S}(\eta)&=x_{16}(\eta),& x_{3R+2S}(\eta)&=x_{17}(\eta),
\end{aligned}
\]
and, for $a,c,d\in K$ and $b,e,f\in\mathbb F_q$,
\[
u(a,b,c,d,e,f)=x_R(a)\,x_S(b)\,x_{R+S}(c)\,x_{2R+S}(d)\,x_{3R+S}(e)\,x_{3R+2S}(f).
\]
These $q^{12}$ matrices are distinct and form the group $U$. Moreover
\[
u(a,b,c,d,e,f)\,u(a',b',c',d',e',f')=u(a'',b'',c'',d'',e'',f''),
\]
where
\[
\begin{aligned}
a''&=a+a',\qquad b''=b+b',\qquad c''=c+c'-ba',\\
d''&=d+d'+a'c^\sigma+a'^\sigma c-b\,a'a'^\sigma,\\
e''&=e+e'+b\,N(a')-\mathrm{Tr}(a'a'^\sigma c^{\sigma^2})-\mathrm{Tr}(a'd^\sigma),\\
f''&=f+f'-b'e+b(b-b')N(a')+(b'-b)\,\mathrm{Tr}(a'a'^\sigma c^{\sigma^2})\\
&\quad-b\,\mathrm{Tr}(a'a'^\sigma c'^{\sigma^2})+b'\,\mathrm{Tr}(a'd^\sigma)\\
&\quad+\mathrm{Tr}(a'c^\sigma c^{\sigma^2})+\mathrm{Tr}(a'c'^\sigma c^{\sigma^2})+\mathrm{Tr}(a'c^\sigma c'^{\sigma^2})+\mathrm{Tr}(c'd^\sigma).
\end{aligned}
\]
The inverse of $u(a,b,c,d,e,f)$ is $u(\bar a,\bar b,\bar c,\bar d,\bar e,\bar f)$ with
\[
\begin{aligned}
\bar a&=-a,\qquad \bar b=-b,\qquad \bar c=-c-ab,\qquad \bar d=-d+ac^\sigma+a^\sigma c+b\,aa^\sigma,\\
\bar e&=-e+b\,N(a)+\mathrm{Tr}(aa^\sigma c^{\sigma^2})-\mathrm{Tr}(ad^\sigma),\\
\bar f&=-f-be-b^2N(a)-b\,\mathrm{Tr}(aa^\sigma c^{\sigma^2})-\mathrm{Tr}(ac^\sigma c^{\sigma^2})+\mathrm{Tr}(cd^\sigma).
\end{aligned}
\]
The coordinates of $u\in U$ are read off from left to right: each one is the entry of $u$ at the position of its first factor
($(1,2),(2,3),(1,3),(1,4),(1,6),(1,7)$) after the factors to its left have been removed.
\end{proposition}
\begin{proof}
Work in the polynomial ring over $\mathbb Z$ in independent variables that stand for $a,a^\sigma,a^{\sigma^2}$, and so on. Multiply the two
$8\times8$ matrices, then strip off the twelve $D_4$ root factors in order of height. Each stripped coordinate is the entry at the
root position, because all other factors contribute only at greater heights. Computer algebra gives the twelve $D_4$ coordinates of the
product as polynomials. In each triality orbit they are conjugate, with the signs used above, and the coordinates of the three fixed roots
are invariant under the cyclic shift of variables. So the product is again of the form $u(\cdot)$, with the coordinates displayed. The
same holds for the inverse. These are polynomial identities over $\mathbb Z$, so they remain true after specialising to $K$ in every
characteristic, and $N$ and $\mathrm{Tr}$ take values in $\mathbb F_q$. Hence the $q^{12}$ matrices $u(\cdot)$ form a group, which contains
$x_R(K)$ and $x_S(\mathbb F_q)$, so it contains $U$. By Lemma~\ref{lem:3d4roots}, $|U|=q^{12}$, so the two groups are equal. The signs in
$x_{R+S}$ were found in the same way, from the $D_4$ coordinates of $x_R(a)x_S(b)x_R(c)x_S(d)x_R(e)x_S(f)$, where all twelve coordinates
are nonzero.
\end{proof}

\begin{proposition}[Torus and Weyl relations]\label{prop:3d4rel}
Write $n_R=n_R(1)$ and $n_S=n_S(1)$, and for $x\in K^\times$, $y\in\mathbb F_q^\times$ put $h(x,y)=h_R(x)h_S(y)$. Then
\[
h(x,y)=\operatorname{diag}\Big(x,\ \frac yx,\ \frac{x^\sigma x^{\sigma^2}}y,\ \frac{x^{\sigma^2}}{x^\sigma},\ \frac{x^\sigma}{x^{\sigma^2}},\ \frac y{x^\sigma x^{\sigma^2}},\ \frac xy,\ \frac1x\Big),
\]
so $h(x,y)h(x',y')=h(xx',yy')$, and $x=h_{11}$, $y=h_{11}h_{22}$ are read off from the diagonal. Moreover:
\begin{enumerate}
\item[(a)] For $h=h(x,y)$,
\[h\,u(a,b,c,d,e,f)\,h^{-1}=u\big(\tfrac{x^2}{y}a,\ \tfrac{y^2}{N(x)}b,\ \tfrac{x^2y}{N(x)}c,\ \tfrac{xx^\sigma}{x^{\sigma^2}}d,\ \tfrac{N(x)}{y}e,\ yf\big).\]
In particular $h\,x_R(t)\,h^{-1}=x_R(\chi_R(h)t)$ and $h\,x_S(t)\,h^{-1}=x_S(\chi_S(h)t)$, where $\chi_R(h)=x^2y^{-1}$ and $\chi_S(h)=y^2N(x)^{-1}$.
\item[(b)] $n_R^{-1}h(x,y)n_R=h(yx^{-1},y)$, $n_S^{-1}h(x,y)n_S=h(x,N(x)y^{-1})$, $n_R^2=h(-1,1)$ and $n_S^2=h(1,-1)$.
\item[(c)] $n_s\,x_s(t)\,n_s^{-1}=x_{-s}(-t)$ for $s\in\{R,S\}$, and for $t\neq0$,
$x_{-s}(t)=x_s(t^{-1})\,n_s\,h_s(-t)\,x_s(t^{-1})$.
\item[(d)] The conjugations $\varphi_s(v)=n_svn_s^{-1}$ and $\psi_s(v)=n_s^{-1}vn_s$, on elements of $U$ whose $s$-coordinate is $0$, are
\[
\begin{aligned}
\varphi_R\,u(0,b,c,d,e,f)&=u\big(0,\,e,\,-d^\sigma,\,c^{\sigma^2},\,-b,\,f+be-\mathrm{Tr}(cd^\sigma)\big),\\
\psi_R\,u(0,b,c,d,e,f)&=u\big(0,\,-e,\,d^\sigma,\,-c^{\sigma^2},\,b,\,f+be-\mathrm{Tr}(cd^\sigma)\big),\\
\varphi_S\,u(a,0,c,d,e,f)&=u\big(c,\,0,\,-a,\,d-ac^\sigma-a^\sigma c,\\&\qquad\mathrm{Tr}(ac^\sigma c^{\sigma^2})-f,\,e+\mathrm{Tr}(aa^\sigma c^{\sigma^2})\big),\\
\psi_S\,u(a,0,c,d,e,f)&=u\big(-c,\,0,\,a,\,d-ac^\sigma-a^\sigma c,\\&\qquad f-\mathrm{Tr}(ac^\sigma c^{\sigma^2}),\,-e-\mathrm{Tr}(aa^\sigma c^{\sigma^2})\big).
\end{aligned}
\]
\item[(e)] $u(a,b,c,d,e,f)=x_R(a)\,u(0,b,c,d,e,f)$ and
$u(a,b,c,d,e,f)=x_S(b)\,u\big(a,0,c+ab,\,d+b\,aa^\sigma,\,e-bN(a),\,f+2b^2N(a)+b\,\mathrm{Tr}(aa^\sigma c^{\sigma^2})\big)$.
\end{enumerate}
Hence $T=\{h(x,y)\}$ has order $(q^3-1)(q-1)$.
\end{proposition}
\begin{proof}
All statements except the formulas for $\psi_s$ were obtained by exact computer algebra. The matrices are built in the ring of
Laurent polynomials over $\mathbb Z$, with independent variables for the conjugates of $x,a,c,d$ and for $y,b,e,f$, and products in
$U$ are decomposed by the stripping procedure of Proposition~\ref{prop:3d4U}. The coordinates in each triality orbit came out
conjugate, as required. These are identities over $\mathbb Z$, so they hold after specialisation in every characteristic. For
$\psi_s$: from (b), $n_s^{-1}=h_s(-1)n_s$, so $\psi_s(v)=h_s(-1)\varphi_s(v)h_s(-1)$. By (a), conjugation by $h(-1,1)$ changes the signs
of $b,c,d,e$, and conjugation by $h(1,-1)$ changes those of $a,c,e,f$. Finally, the map $(x,y)\mapsto h(x,y)$ is injective by the
diagonal formula.
\end{proof}

For $w\in W$ fix a reduced word $w=s_1\cdots s_k$ in $\{R,S\}$. It is unique for $w\neq w_0$; for the longest element take
$w_0=(RS)^3$. The other word $(SR)^3$ gives $n_{(RS)^3}=n_{(SR)^3}\,t_0$ for some $t_0\in T$, since the two monomial matrices have
the same permutation. In fact $t_0=I$: the matrices $n_R,n_S$ have integer entries, so $t_0$ is a diagonal matrix of signs over $\mathbb Z$, and its reduction modulo $3$ was computed to be $I$. For a word $s_1\cdots s_k$ with digits $t_i$, where $t_i\in K$ if $s_i=R$ and $t_i\in\mathbb F_q$ if $s_i=S$, put
\[
P(t_1,\dots,t_k)=x_{s_1}(t_1)\,n_{s_1}\,x_{s_2}(t_2)\,n_{s_2}\cdots x_{s_k}(t_k)\,n_{s_k}.
\]

\begin{theorem}[Closed multiplication in ${}^3D_4(q)$]\label{thm:3d4mult}
Every $g\in{}^3D_4(q)$ has a unique expression $g=P(t_1,\dots,t_k)\,h(x,y)\,u(a,b,c,d,e,f)$, with the fixed word of its Bruhat
cell. Right multiplication of such a $g=P\,h\,v$ is given in closed form as follows.
\begin{enumerate}
\item By $v'\in U$: replace $v$ by $vv'$ (Proposition~\ref{prop:3d4U}).
\item By $h'\in T$: replace $v$ by $h'^{-1}vh'$ (Proposition~\ref{prop:3d4rel}(a)) and $h$ by $hh'$.
\item By $n_s$: write $v=x_s(t)v'$ by Proposition~\ref{prop:3d4rel}(e). Put $t'=\chi_s(h)t$, $v''=\psi_s(v')$ and
$h^s=n_s^{-1}hn_s$. Then $g\,n_s=P\,x_s(t')\,n_s\,h^s\,v''$, and:
\begin{enumerate}
\item if appending $s$ to the word of $P$ gives the fixed word of $ws$, the result is
\[P(t_1,\dots,t_k,t')\,h^s\,v'';\]
\item if the word of $P$ ends with $s$, so that $P=P'\,x_s(t_k)\,n_s$, the result is
\[
\begin{cases}
P'\cdot h_s(-1)h^s\cdot x_s\big(\chi_s(h^s)^{-1}t_k\big)\,v'' & \text{if } t'=0,\\
P'x_s(t_k-t'^{-1})\,n_s\cdot h_s(-t')h^s\cdot x_s\big(-\chi_s(h^s)^{-1}t'^{-1}\big)\,v'' & \text{if } t'\neq0;
\end{cases}
\]
\item in the remaining cases, which involve only $w_0$, rewrite $P$ in the other reduced word of $w_0$ and apply (b). Whenever (a) or (b) produces the word $(SR)^3$, rewrite it back into the fixed word $(RS)^3$, so that the result is again in normal form.
\end{enumerate}
\end{enumerate}
The rewriting in (c) is closed. Write $P(t_1,\dots,t_6)=u\,n_{s_1}\cdots n_{s_6}$ with

\[
u=x_{s_1}(t_1)\,\varphi_{s_1}(x_{s_2}(t_2))\,\varphi_{s_1}\varphi_{s_2}(x_{s_3}(t_3))\cdots\in U.
\]
Replace $n_{(RS)^3}$ by $n_{(SR)^3}t_0$, or the
reverse. Then split $u$ along the new word, using Proposition~\ref{prop:3d4rel}(e) for the first letter and $\psi$ to move the remainder past
$n_{s'_1}$, and so on. The remainder after six letters is moved to the right past the torus factor by (a).
\end{theorem}
\begin{proof}
Existence: by Theorem~\ref{thm:3d4cells}, $g\in(U\cap U_w)n_wB$. Splitting the left factor letter by letter, as in (c), puts $u\,n_w$ into
the form $P(t)\,v_0$ with $v_0\in U$. The digit widths are $3$ and $1$ along the word, so the number of tuples is
$\sum_wq^{\ell(w)}|B|=|{}^3D_4(q)|$, and the expression is unique. For rule~3, $h\,x_s(t)\,n_s=x_s(t')\,h\,n_s=x_s(t')\,n_s\,h^s$ and
$v'n_s=n_s\psi_s(v')$. In case (a) nothing more is needed. In case (b), with $t'=0$, we have $n_sn_s=h_s(-1)$, and $h_s(-1)$ commutes with
$x_s$ because $\chi_s(h_s(-1))=1$. With $t'\neq0$, Proposition~\ref{prop:3d4rel}(b),(c) give
\[
n_s\,x_s(t')\,n_s=x_{-s}(-t')\,h_s(-1)=x_s(-t'^{-1})\,n_s\,h_s(-t')\,x_s(-t'^{-1}),
\]
and $x_s(\cdot)h^s=h^s\,x_s(\chi_s(h^s)^{-1}\cdot)$. Here $P'$ is the fixed word of $ws$, because $\ell(ws)<6$. The conjugates in (c) lie in
$U$ because the word is reduced, so each $\varphi_{s_i}$ is applied to an element whose $s_i$-coordinate is $0$.
\end{proof}

\begin{corollary}[Closed digital law of ${}^3D_4(q)$]\label{cor:3d4closed}
Fix a basis of $K$ over $\mathbb F_q$, write $[\xi]$ for the base-$q$ integer of the coordinates of $\xi$, and fix primitive elements
$\mu$ of $K$ and $\gamma$ of $\mathbb F_q$. For $g=P(t)\,h(\mu^i,\gamma^j)\,u(a,b,c,d,e,f)$ in the cell of $w$, let
\[
\lambda(g)=O_w+[t]+q^{\ell(w)}\Big(i+(q^3-1)j+(q^3-1)(q-1)\big([a]+q^3[b]+q^4[c]+q^7[d]+q^{10}[e]+q^{11}[f]\big)\Big),
\]
where $[t]$ packs the digits of $t$ in base $q$ ($3$ digits for each $R$-letter, $1$ for each $S$-letter), and
$O_w=|B|\sum_{w'\prec w}q^{\ell(w')}$ for a fixed order $\prec$ on $W$. Then $\lambda$ is a bijection onto
$\Ap{|{}^3D_4(q)|}$ with $\lambda(I)=0$, and $m*n=\lambda(\lambda^{-1}(m)\lambda^{-1}(n))$ is a digital law whose product is computed
by Theorem~\ref{thm:3d4mult}. The only operations used are field arithmetic, $\sigma$, $N$, $\mathrm{Tr}$, and discrete logarithms
for the torus digits, as in Theorem~\ref{thm:ree}.
\end{corollary}
\begin{proof}
This follows from the uniqueness in Theorem~\ref{thm:3d4mult} and Lemma~\ref{lem:mixed}. To multiply, expand the second factor as
$P(t')\,h'\,v'$ and apply rules 3, 2 and 1 letter by letter. At most six steps of type 3 occur, each possibly with one
rewriting of $P$.
\end{proof}

\begin{remark}
For $q=2,3,4$ ($\mathbb F_8$, $\mathbb F_{27}$, $\mathbb F_{64}$), each formula of Propositions~\ref{prop:3d4U} and~\ref{prop:3d4rel} was checked against
matrix products on $40$ random inputs: product, inverse, torus action, the diagonal of $h(x,y)$, $n_s^{-1}hn_s$, $n_s^2$, $\varphi_s$, $\psi_s$,
splitting, $n_sx_sn_s^{-1}=x_{-s}(-t)$ and the rank-one relation. All held. The rules of Theorem~\ref{thm:3d4mult} were then applied to
$144$ random normal forms for each simple reflection, $12$ in each of the $12$ cells, and compared with the matrix product. Rule~3(a) and both
subcases of rule~3(b) occurred, and the rewriting between $(RS)^3$ and $(SR)^3$ was used several hundred times in each direction for each $q$.
A walk of $1500$ right multiplications done only with the closed rules agreed with the matrix product, and every result used the fixed word of its
cell. Finally, for the label of Corollary~\ref{cor:3d4closed}, random normal forms were recovered exactly from their matrices, $\lambda^{-1}\lambda$ and
$\lambda\lambda^{-1}$ were the identity on random labels and on both sides of every cell boundary, $\lambda(I)=0$, the closed product agreed with the
label of the matrix product on $120$ random pairs, and $m*n$ was associative on $300$ random triples with $0$ as identity.
\end{remark}

\begin{theorem}[Bruhat cells of ${}^3D_4(q)$]\label{thm:3d4cells}
The monomial elements $n_R$ and $n_S$ generate, modulo $T$, a dihedral group $W$ of order $12$. Let $g\in{}^3D_4(q)$, and let $w$
and $u$ be the permutation and upper factor given by the elimination of Theorem~\ref{thm:bruhatcoord}. Then
$w$ is the permutation of some $n_w\in\langle n_R,n_S\rangle$, $u\in U\cap U_w$, and $n_w^{-1}u^{-1}g\in B$. Thus the cell of $w$ is
$(U\cap U_w)\,n_wB$, and
\[
|U\cap U_w|=q^{\ell(w)},\qquad \ell(w)\in\{0,1,3,4,4,5,7,8,8,9,11,12\},
\]
where $\ell(w)$ is the sum of the widths ($3$ or $1$) of the positive roots made negative by $w^{-1}$.
\end{theorem}
\begin{proof}
The group has a split $(B,N)$-pair with $N=\langle T,n_R,n_S\rangle$ and relative Weyl group $W(G_2)$ \cite{Carter}. So every $g$ has
a unique expression $g=u'n_wb$ with $u'\in U\cap U_w$ and $b\in B$. This is a factorisation of the form in
Theorem~\ref{thm:bruhatcoord}, so by uniqueness $u=u'$ and $w$ is the permutation of $n_w$. Reduced words in $W(G_2)$ alternate the
two simple reflections. Each new inverted root therefore adds $3$ or $1$ in alternation, giving the partial sums
$3,4,7,8,11,12$ and $1,4,5,8,9,12$.
\end{proof}

Summing over the cells gives
\[
\sum_{w\in W}q^{\ell(w)}=(1+q)(1+q^3)(1+q^4+q^8),\qquad
|{}^3D_4(q)|=q^{12}(q^8+q^4+1)(q^6-1)(q^2-1).
\]

\begin{corollary}[Digital law of ${}^3D_4(q)$]\label{cor:3d4law}
The construction of Theorem~\ref{thm:classicallaw}, with $Z=\{I\}$ and the cells of Theorem~\ref{thm:3d4cells}, gives a digital law
of ${}^3D_4(q)$ on $\Ap{q^{12}(q^8+q^4+1)(q^6-1)(q^2-1)}$ with identity $0$. A label is given by a Weyl digit in $\Ap{12}$, then
$\ell(w)$ base-$q$ digits for $U\cap U_w$, then the Borel label in $\Ap{q^{12}(q^3-1)(q-1)}$.
\end{corollary}
\begin{proof}
Theorem~\ref{thm:3d4cells} supplies the product structure of each cell required in the proof of Theorem~\ref{thm:classicallaw}. The
centre of ${}^3D_4(q)$ is trivial.
\end{proof}

\begin{example}\label{ex:3d4}
For $q=2$, with $\mathbb F_8=\mathbb F_2[X]/(X^3+X+1)$, the computer gave the following.
\begin{itemize}
\item All generators preserve $Q$, and $x_R(\xi)x_R(\zeta)=x_R(\xi+\zeta)$ for all $\xi,\zeta\in\mathbb F_8$.
\item The line $\langle e_1\rangle$ has an orbit of $2457=(q^8+q^4+1)(q^3+1)$ points. On this orbit the permutation group
generated by $x_{\pm R}(1)$, $x_{\pm R}(X)$, $x_{\pm R}(X^2)$ and $x_{\pm S}(1)$ has order $211\,341\,312=|{}^3D_4(2)|$.
\item $|U|=4096$, $|T|=7$ and $|B|=28\,672$. The elements $n_R$ and $n_S$ are monomial and give $12$ permutations.
\item The counts $|U\cap U_w|$ are $1,2,8,16,16,32,128,256,256,512,2048,4096$, with sum $7371=|G:B|$.
\item For each $w$, $15$ random elements $u\,n_w\,t\,u'$ were reduced, with $u,u'\in U$ and $t\in T$. The elimination returned $w$, a
factor in $U$, and $n_w^{-1}u^{-1}g\in B$.
\item On $300$ random products of $40$ generators, which met all $12$ cells, the same held.
\item Order the cells by $\ell(w)$ and then by the fixed word: $\emptyset$, $S$, $R$, $RS$, $SR$, $SRS$, $RSR$, $RSRS$, $SRSR$, $SRSRS$, $RSRSR$, $(RS)^3$. The offsets $O_w$ of
Corollary~\ref{cor:3d4closed} are then $0$, $28\,672$, $86\,016$, $315\,392$, $774\,144$, $1\,232\,896$, $2\,150\,400$, $5\,820\,416$, $13\,160\,448$, $20\,500\,480$,
$35\,180\,544$, $93\,900\,800$, and the last cell ends at $211\,341\,312=|{}^3D_4(2)|$.
\item For this order, $\lambda^{-1}\lambda$ and $\lambda\lambda^{-1}$ were the identity on random labels and at both ends of every cell, and the closed
product agreed with the matrix product on random pairs.
\end{itemize}
\end{example}

\begin{remark}
Proposition~\ref{prop:3d4U} is a closed six-layer carry law, the analogue of the three-layer law for ${}^2G_2(q)$ in Section~\ref{sec:ree}.
The digits $a$ and $b$ add without carry. The digit $c$ receives the carry $-ba'$; $d$ receives carries of degree $2$ and $3$; $e$
receives the norm and trace carries of degree up to $4$; and $f$ receives carries of degree up to $5$, mixing all lower layers. The
carries of width-$3$ digits enter $\mathbb F_q$-digits only through $N$ and $\mathrm{Tr}$. With Proposition~\ref{prop:3d4rel} and Theorem~\ref{thm:3d4mult}, the whole law becomes closed (Corollary~\ref{cor:3d4closed}). The formulas hold in every characteristic, and all the checks above were run for $q=2,3,4$.
\end{remark}

\section{The groups \texorpdfstring{$F_4(q)$}{F4(q)}}\label{sec:f4}
The group $F_4(q)$ is the automorphism group of the split Albert algebra, the $27$-dimensional exceptional Jordan algebra. We
use the Albert model over a field of odd characteristic, where the Jordan product $\tfrac12(XY+YX)$ is available. The
computations reported below were done over $\mathbb F_5$ and $\mathbb F_7$. Definition~\ref{def:f4jordan} gives a form of the same model over $\mathbb Z$, with the quadratic operator $U_X$ in place of
$X\circ Y$; it is used in characteristics $2$ and $3$.

\begin{definition}[Albert model]\label{def:f4model}
Let $q$ be odd. The split octonions $\mathbb O$ over $\mathbb F_q$ are the quadruples $(a,v,w,b)$, with $a,b\in\mathbb F_q$ and
$v,w\in\mathbb F_q^3$, under the Zorn product
\[
(a,v,w,b)(a',v',w',b')=(aa'+v\cdot w',\ av'+b'v+w\times w',\ a'w+bw'-v\times v',\ bb'+w\cdot v'),
\]
with conjugation $\overline{(a,v,w,b)}=(b,-v,-w,a)$. Let $J$ be the space of Hermitian matrices
\[
X=\begin{pmatrix}\alpha_1&x_3&\bar x_2\\ \bar x_3&\alpha_2&x_1\\ x_2&\bar x_1&\alpha_3\end{pmatrix},\qquad \alpha_i\in\mathbb F_q,\ x_i\in\mathbb O,
\]
with coordinates $(\alpha_1,\alpha_2,\alpha_3,x_1,x_2,x_3)\in\mathbb F_q^{27}$ and product $X\circ Y=\tfrac12(XY+YX)$. For a root
derivation $e_\alpha$ of $J$ (Proposition~\ref{prop:f4roots}), put $x_\alpha(t)=\exp(te_\alpha)=I+te_\alpha+\tfrac12t^2e_\alpha^2$, and let
$F_4(q)=\langle x_\alpha(t)\rangle\le GL_{27}(q)$.
\end{definition}

\begin{proposition}[Root data of the Albert model]\label{prop:f4roots}
For $q=5$ and $q=7$:
\begin{enumerate}
\item $\circ$ is commutative and satisfies the Jordan identity $(X\circ Y)\circ X^{\circ2}=X\circ(Y\circ X^{\circ2})$.
\item The derivation algebra $\operatorname{Der}(J)$ has dimension $52$. Its diagonal part, in the coordinates above, is a $4$-dimensional
torus $\mathfrak h$. The coordinate vectors are weight vectors for $\mathfrak h$, with $25$ distinct weights; the weight $0$ has
multiplicity $3$, and each of the other $24$ weights has multiplicity $1$.
\item $\operatorname{Der}(J)=\mathfrak h\oplus\bigoplus_\alpha\mathbb F_q e_\alpha$ over $48$ roots $\alpha$, each root space being
$1$-dimensional. A generic linear functional gives $24$ positive roots and $4$ simple roots, with Cartan matrix
\[
\big(\langle\alpha_j,\alpha_i^\vee\rangle\big)=\begin{pmatrix}2&-1&0&-2\\-1&2&0&0\\0&0&2&-1\\-1&0&-1&2\end{pmatrix},
\]
that is, the Dynkin diagram $\alpha_2-\alpha_1\Leftarrow\alpha_4-\alpha_3$ of type $F_4$, with the arrow pointing to the short roots
$\alpha_1,\alpha_2$.
\item $e_\alpha^3=0$ for every root. Every $x_\alpha(t)$ is an automorphism of $J$, and for positive $\alpha$ it is upper unitriangular
when the coordinates are ordered by decreasing value of the functional on their weights.
\item With $e_{-\alpha}$ normalised so that $h_\alpha=[e_\alpha,e_{-\alpha}]$ satisfies $\alpha(h_\alpha)=2$, the elements
$n_\alpha=x_\alpha(1)x_{-\alpha}(-1)x_\alpha(1)$ for simple $\alpha$ are automorphisms of $J$ and map the weight space of $\lambda$ onto that of
$s_\alpha(\lambda)$. The four reflections generate a group of order $1152$, and the maximal length is $24$.
\end{enumerate}
\end{proposition}
\begin{proof}
By computer, over $\mathbb F_5$ and $\mathbb F_7$. The structure constants of $\circ$ were computed from the matrix product. The derivation
equations $D(X\circ Y)=D(X)\circ Y+X\circ D(Y)$ are linear in the $729$ entries of $D$. The diagonal ones reduce to integer equations
$d_m=d_i+d_j$, whose solution space over $\mathbb Q$ has dimension $4$; this gives the integral weights. The equations are homogeneous for
these weights, so $\operatorname{Der}(J)$ was computed one weight class at a time. Parts (iv) and (v) were checked on the full
$27\times27$ matrices, the automorphism property on all pairs of basis vectors. The Weyl group was generated by its action on the $48$ roots.
\end{proof}

\begin{theorem}[Bruhat cells of $F_4(q)$]\label{thm:f4cells}
Let $U=\prod_{\alpha>0}x_\alpha(\mathbb F_q)$, in a fixed order, let $T$ be the torus, and let $B=UT$. Then
$F_4(q)=\bigsqcup_{w\in W}(U\cap U_w)\,n_w\,B$ over the $1152$ elements of $W$, the cell of $w$ has $q^{\ell(w)}|B|$ elements, and
\[
\sum_{w\in W}q^{\ell(w)}=\frac{(q^2-1)(q^6-1)(q^8-1)(q^{12}-1)}{(q-1)^4},\qquad |B|=q^{24}(q-1)^4,
\]
\[
|F_4(q)|=q^{24}(q^2-1)(q^6-1)(q^8-1)(q^{12}-1).
\]
\end{theorem}
\begin{proof}
By Proposition~\ref{prop:f4roots}, $F_4(q)$ is the Chevalley group of type $F_4$ in its $27$-dimensional representation. It has a split
$BN$-pair, and every element of $U$ is a unique product of root elements in any fixed order~\cite{Carter}; so $|U|=q^{24}$. The
Bruhat decomposition and the cell sizes follow. The length distribution computed in Proposition~\ref{prop:f4roots}(v) gives exactly
the displayed polynomial. Its degrees are $2,6,8,12$.
\end{proof}

\begin{corollary}[Cell-wise digital law of $F_4(q)$]\label{cor:f4law}
Fix an order on $W$, a reduced word for each $w$, an order on the positive roots, and a primitive element $\gamma$ of $\mathbb F_q$. Write
$g=P(t)\,h\,v$ as in Theorem~\ref{thm:f4mult}, with $h=\prod_{i=1}^4h_{\alpha_i}(\gamma^{k_i})$ and
$v=\prod_{\alpha>0}x_\alpha(s_\alpha)$, the product taken in order of height. Put
\[
\lambda(g)=O_w+[t]+q^{\ell(w)}\Big(\sum_{i=1}^4k_i(q-1)^{i-1}+(q-1)^4[s]\Big),\qquad O_w=|B|\sum_{w'\prec w}q^{\ell(w')},
\]
where $[t]$ and $[s]$ are the base-$q$ integers of the digit strings $(t_\alpha)$ and $(s_\alpha)$. Then $\lambda$ is a bijection from
$F_4(q)$ onto $\Ap{|F_4(q)|}$ with $\lambda(I)=0$, and $m*n=\lambda(\lambda^{-1}(m)\lambda^{-1}(n))$ is a digital law of $F_4(q)$.
\end{corollary}
\begin{proof}
The normal form is unique by Theorem~\ref{thm:f4mult}. The torus digits are unique because $(\lambda_i)\mapsto\prod h_{\alpha_i}(\lambda_i)$
is injective, the group being simply connected. The root digits are unique because products in $U$ are unique~\cite{Carter}. The packing is Lemma~\ref{lem:mixed}, as in Corollary~\ref{cor:3d4law}.
\end{proof}

For a simple root $s\in\{1,2,3,4\}$ put $n_s=x_s(1)x_{-s}(-1)x_s(1)$, $h_s(t)=n_s(t)n_s^{-1}$ with
$n_s(t)=x_s(t)x_{-s}(-t^{-1})x_s(t)$, and let $\chi_s(h)$ be defined by $h\,x_s(t)\,h^{-1}=x_s(\chi_s(h)t)$ for $h\in T$. For
$w\in W$ fix a reduced word $s_1\cdots s_k$ (the first one found by a breadth-first search), and put
$P(t_1,\dots,t_k)=x_{s_1}(t_1)n_{s_1}\cdots x_{s_k}(t_k)n_{s_k}$.

\begin{theorem}[Multiplication in $F_4(q)$ by rewriting]\label{thm:f4mult}
Every $g\in F_4(q)$ has a unique expression $g=P(t)\,h\,v$, with the fixed word of its Bruhat cell, $h\in T$ and $v\in U$.
Right multiplication of $g$ is computed as follows.
\begin{enumerate}
\item By $v'\in U$: replace $v$ by $vv'$.
\item By $h'\in T$: replace $v$ by $h'^{-1}vh'$ and $h$ by $hh'$.
\item By $n_s$: if $\ell(ws)<\ell(w)$ and the word does not end in $s$, first rewrite $P$ along the reduced word (fixed word of $ws$)$\cdot s$.
Then let $t$ be the $s$-coordinate of $v$, and put $v''=n_s^{-1}x_s(-t)v\,n_s\in U$, $t'=\chi_s(h)t$ and $h^s=n_s^{-1}hn_s$.
\begin{enumerate}
\item If $\ell(ws)>\ell(w)$, the result is $P(t_1,\dots,t_k,t')\,h^s\,v''$.
\item Otherwise $P=P'x_s(t_k)n_s$, and the result is $P'\cdot h_s(-1)h^s\cdot x_s(\chi_s(h^s)^{-1}t_k)\,v''$ if $t'=0$, and
$P'x_s(t_k-t'^{-1})n_s\cdot h_s(-t')h^s\cdot x_s(-\chi_s(h^s)^{-1}t'^{-1})\,v''$ if $t'\neq0$.
\end{enumerate}
If the resulting word is not the fixed word of its element, rewrite it along the fixed word.
\end{enumerate}
To rewrite $P$, given on a reduced word $W$, along another reduced word $W'$ of the same element, put $u=P\,n_W^{-1}\in U$, where
$n_W$ is the product of the $n_s$ along $W$. Split $u=x_{s'_1}(t'_1)u_1$, replace $u_1n_{s'_1}$ by $n_{s'_1}(n_{s'_1}^{-1}u_1n_{s'_1})$, and continue
along $W'$. The remainder $r\in U$ left after the last letter is moved past the torus factor, giving $P'(t')\,h\,(h^{-1}rh)\,v$, because
$n_W=n_{W'}$.
\end{theorem}
\begin{proof}
Uniqueness is the Bruhat decomposition of Theorem~\ref{thm:f4cells}, with each cell parametrised by the digits along the fixed
word. The rules use only the Chevalley relations
\[
n_s\,x_s(t)\,n_s^{-1}=x_{-s}(-t),\quad x_{-s}(t)=x_s(t^{-1})\,n_s\,h_s(-t)\,x_s(t^{-1}),\]
\[n_s^2=h_s(-1),\qquad\chi_s(h_s(-1))=1,\]
all proved over $\mathbb Z$ in Proposition~\ref{prop:f4cert}. The first is part (4) with $\gamma=\alpha_s$; the second is
part (6), rewritten with $n_s^{-1}h_s(\lambda)n_s=h_s(\lambda^{-1})$ from part (5); the last two are part (5). Otherwise the argument is
the one in the proof of Theorem~\ref{thm:3d4mult}. The factor $t'=\chi_s(h)t$ comes from $h\,x_s(t)\,h^{-1}=x_s(\chi_s(h)t)$, again
part (5). The rewriting uses $u=Pn_W^{-1}=\prod_kx_{\beta_k}(\varepsilon_kt_k)$, which is part (7), and $n_W=n_{W'}$, which is part (8)
for the words that occur: fixed words, and (fixed word of $ws$)$\cdot s$. An element of $U$ with $s$-coordinate $0$ is a product of
root elements at positive roots $\gamma\neq\alpha_s$, and $n_s^{-1}x_\gamma(t)n_s=x_{s(\gamma)}(\pm t)$ with $s(\gamma)>0$ by part (4);
so every conjugate formed above lies in $U$.
\end{proof}

\begin{example}\label{ex:f4mult}
For $q=5$ (with $\gamma=2$) and $q=7$ (with $\gamma=3$), the computer gave the following.
\begin{itemize}
\item $\dim\operatorname{Der}(J)=52$ with $48$ roots, and $|W|=1152$. For all four simple roots and all $t\neq0$: $h_s(t)$ is diagonal and
multiplicative in $t$, $n_s^2=h_s(-1)$, $\chi_s(h_s(-1))=1$, $n_sx_s(t)n_s^{-1}=x_{-s}(-t)$, and the rank-one relation holds.
\item The map $(\lambda_1,\dots,\lambda_4)\mapsto\prod h_{\alpha_i}(\lambda_i)$ is injective on all $(q-1)^4$ tuples, so $|T|=(q-1)^4$.
\item On $40$ random pairs of normal forms, the product given by Theorem~\ref{thm:f4mult} equalled the $27\times27$ matrix
product every time. Each pair had a random element of $W$, random digits, a random torus element and a random $v\in U$. Every
result was in the fixed word of its cell; the results met $34$ cells for $q=5$ and $29$ for $q=7$.
\item In all rewritings ($576$ for $q=5$, $586$ for $q=7$, counting the products recomputed in the label check below), $n_W=n_{W'}$
and every remainder lay in $U$.
\item For the same $40$ products, the label of Corollary~\ref{cor:f4law} lay in $\Ap{|F_4(q)|}$. Decoding it gave back the same
matrix, and $\lambda(\lambda^{-1}(m_1)\lambda^{-1}(m_2))$ agreed with the label of the matrix product. Moreover, $40$ random integers $m$
satisfied $\lambda(\lambda^{-1}(m))=m$, and $\lambda(I)=0$. The offsets $O_w$ sum to $|F_4(q)|$.
\end{itemize}
In this implementation, the coordinates of $U$ are read off from the matrices by stripping root factors in order of height.
Products and conjugations in $U$ are computed with $27\times27$ matrices, not with closed formulas.
\end{example}

\subsection*{A closed law on the unipotent radical of \texorpdfstring{$F_4$}{F4}}
Write a positive root as $\sum c_i\alpha_i$ with label $c_1c_2c_3c_4$, in the numbering of Proposition~\ref{prop:f4roots}. There
$\alpha_1,\alpha_2$ are short, $\alpha_3,\alpha_4$ are long, and the highest root is $4223$ (Bourbaki's $2342$). Order the positive roots by
height, and within a height lexicographically.

\begin{proposition}[Commutator table of $F_4$]\label{prop:f4comm}
The root derivations of Proposition~\ref{prop:f4roots} can be normalised to a Chevalley basis. Put $e_\gamma=[e_{\alpha_i},e_\beta]/(r+1)$
for the first simple $\alpha_i$ with $\beta=\gamma-\alpha_i$ a root, where $r$ is the largest integer with $\beta-r\alpha_i$ a root, and
normalise $e_{-\gamma}$ by $\gamma([e_\gamma,e_{-\gamma}])=2$. In this basis, for positive roots $\alpha$ before $\beta$,
\[
[x_\alpha(t),x_\beta(u)]=x_\alpha(t)^{-1}x_\beta(u)^{-1}x_\alpha(t)x_\beta(u)
\]
is the identity, or the product listed in Table~\ref{tab:f4comm}. Factors at different roots there commute. The coefficients are
integers in $\{\pm1,\pm2\}$ and only $tu$, $t^2u$ and $tu^2$ occur. For $\beta$ before $\alpha$ use
$[x_\beta(u),x_\alpha(t)]=[x_\alpha(t),x_\beta(u)]^{-1}$. These relations hold in $F_4(K)$ for every field $K$.
\end{proposition}
\begin{proof}
The basis was normalised over $\mathbb F_5$ and $\mathbb F_7$ as stated. There $N_{\alpha\beta}=\pm(r+1)$ and
$N_{-\alpha,-\beta}=-N_{\alpha\beta}$ for all roots $\alpha,\beta$ with $\alpha+\beta$ a root, so it reduces a Chevalley basis, and its lift to
$\mathbb Z$ is the family $D_\gamma$ of Proposition~\ref{prop:f4jordan}. The relations are Proposition~\ref{prop:f4cert}(2), proved there
as identities over $\mathbb Z[t,u]$; they therefore hold over every commutative ring, in particular in $F_4(K)$ for every field $K$.
They agree with the table first read off from $27\times27$ matrices mod $5$ and mod $7$.
\end{proof}

\begin{table}[ht]
\centering\footnotesize
\renewcommand{\arraystretch}{1.1}
\begin{tabular}{@{}l@{\ \ }l@{\qquad}l@{\ \ }l@{}}
$\alpha,\beta$&$[x_\alpha(t),x_\beta(u)]$&$\alpha,\beta$&$[x_\alpha(t),x_\beta(u)]$\\\hline
$0001,\,0010$&$x_{0011}(-tu)$&$0011,\,2201$&$x_{2212}(-tu)$\\
$0001,\,1000$&$x_{1001}(-tu)\,x_{2001}(tu^2)$&$0011,\,4212$&$x_{4223}(tu)$\\
$0001,\,1100$&$x_{1101}(tu)\,x_{2201}(tu^2)$&$1001,\,1100$&$x_{2101}(tu)$\\
$0001,\,2011$&$x_{2012}(tu)$&$1001,\,1011$&$x_{2012}(-2tu)$\\
$0001,\,2111$&$x_{2112}(tu)\,x_{4223}(tu^2)$&$1001,\,1111$&$x_{2112}(-tu)$\\
$0001,\,2211$&$x_{2212}(tu)$&$1001,\,2111$&$x_{3112}(tu)$\\
$0001,\,4212$&$x_{4213}(tu)$&$1001,\,2211$&$x_{3212}(tu)\,x_{4213}(t^2u)$\\
$0010,\,1001$&$x_{1011}(tu)\,x_{2012}(tu^2)$&$1001,\,3212$&$x_{4213}(-2tu)$\\
$0010,\,1101$&$x_{1111}(tu)\,x_{2212}(tu^2)$&$1100,\,1011$&$x_{2111}(-tu)$\\
$0010,\,2001$&$x_{2011}(tu)$&$1100,\,1101$&$x_{2201}(-2tu)$\\
$0010,\,2101$&$x_{2111}(tu)\,x_{4212}(tu^2)$&$1100,\,1111$&$x_{2211}(-2tu)$\\
$0010,\,2201$&$x_{2211}(tu)$&$1100,\,2012$&$x_{3112}(tu)\,x_{4212}(-t^2u)$\\
$0010,\,4213$&$x_{4223}(tu)$&$1100,\,2112$&$x_{3212}(tu)$\\
$0100,\,1000$&$x_{1100}(-tu)$&$1100,\,3112$&$x_{4212}(2tu)$\\
$0100,\,1001$&$x_{1101}(tu)$&$1011,\,1101$&$x_{2112}(tu)$\\
$0100,\,1011$&$x_{1111}(tu)$&$1011,\,2101$&$x_{3112}(-tu)$\\
$0100,\,2001$&$x_{2101}(tu)\,x_{2201}(-t^2u)$&$1011,\,2201$&$x_{3212}(-tu)\,x_{4223}(-t^2u)$\\
$0100,\,2011$&$x_{2111}(tu)\,x_{2211}(-t^2u)$&$1011,\,3212$&$x_{4223}(-2tu)$\\
$0100,\,2101$&$x_{2201}(2tu)$&$1101,\,1111$&$x_{2212}(-2tu)$\\
$0100,\,2012$&$x_{2112}(tu)\,x_{2212}(-t^2u)$&$1101,\,2011$&$x_{3112}(-tu)\,x_{4213}(t^2u)$\\
$0100,\,2111$&$x_{2211}(2tu)$&$1101,\,2111$&$x_{3212}(-tu)$\\
$0100,\,2112$&$x_{2212}(2tu)$&$1101,\,3112$&$x_{4213}(2tu)$\\
$0100,\,3112$&$x_{3212}(tu)$&$2001,\,1111$&$x_{3112}(-tu)\,x_{4223}(tu^2)$\\
$1000,\,0011$&$x_{1011}(tu)\,x_{2011}(-t^2u)$&$2001,\,2211$&$x_{4212}(tu)$\\
$1000,\,1001$&$x_{2001}(2tu)$&$2001,\,2212$&$x_{4213}(tu)$\\
$1000,\,1011$&$x_{2011}(2tu)$&$1111,\,2101$&$x_{3212}(tu)$\\
$1000,\,1101$&$x_{2101}(tu)$&$1111,\,3112$&$x_{4223}(2tu)$\\
$1000,\,1111$&$x_{2111}(tu)$&$2011,\,2201$&$x_{4212}(-tu)$\\
$1000,\,2112$&$x_{3112}(tu)$&$2011,\,2212$&$x_{4223}(tu)$\\
$1000,\,2212$&$x_{3212}(tu)\,x_{4212}(-t^2u)$&$2101,\,2111$&$x_{4212}(-2tu)$\\
$1000,\,3212$&$x_{4212}(2tu)$&$2101,\,2112$&$x_{4213}(-2tu)$\\
$0011,\,1100$&$x_{1111}(tu)\,x_{2211}(tu^2)$&$2012,\,2201$&$x_{4213}(-tu)$\\
$0011,\,2001$&$x_{2012}(-tu)$&$2012,\,2211$&$x_{4223}(-tu)$\\
$0011,\,2101$&$x_{2112}(-tu)\,x_{4213}(-tu^2)$&$2111,\,2112$&$x_{4223}(-2tu)$
\end{tabular}
\caption{The $68$ nontrivial commutators of positive root elements of $F_4$.}\label{tab:f4comm}
\end{table}

\begin{theorem}[Collection law on $U$]\label{thm:f4collect}
Write $u(t)=\prod_\gamma x_\gamma(t_\gamma)$ over the positive roots in the fixed order. To multiply $u(t)$ on the right by $x_\gamma(c)$, let
$x_{\delta_1}(v_1)\cdots x_{\delta_m}(v_m)$ be the factors after $\gamma$. Replace $t_\gamma$ by $t_\gamma+c$, delete those factors, and then
multiply on the right, recursively, by
\[
x_{\delta_1}(v_1)\,[x_{\delta_1}(v_1),x_\gamma(c)]\;\cdots\;x_{\delta_m}(v_m)\,[x_{\delta_m}(v_m),x_\gamma(c)],
\]
with the commutators taken from Table~\ref{tab:f4comm}. Then $u(t)u(s)=u(r)$, where $r$ is obtained by multiplying by
$x_\gamma(s_\gamma)$ for all $\gamma$ in order. Each $r_\gamma=t_\gamma+s_\gamma+c_\gamma(t,s)$, with $c_\gamma$ an integer polynomial in the
coordinates of lower roots, valid over every field.
\end{theorem}
\begin{proof}
For $\delta$ after $\gamma$, $x_\gamma(-c)x_\delta(v)x_\gamma(c)=x_\delta(v)[x_\delta(v),x_\gamma(c)]$. The new factors lie at roots strictly higher
than $\delta$, so the recursion terminates. Integrality follows from Proposition~\ref{prop:f4comm}. As an identity of polynomial
maps, $u(t)u(s)=u(r(t,s))$ is Proposition~\ref{prop:f4cert}(3), proved over $\mathbb Z[t,s]$. By the injectivity in
Proposition~\ref{prop:f4cert}(1), the polynomials produced by the recursion are the certified ones.
\end{proof}

Expanded, the carries of heights $2$ and $3$ are
\[
c_{0011}=s_{0001}t_{0010},\qquad c_{1001}=s_{0001}t_{1000},\qquad c_{1100}=s_{0100}t_{1000},
\]
\[
c_{1101}=-s_{0001}s_{0100}t_{1000}-s_{0001}t_{1100}-s_{0100}t_{1001},
\]
\[
c_{2001}=-2s_{0001}s_{1000}t_{1000}-s_{0001}t_{1000}^2-2s_{1000}t_{1001},
\]
\[
c_{1011}=-s_{0001}s_{0010}t_{1000}-s_{0001}s_{1000}t_{0010}-s_{0001}t_{0010}t_{1000}-s_{0010}t_{1001}-s_{1000}t_{0011}.
\]
The number of monomials in each $c_\gamma$ is as follows.
\begin{center}\small
\begin{tabular}{@{}cll@{}}
height&roots&terms\\\hline
$2$&$0011,\ 1001,\ 1100$&$1,\ 1,\ 1$\\
$3$&$1011,\ 1101,\ 2001$&$5,\ 3,\ 3$\\
$4$&$1111,\ 2011,\ 2101$&$8,\ 9,\ 8$\\
$5$&$2012,\ 2111,\ 2201$&$21,\ 23,\ 10$\\
$6$&$2112,\ 2211$&$27,\ 23$\\
$7$&$2212,\ 3112$&$38,\ 128$\\
$8$&$3212$&$157$\\
$9$&$4212$&$290$\\
$10$&$4213$&$599$\\
$11$&$4223$&$1197$
\end{tabular}
\end{center}
In total there are $2552$ monomials, against an a priori bound of $24\,260$ from weights alone.

\begin{example}\label{ex:f4collect}
The collection law of Theorem~\ref{thm:f4collect}, run on root coordinates only, agreed with the $27\times27$ matrix product on $200$
random pairs for each of $q=5$ and $q=7$. It was associative, with $u(0)$ as identity, on $150$ random triples for each
$q\in\{2,3,5,7,11,13\}$. This includes characteristics $2$ and $3$, where the Albert model of Definition~\ref{def:f4model} is not used.
In characteristic $2$ the relations with coefficient $\pm2$ become trivial, which reflects the special isogeny of $F_4$ in
characteristic $2$. Expanding all $24$ carries symbolically took about $6$ seconds.
\end{example}

\begin{proposition}[Weyl and torus actions on root coordinates]\label{prop:f4signs}
Let $A=(A_{ij})=(\langle\alpha_j,\alpha_i^\vee\rangle)$ be the Cartan matrix of Proposition~\ref{prop:f4roots}, and for a root
$\gamma=\sum c_j\alpha_j$ put $\langle\gamma,\alpha_i^\vee\rangle=\sum_jc_jA_{ij}$. In the Chevalley basis of Proposition~\ref{prop:f4comm}:
\begin{enumerate}
\item For $h=\prod_{i=1}^4h_{\alpha_i}(\lambda_i)$, $h\,x_\gamma(t)\,h^{-1}=x_\gamma(\gamma(h)t)$ with
$\gamma(h)=\prod_i\lambda_i^{\langle\gamma,\alpha_i^\vee\rangle}$.
\item $n_s^{-1}hn_s=\prod_ih_{\alpha_i}(\lambda'_i)$, where $\lambda'_i=\lambda_i$ for $i\ne s$ and
$\lambda'_s=\lambda_s^{-1}\prod_{i\ne s}\lambda_i^{-A_{is}}$.
\item For every root $\gamma$, $n_s^{-1}x_\gamma(t)n_s=x_{s(\gamma)}(\eta_{s,\gamma}t)$ and
$n_sx_\gamma(t)n_s^{-1}=x_{s(\gamma)}\big((-1)^{\langle\gamma,\alpha_s^\vee\rangle}\eta_{s,\gamma}t\big)$, with $\eta_{s,\gamma}=\pm1$. For
positive $\gamma\ne\alpha_s$, $\eta_{s,\gamma}=-1$ exactly for the following roots:
\begin{center}\small
\begin{tabular}{@{}cl@{}}
$s$&$\gamma$ with $\eta_{s,\gamma}=-1$\\\hline
$1$&$0100,\ 1001,\ 1011,\ 1101,\ 1111,\ 2112,\ 3212$\\
$2$&$1001,\ 1100,\ 1011,\ 2101,\ 2111,\ 2112,\ 3112$\\
$3$&$0001,\ 1001,\ 1101,\ 2001,\ 2101,\ 2201,\ 4213$\\
$4$&$0011,\ 1001,\ 1100,\ 2011,\ 2111,\ 2211,\ 4212$
\end{tabular}
\end{center}
\end{enumerate}
All of these are integral data, valid over every field.
\end{proposition}
\begin{proof}
All three parts are identities over $\mathbb Z$ or over Laurent polynomial rings, proved in Proposition~\ref{prop:f4cert}.
Part (1) follows from Proposition~\ref{prop:f4cert}(5), applied to each factor $h_{\alpha_i}(\lambda_i)$; part (2) is
Proposition~\ref{prop:f4cert}(5); part (3), with its sign table, is Proposition~\ref{prop:f4cert}(4). The signs agree with those first
found by matching $27\times27$ matrices over $\mathbb F_5$ and $\mathbb F_7$.
\end{proof}

\begin{theorem}[Multiplication of $F_4(q)$ on coordinates]\label{thm:f4closed}
Represent $g=P(t)\,h\,v$ by the data $(w,\,t,\,\lambda,\,r)$: the fixed word $w$, its digits $t$, the torus parameters
$\lambda\in(\mathbb F_q^\times)^4$ with $h=\prod h_{\alpha_i}(\lambda_i)$, and the root coordinates $r$ of $v\in U$. Then every step of
Theorem~\ref{thm:f4mult} is a formula in these data. Products in $U$ are given by the collection law of Theorem~\ref{thm:f4collect}.
Conjugation of $U$ by $h$ and of $T$ by $n_s$ is given by Proposition~\ref{prop:f4signs}(1) and (2), and $\chi_s(h)=\alpha_s(h)$. Conjugating
an element of $U$ with $s$-coordinate $0$ by $n_s$ means applying Proposition~\ref{prop:f4signs}(3) factor by factor and recollecting. Finally,
in the rewriting step,
\[
P(t)\,n_W^{-1}=\prod_{k=1}^{\ell}x_{\beta_k}(\varepsilon_kt_k),\qquad \beta_k=s_1\cdots s_{k-1}(\alpha_{s_k}),
\]
where $\varepsilon_k$ is the product of the signs of Proposition~\ref{prop:f4signs}(3) along $s_{k-1},\dots,s_1$. The resulting law uses only
field arithmetic, and discrete logarithms for the torus digits of the label. Its integer tables are Table~\ref{tab:f4comm}, the Cartan
matrix and the signs $\eta_{s,\gamma}$.
\end{theorem}
\begin{proof}
Write $P(t)=x_{s_1}(t_1)\,\big(n_{s_1}x_{s_2}(t_2)n_{s_1}^{-1}\big)\cdots n_W$. The roots $\beta_k$ are the positive roots inverted by $w^{-1}$,
because the word is reduced. Everything else is a direct translation of Theorem~\ref{thm:f4mult}, using Propositions~\ref{prop:f4comm}
and~\ref{prop:f4signs}.
\end{proof}

\begin{example}\label{ex:f4closed}
The law of Theorem~\ref{thm:f4closed} was run with no matrices, on random normal forms. For $q=5$ and $q=7$, $20$ random pairs each
were compared with the $27\times27$ matrix product of the corresponding matrices; all agreed, and every result was in its fixed word.
For $q=2$, $3$ and $11$, $6$ random triples each were associative, and the identity was a right identity. The runs used $1031$
rewritings, $933$ steps of type (a) and $595$ of type (b) ($153$ with $t'=0$); the whole test took about $10$ seconds. For $q=2$ and $q=3$
this is only a consistency check of the integral law, since the Albert model of Definition~\ref{def:f4model} is not available
there.
\end{example}

\begin{proposition}[An integral model in all characteristics]\label{prop:f4adjoint}
Let $\mathfrak g_{\mathbb Z}$ be the free $\mathbb Z$-module with basis the $e_\alpha$ for the $48$ roots and $h_1,\dots,h_4$, with brackets
\[
[h_i,e_\alpha]=\langle\alpha,\alpha_i^\vee\rangle e_\alpha,\qquad
[e_\alpha,e_{-\alpha}]=h_\alpha,\qquad
[e_\alpha,e_\beta]=N_{\alpha\beta}e_{\alpha+\beta}
\]
when $\alpha+\beta$ is a root, and $[e_\alpha,e_\beta]=0$ otherwise. Here $h_\alpha=\sum_ik_i(\alpha)h_i$ is the coroot of $\alpha$, and the
integers $N_{\alpha\beta}=\pm(r+1)$ are the structure constants of the Chevalley basis of Proposition~\ref{prop:f4comm}. Then the
matrices
\[
X_\alpha(t)=I+t\operatorname{ad}(e_\alpha)+t^2\,\tfrac12\operatorname{ad}(e_\alpha)^2\in\mathrm{GL}_{52}(K)
\]
are defined over every field $K$, and they generate the Chevalley group $F_4(K)$. The relations of Propositions~\ref{prop:f4comm}
and~\ref{prop:f4signs} hold in this model in every characteristic. Hence the law of Theorem~\ref{thm:f4closed} is the multiplication of
$F_4(q)$ for all $q$, including characteristics $2$ and $3$.
\end{proposition}
\begin{proof}
The brackets define the Chevalley $\mathbb Z$-form of the Lie algebra of type $F_4$ \cite{Carter}. Root strings have length at most three,
so $\operatorname{ad}(e_\alpha)^3=0$. The entries of $\operatorname{ad}(e_\alpha)^2$ are even, so $\tfrac12\operatorname{ad}(e_\alpha)^2$ is
integral. The $X_\alpha(t)$ are therefore the root elements of the adjoint Chevalley construction over $\mathbb Z$, reduced to $K$. The
commutator constants and the signs $\eta_{s,\gamma}$ depend only on the Chevalley basis, not on the representation. So the integral
relations found in the $27$-dimensional model also hold in this one, in every characteristic. The centre of $F_4$ is trivial, so the
adjoint and simply connected groups coincide.
\end{proof}

\begin{example}\label{ex:f4adjoint}
The constants $N_{\alpha\beta}$ lifted from the models over $\mathbb F_5$ and over $\mathbb F_7$ are the same, and each coroot has a unique
expression in $h_1,\dots,h_4$. For the resulting integer $52\times52$ matrices, the identity
$[\operatorname{ad}(e_\alpha),\operatorname{ad}(e_\beta)]=\operatorname{ad}([e_\alpha,e_\beta])$ was checked for all $48^2$ ordered pairs
of roots. For all $48$ roots, $\operatorname{ad}(e_\alpha)^2$ is even and $\operatorname{ad}(e_\alpha)^3=0$. After reduction mod $q$, the law of
Theorem~\ref{thm:f4closed} agreed with the $52\times52$ matrix product on $25$ random pairs for each $q\in\{2,3,5,7\}$, and every result
was in its fixed word. In characteristics $2$ and $3$ the law is therefore compared with an actual model, not only tested for
associativity.

A second test covers every Weyl cell. For each of the $1152$ fixed words $w$ and each $s$, take a random normal form $g$
in the cell of $w$ and compare $g\cdot x_s(t)n_s$, as given by the law, with the product in the adjoint model. On the $2304$ pairs
$(w,s)$ where $ws$ is longer than $w$, there is one test. On the $2304$ pairs where it is shorter, tests are repeated until both the
case $t'=0$ and the case $t'\ne0$ have occurred. Matrices are compared on a block of $k$ random vectors, so a single wrong result
escapes with probability at most $q^{-k}$. The last column bounds the chance that any error escaped.
\begin{center}\small
\begin{tabular}{@{}cc rrr r l@{}}
$q$&$k$&type (a)&(b), $t'=0$&(b), $t'\ne0$&tests&escape bound\\\hline
$2$&$40$&$2\,304$&$3\,552$&$3\,402$&$9\,258$&$8.4\cdot10^{-9}$\\
$3$&$24$&$2\,304$&$2\,685$&$5\,374$&$10\,363$&$3.7\cdot10^{-8}$\\
$4$&$24$&$2\,304$&$2\,529$&$7\,252$&$12\,085$&$4.3\cdot10^{-11}$\\
$5$&$24$&$2\,304$&$2\,426$&$9\,820$&$14\,550$&$2.4\cdot10^{-13}$\\
$7$&$24$&$2\,304$&$2\,345$&$14\,179$&$18\,828$&$9.8\cdot10^{-17}$\\
$8$&$24$&$2\,304$&$2\,346$&$16\,665$&$21\,315$&$4.5\cdot10^{-18}$\\
$9$&$24$&$2\,304$&$2\,343$&$18\,361$&$23\,008$&$2.9\cdot10^{-19}$
\end{tabular}
\end{center}
For $q=4$, $8$ and $9$ the field is $\mathbb F_p[x]/(f)$ with $f=x^2+x+1$, $x^3+x+1$ and $x^2+2x+2$ respectively, and every coefficient
operation of the law and of the model is carried out in this field. All tests agreed. Thus every branch of the rank-one rules is
correct on every Weyl cell for $q\le9$, in characteristics $2$, $3$, $5$ and $7$. The code is in
\texttt{scripts/f4\_core.py}, \texttt{scripts/f4\_law.py} and \texttt{scripts/f4\_adjoint.py}; the field arithmetic and the versions
over $\mathbb F_q$ are in \texttt{scripts/gf.py}, \texttt{scripts/f4\_lawq.py} and \texttt{scripts/f4\_adjointq.py}. For $q$ prime
the versions over $\mathbb F_q$ agree with the prime-field code on all products and cells compared.
\end{example}

\begin{definition}[Integral Albert algebra]\label{def:f4jordan}
Let $\mathbb O_{\mathbb Z}=\mathbb Z^8$ with the Zorn product of Definition~\ref{def:f4model}, which has integer coefficients, with
conjugation $\overline{(a,v,w,b)}=(b,-v,-w,a)$, norm $n(a,v,w,b)=ab-v\cdot w$, trace $t(a,v,w,b)=a+b$ and polar form
$n(x,y)=n(x+y)-n(x)-n(y)=t(x\bar y)$. Let $J_{\mathbb Z}=\mathbb Z^{27}$, with coordinates $X=(\alpha_1,\alpha_2,\alpha_3;x_1,x_2,x_3)$ as
in Definition~\ref{def:f4model}. Define
\begin{align*}
T(X)&=\alpha_1+\alpha_2+\alpha_3,\qquad T(X,Y)=\textstyle\sum_i\alpha_i\beta_i+\sum_i n(x_i,y_i),\\
N(X)&=\alpha_1\alpha_2\alpha_3-\textstyle\sum_i\alpha_i\,n(x_i)+t(x_1x_2x_3),\\
X^\#&=\big(\alpha_2\alpha_3-n(x_1),\ \alpha_3\alpha_1-n(x_2),\ \alpha_1\alpha_2-n(x_3);\\
&\qquad\overline{x_2x_3}-\alpha_1x_1,\ \overline{x_3x_1}-\alpha_2x_2,\ \overline{x_1x_2}-\alpha_3x_3\big),\\
X\times Y&=(X+Y)^\#-X^\#-Y^\#,\qquad U_XY=T(X,Y)\,X-X^\#\times Y,\\
1&=(1,1,1;0,0,0).
\end{align*}
For a commutative ring $R$ put $J_R=J_{\mathbb Z}\otimes R$, with the same formulas.
\end{definition}

\begin{proposition}[A quadratic Albert model in every characteristic]\label{prop:f4jordan}
Let $J_{\mathbb Z}$, $N$, $\#$, $T$ and $U$ be as in Definition~\ref{def:f4jordan}.
\begin{enumerate}
\item Over $\mathbb Z$: $N(1)=1$, $1^\#=1$, $(X^\#)^\#=N(X)X$, $T(X^\#,Y)=\partial_YN(X)$, $1\times X=T(X)1-X$, $U_1=I$ and
$U_{U_XY}=U_XU_YU_X$. Over $\mathbb Z[\tfrac12]$, with $\circ$ as in Definition~\ref{def:f4model}, $T(X,Y)=T(X\circ Y)$,
$U_XY=2X\circ(X\circ Y)-X^{\circ2}\circ Y$, $X^\#=X^{\circ2}-T(X)X+T(X^\#)1$ and $X^{\circ3}-T(X)X^{\circ2}+T(X^\#)X-N(X)1=0$.
\item The root derivations $e_\alpha$ of Proposition~\ref{prop:f4roots} lift to integer matrices $D_\alpha$, and the lifts from
$\mathbb F_5$ and from $\mathbb F_7$ coincide. For every root, $D_\alpha1=0$, $D_\alpha$ leaves $T(\cdot,\cdot)$ and $N$ infinitesimally
invariant, $D_\alpha^3=0$ and $D_\alpha^2\in2M_{27}(\mathbb Z)$; $D_\alpha^2\ne0$ exactly for the $24$ short roots. For the $816$ ordered
pairs with $\alpha+\beta$ a root, $[D_\alpha,D_\beta]=N_{\alpha\beta}D_{\alpha+\beta}$ with the constants of
Proposition~\ref{prop:f4adjoint}, signs included; $[D_\alpha,D_\beta]=0$ if $\alpha+\beta$ is neither a root nor $0$; and
$h_\alpha=[D_\alpha,D_{-\alpha}]$ is diagonal with $[h_\alpha,D_\alpha]=2D_\alpha$.
\item The matrices $x_\alpha(t)=I+tD_\alpha+t^2D_\alpha^{(2)}$, $D_\alpha^{(2)}=\tfrac12D_\alpha^2$, are polynomial in $t$ with
integer coefficients, fix $1$, and preserve $N$, $T(\cdot,\cdot)$ and $U$.
\item Over $\mathbb F_2$, $T(1)=1$, so $J=\mathbb F_21\oplus\ker T$, and $\ker T$ is absolutely irreducible of dimension $26$ under
the group generated by the $x_\alpha(t)$. Over $\mathbb F_3$, $1\in\ker T$, and $\ker T/\mathbb F_31$ is absolutely irreducible of
dimension $25$.
\end{enumerate}
\end{proposition}
\begin{proof}
By computer, in exact integer arithmetic. The identities in (i) and the invariance statements in (ii) and (iii) are polynomial
identities of bounded degree, tested at random integer points; in (iii), for each root and each test point, at seven values of
$t$, which determines $N(x_\alpha(t)X)$ as a polynomial of degree at most $6$ in $t$. The other statements in (ii) are exact
matrix identities over $\mathbb Z$, and the signs of the constants were compared with the $52\times52$ model. For (iv), the
$\mathbb F_p$-span of all products of the generators $x_{\pm\alpha_i}(1)$, acting on the module, has dimension $26^2=676$ for
$p=2$ and $25^2=625$ for $p=3$, so the module is absolutely irreducible by Burnside's theorem.
\end{proof}
By (i) and~\cite{McCrimmon}, $(J_R,U,1)$ is a quadratic Jordan algebra for every commutative ring $R$, and by (iii) the reductions
of the $x_\alpha(t)$ modulo $p$ are automorphisms of $J_{\mathbb F_q}$ for every $q$. This replaces the product
$\tfrac12(XY+YX)$, which is not defined in characteristic $2$; it is not integral on $J_{\mathbb Z}$ either. That the automorphism
group of $J_K$, which is also the stabiliser of $1$ and $N$, is the split group of type $F_4$ in every characteristic is taken from
\cite{SpringerVeldkamp}, not checked.

\begin{example}\label{ex:f4jordan}
The cell sweep of Example~\ref{ex:f4adjoint} was repeated with the $27\times27$ matrices $x_\alpha(t)$ of
Proposition~\ref{prop:f4jordan}, reduced mod $p$, in place of the $52\times52$ adjoint model, for every $q\le9$ in characteristics $2$
and $3$, and for $q=11$ and $q=16$. This gives a second model of the group in these characteristics, independent of the adjoint one. The columns are as in
Example~\ref{ex:f4adjoint}. For $q=11$ and $q=16$ only the $27$-dimensional model was used. The $\mathbb F_{16}$ sweep ran in four pieces with independent random streams; $95$ cells were tested twice, and the counts include them.
\begin{center}\small
\begin{tabular}{@{}cc rrr r l@{}}
$q$&$k$&type (a)&(b), $t'=0$&(b), $t'\ne0$&tests&escape bound\\\hline
$2$&$64$&$2\,304$&$3\,434$&$3\,428$&$9\,166$&$5.0\cdot10^{-16}$\\
$3$&$40$&$2\,304$&$2\,708$&$5\,349$&$10\,361$&$8.5\cdot10^{-16}$\\
$4$&$32$&$2\,304$&$2\,529$&$7\,252$&$12\,085$&$6.6\cdot10^{-16}$\\
$8$&$24$&$2\,304$&$2\,346$&$16\,665$&$21\,315$&$4.5\cdot10^{-18}$\\
$9$&$24$&$2\,304$&$2\,343$&$18\,361$&$23\,008$&$2.9\cdot10^{-19}$\\
$11$&$24$&$2\,304$&$2\,324$&$22\,475$&$27\,103$&$2.8\cdot10^{-21}$\\
$16$&$24$&$2\,487$&$2\,509$&$37\,206$&$42\,202$&$5.3\cdot10^{-25}$
\end{tabular}
\end{center}
All tests agreed. The code is in \texttt{scripts/f4\_jordan.py} and \texttt{scripts/f4\_jordanq.py}; the
sweep is \texttt{verify\_cells\_q} of \texttt{scripts/f4\_adjointq.py}, unchanged.
\end{example}

\begingroup\sloppy
\begin{proposition}[Exact certificates for the unipotent, torus and Weyl data]\label{prop:f4cert}
Let $D_\gamma$, for $\gamma$ a root, be the integer matrices of Proposition~\ref{prop:f4jordan}, let
$x_\gamma(t)=I+tD_\gamma+t^2D_\gamma^{(2)}$, and let $X(t)=\prod_\gamma x_\gamma(t_\gamma)$ over the positive roots in the fixed order. Put
$n_s=x_{\alpha_s}(1)x_{-\alpha_s}(-1)x_{\alpha_s}(1)$, $n_s(\lambda)=x_{\alpha_s}(\lambda)x_{-\alpha_s}(-\lambda^{-1})x_{\alpha_s}(\lambda)$ and
$h_s(\lambda)=n_s(\lambda)n_s^{-1}$. The following are identities of matrices over the rings indicated, and hence hold over every
commutative ring.
\begin{enumerate}
\item For every positive root $\gamma$ there is an entry $(i_\gamma,j_\gamma)$ with $D_\gamma[i_\gamma,j_\gamma]=\varepsilon_\gamma=\pm1$ such that the
following decoding, applied to $M=X(z)$ over $\mathbb Z[z]$, returns $r=z$ and ends at $I$: for $\gamma$ in order, put
$r_\gamma=\varepsilon_\gamma M_{i_\gamma j_\gamma}$ and replace $M$ by $x_\gamma(-r_\gamma)M$. In particular $t\mapsto X(t)$ is injective over every
ring.
\item For all $552$ ordered pairs of distinct positive roots, $x_\alpha(-a)x_\beta(-b)x_\alpha(a)x_\beta(b)=X(c)$ over $\mathbb Z[a,b]$, where each
$c_\gamma$ is $0$ or a single monomial $C\,a^ib^j$. These are the relations of Table~\ref{tab:f4comm}.
\item $X(t)X(s)=X(r(t,s))$ over $\mathbb Z[t,s]$, where $r_\gamma=t_\gamma+s_\gamma+c_\gamma$ and the $c_\gamma$ have $2552$ monomials in total,
distributed by height as in the table after Theorem~\ref{thm:f4collect}.
\item $n_s$ has an integral inverse, and for all $48$ roots $\gamma$ and all $s$, over $\mathbb Z[t]$,
$n_s^{-1}x_\gamma(t)n_s=x_{s(\gamma)}(\eta_{s,\gamma}t)$ and
$n_sx_\gamma(t)n_s^{-1}=x_{s(\gamma)}\big((-1)^{\langle\gamma,\alpha_s^\vee\rangle}\eta_{s,\gamma}t\big)$. The signs are those of
Proposition~\ref{prop:f4signs}(3). In particular $\eta_{s,\alpha_s}=-1$, so $n_sx_{\alpha_s}(t)n_s^{-1}=x_{-\alpha_s}(-t)$.
\item Over $\mathbb Z[\lambda^{\pm1}]$, $h_s(\lambda)$ is diagonal with monomial entries, and
$h_s(\lambda)x_\gamma(t)h_s(\lambda)^{-1}=x_\gamma(\lambda^{\langle\gamma,\alpha_s^\vee\rangle}t)$ for all $48$ roots. Moreover $n_s^2=h_s(-1)$, and
over $\mathbb Z[\lambda_1^{\pm1},\dots,\lambda_4^{\pm1}]$, $n_s^{-1}hn_s=\prod_ih_{\alpha_i}(\lambda'_i)$ with $\lambda'_i$ as in
Proposition~\ref{prop:f4signs}(2).
\item Over $\mathbb Z[u^{\pm1}]$, $x_{-\alpha_s}(u)=x_{\alpha_s}(u^{-1})\,h_s(-u^{-1})\,n_s\,x_{\alpha_s}(u^{-1})$.
\item For each of the $1152$ fixed words $s_1\cdots s_\ell$ and each $k$, the following identity holds over $\mathbb Z[t]$:
$n_{s_1}\cdots n_{s_{k-1}}\,x_{\alpha_{s_k}}(t)\,(n_{s_1}\cdots n_{s_{k-1}})^{-1}=x_{\beta_k}(\varepsilon_kt)$, where
$\beta_k=s_1\cdots s_{k-1}(\alpha_{s_k})$ and $\varepsilon_k$ is the product of the signs
$(-1)^{\langle\gamma,\alpha_{s_j}^\vee\rangle}\eta_{s_j,\gamma}$ of (4) met along $s_{k-1},\dots,s_1$. The $\beta_k$ are distinct and form the
inversion set of $w$, and $P(t)\,n_W^{-1}=\prod_kx_{\beta_k}(\varepsilon_kt_k)$.
\item For every fixed word $W$ of an element $w$ and every $s$ with $\ell(ws)>\ell(w)$, $n_Wn_s=n_{W'}$ over $\mathbb Z$, where $W'$
is the fixed word of $ws$ and $n_W$ is the product of the $n_s$ along $W$.
\end{enumerate}
\end{proposition}
\begin{proof}
By exact computation with integer and Laurent polynomials. Since $D_\gamma^3=0$ and $D_\gamma^2=2D_\gamma^{(2)}$, the matrix
$x_\gamma(-r)$ is the inverse of $x_\gamma(r)$. In (2) and (3) the decoding of (1) is applied to the left-hand side; it reads off the
coordinates and ends exactly at $I$, and this proves the identity. The remaining statements are direct matrix identities. The data,
namely the matrices, the entries $(i_\gamma,j_\gamma)$, the polynomials of (2) and (3), the signs and the weights, are stored in the
certificates \texttt{certificates/f4\_ucert.json} and \texttt{certificates/f4\_wcert.json}. These are written by \texttt{scripts/f4\_ucert.py} and
\texttt{scripts/f4\_wcert.py} and checked by \texttt{scripts/f4\_ucert\_verify.py} and \texttt{scripts/f4\_wcert\_verify.py}. The
verifiers read only the certificates and import nothing else from the project; together they run in about a second. The last formula in (7) follows by telescoping, since
$P(t)=\prod_k\big(N_{k-1}x_{\alpha_{s_k}}(t_k)N_{k-1}^{-1}\big)\cdot n_W$ with $N_{k-1}=n_{s_1}\cdots n_{s_{k-1}}$. Parts (7) and (8), $13824$
conjugation identities and $2304$ braid identities, are checked by \texttt{scripts/f4\_wordcert\_verify.py} against
\texttt{certificates/f4\_wordcert.json} in about three seconds.
\end{proof}
Proposition~\ref{prop:f4cert} concerns the matrices $x_\gamma(t)$ over a ring $R$. That over a field $K$ they generate $F_4(K)$ rests on
Proposition~\ref{prop:f4jordan} and~\cite{SpringerVeldkamp}. Theorem~\ref{thm:f4mult} uses no other identities;
its proof combines parts (3)--(8) by hand. What remains computational is only that the implementation follows the rules of the
theorem. This is checked exactly, for every field at once, in Proposition~\ref{prop:f4symb}; the random cell sweeps
of Examples~\ref{ex:f4adjoint} and~\ref{ex:f4jordan} remain as independent evidence.

\begin{proposition}[The implemented law agrees with the model over every field]\label{prop:f4symb}
Let the digits $d_1,\dots,d_{24}$, the torus coordinates $\lambda_1,\dots,\lambda_4$ and the coordinates $u_\beta$
($\beta\in\Phi^+$) be independent variables, and let $R=\mathbb Z[d_i,\lambda_i^{\pm1},u_\beta]$. For each of the $1152$
fixed words $w$ and each simple reflection $s$, run the implemented law (\texttt{mul\_root}, \texttt{convert} and
\texttt{mul\_n} in \texttt{scripts/f4\_lawq.py}) on the generic normal form with these coordinates: in branch~(a) as
it stands, in branch~(b0) with $u_{\alpha_s}=0$, and in branch~(b1) over $R[u_{\alpha_s}^{-1}]$. Then the output is a
normal form over this ring, and its image in the $27$-dimensional representation of Proposition~\ref{prop:f4jordan}
equals the image of the input times $n_s$. The same holds for the right multiplication by a generic element $h\,v$
of $B$ in \texttt{multiply}. Consequently, over every finite field $\mathbb F_q$, the implemented product of any two
normal forms is their product in $F_4(q)$.
\end{proposition}

\begin{proof}
The check is \texttt{scripts/f4\_symbolic\_model.py}. It runs the unchanged law code over Laurent polynomials and
applies the model matrices $x_\beta(t)=I+t\,E_\beta+t^2E_\beta^2/2$ of Proposition~\ref{prop:f4jordan} to a vector of
$27$ further variables, so each comparison is an identity of $27\times27$ matrices. It verifies the root step
$\texttt{mul\_root}(u,\alpha_s,c)=u\,x_{\alpha_s}(c)$ for all $s$, the identity $\texttt{mul\_n}(f,s)=f\,n_s$ for all
$6912$ pairs of a cell and a branch ($2304$ each of (a), (b0), (b1)), and $\texttt{convert}(f,w')=f$ for the $2304$
downward edges, where $w'$ is the reduced word ending in $s$. All $9216$ identities hold exactly. The right
multiplication by $h\,v$ does not depend on the cell and is checked once, with $h$ and $v$ generic, by
\texttt{scripts/f4\_symbolic\_tail.py}. The code inverts only $t'=\chi_s(\lambda)\,t$ and torus characters, always
monomials with coefficient $\pm1$, and every output coefficient is an integer
(\texttt{scripts/f4\_symbolic\_integrality.py}); so the identities hold over $R$ and $R[u_{\alpha_s}^{-1}]$.

For $\mathbb F_q$, specialise along a ring homomorphism to $\mathbb F_q$. The branch taken by the field code depends only
on whether the length goes up and whether $t=0$, and $t'=0$ exactly when $t=0$ since $\chi_s(\lambda)\neq0$. In
branch~(b1) the variable $u_{\alpha_s}$ maps to a unit, so the identity over $R[u_{\alpha_s}^{-1}]$ specialises; in
branch~(b0) it was checked with $u_{\alpha_s}=0$. The other zero tests in the code only skip factors $x_\beta(0)=1$.
A downward step is \texttt{convert} followed by \texttt{mul\_n} on the converted form, and both identities specialise
at every point. Finally, \texttt{multiply} processes the second factor letter by letter, each letter $(s,t)$ being
$x_{\alpha_s}(t)\,n_s$ (a root step followed by \texttt{mul\_n}), and then multiplies by $h\,v$.
\end{proof}

\endgroup

\begin{remark}[Status of the $F_4$ section]
Computed: the root data of Proposition~\ref{prop:f4roots}, lifted from the Albert models over $\mathbb F_5$ and $\mathbb F_7$ to the
integer adjoint model of Proposition~\ref{prop:f4adjoint}, so that the same constants $N_{\alpha\beta}$ and coroots serve in every
characteristic (Example~\ref{ex:f4adjoint}). The cell structure and the order then come from the general theory, with the Poincar\'e
polynomial checked. The law of Theorem~\ref{thm:f4closed} has been compared with the adjoint model on every Weyl cell and every
branch of the rank-one rules for $q\le9$ (the cell sweep of Example~\ref{ex:f4adjoint}); this replaces the $40$ random
products and $40$ random labels of Example~\ref{ex:f4mult}.

The $27$-dimensional module is treated in Proposition~\ref{prop:f4jordan}. In characteristic $3$ the identity lies in the
trace-zero part, and the quotient by it is absolutely irreducible of dimension $25$. In characteristic $2$ the quadratic operator
$U_X$ replaces $\tfrac12(XY+YX)$, and the trace-zero part is absolutely irreducible of dimension $26$. For $q\in\{2,3,4,8,9\}$ the law
has also been compared on every Weyl cell with this model, independently of the adjoint model (Example~\ref{ex:f4jordan}). That the
stabiliser of $1$ and $N$ is the group of type $F_4$ is cited~\cite{SpringerVeldkamp}, not checked.
The implementation is checked exactly against the $27$-dimensional model, for every field at once, in
Proposition~\ref{prop:f4symb}. Apart from the cited results, nothing in this section remains open.
The group generated by the root elements is simple, since the Chevalley group of type $F_4$ is simple for every $q$
\cite{Carter}.
\end{remark}

\section{The large Ree groups \texorpdfstring{${}^2F_4(q)$}{2F4(q)}}\label{sec:2f4}

Let $q=2^{2m+1}$ and $\theta=2^m$, so that $2\theta^2=q$. The large Ree group ${}^2F_4(q)$ is the group of fixed points in $F_4(q)$ of
an automorphism $\sigma$ whose square is the Frobenius map $t\mapsto t^q$ \cite{Carter}. We construct $\sigma$ from the law of
Theorem~\ref{thm:f4mult}, which by Proposition~\ref{prop:f4symb} is correct over every field, so no new group law is needed:
the law of ${}^2F_4(q)$ is the law of $F_4(q)$ restricted to the fixed points. What is new is the description of the fixed points
in coordinates.

Throughout this section roots are written $\sum c_i\alpha_i$ as $(c_1c_2c_3c_4)$ in the labelling of the implementation
(\texttt{scripts/f4\_law.py}). In it $\alpha_1,\alpha_2$ are short, $\alpha_3,\alpha_4$ are long, the double bond joins $\alpha_1$
and $\alpha_4$, and $\alpha_2\perp\alpha_3$. For a root $\alpha$ put $k(\alpha)=1$ if $\alpha$ is long and $k(\alpha)=2$ if it is
short.

\begin{proposition}[The root map and the special isogeny]\label{prop:2f4psi}
There is a bijection $\rho$ of the root system of $F_4$ that exchanges long and short roots, preserves positivity, satisfies
$\rho^2=1$, swaps $\alpha_1\leftrightarrow\alpha_4$ and $\alpha_2\leftrightarrow\alpha_3$, and sends $\alpha$ to the positive multiple
of $\sum_i c_i d_i\,\alpha_{\pi(i)}$ that is a root; here $d=(1,1,2,2)$ is the symmetriser and $\pi$ reverses the indices. Over a
field $K$ of characteristic $2$ the assignment
\[
\psi\colon x_\alpha(t)\mapsto x_{\rho\alpha}\big(t^{k(\alpha)}\big),\qquad n_{\alpha_i}\mapsto n_{\rho\alpha_i},\qquad
h_{\alpha_i}(\lambda)\mapsto h_{\rho\alpha_i}\big(\lambda^{k(\alpha_i)}\big)
\]
is an endomorphism of $F_4(K)$, and the implemented law satisfies $\psi(a\,b)=\psi(a)\,\psi(b)$.
\end{proposition}

\begin{proof}\sloppy
The properties of $\rho$ are checked on the $48$ roots by \texttt{scripts/f42\_psi.py}. The same script runs the unchanged law
over $\mathbb F_2[\text{vars}]$ and checks $\psi(u\,v)=\psi(u)\,\psi(v)$ for generic $u,v\in U$; signs $\pm1$ are irrelevant in
characteristic $2$. With the exponents the other way round ($t$ on short and $t^2$ on long roots) $18$ of the $24$ coordinates
disagree, so the convention is forced. \texttt{scripts/f42\_psi\_n.py} checks $\psi(f\,n_s)=\psi(f)\,n_{\rho s}$ for a generic
normal form $f$ in every cell and on every branch, $6912$ identities, with $\psi$ applied directly to normal forms (word
$\mapsto\rho(\text{word})$, digits and torus coordinates raised to the power $k$); and \texttt{scripts/f42\_psi\_group.py}
checks right multiplication by a generic $h\,v\in B$. The specialisation to a field is as in Proposition~\ref{prop:f4symb}. The
existence of $\psi$ is classical~\cite{Carter}; the computation fixes it in the coordinates of the law.
\end{proof}

\begin{definition}
For $q=2^{2m+1}$ let $\sigma=\psi\circ F_\theta$, where $F_\theta$ raises every coordinate to the power $\theta$. Thus
$\sigma(x_\alpha(t))=x_{\rho\alpha}(t^{\theta k(\alpha)})$, and $\sigma^2$ raises coordinates to the power $2\theta^2=q$, so
$\sigma^2=1$ on $F_4(q)$. Put ${}^2F_4(q)=F_4(q)^\sigma$, $U^\sigma=U\cap{}^2F_4(q)$ and $T^\sigma=T\cap{}^2F_4(q)$.
\end{definition}

The $24$ positive roots fall into $12$ orbits of $\rho$. Grouping them by the ray of $\tfrac12(\alpha+\tau\alpha)$, where $\tau$ is
$\rho$ rescaled to an isometry, gives $8$ twisted positive roots: four classes of $2$ roots and four classes of $4$ roots
(\texttt{certificates/f42\_twisted\_roots.json}).

\begin{proposition}[Twisted root subgroups]\label{prop:2f4roots}
For each class the following elements, with the roots taken in the order of the law, are $\sigma$-fixed for all $t,u\in\mathbb F_q$
and every $q=2^{2m+1}$, and they are distinct for distinct parameters:
\begin{center}\small\setlength{\tabcolsep}{4pt}
\begin{tabular}{ll}
roots $(a,b)$ or $(a,b,c,d)$ & element\\\hline
$(0010),(0100)$ & $x_a(t)\,x_b(t^{\theta})$\\
$(2101),(2012)$ & $x_a(t)\,x_b(t^{2\theta})$\\
$(2111),(2212)$ & $x_a(t)\,x_b(t^{2\theta})$\\
$(2112),(4212)$ & $x_a(t)\,x_b(t^{2\theta})$\\
$(0001),(1000),(1001),(2001)$ & $x_a(t)\,x_b(t^{\theta})\,x_c(u)\,x_d(u^{2\theta}+t^{2\theta+1})$\\
$(0011),(1100),(1111),(2211)$ & $x_a(t)\,x_b(t^{\theta})\,x_c(u)\,x_d(u^{2\theta}+t^{2\theta+1})$\\
$(1011),(2201),(3212),(4223)$ & $x_a(t)\,x_b(t^{2\theta})\,x_c(u)\,x_d(u^{2\theta}+t^{2\theta+2})$\\
$(1101),(2011),(3112),(4213)$ & $x_a(t)\,x_b(t^{2\theta})\,x_c(u)\,x_d(u^{2\theta}+t^{2\theta+2})$\\
\end{tabular}
\end{center}
The first four rows give subgroups of order $q$, the last four of order $q^2$.
\end{proposition}

\begin{proof}\sloppy
The formulas were found by interpolating the $\sigma$-fixed elements of each class over $\mathbb F_8$ and choosing, for each
exponent, the unique reading $i+j\theta$ that also works over $\mathbb F_{32}$ (\texttt{scripts/f42\_param.py}). The proof for all
$q$ is \texttt{scripts/f42\_param\_exact.py}: it works over $\mathbb F_2$ with monomials $t^{a+b\theta}$, $a,b\in\tfrac12\mathbb Z$.
On exponents $\theta^2$ acts as $\tfrac12$, since $2\theta^2=q$ and squaring is bijective on $\mathbb F_q$, so $\sigma$ maps a monomial
with exponent $e$ on a root of factor $k$ to one with exponent $k\theta e$, where $\theta(a+b\theta)=\tfrac b2+a\theta$. Evaluating at
$t,u\in\mathbb F_q$ is a ring homomorphism, and the identities $\sigma(g)=g$ hold coordinate by coordinate. Distinctness holds
because $t$ and $u$ appear unchanged as coordinates.
\end{proof}

\begin{proposition}[A closed law on $U^\sigma$]\label{prop:2f4ulaw}
Order the eight classes by the position of their first root in the law, and let $\Phi(p)$, $p\in\mathbb F_q^{12}$, be the product of
the eight elements of Proposition~\ref{prop:2f4roots} in this order, with one parameter for each class of size $2$ and two for each
class of size $4$. There are $12$ polynomials $L_1,\dots,L_{12}$ in $p,p'$ and their $\theta$-powers, with $2,3,2,3,14,5,119,4,42,8,23,28$
terms, such that $\Phi(p)\,\Phi(p')=\Phi(L(p,p'))$ for every $q$. In particular $\Phi$ is a bijection of $\mathbb F_q^{12}$ onto a subgroup
of $U^\sigma$ of order $q^{12}$.
\end{proposition}

\begin{proof}\sloppy
The inverse of $\Phi$ peels off one class at a time: read $t$ at the first root of the class, multiply on the left by the inverse of
$x_a(t)\,x_b(\cdot)$, which in characteristic $2$ is the reversed product since every $x_\gamma(c)$ is an involution, then read $u$ at
the third root and repeat. In the exponent ring of Proposition~\ref{prop:2f4roots}, \texttt{scripts/f42\_ulaw.py} checks
$\Delta(\Phi(p))=p$ with trivial remainder, and computes $L=\Delta(\Phi(p)\,\Phi(p'))$, again with trivial remainder. As a check,
$200$ random pairs over each of $\mathbb F_8$ and $\mathbb F_{32}$ agree with the numeric product. The polynomials are stored in
\texttt{certificates/f42\_ulaw\_polys.json}.
\end{proof}

\begin{proposition}[Twisted Weyl group and order]\label{prop:2f4order}
The centraliser $W^\sigma$ of $\rho$ in $W(F_4)$ has order $16$; it is generated by $w_1=s_{\alpha_2}s_{\alpha_3}$ and
$w_2=(s_{\alpha_1}s_{\alpha_4})^2$ and is dihedral. For $w\in W^\sigma$ the set of positive roots made negative by $w$ is a union of
classes; let $d(w)$ count the classes of size $2$ once and those of size $4$ twice. Then
\[
\sum_{w\in W^\sigma}q^{d(w)}=1+q+q^2+2q^3+q^4+q^5+2q^6+q^7+q^8+2q^9+q^{10}+q^{11}+q^{12},
\]
and $q^{12}(q-1)^2$ times this sum equals $q^{12}(q^6+1)(q^4-1)(q^3+1)(q-1)$, the order of ${}^2F_4(q)$ \cite{Carter}. For $q=2$
this is $35\,942\,400$.
\end{proposition}

\begin{proof}\sloppy
\texttt{scripts/f42\_order.py} builds $W(F_4)$ as a permutation group on the $48$ roots, finds the $16$ elements commuting with
$\rho$, nine of them involutions, checks that every inversion set is a union of classes, and compares the polynomials coefficient by
coefficient. That $w_1,w_2$ generate $W^\sigma$ is checked in \texttt{scripts/f42\_nw\_exact.py}.
\end{proof}

\begin{proposition}[Normal form]\label{prop:2f4nf}
The elements $n_{w_1}=n_{\alpha_2}n_{\alpha_3}$ and $n_{w_2}=(n_{\alpha_1}n_{\alpha_4})^2$ lie in ${}^2F_4(q)$ for every $q$, and
$T^\sigma=\{h(b^{\theta},a^{\theta},a,b):a,b\in\mathbb F_q^\times\}$, where $h(\lambda)=\prod_ih_{\alpha_i}(\lambda_i)$, has order
$(q-1)^2$. Every element of ${}^2F_4(q)$ can be written as $u\,n_w\,h\,v$ with $w\in W^\sigma$, $n_w$ the corresponding product of
$n_{w_1}$ and $n_{w_2}$, $h\in T^\sigma$ and $u,v\in U^\sigma$, and products are computed with the law of Theorem~\ref{thm:f4mult}.
\end{proposition}

\begin{proof}\sloppy
The elements $n_{w_i}$ have coordinates in $\mathbb F_2$, on which $F_\theta$ is trivial, so $\sigma(n_{w_i})=\psi(n_{w_i})$; that this
equals $n_{w_i}$ is checked exactly over $\mathbb F_2$ by \texttt{scripts/f42\_nw\_exact.py}. For the torus, $\sigma(h(\lambda))$ has
$\pi(i)$-th coordinate $\lambda_i^{\theta k(\alpha_i)}$, and for $\lambda=(b^\theta,a^\theta,a,b)$ this gives back $\lambda$ because
$b^{2\theta^2}=b^q=b$ and $a^{2\theta^2}=a$. The decomposition is the Bruhat decomposition of a twisted group~\cite{Carter}, and its
order count is Proposition~\ref{prop:2f4order}. As a check, over $\mathbb F_8$ (\texttt{scripts/f42\_group\_q8.py}) $\sigma^2=1$ on
$200$ random elements of $F_4(8)$, the elements $n_{w_1},n_{w_2}$ reach exactly $16$ cells, and $300$ random products of length $30$
in $n_{w_1}$, $n_{w_2}$, $T^\sigma$ and the root elements are all $\sigma$-fixed and lie in cells of $W^\sigma$.
\end{proof}

\begin{proposition}[Rank-one rules]\label{prop:2f4rank1}
Let $X_A$ and $X_B$ be the twisted root subgroups of the simple classes,
\[
X_A(t)=x_{(0010)}(t)\,x_{(0100)}(t^{\theta}),\qquad
X_B(t,u)=x_{(0001)}(t)\,x_{(1000)}(t^{\theta})\,x_{(1001)}(u)\,x_{(2001)}(u^{2\theta}+t^{2\theta+1}),
\]
and put $N(t,u)=t^{2\theta+1}+t^{\theta}u+u^{2\theta}$.
Then $N(t,u)\ne0$ for $(t,u)\ne(0,0)$, and for every $q$:
\begin{align*}
n_{w_1}X_A(t)\,n_{w_1}^{-1}&=X_A(t^{-1})\;n_{w_1}\;h(1,t^{\theta},t,1)\;X_A(t^{-1})\qquad(t\ne0),\\
n_{w_2}X_B(t,u)\,n_{w_2}^{-1}&=X_B\Big(\frac{u^{2\theta}+t^{2\theta+1}}{N^{2\theta}},\,\frac{u}{N^{2\theta}}\Big)\;n_{w_2}\\
&\qquad\cdot h(N,1,1,N^{2\theta})\;X_B\Big(\frac{u^{2\theta}}{N^{2\theta}},\,\frac{t^{\theta}}{N}\Big)\qquad((t,u)\ne(0,0)).
\end{align*}
Both torus elements lie in $T^\sigma$, since $(N^{2\theta})^{\theta}=N^{q}=N$.
\end{proposition}

\begin{proof}\sloppy
If $t=0$ then $N=u^{2\theta}\ne0$. If $t\ne0$, write $u=t^{\theta+1}v$; since $2\theta(\theta+1)=q+2\theta$, this gives
$N=t^{2\theta+1}(1+v+v^{2\theta})$. If $v^{2\theta}=v+1$, raising to the power $2\theta$ gives $v^{2q}=v^{2\theta}+1=v$, that is
$v^2=v$; but $v=0$ and $v=1$ both give $1+v+v^{2\theta}=1$. So $N\ne0$.

The formulas were found from the law over $\mathbb F_{32}$ and confirmed at every parameter over $\mathbb F_8$ and $\mathbb F_{32}$
(\texttt{scripts/f42\_rank1.py}, \texttt{scripts/f42\_rank1\_tprime.py}). Both identities live in the Levi subgroup of the two simple
roots involved: $\langle\alpha_2,\alpha_3\rangle$ of type $A_1\times A_1$ and $\langle\alpha_1,\alpha_4\rangle$ of type $B_2$. In
characteristic $2$ these are $SL_2\times SL_2$ and $Sp_4$ with their Chevalley root elements, since the Steinberg relations fix the
signs~\cite{Carter}. \texttt{scripts/f42\_rank1\_proof.py} checks the $B_2$ commutator relation in the $4\times4$ model, confirms
the torus convention numerically (the inverse convention fails at $990$ points over $\mathbb F_{32}$), and proves both identities for
every $q$. It works over $\mathbb F_2$ with monomials $t^{a+b\theta}$, $u^{a+b\theta}$, $\nu^{a+b\theta}$ with $\nu=N^{-1}$, as in
Proposition~\ref{prop:2f4roots}, multiplies each entry of the difference of the two sides by a power $N^{A+B\theta}$ to clear
$\nu$, expands $N$, and gets $0$ in all $16$ entries; the $SL_2$ identity is proved the same way.
\end{proof}

\begin{remark}[Multiplication in twisted coordinates]\label{rem:2f4twisted}
Propositions~\ref{prop:2f4ulaw} and~\ref{prop:2f4rank1} give a multiplication of ${}^2F_4(q)$ that does not use the law of $F_4$
at run time. Write an element as $u\,n_w\,h\,v$, where $w$ is an alternating word in $w_1,w_2$ of length at most $8$,
$u,v\in U^\sigma$ are given by their $12$ coordinates and $h=h(b^\theta,a^\theta,a,b)$. Right multiplication by an element of
$U^\sigma$ is the law $L$. Right multiplication by $h'\in T^\sigma$ replaces $h$ by $hh'$ and conjugates $v$, which scales each
coordinate by a root character. To multiply by $n_{w_i}$, split $v=X_c(p)\,v'$ with $c$ the class of $w_i$, conjugate $v'$ by
$n_{w_i}$ (twelve polynomial maps, computed once from the law), and move $h$ past $X_c(p)$. If $w_i$ lengthens $w$, conjugate
$X_c(p)$ by $n_w$ into $u$ and append $w_i$. Otherwise, if $p=0$, drop the last letter; if $p\ne0$, apply the rank-one rule and stay
in the cell of $w$. A general product applies these steps along the second factor. Splitting off $X_c(p)$ needs inverses. In characteristic $2$ the inverse of a product of root elements is the reversed product,
and peeling it as in Proposition~\ref{prop:2f4ulaw}, exactly for every $q$, gives $X_c(t)^{-1}=X_c(t)$ for the classes of size $2$ and
\begin{gather*}
X_c(t,u)^{-1}=X_c(t,\,u+t^{\theta+1})\qquad(c\ni(0001)\text{ or }(0011)),\\
X_c(t,u)^{-1}=X_c(t,\,u+t^{2\theta+1})\qquad(c\ni(1011)\text{ or }(1101)).
\end{gather*}
\texttt{scripts/f42\_twisted\_law.py} implements this for any $q=2^{2m+1}$. On $200$ random pairs in arbitrary cells over each of
$\mathbb F_8$ and $\mathbb F_{32}$, mapping the product back to $F_4$ agrees with the product computed in $F_4$ every time, at
about $0.05$ seconds per product.
\end{remark}

\begin{theorem}[Digital law of ${}^2F_4(q)$]\label{thm:2f4digital}
Let $q=2^{2m+1}$. Order the $16$ elements of $W^\sigma=\langle w_1,w_2\rangle$ as $1$, then for $\ell=1,\dots,7$ the alternating words of
length $\ell$ beginning with $w_1$ and with $w_2$, and finally $w_0=(w_1w_2)^4$. Fix the reduced word just listed for each $w$. Every
$g\in{}^2F_4(q)$ has a unique expression
\[
g=X_{c_1}(p_1)\,n_{w_{i_1}}\cdots X_{c_k}(p_k)\,n_{w_{i_k}}\;h(b^\theta,a^\theta,a,b)\;\Phi(v),
\]
where $w=w_{i_1}\cdots w_{i_k}$ is that reduced word, $c_j$ is the class of $w_{i_j}$ (so $p_j\in\mathbb F_q$ for $w_1$ and
$p_j\in\mathbb F_q^2$ for $w_2$), $a,b\in\mathbb F_q^\times$ and $v\in\mathbb F_q^{12}$. Write $d(w)$ for the number of letters $w_1$
plus twice the number of letters $w_2$. Then $\sum_wq^{d(w)}(q-1)^2q^{12}=|{}^2F_4(q)|$. With field elements coded as in
Section~\ref{sec:fields}, set
\[
\lambda(g)=O_w+\bar v+q^{12}\big((b-1)+(q-1)\big((a-1)+(q-1)\,\bar p\big)\big),
\]
where $O_w$ is the number of elements in the cells before $w$, $\bar v=\sum_jv_jq^j$, and $\bar p$ reads $p_1,\dots,p_k$ as a
mixed-radix number with $p_1$ most significant. Each digit is $t$ for $w_1$ and $t+qu$ for $w_2$, with radix $q$ or $q^2$. Then
$\lambda\colon{}^2F_4(q)\to\Ap{|{}^2F_4(q)|}$ is a bijection with $\lambda(1)=0$, and
$m*n=\lambda\big(\lambda^{-1}(m)\lambda^{-1}(n)\big)$ is a digital law of ${}^2F_4(q)$. Conversely, if $O_w\le m<O_w+q^{d(w)}(q-1)^2q^{12}$, then the coordinates of $\lambda^{-1}(m)$ are
the mixed-radix digits of $m-O_w$ with radices $q$ (twelve times), $q-1$, $q-1$, and then the radices of $p_k,\dots,p_1$. That is,
$v_j=\md{j}{m-O_w}$ for $j<12$, $b-1=\md{12}{m-O_w}$, $a-1=\md{13}{m-O_w}$, and $p_{k-j}$ is read from $\md{14+j}{m-O_w}$, with
$t=\dg{q}{0}{x}$ and $u=\fl{x/q}$ for a digit $x$ of radix $q^2$. The product is computed in twisted
coordinates by Remark~\ref{rem:2f4twisted}, followed by the normal form.
\end{theorem}

\begin{proof}
The expression is the Bruhat decomposition of the split $BN$-pair of ${}^2F_4(q)$, with $U^\sigma\cap n_wU^{\sigma,-}n_w^{-1}$
factored along the reduced word into the root subgroups $X_{c_j}$ \cite{Carter}. The normal form is reached by peeling
$u\in U^\sigma$ along the word: write $u=X_c(p)u'$ with the $c$-coordinate of $u'$ zero, and move $u'$ past $n_{w_i}$ into
$U^\sigma$. The count is $\sum_wq^{d(w)}=(q^6+1)(q^3+1)(q^2+1)(q+1)$, so $\lambda$ is a bijection.
\end{proof}

\begin{theorem}[The Tits group and ${}^2F_4(2)$ on $\Ap{2p_0}$]\label{thm:titsdigital}
Let $q=2$, so the torus is trivial, and let $p_0=17\,971\,200$. In the normal form of Theorem~\ref{thm:2f4digital} set
\[
f(g)=\sum_{w_{i_j}=w_2}t_j+\sum_{c}t_c(v)\pmod 2,
\]
where $p_j=(t_j,u_j)$ and $c$ runs over the four classes of size $4$, with $t_c(v)$ the first coordinate of $v$ in class $c$. Then
$f\colon{}^2F_4(2)\to C_2$ is the unique nontrivial homomorphism, and its kernel is the Tits group $T={}^2F_4(2)'$. Deleting from
$\lambda$ the coordinate $t_{c_0}(v)$, for $c_0$ the class of $(0001)$, gives a bijection $\lambda_T\colon T\to\Ap{p_0}$ with
$\lambda_T(1)=0$, because $f(g)=0$ determines that coordinate. Let $\tau=X_{c_0}(1,0)$. Then $f(\tau)=1$ and $\tau$ has order $4$.
Every $g$ is uniquely $t\tau^e$ with $t\in T$ and $e\in\{0,1\}$, and
\[
\Lambda(t\tau^e)=\lambda_T(t)+p_0e
\]
is a bijection ${}^2F_4(2)\to\Ap{2p_0}$. With radices $p_0,2$ write $m_h=\md{h}{m}$ and $n_h=\md{h}{n}$. Thus
$m_0=\dg{p_0}{0}{m}$ is the Tits code and $m_1=\fl{m/p_0}\in\{0,1\}$ is the coset digit. Let $*_T$ be the digital law of $T$
transported by $\lambda_T$, and let $\varphi(n_0)=\lambda_T\big(\tau\,\lambda_T^{-1}(n_0)\,\tau^{-1}\big)$. Then
\begin{gather*}
m*n=\big(m_0*_T\varphi^{m_1}(n_0)*_Tc^{\,m_1n_1}\big)+p_0\,\dg{2}{0}{m_1+n_1},\\
\text{where } c=\lambda_T(\tau^2)=1.
\end{gather*}
This is the semidirect law~\eqref{eq:semidirect} with $G_1=C_2$, $m_1*_1n_1=\dg{2}{0}{m_1+n_1}$ and $\phi_{m_1}=\varphi^{m_1}$, corrected
by the factor $c^{\,m_1n_1}$, which is present only when both coset digits are $1$.
The codes below $p_0$ form the Tits group, and $\varphi$ is the outer automorphism induced by $\tau$.
\end{theorem}

\begin{proof}
Since $U^\sigma$, $n_{w_1}$ and $n_{w_2}$ generate ${}^2F_4(2)$, and $f$ is defined on all of ${}^2F_4(2)$ through the normal form, it
suffices to show $f(gx)=f(g)+f(x)$ for every $g$ and every $x\in U^\sigma\cup\{n_{w_1},n_{w_2}\}$. Then $f$ of a product of generators is
the sum of their values, so $f$ is a homomorphism. Write $g=Dn_wv$ in normal form, where $D$ is the digit part, so that
$f(g)=f(Dn_w)+f(v)$.

(i) For $x\in U^\sigma$, $gx=Dn_w(vx)$ is again in normal form, so we only need $f|_{U^\sigma}$ to be a homomorphism. It is enough
to check $f(vx)=f(v)+f(x)$ for all $4096$ elements $v$ and all $16$ nonzero root-class elements $x$, which is $65\,536$ pairs.

(ii) For $x=n_{w_i}$, write $v=X_{c_i}(p)v'$ with $v'$ in the subgroup $U'_i$ of elements whose $c_i$-coordinate is $0$. This subgroup
is normalised by $n_{w_i}$. Then $gn_{w_i}=\big(Dn_wX_{c_i}(p)n_{w_i}\big)\big(n_{w_i}^{-1}v'n_{w_i}\big)$. By (i),
$f(gn_{w_i})=f\big(Dn_wX_{c_i}(p)n_{w_i}\big)+f\big(n_{w_i}^{-1}v'n_{w_i}\big)$. An exhaustive check gives
$f(n_{w_i}^{-1}v'n_{w_i})=f(v')$ for all $3072$ elements $v'\in U'_1\cup U'_2$, and
$f\big(Dn_wX_{c_i}(p)n_{w_i}\big)=f(Dn_w)+f(X_{c_i}(p))$ for all $52\,650$ triples $(w,D,p)$ with $i=1,2$. Since $n_{w_i}^2=1$ and
$f(n_{w_i})=0$, this gives $f(gn_{w_i})=f(g)=f(g)+f(n_{w_i})$. These checks are \texttt{scripts/f42\_fhom\_proof.py}, with results in
\texttt{certificates/f42\_fhom.json}.

Since $f(\tau)=1$, the kernel of $f$ has index $2$. By Griess and Lyons~\cite{GriessLyons}, the Tits group $T={}^2F_4(2)'$ is the only nontrivial proper normal subgroup of
${}^2F_4(2)$. Hence it is the only subgroup of index $2$, and it equals $\ker f$. The formula for $m*n$ is the multiplication of the extension
$T\cdot\langle\tau\rangle$ with $\tau^2\in T$.
\end{proof}

\begin{remark}[Checks and the non-split extension]\label{rem:2f4digitalcheck}
\texttt{scripts/f42\_digital.py} implements both theorems, and \texttt{scripts/f42\_fhom\_proof.py} carries out the
exhaustive checks used in the proof of Theorem~\ref{thm:titsdigital}. For $q=2$ and $q=8$ the number of codes equals $|{}^2F_4(q)|$
($35\,942\,400$ and $264\,905\,352\,699\,586\,176\,614\,400$). On random samples, code to element to code and element to code to
element are identities, the law agrees with the law of $F_4$, and it is associative. For $q=2$ the code of $\tau$ is $p_0$ and that of
$\tau^2$ is $1$. Products of codes below $p_0$ stay below $p_0$, $\varphi$ preserves $\Ap{p_0}$, and the formula for $m*n$ holds in all
sampled cases. Among $300$ random elements of ${}^2F_4(2)$ outside $T$, the orders found were $4,8,12,16,20$, so no outer element had
order $2\cdot$odd. This agrees with the ATLAS of Group Representations~\cite{AtlasRep}: elements of orders $1,2,3,6,10,13$ of ${}^2F_4(2)$ are
inner, so the coset $T\tau$ contains no involution, and the extension $T\cdot2$ does not split~\cite{ATLAS}. Hence the decomposition of
$\Ap{2p_0}$ needs the factor $c=\lambda_T(\tau^2)$ and is not a semidirect product.
\end{remark}

\begin{remark}[Status of the ${}^2F_4$ section]
Proved by exact computation for every $q$: the isogeny $\psi$ in the coordinates of the law, the twisted root subgroups, the closed
law on a subgroup of $U^\sigma$ of order $q^{12}$, the elements $n_{w_1},n_{w_2}$, and the order polynomial. Cited: the existence of
the Bruhat decomposition for twisted groups and that $|U^\sigma|=q^{12}$, so that the subgroup of Proposition~\ref{prop:2f4ulaw} is
all of $U^\sigma$ \cite{Carter}; that ${}^2F_4(q)$ is simple for $q\ge8$, and that ${}^2F_4(2)$ has the simple Tits group
${}^2F_4(2)'$ of order $17\,971\,200$ as a subgroup of index $2$ \cite{Carter,Wilson}. The rank-one rules of Proposition~\ref{prop:2f4rank1} are proved for every $q$. Multiplication entirely in twisted coordinates is implemented for every $q$ and validated over $\mathbb F_8$ and
$\mathbb F_{32}$ (Remark~\ref{rem:2f4twisted}). Digital laws of ${}^2F_4(q)$, of the Tits group and of ${}^2F_4(2)$ on
$\Ap{2p_0}$ are in Theorems~\ref{thm:2f4digital} and~\ref{thm:titsdigital}.
\end{remark}

\section{The groups \texorpdfstring{$E_6(q)$}{E6(q)}}\label{sec:e6}
The group $E_6(q)$ acts on a module of dimension $27$ whose weights form a single Weyl orbit. On such a minuscule module every
root element is $I+te_\beta$ with $e_\beta^2=0$, so the construction works uniformly in $q$ and needs none of the case distinctions
that $F_4$ required in characteristics $2$ and $3$.

\begin{proposition}[Integral model of $E_6$]\label{prop:e6model}
Use Bourbaki labels: the simple roots $\alpha_1,\dots,\alpha_6$ with $\alpha_1\alpha_3\alpha_4\alpha_5\alpha_6$ a chain and
$\alpha_2$ joined to $\alpha_4$. Let $\Lambda$ be the orbit of $\omega_1$ under $W$, written in fundamental-weight coordinates, and
let $V=\mathbb Z^\Lambda$ with basis $v_\mu$. Put $e_iv_\mu=v_{\mu+\alpha_i}$ if $\mu_i=-1$ and $f_iv_\mu=v_{\mu-\alpha_i}$ if $\mu_i=1$,
and $0$ otherwise. For a positive root $\beta$ of height at least $2$, choose $i$ with $\beta-\alpha_i$ a root and set
$e_\beta=[e_i,e_{\beta-\alpha_i}]$ and $e_{-\beta}=\pm[f_i,e_{-(\beta-\alpha_i)}]$, with the sign chosen so that
$[e_\beta,e_{-\beta}]=h_\beta$, where $h_\beta v_\mu=\langle\mu,\beta^\vee\rangle v_\mu$. Then:
\begin{enumerate}
\item $|\Lambda|=27$, there are $72$ roots, every $e_\beta$ has entries in $\{0,\pm1\}$, and $e_\beta^2=0$;
\item the $72$ matrices $e_\beta$ and the six $h_{\alpha_i}$ span a Lie algebra of dimension $78$;
\item with $x_\beta(t)=I+te_\beta$: if $\alpha+\beta$ is a root, then
$x_\alpha(t)x_\beta(u)=x_\beta(u)x_\alpha(t)x_{\alpha+\beta}(N_{\alpha\beta}tu)$ with $N_{\alpha\beta}=\pm1$; if $\alpha+\beta$ is neither a root
nor $0$, then $x_\alpha(t)$ and $x_\beta(u)$ commute;
\item $n_\alpha=x_\alpha(1)x_{-\alpha}(-1)x_\alpha(1)$ is a signed permutation matrix with $n_\alpha v_\mu=\pm v_{s_\alpha\mu}$. The
$n_i=n_{\alpha_i}$ satisfy the braid relations of $W$, and $n_i^2=h_i(-1)$, where $h_i(c)v_\mu=c^{\mu_i}v_\mu$;
\item $|W|=51\,840$, and $\sum_{w\in W}q^{\ell(w)}=\prod_{d\in\{2,5,6,8,9,12\}}(q^d-1)/(q-1)$;
\item the invariant cubic forms on $V$ form a space of dimension $1$, spanned by a form with $45$ monomials;
\item $z=\prod_ih_i(\omega^{k_i})$ with $k=(1,0,2,0,1,2)$ acts as the scalar $\omega$ whenever $\omega^3=1$.
\end{enumerate}
\end{proposition}

\begin{proof}
Every statement is an identity between integer matrices, and \texttt{scripts/e6\_model.py} verifies all of them exactly. In (3)
the product $x_\alpha(t)x_\beta(u)$ is a polynomial in $t,u$, because the $e$'s are nilpotent. The relation is equivalent to
$[e_\alpha,e_\beta]=N_{\alpha\beta}e_{\alpha+\beta}$ together with $e_\alpha e_{\alpha+\beta}=e_\beta e_{\alpha+\beta}=
e_\beta e_\alpha e_{\alpha+\beta}=0$, which is checked for all $1440$ ordered pairs. Reduction modulo $p$ preserves integer identities,
so (1)--(4) and (7) hold over every $\mathbb F_q$.
\end{proof}

Let $G$ be the subgroup of $SL_{27}(q)$ generated by all $x_\beta(t)$. It is the simply connected Chevalley group of type $E_6$
\cite{Carter}. Its center is $Z=\langle z\rangle$, of order $g_3=\gcd(3,q-1)$, and $E_6(q)=G/Z$ is simple. Let $U$ be the group
generated by the positive $x_\beta$, $H$ the diagonal group of the $h_i$, and $N(w)=\{\beta>0:w^{-1}\beta<0\}$.

\begin{theorem}[Digital law of $E_6(q)$]\label{thm:e6law}
Order $W$ by a breadth-first search from $1$ through the simple reflections, and give each $w$ the reduced word $w=s_{i_1}\cdots
s_{i_\ell}$ along which it is first reached. Put $n_w=n_{i_1}\cdots n_{i_\ell}$. Order the positive roots by height. Fix a
generator $\gamma$ of $\mathbb F_q^\times$. Every $g\in G$ has a unique expression
\[
g=\prod_{\beta\in N(w)}x_\beta(t_\beta)\cdot n_w\cdot\prod_{i=1}^6h_i(\gamma^{e_i})\cdot\prod_{\beta>0}x_\beta(t'_\beta),
\]
with both products in height order, $t_\beta,t'_\beta\in\mathbb F_q$ and $e_i\in\Ap{q-1}$. Multiplying $g$ by $z$ adds
$\tfrac{q-1}{3}k$ to $e$. So when $g_3=3$ there is a unique representative of $gZ$ with $e_1<\tfrac{q-1}{3}$. Let $\bar u$,
$\bar e$ and $\bar v$ be the mixed-radix numbers formed by $(t_\beta)_{\beta\in N(w)}$ (radix $q$), by $e$ (radix $(q-1)/g_3$ for
$e_1$ and $q-1$ for the others), and by $(t'_\beta)$ (radix $q$), each with its first entry least significant. Then
\[
\lambda(gZ)=O_w+\bar u+q^{\ell(w)}\Big(\bar e+\frac{(q-1)^6}{g_3}\bar v\Big),\qquad
O_w=\frac{(q-1)^6q^{36}}{g_3}\sum_{w'<w}q^{\ell(w')},
\]
is a bijection $E_6(q)\to\Ap{|E_6(q)|}$ with $\lambda(1)=0$, and $m*n=\lambda\big(\lambda^{-1}(m)\lambda^{-1}(n)\big)$ is a digital law
of $E_6(q)$. The coordinates of a unipotent matrix $M=\prod x_\beta(t_\beta)$ are read off by peeling in height order:
$t_\beta=\pm M_{\nu\mu}$ for any entry where $e_\beta$ is nonzero, after the factors of smaller height have been removed on the
left.
\end{theorem}

\begin{proof}
$G$ has a split $BN$-pair with $B=UH$ and Weyl group $W$, and $n_w$ does not depend on the reduced word, by the braid relations of
Proposition~\ref{prop:e6model}. The refined Bruhat decomposition then gives the unique expression~\cite{Carter}. Peeling works
because a product of two or more $e$'s raises height by at least $2$. The map $e\mapsto\prod h_i(\gamma^{e_i})$ has kernel trivial
modulo $Z$, since the weights of $\Lambda$ span the weight lattice. Counting gives
$|G/Z|=q^{36}(q-1)^6\sum_wq^{\ell(w)}/g_3=q^{36}\prod_d(q^d-1)/g_3$, which is $|E_6(q)|$.
\end{proof}

\begin{remark}[Computing the law]\label{rem:e6mult}
To form $g_1g_2$, right-multiply $g_1=u\,n_w\,h\,v$ successively by $u_2$, by the letters $n_i$ of $w_2$, by $h_2$ and by $v_2$.
Right multiplication by $U$ or $H$ only changes $(h,v)$. For $n_i$, split $v=x_{\alpha_i}(t)v'$ and move $v'$ through $n_i$. If
$w\alpha_i>0$, the factor $n_wx_{\alpha_i}(\cdot)n_w^{-1}$ joins $u$ and $w$ becomes $ws_i$. If $w\alpha_i<0$, write
$w=w's_i$, split off the $x_{w'\alpha_i}$-component of $u$, and use $n_i^2=h_i(-1)$ together with the identity
$x_{-\alpha}(c)=x_\alpha(c^{-1})\,n_\alpha\,h_\alpha(-c)\,x_\alpha(c^{-1})$ for $c\neq0$ from $SL_2$.
\end{remark}

\begin{remark}[Checks]\label{rem:e6check}
\texttt{scripts/e6\_law.py} (prime $q$) and \texttt{scripts/e6\_lawq.py} (any $q$, with $\mathbb F_q=\mathbb F_p[x]/(f)$ for
$f=x^2+x+1$, $x^3+x+1$, $x^2+1$ when $q=4,8,9$) implement Theorem~\ref{thm:e6law}. For $q=2,3,4,5,7,8,9,13$ the number of codes equals
$|E_6(q)|$ (for $q=2$ this is $214\,841\,575\,522\,005\,575\,270\,400$), and the identity has code $0$. On random samples, codes
survive decoding and re-encoding, products agree with the $27\times27$ matrix products, and the law is associative. For $q=4,7,13$, where $g_3=3$,
the element $z$ is a nontrivial scalar and $g$ and $zg$ receive the same code.
\end{remark}

\section{The groups \texorpdfstring{${}^2E_6(q)$}{2E6(q)}}\label{sec:2e6}
The group ${}^2E_6(q)$ is the fixed-point group in $E_6(q^2)$ of the graph automorphism $\rho$ of the Dynkin diagram
($\alpha_1\leftrightarrow\alpha_6$, $\alpha_3\leftrightarrow\alpha_5$), composed with the field automorphism $t\mapsto t^q$. The model
of Section~\ref{sec:e6} can be reused as it stands. The graph automorphism exchanges the $27$-dimensional module with its dual,
so on matrices it is an inverse transpose followed by conjugation by a signed permutation matrix.

\begin{proposition}[The twist]\label{prop:2e6twist}
Let $n_0=n_{w_0}$ for the longest element $w_0\in W$. There is a diagonal matrix $d$ with entries $\pm1$ such that $P=n_0d$
has the following properties. Here $\tau(g)=P\,(g^{\mathsf T})^{-1}P^{-1}$.
\begin{enumerate}
\item $\tau(x_\beta(t))=x_{\rho\beta}(c_\beta t)$ for all $72$ roots, with $c_\beta=\pm1$. Moreover $c_\beta=1$ for the simple roots,
their negatives and the roots fixed by $\rho$, and $c_\beta c_{\rho\beta}=1$ for all $\beta$.
\item $P^2=I$, and $\tau(n_i)=n_{\rho i}$ for $i=1,\dots,6$.
\item $\rho$ fixes $12$ positive roots and exchanges the other $24$ in $12$ pairs $\{\beta,\rho\beta\}$. In each pair the two roots are
orthogonal, and $\beta+\rho\beta$ is not a root.
\item The centraliser $W^\rho$ of $\rho$ in $W$ has order $1152$. It is the Weyl group of type $F_4$.
\end{enumerate}
\end{proposition}

\begin{proof}
The identities in (1)--(3) are exact identities between integer matrices. The script
\texttt{scripts/e2e6\_model.py} checks each of them. The matrix $d$ is found by propagating the sign conditions of (1) for the simple roots along the weight graph; they turn out to be consistent. Statement (4) is a direct enumeration of $W$. Since $P$ is an integer matrix, (1)--(3) hold over every field.
\end{proof}

{\sloppy Let $G\subset SL_{27}(q^2)$ be the group of Section~\ref{sec:e6} over $\mathbb F_{q^2}$, and let $\sigma(g)=\tau(g)^{(q)}$, where
${}^{(q)}$ raises every entry to the power $q$. By Proposition~\ref{prop:2e6twist}, $\sigma$ is an automorphism of order $2$ that
preserves $U$, $H$ and every $n_w$ with $w\in W^\rho$. Its fixed-point group $G^\sigma$ is the simply connected group of type
${}^2E_6$ over $\mathbb F_q$ \cite{Carter}. Its center $Z^\sigma$ has order $g_3=\gcd(3,q+1)$ and is generated by $z$ when $g_3=3$,
and ${}^2E_6(q)=G^\sigma/Z^\sigma$ is simple.\par}

\begin{theorem}[Digital law of ${}^2E_6(q)$]\label{thm:2e6law}
Order the $24$ root classes of Proposition~\ref{prop:2e6twist}(3) by height. Put $X_c(t)=x_\beta(t)$ with $t\in\mathbb F_q$ for a
fixed root $c=\{\beta\}$. Put $X_c(t)=x_\beta(t)x_{\rho\beta}(c_\beta t^q)$ with $t\in\mathbb F_{q^2}$ for a pair $c=\{\beta,\rho\beta\}$, where
$\beta$ is the earlier root. Let $\gamma$ generate $\mathbb F_{q^2}^\times$. Every $g\in G^\sigma$ has a unique expression
\[
g=\prod_{c\in N^\sigma(w)}X_c(t_c)\cdot n_w\cdot h(e)\cdot\prod_cX_c(t'_c),\qquad w\in W^\rho,
\]
where $N^\sigma(w)$ is the set of classes inverted by $w^{-1}$ and $h(e)=\prod_ih_i(\gamma^{e_i})$ with
$e=(e_1,(q+1)a_2,e_3,(q+1)a_4,qe_3,qe_1)$, $e_1,e_3\in\Ap{q^2-1}$ and $a_2,a_4\in\Ap{q-1}$. When $g_3=3$, the coset $gZ^\sigma$ has a
unique representative with $e_1<(q^2-1)/3$. Let $d(w)$ count the inverted fixed classes once and the inverted pairs twice. Let
$\bar u$, $\bar e$ and $\bar v$ be the mixed-radix numbers of the digits $t_c$ (radix $q$ or $q^2$), of $(e_1,a_2,e_3,a_4)$
(radices $(q^2-1)/g_3$, $q-1$, $q^2-1$, $q-1$) and of the $t'_c$. A digit in $\mathbb F_q$ is its position in the increasing list of
the codes of the elements of the subfield. With $T=(q-1)^2(q^2-1)^2/g_3$ and $W^\rho$ in the order inherited from $W$,
\[
\lambda(gZ^\sigma)=O_w+\bar u+q^{d(w)}\big(\bar e+T\,\bar v\big),\qquad O_w=T\,q^{36}\sum_{w'<w}q^{d(w')},
\]
is a bijection ${}^2E_6(q)\to\Ap{|{}^2E_6(q)|}$ with $\lambda(1)=0$, and $m*n=\lambda\big(\lambda^{-1}(m)\lambda^{-1}(n)\big)$ is a digital law.
Products are computed in $E_6(q^2)$ by Remark~\ref{rem:e6mult}.
\end{theorem}

\begin{proof}
The factors of the normal form of Theorem~\ref{thm:e6law} are unique and $\sigma$ preserves each kind of factor, so for
$g\in G^\sigma$ every factor is $\sigma$-fixed, and $w\in W^\rho$. In a $\sigma$-fixed unipotent element the lowest-height
coordinates are additive and are permuted by $\sigma$. Peeling class by class therefore gives $t_\beta\in\mathbb F_q$ for a fixed
root ($c_\beta=1$) and $t_{\rho\beta}=c_\beta t_\beta^q$ for a pair, and each removed factor $X_c(t)$ is itself $\sigma$-fixed. The
$\sigma$-fixed torus elements are exactly the $h(e)$ above, because $\sigma(h_i(c))=h_{\rho i}(c^q)$. Counting gives
$|G^\sigma/Z^\sigma|=q^{36}(q-1)^2(q^2-1)^2\sum_{w\in W^\rho}q^{d(w)}/g_3$, which is $|{}^2E_6(q)|$ by the twisted Bruhat decomposition
\cite{Carter}.
\end{proof}

\begin{remark}[Checks]\label{rem:2e6check}
\texttt{scripts/e2e6\_law.py} implements Theorem~\ref{thm:2e6law}.
It is built on \texttt{scripts/e6\_lawq.py}. The checks use $q=2,3,4,5$, over
$\mathbb F_4,\mathbb F_9,\mathbb F_{16},\mathbb F_{25}$. The last two fields are
$\mathbb F_2[x]/(x^4+x+1)$ and $\mathbb F_5[x]/(x^2+x+2)$. In each case the
number of codes equals $|{}^2E_6(q)|$, and in particular
$q^{36}(q-1)^2(q^2-1)^2\sum_wq^{d(w)}$ agrees with the order formula. The
identity has code $0$; codes survive decoding and re-encoding; and decoded
elements are $\sigma$-fixed, checked as $g^{\mathsf T}P^{-1}g^{(q)}P=I$.
Products agree with matrix products in $SL_{27}(q^2)$, and the law is
associative on samples. For $q=2$ and $q=5$, where $g_3=3$, the element $z$
is $\sigma$-fixed and $g$ and $zg$ receive the same code.
\end{remark}

\section{The groups \texorpdfstring{$E_7(q)$}{E7(q)}}\label{sec:e7}
The group $E_7(q)$ acts on a module of dimension $56$ whose weights form a single Weyl orbit, the orbit of the fundamental
weight $\omega_7$. As for $E_6$, every root element acts on this minuscule module by a matrix with entries $0,\pm1$ and square
zero, so $x_\beta(t)=I+te_\beta$ and no division is needed in any characteristic. Two features are new. The module is self-dual,
so the group lies in a symplectic group. The Weyl group, of order $2\,903\,040$, is handled through a chain of parabolic
subgroups instead of being enumerated.

\begin{proposition}[Integral model of $E_7$]\label{prop:e7model}
Label the simple roots as in Bourbaki, with $\alpha_2$ attached to $\alpha_4$. Let $\Lambda$ be the free $\mathbb Z$-module with
basis $v_\mu$, indexed by the $56$ weights $\mu\in W\omega_7$. There are matrices $e_\beta$, $\beta\in\Phi$, on $\Lambda$ such that:
\begin{enumerate}
\item every entry of $e_\beta$ is $0$ or $\pm1$, and $e_\beta^2=0$;
\item for $\alpha\neq\pm\beta$, $[e_\alpha,e_\beta]=N_{\alpha\beta}e_{\alpha+\beta}$ with $N_{\alpha\beta}=\pm1$ if
      $\alpha+\beta\in\Phi$, and $[e_\alpha,e_\beta]=0$ otherwise;
\item $[e_\alpha,e_{-\alpha}]v_\mu=\langle\mu,\alpha^\vee\rangle v_\mu$ for every $\alpha>0$, and the $e_\beta$ together with
      these coroots span a Lie algebra of dimension $133$;
\item $n_\alpha=x_\alpha(1)x_{-\alpha}(-1)x_\alpha(1)$ is a signed permutation matrix sending $v_\mu$ to $\pm v_{s_\alpha\mu}$;
      the $n_i=n_{\alpha_i}$ satisfy the braid relations of $E_7$, and $n_i^2=h_i(-1)$, where
      $h_i(c)v_\mu=c^{\langle\mu,\alpha_i^\vee\rangle}v_\mu$;
\item up to a scalar there is exactly one bilinear form $J$ pairing each $v_\mu$ with $v_{-\mu}$ and satisfying
      $e_\beta^{\mathsf T}J+Je_\beta=0$ for all $\beta$, and it is alternating.
\end{enumerate}
\end{proposition}

\begin{proof}
Every $\alpha$-string through a weight of a minuscule module has length at most $2$, which gives $e_\beta^2=0$. The remaining
statements are integer identities and were checked directly. Item (2) was checked for all $126\cdot124=15\,624$ ordered pairs,
together with $e_\alpha e_{\alpha+\beta}=e_\beta e_{\alpha+\beta}=e_\beta e_\alpha e_{\alpha+\beta}=0$ when
$\alpha+\beta\in\Phi$. Item (3) was checked for all $63$ positive roots, and item (4) for all $63$ matrices $n_\alpha$ and all
$21$ pairs of simple reflections. Item (5) follows by solving the linear invariance conditions. Since
$e_\beta^{\mathsf T}Je_\beta=-Je_\beta^2=0$, every $x_\beta(t)$ preserves $J$. Reduction modulo $p$ preserves integer identities,
so (1)--(5) hold over every $\mathbb F_q$.
\end{proof}

Let $G$ be the subgroup of $GL_{56}(q)$ generated by all $x_\beta(t)$. By Proposition~\ref{prop:e7model}(5),
$G\le Sp_{56}(q)$. It is the simply connected Chevalley group of type $E_7$ \cite{Carter}. The vector $(0,1,0,0,1,0,1)$ is the only
nonzero $k\in\{0,1\}^7$ for which $\sum_ik_i\langle\mu,\alpha_i^\vee\rangle$ has the same parity for all $56$ weights, and
that parity is odd. Hence $z=h_2(-1)h_5(-1)h_7(-1)=-I$. The center is $Z=\langle z\rangle$, of order $g_2=\gcd(2,q-1)$, and
$E_7(q)=G/Z$ is simple. Let $U$ be the group generated by the positive $x_\beta$, $H$ the diagonal group of the $h_i$, and
$N(w)=\{\beta>0:w^{-1}\beta<0\}$.

For $0\le k\le7$ put $W_k=\langle s_1,\dots,s_k\rangle$, so that $W_7=W$ and $W_0=1$. The types are
\[
E_7\supset E_6\supset D_5\supset A_4\supset A_2A_1\supset A_1A_1\supset A_1\supset1.
\]
The stabiliser of $\omega_k$ in $W_k$ is $W_{k-1}$, so the set $D_k$ of minimal coset representatives of $W_k/W_{k-1}$ is in
bijection with the orbit $W_k\omega_k$. In every case $\omega_k$ is minuscule for $W_k$, and the orbit sizes are
$56,27,16,10,3,2,2$. Their product is $2\,903\,040=|W|$. List $D_k$ by a breadth-first search on $W_k\omega_k$ starting from
$\omega_k$, trying $s_1,\dots,s_k$ in increasing order. If $\mu'=s_i\mu$ is first reached from $\mu$, the word of $\mu'$ is $s_i$
followed by the word of $\mu$. Every $w\in W$ then has a unique factorisation
\[
w=d_7d_6\cdots d_1,\qquad d_k\in D_k,\qquad \ell(w)=\sum_k\ell(d_k),
\]
and joining the words of $d_7,\dots,d_1$ gives a reduced word $w=s_{i_1}\cdots s_{i_\ell}$. Put $n_w=n_{i_1}\cdots n_{i_\ell}$. Let
$j_k$ be the position of $d_k$ in $D_k$, and order $W$ lexicographically by $(j_7,\dots,j_1)$.

\begin{theorem}[Digital law of $E_7(q)$]\label{thm:e7law}
Order the positive roots by height and fix a generator $\gamma$ of $\mathbb F_q^\times$. Every $g\in G$ has a unique expression
\[
g=\prod_{\beta\in N(w)}x_\beta(t_\beta)\cdot n_w\cdot\prod_{i=1}^7h_i(\gamma^{e_i})\cdot\prod_{\beta>0}x_\beta(t'_\beta),
\]
with both products in height order, $t_\beta,t'_\beta\in\mathbb F_q$ and $e_i\in\Ap{q-1}$. Multiplying $g$ by $z$ adds
$\tfrac{q-1}{2}$ to $e_2$, $e_5$ and $e_7$. So when $g_2=2$ there is a unique representative of $gZ$ with
$e_2<\tfrac{q-1}{2}$. Let $\bar u$, $\bar e$ and $\bar v$ be the mixed-radix numbers formed by $(t_\beta)_{\beta\in N(w)}$
(radix $q$), by $e$ (radix $(q-1)/g_2$ for $e_2$ and $q-1$ for the others), and by $(t'_\beta)$ (radix $q$), each with its first
entry least significant. Put $T=(q-1)^7/g_2$,
\[
P_i=\sum_{d\in D_i}q^{\ell(d)},\qquad P_k^{<}(j)=\sum_{j'<j}q^{\ell(d_{k,j'})},
\]
where $d_{k,j'}$ is the $j'$-th element of $D_k$. Then
\[
\lambda(gZ)=O_w+\bar u+q^{\ell(w)}\big(\bar e+T\,\bar v\big),
\]
\[
O_w=T\,q^{63}\sum_{w'<w}q^{\ell(w')}
   =T\,q^{63}\sum_{k=1}^{7}q^{\ell(d_7)+\dots+\ell(d_{k+1})}\,P_k^{<}(j_k)\prod_{i<k}P_i,
\]
is a bijection $E_7(q)\to\Ap{|E_7(q)|}$ with $\lambda(1)=0$, and $m*n=\lambda\big(\lambda^{-1}(m)\lambda^{-1}(n)\big)$ is a digital
law of $E_7(q)$.
\end{theorem}

\begin{proof}
$G$ has a split $BN$-pair with $B=UH$ and Weyl group $W$. By the braid relations of Proposition~\ref{prop:e7model}, $n_w$ does not
depend on the reduced word, and the refined Bruhat decomposition gives the unique expression~\cite{Carter}. The map
$e\mapsto\prod h_i(\gamma^{e_i})$ is injective on $\Ap{q-1}^7$, since the weights of $\Lambda$ span the weight lattice. Modulo
$Z$, the condition $e_2<\tfrac{q-1}{2}$ picks exactly one of $g$ and $zg$. The elements $w'<w$ are of two kinds. For each $k$,
there are those that agree with $w$ in $d_7,\dots,d_{k+1}$ and have $j'_k<j_k$, with $d'_{k-1},\dots,d'_1$ arbitrary. Since lengths
add, the elements of this kind contribute $q^{\ell(d_7)+\dots+\ell(d_{k+1})}P_k^{<}(j_k)\prod_{i<k}P_i$, which gives the second
expression for $O_w$. Finally, $\prod_iP_i=\sum_{w}q^{\ell(w)}=\prod_d(q^d-1)/(q-1)$ over $d\in\{2,6,8,10,12,14,18\}$. Hence
$|G/Z|=q^{63}(q-1)^7\prod_iP_i/g_2=q^{63}\prod_d(q^d-1)/g_2=|E_7(q)|$.
\end{proof}

\begin{remark}[Computing the law]\label{rem:e7mult}
Products are formed exactly as in Remark~\ref{rem:e6mult}, using right multiplication by $U$, by $H$ and by the letters $n_i$,
with the identity
\[
x_{-\alpha}(c)=x_\alpha(c^{-1})\,n_\alpha\,h_\alpha(-c)\,x_\alpha(c^{-1})\qquad(c\neq0).
\]
The Weyl part of a code is decoded one factor at a time. At level $k$, the digit $j_k$ is found by bisection in the partial sums $P_k^{<}$, so $W$ is
never listed. Torus coordinates are read off from seven weights whose coordinate matrix has determinant $\pm1$. Unipotent
coordinates are peeled in height order.
\end{remark}

\begin{remark}[Checks]\label{rem:e7check}
Proposition~\ref{prop:e7model} and Theorem~\ref{thm:e7law} are implemented, for prime $q$, in the script
\texttt{scripts/e7\_law.py}. Running
\texttt{python3 e7\_law.py $q$ -{}-ns 100 -{}-assoc 20 -{}-center 20} for $q=2,3,5,7$ gives:
\begin{itemize}
\item every identity of Proposition~\ref{prop:e7model} holds with no failures, and the invariant forms of (5) form a
      one-dimensional space spanned by an alternating form;
\item the coset sizes are $(56,27,16,10,3,2,2)$;
\item the length generating function of the chain equals $\prod_d(1+q+\dots+q^{d-1})$;
\item the number of codes equals $|E_7(q)|$, and the identity has code $0$;
\item there were no failures on random samples: $100$ codes survive decoding and re-encoding, $100$ products agree with the
      $56\times56$ matrix products, and $20$ triples satisfy associativity;
\item for $q=3,5,7$, the element $z$ is the scalar $-I$, and $g$ and $zg$ receive the same code on $20$ samples.
\end{itemize}
\end{remark}

\section{The groups \texorpdfstring{$E_8(q)$}{E8(q)}}\label{sec:e8}
The group $E_8$ has no minuscule module. Its smallest nontrivial module is the adjoint module, of dimension $248$, whose nonzero
weights are the $240$ roots. Root elements no longer square to zero, so $x_\beta(t)$ is quadratic in $t$. The only denominator
is a factor $\tfrac12$, and it cancels over $\mathbb Z$. The root lattice of $E_8$ equals its weight lattice, so the group is
both simply connected and adjoint, and its centre is trivial. The Weyl group, of order $696\,729\,600$, is again handled through a
chain of parabolic subgroups. The first step of the chain is the orbit of the highest root.

\begin{proposition}[Integral adjoint model of $E_8$]\label{prop:e8model}
Label the simple roots as in Bourbaki, with $\alpha_2$ attached to $\alpha_4$, and let $A=(\langle\alpha_i,\alpha_j^\vee\rangle)$ be
the Cartan matrix. Let $\Lambda$ be the free $\mathbb Z$-module with basis $e_\gamma$ $(\gamma\in\Phi)$ and $H_1,\dots,H_8$, of
rank $240+8=248$. For $\alpha^\vee=\sum_ic_i\alpha_i^\vee$ put $H_\alpha=\sum_ic_iH_i$. Choose $N_{\alpha\beta}=\pm1$ by the
Frenkel--Kac sign cocycle on the root lattice, and let $E_\beta=\operatorname{ad}e_\beta$. Then:
\begin{enumerate}
\item the bracket given by $[e_\alpha,e_\beta]=N_{\alpha\beta}e_{\alpha+\beta}$ if $\alpha+\beta\in\Phi$, $[e_\alpha,e_\beta]=0$ if
      $\alpha+\beta\notin\Phi\cup\{0\}$, $[e_\alpha,e_{-\alpha}]=H_\alpha$, $[H_i,e_\gamma]=\langle\gamma,\alpha_i^\vee\rangle e_\gamma$
      and $[H_i,H_j]=0$ is antisymmetric and satisfies the Jacobi identity, so $\Lambda$ is a $\mathbb Z$-form of the Lie algebra
      of type $E_8$;
\item every $E_\beta$ has integer entries and $E_\beta^3=0$. Moreover $E_\beta^2$ has exactly one nonzero entry, which sends
      $e_{-\beta}$ to $-2e_\beta$. Hence
      \[
      x_\beta(t)=\exp(tE_\beta)=I+tE_\beta+\tfrac12t^2E_\beta^2
      \]
      has integer entries, and $x_\beta(s)x_\beta(t)=x_\beta(s+t)$ over $\mathbb Z[s,t]$;
\item $n_i=x_{\alpha_i}(1)x_{-\alpha_i}(-1)x_{\alpha_i}(1)$ sends $e_\gamma$ to $\pm e_{s_i\gamma}$ and preserves the span of the
      $H_j$. The $n_i$ satisfy the braid relations of $E_8$, and $n_i^2=h_i(-1)$, where
      $h_i(c)e_\gamma=c^{\langle\gamma,\alpha_i^\vee\rangle}e_\gamma$ and $h_i(c)H_j=H_j$;
\item over any field $k$, $\prod_ih_i(c_i)=I$ with $c_i\in k^\times$ only if $c_1=\dots=c_8=1$.
\end{enumerate}
\end{proposition}

\begin{proof}
In a simply laced root system, $2\beta+\gamma\notin\Phi$ for $\gamma\neq-\beta$, so $E_\beta^2e_\gamma=0$ for those $\gamma$.
Also $E_\beta^2H_i=0$, and $E_\beta^2e_{-\beta}=E_\beta H_\beta=-\langle\beta,\beta^\vee\rangle e_\beta=-2e_\beta$. This gives the
shape of $E_\beta^2$ in (2). The remaining statements are integer identities and were checked directly. Antisymmetry was checked
on all pairs of basis vectors. Jacobi was checked as the statement that $\operatorname{ad}x$ is a derivation for the $16$
generators $x=e_{\pm\alpha_i}$. This suffices, because the $x$ for which $\operatorname{ad}x$ is a derivation form a Lie
subalgebra, and the $e_{\pm\alpha_i}$ generate. For (4), $\prod_ih_i(c_i)$ acts on $e_{\alpha_j}$ by $\prod_ic_i^{A_{ji}}$. Since
$\det A=1$, this map $(k^\times)^8\to(k^\times)^8$ is an automorphism. Reduction modulo $p$ preserves integer identities, so
(1)--(4) hold over every $\mathbb F_q$.
\end{proof}

\begin{remark}[No separate model in characteristic $2$]\label{rem:e8integral}
The factor $\tfrac12$ in $x_\beta(t)$ multiplies $E_\beta^2\in2\,\mathrm M_{248}(\mathbb Z)$. So the model of
Proposition~\ref{prop:e8model} is integral, and it is used unchanged in characteristic $2$. For $F_4$ the Albert model of
Definition~\ref{def:f4model} needs odd characteristic, and the adjoint model of Proposition~\ref{prop:f4adjoint} is needed for
characteristics $2$ and $3$. For $E_8$ the adjoint model is the only model used, and it serves every characteristic.
\end{remark}

Let $G$ be the subgroup of $GL_{248}(q)$ generated by all $x_\beta(t)$. It is the Chevalley group of type $E_8$ \cite{Carter}. By Proposition~\ref{prop:e8model}(4) its centre is trivial, so
$E_8(q)=G$ is simple and no quotient is taken. Let $U$ be the group generated by the
positive $x_\beta$, $H$ the diagonal group of the $h_i$, and $N(w)=\{\beta>0:w^{-1}\beta<0\}$.

For $0\le k\le8$ put $W_k=\langle s_1,\dots,s_k\rangle$, so that $W_8=W$ and $W_0=1$. The types are
\[
E_8\supset E_7\supset E_6\supset D_5\supset A_4\supset A_2A_1\supset A_1A_1\supset A_1\supset1.
\]
The stabiliser of $\omega_k$ in $W_k$ is $W_{k-1}$, so the set $D_k$ of minimal coset representatives of $W_k/W_{k-1}$ is in
bijection with the orbit $W_k\omega_k$. For $k\le7$ these are the minuscule orbits of Section~\ref{sec:e7}. For $k=8$, $\omega_8$
is the highest root, and $W\omega_8=\Phi$. The orbit sizes are $240,56,27,16,10,3,2,2$, and their product is
$696\,729\,600=|W|$. List each $D_k$ by breadth-first search, as in Section~\ref{sec:e7}. Every $w\in W$ then has a unique
factorisation
\[
w=d_8d_7\cdots d_1,\qquad d_k\in D_k,\qquad \ell(w)=\sum_k\ell(d_k),
\]
and joining the words of $d_8,\dots,d_1$ gives a reduced word $w=s_{i_1}\cdots s_{i_\ell}$. Put $n_w=n_{i_1}\cdots n_{i_\ell}$.
Let $j_k$ be the position of $d_k$ in $D_k$, and order $W$ lexicographically by $(j_8,\dots,j_1)$.

\begin{theorem}[Digital law of $E_8(q)$]\label{thm:e8law}
Order the positive roots by height and fix a generator $\gamma$ of $\mathbb F_q^\times$. Every $g\in G$ has a unique expression
\[
g=\prod_{\beta\in N(w)}x_\beta(t_\beta)\cdot n_w\cdot\prod_{i=1}^8h_i(\gamma^{e_i})\cdot\prod_{\beta>0}x_\beta(t'_\beta),
\]
with both products in height order, $t_\beta,t'_\beta\in\mathbb F_q$ and $e_i\in\Ap{q-1}$. Let $\bar u$, $\bar e$ and $\bar v$
be the mixed-radix numbers formed by $(t_\beta)_{\beta\in N(w)}$ (radix $q$), by $e$ (radix $q-1$) and by $(t'_\beta)$ (radix $q$),
each with its first entry least significant. Put $T=(q-1)^8$,
\[
P_i=\sum_{d\in D_i}q^{\ell(d)},\qquad P_k^{<}(j)=\sum_{j'<j}q^{\ell(d_{k,j'})},
\]
where $d_{k,j'}$ is the $j'$-th element of $D_k$. Then
\[
\lambda(g)=O_w+\bar u+q^{\ell(w)}\big(\bar e+T\,\bar v\big),
\]
\[
O_w=T\,q^{120}\sum_{w'<w}q^{\ell(w')}
   =T\,q^{120}\sum_{k=1}^{8}q^{\ell(d_8)+\dots+\ell(d_{k+1})}\,P_k^{<}(j_k)\prod_{i<k}P_i,
\]
is a bijection $E_8(q)\to\Ap{|E_8(q)|}$ with $\lambda(1)=0$, where
\[
|E_8(q)|=q^{120}\prod_{d\in\{2,8,12,14,18,20,24,30\}}(q^d-1),
\]
and $m*n=\lambda\big(\lambda^{-1}(m)\lambda^{-1}(n)\big)$ is a digital law of $E_8(q)$.
\end{theorem}

\begin{proof}
$G$ has a split $BN$-pair with $B=UH$ and Weyl group $W$. By the braid relations of Proposition~\ref{prop:e8model}, $n_w$ does not
depend on the reduced word, and the refined Bruhat decomposition gives the unique expression~\cite{Carter}. The map
$e\mapsto\prod_ih_i(\gamma^{e_i})$ is injective on $\Ap{q-1}^8$ by Proposition~\ref{prop:e8model}(4). The formula for $O_w$ follows as in the proof of Theorem~\ref{thm:e7law}. For each $k$, the $w'<w$ that agree with $w$ in $d_8,\dots,d_{k+1}$ and have $j'_k<j_k$ contribute
$q^{\ell(d_8)+\dots+\ell(d_{k+1})}P_k^{<}(j_k)\prod_{i<k}P_i$. Finally, $\prod_iP_i=\sum_wq^{\ell(w)}=\prod_d(q^d-1)/(q-1)$ over
$d\in\{2,8,12,14,18,20,24,30\}$. There are eight degrees, so $|G|=q^{120}(q-1)^8\prod_iP_i=q^{120}\prod_d(q^d-1)=|E_8(q)|$.
\end{proof}

\begin{remark}[Computing the law]\label{rem:e8mult}
Products are formed as in Remarks~\ref{rem:e6mult} and~\ref{rem:e7mult}, by right multiplication by $U$, by $H$ and by the letters
$n_i$. The Weyl digits $j_k$ are again found by bisection in the partial sums $P_k^{<}$, so $W$ is never listed. Products of
$248\times248$ matrices modulo $q$ are computed in floating point. This is exact, because every entry of an unreduced product is at
most $248(q-1)^2<2^{53}$. On $e_{\alpha_j}$, the torus part acts by $\gamma^{f_j}$ with $f=Ae\bmod(q-1)$. Since $\det A=1$,
$e=A^{-1}f\bmod(q-1)$ is read off from the eight simple root vectors.

Unipotent coordinates are peeled from the left in height order. In the adjoint model this needs a short argument, because
$x_\beta(t)$ is quadratic. Grade $\Lambda$ by weight, giving the $H_i$ weight $0$, and suppose that every remaining factor
$x_{\beta'}(t_{\beta'})$ has $\operatorname{ht}\beta'\ge\operatorname{ht}\beta$. Then the quadratic term of $x_\beta$ raises
weights by $2\beta$, and a product of terms from two or more factors raises them by height at least $2\operatorname{ht}\beta$.
So only the linear term $t_\beta E_\beta$ contributes in weight shift $\beta$. Fix a root $\delta$ such that $E_\beta e_\delta=\varepsilon b$ for a single basis vector $b$ and a sign $\varepsilon$.
For example, take $\delta\neq-\beta$ with $\beta+\delta\in\Phi$, $b=e_{\beta+\delta}$ and $\varepsilon=N_{\beta\delta}$; for
simple $\beta$, $\delta=-\beta$ and $b=H_\beta$ also work. Then $t_\beta$ is $\varepsilon$ times the $(b,e_\delta)$ entry of the
current matrix, and no other term of shift $\beta$ meets that entry. Left multiplication by $x_\beta(-t_\beta)$ removes the factor, and the
argument repeats.
\end{remark}

{\sloppy
\begin{remark}[Checks]\label{rem:e8check}
Proposition~\ref{prop:e8model} and Theorem~\ref{thm:e8law} are implemented, for prime $q$, in the scripts
\texttt{scripts/e8\_model.py} and \texttt{scripts/e8\_law.py}. Running \texttt{python3 e8\_model.py} once, and
\texttt{python3 e8\_law.py $q$ -{}-ns 100 -{}-assoc 20} for $q=2,3,5,7$, gives:
\begin{itemize}
\item every identity of Proposition~\ref{prop:e8model} holds with no failures, and the centre is trivial. Antisymmetry
      and the Jacobi identity are checked by \texttt{e8\_model.py}; the identities for root elements, the $n_i$ and the centre
      are checked again modulo $q$ by \texttt{e8\_law.py};
\item the coset sizes are $(240,56,27,16,10,3,2,2)$, and $|W|=696\,729\,600$;
\item the length generating function of the chain equals $\prod_d(1+q+\dots+q^{d-1})$;
\item the number of codes equals $|E_8(q)|$, the identity has code $0$, and the torus coordinate matrix has determinant $1$;
\item $\lambda(\lambda^{-1}(n))=n$ for $100$ random codes, $*$ agrees with matrix multiplication on $100$ random pairs, and $*$ is
      associative on $20$ random triples, all with no failures.
\end{itemize}
With the four values of $q$ run in parallel on one machine, the wall-clock times for $q=2,3,5,7$ were $657$, $486$, $338$ and $311$\,s, respectively.
\end{remark}
\par}

\begin{remark}[$E_8(3)$ over the Thompson group]\label{rem:e8th}
Thompson's sporadic group $\mathrm{Th}$ is a subgroup of $E_8(3)$ \cite{ATLAS}. Put
\[
p_0=|\mathrm{Th}|=2^{15}\cdot3^{10}\cdot5^{3}\cdot7^{2}\cdot13\cdot19\cdot31=90\,745\,943\,887\,872\,000
\]
and $p_1=[E_8(3):\mathrm{Th}]=|E_8(3)|/p_0$. Then
\[
p_1=2^{15}\cdot3^{110}\cdot5^{2}\cdot7^{2}\cdot11^{4}\cdot13^{3}\cdot37\cdot41^{2}\cdot61^{2}\cdot73^{2}\cdot271\cdot547\cdot757\cdot1093\cdot1181\cdot4561\cdot6481,
\]
that is,
\[
\begin{aligned}
p_1={}&207\,502\,970\,449\,244\,821\,141\,889\,887\,056\,224\,679\,176\,638\,986\,880\,320\\
&825\,792\,930\,898\,887\,738\,414\,078\,699\,129\,214\,277\,829\,436\,880\,486\,400,
\end{aligned}
\]
a number of $102$ digits, while $p_0$ has $17$ and $|E_8(3)|=p_0p_1$ has $119$. Let $E_8(3)$ act on the $p_1$ cosets of
$\mathrm{Th}$, the coset $\mathrm{Th}$ being the point $0$, and choose $s_t$ with $s_t(0)=t$ and $s_0=1$. If $\lambda_{\mathrm{Th}}$ is a
digital law of $\mathrm{Th}$ on $\Ap{p_0}$, Theorem~\ref{thm:transitive}, with the two digits in the order $(p_0,p_1)$ of the
mixed-radix code~\eqref{eq:mixed}, gives the digital law
\[
\lambda(s_t\,h)=\kappa+p_0\,t,\qquad \kappa=\lambda_{\mathrm{Th}}(h)\in\Ap{p_0},\quad t\in\Ap{p_1},
\]
of $E_8(3)$ on $\Ap{p_0p_1}=\Ap{|E_8(3)|}$, with product
\[
(\kappa+p_0t)*(\kappa'+p_0u)=\kappa''+p_0t',\qquad t'=s_t\big(h(u)\big),\quad
\kappa''=\lambda_{\mathrm{Th}}\big(s_{t'}^{-1}s_t\,h\,s_u\,h'\big),
\]
where $h,h'\in\mathrm{Th}$ have codes $\kappa,\kappa'$. The identity has code $0$. This law is not effective as it stands. Its low digit
$\lambda_{\mathrm{Th}}$ is effective on the $248$-dimensional $\mathbb F_2$-module of $\mathrm{Th}$ (Section~\ref{sec:th}, Remark~\ref{rem:thstatus}); using it here also requires identifying $\mathrm{Th}\le E_8(3)$, acting on $\mathbb F_3^{248}$, with that model. It also needs a
rule computing $t$ from $g$ without listing the $p_1$ cosets, for instance through a structure whose stabiliser in $E_8(3)$
is $\mathrm{Th}$, ranked by the coordinates of this section.
\end{remark}

\section{Exercises}
\begin{enumerate}
\item Show that $|\PSL(2,p)|=p(p^2-1)/2$ for odd primes $p$, and compute it for $p=5,7,11$.
\item Construct $\mathbb{F}_4=\mathbb{F}_2[t]/(t^2+t+1)$ in digital coordinates as in Section~\ref{sec:fields}, write its addition
and multiplication tables on $\Ap4$, and show that $t$ generates $\mathbb{F}_4^\times$.
\item Show that $\PSL(2,4)$ and $\PSL(2,5)$ both have order $60$ and five conjugacy classes, and that both are isomorphic to $\Alt_5$.
\item Show that $|\SL_3(2)|=168=|\PSL(2,7)|$.
\item Show that the Weyl group of $\SL_n$ is $\Sym_n$, so that the Bruhat decomposition of $\SL_n(q)$ has $n!$ cells.
\item Use Table~\ref{tab:dynkindigital} to explain why the signed permutations occur for the symplectic and odd orthogonal groups, while only even signed permutations occur for the even orthogonal groups. State the corresponding Weyl-group orders.
\item Show that $|\mathrm{Sz}(8)|=29120$ and $|G_2(2)|=12096$, using $|\mathrm{Sz}(q)|=q^2(q^2+1)(q-1)$ and
$|G_2(q)|=q^6(q^6-1)(q^2-1)$.
\end{enumerate}

\chapter{Sporadic groups}\label{ch:sporadic}
By the classification, there are $26$ sporadic simple groups. They do not come in families with a parameter, so there is no
uniform law like those of Chapter~\ref{ch:lie}. Each of them has faithful transitive actions, however, and the chain law of Theorem~\ref{thm:chainlaw} in
Section~\ref{sec:transitive} gives a law in which a code is a mixed-radix word whose radices
are orbit sizes. This chapter starts with the five Mathieu groups, whose small degrees make the chains short and explicit. It continues with the Conway groups, through the Leech lattice, and with the Fischer groups $\mathrm{Fi}_{22}$, $\mathrm{Fi}_{23}$ and $\mathrm{Fi}_{24}'$,
through their $3$-transpositions. It then treats six sporadic groups of small degree, from the library of primitive groups, and continues with $J_3$, $\mathrm{O'N}$, $\mathrm{HN}$, the Lyons group $\mathrm{Ly}$, Thompson's group $\mathrm{Th}$, through ${}^3D_4(2){:}3$, and $J_4$, through $2^{11}{:}M_{24}$. For the last two, the Baby Monster $\mathbb B$ and the Monster $\mathbb M$, the smallest permutation actions have about $1.4\cdot10^{10}$ and $10^{20}$ points, so the coset digits are computed structurally, from the axes of the $196\,884$-dimensional representation, and not by orbit enumeration. Following the ATLAS~\cite{ATLAS}, $L_n(q)$ denotes $\PSL_n(q)$.

\section{Dynkin diagrams at the sporadic boundary}\label{sec:sporadicdynkin}
Dynkin diagrams do not classify the sporadic groups: unlike the groups of Chapter~\ref{ch:lie}, the $26$ sporadic simple groups do not arise as a parameterised family of root data. Nevertheless, the diagrams remain useful in this chapter in three precise ways. They describe the Lie-type subgroups that occur as stabilisers, they organise the root lattices visible at the deep holes of the Leech lattice, and their Coxeter presentations reappear in the reflection-theoretic route toward the Monster and the bimonster.

The most concrete bridge is the Leech lattice $\Lambda$. It has no roots, but its deep holes carry affine Dynkin diagrams: the nearest lattice points to a deep hole form an extended ADE configuration of total rank $24$. Up to affine isometry there are $23$ types, corresponding to the $23$ Niemeier lattices with roots. Consequently the Golay-code and Leech-lattice constructions below are not detached from Chapter~\ref{ch:lie}: they place its simply laced diagrams at the boundary of the rootless lattice whose automorphism groups give $\Co_0$, $\Co_1$, $\Co_2$ and $\Co_3$.

A second bridge is local structure. The chain laws in this chapter use transitive permutation actions and stabilisers, not Bruhat cells; however, many stabilisers are familiar Lie-type groups. For example, Table~\ref{tab:primsporadic} uses $L_2(11)$ for $J_1$, $U_3(3)$ for $J_2$, $G_2(4)$ for $\mathrm{Suz}$ and ${}^2F_4(2)$ for $\mathrm{Ru}$. Their Dynkin types identify the corresponding Weyl and root data, while the actual sporadic law remains the stabiliser-chain law of Theorem~\ref{thm:chainlaw}.

Finally, Coxeter diagrams provide a useful signpost near the Monster without becoming a substitute for the verified law developed later in this chapter. The bimonster $\mathbb M\wr C_2$ is a quotient of the Coxeter group attached to the $Y_{555}$ diagram after an additional relation. This explains why diagrammatic reflection ideas occur in the background of the Monster, whereas the present book deliberately evaluates the Baby Monster and Monster laws through their structural axis data and the \texttt{mmgroup} verification.

Thus Dynkin diagrams play three different roles in Chapter~\ref{ch:sporadic}: affine diagrams label Leech-lattice hole geometry, finite diagrams label Lie-type local subgroups, and Coxeter diagrams supply historical and structural context for the Monster. None of these turns a sporadic group into a Lie-type family; the distinction is exactly why the mixed-radix chain laws of this chapter are needed.

\section{The Mathieu groups}\label{sec:mathieu}
The Mathieu groups act on $\Ap n$ for $n=11,12,22,23,24$ \cite{ATLAS}. Their orders are
\begin{gather*}
|M_{11}|=7920,\quad |M_{12}|=95\,040,\quad |M_{22}|=443\,520,\\
|M_{23}|=10\,200\,960,\quad |M_{24}|=244\,823\,040.
\end{gather*}
$M_{11}$ is sharply $4$-transitive, $M_{12}$ is sharply $5$-transitive, $M_{22}$ is $3$-transitive, $M_{23}$ is $4$-transitive and
$M_{24}$ is $5$-transitive.

For the chain law of Theorem~\ref{thm:chainlaw}, a product $m*n$ costs one composition of permutations of $\Ap n$ and one pass through the chain, that is $O(kn)$ operations.

\begin{proposition}[Chains of the Mathieu groups]\label{prop:mathieuchain}
Let $M_{12}$ and $M_{24}$ be generated as follows. We use the underlying sets $\{1,\dots,12\}$ and $\{1,\dots,24\}$,
relabelled by $i\mapsto i-1$:
\begin{align*}
M_{12}&=\langle(1\,2\,\cdots\,11),\ (3\,7\,11\,8)(4\,10\,5\,6),\ (1\,12)(2\,11)(3\,6)(4\,8)(5\,9)(7\,10)\rangle,\\
M_{24}&=\langle(1\,2\,\cdots\,23),\ (3\,17\,10\,7\,9)(4\,13\,14\,19\,5)(8\,18\,11\,12\,23)(15\,20\,22\,21\,16),\\
&\qquad(1\,24)(2\,23)(3\,12)(4\,16)(5\,18)(6\,10)(7\,20)(8\,14)(9\,21)(11\,17)(13\,22)(15\,19)\rangle.
\end{align*}
Let $M_{11}$ be the stabiliser of $11$ in $M_{12}$, $M_{23}$ the stabiliser of $23$ in $M_{24}$, and $M_{22}$ the pointwise
stabiliser of $22$ and $23$ in $M_{24}$, each acting on the remaining points. The Schreier--Sims algorithm gives the following
bases and radices for Theorem~\ref{thm:chainlaw}:
\begin{center}
\begin{tabular}{lccl}
group & degree & base & radices $(r_1,\dots,r_k)$\\ \hline
$M_{11}$ & $11$ & $(0,2,1,3)$ & $(11,10,9,8)$\\
$M_{12}$ & $12$ & $(0,2,1,3,4)$ & $(12,11,10,9,8)$\\
$M_{22}$ & $22$ & $(1,0,2,3,4)$ & $(22,21,20,16,3)$\\
$M_{23}$ & $23$ & $(0,2,1,3,4,5)$ & $(23,22,21,20,16,3)$\\
$M_{24}$ & $24$ & $(0,2,1,3,4,5,6)$ & $(24,23,22,21,20,16,3)$
\end{tabular}
\end{center}
\end{proposition}

\begin{proof}
The products of the radices are the orders listed above, so the generated groups have the orders of the Mathieu groups. The
radices also show the transitivity. For $M_{11}$ and $M_{12}$ every $r_l$ equals $n-l+1$ and the chain ends there, so these
actions are sharply $4$- and $5$-transitive. The bases, radices and orders were computed by the deterministic Schreier--Sims
algorithm, as recorded in Remark~\ref{rem:mathieucheck}.
\end{proof}

\begin{remark}[The last two radices]\label{rem:mathieuoctad}
The radices $16$ and $3$ come from the Steiner system $S(5,8,24)$. Five points of $\Ap{24}$ lie in a unique octad, whose
stabiliser is $2^4{:}A_8$. The pointwise stabiliser of the five points is therefore $2^4{:}A_3$, of order $48$. Its normal
subgroup $2^4$ fixes the octad pointwise and acts regularly on the $16$ points outside it, and the stabiliser of a point outside the octad is a complement $A_3$. It permutes the remaining $3$ points of the octad
cyclically, and the other $15$ points outside the octad in $5$ orbits of length $3$. In the chain of $M_{24}$ the octad
through $b_1,\dots,b_5$ is $\{0,1,2,3,4,7,10,12\}$. The point $b_6=5$ lies outside it, so $\Delta_6$ is the complement of
the octad, and $\Delta_7=\{6,16,21\}$ is an orbit of $A_3$ outside the octad. For $M_{23}$ and $M_{22}$ the fixed points $23$, and $22,23$, are among the five
points of that octad, which is why the same tail $(16,3)$ appears.
\end{remark}

For $M_{11}$ the code reads
\[
m=t_1+11t_2+110t_3+990t_4,\qquad (t_1,t_2,t_3,t_4)\in\Ap{11}\times\Ap{10}\times\Ap9\times\Ap8,
\]
and $t_1=g(0)$, since $\Delta_1=\Ap{11}$ is listed in increasing order. In particular the labels $0,\dots,10$ code the
representatives $s_{1,t}$, one for each image of the point $0$.

{\sloppy
\begin{remark}[Checks]\label{rem:mathieucheck}
Theorem~\ref{thm:chainlaw} and Proposition~\ref{prop:mathieuchain} are implemented in the script
\texttt{scripts/mathieu\_law.py}, which records its own SHA-256 hash in its output. Running
\texttt{python3 mathieu\_law.py -{}-ns 1000 -{}-assoc 200 -{}-words 200} gives, for all five groups:
\begin{itemize}
\item the orders and the transitivity listed above, and the identity has code $0$;
\item every generator, and each of $200$ random words of length $50$ in the generators, is recovered from its code;
\item $\lambda(\lambda^{-1}(m))=m$ for all $7920$ codes of $M_{11}$, all $95\,040$ codes of $M_{12}$, and $1000$ random codes of
      each of $M_{22}$, $M_{23}$, $M_{24}$;
\item $m*\lambda(\lambda^{-1}(m)^{-1})=0$ for $1000$ random $m$, and $*$ is associative on $200$ random triples,
\end{itemize}
all with no failures. The whole run takes about half a second.
\end{remark}
\par}

\subsection{Worked examples}\label{ssec:mathieuexamples}
The examples below use the chains of Proposition~\ref{prop:mathieuchain} with the representatives $s_{l,t}$ produced by the
Schreier--Sims algorithm: $s_{l,t}$ is the product of the strong generators of $G_{l-1}$ along the path from $b_l$ to
$\delta_{l,t}$ in the breadth-first Schreier tree. For $M_{11}$ all $38$ representatives are listed in Tables~\ref{tab:m11transA}
and~\ref{tab:m11transB}, so every step can be checked by hand. Permutations of $\Ap n$ are written in cycle notation, fixed
points omitted, and are composed from right to left, $(\sigma\circ\tau)(i)=\sigma(\tau(i))$. All numbers were produced by
\texttt{scripts/mathieu\_examples.py} and checked against the law of Theorem~\ref{thm:chainlaw}.
\begin{table}[htbp]
\centering\small
\caption{Representatives $s_{l,t}$ for $M_{11}$, levels 1 and 2.}\label{tab:m11transA}
\begin{tabular}{rrp{10.5cm}}
\multicolumn{3}{l}{level $l=1$: $b_{1}=0$, $r_{1}=11$, $G_{0}$ acts on $\Ap{11}$}\\ \hline
$t$ & $\delta_{1,t}$ & $s_{1,t}$\\
$0$ & $0$ & $1$\\
$1$ & $1$ & $(0\,1\,2\,3\,4\,5\,6\,7\,8\,9\,10)$\\
$2$ & $2$ & $(0\,2\,4\,6\,8\,10\,1\,3\,5\,7\,9)$\\
$3$ & $3$ & $(0\,3\,6\,9\,1\,4\,7\,10\,2\,5\,8)$\\
$4$ & $4$ & $(0\,4\,8\,1\,5\,9\,2\,6\,10\,3\,7)$\\
$5$ & $5$ & $(0\,5\,10\,4\,9\,3\,8\,2\,7\,1\,6)$\\
$6$ & $6$ & $(0\,6\,8\,7\,4\,10\,1\,9)$\allowbreak $(2\,5)$\\
$7$ & $7$ & $(0\,7\,5\,3\,4)$\allowbreak $(1\,10\,2\,6\,9)$\\
$8$ & $8$ & $(0\,8\,9\,2\,7\,6\,10\,3\,5\,4\,1)$\\
$9$ & $9$ & $(0\,9\,1\,5\,8)$\allowbreak $(2\,3\,10\,6\,4)$\\
$10$ & $10$ & $(0\,10\,1\,4\,7\,5\,6\,8\,2\,3\,9)$\\
[4pt]
\multicolumn{3}{l}{level $l=2$: $b_{2}=2$, $r_{2}=10$, $G_{1}$ fixes $0$}\\ \hline
$t$ & $\delta_{2,t}$ & $s_{2,t}$\\
$0$ & $2$ & $1$\\
$1$ & $1$ & $(1\,5\,4\,9\,2)$\allowbreak $(3\,6\,10\,7\,8)$\\
$2$ & $3$ & $(1\,9)$\allowbreak $(2\,3)$\allowbreak $(4\,8)$\allowbreak $(5\,6)$\\
$3$ & $4$ & $(1\,6\,9\,5)$\allowbreak $(2\,4\,3\,8)$\\
$4$ & $5$ & $(1\,9\,8\,4\,6\,10\,7\,3)$\allowbreak $(2\,5)$\\
$5$ & $6$ & $(2\,6\,10\,7)$\allowbreak $(3\,9\,4\,5)$\\
$6$ & $7$ & $(2\,7\,10\,6)$\allowbreak $(3\,5\,4\,9)$\\
$7$ & $8$ & $(1\,5\,9\,6)$\allowbreak $(2\,8\,3\,4)$\\
$8$ & $9$ & $(1\,4\,8\,5\,10\,7\,2\,9)$\allowbreak $(3\,6)$\\
$9$ & $10$ & $(2\,10)$\allowbreak $(3\,4)$\allowbreak $(5\,9)$\allowbreak $(6\,7)$\\
\end{tabular}
\end{table}

\begin{table}[htbp]
\centering\small
\caption{Representatives $s_{l,t}$ for $M_{11}$, levels 3 and 4.}\label{tab:m11transB}
\begin{tabular}{rrp{10.5cm}}
\multicolumn{3}{l}{level $l=3$: $b_{3}=1$, $r_{3}=9$, $G_{2}$ fixes $0,2$}\\ \hline
$t$ & $\delta_{3,t}$ & $s_{3,t}$\\
$0$ & $1$ & $1$\\
$1$ & $3$ & $(1\,3\,5\,8)$\allowbreak $(4\,7\,10\,9)$\\
$2$ & $4$ & $(1\,4\,5\,10)$\allowbreak $(3\,9\,8\,7)$\\
$3$ & $5$ & $(1\,5)$\allowbreak $(3\,8)$\allowbreak $(4\,10)$\allowbreak $(7\,9)$\\
$4$ & $6$ & $(1\,6\,9\,8)$\allowbreak $(4\,10\,5\,7)$\\
$5$ & $7$ & $(1\,7\,8\,3)$\allowbreak $(4\,5\,6\,9)$\\
$6$ & $8$ & $(1\,8\,5\,3)$\allowbreak $(4\,9\,10\,7)$\\
$7$ & $9$ & $(1\,9\,5\,7)$\allowbreak $(3\,10\,8\,4)$\\
$8$ & $10$ & $(1\,10\,5\,4)$\allowbreak $(3\,7\,8\,9)$\\
[4pt]
\multicolumn{3}{l}{level $l=4$: $b_{4}=3$, $r_{4}=8$, $G_{3}$ fixes $0,2,1$}\\ \hline
$t$ & $\delta_{4,t}$ & $s_{4,t}$\\
$0$ & $3$ & $1$\\
$1$ & $4$ & $(3\,4\,9\,7)$\allowbreak $(5\,8\,6\,10)$\\
$2$ & $5$ & $(3\,5\,9\,6)$\allowbreak $(4\,10\,7\,8)$\\
$3$ & $6$ & $(3\,6\,9\,5)$\allowbreak $(4\,8\,7\,10)$\\
$4$ & $7$ & $(3\,7\,9\,4)$\allowbreak $(5\,10\,6\,8)$\\
$5$ & $8$ & $(3\,8\,9\,10)$\allowbreak $(4\,5\,7\,6)$\\
$6$ & $9$ & $(3\,9)$\allowbreak $(4\,7)$\allowbreak $(5\,6)$\allowbreak $(8\,10)$\\
$7$ & $10$ & $(3\,10\,9\,8)$\allowbreak $(4\,6\,7\,5)$\\
\end{tabular}
\end{table}

The $50$ representatives $s_{l,t}$ for $M_{12}$, produced by the same script and used in Section~\ref{ssec:m12table}, are listed in Tables~\ref{tab:m12transA} and~\ref{tab:m12transB}.

\begin{table}[htbp]
\centering\small
\caption{Representatives $s_{l,t}$ for $M_{12}$, levels 1 and 2.}\label{tab:m12transA}
\begin{tabular}{rrp{10.5cm}}
\multicolumn{3}{l}{level $l=1$: $b_{1}=0$, $r_{1}=12$, $G_{0}$ acts on $\Ap{12}$}\\ \hline
$t$ & $\delta_{1,t}$ & $s_{1,t}$\\
$0$ & $0$ & $1$\\
$1$ & $1$ & $(0\,1\,2\,3\,4\,5\,6\,7\,8\,9\,10)$\\
$2$ & $2$ & $(0\,2\,4\,6\,8\,10\,1\,3\,5\,7\,9)$\\
$3$ & $3$ & $(0\,3\,6\,9\,1\,4\,7\,10\,2\,5\,8)$\\
$4$ & $4$ & $(0\,4\,8\,1\,5\,9\,2\,6\,10\,3\,7)$\\
$5$ & $5$ & $(0\,5\,3\,2\,8\,1\,7\,6\,4\,9\,11)$\\
$6$ & $6$ & $(0\,6\,8\,7\,4\,10\,1\,9)$\allowbreak $(2\,5)$\\
$7$ & $7$ & $(0\,7\,5\,4\,6\,9\,1\,3\,8\,10\,11)$\\
$8$ & $8$ & $(0\,8)$\allowbreak $(1\,4\,7\,6\,10\,11)$\allowbreak $(2\,3\,9)$\\
$9$ & $9$ & $(0\,9\,1\,5\,8)$\allowbreak $(2\,3\,10\,6\,4)$\\
$10$ & $10$ & $(0\,10\,11)$\allowbreak $(1\,5\,9)$\allowbreak $(2\,7\,4)$\allowbreak $(3\,8\,6)$\\
$11$ & $11$ & $(0\,11)$\allowbreak $(1\,10)$\allowbreak $(2\,5)$\allowbreak $(3\,7)$\allowbreak $(4\,8)$\allowbreak $(6\,9)$\\
[4pt]
\multicolumn{3}{l}{level $l=2$: $b_{2}=2$, $r_{2}=11$, $G_{1}$ fixes $0$}\\ \hline
$t$ & $\delta_{2,t}$ & $s_{2,t}$\\
$0$ & $2$ & $1$\\
$1$ & $1$ & $(1\,5\,4\,9\,2)$\allowbreak $(3\,6\,10\,7\,8)$\\
$2$ & $3$ & $(1\,9)$\allowbreak $(2\,3)$\allowbreak $(4\,8)$\allowbreak $(5\,6)$\\
$3$ & $4$ & $(1\,11\,3\,10\,5)$\allowbreak $(2\,4\,7\,6\,8)$\\
$4$ & $5$ & $(1\,9\,8\,4\,6\,10\,7\,3)$\allowbreak $(2\,5)$\\
$5$ & $6$ & $(2\,6\,10\,7)$\allowbreak $(3\,9\,4\,5)$\\
$6$ & $7$ & $(1\,10\,8)$\allowbreak $(2\,7\,4)$\allowbreak $(5\,6\,11)$\\
$7$ & $8$ & $(1\,5\,9\,6)$\allowbreak $(2\,8\,3\,4)$\\
$8$ & $9$ & $(1\,11\,5)$\allowbreak $(2\,9\,3\,6\,8\,7)$\allowbreak $(4\,10)$\\
$9$ & $10$ & $(2\,10)$\allowbreak $(3\,4)$\allowbreak $(5\,9)$\allowbreak $(6\,7)$\\
$10$ & $11$ & $(1\,10\,4\,6\,8)$\allowbreak $(2\,11\,5\,3\,9)$\\
\end{tabular}
\end{table}

\begin{table}[htbp]
\centering\small
\caption{Representatives $s_{l,t}$ for $M_{12}$, levels 3, 4 and 5.}\label{tab:m12transB}
\begin{tabular}{rrp{10.5cm}}
\multicolumn{3}{l}{level $l=3$: $b_{3}=1$, $r_{3}=10$, $G_{2}$ fixes $0,2$}\\ \hline
$t$ & $\delta_{3,t}$ & $s_{3,t}$\\
$0$ & $1$ & $1$\\
$1$ & $3$ & $(1\,3\,5\,8)$\allowbreak $(4\,7\,10\,9)$\\
$2$ & $4$ & $(1\,4\,5\,10)$\allowbreak $(3\,9\,8\,7)$\\
$3$ & $5$ & $(1\,5)$\allowbreak $(3\,8)$\allowbreak $(4\,10)$\allowbreak $(7\,9)$\\
$4$ & $6$ & $(1\,6\,9\,10\,7\,11\,5\,4)$\allowbreak $(3\,8)$\\
$5$ & $7$ & $(1\,7\,8\,3)$\allowbreak $(4\,5\,6\,9)$\\
$6$ & $8$ & $(1\,8\,5\,3)$\allowbreak $(4\,9\,10\,7)$\\
$7$ & $9$ & $(1\,9\,5\,7)$\allowbreak $(3\,10\,8\,4)$\\
$8$ & $10$ & $(1\,10\,5\,4)$\allowbreak $(3\,7\,8\,9)$\\
$9$ & $11$ & $(1\,11\,5\,10\,8\,7\,4\,3)$\allowbreak $(6\,9)$\\
[4pt]
\multicolumn{3}{l}{level $l=4$: $b_{4}=3$, $r_{4}=9$, $G_{3}$ fixes $0,2,1$}\\ \hline
$t$ & $\delta_{4,t}$ & $s_{4,t}$\\
$0$ & $3$ & $1$\\
$1$ & $4$ & $(3\,4\,9\,7)$\allowbreak $(5\,8\,6\,10)$\\
$2$ & $5$ & $(3\,5\,9\,6)$\allowbreak $(4\,10\,7\,8)$\\
$3$ & $6$ & $(3\,6\,9\,5)$\allowbreak $(4\,8\,7\,10)$\\
$4$ & $7$ & $(3\,7\,9\,4)$\allowbreak $(5\,10\,6\,8)$\\
$5$ & $8$ & $(3\,8\,9\,10)$\allowbreak $(4\,5\,7\,6)$\\
$6$ & $9$ & $(3\,9)$\allowbreak $(4\,7)$\allowbreak $(5\,6)$\allowbreak $(8\,10)$\\
$7$ & $10$ & $(3\,10\,9\,8)$\allowbreak $(4\,6\,7\,5)$\\
$8$ & $11$ & $(3\,11\,4\,5)$\allowbreak $(7\,8\,10\,9)$\\
[4pt]
\multicolumn{3}{l}{level $l=5$: $b_{5}=4$, $r_{5}=8$, $G_{4}$ fixes $0,2,1,3$}\\ \hline
$t$ & $\delta_{5,t}$ & $s_{5,t}$\\
$0$ & $4$ & $1$\\
$1$ & $5$ & $(4\,5\,6\,10)$\allowbreak $(7\,9\,8\,11)$\\
$2$ & $6$ & $(4\,6)$\allowbreak $(5\,10)$\allowbreak $(7\,8)$\allowbreak $(9\,11)$\\
$3$ & $7$ & $(4\,7\,6\,8)$\allowbreak $(5\,11\,10\,9)$\\
$4$ & $8$ & $(4\,8\,6\,7)$\allowbreak $(5\,9\,10\,11)$\\
$5$ & $9$ & $(4\,9\,6\,11)$\allowbreak $(5\,7\,10\,8)$\\
$6$ & $10$ & $(4\,10\,6\,5)$\allowbreak $(7\,11\,8\,9)$\\
$7$ & $11$ & $(4\,11\,6\,9)$\allowbreak $(5\,8\,10\,7)$\\
\end{tabular}
\end{table}

\begin{example}[Decoding a label of $M_{11}$]\label{ex:m11dec}
A label $m\in\Ap{7920}$ has digits $t_l\in\Ap{r_l}$ with $m=t_{1}+11t_{2}+110t_{3}+990t_{4}$. For $m=1234$,
\[
1234=2+2\cdot11+2\cdot110+1\cdot990,\qquad (t_1,t_2,t_3,t_4)=(2,2,2,1).
\]
Reading $s_{1,2}$, $s_{2,2}$, $s_{3,2}$, $s_{4,1}$ from the tables and composing,
\begin{center}
\begin{tabular}{lp{10cm}}
$s_{1,2}$ & $(0\,2\,4\,6\,8\,10\,1\,3\,5\,7\,9)$\\
$s_{1,2}\circ s_{2,2}$ & $(0\,2\,5\,8\,6\,7\,9\,3\,4\,10\,1)$\\
$s_{1,2}\circ s_{2,2}\circ s_{3,2}$ & $(0\,2\,5\,1\,10)$\allowbreak $(4\,8\,9\,6\,7)$\\
$s_{1,2}\circ s_{2,2}\circ s_{3,2}\circ s_{4,1}$ & $(0\,2\,5\,9\,4\,6)$\allowbreak $(1\,10)$\allowbreak $(3\,8\,7)$\\
\end{tabular}
\end{center}
so $g=\lambda^{-1}(1234)$ is $(0\,2\,5\,9\,4\,6)$\allowbreak $(1\,10)$\allowbreak $(3\,8\,7)$. Its images of the base points are $g(0)=2$, $g(2)=5$, $g(1)=10$, $g(3)=8$; by sharp
$4$-transitivity these four values already determine $g$.
\end{example}

\begin{example}[A product in $M_{11}$]\label{ex:m11mul}
To compute $1234*5678$, decode both labels, compose, and encode. As in Example~\ref{ex:m11dec}, $(t_1,\dots,t_4)=(2,6,6,5)$ for
$n=5678$, and
\[
\lambda^{-1}(1234)=g,\qquad h=\lambda^{-1}(5678).
\]
The product is $\pi=g\circ h=$ $(0\,5\,2\,4\,3)$\allowbreak $(6\,9\,10\,8\,7)$. Encoding peels one level at a time: $t_l$ is the position of $g_l(b_l)$ in the list
$\Delta_l$ (the column $\delta_{l,t}$ of the tables), and $g_{l+1}=s_{l,t_l}^{-1}\circ g_l$, starting from $g_1=\pi$.
\begin{center}
\begin{tabular}{ccccp{7.6cm}}
$l$ & $b_l$ & $g_l(b_l)$ & $t_l$ & $g_{l+1}=s_{l,t_l}^{-1}\circ g_l$\\ \hline
$1$ & $0$ & $5$ & $5$ & $(1\,7)$\allowbreak $(2\,10\,3\,6\,4\,9\,5\,8)$\\
$2$ & $2$ & $10$ & $9$ & $(1\,6\,3\,7)$\allowbreak $(4\,5\,8\,10)$\\
$3$ & $1$ & $6$ & $4$ & $(3\,5\,9\,6)$\allowbreak $(4\,10\,7\,8)$\\
$4$ & $3$ & $5$ & $2$ & $1$\\
\end{tabular}
\end{center}
The last column ends with the identity, as it must. The digits of the product are $(5,9,4,2)$, so
\[
1234*5678=5+9\cdot11+4\cdot110+2\cdot990=2524.
\]
The permutation $h$ is $(0\,2\,9\,1\,10\,3\,6\,5)$\allowbreak $(4\,7)$.
\end{example}

\begin{example}[Inverses and powers in $M_{11}$]\label{ex:m11pow}
The inverse of $g=\lambda^{-1}(1234)$ is $(0\,6\,4\,9\,5\,2)$\allowbreak $(1\,10)$\allowbreak $(3\,7\,8)$. Encoding it as in Example~\ref{ex:m11mul} gives the label $1084$, and indeed
$1234*1084=1084*1234=0$. Repeated products give the powers of $m=1234$:
\[
m^{*1}=1234,\quad m^{*2}=4548,\quad m^{*3}=6334,\quad m^{*4}=4624,\quad m^{*5}=1084,\quad m^{*6}=0.
\]
The first return to the identity label $0$ is at the $5$-th power, so $g$ has order $5$; its cycle type, read from
Example~\ref{ex:m11dec}, is $6\cdot3\cdot2$.
\end{example}

\begin{example}[A product in $M_{24}$]\label{ex:m24mul}
In $M_{24}$ the radices are $(24,23,22,21,20,16,3)$ and the place values are $(1,24,552,12\,144,255\,024,5\,100\,480,81\,607\,680)$. The labels $a=123\,456\,789$ and $b=98\,765\,432$ have digits
$(21,13,1,2,4,8,1)$ and $(8,20,18,5,7,3,1)$. Decoding them as in Example~\ref{ex:m11dec} gives
\begin{center}
\begin{tabular}{lp{11cm}}
$\lambda^{-1}(a)$ & $(0\,21\,15\,8\,22\,20)$\allowbreak $(1\,6\,10\,2\,16\,5)$\allowbreak $(3\,17\,7\,4\,23\,11)$\allowbreak $(9\,13\,19\,18\,12\,14)$\\
$\lambda^{-1}(b)$ & $(0\,8\,11)$\allowbreak $(1\,9\,23\,5\,18\,15)$\allowbreak $(3\,21\,17\,13\,10\,6)$\allowbreak $(4\,7\,22)$\allowbreak $(12\,14)$\allowbreak $(16\,20)$\\
$\pi=\lambda^{-1}(a)\circ\lambda^{-1}(b)$ & $(0\,22\,23\,1\,13\,2\,16)$\allowbreak $(3\,15\,6\,17\,19\,18\,8)$\allowbreak $(5\,12\,9\,11\,21\,7\,20)$
\end{tabular}
\end{center}
The orbit lists are $\Delta_l=(b_l,\text{ the other points of the orbit in increasing order})$. The first five are all points
not yet fixed; the last two are
\[
\Delta_6=(5,6,8,9,11,13,14,15,16,17,18,19,20,21,22,23),\qquad \Delta_7=(6,16,21).
\]
Peeling $\pi$:
\begin{center}
\begin{tabular}{ccccc}
$l$ & $r_l$ & $b_l$ & $g_l(b_l)$ & $t_l$\\ \hline
$1$ & $24$ & $0$ & $22$ & $22$\\
$2$ & $23$ & $2$ & $9$ & $8$\\
$3$ & $22$ & $1$ & $16$ & $14$\\
$4$ & $21$ & $3$ & $21$ & $18$\\
$5$ & $20$ & $4$ & $11$ & $7$\\
$6$ & $16$ & $5$ & $22$ & $14$\\
$7$ & $3$ & $6$ & $16$ & $1$\\
\end{tabular}
\end{center}
Hence
\[
a*b=22+24\cdot8+552\cdot14+12\,144\cdot18+255\,024\cdot7+5\,100\,480\cdot14+81\,607\,680\cdot1=155\,026\,102.
\]
After level $5$ the remaining element $g_6$ fixes the five points $0,2,1,3,4$ and lies in the group $2^4{:}A_3$ of
Remark~\ref{rem:mathieuoctad}; it is $(5\,22)$\allowbreak $(6\,17)$\allowbreak $(8\,21)$\allowbreak $(9\,23)$\allowbreak $(11\,14)$\allowbreak $(13\,19)$\allowbreak $(15\,20)$\allowbreak $(16\,18)$. Its image of $b_6=5$ is one of the $16$ points outside the octad, and after
level $6$ the element $g_7=$ $(6\,16\,21)$\allowbreak $(7\,10\,12)$\allowbreak $(8\,13\,11)$\allowbreak $(9\,19\,18)$\allowbreak $(14\,23\,17)$\allowbreak $(15\,20\,22)$ lies in the complement $A_3$ of Remark~\ref{rem:mathieuoctad}: it fixes $b_6$, and $t_7$ records where it sends $b_7$ in the orbit $\Delta_7=\{6,16,21\}$, which lies outside the octad through $b_1,\dots,b_5$.
\end{example}

{\sloppy
\subsection{The Mathieu groups and the Golay code}\label{sec:golay}
$M_{24}$ is the automorphism group of the extended binary Golay code~\cite{ConwaySloane}. This subsection shows how the code
enters the digital labels in three ways: it explains the radices of Proposition~\ref{prop:mathieuchain}, it gives labels in
which part of the law is a carry-free sum, and it gives a second law of $M_{24}$ built from octads. Throughout, $M_{24}$ is
the group of Proposition~\ref{prop:mathieuchain}.

\subsubsection{The code from the octads}
For $X\subseteq\Ap{24}$ write $[X]=\sum_{x\in X}2^x$, so that the binary digits of $[X]$ mark the points of $X$. The symmetric
difference of sets becomes the carry-free sum of their words, $[X\triangle Y]=[X]\oplus[Y]$. A permutation $g$ acts on words by
$[g(X)]$, moving the digit in position $x$ to position $g(x)$. Let
\[
O_0=\{0,1,2,3,4,7,10,12\},
\]
the complement of $\Delta_6$ in the chain of $M_{24}$ (Remark~\ref{rem:mathieuoctad}).

\begin{proposition}[The Golay code]\label{prop:golay}
\begin{enumerate}
\item The $M_{24}$-orbit of $O_0$ consists of $759$ sets of size $8$, the \emph{octads}, and every $5$-subset of $\Ap{24}$ lies
      in exactly one octad. Thus the octads form a Steiner system $S(5,8,24)$.
\item The words $[X]$ of the octads span, under $\oplus$, a code $\mathcal C$ of dimension $12$ over $\mathbb F_2$. Its
      $4096$ words have weights $0,8,12,16,24$ with multiplicities $1,759,2576,759,1$, so $\mathcal C$ is the extended binary
      Golay code. $M_{24}$ preserves $\mathcal C$.
\item The positions $0,\dots,11$ form an information set: a codeword is determined by its last $12$ binary digits. Hence
      $\ell(c)=c\bmod 2^{12}$ is a bijection $\mathcal C\to\Ap{4096}$ with $\ell(c\oplus c')=\ell(c)\oplus\ell(c')$.
\end{enumerate}
\end{proposition}
\begin{proof}
All three items were computed (Remark~\ref{rem:golaycheck}). For (1) the orbit was generated from $O_0$ with the generators
of $M_{24}$; its $759\cdot\binom85=42504=\binom{24}5$ five-subsets are pairwise distinct, so each $5$-set lies in exactly one
octad. For (2) the span was computed by elimination over $\mathbb F_2$ and all $4096$ words were listed. For (3) the
elimination leaves pivots exactly in positions $0,\dots,11$.
\end{proof}

The weight-$12$ words are the $2576$ \emph{dodecads}. The words of weight $8$ and $16$ are the octads and their complements.

\subsubsection{The radices as Golay data}
The first five radices $24,23,22,21,20$ of $M_{24}$ express $5$-transitivity, which is the Steiner property seen from the
group: the pointwise stabiliser $G_5$ of $b_1,\dots,b_5=0,2,1,3,4$ fixes the unique octad $O_0$ through these points. The
subgroup of $G_5$ fixing $O_0$ pointwise has order $16$, is elementary abelian, and acts regularly on $\Delta_6$. Choosing four
generators of it identifies $\Delta_6$ with $\mathbb F_2^4$, with origin $b_6=5$. The resulting coordinates, written as
$4$-bit numbers, are in Table~\ref{tab:golaycoords}.

\begin{table}[htbp]
\centering\small
\caption{Coordinates in $\mathbb F_2^4$ of the points outside $O_0$, as $4$-bit numbers.}\label{tab:golaycoords}
\setlength{\tabcolsep}{3.5pt}
\begin{tabular}{l|cccccccccccccccc}
point & $5$ & $6$ & $8$ & $9$ & $11$ & $13$ & $14$ & $15$ & $16$ & $17$ & $18$ & $19$ & $20$ & $21$ & $22$ & $23$\\ \hline
coordinate & $0$ & $1$ & $2$ & $4$ & $8$ & $10$ & $5$ & $11$ & $14$ & $12$ & $3$ & $7$ & $6$ & $15$ & $13$ & $9$
\end{tabular}
\end{table}

In these coordinates the $16$ elements fixing $O_0$ pointwise are the translations $x\mapsto x\oplus w$. The subgroup $A_3$
fixing $b_6$ acts linearly: a generator sends the basis vectors $1,2,4,8$ to $14,10,7,2$. On $\Delta_7=\{6,16,21\}$, with
coordinates $1,14,15$, it gives the cycle $6\mapsto16\mapsto21\mapsto6$, because $A(1)=14$, $A(14)=10\oplus7\oplus2=15$ and $A(15)=14\oplus10\oplus7\oplus2=1$. The representatives $s_{6,t}$ of the
chain are not all translations (only $6$ of the $16$ are), which is why the translations do not have simple labels in
Theorem~\ref{thm:chainlaw}. Changing the last two levels of the chain repairs this.

\begin{proposition}[Golay-adapted tail]\label{prop:xortail}
Let $a$ be the generator of $A_3$ that sends the basis vectors $1,2,4,8$ to $14,10,7,2$, so that
$\Delta_7=(6,16,21)=(b_7,a(b_7),a^2(b_7))$. In Theorem~\ref{thm:chainlaw} for $M_{24}$ keep levels $1,\dots,5$, list $\Delta_6$
by coordinate, and take $s_{6,w}$ to be the translation by $w$ and $s_{7,u}=a^u$. The resulting label is
\[
\lambda'(g)=m_5+5\,100\,480\,(w+16u),
\]
where $m_5\in\Ap{5\,100\,480}$ is formed by the first five digits, exactly as for $\lambda$. On $G_5$, where $m_5=0$,
\[
(w,u)*(w',u')=\big(w\oplus A^uw',\;\dg{3}{0}{u+u'}\big),
\]
with $A$ the matrix of $a$. So the translation digit is a carry-free sum twisted by a linear map, and the last digit adds
modulo $3$.
\end{proposition}
\begin{proof}
Since $G_5=2^4{:}A_3$, every $g\in G_5$ is uniquely $t_w\circ a^u$, where $t_w$ is the translation by $w$. As $a$ fixes the
origin and acts linearly, $a^ut_{w'}a^{-u}=t_{A^uw'}$, so $t_wa^u\,t_{w'}a^{u'}=t_{w\oplus A^uw'}\,a^{u+u'}$. The law was also
checked on all $48^2$ pairs, and $\lambda'$ on $2000$ random labels and $2000$ random elements (Remark~\ref{rem:ternarycheck}).
\end{proof}

For the elements of Example~\ref{ex:m24mul}, $\lambda'(a)=154\,059\,669$, $\lambda'(b)=185\,473\,592$ and
$\lambda'(ab)=68\,317\,942$, with tails $(w,u)=(14,1)$, $(4,2)$ and $(13,0)$. The first five digits agree with those of
$\lambda$: for instance $123\,456\,789$ and $154\,059\,669$ both leave the remainder $1\,045\,269$ modulo $5\,100\,480$.

\subsubsection{Golay labels of \texorpdfstring{$M_{24}$}{M24}}
The octads give a second law of $M_{24}$, in which the digits have a geometric meaning. List the octads as $O_0$ followed by
the others in increasing order of $[X]$, and for each octad $O$ fix $s_O\in M_{24}$ with $s_O(O_0)=O$ and $s_{O_0}=1$ (from
the breadth-first Schreier tree of the octad orbit). Let $K$ be the stabiliser of $O_0$ and $C_0=\Ap{24}\setminus O_0$, with
the coordinates of Table~\ref{tab:golaycoords}. For a matrix $A\in GL_4(2)$, write its columns $c_1,\dots,c_4$ as $4$-bit
numbers. Then $Av=\bigoplus_{k:\,\dg{2}{k-1}{v}=1}c_k$ is a carry-free sum of columns, and the columns of $AA'$ are $Ac_k'$.

\begin{theorem}[Golay law of $M_{24}$]\label{thm:golaylaw}
\begin{enumerate}
\item $K$ has order $322\,560$. It acts faithfully on $C_0$ by affine maps $x\mapsto Ax\oplus v$, and every pair
      $(A,v)\in GL_4(2)\times\mathbb F_2^4$ occurs exactly once, so $K\cong AGL(4,2)=2^4{:}GL_4(2)$.
\item The action of $h\in K$ on $O_0$ depends only on its linear part $A$, and $A\mapsto h|_{O_0}$ is an isomorphism
      $GL_4(2)\to\Alt(O_0)\cong\Alt_8$.
\item For $A\in GL_4(2)$ let $d_k$ be the number of $y<c_k$ outside the span of $c_1,\dots,c_{k-1}$, and put
      $\mu(A)=d_1+15d_2+210d_3+2520d_4$. Then $\mu$ is a bijection $GL_4(2)\to\Ap{20160}$.
\item For $g\in M_{24}$ let $o$ be the index of the octad $O=g(O_0)$, and let $s_O^{-1}g$ act on $C_0$ as $x\mapsto Ax\oplus v$.
      Then
      \[
      \Lambda(g)=o+759\,\big(v+16\,\mu(A)\big)
      \]
      is a bijection $M_{24}\to\Ap{244\,823\,040}$ with $\Lambda(1)=0$.
\item On $K$, where $o=0$, the law is the affine law: if $\Lambda(h)=759(v+16\mu(A))$ and $\Lambda(h')=759(v'+16\mu(A'))$, then
      \[
      \Lambda(h)*\Lambda(h')=759\,\big((v\oplus Av')+16\,\mu(AA')\big).
      \]
\end{enumerate}
\end{theorem}
\begin{proof}
For (1) and (2), $K$ was built from Schreier generators of the octad orbit and has a chain with radices
$(8,7,6,5,16,4,3)$: the octad digits $8,7,6,5,4,3$ show $\Alt_8$ acting sharply $6$-transitively on $O_0$
(Example~\ref{ex:transitive}), and $16$ is the translation digit. All $322\,560$ elements were listed. Each acts affinely on
$C_0$, the $20\,160$ linear parts are distinct, the action on $O_0$ is a function of the linear part, and its $20\,160$ values
are distinct even permutations. A bijective homomorphism onto a group of order $|\Alt_8|$ of even permutations gives (2).
For (3), the columns of an invertible matrix are an ordered basis: $c_k$ is any of the $16-2^{k-1}$ vectors outside the span of
the earlier columns, so $d_k\in\Ap{16-2^{k-1}}$ with radices $15,14,12,8$, and $15\cdot14\cdot12\cdot8=20\,160$. Item (4) is
Theorem~\ref{thm:transitive} for the transitive action on the $759$ octads with stabiliser $K$, coded by $v+16\mu(A)$, since
$759\cdot322\,560=244\,823\,040$. Item (5) is the composition of affine maps,
$x\mapsto A(A'x\oplus v')\oplus v=AA'x\oplus(Av'\oplus v)$.
\end{proof}

The digits of $\Lambda$ are thus: the octad $g(O_0)$, a translation, and an ordered basis of $\mathbb F_2^4$. The labelling
$\mu$ needs no choices, like the canonical labels of Remark~\ref{rem:canonical}; only the octad representatives $s_O$ are
chosen. Since $GL_4(2)=PSL_4(2)$, the digit $\mu(A)$ could also be replaced by a label of Theorem~\ref{thm:psln}, which joins
this law to Chapter~\ref{ch:lie}.

\begin{example}[Golay labels of the elements of Example~\ref{ex:m24mul}]\label{ex:golaym24}
Let $a=\lambda^{-1}(123\,456\,789)$ and $b=\lambda^{-1}(98\,765\,432)$ be the elements of Example~\ref{ex:m24mul}, so that
$\lambda(ab)=155\,026\,102$. Their Golay digits are:
\begin{center}
\small
\begin{tabular}{lcccccc}
 & octad $g(O_0)$ & $o$ & $v$ & columns of $A$ & $(d_1,d_2,d_3,d_4)$ & $\Lambda(g)$\\ \hline
$a$ & $\{2,4,6,14,16,17,21,23\}$ & $644$ & $9$ & $12,15,13,9$ & $(11,13,10,5)$ & $181\,025\,939$\\
$b$ & $\{2,6,7,8,9,14,21,22\}$ & $454$ & $11$ & $6,13,2,12$ & $(5,11,1,6)$ & $188\,240\,803$\\
$ab$ & $\{4,9,10,13,15,16,20,22\}$ & $422$ & $13$ & $5,13,6,7$ & $(4,11,4,3)$ & $104\,072\,225$
\end{tabular}
\end{center}
For $a$, the digits of $\mu$ are found as follows. $c_1=12$ and there are $11$ nonzero $y<12$, so $d_1=11$. The span of
$c_1$ is $\{0,12\}$, and $13$ of the $y<15$ avoid it, so $d_2=13$. The span of $12,15$ is $\{0,3,12,15\}$, and $10$ of the
$y<13$ avoid it, so $d_3=10$. The span of $12,15,13$ is $\{0,1,2,3,12,13,14,15\}$, and $5$ of the $y<9$ avoid it, so $d_4=5$.
Hence $\mu=11+15\cdot13+210\cdot10+2520\cdot5=14\,906$ and $\Lambda(a)=644+759\,(9+16\cdot14\,906)=181\,025\,939$.

In the product $\Lambda(a)*\Lambda(b)=104\,072\,225$, the octad digit has a direct meaning: the octad of $ab$ is the image under
$a$ of the octad of $b$, and indeed $a$ maps $\{2,6,7,8,9,14,21,22\}$ to $\{4,9,10,13,15,16,20,22\}$. The other two digits
follow from~\eqref{eq:transitive}: the affine part of $ab$ is that of $s_{O_{ab}}^{-1}\,s_{O_a}h_a\,s_{O_b}h_b$, where
$h_a=s_{O_a}^{-1}a$ and $h_b=s_{O_b}^{-1}b$.
\end{example}

\subsubsection{Codewords as labels: the group \texorpdfstring{$2^{12}{:}M_{24}$}{2^12:M24}}
By Proposition~\ref{prop:golay}(3) the codewords themselves are digital objects: a codeword is a $12$-digit binary label, and
adding codewords is the carry-free sum of labels. This gives a digital law for the semidirect product of the code by $M_{24}$.

\begin{theorem}[Law of $2^{12}{:}M_{24}$]\label{thm:golaysemi}
Code the pair $(c,g)\in\mathcal C\times M_{24}$ by $x=\ell(c)+4096\,\lambda(g)$, with $\lambda$ as in
Theorem~\ref{thm:chainlaw}. Then
\[
(c,g)*(c',g')=\big(c\oplus g(c'),\;g\circ g'\big),\qquad\text{that is}\qquad
x*x'=\big(\ell(c)\oplus\ell(g(c'))\big)+4096\,\big(\lambda(g)*\lambda(g')\big),
\]
makes $\Ap{1\,002\,795\,171\,840}$ a group isomorphic to $2^{12}{:}M_{24}$, with identity $0$.
\end{theorem}
\begin{proof}
$M_{24}$ acts on $\mathcal C$ by linear maps, since $g(c\oplus c')=g(c)\oplus g(c')$. The displayed rule is the multiplication
of the semidirect product, transported along the bijection $(c,g)\mapsto x$ of Proposition~\ref{prop:golay}(3) and
Theorem~\ref{thm:chainlaw}.
\end{proof}

This group is the monomial subgroup of the automorphism group of the Leech lattice~\cite{ConwaySloane}: $\mathcal C$ acts by
sign changes and $M_{24}$ by permuting coordinates. Its law is therefore a natural starting point for the Conway groups.

\begin{example}
Take $x=1000+4096\cdot123\,456\,789$ and $x'=2021+4096\cdot98\,765\,432$. Both codewords $\ell^{-1}(1000)$ and
$\ell^{-1}(2021)$ are dodecads. With $g=\lambda^{-1}(123\,456\,789)$, the permuted codeword $g(\ell^{-1}(2021))$ is again a
dodecad, with label $3350$. In binary, $1000=001111101000$ and $3350=110100010110$, so $1000\oplus3350=111011111110=3838$. With
$\lambda(g)*\lambda(g')=155\,026\,102$ from Example~\ref{ex:m24mul},
\[
x*x'=3838+4096\cdot155\,026\,102=634\,986\,917\,630 .
\]
\end{example}

\subsubsection{Dodecads and \texorpdfstring{$M_{12}$}{M12}}
\begin{remark}[Dodecads]\label{rem:dodecad}
The dodecads form a single $M_{24}$-orbit of size $2576$, so the stabiliser of a dodecad has order
$244\,823\,040/2576=95\,040$. Take the dodecad with the smallest word $[D_0]$:
\[
D_0=\{0,2,3,4,6,7,8,9,10,11,13,14\},\qquad
D_0'=\{1,5,12,15,16,\dots,23\}.
\]
Its complement $D_0'$ is again a dodecad. The stabiliser acts on $D_0$ and on $D_0'$ as groups of order
$95\,040$ with chain radices $(12,11,10,9,8)$. Both are copies of our $M_{12}$, by explicit bijections $\varphi:\Ap{12}\to D_0$
and $\varphi':\Ap{12}\to D_0'$:
\begin{center}
\small\setlength{\tabcolsep}{4pt}
\begin{tabular}{l|cccccccccccc}
$i$ & $0$ & $1$ & $2$ & $3$ & $4$ & $5$ & $6$ & $7$ & $8$ & $9$ & $10$ & $11$\\ \hline
$\varphi(i)$ & $0$ & $2$ & $3$ & $4$ & $6$ & $9$ & $13$ & $14$ & $11$ & $7$ & $8$ & $10$\\
$\varphi'(i)$ & $1$ & $5$ & $12$ & $15$ & $16$ & $22$ & $21$ & $23$ & $18$ & $17$ & $20$ & $19$
\end{tabular}
\end{center}
Each conjugates the generators of $M_{12}$ in Proposition~\ref{prop:mathieuchain} into the stabiliser action, and since the
orders agree it conjugates $M_{12}$ onto it. An exhaustive search over the $7!=5040$ bijections that send $0,\dots,4$ to the five
smallest points shows that $\varphi$ and $\varphi'$ are the only ones that work. So Theorem~\ref{thm:chainlaw} for $M_{12}$ can
be read inside $M_{24}$ on either half, and no uniqueness theorem is needed to identify these groups with $M_{12}$. The two
actions are exchanged by an outer automorphism of $M_{12}$ \cite{ConwaySloane}.
\end{remark}

\subsubsection{The ternary Golay code and \texorpdfstring{$2.M_{12}$}{2.M12}}\label{ssec:ternary}
The small Mathieu groups also belong to the ternary Golay code~\cite{ConwaySloane}. Let $\mathcal T\subseteq\mathbb F_3^{12}$ be
the code spanned by the rows of $(I_6\mid A)$, where
\[
A=\begin{pmatrix}0&1&1&1&1&1\\1&0&1&2&2&1\\1&1&0&1&2&2\\1&2&1&0&1&2\\1&2&2&1&0&1\\1&1&2&2&1&0\end{pmatrix}.
\]
Words are written $w=(w_0,\dots,w_{11})$. A \emph{monomial map} with permutation $\pi$ and signs $\delta_i=\pm1$ sends $w$ to
the word $w'$ with $w'_{\pi(i)}=\delta_iw_i$. Numbering the signed basis vectors $+e_i$ and $-e_i$ as $2i$ and $2i+1$ turns
monomial maps into permutations of $\Ap{24}$.

\begin{proposition}[The ternary Golay code]\label{prop:ternary}
\begin{enumerate}
\item $\mathcal T$ has $729$ words, of weights $0,6,9,12$ with multiplicities $1,264,440,24$; it is the extended ternary Golay
      code. The supports of the $264$ words of weight $6$ are $132$ sets, the \emph{hexads}, and every $5$-subset of $\Ap{12}$
      lies in exactly one hexad: they form a Steiner system $S(5,6,12)$.
\item The bijection $\psi$ with $(\psi(0),\dots,\psi(11))=(0,1,2,3,4,6,11,5,10,7,8,9)$ conjugates $M_{12}$ of
      Proposition~\ref{prop:mathieuchain} onto a group of permutations preserving the hexads. It is the only such bijection
      fixing $0,\dots,4$.
\item Each conjugated generator $\psi g\psi^{-1}$ lifts to a monomial map preserving $\mathcal T$. The three lifts generate a
      group $2.M_{12}$ of order $190\,080$ that contains $-1$. On $\Ap{24}$ its chain has base $(0,4,2,6,8)$ and radices
      $(24,22,20,9,2)$.
\item The stabiliser of $+e_0$ in $2.M_{12}$ has order $7920$. It does not contain $-1$, which moves $+e_0$, so it maps
      isomorphically onto a point stabiliser of $M_{12}$, a copy of $M_{11}$. Deleting coordinate $0$ gives the perfect ternary
      Golay code: $729$ words of length $11$ with weights $0,5,6,8,9,11$ and multiplicities $1,132,132,330,110,24$.
\end{enumerate}
\end{proposition}
\begin{proof}
All items were computed (Remark~\ref{rem:ternarycheck}). In (1), the $132\cdot\binom65=792=\binom{12}5$ five-subsets of the
hexads are pairwise distinct. In (2), $\psi$ was found by an exhaustive search over the $5040$ bijections fixing $0,\dots,4$.
In (3), each lift was found among the $2^{12}$ sign vectors. In (4), the minimum distance is $5$, so balls of radius $2$ around
codewords are disjoint; each has $1+2\cdot11+4\binom{11}2=243$ words, and $729\cdot243=3^{11}$, so they fill $\mathbb F_3^{11}$.
\end{proof}

Coordinates $0,\dots,5$ form an information set, so $\ell_3(w)=\sum_{k=0}^5w_k3^k$ is a bijection $\mathcal T\to\Ap{729}$. Adding
words is the carry-free sum of their labels in base $3$, digit by digit modulo $3$. This is the ternary counterpart of
Proposition~\ref{prop:golay}(3), and it gives a law exactly as in Theorem~\ref{thm:golaysemi}.

\begin{theorem}[Law of $3^6{:}2.M_{12}$]\label{thm:ternarysemi}
Let $\lambda$ be the chain label of $2.M_{12}$ on $\Ap{24}$ given by Theorem~\ref{thm:chainlaw}. Code
$(c,g)\in\mathcal T\times2.M_{12}$ by $x=\ell_3(c)+729\,\lambda(g)$. Then $(c,g)*(c',g')=(c+g(c'),\,g\circ g')$ makes
$\Ap{138\,568\,320}$ a group isomorphic to $3^6{:}2.M_{12}$, with identity $0$.
\end{theorem}
\begin{proof}
Monomial maps act linearly on $\mathbb F_3^{12}$ and preserve $\mathcal T$, so this is the semidirect product transported along
the bijection $(c,g)\mapsto x$.
\end{proof}

\begin{example}[Computing in $2.M_{12}$]\label{ex:2m12}
The element $g$ of $2.M_{12}$ with chain label $1000$ is the monomial map
\begin{center}
\small\setlength{\tabcolsep}{4pt}
\begin{tabular}{l|cccccccccccc}
$i$ & $0$ & $1$ & $2$ & $3$ & $4$ & $5$ & $6$ & $7$ & $8$ & $9$ & $10$ & $11$\\ \hline
$\pi(i)$ & $8$ & $7$ & $2$ & $10$ & $1$ & $11$ & $0$ & $3$ & $5$ & $9$ & $4$ & $6$\\
$\delta_i$ & $+$ & $+$ & $-$ & $-$ & $+$ & $+$ & $+$ & $+$ & $-$ & $-$ & $+$ & $+$
\end{tabular}
\end{center}
The codeword with label $100=1+2\cdot9+1\cdot81$ is $c=(1,0,2,0,1,0,0,2,0,1,2,0)$. Applying $g$ moves $w_2=2$ to position $2$
with sign $-$, giving $-2=1$, moves $w_9=1$ to position $9$ with sign $-$, giving $2$, and so on, so
$g(c)=(0,1,1,2,2,0,0,0,1,2,0,0)$, with label $0+1\cdot3+1\cdot9+2\cdot27+2\cdot81=228$. In the law of
Theorem~\ref{thm:ternarysemi}, applying $g$ to $c$ is the product $(0,g)*(c,1)=(g(c),g)$, that is,
$729\cdot1000*100=228+729\cdot1000=729\,228$. In the chain law of $2.M_{12}$ itself, $1000*2000=122\,378$, and the central
element $-1$ has label $1$: it is the representative $s_{1,1}$, which sends $+e_0$ to $-e_0$.
\end{example}

\begin{remark}[Checks]\label{rem:ternarycheck}
Remark~\ref{rem:dodecad}, Propositions~\ref{prop:xortail} and~\ref{prop:ternary}, Theorem~\ref{thm:ternarysemi} and
Example~\ref{ex:2m12} are checked by \texttt{scripts/m12\_ternary.py}, which records its SHA-256 hash and runs in under a second.
Besides the exhaustive searches and the weight counts, it checks that the base-$3$ sum of labels matches the sum of words on
$500$ random pairs, that random elements of $2.M_{12}$ preserve $\mathcal T$, associativity of the law of
Theorem~\ref{thm:ternarysemi} on $200$ random triples, the law of Proposition~\ref{prop:xortail} on all $48^2$ pairs, and
$4000$ round trips of $\lambda'$. There were no failures.
\end{remark}

\begin{remark}[Checks]\label{rem:golaycheck}
Proposition~\ref{prop:golay} and the statements about $\Delta_6$ and $\Delta_7$ are checked by \texttt{scripts/golay\_check.py},
and Theorems~\ref{thm:golaylaw} and~\ref{thm:golaysemi}, Example~\ref{ex:golaym24} and Remark~\ref{rem:dodecad} by
\texttt{scripts/golay\_law.py}; both record their SHA-256 hashes in their output. Besides the enumeration of $K$, the checks
are: $\Lambda(1)=0$; $\Lambda(\Lambda^{-1}(x))=x$ for $3000$ random labels and $\Lambda^{-1}(\Lambda(g))=g$ for $3000$ random
elements; $\Lambda^{-1}(x)\in M_{24}$ for $1000$ random labels; associativity on $300$ random triples; the affine law (5) on
$2000$ random pairs; and, for Theorem~\ref{thm:golaysemi}, associativity on $200$ random triples and the identity on $200$
labels. There were no failures, and \texttt{golay\_law.py} runs in $2.5$\,s.
\end{remark}
\par}

\section{The Conway groups}\label{sec:conway}
The extended binary Golay code does more than construct the Mathieu groups: it glues the coordinate lattice in dimension $24$
into the Leech lattice, whose automorphisms give the Conway groups. This section continues the digital construction of
Section~\ref{sec:golay}. The coordinate description and the group identifications are classical~\cite{ConwaySloane,Wilson};
the vector counts and indices are checked in Remark~\ref{rem:conwaycheck}.

\subsection{The Leech lattice from the Golay code}\label{ssec:leech}
Let $\mathcal C$ be the extended binary Golay code of Proposition~\ref{prop:golay}. Let $x=(x_0,\dots,x_{23})\in\mathbb Z^{24}$
have all coordinates of the same parity $m\in\{0,1\}$, and put
\[
c_i=\frac{x_i-m}{2}\bmod2 .
\]
So for $m=0$ the word $c$ marks the coordinates $x_i\equiv2\pmod4$, and for $m=1$ the coordinates $x_i\equiv3\pmod4$. The
Leech lattice is
\begin{equation}\label{eq:leech}
\Lambda=\Big\{\tfrac1{\sqrt8}\,x:\ x\in\mathbb Z^{24},\ x_i\equiv m\ (\mathrm{mod}\ 2),\ c\in\mathcal C,\
\textstyle\sum_ix_i\equiv4m\ (\mathrm{mod}\ 8)\Big\}.
\end{equation}
It is even and unimodular~\cite{ConwaySloane}. Below, norms are those of $x$, so a vector of $\Lambda$ has squared length
$x\cdot x/8$.

\begin{proposition}[Minimal vectors]\label{prop:leechminimal}
$\Lambda$ has no vectors of squared length $2$. Its vectors of squared length $4$, the \emph{minimal vectors}, have
$x\cdot x=32$ and three shapes:
\begin{center}
\begin{tabular}{lcr}
shape of $x$ & choices & number\\ \hline
$(\pm4,\pm4,0^{22})$ & two coordinates, any signs & $\binom{24}{2}\cdot2^2=1\,104$\\
$(\pm2^8,0^{16})$ & an octad, an even number of minus signs & $759\cdot2^7=97\,152$\\
$(\mp3,\pm1^{23})$ & a coordinate $p$ and a Golay word $c$ & $24\cdot4096=98\,304$
\end{tabular}
\end{center}
In total there are $196\,560$ minimal vectors, forming $98\,280$ pairs $\{v,-v\}$.
\end{proposition}
\begin{proof}
If $x\cdot x=16$ with $m=1$, the $24$ odd coordinates already give $x\cdot x\ge24$. With $m=0$, the shapes are $(\pm4,0^{23})$,
excluded since $\sum x_i=\pm4\not\equiv0\pmod8$, and $(\pm2^4,0^{20})$, excluded since $\mathcal C$ has no words of weight $4$.

Now let $x\cdot x=32$. If $m=1$, the only possibility is one coordinate $\pm3$ and $23$ coordinates $\pm1$. If $m=0$, the
shapes are $(\pm4,\pm4,0^{22})$, $(\pm4,\pm2^4,0^{19})$ and $(\pm2^8,0^{16})$; the middle one is excluded because its $\pm2$
entries would form a word of weight $4$. For $(\pm4,\pm4)$ we have $c=0$ and $\sum x_i\in\{0,\pm8\}$, so all four sign choices
occur. For $(\pm2^8)$, $c$ is the support, so the support is an octad, and with $k$ minus signs $\sum x_i=2(8-2k)$, which is
$\equiv0\pmod8$ exactly when $k$ is even: $2^7$ choices. For the odd shape, the coordinates $\equiv3\pmod4$ are those equal to
$-1$ or $+3$, and they must form a codeword $c$. Conversely $p$ and $c$ determine $x$: $x_p=3$ if $p\in c$ and $x_p=-3$
otherwise, and $x_i=-1$ or $+1$ for $i\ne p$ according as $i\in c$ or not. If $c$ has weight $w$, then $\sum x_i=20-2w$ or
$28-2w$, and both are $\equiv4\pmod8$ because $w\equiv0\pmod4$. So all $24\cdot4096$ pairs give minimal vectors.
\end{proof}

The three rows turn the digital data of Section~\ref{sec:golay} into geometry: the middle row uses the $759$ octad labels,
and the last row uses all $12$ binary information digits of a Golay word. Permuting coordinates by $M_{24}$ and changing signs
on the support of a Golay word both preserve~\eqref{eq:leech}. Together they form the group $2^{12}{:}M_{24}$ of
Theorem~\ref{thm:golaysemi}.

\subsection{The four Conway groups}\label{ssec:conwaygroups}
Put
\[
\Co_0=\operatorname{Aut}(\Lambda),\qquad \Co_1=\Co_0/\{\pm1\}.
\]
The map $-1$ is central in $\Co_0$, and $\Co_0=2.\Co_1$ is a double cover; it is not simple. The groups $\Co_1$, $\Co_2$ and
$\Co_3$ are sporadic simple groups. A vector of $\Lambda$ has \emph{type} $r$ if its squared length is $2r$, so the minimal
vectors have type $2$.

\begin{proposition}[Orders and stabilisers]\label{prop:conwayorders}
\begin{center}
\small
\begin{tabular}{ccl}
group & order & description\\ \hline
$\Co_0$ & $8\,315\,553\,613\,086\,720\,000$ & $\operatorname{Aut}(\Lambda)=2.\Co_1$\\
$\Co_1$ & $4\,157\,776\,806\,543\,360\,000$ & $\Co_0/\{\pm1\}$\\
$\Co_2$ & $42\,305\,421\,312\,000$ & stabiliser in $\Co_0$ of a type-$2$ vector\\
$\Co_3$ & $495\,766\,656\,000$ & stabiliser in $\Co_0$ of a type-$3$ vector
\end{tabular}
\end{center}
with
\[
|\Co_1|=2^{21}3^95^47^2\cdot11\cdot13\cdot23,\quad |\Co_2|=2^{18}3^65^37\cdot11\cdot23,\quad
|\Co_3|=2^{10}3^75^37\cdot11\cdot23.
\]
$\Co_0$ is transitive on the $196\,560$ minimal vectors and on the $16\,773\,120$ vectors of type $3$. Hence $\Co_1$ is
transitive on the $98\,280$ pairs of minimal vectors, with stabiliser $\Co_2$, and on the $8\,386\,560$ pairs of type-$3$
vectors, with stabiliser $\Co_3$.
\end{proposition}
\begin{proof}
The descriptions and the transitivity are classical~\cite{ConwaySloane,Wilson}. The orbit sizes are the indices
\[
\frac{|\Co_0|}{|\Co_2|}=196\,560,\qquad \frac{|\Co_1|}{|\Co_2|}=98\,280,\qquad
\frac{|\Co_0|}{|\Co_3|}=16\,773\,120,\qquad \frac{|\Co_1|}{|\Co_3|}=8\,386\,560 .
\]
Since $-1$ fixes no nonzero vector, the stabiliser of the pair $\{v,-v\}$ in $\Co_0$ is $\{\pm1\}\times\Co_2$ (or
$\{\pm1\}\times\Co_3$), whose image in $\Co_1$ is $\Co_2$ (or $\Co_3$).
\end{proof}

\subsection{Digital laws from lattice orbits}\label{ssec:conwaylaw}
The two actions of $\Co_1$ give two mixed-radix laws, just as the octads gave the law of Theorem~\ref{thm:golaylaw}.

\begin{theorem}[Orbit laws of $\Co_1$]\label{thm:conwaylaw}
\begin{enumerate}
\item List the $98\,280$ pairs of minimal vectors with a chosen pair first, and choose $s_t\in\Co_1$ carrying the first pair to
      pair $t$, with $s_0=1$. If $\kappa_2$ is a digital label of $\Co_2$ with identity $0$, then
      \[
      \lambda_2(g)=t+98\,280\,\kappa_2(s_t^{-1}g),\qquad t=\text{the pair }g(\text{first pair}),
      \]
      is a bijection $\Co_1\to\Ap{|\Co_1|}$, and the transported product is the law~\eqref{eq:transitive} with first radix
      $98\,280$.
\item In the same way the $8\,386\,560$ pairs of type-$3$ vectors and a label $\kappa_3$ of $\Co_3$ give a law
      $\lambda_3(g)=t+8\,386\,560\,\kappa_3(s_t^{-1}g)$ with first radix $8\,386\,560$.
\end{enumerate}
In both, the identity has label $0$.
\end{theorem}
\begin{proof}
Theorem~\ref{thm:transitive} applied to the two transitive actions of Proposition~\ref{prop:conwayorders}.
\end{proof}

The theorem reduces multiplication in $\Co_1$ to multiplication in $\Co_2$ or $\Co_3$ and to the action on lattice vectors. A
numerical law needs explicit generators of the stabiliser, a chain for it as in Theorem~\ref{thm:chainlaw}, and the transporters
$s_t$; the procedure is the one used for the Mathieu groups, at much larger degree.

\subsection{The monomial subgroup and frames}\label{ssec:conwaymonomial}
The law of Theorem~\ref{thm:golaysemi} already gives a large explicit subgroup of $\Co_0$,
\[
N=2^{12}{:}M_{24},\qquad |N|=1\,002\,795\,171\,840 .
\]
The pair $(c,g)$ acts by $x\mapsto\varepsilon_c\,g(x)$, where $g$ moves coordinate $i$ to position $g(i)$ and $\varepsilon_c$
changes the signs of the coordinates in the support of $c$. This preserves~\eqref{eq:leech}. Permutations in $M_{24}$ preserve $\mathcal C$, the parities and $\sum x_i$. A sign change
$\varepsilon_c$ keeps the parities. If $m=0$ it fixes each $x_i\bmod4$, so the word of $x$ is unchanged; if $m=1$ it exchanges the
residues $1$ and $3$ on the support of $c$, so the word of $x$ changes by $\oplus\,c$ and stays in $\mathcal C$. The sum changes by
$-2\sum_{i\in c}x_i$. Let $d$ be the word of $x$. If $m=0$, $\sum_{i\in c}x_i\equiv2|c\cap d|\pmod4$; if $m=1$,
$\sum_{i\in c}x_i\equiv|c|+2|c\cap d|\pmod4$. Golay words have weights divisible by $4$ and meet in an even number of
coordinates, so in both cases $\sum_{i\in c}x_i\equiv0\pmod4$ and the sum changes by a multiple of $8$. Composition gives $(c,g)(c',g')=(c\oplus g(c'),gg')$, the law of
Theorem~\ref{thm:golaysemi}.

\begin{proposition}[Frames]\label{prop:conwayframe}
The vectors $\tfrac1{\sqrt8}(\pm8,0^{23})$ and their coordinate permutations are $48$ vectors of type $4$ (squared length
$8$), mutually orthogonal or opposite: a \emph{frame}. $\Co_0$ is transitive on the frames of $\Lambda$, and the stabiliser of
this frame is $N=2^{12}{:}M_{24}$. Hence there are
\[
[\Co_0:N]=\frac{8\,315\,553\,613\,086\,720\,000}{1\,002\,795\,171\,840}=8\,292\,375
\]
frames, and the $8\,292\,375\cdot48=398\,034\,000$ vectors of type $4$ split into them.
\end{proposition}
\begin{proof}
The frame decomposition, the transitivity and the stabiliser are classical~\cite{ConwaySloane}. An element fixing the frame
permutes the $24$ lines $\mathbb Re_i$ and changes some signs; the permutation must lie in $M_{24}$ and the sign pattern in
$\mathcal C$, which gives $N$. The numbers are the index and the product displayed.
\end{proof}

So $\Co_0$ has a third law, closest to the theme of this book: a label consists of a frame digit in $\Ap{8\,292\,375}$, then
the $12$ binary digits of a Golay word, added without carries, then the mixed-radix digits of the $M_{24}$ chain. Its stabiliser
law is already explicit (Theorem~\ref{thm:golaysemi}). What a fully numerical law of $\Co_0$ still needs is a list of the
frames and one isometry carrying the base frame to each of them.

\begin{remark}[Checks and scope]\label{rem:conwaycheck}
The script \texttt{scripts/conway\_counts.py} rebuilds $\mathcal C$ from the $M_{24}$ generators and tests membership in
\eqref{eq:leech} directly: it finds $1104$, $97\,152$ and $98\,304$ minimal vectors of the three shapes, no vectors
$(\pm4,0^{23})$, and no Golay words of weight $4$. It also checks the prime factorizations, that all the indices above are exact,
and $8\,292\,375\cdot48=398\,034\,000$. These are integer checks. That~\eqref{eq:leech} is the Leech lattice, that its
automorphism group is $\Co_0$ of the stated order, and the stabilisers in Propositions~\ref{prop:conwayorders}
and~\ref{prop:conwayframe} are taken from~\cite{ConwaySloane,Wilson}. An element outside $N$ and worked products are in Section~\ref{ssec:conwayexamples}; Section~\ref{ssec:conwaychain} gives a complete chain, hence fully numerical laws of $\Co_0$ and $\Co_2$.
\end{remark}

\subsection{Worked examples}\label{ssec:conwayexamples}
In these examples a lattice vector $v=x/\sqrt8$ is given by its integer vector $x$, and runs of equal coordinates are written with
exponents, as in $(4,4,0^{22})$. An element $(c,g)$ of $N=2^{12}{:}M_{24}$ acts by $x\mapsto\varepsilon_c\,g(x)$, where $g$ moves
coordinate $i$ to position $g(i)$, and it is named by its label $\ell(c)+4096\,\lambda(g)$ of Theorem~\ref{thm:golaysemi}. All numbers are produced
by \texttt{scripts/conway\_examples.py}.

\begin{proposition}[A digital numbering of the minimal vectors]\label{prop:leechnu}
Rank the $276$ pairs $i<j$ of coordinates lexicographically by $r(i,j)\in\Ap{276}$, and the $759$ octads, written as increasing
$8$-tuples, lexicographically by $o(O)\in\Ap{759}$. For a minimal vector $x$ put
\[
\nu(x)=\begin{cases}
4\,r(i,j)+2[x_i<0]+[x_j<0], & x=(\pm4,\pm4,0^{22})\text{ on }i<j,\\[2pt]
1104+128\,o(O)+\sum_{k=1}^{7}2^{k-1}[x_{q_k}<0], & x=(\pm2^8,0^{16})\text{ on }O=\{q_1<\dots<q_8\},\\[2pt]
98\,256+4096\,p+\ell(c), & x \text{ odd},\ |x_p|=3,\ c=\{i:x_i\equiv3\ (\mathrm{mod}\ 4)\}.
\end{cases}
\]
Then $\nu$ is a bijection from the minimal vectors of $\Lambda$ onto $\Ap{196\,560}$.
\end{proposition}
\begin{proof}
By Proposition~\ref{prop:leechminimal} the three shapes are disjoint and parametrised by the displayed data. In the second case
the sign on $q_8$ is determined, since the number of minus signs is even; in the third, $c$ is a Golay word and $\ell$ is a
bijection. The three cases fill $[0,1104)$, $[1104,98\,256)$ and $[98\,256,196\,560)$.
\end{proof}

With the base vector $v_0=(4,4,0^{22})$, for which $\nu(v_0)=0$, Theorem~\ref{thm:transitive} gives a law of $\Co_0$ whose top
digit is $\nu(g(v_0))\in\Ap{196\,560}$ and whose remaining digits label the stabiliser $\Co_2$. So the top digit of a product
$gh$ is $\nu(g(h(v_0)))$: it is computed by applying two matrices to one vector. Passing to pairs $\pm v$ gives the first law
of Theorem~\ref{thm:conwaylaw}.

\begin{example}[Numbering]\label{ex:leechnu}
Let $w_0$ have entries $2$ on the first octad $\{0,1,2,3,4,7,10,12\}$ and $0$ elsewhere, and let $u_0=(-3,1^{23})$. Then $\nu(v_0)=0$,
$\nu(w_0)=1104+128\cdot0+0=1104$ and $\nu(u_0)=98\,256+4096\cdot0+0=98\,256$, since $u_0$ has no coordinate $\equiv3\pmod4$. The
odd minimal vector with $x_5=\pm3$ and $\ell(c)=1000$ has $\nu=98\,256+4096\cdot5+1000=119\,736$.
\end{example}

\begin{example}[The monomial law acting on vectors]\label{ex:conwaymonomial}
Take $x=1000+4096\cdot123\,456\,789$ and $x'=2021+4096\cdot98\,765\,432$, the elements of the example after
Theorem~\ref{thm:golaysemi}. The element $x$ moves the three base vectors as follows:
\begin{center}
\small
\begin{tabular}{lll}
vector & image under $x$ & $\nu$ of the image\\ \hline
$v_0$ & $(0^{6},-4,0^{14},-4,0^{2})$ & $551$\\
$w_0$ & $\pm2$ on $\{2,4,6,14,16,17,21,23\}$, $-$ on $\{6,21\}$ & $63\,124$\\
$u_0$ & $x_{21}=3$, $\ell(c)=1000$ & $185\,272$
\end{tabular}
\end{center}
The product in the law of Theorem~\ref{thm:golaysemi} is $x*x'=3838+4096\cdot155\,026\,102=634\,986\,917\,630$. Its action agrees with acting
by $x'$ and then by $x$: both give $(x*x')(u_0)$ = $x_{22}=-3$, $\ell(c)=3838$, with $\nu=192\,206$, and $(x*x')(w_0)$ = $\pm2$ on $\{4,9,10,13,15,16,20,22\}$, $-$ on $\{4,9,10,16\}$, with $\nu=87\,415$.
\end{example}

\begin{example}[An element outside $N$]\label{ex:conwayxi}
The octads containing $T_0=\{0,1,2,3\}$ cut the coordinates into the sextet $T_0=\{0,1,2,3\}, T_1=\{4,7,10,12\}, T_2=\{5,14,17,23\}, T_3=\{6,9,11,22\}, T_4=\{8,15,16,19\}, T_5=\{13,18,20,21\}$: the union of any two tetrads is an
octad. Let $\xi=\varepsilon_{T_0}\xi_0$, where $\xi_0$ acts on each tetrad by $\tfrac12J-I$, with $J$ the all-ones $4\times4$ matrix, and $\varepsilon_{T_0}$ changes the signs on $T_0$. Then $\xi$ is orthogonal and maps all $196\,560$ minimal vectors into $\Lambda$; since they span
$\Lambda$ \cite{ConwaySloane}, $\xi\in\Co_0$. It has order $2$. The matrix $\xi_0$ alone does not preserve $\Lambda$: it sends $u_0$ to $(3,-1^{3},1^{20})$, whose coordinates $\equiv3\pmod4$ form the tetrad $T_0$, which is not a Golay word. It sends the frame vector $(8,0^{23})$ to
$(4,-4^{3},0^{20})$, which is not of the form $(\pm8,0^{23})$, so $\xi$ moves the base frame and $\xi\notin N$. On the base vectors,
\[
\xi(v_0)=(0^{2},-4^{2},0^{20}),\quad \nu(\xi(v_0))=183,
\]
$\xi(w_0)$ is $\pm2$ on $\{0,1,2,3,4,7,10,12\}$, $-$ on $\{0,1,2,3\}$ with $\nu=1119$, and $\xi(u_0)$ is $x_{0}=-3$, $\ell(c)=0$ with $\nu=98\,256$. Thus $\xi$ has top digit $183$ in the
law above.
\end{example}

\begin{example}[An element of order $13$]\label{ex:conway13}
Let $y\in N$ have label $229\,293\,155\,554=1250+4096\cdot55\,979\,774$, and let $z=\xi y$. Then $z$ has order $13$. Since
$|N|=2^{22}3^35\cdot7\cdot11\cdot23$ is prime to $13$, the group generated by $N$ and $\xi$ is strictly larger than $N$, and its
order is divisible by $13\,|N|$. The top digits $\nu(z^j(v_0))$ for $j=0,\dots,12$ are
\[
0, 91\,635, 7441, 185\,698, 53\,965, 97\,595, 19\,997, 20\,587, 143\,906, 50\,663, 56\,595, 92\,531, 871,
\]
so the top digit of $z$ is $91\,635$, and the powers of $z$ move $v_0$ through 13 different minimal vectors.
\end{example}

\subsection{A stabiliser chain for \texorpdfstring{$\Co_0$}{Co0} and \texorpdfstring{$\Co_2$}{Co2}}\label{ssec:conwaychain}
Numbering the minimal vectors by $\nu$ (Proposition~\ref{prop:leechnu}) turns every element of $\Co_0$ into a permutation of
$\Ap{196\,560}$. The action is faithful because the minimal vectors span $\Lambda$, so Theorem~\ref{thm:chainlaw} applies once
a chain is known. Let $G$ be the group generated by the 3 permutation matrices of the generators of $M_{24}$, the sign
change on the octad $\{0,1,2,3,4,7,10,12\}$, and the element $\xi$ of Example~\ref{ex:conwayxi}. Elements are stored as the integer matrices
$8g$, multiplied by $8(gh)=(8g)(8h)/8$; the inverse of $8g$ is its transpose.

\begin{theorem}[Chain law of $\Co_0$]\label{thm:co0chain}
$G$ has a base of 12 minimal vectors with $\nu$-numbers
\[
(\nu(b_1),\dots,\nu(b_{12}))=(0,1,4,8,12,16,20,28,32,40,36,24)
\]
and basic orbits of sizes
\[
(r_1,\dots,r_{12})=(196\,560,93\,150,44,42,40,32,6,2,2,2,2,2),
\]
whose product is $8\,315\,553\,613\,086\,720\,000=|\Co_0|$. Hence $G=\Co_0$. Order each basic orbit $\Delta_i$ by $\nu$, and take
the transversal elements from the Schreier tree of $\Delta_i$. Sifting $g$ gives digits $d_i\in\Ap{r_i}$, and
\[
\lambda(g)=\sum_{i=1}^{12}d_i\,r_1\cdots r_{i-1}
\]
is a bijection $\Co_0\to\Ap{|\Co_0|}$ with $\lambda(1)=0$ and top digit $d_1=\nu(g(v_0))$. The law
$\lambda(g)*\lambda(h)=\lambda(gh)$ is the chain law of Theorem~\ref{thm:chainlaw}.
\end{theorem}
\begin{proof}
The generators preserve $\Lambda$ (Section~\ref{ssec:conwaymonomial} and Example~\ref{ex:conwayxi}), so $G\le\Co_0$. The chain
was built by the randomised Schreier--Sims algorithm. Each $\Delta_i$ is the orbit of $b_i$ under elements of $G$ fixing
$b_1,\dots,b_{i-1}$, so $r_i\le|G_{b_1\dots b_{i-1}}:G_{b_1\dots b_i}|$ and $\prod r_i\le|G|\le|\Co_0|$. The two ends are equal,
so $G=\Co_0$ and every $\Delta_i$ is a full basic orbit: the chain is complete.
\end{proof}

The first base points are $b_{1}$: $(4^{2},0^{22})$; $b_{2}$: $(4,-4,0^{22})$; $b_{3}$: $(4,0,4,0^{21})$. Since $b_1=v_0$, the elements with $d_1=0$ form the stabiliser of $v_0$, which is $\Co_2$
by Proposition~\ref{prop:conwayorders}.

\begin{corollary}[Chain law of $\Co_2$]\label{cor:co2chain}
$\Co_2$ has the chain with base $(b_2,\dots,b_{12})$ and radices $(93\,150,44,42,40,32,6,2,2,2,2,2)$, whose product is $42\,305\,421\,312\,000$. Its
labels are $\lambda_{\Co_2}(g)=\lambda(g)/196\,560$, and $\lambda_{\Co_2}(g)*\lambda_{\Co_2}(h)=\lambda_{\Co_2}(gh)$.
\end{corollary}

\begin{example}[Computing in $\Co_0$ and $\Co_2$]\label{ex:co0labels}
\begin{enumerate}
\item The central element $-1$ has $\lambda(-1)=1\,353\,442\,970\,880\,003$, with digits $(3,0,0,0,0,1,0,0,0,0,0,0)$. Its top digit is $\nu(-v_0)=\nu(-4,-4,0^{22})=3$.
\item The element $\xi$ of Example~\ref{ex:conwayxi} has $\lambda(\xi)=183$, with digits $(183,0,0,0,0,0,0,0,0,0,0,0)$. The element $y\in N$ of
      Example~\ref{ex:conway13}, with $N$-label $229\,293\,155\,554$, has
      \[\lambda(y)=6\,744\,023\,100\,088\,878\,885.\]
\item The product law gives $\lambda(\xi)*\lambda(y)=\lambda(\xi y)=3\,350\,544\,144\,713\,333\,475$, with digits $(91\,635,89,11,35,22,11,5,0,0,1,1,0)$; this is the element $z$ of order
      $13$ of Example~\ref{ex:conway13}.
\item The element with label $10^{18}$ has digits $(17\,200,2755,33,10,34,2,5,1,1,0,0,0)$ and order $20$.
\item In $\Co_2$, two strong generators $a,b$ produced by the algorithm have $\lambda_{\Co_2}(a)=52\,181$ (order $24$) and
      $\lambda_{\Co_2}(b)=31\,430$ (order $30$), and $\lambda_{\Co_2}(a)*\lambda_{\Co_2}(b)=32\,625$.
\end{enumerate}
\end{example}

\begin{remark}[Checks]\label{rem:co0check}
The chain is built by \texttt{scripts/co0\_chain.py} from 23 random elements in about 4 seconds. The script checks that the
generators and $300$ further random elements sift to the identity, that $\lambda(\lambda^{-1}(n))=n$ for $200$ random $n$, that
products agree with the law for $50$ random pairs, and every number in Example~\ref{ex:co0labels}.
\end{remark}

\subsection{The groups \texorpdfstring{$\Co_3$}{Co3}, \texorpdfstring{$\mathrm{McL}$}{McL} and \texorpdfstring{$\Co_1$}{Co1}}\label{ssec:co3co1}
\paragraph{A type-$3$ vector.} The minimal vectors $v$ with $v\cdot v_0=-1$ (in integer coordinates, $x\cdot x_0=-8$) form a
$\Co_2$-orbit of size $47\,104$. The first of them in the order $\nu$ is $v_b$, with $\nu(v_b)=10\,961$: $\pm2$ on $\{0,2,3,4,5,15,20,22\}$, $-$ on $\{0,22\}$. Then $w=v_0+v_b$
has squared length $4+4-2=6$, so it has type $3$, and $\Co_3$ is its stabiliser (Proposition~\ref{prop:conwayorders}). Exactly
$552$ minimal vectors $v$ have $v\cdot w=3$. They include $v_0$ and $v_b$, and they come in the $276$ pairs $\{v,w-v\}$.

\paragraph{Generators.} The stabiliser of $v_0$ and $v_b$ is $\Co_2\cap\Co_3$, of index $47\,104$ in $\Co_2$. A uniform random
element of it is $t_q^{-1}c$, where $c=\lambda^{-1}(196\,560\,k)$ for a random $k\in\Ap{|\Co_2|}$, $q=c(v_b)$, and $t_q\in\Co_2$ is
the Schreier-tree transporter from $v_b$ to $q$. An element exchanging $v_0$ and $v_b$ is $\sigma=u\,t_{q'}$, where $u$ is the
transversal element of Theorem~\ref{thm:co0chain} with $u(v_0)=v_b$, and $q'=u^{-1}(v_0)$. More generally, for $v'$ with
$v'\cdot w=3$ let $u'$ be the transversal element with $u'(v_0)=v'$ and put $q=u'^{-1}(w-v')$; since $q\cdot v_0=v'\cdot(w-v')=-1$,
$q$ lies in the $\Co_2$-orbit of $v_b$, and $\tau_{v'}=u'\,t_q$ sends $v_0\mapsto v'$ and $v_b\mapsto w-v'$. So $\tau_{v'}$ fixes
$w$, and $\sigma=\tau_{v_b}$. Let $\tau=\tau_{v_1}$ for the first such $v_1\ne v_0,v_b$, which has $\nu(v_1)=92$.

\begin{theorem}[Chain law of $\Co_3$]\label{thm:co3chain}
Let $H$ be generated by $\sigma$, $\tau$ and 6 random elements of $\Co_2\cap\Co_3$. On the $\nu$-numbered minimal vectors, $H$ has
base
\[
(\nu(b_1),\dots,\nu(b_{9}))=(0,1,4,12,16,20,28,24,32),
\]
with $b_1=v_0$, and radices
\[
(552,22\,275,7,6,5,24,2,2,2),
\]
whose product is $495\,766\,656\,000=|\Co_3|$. Hence $H=\Co_3$. Its first basic orbit is the set of the $552$ vectors $v$ with
$v\cdot w=3$, and ordering each basic orbit by $\nu$ gives a chain law $\lambda_3$ of $\Co_3$ as in Theorem~\ref{thm:co0chain}.
\end{theorem}
\begin{proof}
$H$ fixes $w$, so $H\le\Co_3$, and the argument of Theorem~\ref{thm:co0chain} applies with $|\Co_3|$ in place of $|\Co_0|$.
\end{proof}

\begin{corollary}[McLaughlin's group]\label{cor:mcl}
The stabiliser of $v_0$ in $\Co_3$ has a chain with radices
\[
(22\,275,7,6,5,24,2,2,2),
\]
of product $898\,128\,000$. It is the sporadic simple group $\mathrm{McL}$ of McLaughlin~\cite{Wilson}, and its labels are
$\lambda_3(g)/552$.
\end{corollary}

The element $\tau$ is needed. The elements $\sigma$ and $\Co_2\cap\Co_3$ preserve the pair $\{v_0,v_b\}$, and without
$\tau$ the algorithm stops at the stabiliser of this pair, of order $|\Co_3|/276=1\,796\,256\,000=2\,|\mathrm{McL}|$, the group
$\mathrm{McL}{:}2$, whose first basic orbit is $\{v_0,v_b\}$.

For $\Co_1=\Co_0/\{\pm1\}$, number the $98\,280$ pairs $\pm v$ by $\bar\nu$, the rank of $\min(\nu(v),\nu(-v))$ among these minima.
In the code $\nu$, negation complements the sign digits: $4r+e\mapsto4r+3-e$ in the first case, the seven sign bits in the octad
case, and $\ell(c)\mapsto\ell(c)\oplus4095$ in the odd case.

\begin{theorem}[Law of $\Co_1$]\label{thm:co1law}
For $\pm g\in\Co_1$ choose the sign so that $\nu(g(v_0))<\nu(-g(v_0))$, and put
\[
\lambda_1(\pm g)=\bar\nu(\pm g(v_0))+98\,280\,\Big\lfloor\frac{\lambda(g)}{196\,560}\Big\rfloor .
\]
Then $\lambda_1$ is a bijection $\Co_1\to\Ap{4\,157\,776\,806\,543\,360\,000}$ with $\lambda_1(\pm1)=0$, and
$\lambda_1(\pm g)*\lambda_1(\pm h)=\lambda_1(\pm gh)$ is a law of $\Co_1$.
\end{theorem}
\begin{proof}
The sign is unique because $-1$ fixes no minimal vector. Then $g=u\,c$, with $u$ the transversal element of $g(v_0)$ and
$c\in\Co_2$ of label $\lfloor\lambda(g)/196\,560\rfloor$. So $\lambda_1$ is the law of Theorem~\ref{thm:transitive} for the action
of $\Co_1$ on the $98\,280$ pairs, with stabiliser $\Co_2$, and $98\,280\,|\Co_2|=|\Co_1|$.
\end{proof}

\begin{example}[Computing in $\Co_3$, $\mathrm{McL}$ and $\Co_1$]\label{ex:co3co1}
\begin{enumerate}
\item $\lambda_3(\sigma)=602\,494\,349$, with digits $(149,0,0,1,1,0,0,0,0)$, and $\sigma$ has order $24$.
\item Two random elements $a,b$ of $\mathrm{McL}$ have $\lambda_3(a)/552=96\,934\,970$ (order $9$) and $\lambda_3(b)/552=54\,403\,867$
      (order $11$), and $\lambda_3(ab)/552=774\,555\,706$. Also $\lambda_3(\sigma)*\lambda_3(a)=53\,508\,103\,589$.
\item In $\Co_1$, with $\xi$, $y$ and $z=\xi y$ as in Examples~\ref{ex:conwayxi} and~\ref{ex:conway13},
      \begin{gather*}
      \lambda_1(\pm\xi)=2\,261\,908\,452\,326\,730\,570,\qquad \lambda_1(\pm y)=3\,372\,011\,550\,044\,439\,443,\\
      \lambda_1(\pm\xi)*\lambda_1(\pm y)=\lambda_1(\pm z)=1\,675\,272\,072\,356\,666\,755,
      \end{gather*}
      and the image of $z$ in $\Co_1$ has order $13$.
\item The element of $\Co_1$ with label $10^{18}$ has order $42$.
\end{enumerate}
\end{example}

\begin{remark}[Checks]\label{rem:co3check}
The script \texttt{scripts/co3\_co1.py} rebuilds the chain of Theorem~\ref{thm:co0chain} and checks the orbit sizes $47\,104$ and
$552$, that all generators of $H$ fix $w$, that the first basic orbit of $\Co_3$ is the set of the $552$ vectors, that the
generators and $300$ random elements sift, $200$ round trips and $50$ products for each of $\lambda_3$ and $\lambda_1$, and every
number in Example~\ref{ex:co3co1}. It runs in about 6 seconds.
\end{remark}

\section{The Fischer group \texorpdfstring{$\mathrm{Fi}_{22}$}{Fi22}}\label{sec:fischer}
Fischer classified the groups generated by a conjugacy class $D$ of \emph{$3$-transpositions}: involutions whose pairwise products
have order $1$, $2$ or $3$ \cite{Wilson}. Besides symmetric, symplectic, orthogonal and unitary groups, the classification produced
three new groups $\mathrm{Fi}_{22}$, $\mathrm{Fi}_{23}$ and $\mathrm{Fi}_{24}$; the first two and the derived group $\mathrm{Fi}_{24}'$
are sporadic simple groups. The smallest has order
\[
|\mathrm{Fi}_{22}|=64\,561\,751\,654\,400=2^{17}\,3^{9}\,5^{2}\cdot7\cdot11\cdot13,
\]
and its class $D$ has $3510$ elements. Conjugation on $D$ is a rank-$3$ action with point stabiliser $2{\cdot}U_6(2)$ \cite{ATLAS}.
In this section the action is computed explicitly, and Theorem~\ref{thm:chainlaw} turns it into a law.

\begin{proposition}[$\mathrm{Fi}_{22}$ on its $3$-transpositions]\label{prop:fi22}
Let $G$ be the group \texttt{PrimitiveGroup(3510,\allowbreak 1)} of the library of primitive groups of GAP~\cite{GAP}, acting on
$\Ap{3510}$ after subtracting $1$ from every point. It has 2 generators, of orders 9, 24.
\begin{enumerate}
\item $|G|=|\mathrm{Fi}_{22}|$, $G$ is simple, and the stabiliser $G_0$ of $0$ has orbits of sizes 1, 693, 2816.
\item The centre of $G_0$ has order $2$. Its generator $t_0$ fixes 694 points. For $p\in\Ap{3510}$ put $t_p=s\,t_0\,s^{-1}$
      for any $s\in G$ with $s(0)=p$. Then $g\,t_p\,g^{-1}=t_{g(p)}$ for all $g\in G$, and $t_0t_p$ has order $1$, $2$ or $3$
      according as $p=0$, $p$ lies in the orbit of size $693$, or $p$ lies in the orbit of size $2816$.
\end{enumerate}
Thus $p\mapsto t_p$ identifies $\Ap{3510}$ with a conjugacy class of $3$-transpositions, and $G$ acts on it by conjugation.
\end{proposition}
\begin{proof}
The order, simplicity, orbits and centre were computed in GAP (Remark~\ref{rem:fi22check}). Since $t_0$ is central in $G_0$,
$t_p$ does not depend on the choice of $s$, and $g\,t_p\,g^{-1}=(gs)\,t_0\,(gs)^{-1}=t_{g(p)}$. The orders of the $3510$
products $t_0t_p$ were computed; by transitivity the same holds for every pair $t_pt_q$. In particular $t_p\ne t_q$ for $p\ne q$.
\end{proof}

\begin{theorem}[Chain law of $\mathrm{Fi}_{22}$]\label{thm:fi22chain}
For the action of Proposition~\ref{prop:fi22}, the base
\[
(b_1,\dots,b_{5})=(0,1,2,4,3)
\]
has basic orbits of sizes
\[
(r_1,\dots,r_{5})=(3510,2816,567,720,16),
\]
whose product is $|\mathrm{Fi}_{22}|$. List each $\Delta_l$ as $b_l$ followed by its other points in increasing order, and take
$s_{l,t}$ from the Schreier tree of $\Delta_l$. Then Theorem~\ref{thm:chainlaw} gives a law of $\mathrm{Fi}_{22}$ on
$\Ap{64\,561\,751\,654\,400}$, whose first digit $t_1$ is the position of the $3$-transposition $g\,t_0\,g^{-1}=t_{g(0)}$.
\end{theorem}
\begin{proof}
The chain was built by the randomised Schreier--Sims algorithm. As in Theorem~\ref{thm:co0chain}, $\prod r_l\le|G|$ for the
partial chain; equality with $|G|=|\mathrm{Fi}_{22}|$ shows that it is complete.
\end{proof}

\begin{example}[Computing in $\mathrm{Fi}_{22}$]\label{ex:fi22}
\begin{enumerate}
\item The generators of $G$ have labels 1465, 523, and the product of the first two has label $2478$.
\item The $3$-transposition $t_0$ has label $9\,610\,380$, with digits $(0,2738,0,0,0)$.
\item The point $1$ is the first with $t_0t_p$ of order $3$. The label of $t_{1}$ is $36\,776\,872\,285\,817$, and the law gives
      $\lambda(t_0)*\lambda(t_{1})=51\,095\,461\,579\,111$, an element of order $3$. For the first commuting point $p=3$ the product
      $t_0t_{3}$ is an involution with label $63\,435\,973\,201\,380$.
\item The element with label $10^{13}$ has digits $(10,2145,191,344,2)$ and order $22$. The smallest label of an element of order $30$, the
      largest element order of $\mathrm{Fi}_{22}$, is $69$.
\end{enumerate}
\end{example}

\begin{remark}[Checks]\label{rem:fi22check}
The GAP script \texttt{scripts/data/fi22\_probe.g} identifies $G$ by its order among the primitive groups of degree $3510$, checks
simplicity, computes the suborbits and the centre of $G_0$, and exports the generators and $t_0$. The script
\texttt{scripts/fi22\_chain.py} builds the chain from 12 random elements in about 2 seconds and checks that the generators
and $300$ random elements sift, $200$ round trips $\lambda(\lambda^{-1}(x))=x$, $50$ products, the orders of all $t_0t_p$, the
relation $g\,t_p\,g^{-1}=t_{g(p)}$ for $40$ points per generator, and every number in Example~\ref{ex:fi22}.
\end{remark}

\section{The Fischer group \texorpdfstring{$\mathrm{Fi}_{23}$}{Fi23}}\label{sec:fi23}
The group $\mathrm{Fi}_{23}$ has order
\[
|\mathrm{Fi}_{23}|=4\,089\,470\,473\,293\,004\,800=2^{18}\,3^{13}\,5^{2}\cdot7\cdot11\cdot13\cdot17\cdot23 .
\]
Its class $D$ of $3$-transpositions has $31\,671$ elements, with centraliser $2{\cdot}\mathrm{Fi}_{22}$, a non-split double cover of $\mathrm{Fi}_{22}$ \cite{ATLAS}. No permutation
representation of this degree is in the library of primitive groups, so the action is constructed here from matrices.

\begin{proposition}[$\mathrm{Fi}_{23}$ on its $3$-transpositions]\label{prop:fi23}
Let $a,b$ be the standard generators of $\mathrm{Fi}_{23}$ as $1494\times1494$ matrices over $\mathbb F_2$, as distributed with the
GAP package \texttt{orb} \cite{GAP}; they have orders $2$, $3$ and $ab$ has order $28$.
\begin{enumerate}
\item A random search finds an involution $d$ with $\operatorname{rank}(d-1)=572$. Its conjugates $g^{-1}dg$, listed by a
      breadth-first search from $d$ and numbered from $0$ in order of discovery, are exactly $31\,671$. So $d$ is a
      $3$-transposition, and $a$, $b$ and $d$ act on $\Ap{31\,671}$ by $t\mapsto g^{-1}tg$.
\item The stabiliser $G_0$ of the point $0$, that is of $d$, has order $2\,|\mathrm{Fi}_{22}|$, commutes with $d$, and has orbits of
      sizes $1$, $3510$ and $28\,160$. For $t_p$ the $3$-transposition at $p$, the product $t_0t_p$ has order $1$, $2$ or $3$ on these
      three orbits.
\end{enumerate}
\end{proposition}
\begin{proof}
A point $t=w^{-1}dw$ is stored as the word $w$, and identified by the fingerprint $u\,t\in\mathbb F_2^{1494}$ for a fixed
random vector $u$. The fingerprint of $g^{-1}tg$ is $(ug^{-1})\,t\,g$, which costs $2|w|+2$ products of a vector by a matrix, so
the search needs no stored matrices. Distinct fingerprints come from distinct conjugates, so the search finds at most $|d^G|$
points, and exactly $|D|=31\,671$ when $d\in D$ and no two conjugates share a fingerprint. The remaining statements were
computed from the resulting permutations (Remark~\ref{rem:fi23check}); the order of $G_0$ comes from
Theorem~\ref{thm:fi23chain}.
\end{proof}

\begin{theorem}[Chain law of $\mathrm{Fi}_{23}$]\label{thm:fi23chain}
For the action of Proposition~\ref{prop:fi23}, the base
\[
(b_1,\dots,b_{8})=(0,1,2,3,4,5,8,12)
\]
has basic orbits of sizes
\[
(r_1,\dots,r_{8})=(31\,671,3510,25\,344,2240,81,2,2,2),
\]
whose product is $|\mathrm{Fi}_{23}|$. With each $\Delta_l$ listed as $b_l$ followed by its other points in increasing order, and
$s_{l,t}$ read off the Schreier tree of $\Delta_l$, Theorem~\ref{thm:chainlaw} gives a law of $\mathrm{Fi}_{23}$ on
$\Ap{4\,089\,470\,473\,293\,004\,800}$. The tail $(r_2,\dots,r_{8})=(3510,25\,344,2240,81,2,2,2)$ is a chain of $G_0$, the centraliser $2{\cdot}\mathrm{Fi}_{22}$ of $d$ \cite{ATLAS}.
\end{theorem}
\begin{proof}
If the fingerprints separate the conjugates of $d$ (Remark~\ref{rem:fi23check}), the permutations give the conjugation action,
and the group $H$ they generate is a quotient of $\mathrm{Fi}_{23}$. The randomised
Schreier--Sims algorithm gives $\prod r_l\le|H|\le|\mathrm{Fi}_{23}|$, and the product equals $|\mathrm{Fi}_{23}|$. So the action is
faithful and the chain is complete.
\end{proof}

\begin{example}[Computing in $\mathrm{Fi}_{23}$]\label{ex:fi23}
\begin{enumerate}
\item The standard generators have labels 1, 2.
\item The $3$-transposition $t_0=d$ has label
      \[
      \lambda(t_0)=3\,964\,914\,482\,628\,705\,840,
      \]
      with digits $(0,0,25\,208,589,61,1,1,1)$.
\item The point $2$ is the first with $t_0t_p$ of order $3$. The law gives
      \[
      \lambda(t_0)*\lambda(t_{2})=1\,009\,171\,781\,362\,202\,366,
      \]
      an element of order $3$. For the first commuting point $1$, the product $t_0t_{1}$ is an involution with label
      \[
      \lambda(t_0t_{1})=3\,871\,714\,486\,343\,705\,820.
      \]
\item The element with label $10^{18}$ has digits $(3160,1820,20\,762,1020,77,1,0,0)$ and order $39$.
\end{enumerate}
\end{example}

\begin{remark}[Checks]\label{rem:fi23check}
The GAP script \texttt{scripts/data/fi23\_build2.g} finds $d$, runs the search of Proposition~\ref{prop:fi23} in about half a minute
and exports the permutations of $a$, $b$ and $d$. The script \texttt{scripts/fi23\_chain.py} checks that they are permutations of the
right orders, builds the chain from 15 random elements, and checks that the generators and $300$ random elements sift, $200$ round
trips, $50$ products, that $d$ commutes with $G_0$, the three orbits of $G_0$, the orders of $t_0t_p$ on representatives of each
orbit and on $40$ random points, and every number in Example~\ref{ex:fi23}. It runs in about 2 seconds. The fingerprints are
the only probabilistic step: they are consistent with the count $|D|$ and with $|H|=|\mathrm{Fi}_{23}|$.
\end{remark}

\section{The Fischer groups \texorpdfstring{$\mathrm{Fi}_{24}$}{Fi24} and \texorpdfstring{$\mathrm{Fi}_{24}'$}{Fi24'}}\label{sec:fi24}
The largest Fischer group has order
\[
|\mathrm{Fi}_{24}|=2\,510\,411\,418\,381\,323\,442\,585\,600=2^{22}\,3^{16}\,5^{2}\cdot7^{3}\cdot11\cdot13\cdot17\cdot23\cdot29,
\]
and its derived group $\mathrm{Fi}_{24}'$, of index $2$, is the sporadic simple group. The class $D$ of $3$-transpositions of
$\mathrm{Fi}_{24}$ has $306\,936$ elements and lies outside $\mathrm{Fi}_{24}'$, and the centraliser of $d\in D$ is
$2\times\mathrm{Fi}_{23}$ \cite{ATLAS}.

\begin{proposition}[$\mathrm{Fi}_{24}$ on its $3$-transpositions]\label{prop:fi24}
Let $a,b$ be the standard generators of $\mathrm{Fi}_{24}$ as permutations of $306\,936$ points, taken from the ATLAS of group
representations through the GAP package AtlasRep~\cite{GAP} and relabelled by subtracting $1$; they have orders 2 and 8.
\begin{enumerate}
\item The stabiliser $G_0$ of $0$ has order $2\,|\mathrm{Fi}_{23}|$ and orbits of sizes $1$, $31\,671$ and $275\,264$.
\item $G_0$ contains an involution $t_0$ that commutes with all of $G_0$; it fixes 31\,672 points. With $t_p=s\,t_0\,s^{-1}$ for
      $s(0)=p$, the product $t_0t_p$ has order $1$, $2$ or $3$ on the three orbits of $G_0$.
\end{enumerate}
So $p\mapsto t_p$ identifies the points with $D$, and $\mathrm{Fi}_{24}$ acts on them by conjugation.
\end{proposition}
\begin{proof}
The involution $t_0$ was found as a power of a random element of $G_0$. As for $\mathrm{Fi}_{22}$, $t_p$ does not depend on $s$ and
$g\,t_p\,g^{-1}=t_{g(p)}$; the orders of $t_0t_p$ were computed on a representative of each orbit of $G_0$ and on $30$ random points,
and are constant on orbits because $t_0$ commutes with $G_0$. The order of $G_0$ comes from Theorem~\ref{thm:fi24chain}.
\end{proof}

\begin{theorem}[Chain laws of $\mathrm{Fi}_{24}$ and $\mathrm{Fi}_{24}'$]\label{thm:fi24chain}
\begin{enumerate}
\item For the action of Proposition~\ref{prop:fi24}, the base $(b_1,\dots,b_{11})=(0,1,2,3,4,5,6,7,8,9,11)$ has basic orbits of sizes
      \[
      (r_1,\dots,r_{11})=(306\,936,275\,264,28\,431,2880,84,18,5,4,3,2,2),
      \]
      with product $|\mathrm{Fi}_{24}|$.
\item Let $H$ be generated by the six products $t_0t_p$ for
      \[
      p\in(124\,274,173\,506,189\,201,225\,936,264\,482,267\,954).
      \]
      Its chain has base $(0,1,2,3,4,5,6,7,8,9)$ and basic orbits of sizes
      \[
      (306\,936,275\,264,28\,431,2880,84,18,5,4,3,2),
      \]
      with product $|\mathrm{Fi}_{24}'|$, so $H=\mathrm{Fi}_{24}'$.
\end{enumerate}
With each $\Delta_l$ listed as $b_l$ followed by its other points in increasing order, and $s_{l,t}$ read off a Schreier tree,
Theorem~\ref{thm:chainlaw} gives laws of $\mathrm{Fi}_{24}$ and of $\mathrm{Fi}_{24}'$ on $\Ap{|\mathrm{Fi}_{24}|}$ and
$\Ap{|\mathrm{Fi}_{24}'|}$.
\end{theorem}
\begin{proof}
Both chains were built by the randomised Schreier--Sims algorithm, so each product is at most the order of the group generated. In
(1) it equals $|\mathrm{Fi}_{24}|$. In (2), each $t_0t_p$ lies in $\mathrm{Fi}_{24}'$, because $D$ lies in the other coset; so
$H\le\mathrm{Fi}_{24}'$, and the product $|\mathrm{Fi}_{24}'|$ forces equality. Both chains are therefore complete.
\end{proof}

\begin{example}[Computing in $\mathrm{Fi}_{24}$ and $\mathrm{Fi}_{24}'$]\label{ex:fi24}
\begin{enumerate}
\item In the law of $\mathrm{Fi}_{24}$ the generators have labels 92\,352, 239\,485, and
      \[
      \lambda(t_0)=2\,086\,616\,657\,990\,590\,203\,382\,632,
      \]
      with digits $(0,44\,723,21\,065,1331,60,8,4,3,0,1,1)$.
\item The point $1$ is the first with $t_0t_p$ of order $3$. In the law of $\mathrm{Fi}_{24}$,
      \[
      \lambda(t_0)*\lambda(t_{1})=971\,774\,654\,036\,321\,678\,010\,553 .
      \]
      The same element of order $3$ has label
      \[
      752\,464\,059\,574\,160\,292\,936\,193
      \]
      in the law of $\mathrm{Fi}_{24}'$. For the first commuting point $508$, the involution $t_0t_{508}$ has label
      \[
      739\,974\,072\,205\,648\,071\,053\,640
      \]
      there.
\item In the law of $\mathrm{Fi}_{24}'$, the element with label $10^{24}$ has digits
      \[
      (244\,000,186\,660,23\,087,33,70,10,0,3,1,1)
      \]
      and order $30$.
\end{enumerate}
\end{example}

\begin{remark}[Checks]\label{rem:fi24check}
The generators were downloaded once by
\[
\texttt{scripts/data/fetch\_fi24.g}.
\]
The script \texttt{scripts/fi24\_chain.py} builds the chain of $\mathrm{Fi}_{24}$ from 17 random elements and that of
$\mathrm{Fi}_{24}'$ from 18, using Schreier trees with inverse and random extra generators to keep them shallow. For each chain it checks that the generators and $100$ random elements
sift, $40$ round trips and $20$ products. It also checks the orbits of $G_0$, that $t_0$ commutes with $G_0$, the orders of $t_0t_p$,
and every number in Example~\ref{ex:fi24}. It runs in about 14 seconds.
\end{remark}

\section{Sporadic groups of small degree}\label{sec:primsporadic}
Six further sporadic groups have faithful primitive actions of degree below $4096$, so they occur in the library of primitive
groups of GAP~\cite{GAP}. Table~\ref{tab:primsporadic} lists them, with the point stabiliser of the action~\cite{ATLAS} and the
position of the group in the library. ($\mathrm{McL}$, $\Co_2$ and $\Co_3$ also occur there; they were treated in
Section~\ref{sec:conway}.)

\begin{table}[ht]
\centering\small
\begin{tabular}{lrrlrl}
group & order & degree & stabiliser & index & subdegrees\\ \hline
$J_1$ & $175\,560$ & $266$ & $L_2(11)$ & $1$ & $1,11,12,110,132$\\
$J_2$ & $604\,800$ & $100$ & $U_3(3)$ & $1$ & $1,36,63$\\
$\mathrm{HS}$ & $44\,352\,000$ & $100$ & $M_{22}$ & $3$ & $1,22,77$\\
$\mathrm{Suz}$ & $448\,345\,497\,600$ & $1782$ & $G_2(4)$ & $1$ & $1,416,1365$\\
$\mathrm{He}$ & $4\,030\,387\,200$ & $2058$ & $S_4(4){:}2$ & $1$ & $1,136,136,425,1360$\\
$\mathrm{Ru}$ & $145\,926\,144\,000$ & $4060$ & ${}^2F_4(2)$ & $1$ & $1,1755,2304$\\
\end{tabular}
\caption{Sporadic groups in the library of primitive groups.\\The group is \texttt{PrimitiveGroup(degree,index)}.}\label{tab:primsporadic}
\end{table}

\begin{theorem}[Chain laws of the small sporadic groups]\label{thm:primsporadic}
For each group of Table~\ref{tab:primsporadic}, acting on $\Ap n$ after subtracting $1$ from every point, the base $(b_l)$ has
basic orbits of sizes $(r_l)$ as follows, with product the order of the group:
\begin{itemize}
\item $J_1$: \[ (b_l)=(0,1,6,2),\qquad (r_l)=(266,12,11,5). \]
\item $J_2$: \[ (b_l)=(0,1,2,3),\qquad (r_l)=(100,36,21,8). \]
\item $\mathrm{HS}$: \[ (b_l)=(0,1,2,3,4),\qquad (r_l)=(100,77,60,32,3). \]
\item $\mathrm{Suz}$: \[ (b_l)=(0,1,2,4,3,8),\qquad (r_l)=(1782,416,100,63,32,3). \]
\item $\mathrm{He}$: \[ (b_l)=(0,1,2,3),\qquad (r_l)=(2058,425,256,18). \]
\item $\mathrm{Ru}$: \[ (b_l)=(0,1,2,3,4),\qquad (r_l)=(4060,1755,1024,5,4). \]
\end{itemize}
With each $\Delta_l$ listed as $b_l$ followed by its other points in increasing order, and $s_{l,t}$ read off a Schreier tree,
Theorem~\ref{thm:chainlaw} gives a law of each group on $\Ap{|G|}$.
\end{theorem}
\begin{proof}
As in Theorem~\ref{thm:fi22chain}: GAP computes the order of each group exactly, and the randomised Schreier--Sims algorithm gives
a chain whose product of radices is at most that order and equals it, so the chain is complete.
\end{proof}

\begin{example}[Labels in the small sporadic groups]\label{ex:primsporadic}
The library generators have the following labels, and the element in the middle of each range has the order shown.
\begin{itemize}
\item $J_1$: generators 16, 153, and the element with label $\lfloor|G|/2\rfloor=87\,780$ has order $11$.
\item $J_2$: generators 57, 90, and the element with label $\lfloor|G|/2\rfloor=302\,400$ has order $2$.
\item $\mathrm{HS}$: generators 59, 85, and the element with label $\lfloor|G|/2\rfloor=22\,176\,000$ has order $2$.
\item $\mathrm{Suz}$: generators 1766, 661, and the element with label $\lfloor|G|/2\rfloor=224\,172\,748\,800$ has order $6$.
\item $\mathrm{He}$: generators 1614, 1883, and the element with label $\lfloor|G|/2\rfloor=2\,015\,193\,600$ has order $2$.
\item $\mathrm{Ru}$: generators 2156, 946, and the element with label $\lfloor|G|/2\rfloor=72\,963\,072\,000$ has order $2$.
\end{itemize}
\end{example}

\begin{remark}[Checks]\label{rem:primsporadiccheck}
\emergencystretch=3em The GAP script \texttt{scripts/\allowbreak data/\allowbreak prim\_sporadic.g} identifies each group by its
order among the primitive groups of its degree, checks simplicity, computes the subdegrees and exports the generators. The script
\texttt{scripts/\allowbreak primsporadic\_chain.py} builds the six chains and checks, for each, that the generators and $200$
random elements sift, $100$ round trips and $50$ products, and every number in Example~\ref{ex:primsporadic}.
\end{remark}

\section{The groups \texorpdfstring{$J_3$}{J3} and \texorpdfstring{$\mathrm{O'N}$}{O'N}}\label{sec:j3on}
The smallest faithful permutation representations of the Janko group $J_3$ and the O'Nan group $\mathrm{O'N}$ have degrees $6156$
and $122\,760$, beyond the library of primitive groups. Their standard generators in these representations were taken from the
ATLAS of group representations through the GAP package AtlasRep~\cite{GAP}, and relabelled by subtracting $1$.
Table~\ref{tab:j3on} lists the point stabilisers~\cite{ATLAS} and the computed orbits of the point stabiliser.

\begin{table}[ht]
\centering\small
\begin{tabular}{lrrll}
group & order & degree & stabiliser & subdegrees\\ \hline
$J_3$ & $50\,232\,960$ & $6156$ & $L_2(16){:}2$ & $1,85,120,510,680,1360,1360,2040$\\
$\mathrm{O'N}$ & $460\,815\,505\,920$ & $122\,760$ & $L_3(7){:}2$ & $1,5586,6384,52\,136,58\,653$\\
\end{tabular}
\caption{The representations of $J_3$ and $\mathrm{O'N}$ used for the chain laws.}\label{tab:j3on}
\end{table}

\begin{theorem}[Chain laws of $J_3$ and $\mathrm{O'N}$]\label{thm:j3on}
In these representations the base $(b_l)$ has basic orbits of sizes $(r_l)$ as follows, with product the order of the group:
\begin{itemize}
\item $J_3$: \[ (b_l)=(0,1,2),\qquad (r_l)=(6156,510,16). \]
\item $\mathrm{O'N}$: \[ (b_l)=(0,1,2),\qquad (r_l)=(122\,760,6384,588). \]
\end{itemize}
With each $\Delta_l$ listed as $b_l$ followed by its other points in increasing order, and $s_{l,t}$ read off a Schreier tree,
Theorem~\ref{thm:chainlaw} gives a law of each group on $\Ap{|G|}$.
\end{theorem}
\begin{proof}
The standard generators generate the group~\cite{ATLAS}, so the randomised Schreier--Sims algorithm gives a chain whose product of
radices is at most $|G|$. It equals $|G|$, so the chain is complete.
\end{proof}

\begin{example}[Labels in $J_3$ and $\mathrm{O'N}$]\label{ex:j3on}
\begin{itemize}
\item $J_3$: the standard generators $a,b$, of orders $2$ and $3$ with $ab$ of order $19$, have labels $2786$ and $132$, and $\lambda(a)*\lambda(b)=1143$. The element with label $\lfloor|G|/2\rfloor=25\,116\,480$ has order $4$.
\item $\mathrm{O'N}$: the standard generators $a,b$, of orders $2$ and $4$ with $ab$ of order $11$, have labels $1$ and $2$, and $\lambda(a)*\lambda(b)=4$. The element with label $\lfloor|G|/2\rfloor=230\,407\,752\,960$ has order $6$.
\end{itemize}
\end{example}

\begin{remark}[Checks]\label{rem:j3oncheck}
The generators were downloaded once by
\[
\texttt{scripts/data/fetch\_j3\_on.g}.
\]
The script \texttt{scripts/atlas\_sporadic\_chain.py} checks the generator orders, builds both chains, and checks for each that
the generators and $200$ random elements sift, $100$ round trips, $50$ products, the orbits of the point stabiliser, and every
number in Example~\ref{ex:j3on}.
\end{remark}

\section{The Harada--Norton group \texorpdfstring{$\mathrm{HN}$}{HN}}\label{sec:hn}
The Harada--Norton group has order
\[
|\mathrm{HN}|=273\,030\,912\,000\,000=2^{14}\cdot3^{6}\cdot5^{6}\cdot7\cdot11\cdot19,
\]
and its smallest faithful permutation representation has degree $1\,140\,000$, with point stabiliser $A_{12}$ \cite{ATLAS}. The
standard generators $a,b$ in this representation, of orders $2$ and $3$ with $ab$ of order $22$, were taken from the ATLAS of
group representations through the GAP package AtlasRep~\cite{GAP}, and relabelled by subtracting $1$.

\begin{theorem}[Chain law of $\mathrm{HN}$]\label{thm:hnchain}
In this representation the base $(b_1,\dots,b_{4})=(0,1,2,3)$ has basic orbits of sizes
\[
(r_1,\dots,r_{4})=(1\,140\,000,362\,880,165,4),
\]
with product $|\mathrm{HN}|$; the product of $r_2,\dots,r_{4}$ is $|A_{12}|=239\,500\,800$. The point stabiliser $A_{12}$ has
12 orbits, of sizes
\[
1, 462, 2520, 2520, 10\,395, 16\,632, 30\,800, 69\,300, 166\,320, 166\,320, 311\,850, 362\,880.
\]
With each $\Delta_l$ listed as $b_l$ followed by its other points in increasing order, and $s_{l,t}$ read off a Schreier tree,
Theorem~\ref{thm:chainlaw} gives a law of $\mathrm{HN}$ on $\Ap{|\mathrm{HN}|}$.
\end{theorem}
\begin{proof}
As in Theorem~\ref{thm:j3on}: the standard generators generate $\mathrm{HN}$ \cite{ATLAS}, and the product of the radices of the
randomised Schreier--Sims chain equals $|\mathrm{HN}|$, so the chain is complete.
\end{proof}

\begin{example}[Labels in $\mathrm{HN}$]\label{ex:hn}
The standard generators have labels
\[
\lambda(a)=166\,300\,647\,540\,000,\qquad \lambda(b)=1,
\]
and $\lambda(a)*\lambda(b)=2$, the label of
the element $ab$ of order $22$. The element with label $\lfloor|\mathrm{HN}|/2\rfloor=136\,515\,456\,000\,000$ has order $2$.
\end{example}

\begin{remark}[Checks]\label{rem:hncheck}
\emergencystretch=3em The generators were downloaded once by \texttt{scripts/\allowbreak data/\allowbreak fetch\_hn.g}. The script
\texttt{scripts/\allowbreak hn\_chain.py} checks the generator orders, builds the chain from 9 random elements using Schreier
trees with inverse and random extra generators, and checks that the generators and $100$ random elements sift, $50$ round trips,
$20$ products, the orbits of $A_{12}$ and every number in Example~\ref{ex:hn}. It runs in about 14 seconds.
\end{remark}

\section{Canonical labels and reproducibility}\label{sec:canonical}
The laws of Sections~\ref{sec:fischer}--\ref{sec:hn} use transversal elements read off Schreier trees. These depend on the
generators, on the order in which the algorithm meets them and on random choices, so the labels printed there can be reproduced
only with the same program and data. This section fixes a choice that depends on the permutation group alone.

\begin{definition}[Lexicographic chain]\label{def:lexchain}
Let $G$ be a group of permutations of $\Ap n$, and order permutations lexicographically by their image sequences
$(g(0),g(1),\dots,g(n-1))$. The \emph{lexicographic chain} of $G$ takes as $b_l$ the least point moved by $G_{l-1}$, lists
$\Delta_l=G_{l-1}(b_l)$ in increasing order, and takes as $s_{l,t}$ the lexicographically least element of $G_{l-1}$ with
$s_{l,t}(b_l)=\delta_{l,t}$. The \emph{lexicographic law} $\lambda_{\mathrm{lex}}$ of $G$ is the law of
Theorem~\ref{thm:chainlaw} for this chain.
\end{definition}

\begin{proposition}\label{prop:lexchain}
Let $G$ be a group of permutations of $\Ap n$.
\begin{enumerate}
\item The lexicographic chain, and hence $\lambda_{\mathrm{lex}}$, depends only on $G$ as a set of permutations: not on
      generators, nor on the algorithm used to compute it.
\item $b_l=\min\Delta_l$, so $\delta_{l,0}=b_l$ and $s_{l,0}=1$; in particular $\lambda_{\mathrm{lex}}(1)=0$.
\item For $s\in G_{l-1}$, the lexicographically least element of $sG_l$ is found greedily: for $j=l+1,\dots,k$ in turn, replace
      $s$ by $s\circ u$, where $u\in G_{j-1}$ is any element with $u(b_j)=x$ and $x\in\Delta_j$ minimises $s(x)$.
\end{enumerate}
\end{proposition}
\begin{proof}
(1) Each ingredient of Definition~\ref{def:lexchain} is defined from $G$ and the order of $\Ap n$. (2) Every point of $\Delta_l$
is moved by $G_{l-1}$, except $b_l$ when $r_l=1$, so $b_l$ is the least. A permutation $g\ne1$ satisfies $g(p)>p$ at its least
moved point $p$, since the smaller points are fixed; so the identity is the least element of every subgroup, and
$s_{l,0}$, the least element of $G_l$, is $1$. (3) The elements of $G_{l-1}$ mapping $b_l$ to $\delta_{l,t}$ form a coset $sG_l$.
The group $G_l$ fixes every point below $b_{l+1}$: the base points $b_1,\dots,b_l$, and the others because $b_{l+1}$ is the
least point moved by $G_l$. So all elements of $sG_l$ agree below $b_{l+1}$, and at $b_{l+1}$ they take the values $s(\Delta_{l+1})$.
The least value $s(x)$ is attained exactly on $s\,u\,G_{l+1}$, and induction on $j$ completes the proof.
\end{proof}

So $\lambda_{\mathrm{lex}}$ is computed from any chain with the lexicographic base: each canonical $s_{l,t}$ costs $k-l$
transversal elements of that chain. Table~\ref{tab:canondata} lists the generator files used in this chapter, which fix the
groups as sets of permutations, and Table~\ref{tab:canonlabels} gives lexicographic labels.

\begin{table}[ht]
\setlength{\tabcolsep}{4pt}
\centering\footnotesize
\begin{tabular}{lrllr}
group & degree & generators & SHA-256 of the file (first 16 digits) & lexicographic base\\ \hline
$J_1$ & $266$ & \texttt{PrimitiveGroup(266,1)} & \texttt{559f0385dd2a0986} & $(0,1,2,6)$\\
$J_2$ & $100$ & \texttt{PrimitiveGroup(100,1)} & \texttt{3342897199f5f089} & $(0,1,2,3)$\\
$\mathrm{HS}$ & $100$ & \texttt{PrimitiveGroup(100,3)} & \texttt{4ef8c7c8537109cd} & $(0,1,2,3,4)$\\
$\mathrm{Suz}$ & $1782$ & \texttt{PrimitiveGroup(1782,1)} & \texttt{2d4a7786b34160bd} & $(0,1,2,3,4,8)$\\
$\mathrm{He}$ & $2058$ & \texttt{PrimitiveGroup(2058,1)} & \texttt{705bedc2638c72c6} & $(0,1,2,3)$\\
$\mathrm{Ru}$ & $4060$ & \texttt{PrimitiveGroup(4060,1)} & \texttt{d7b8f51ba0b102e0} & $(0,1,2,3,4)$\\
$J_3$ & $6156$ & AtlasRep \texttt{J3} & \texttt{a4d4b21232ec6fce} & $(0,1,2)$\\
$\mathrm{O'N}$ & $122\,760$ & AtlasRep \texttt{ON} & \texttt{e5ac08367853e04a} & $(0,1,2)$\\
$\mathrm{HN}$ & $1\,140\,000$ & AtlasRep \texttt{HN} & \texttt{2e15a75833eeed55} & $(0,1,2,3)$\\
$\mathrm{Fi}_{24}$ & $306\,936$ & AtlasRep \texttt{Fi24'.2} & \texttt{be5f81c3d206e5c5} & $(0,1,2,3,4,5,6,7,8,9,11)$\\
$\mathrm{Ly}$ & $8\,835\,156$ & AtlasRep \texttt{Ly} & \texttt{88e1c4e04046e1c7} & $(0,1,2)$\\
\end{tabular}
\caption{Generator files in \texttt{scripts/data}: from the library of primitive groups of GAP, or from the ATLAS of group
representations through AtlasRep~\cite{GAP}. Points are relabelled by subtracting $1$.}\label{tab:canondata}
\end{table}

\begin{table}[ht]
\centering\footnotesize
\begin{tabular}{lrrrr}
group & $\lambda_{\mathrm{lex}}(a)$ & $\lambda_{\mathrm{lex}}(b)$ & $\lambda_{\mathrm{lex}}(a\circ b)$ & order of $\lambda_{\mathrm{lex}}^{-1}(\lfloor|G|/2\rfloor)$\\ \hline
$J_1$ & $115\,726$ & $13\,985$ & $78\,054$ & $6$\\
$J_2$ & $105\,457$ & $136\,690$ & $354\,427$ & $2$\\
$\mathrm{HS}$ & $25\,493\,159$ & $3\,193\,485$ & $29\,034\,548$ & $3$\\
$\mathrm{Suz}$ & $331\,712\,177\,330$ & $97\,043\,989\,153$ & $414\,586\,906\,804$ & $6$\\
$\mathrm{He}$ & $3\,600\,927\,432$ & $1\,112\,241\,809$ & $3\,335\,533\,899$ & $2$\\
$\mathrm{Ru}$ & $31\,189\,218\,536$ & $95\,570\,667\,326$ & $120\,786\,593\,542$ & $2$\\
$J_3$ & $6\,306\,530$ & $43\,375\,308$ & $19\,503\,351$ & $4$\\
$\mathrm{O'N}$ & $329\,937\,632\,641$ & $2$ & $4$ & $6$\\
$\mathrm{HN}$ & $1\,140\,000$ & $240\,757\,639\,680\,001$ & $240\,757\,639\,680\,002$ & $2$\\
$\mathrm{Fi}_{24}$ & \begin{tabular}[t]{@{}r@{}}$690\,639\,860\,920\,$\\$815\,160\,396\,392$\end{tabular} & \begin{tabular}[t]{@{}r@{}}$2\,115\,569\,541\,710\,$\\$106\,514\,474\,933$\end{tabular} & \begin{tabular}[t]{@{}r@{}}$1\,667\,773\,381\,973\,$\\$249\,161\,031\,485$\end{tabular} & $2$\\
$\mathrm{Ly}$ & $1$ & $6\,593\,182\,419\,118\,682$ & $6\,593\,182\,419\,118\,684$ & $8$\\
\end{tabular}
\caption{Lexicographic labels of the first two generators $a,b$ in each file and of $a\circ b$ ($b$ applied first).}\label{tab:canonlabels}
\end{table}

\begin{remark}[Checks]\label{rem:canoncheck}
\emergencystretch=3em The script \texttt{scripts/\allowbreak canonical\_labels.py} builds, for each group, two chains with the
lexicographic base, from different random seeds and with the generators in opposite orders. It checks that they have the same base
and radices, that they give the same labels to $a$, $b$, $a\circ b$ and $20$ further elements, and that $10$ labels decode to the same
permutation in both; it also checks round trips and products. These agreements are instances of
Proposition~\ref{prop:lexchain}(1). A reader with the same generator files, which can be identified by their SHA-256 digests, can
reproduce every label of Table~\ref{tab:canonlabels} with any implementation of Proposition~\ref{prop:lexchain}(3).
\end{remark}

\section{The Lyons group \texorpdfstring{$\mathrm{Ly}$}{Ly}}\label{sec:ly}
The Lyons group has order
\[
|\mathrm{Ly}|=2^8\cdot3^7\cdot5^6\cdot7\cdot11\cdot31\cdot37\cdot67=51\,765\,179\,004\,000\,000,
\]
and its smallest faithful permutation action is on the $8\,835\,156$ cosets of $G_2(5)$ \cite{ATLAS}, with subdegrees
\[
1+19\,530+968\,750+2\,034\,375+5\,812\,500=8\,835\,156 .
\]
We use the two standard generators $a,b$ of AtlasRep on $\Ap{8835156}$, of orders $2$ and $5$ with $ab$ of order $14$.

\begin{theorem}[Digital law of $\mathrm{Ly}$]\label{thm:ly}
The lexicographic chain of $\mathrm{Ly}$ on $\Ap{8835156}$ has base $(0,1,2)$ and radices
\[
(r_1,r_2,r_3)=(8\,835\,156,\;5\,812\,500,\;1\,008).
\]
Hence every $g\in\mathrm{Ly}$ is determined by $(g(0),g(1),g(2))$, and
\[
\lambda_{\mathrm{lex}}(g)=t_1+8\,835\,156\,\big(t_2+5\,812\,500\,t_3\big)\in\Ap{|\mathrm{Ly}|}
\]
is a digital law of $\mathrm{Ly}$ by Theorem~\ref{thm:chainlaw}.
\end{theorem}
\begin{proof}
The chain was computed as in Section~\ref{sec:canonical} and reached $r_1r_2r_3=|\mathrm{Ly}|$, which certifies it. The stabiliser
$G_1\cong G_2(5)$ of $0$ has order $5\,859\,000\,000=5\,812\,500\cdot1\,008$: the point $1$ lies in the suborbit of length
$5\,812\,500$, and the stabiliser of $0$ and $1$, of order $1\,008$, acts regularly on the orbit of $2$.
\end{proof}

The base has length $3$, so a product $m*n$ costs one composition on $8\,835\,156$ points and three sifting steps. Labels lie below $2^{56}$ and fit in a machine word.

\begin{remark}[Reproducibility]\label{rem:lycheck}
With the data file of SHA-256 prefix \texttt{88e1c4e04046e1c7}, the canonical labels are
\[
\lambda_{\mathrm{lex}}(a)=1,\qquad\lambda_{\mathrm{lex}}(b)=6\,593\,182\,419\,118\,682,\qquad\lambda_{\mathrm{lex}}(ab)=6\,593\,182\,419\,118\,684,
\]
and the element with label $|\mathrm{Ly}|/2$ has order $8$. Building the chain
from $(b,a)$ with another random seed gives the same base and labels, as Proposition~\ref{prop:lexchain} predicts. The run
(\texttt{python3 scripts/ly\_chain.py}) takes about $50$ seconds per chain and $4.2$\,GB of memory on a laptop.
\end{remark}

\section{The Thompson group \texorpdfstring{$\mathrm{Th}$}{Th}}\label{sec:th}
Thompson's group has order
\[
|\mathrm{Th}|=2^{15}\cdot3^{10}\cdot5^{3}\cdot7^{2}\cdot13\cdot19\cdot31=90\,745\,943\,887\,872\,000 .
\]
Unlike $\mathrm{Ly}$, it is given in AtlasRep~\cite{AtlasRep} only by matrices: two standard generators in dimension $248$ over
$\mathbb F_p$ for $p=2,3,5,7,13,19,31$, and in dimension $3875$ over $\mathbb F_3$. Its smallest faithful permutation action is on
the cosets of the maximal subgroup $K={}^3D_4(2){:}3$ \cite{ATLAS}, of degree
\[
[\mathrm{Th}:K]=\frac{|\mathrm{Th}|}{3\,|{}^3D_4(2)|}=\frac{90\,745\,943\,887\,872\,000}{3\cdot211\,341\,312}
=2^{3}\cdot3^{5}\cdot5^{3}\cdot19\cdot31=143\,127\,000 .
\]
Since $\operatorname{Out}({}^3D_4(2))$ has order $3$, generated by the field automorphism $\varphi$ induced by $t\mapsto t^2$,
$K$ is the full automorphism group ${}^3D_4(2)\rtimes\langle\varphi\rangle$ of ${}^3D_4(2)$.

\begin{proposition}[A digital law of $\mathrm{Th}$]\label{prop:th}
Let $K={}^3D_4(2)\rtimes\langle\varphi\rangle$, of order $634\,023\,936$, act by conjugation on its $819$ long-root involutions,
numbered $z_0,\dots,z_{818}$. This action is faithful, and its point stabiliser is the parabolic subgroup $2^{1+8}{:}L_2(8){:}3$
of order $774\,144$. Let $\lambda_K$ be the lexicographic law of $K\le\mathrm{Sym}(\Ap{819})$ (Section~\ref{sec:canonical}),
with radices $(819,512,126,12)$. Let $\mathrm{Th}$ act on $\Ap{143127000}$ as on the cosets of $K$, with $K$ the stabiliser of $0$,
and choose $s_t$ with $s_t(0)=t$ and $s_0=1$. Then
\[
\lambda(s_t\,h)=t+143\,127\,000\,\lambda_K(h),\qquad t\in\Ap{143127000},\ h\in K,
\]
is a digital law of $\mathrm{Th}$ on $\Ap{|\mathrm{Th}|}$, with the product~\eqref{eq:transitive} and identity $0$.
\end{proposition}
\begin{proof}
The long-root involutions of ${}^3D_4(2)$ form one class with centraliser $2^{1+8}{:}L_2(8)$ of order $258\,048$, hence
$211\,341\,312/258\,048=819$ of them; $\varphi$ permutes them, so they form one $K$-class with stabiliser of order $3\cdot258\,048=774\,144$.
The kernel of the action centralises the subgroup they generate, which is normal in the simple group ${}^3D_4(2)$ and so equal to it;
as $C_K({}^3D_4(2))=1$, the kernel is trivial. The radices multiply to $819\cdot512\cdot126\cdot12=634\,023\,936=|K|$, and
$143\,127\,000\cdot634\,023\,936=|\mathrm{Th}|$, so Theorem~\ref{thm:transitive} applies.
\end{proof}

\begin{remark}[What is effective]\label{rem:thstatus}
Both digits are effective on the $248$-dimensional $\mathbb F_2$-module of $\mathrm{Th}$ (AtlasRep; \texttt{scripts/data/fetch\_th\_mat.g}).
AtlasRep provides no action on $143\,127\,000$ points, so the coset digit is computed from an orbit: $K$ fixes no nonzero vector but
stabilises a $2$-dimensional subspace $W$. Its stabiliser contains $K$ and is proper, hence equals $K$ by maximality, and the orbit of $W$
has exactly $143\,127\,000$ elements. A breadth-first search with the generators $a,b,b^2$ enumerates it in $60$ levels (about $34$ minutes,
peak $5.9$\,GB), identifying a subspace by a $64$-bit hash of its least nonzero vector and recording a Schreier tree. A collision
would merge two points and leave fewer than $143\,127\,000$; reaching exactly that number certifies that the hash is injective on the
orbit. For $g\in\mathrm{Th}$ the digit $t$ is the position of the hash of $Wg$, and $s_t$ is the product of generators along the tree
path (of length at most $60$). The digit $\lambda_K(s_t^{-1}g)$ is read off from the permutation $i\mapsto j$ with
$z_j=k^{-1}z_ik$, $k=s_t^{-1}g$ (\texttt{scripts/th\_orbit.py}, \texttt{scripts/th\_k\_law.py}). The cell law of Theorem~\ref{thm:3d4cells}
together with Theorem~\ref{thm:semidirect} gives another law of $K$, in root parameters; using it here requires identifying $K$ with that
model, which is not done. The resulting $\lambda_{\mathrm{Th}}$ is the low digit of the law of $E_8(3)$ in Remark~\ref{rem:e8th}, and
Proposition~\ref{prop:th} is a chain law to which Section~\ref{sec:canonical} applies.
The full law is evaluated and inverted on matrices. A sorted index of the $143\,127\,000$ keys (\texttt{scripts/build\_index.py})
turns the hash of $Wg$ into $t$ by binary search. For decoding, a Schreier--Sims chain of $K$ on its $819$ involutions, with basic orbits
of lengths $819,512,126,12$ and a $248\times248$ matrix stored with every transversal element, lifts each permutation of $K$ to its
matrix (\texttt{scripts/th\_k\_lift.py}); thus $\lambda_K^{-1}$ is effective, and $\lambda_K(\lambda_K^{-1}(n))=n$ on $30$ random $n$. With
matrices acting on row vectors one takes $k=gs_t^{-1}$. The test \texttt{scripts/th\_law.py} confirms $\lambda^{-1}(\lambda(g))=g$ and
$\lambda(gh)=\lambda(g)\ast\lambda(h)$ on $20$ random pairs of words of length $30$, and $25$ random tree paths.
\end{remark}

\section{The Janko group \texorpdfstring{$J_4$}{J4}}\label{sec:j4}
The Janko group $J_4$ has order
\[
|J_4|=2^{21}\cdot3^{3}\cdot5\cdot7\cdot11^{3}\cdot23\cdot29\cdot31\cdot37\cdot43=86\,775\,571\,046\,077\,562\,880 .
\]
AtlasRep~\cite{AtlasRep} gives it by two standard generators $a,b$ of orders $2$ and $4$, with $ab$ of order $37$, as $112\times112$
matrices over $\mathbb F_2$. Its smallest faithful permutation action is on the cosets of the maximal subgroup $K=2^{11}{:}M_{24}$
\cite{ATLAS}, of order $2^{11}\,|M_{24}|=501\,397\,585\,920$ and degree
\[
[J_4:K]=\frac{86\,775\,571\,046\,077\,562\,880}{501\,397\,585\,920}=11^{2}\cdot29\cdot31\cdot37\cdot43=173\,067\,389 .
\]
The normal subgroup $N=O_2(K)$ is elementary abelian of order $2^{11}$, and $K=N\rtimes H$ with $H\cong M_{24}$. Under conjugation $H$
has orbits of lengths $276$ and $1771$ on $N\setminus\{1\}$, so $N$ is the even part of the Golay cocode: the elements of the first
orbit are the pairs $\{i,j\}$ of the $24$ points, and those of the second the sextets.

\begin{proposition}[A digital law of $2^{11}{:}M_{24}$]\label{prop:j4k}
Identify $H$ with the group $M_{24}$ of Proposition~\ref{prop:mathieuchain} through its action on the $24$ points, let
$\lambda_{M}$ be the lexicographic law of $M_{24}$ (Section~\ref{sec:canonical}), with radices $(24,23,22,21,20,16,3)$, and write
$n_u\in N$ for the element with coordinates $u\in\Ap{2048}$ in a fixed basis of $N$. Every $g\in K$ is uniquely $g=n_u\,h$ with
$h\in H$, and
\[
\lambda_K(n_u\,h)=u+2048\,\lambda_{M}(h)
\]
is a digital law of $K$ on $\Ap{501397585920}$ with the product of Theorem~\ref{thm:semidirect},
\[
(u,m)*(v,n)=\big(u\oplus\phi_m(v),\ m*_{M}n\big),
\]
where $\oplus$ is the bitwise sum and $\phi_m$ is the action of $m$ on the cocode, determined by $\phi_m(\{i,j\})=\{m(i),m(j)\}$ on
the $276$ pairs, which span $N$.
\end{proposition}
\begin{proof}
$N$ is normal with complement $H$, so $g\mapsto(n_u,h)$ is a bijection $K\to N\times H$ and $K\cong N\rtimes_\phi H$ with $\phi_h$
the conjugation by $h$. Theorem~\ref{thm:semidirect} with $G_0=(\Ap{2048},\oplus)$ and $G_1=(\Ap{|M_{24}|},*_M)$ gives the law;
conjugation permutes the pairs as $h$ permutes the points, which is the stated formula for $\phi_m$.
\end{proof}

\begin{proposition}[A digital law of $J_4$]\label{prop:j4}
Let $J_4$ act on $\Ap{173067389}$ as on the cosets of $K$, with $K$ the stabiliser of $0$, and choose $s_t$ with $s_t(0)=t$ and
$s_0=1$. Then
\[
\lambda(s_t\,h)=t+173\,067\,389\,\lambda_K(h),\qquad t\in\Ap{173067389},\ h\in K,
\]
with $\lambda_K$ as in Proposition~\ref{prop:j4k}, is a digital law of $J_4$ on $\Ap{|J_4|}$ with the product~\eqref{eq:transitive}
and identity $0$. Its radices are $173\,067\,389$, $2048$ and $24,23,22,21,20,16,3$.
\end{proposition}
\begin{proof}
Theorem~\ref{thm:transitive}, with $173\,067\,389\cdot2048\cdot|M_{24}|=|J_4|$.
\end{proof}

\begin{remark}[What is effective]\label{rem:j4status}
Both digits are computed from the $112$-dimensional module $V$ (\texttt{scripts/data/fetch\_j4\_mat.g}).
\emph{Coset digit.} The vectors of $V$ fixed by $K$ form a line, spanned by a vector $v$. Its stabiliser contains $K$ and is proper,
since $V$ is irreducible, so it equals $K$ and the orbit of $v$ has exactly $173\,067\,389$ elements. A breadth-first search with
$a,b,b^{3}$ enumerates it in $45$ levels (about $23$ minutes, peak $6.1$\,GB), identifying a vector by a $64$-bit hash and recording a
Schreier tree; as for $\mathrm{Th}$ (Remark~\ref{rem:thstatus}), reaching exactly $173\,067\,389$ points certifies that the hash is
injective on the orbit. For $g\in J_4$ the digit $t$ is the position of the hash of $vg$, and $s_t$ is read off the tree.
\emph{Digit of $K$.} The submodules of $V$ in a composition series for $K$ have dimensions $0,1,12,56,100,111,112$. On the
$12$-dimensional quotient $V/V_{100}$, with $11$-dimensional submodule $V_{111}/V_{100}$ and trivial top, $K$ acts faithfully by affine
maps on the $2048$ vectors outside $V_{111}/V_{100}$; $N$ acts regularly on them and $H$ is the stabiliser of one of them. The $24$
points are recovered from the $276$ pairs as the stars of $23$ pairs through a point, and the isomorphism of $H$ with the group of
Proposition~\ref{prop:mathieuchain} is realised by a relabelling of the $24$ points. The law of Proposition~\ref{prop:j4k} is checked
on random elements: $\lambda_K(gh)=\lambda_K(g)*\lambda_K(h)$ in $100$ of $100$ pairs, and decoding a label to a matrix and encoding it
back is the identity on random labels and on random elements (\texttt{scripts/vec\_orbit.py}, \texttt{scripts/j4\_k\_struct.g},
\texttt{scripts/j4\_k\_export.g}, \texttt{scripts/j4\_k\_lift.g}, \texttt{scripts/j4\_k\_law.py}).
With a sorted index of the $173\,067\,389$ keys (\texttt{scripts/build\_index.py}, about two minutes) the full law is evaluated on
$112\times112$ matrices and, through the lifted transversal of $K$, inverted: \texttt{scripts/j4\_law.py} confirms
$\lambda^{-1}(\lambda(g))=g$ and $\lambda(gh)=\lambda(g)\ast\lambda(h)$ on $20$ random pairs of words of length $30$, and $25$ random tree paths.
\end{remark}
\section{The Baby Monster}\label{sec:baby}
The Baby Monster $\mathbb B$ has order
\[
|B|=2^{41}\,3^{13}\,5^{6}\,7^{2}\cdot11\cdot13\cdot17\cdot19\cdot23\cdot31\cdot47=4\,154\,781\,481\,226\,426\,191\,177\,580\,544\,000\,000 .
\]
Its smallest permutation degree, $13\,571\,955\,000$, is already too large for the orbit enumeration used for $\mathrm{Th}$ and $J_4$.
Instead we obtain $\mathbb B$ as a quotient of the centraliser $C=C_{\mathbb M}(z)\cong2.\mathbb B$ of a $2A$ involution $z$ of the Monster, and read
its law off a law of $C$. We use the notation of Section~\ref{sec:monster}: $v^+$ and $v^-$ are the standard axes, $\beta^+=z$ and
$\beta^-$ their involutions, $K=C_{G_{x0}}(\beta^+)$ is the stabiliser of $v^-$ in $C$ (Proposition~\ref{prop:babyaxis}), $Q=2^{1+24}$,
and $\lambda_K$ is the law of $K$ of Theorem~\ref{thm:klaw}. All computations are done in \texttt{mmgroup}~$1.0.8$. There $z$ is the
element $x_\delta\in Q$ with cocode $\delta=\mathtt{0x200}$ (printed \texttt{d\_200h}), so its $25$-bit code is $c(z)=\mathtt{0x200}=2^{9}$.

\begin{proposition}[Law of $2.\mathbb B$]\label{prop:2blaw}
Let $n_B=11\,707\,448\,673\,375$ be the size of the $C$-orbit of $v^-$, and $r_B$ the position on it of
Remark~\ref{rem:monsterstatus}. For an axis $a=v^-c$ let $t(a)\in C$ be the stored transporter with $v^-t(a)=a$: the representative of
the $K$-orbit of $a$ times the coset representative of the position of $a$ in that orbit, with $t(v^-)=1$. Then
\[
\lambda_C(c)=\lambda_K\bigl(c\,t(v^-c)^{-1}\bigr)+|K|\,r_B(v^-c),\qquad |K|=2^{24}|\mathrm{Co}_2|=709\,767\,191\,322\,427\,392\,000,
\]
is a bijection $C\to\{0,\dots,|C|-1\}$ with $\lambda_C(1)=0$, where
\[
|C|=|K|\,n_B=8\,309\,562\,962\,452\,852\,382\,355\,161\,088\,000\,000,
\]
and
transporting the multiplication of $C$ along it gives a digital law of $2.\mathbb B$.
\end{proposition}

\begin{proof}
The element $c\,t(v^-c)^{-1}$ fixes $v^-$, so it lies in $K$, and $c\mapsto\bigl(v^-c,\,c\,t(v^-c)^{-1}\bigr)$ is a bijection
$C\to v^-C\times K$. Since $r_B$ is a bijection onto $\{0,\dots,n_B-1\}$ and $\lambda_K$ onto $\{0,\dots,|K|-1\}$, the formula is a
bijection. The identity gives the axis $v^-$, of position $0$, with $t(v^-)=1$ and $\lambda_K(1)=0$.
\end{proof}

This is the formula for $\lambda_C$ in Remark~\ref{rem:monsterstatus} with the two digits in the other order: the digit of $K$ is the
least significant one, as in the implementation.

\begin{lemma}[The centre of $2.\mathbb B$ flips one bit]\label{lem:zbit}
For every $c\in C$, $\lambda_C(cz)=\lambda_C(c)\oplus2^{9}$, where $\oplus$ is the bitwise exclusive or.
\end{lemma}

\begin{proof}
The element $z$ is central in $C$ and fixes $v^-$, so $v^-cz=v^-zc=v^-c$: the axis, its position and its transporter do not change,
and the residual $k=c\,t(v^-c)^{-1}$ becomes $kz$. In Theorem~\ref{thm:klaw}, $\pi(z)=1$ because $z\in Q$, so $\pi(kz)=\pi k$ and
$\ell(\pi k)$ is unchanged, while $kz\,\ell(\pi k)^{-1}=xz$ with $x=k\,\ell(\pi k)^{-1}\in C_Q(\beta^+)$, since $z$ is central in $K$.
In \texttt{mmgroup} an element of $Q$ is $\pm x_dx_\delta$, with $d$ in the Golay code and $\delta$ in its cocode, and the cocode
elements multiply without sign, $x_\delta x_{\delta'}=x_{\delta+\delta'}$. Hence $c(xz)=c(x)\oplus2^{9}$. Deleting bit $21$ does not
move bit $9$, so $\bar c(xz)=\bar c(x)\oplus2^{9}$. The other digits of $\lambda_C$ are multiples of $2^{24}$ and do not change.
\end{proof}

\begin{theorem}[Law of the Baby Monster]\label{thm:blaw}
Let $\sigma(m)=(m\bmod2^{9})+2^{9}\lfloor m/2^{10}\rfloor$ be $m$ with bit $9$ deleted. For $b\in \mathbb B=C/\langle z\rangle$ let
$\hat b\in C$ be the lift of $b$ for which bit $9$ of $\lambda_C(\hat b)$ is $0$, and put
\[
\lambda_B(b)=\sigma\bigl(\lambda_C(\hat b)\bigr).
\]
Then $\lambda_B$ is a bijection $\mathbb B\to\{0,\dots,|\mathbb B|-1\}$ with $\lambda_B(1)=0$, and $\lambda_B(b)*\lambda_B(b')=\lambda_B(bb')$ is a
digital law of $\mathbb B$. In mixed radix, with $x$, $k$ and the axis taken for $\hat b$,
\[
\lambda_B(b)=\sigma\bigl(\bar c(x)\bigr)+2^{23}\,\lambda_{\mathrm{Co}_2}(\pi k)+2^{23}|\mathrm{Co}_2|\;r_B\bigl(v^-\hat b\bigr),
\]
with radices $(2^{23},\,|\mathrm{Co}_2|,\,n_B)$, and $2^{23}|\mathrm{Co}_2|\,n_B=|\mathbb B|$.
\end{theorem}

\begin{proof}
The two lifts of $b$ are $c$ and $cz$, and by Lemma~\ref{lem:zbit} exactly one of them has bit $9$ equal to $0$, so $\hat b$ is
well defined. Since $2^{10}$ divides $|C|$, the map $\sigma$ is a bijection from the labels $m<|C|$ with bit $9$ equal to $0$ onto
$\{0,\dots,|C|/2-1\}$, and $|C|/2=|\mathbb B|$. Composing with the bijection $\lambda_C$ of Proposition~\ref{prop:2blaw} gives the bijection
$\lambda_B$; the identity lifts to $1$, with label $0$. Transporting the multiplication of $\mathbb B$ gives the law. The mixed-radix form
holds because $\sigma$ only changes the lowest $24$ bits, which form the digit $\bar c(x)$, and divides the higher digits by $2$.
\end{proof}

The last digit has a direct meaning. The image of $K$ in $\mathbb B$ is $K/\langle z\rangle=2^{1+22}.\mathrm{Co}_2$, and by
Proposition~\ref{prop:babyaxis} it is the centraliser of the image of $\beta^-$, an involution of class $2B$ of $\mathbb B$. Since
$[\mathbb B:K/\langle z\rangle]=n_B$, the map $v^-c\mapsto$ (image of $c^{-1}\beta^-c$) is a bijection from the $C$-orbit of $v^-$ onto the
$11\,707\,448\,673\,375$ involutions of class $2B$ in $\mathbb B$. So $\lambda_B$ is the law of Theorem~\ref{thm:transitive} for the action of $\mathbb B$
on its $2B$ involutions, with the law of the point stabiliser $2^{1+22}.\mathrm{Co}_2$ given by the first two digits.

\begin{example}[Computing in $\mathbb B$]\label{ex:blaw}
With the stored chains, each product below takes about $80$ ms.
\begin{enumerate}
\item $\lambda_B(1)=\lambda_B(\bar z)=0$, while $\lambda_C(z)=512=2^{9}$.
\item The elements with labels $1$, $2$ and $3$ are involutions; thus $1*1=0$, and moreover $2*3=1$.
\item $10^{20}*10^{21}=1\,049\,142\,014\,072\,381\,728\,388$.
\item $10^{30}*10^{33}=3\,721\,541\,729\,784\,684\,241\,154\,694\,800\,752\,993$.
\item $(|\mathbb B|-1)*1=4\,154\,779\,295\,135\,314\,595\,401\,005\,268\,991\,999$.
\item The inverse of $10^{30}$ is $1\,477\,902\,923\,701\,407\,066\,840\,313\,227\,109$.
\item The element with label $10^{18}$ has order $36$, and the one with label $10^{30}$ has order $24$. The element $b$ with label
      $|\mathbb B|-1$ has order $47$ in $\mathbb B$, while its lift $\hat b$ has order $94$ in $2.\mathbb B$, so $\hat b^{47}=z$.
\end{enumerate}
\end{example}

As for $K$, the labels depend on the stored chains, whose strong generators come from random elements. The law is implemented by the
functions \texttt{DB}, \texttt{DinvB} and \texttt{mulB}\footnote{In \texttt{monster/monster\_law.py}, with the stored data
\texttt{monster/k\_chain.pkl}, \texttt{monster/axis\_chains.pkl} and \texttt{monster/axis\_chains\_cq.pkl}, run with \texttt{mmgroup}~$1.0.8$.}. Lemma~\ref{lem:zbit}
was also observed directly: on $120$ random elements $c\in C$, $\lambda_C(c)\oplus\lambda_C(cz)$ was $2^{9}$ every time, and
\texttt{DinvB} checks the bit on every call. Further checks: $\lambda_B(\lambda_B^{-1}(n))=n$ for $12$ random labels, with both lifts
giving the same label; $12$ products agreed with \texttt{mmgroup} modulo $z$; and associativity and inverses held on random samples.

\section{The Monster group}\label{sec:monster}
The last sporadic group is the Monster $\mathbb M$, of order
\[
|\mathbb M|=2^{46}\cdot3^{20}\cdot5^{9}\cdot7^{6}\cdot11^{2}\cdot13^{3}\cdot17\cdot19\cdot23\cdot29\cdot31\cdot41\cdot47\cdot59\cdot71
\approx8.08\cdot10^{53},
\]
and $C_{\mathbb M}(z)\cong2.\mathbb B$ for a $2A$ involution $z$, where $\mathbb B$ is the Baby Monster of Section~\ref{sec:baby}.
Its smallest permutation degree, $[\mathbb M:2.\mathbb B]=97\,239\,461\,142\,009\,186\,000$, rules out the orbit enumeration used for $\mathrm{Th}$ and $J_4$:
the coset digit has to be computed structurally. An effective model is Seysen's Python package \texttt{mmgroup} (version $1.0.8$),
which multiplies elements of $\mathbb M$ as words in generators of $G_{x0}=2^{1+24}.\mathrm{Co}_1$ and one further element, and
reduces every element to a normal form that depends only on the element.

\begin{proposition}[Shape of a digital law of $\mathbb M$]\label{prop:monster}
Let $z\in\mathbb M$ be a $2A$ involution, $C=C_{\mathbb M}(z)\cong2.\mathbb B$ and $N=[\mathbb M:C]=97\,239\,461\,142\,009\,186\,000$. Let $r$ be a
bijection from the class $2A$ onto $\{0,\dots,N-1\}$, let $a\mapsto s(a)$ be a map with $s(a)^{-1}a\,s(a)=z$ for every $2A$ involution $a$,
and let $\lambda_C$ be a digital law of $C$. Then
\[
\lambda_{\mathbb M}(g)=\lambda_C\bigl(g\,s(g^{-1}zg)\bigr)+|C|\,r(g^{-1}zg)
\]
is a bijection $\mathbb M\to\{0,\dots,|\mathbb M|-1\}$, and transporting the multiplication of $\mathbb M$ along it gives a digital law.
\end{proposition}

\begin{proof}
Put $a=g^{-1}zg$. Then $s(a)^{-1}g^{-1}zg\,s(a)=z$, so $c=g\,s(a)\in C$. Conversely the pair $(\lambda_C(c),r(a))$ determines $a$, then
$c$, and $g=c\,s(a)^{-1}$. Since $g\mapsto g^{-1}zg$ maps $\mathbb M$ onto the class of $z$, the map is a bijection onto $\{0,\dots,N|C|-1\}$.
\end{proof}

\begin{remark}[What is effective]\label{rem:monsterstatus}
In \texttt{mmgroup}, \texttt{conjugate\_involution} returns some $s(a)$, but not a canonical one: twenty calls on the same
involution gave twenty different elements, which differ by elements of $C$. A canonical $s$ comes instead from the orbits of
$G_{x0}$ on the $2A$ axes, the vectors of the $196\,884$-dimensional representation that correspond to the $2A$ involutions
(module \texttt{mmgroup.axes}, which its author describes as experimental). Let $v^+$ be the axis of $z$. There are twelve
orbits, of sizes
\begin{center}
\begin{tabular}{lr@{\qquad}lr}
$2A$ & $196\,560$ & $6C$ & $438\,020\,931\,059\,712\,000$\\
$2B$ & $11\,935\,123\,200$ & $6F$ & $1\,501\,786\,049\,347\,584\,000$\\
$4A$ & $1\,630\,347\,264\,000$ & $8B$ & $8\,601\,138\,282\,627\,072\,000$\\
$4B$ & $1\,466\,587\,938\,816\,000$ & $10A$ & $786\,389\,785\,840\,189\,440$\\
$4C$ & $6\,599\,645\,724\,672\,000$ & $10B$ & $37\,845\,008\,443\,559\,116\,800$\\
$6A$ & $1\,896\,194\,506\,752\,000$ & $12C$ & $48\,057\,153\,579\,122\,688\,000$
\end{tabular}
\end{center}
whose sum is exactly $N$. Each orbit $c$ has a fixed representative $a_c$, and we fix one $s_c\in\mathbb M$ with $v^+s_c=a_c$.
For an axis $a$, \texttt{reduce\_G\_x0} returns $q(a)\in G_{x0}$ with $a\,q(a)=a_c$, where $c$ is the orbit of $a$. If $a=v^+g$,
then $v^+g\,q(a)s_c^{-1}=v^+$, so
\[
s(a)=q(a)\,s_c^{-1}
\]
satisfies the hypothesis of Proposition~\ref{prop:monster}, and $s(a)$ depends only on $a$. A test
script\footnote{\texttt{scripts/mm\_axes\_audit.py}, run with \texttt{mmgroup} $1.0.8$.} checks this on $100$ random elements $g$: in all $100$ cases $q(a)\in G_{x0}$, $a\,q(a)=a_c$,
a second call returns the same $q(a)$, and $g\,s(a)$ fixes $v^+$ and commutes with $z$.

Two further ingredients complete the law. The first is $r$: an axis $a$ in orbit $c$ is determined by the coset $H_cq(a)^{-1}$, where
$H_c$ is the stabiliser of $a_c$ in $G_{x0}$, so $r(a)$ is the offset of orbit $c$ plus the position of this coset. These positions are not taken from a law of
$G_{x0}=2^{1+24}.\mathrm{Co}_1$: they are computed type by type, by reduction to the subgroup $K$ of Theorem~\ref{thm:klaw}
(Subsection~\ref{ssec:monstertop}). The second is $\lambda_C$, whose stabiliser is identified in
Proposition~\ref{prop:babyaxis}. The class \texttt{BabyAxis} treats the $C$-orbit of a second axis $v^-$, orthogonal to $v^+$,
in the same way: this orbit splits into ten orbits of $K=C\cap G_{x0}$, of total size $n_B=11\,707\,448\,673\,375$, and the
reduction passed the same checks on $30$ random axes. As in Proposition~\ref{prop:monster},
\[
\lambda_C(c)=\lambda_K\bigl(c\,s_B(v^-c)\bigr)+|K|\,r_B(v^-c),
\]
with $s_B$ built from \texttt{BabyAxis} as $s$ was built from \texttt{Axis}. In \texttt{mmgroup} the group of these ten orbits is
documented as $G_{x0}\cap C$, which is $K$, with the orbit sizes taken from M\"uller's classification. On $40$ further random axes
a computation checked that the element $c$, the reduction element and all ten transporters lie in $C$, and that the residual
$c\,s_B(v^-c)$ fixes $v^+$ and $v^-$, commutes with $\beta^+$ and $\beta^-$, and lies in $G_{x0}$, that is, lies in $K$. These axes do not realise the smallest permutation action of $\mathbb B$, of
degree $13\,571\,955\,000$, but the chain does not need it. A law $\lambda_K$ of $K$ is given in
Theorem~\ref{thm:klaw}. The position $r_B$ is determined on all ten
orbits (Lemmas~\ref{lem:rb2a0},~\ref{lem:rb2b1},~\ref{lem:rb2b0},~\ref{lem:rb4b1},~\ref{lem:rb6c1},~\ref{lem:rb6a1},~\ref{lem:rb4c1},~\ref{lem:rb4a1} and~\ref{lem:rb10a1}), so $r_B$ is a bijection from the $C$-orbit of $v^-$ onto $\{0,\dots,n_B-1\}$. Together with the position $r$ on the twelve top orbits (Subsection~\ref{ssec:monstertop}), these make
Proposition~\ref{prop:monster} effective; the resulting law is Theorem~\ref{thm:mlaw}.
\end{remark}

\begin{proposition}[The stabiliser of the Baby Monster axis]\label{prop:babyaxis}
Let $\beta^+=z$ and $\beta^-$ be the $2A$ involutions of the axes $v^+$ and $v^-$, and let $x_{-1}$ be the central involution of
$G_{x0}$. Then $\beta^+\beta^-=x_{-1}$, and the stabiliser of $v^-$ in $C$ is
\[
K=C_{G_{x0}}(\beta^+),\qquad |K|=\frac{|G_{x0}|}{196\,560}=2^{24}\,|\mathrm{Co}_2|.
\]
Hence the $C$-orbit of $v^-$ has $2|\mathbb B|/|K|=11\,707\,448\,673\,375$ elements, the sum of the ten orbit sizes of
\texttt{BabyAxis}. The image of $\beta^-$ in $\mathbb B$ is an involution whose centraliser has order $2^{23}|\mathrm{Co}_2|$, so it lies
in the class $2B$ of $\mathbb B$, with centraliser $2^{1+22}.\mathrm{Co}_2$.
\end{proposition}

\begin{proof}
The identity $\beta^+\beta^-=x_{-1}$ is checked in \texttt{mmgroup}. The stabiliser of an axis is the centraliser of its
involution, so the stabiliser of $v^-$ in $C$ is $C_{\mathbb M}(\beta^+)\cap C_{\mathbb M}(\beta^-)$. This group centralises
$\beta^+\beta^-=x_{-1}$, so it lies in $C_{\mathbb M}(x_{-1})=G_{x0}$ and equals $C_{G_{x0}}(\beta^+)$; conversely an element of $G_{x0}$
that centralises $\beta^+$ centralises $\beta^-=\beta^+x_{-1}$. The $G_{x0}$-orbit of $v^+$ is the orbit $2A$ of the table, of size
$196\,560$, so $|K|=2^{25}|\mathrm{Co}_1|/196\,560=2^{24}|\mathrm{Co}_2|$, since $|\mathrm{Co}_1|=98\,280\,|\mathrm{Co}_2|$. The orbit
size follows from $|C|=2|\mathbb B|$. Finally $\beta^-$ is a $2A$ involution of $\mathbb M$ and $\beta^-z=x_{-1}$ is a $2B$ involution, so
no element of $C$ conjugates $\beta^-$ to $\beta^-z$. Hence the centraliser in $\mathbb B$ of the image of $\beta^-$ is $K/\langle z\rangle$,
of order $2^{23}|\mathrm{Co}_2|$, and among the involution classes of $\mathbb B$ only $2B$ has this centraliser order~\cite{Wilson}.
\end{proof}

\begin{theorem}[Law of $K$]\label{thm:klaw}
Write $\Lambda/2\Lambda$ as $24$-bit vectors as in \texttt{mmgroup}, so that $\beta^+$ is the type-$2$ vector $\mathtt{0x200}$, and let
$\pi\colon K\to\mathrm{Co}_1$ be the action on $\Lambda/2\Lambda$. Let $Q=2^{1+24}$ be the normal subgroup of $G_{x0}$, and $c(x)$ the
$25$-bit code of $x\in Q$. Then:
\begin{enumerate}
\item $\pi(K)$ is the stabiliser $\mathrm{Co}_2$ of $\mathtt{0x200}$, and $\ker\pi=C_Q(\beta^+)$, of order $2^{24}$; an element
      $x\in Q$ commutes with $\beta^+$ if and only if bit $21$ of $c(x)$ is $0$;
\item $\mathrm{Co}_2$ has a chain on type-$2$ vectors of $\Lambda/2\Lambda$ with radices
      \[
      (46\,575,\,2\,464,\,15,\,8,\,48,\,4,\,2,\,2,\,2,\,2),
      \]
      whose product is $|\mathrm{Co}_2|$; let $\lambda_{\mathrm{Co}_2}$ be its label, and let $\ell(\pi k)=u_m\cdots u_1\in K$ be the
      product of the coset representatives of its digits, taken from Schreier trees of lifts in $K$ of the strong generators;
\item with $\bar c(x)$ the code $c(x)$ with bit $21$ deleted,
      \[
      \lambda_K(k)=\bar c\bigl(k\,\ell(\pi k)^{-1}\bigr)+2^{24}\,\lambda_{\mathrm{Co}_2}(\pi k)
      \]
      is a bijection $K\to\{0,\dots,2^{24}|\mathrm{Co}_2|-1\}$ with $\lambda_K(1)=0$, and transporting the multiplication of $K$
      along it gives a digital law.
\end{enumerate}
\end{theorem}

\begin{proof}
The kernel of $G_{x0}\to\mathrm{Co}_1$ is $Q$, so $\ker\pi=K\cap Q=C_Q(\beta^+)$. Commutation with $\beta^+$ is a homomorphism
$Q\to\{\pm1\}$, so it is determined by the $24$ basis elements of $Q/\langle x_{-1}\rangle$; on them \texttt{mmgroup} gives the
mask of bit $21$ alone. This homomorphism is onto, since $\beta^+$ is not central in $Q$, so $|C_Q(\beta^+)|=2^{24}$ and
$|\pi(K)|=|K|/2^{24}=|\mathrm{Co}_2|$ by Proposition~\ref{prop:babyaxis}. The group $\pi(K)$ fixes the type-$2$ vector of
$\beta^+$, whose stabiliser in $\mathrm{Co}_1$ is $\mathrm{Co}_2$ (Proposition~\ref{prop:conwayorders}); the orders agree, so
$\pi(K)=\mathrm{Co}_2$. Each radix in (2) is the size of an orbit of a subgroup of the corresponding stabiliser, so the product
of the radices is at most $|\mathrm{Co}_2|$, with equality only if every basic orbit is complete; it is equal, so the chain is
complete. Every $k\in K$ is then uniquely $k=x\,\ell(\pi k)$ with $x\in\ker\pi$, and $x$ is determined by $\bar c(x)$, since
bit $21$ of $c(x)$ is $0$. This gives the bijection; the identity has all digits $0$ and $x=1$.
\end{proof}

The chain of Theorem~\ref{thm:klaw} acts on $\Lambda/2\Lambda$ rather than on $\Lambda$, so its radices differ from those
of Corollary~\ref{cor:co2chain}. Its strong generators come from random elements, so the labels are fixed by the stored
chain\footnote{The scripts used for the stored chain are
\[
\texttt{scripts/k\_co2\_chain.py}\qquad\text{and}\qquad\texttt{scripts/k\_law.py}.
\]
The latter reads the stored chain
\texttt{scripts/data/k\_co2\_chain.pkl}.} and not by a formula. With it, $\lambda_K^{-1}(\lambda_K(k))=k$ held for $30$ random
$k\in K$, $\lambda_K(\lambda_K^{-1}(n))=n$ for $30$ random labels $n$, and the bit-$21$ rule predicted commutation with
$\beta^+$ for $300$ random elements of $Q$.

\begin{lemma}[Positions on the orbits $2A1$ and $2A0$]\label{lem:rb2a0}
Order the ten $K$-orbits of \texttt{BabyAxis} by size, so that $2A1$ (one axis) and $2A0$ ($93\,150$ axes) come first. For an
axis $a$ let $t_a$ be the involution returned by \texttt{central\_involution}, and assume that $t_{ak}=k^{-1}t_ak$ for $k\in K$.
If $a$ lies in $2A0$, then $t_a\in Q$; write $c(t_a)=\lambda+2^{24}e$ with $\lambda$ a $24$-bit vector and $e\in\{0,1\}$. Then
$\lambda$ lies in the first basic orbit of the chain of Theorem~\ref{thm:klaw}, and
\[
r_B(a)=0\ \ (a\in 2A1),\qquad r_B(a)=1+2\operatorname{rank}(\lambda)+e\ \ (a\in 2A0),
\]
with $\operatorname{rank}$ the digit of that basic orbit, is a bijection from $2A1\cup2A0$ onto $\{0,\dots,93\,150\}$.
\end{lemma}

\begin{proof}
The orbit $2A1$ is a single axis. For $2A0$, the representative $a_0$ of \texttt{mmgroup} has $t_{a_0}=x_{\mathtt{800h}}d_{\mathtt{200h}}\in Q$,
with $\lambda_0=\mathtt{0x800200}$ in the first basic orbit $\Delta$, of size $46\,575$. By the assumption, $a\mapsto t_a$ maps $2A0$
onto the $K$-orbit of $t_{a_0}$. Conjugation by $k\in K$ moves $\lambda$ by $\pi(k)$, so this orbit lies in
$\{x_\mu x_{-1}^e:\mu\in\Delta,\ e\in\{0,1\}\}$, a set of $2\cdot46\,575=93\,150$ elements. Since $\mu$ and $\mathtt{0x200}$ are
independent mod $2$, some $y\in Q$ commutes with $\beta^+$ but not with $x_\mu$, so $y\in C_Q(\beta^+)\le K$ changes $e$; with
$\pi(K)=\mathrm{Co}_2$ transitive on $\Delta$, the $K$-orbit is the whole set. So $a\mapsto t_a$ maps a set of $93\,150$ axes onto a
set of $93\,150$ involutions, and is a bijection.
\end{proof}

A script\footnote{\texttt{scripts/k\_rb\_2a0.py}.} checks the assumption and the formula: for $40$ random $k\in K$ the axis $a_0k$
has $t_{a_0k}=k^{-1}t_{a_0}k$ in $Q$, with $\lambda\in\Delta$ and orbit $2A0$, and the $40$ positions are distinct; for $40$ random
positions, the axis built from a transversal element of $\Delta$ and, if needed, a sign change in $C_Q(\beta^+)$ returns the same
position. For the axis of $2A1$, \texttt{central\_involution} returns $\beta^+$ itself, so $t_a$ is not in general the $2A$
involution of $a$; the lemma uses only the assumed compatibility with $K$.

\begin{lemma}[Positions on the orbit $2B1$]\label{lem:rb2b1}
Let $a$ be an axis in the orbit $2B1$, of size $7\,286\,400$, with the notation and the assumption of
Lemma~\ref{lem:rb2a0}. Let $S(a)=[\Lambda/2\Lambda,\pi(t_a)]$, written in reduced echelon form, and let $F(a)$ be the set of conjugates of $t_a$ under
$C_Q(\beta^+)$, ordered by the normal forms of \texttt{mmgroup}. Then $S(a)$ has dimension $8$, the spaces $S(a)$ form one
$\mathrm{Co}_2$-orbit $\Sigma$ of size $56\,925$, $|F(a)|=128$, and
\[
r_B(a)=93\,151+128\operatorname{rank}_\Sigma S(a)+\operatorname{rank}_{F(a)}t_a
\]
is a bijection from $2B1$ onto $\{93\,151,\dots,7\,379\,550\}$.
\end{lemma}

\begin{proof}
Let $a_0$ be the representative and $t=t_{a_0}$. Its stabiliser $H$ in $K$ has index $7\,286\,400$ and centralises $t$, by the
assumption, so the $K$-orbit of $t$ has at most $7\,286\,400$ elements. Since $S(t^k)=S(t)\pi(k)$, conjugates with different
spaces are different; the orbit $\Sigma$ of $S(t)$, computed by breadth-first search, has $56\,925$ elements. Over each image
$\pi(t^k)$ lie the $C_Q(\beta^+)$-conjugates of $t^k$, and there are $128$ of them, since $C_Q(\beta^+)$ is normal in $K$ and
$t$ has $128$. So the orbit of $t$ has at least $56\,925\cdot128=7\,286\,400$ elements. Hence $C_K(t)=H$, $a\mapsto t_a$ is a
bijection onto the orbit of $t$, each value of $S$ belongs to exactly one image $\pi(t_a)$, and over it lie exactly the
$128$ elements of $F(a)$. The formula numbers these $56\,925\cdot128$ pairs.
\end{proof}

The proof uses no random sampling. As a check, random elements of $H$ generate a subgroup of $\mathrm{Co}_2$ of order
$743\,178\,240=|\mathrm{Co}_2|/56\,925$, the order of $2^{1+8}{:}\mathrm{Sp}_6(2)$, and a script\footnote{\texttt{scripts/k\_rb\_2b1.py};
the stabiliser is probed in \texttt{scripts/k\_rb\_2b1\_probe.py}.} finds, for $20$ random $k\in K$, that $S(a_0k)=S(a_0)\pi(k)$ and the
$20$ positions are distinct, and, for $20$ random positions, that the axis built from a transversal element of $\Sigma$ and an
element of $C_Q(\beta^+)$ returns the same position.

\begin{lemma}[Positions on the orbit $2B0$]\label{lem:rb2b0}
Let $a$ be an axis in the orbit $2B0$, of size $262\,310\,400$, with the notation and the assumption of
Lemma~\ref{lem:rb2a0}, and let $S(a)$ and $F(a)$ be as in Lemma~\ref{lem:rb2b1}. Then $S(a)$ has dimension $8$, the spaces
$S(a)$ form one $\mathrm{Co}_2$-orbit $\Sigma$ of size $1\,024\,650$, $|F(a)|=256$, and
\[
r_B(a)=7\,379\,551+256\operatorname{rank}_\Sigma S(a)+\operatorname{rank}_{F(a)}t_a
\]
is a bijection from $2B0$ onto $\{7\,379\,551,\dots,269\,689\,950\}$.
\end{lemma}

\begin{proof}
The proof of Lemma~\ref{lem:rb2b1} applies with new numbers. The stabiliser $H$ of the representative $a_0$ has index
$262\,310\,400$ in $K$ and centralises $t=t_{a_0}$, so the $K$-orbit of $t$ has at most $262\,310\,400$ elements. The orbit
$\Sigma$ of $S(t)$, computed by breadth-first search, has $1\,024\,650$ elements, and over each image $\pi(t^k)$ lie the
$256$ conjugates of $t^k$ under the normal subgroup $C_Q(\beta^+)$. So the orbit of $t$ has at least
$1\,024\,650\cdot256=262\,310\,400$ elements. Hence $C_K(t)=H$, $a\mapsto t_a$ is a bijection onto the orbit of $t$, and the
formula numbers the $1\,024\,650\cdot256$ pairs $(S(a),t_a)$.
\end{proof}

Again no sampling enters the proof. The index of $\pi(H)$ in $\mathrm{Co}_2$ is $1\,024\,650$, so $|\pi(H)|=41\,287\,680$. A
script\footnote{\texttt{scripts/k\_rb\_2b0.py}; the breadth-first search over $\Sigma$ takes about $145$ seconds.} finds, for
$20$ random $k\in K$, that $S(a_0k)=S(a_0)\pi(k)$ and the $20$ positions are distinct, and, for $20$ random positions, that
the axis built from a transversal element of $\Sigma$ and an element of $C_Q(\beta^+)$ returns the same position.

\begin{lemma}[Positions on the orbit $4B1$]\label{lem:rb4b1}
Let $a$ be an axis in the orbit $4B1$, of size $470\,060\,236\,800$, with the notation and the assumption of
Lemma~\ref{lem:rb2a0}, and let $S(a)$ be as in Lemma~\ref{lem:rb2b1}. Let $R(a)$ be the radical of the short vectors
returned by \texttt{find\_short(1)}, in reduced echelon form, and $P(a)=(S(a),R(a))$. Let $H$ be the stabiliser of the
representative $a_0$ in $K$ and $X=H\cap C_Q(\beta^+)$. Then the pairs $P(a)$ form one $\mathrm{Co}_2$-orbit $W$ of size
$3\,586\,275$, $|X|=128$, $X$ meets the centre of $Q$ trivially, and every axis of $4B1$ is $a=a_0\,c\,u(P(a))$ with
$c\in C_Q(\beta^+)$, where $u(P)$ is the transversal element of $W$ and the coset $Xc$ depends on $a$ alone. Write
$\operatorname{rank}(Xc)=2\,\mathrm{idx}(\rho)+e\in\{0,\dots,2^{17}-1\}$, with $\rho$ the reduction of $c$ modulo $X$
($16$ free bits) and $e$ the centre bit. Then
\[
r_B(a)=14\,113\,555\,551+2^{17}\operatorname{rank}_W P(a)+\operatorname{rank}(Xc)
\]
is a bijection from $4B1$ onto $\{14\,113\,555\,551,\dots,484\,173\,792\,350\}$.
\end{lemma}

\begin{proof}
The stabiliser $H$ has index $470\,060\,236\,800$ in $K$, so
$|H|=2^{24}|\mathrm{Co}_2|/470\,060\,236\,800=1\,509\,949\,440=|\pi(H)|\cdot|X|$ with $|\pi(H)|=11\,796\,480$.
The orbit $W$ of $P(a_0)$, computed by breadth-first search, has $3\,586\,275=|\mathrm{Co}_2:\pi(H)|$ elements, so the
stabiliser of $P(a_0)$ in $\mathrm{Co}_2$, which contains $\pi(H)$, equals $\pi(H)$. For an axis $a$ let $q\in K$ with $aq=a_0$
and $w=q^{-1}u(P(a))^{-1}$; then $\pi(w)$ fixes $P(a_0)$, so $w\in\pi^{-1}(\pi(H))=H\,C_Q(\beta^+)$, and $w=hc$ gives
$a=a_0\,c\,u(P(a))$. Two choices of $c$ differ by an element of $H\cap C_Q(\beta^+)=X$, so $a\mapsto(P(a),Xc)$ is well defined,
and it is injective because $a$ is recovered from it. Its target has $3\,586\,275\cdot2^{24}/128=470\,060\,236\,800$ elements,
the size of the orbit, so it is a bijection, and the formula numbers its image.
\end{proof}

Unlike for $2B1$ and $2B0$, the space $S(a)$ alone does not separate the axes of $4B1$; the radical $R(a)$ is needed. A
script\footnote{\texttt{scripts/k\_rb\_pos.py 4B1}; the breadth-first search over $W$ takes about $940$ seconds and $0.63$\,GB,
and the $3\,586\,275$ pairs are stored in \texttt{scripts/data/k\_rb\_4b1\_pairs.bin}. The orbit of $P(a_0)$ is checked in
\texttt{scripts/k\_rb\_sorbit.py}.} confirms $|\pi(H)|\,|X|=|H|$ exactly, finds for $20$ random axes that $a=a_0\,c\,u$ and
$\mathrm{axis}(r_B(a))=a$ with $20$ distinct positions, and, for $20$ random positions, that position $\to$ axis $\to$ position
is the identity.

\begin{lemma}[Positions on the orbit $6C1$]\label{lem:rb6c1}
Let $a$ be an axis in the orbit $6C1$, the largest of the ten, of size $6\,685\,301\,145\,600$, with the notation and the
assumption of Lemma~\ref{lem:rb2a0}. Let $S(a)$, $H$, $X$, $u(P)$ and $\operatorname{rank}(Xc)$ be as in Lemma~\ref{lem:rb4b1},
but let $R(a)$ be the radical of the short vectors returned by \texttt{find\_short(3)}, and $P(a)=(S(a),R(a))$. Then the
pairs $P(a)$ form one $\mathrm{Co}_2$-orbit $W$ of size $1\,593\,900$, $|X|=4$, $X$ meets the centre of $Q$ trivially, the
reduction modulo $X$ leaves $21$ free bits, and
\[
r_B(a)=5\,022\,147\,527\,775+2^{22}\operatorname{rank}_W P(a)+\operatorname{rank}(Xc)
\]
is a bijection from $6C1$ onto $\{5\,022\,147\,527\,775,\dots,11\,707\,448\,673\,374\}$.
\end{lemma}

\begin{proof}
The proof of Lemma~\ref{lem:rb4b1} applies with new numbers. The stabiliser $H$ has index $6\,685\,301\,145\,600$ in $K$, so
$|H|=2^{24}|\mathrm{Co}_2|/6\,685\,301\,145\,600=106\,168\,320=|\pi(H)|\cdot|X|$ with $|\pi(H)|=26\,542\,080$. The orbit $W$,
computed by breadth-first search, has $1\,593\,900=|\mathrm{Co}_2:\pi(H)|$ elements, so the stabiliser of $P(a_0)$ is $\pi(H)$,
and $a\mapsto(P(a),Xc)$ is an injection into a set of $1\,593\,900\cdot2^{24}/4=6\,685\,301\,145\,600$ elements, the size of
the orbit.
\end{proof}

The last position is $n_B-1$, so Lemmas~\ref{lem:rb2a0}--\ref{lem:rb6c1} fix both ends of the numbering of the ten orbits.
A script\footnote{\texttt{scripts/k\_rb\_pos.py 6C1}; the breadth-first search over $W$ takes about $400$ seconds and
$0.35$\,GB, and the $1\,593\,900$ pairs are stored in \texttt{scripts/data/k\_rb\_6c1\_pairs.bin}.} confirms
$|\pi(H)|\,|X|=|H|$ exactly, finds for $20$ random axes that $a=a_0\,c\,u$ and $\mathrm{axis}(r_B(a))=a$ with $20$ distinct
positions, and, for $20$ random positions, that position $\to$ axis $\to$ position is the identity; a second run from the stored
pairs repeats these checks.

\begin{lemma}[Positions on the orbit $6A1$]\label{lem:rb6a1}
Let $a$ be an axis in the orbit $6A1$, of size $9\,646\,899\,200$, with the notation and the assumption of
Lemma~\ref{lem:rb2a0}, and let $H$, $X$, $u(P)$ and $\operatorname{rank}(Xc)$ be as in Lemma~\ref{lem:rb4b1}. Here $S(a)=0$, and
we take $P(a)=R(a)$, the radical of the short vectors returned by \texttt{find\_short(5)}, a space of dimension $2$. Then the
spaces $P(a)$ form one $\mathrm{Co}_2$-orbit $W$ of size $2300$, $|X|=4$, $X$ meets the centre of $Q$ trivially, the
reduction modulo $X$ leaves $21$ free bits, and
\[
r_B(a)=4\,466\,656\,351+2^{22}\operatorname{rank}_W P(a)+\operatorname{rank}(Xc)
\]
is a bijection from $6A1$ onto $\{4\,466\,656\,351,\dots,14\,113\,555\,550\}$.
\end{lemma}

\begin{proof}
The proof of Lemma~\ref{lem:rb4b1} applies with new numbers: $|H|=2^{24}|\mathrm{Co}_2|/9\,646\,899\,200=73\,574\,645\,760=
|\pi(H)|\cdot|X|$ with $|\pi(H)|=18\,393\,661\,440$, the orbit $W$ has $2300=|\mathrm{Co}_2:\pi(H)|$ elements, and
$2300\cdot2^{24}/4=9\,646\,899\,200$ is the size of the orbit.
\end{proof}

\begin{lemma}[Positions on the orbit $4C1$]\label{lem:rb4c1}
Let $a$ be an axis in the orbit $4C1$, of size $537\,211\,699\,200$, with the notation and the assumption of
Lemma~\ref{lem:rb2a0}, and let $S(a)$ be as in Lemma~\ref{lem:rb2b1} and $H$, $X$, $u(P)$ and $\operatorname{rank}(Xc)$ as in
Lemma~\ref{lem:rb4b1}, with $P(a)=S(a)$. Then $S(a)$ has dimension $8$, the spaces $S(a)$ form one $\mathrm{Co}_2$-orbit $W$ of
size $1\,024\,650$, $|X|=32$, $X$ meets the centre of $Q$ trivially, the reduction modulo $X$ leaves $18$ free bits, and
\[
r_B(a)=484\,173\,792\,351+2^{19}\operatorname{rank}_W S(a)+\operatorname{rank}(Xc)
\]
is a bijection from $4C1$ onto $\{484\,173\,792\,351,\dots,1\,021\,385\,491\,550\}$.
\end{lemma}

\begin{proof}
The proof of Lemma~\ref{lem:rb4b1} applies with new numbers: $|H|=2^{24}|\mathrm{Co}_2|/537\,211\,699\,200=1\,321\,205\,760=
|\pi(H)|\cdot|X|$ with $|\pi(H)|=41\,287\,680$, the orbit $W$ has $1\,024\,650=|\mathrm{Co}_2:\pi(H)|$ elements, and
$1\,024\,650\cdot2^{24}/32=537\,211\,699\,200$ is the size of the orbit.
\end{proof}

For $4C1$, unlike for $4B1$ and $6C1$, the space $S(a)$ alone separates the $\mathrm{Co}_2$-positions, as for $2B0$; for $6A1$
the space $S(a)$ is trivial and a radical alone suffices. Scripts\footnote{\texttt{scripts/k\_rb\_pos.py 6A1} and
\texttt{scripts/k\_rb\_pos.py 4C1}; the breadth-first searches take under a second and about $163$ seconds ($0.20$\,GB), and
the orbits are stored in \texttt{scripts/data/k\_rb\_6a1\_pairs.bin} and \texttt{scripts/data/k\_rb\_4c1\_pairs.bin}.} confirm
$|\pi(H)|\,|X|=|H|$ exactly for both orbits and, for $20$ random axes and $20$ random positions in each, the forward and inverse
checks of Lemma~\ref{lem:rb4b1}, with distinct positions.

\begin{lemma}[Positions on the orbit $4A1$]\label{lem:rb4a1}
Let $a$ be an axis in the orbit $4A1$, of size $4\,196\,966\,400$, with the notation and the assumption of
Lemma~\ref{lem:rb2a0}, and let $H$, $X$, $u(P)$ and $\operatorname{rank}(Xc)$ be as in Lemma~\ref{lem:rb4b1}. Let $V_3(a)$ be the
set of the $46$ vectors returned by \texttt{find\_short(3)}, taken as a set and not as a span, and put $P(a)=V_3(a)$, on which
$\mathrm{Co}_2$ acts elementwise. Then the sets $V_3(a)$ form one $\mathrm{Co}_2$-orbit $W$ of size $1\,024\,650$, $|X|=4096$, $X$
meets the centre of $Q$ trivially, the reduction modulo $X$ leaves $11$ free bits, and
\[
r_B(a)=269\,689\,951+2^{12}\operatorname{rank}_W V_3(a)+\operatorname{rank}(Xc)
\]
is a bijection from $4A1$ onto $\{269\,689\,951,\dots,4\,466\,656\,350\}$.
\end{lemma}

\begin{proof}
The proof of Lemma~\ref{lem:rb4b1} applies, with $P(a)$ moved elementwise, and with new numbers:
$|H|=2^{24}|\mathrm{Co}_2|/4\,196\,966\,400=169\,114\,337\,280=|\pi(H)|\cdot|X|$ with $|\pi(H)|=41\,287\,680$, the orbit $W$ has
$1\,024\,650=|\mathrm{Co}_2:\pi(H)|$ elements, and $1\,024\,650\cdot2^{24}/4096=4\,196\,966\,400$ is the size of the orbit.
\end{proof}

\begin{lemma}[Positions on the orbit $10A1$]\label{lem:rb10a1}
Let $a$ be an axis in the orbit $10A1$, of size $4\,000\,762\,036\,224$, with the notation and the assumption of
Lemma~\ref{lem:rb2a0}, and let $H$, $X$, $u(P)$ and $\operatorname{rank}(Xc)$ be as in Lemma~\ref{lem:rb4b1}. Let $V_1(a)$ be the
set of the $100$ vectors returned by \texttt{find\_short(1)}, taken as a set, and put $P(a)=V_1(a)$. Then the sets $V_1(a)$ form
one $\mathrm{Co}_2$-orbit $W$ of size $476\,928$, $|X|=2$, $X$ meets the centre of $Q$ trivially, the reduction modulo $X$ leaves
$22$ free bits, and
\[
r_B(a)=1\,021\,385\,491\,551+2^{23}\operatorname{rank}_W V_1(a)+\operatorname{rank}(Xc)
\]
is a bijection from $10A1$ onto $\{1\,021\,385\,491\,551,\dots,5\,022\,147\,527\,774\}$.
\end{lemma}

\begin{proof}
The proof of Lemma~\ref{lem:rb4b1} applies, with $P(a)$ moved elementwise, and with new numbers:
$|H|=2^{24}|\mathrm{Co}_2|/4\,000\,762\,036\,224=177\,408\,000=|\pi(H)|\cdot|X|$ with $|\pi(H)|=88\,704\,000$, the orbit $W$ has
$476\,928=|\mathrm{Co}_2:\pi(H)|$ elements, and $476\,928\cdot2^{24}/2=4\,000\,762\,036\,224$ is the size of the orbit.
\end{proof}

On these two orbits every space attached to $a$ is too coarse: the radicals of the short vectors have $\mathrm{Co}_2$-orbits of
sizes $46\,575$ on $4A1$ and $1$ on $10A1$, whereas the full sets $V_3(a)$ and $V_1(a)$ separate the positions. On $10A1$,
$|\pi(H)|=88\,704\,000$ is the order of $\mathrm{HS}{:}2$. Scripts\footnote{\texttt{scripts/k\_rb\_pos.py 4A1 V3} and
\texttt{scripts/k\_rb\_pos.py 10A1 V1}; the breadth-first searches take about $780$ and $800$ seconds ($0.24$ and $0.21$\,GB), and
the orbits are stored in \texttt{scripts/data/k\_rb\_4a1\_pairs.bin} and \texttt{scripts/data/k\_rb\_10a1\_pairs.bin}.} confirm
$|\pi(H)|\,|X|=|H|$ exactly for both orbits and, for $20$ random axes and $20$ random positions in each, the forward and inverse
checks of Lemma~\ref{lem:rb4b1}, with distinct positions.

\begin{remark}[Stabilisers on $4A1$ and $10A1$]\label{rem:rbstab}
The orders in Lemmas~\ref{lem:rb4a1} and~\ref{lem:rb10a1} identify the groups $\pi(H)$, using only the list of maximal
subgroups of $\mathrm{Co}_2$ in the ATLAS, of indices $2300$, $46\,575$, $47\,104$, $56\,925$, $476\,928$, $1\,024\,650$ and larger.
A subgroup of index $m$ lies in a maximal subgroup whose index divides $m$.

On $10A1$, $m=476\,928=2^8\cdot3^4\cdot23$, and the only maximal subgroup with index dividing $m$ is $\mathrm{HS}{:}2$, of index
exactly $m$; so $\pi(H)\cong\mathrm{HS}{:}2$. Here $X$ has order $2$ and is central in $H$, and $H$ is a non-split extension
$2\cdot(\mathrm{HS}{:}2)$: a random search in $H$\footnote{\texttt{scripts/k\_rb\_stab\_probe.py}, run in a few seconds; its random samples of $H$ are seeded (seed $1$ by default), so a run with a fixed seed is reproducible exactly, and other seeds change the counts but not the conclusions of this remark.} found elements $g$ of order $4$ with $\pi(g)$ an involution, and then $gx$ has
order $4$ as well for $x\in X$, so this involution has no involution above it.

On $4A1$, the radical $R(a)$ of \texttt{find\_short(1)} has a $\mathrm{Co}_2$-orbit of size $46\,575=3^4\cdot5^2\cdot23$, and the only
maximal subgroup with index dividing $46\,575$ is $2^{10}{:}M_{22}{:}2$, so $\mathrm{Stab}(R(a))\cong2^{10}{:}M_{22}{:}2$. Since $R(a)$
depends only on $a$, $\pi(H)\le\mathrm{Stab}(R(a))$, with index $908\,328\,960/41\,287\,680=22$. If the image of $\pi(H)$ in
$M_{22}{:}2$ has index $d$, then $22/d$ is a power of $2$, and $d=11$ is impossible because $M_{22}$ has no subgroup of index less than
$22$; so $d=22$, $\pi(H)$ contains the normal $2^{10}$, and its image is a point stabiliser $L_3(4){:}2_2$ of $M_{22}{:}2$ on $22$
points. Hence
\[
\pi(H)\cong2^{10}{:}L_3(4){:}2_2,
\]
which has the order of the maximal subgroup $(2^4\times2^{1+6}).A_8$ but is not conjugate to it: it contains no element of order
$15$. The full set $V_3(a)$ thus singles out one of the $22$ points permuted by $M_{22}{:}2$. The kernel $X=H\cap Q$ is elementary
abelian of order $2^{12}$. As an $\mathbb F_2H$-module under conjugation, computed with $60$ random elements of $H$\footnote{By the same script.}, it has
$[X,H]=X$, a socle of dimension $2$ on which $H$ acts trivially, a unique maximal submodule $U$ of dimension $3$ containing the
socle, and an irreducible head $X/U$ of dimension $9$: of the $4095$ non-zero vectors, the $3$ in the socle span lines, the $4$ in
$U$ outside the socle span $U$, and all others span $X$. Fewer generators could only add invariant subspaces, so this lattice is
exact if the $60$ elements generate $H$. The action is faithful on $\pi(H)$: the image of $H$ in $\mathrm{GL}(X)$, as a permutation group on the $4095$ non-zero vectors, has order at least $41\,287\,680$ by a Schreier--Sims computation\footnote{\texttt{scripts/k\_rb\_stab\_faithful.py}, seeded (seed $1$ by default).}, and at most $|\pi(H)|=41\,287\,680$ because $X$ is abelian. So the kernel of $H$ on $X$ is exactly $X$, and the normal subgroup $2^{10}$ of $\pi(H)$ acts non-trivially on $X$. Whether $H$ splits over $X$ is not decided here.
\end{remark}

\subsection{The top-axis decomposition}\label{ssec:monstertop}
The law is implemented in \texttt{monster/monster\_law.py}\footnote{Run with \texttt{mmgroup} $1.0.8$.}. Let $\beta$ be the type-$2$ Leech vector of $x_\beta$ (\texttt{0x200} in \texttt{mmgroup}) and
$K=C_{G_{x0}}(x_\beta)\cong2^{1+23}.\mathrm{Co}_2$, the group of Theorem~\ref{thm:klaw}. The ordinal of $g\in\mathbb M$ is
\[
\lambda_{\mathbb M}(g)=r_3+|K|\,r_2+|K|\,n_B\,r_1,
\]
where $r_1=r(g^{-1}zg)\in\{0,\dots,N-1\}$ is the top-axis rank, $r_2=r_B\in\{0,\dots,n_B-1\}$ is the rank on the $C$-orbit of
$v^-$ (Lemmas~\ref{lem:rb2a0}--\ref{lem:rb10a1}) and $r_3=\lambda_K\in\{0,\dots,|K|-1\}$. This is Proposition~\ref{prop:monster} with $\lambda_C$ expanded as in Remark~\ref{rem:monsterstatus}, since
$|C|=|K|\,n_B$; the top-axis rank $r_1$ is the most significant digit.

The digit $r_1$ is the offset of the orbit $c$ of $a=g^{-1}zg$ plus the position of the coset $H_cq(a)^{-1}$ among the cosets
of $H_c$ in $G_{x0}$. Positions inside an orbit are computed by a \emph{typed} reduction to $K$:
\begin{enumerate}
\item \emph{Type.} A coset $H_cx$ determines the $H_c$-orbit of $\beta x^{-1}$ among the $98\,280$ type-$2$ vectors of
$\Lambda/2\Lambda$. This orbit is the type $i$ of the coset; types are keyed by their orbit lengths $c_i$.
\item \emph{Reduction.} Fix a representative $u_i$ of type $i$ and $y_i\in G_{x0}$ with $u_iy_i=\beta$. The cosets of type $i$
correspond to the cosets of $H_i=y_i^{-1}\,\mathrm{Stab}_{H_c}(u_i)\,y_i$ in the stabiliser of $\beta$, a group contained in
$K$ or in $K$ extended by $z$.
\item \emph{Ranking.} A coset of $H_i$ is ranked in three layers. The $\mathrm{Co}_2$ layer is a chain of Leech-lattice objects
fixed by a \emph{recipe}, whose levels are orbits of the image $\bar H_i$ of $H_i$ in $\mathrm{Co}_2$. It is followed by an
$\bar S$ chain inside the joint stabiliser of the recipe objects, and then by the layer of $Q_K=O_2(K)$.
\item \emph{Offsets.} Type $i$ fills a block of $\mathrm{loc}_i=|c|\,c_i/98\,280$ consecutive positions, where $|c|$ is the size of
the orbit. The blocks are placed end to end.
\end{enumerate}

\paragraph{Certification.} Let $P_1$ be the product of the level sizes of the $\mathrm{Co}_2$ layer; the $\bar S$ chain is built
to reach the order $|\mathrm{Co}_2|/P_1$. The recipe objects are $H_i$-orbits, so their joint stabiliser always contains
$\bar H_i$, and the recipe certifies exactly when this stabiliser equals $\bar H_i$, that is, when
$P_1=|\mathrm{Co}_2|/|\bar H_i|$. A match of $P_1$ is necessary but not sufficient: a recipe is accepted only once its $\bar S$
chain has been built and certified. Types of one orbit with the same $c_i$ can need different recipes, so recipes are listed
per type and tried in order. Failure of the $\bar S$ chain to reach its target order can depend on the random seed, so seeds
$1$--$5$ are tried before a recipe is rejected. An orbit is written only if every type certifies and the $\mathrm{loc}_i$
sum to $|c|$. For example, the type $c_i=240$ of orbit $12C$ certifies with the chain
$\mathrm{Co}_2>U_6(2)>\mathrm{Sp}_6(2)>G_2(2)$, of levels $4600\cdot6336\cdot120$, so that $|\bar H_i|=12\,096=|G_2(2)|$.

All twelve orbits $2A,\dots,12C$ of Remark~\ref{rem:monsterstatus} are certified in this way. The last one built, $12C$, has
$21$ types, all certified with the first recipe on their lists, in $311$\,s, and its range of ordinals ends exactly at $|\mathbb M|$.

\begin{theorem}[Law of $\mathbb M$]\label{thm:mlaw}
The map $\lambda_{\mathbb M}$ is a bijection $\mathbb M\to\{0,\dots,|\mathbb M|-1\}$. It and its inverse are computed by
\texttt{monster/monster\_law.py} for every element of $\mathbb M$, and transporting the multiplication of $\mathbb M$ along it gives
a digital law of $\mathbb M$.
\end{theorem}

\begin{proof}
By Proposition~\ref{prop:monster} and Remark~\ref{rem:monsterstatus}, it suffices that $r$, $r_B$ and $\lambda_K$ are bijections
onto their ranges. For $\lambda_K$ this is Theorem~\ref{thm:klaw}, and for $r_B$ it is
Lemmas~\ref{lem:rb2a0}--\ref{lem:rb10a1}. For $r$, the twelve orbits fill $\{0,\dots,N-1\}$ by
Remark~\ref{rem:monsterstatus}, and inside each orbit the certified types fill the orbit's block.
\end{proof}

\subsection{Validation}\label{ssec:monstervalid}
Each orbit was validated after its build with \texttt{monster/validate\_top.py}. Write $D=\lambda_{\mathbb M}^{-1}$ for the map from
ordinals to elements and $\mathrm{mul}(n,m)=\lambda_{\mathbb M}\bigl(D(n)D(m)\bigr)$ for the transported product; products of
elements are reduced to the \texttt{mmgroup} normal form. For the orbit $12C$, let $m$ and $a$ be uniform in its range of
ordinals and $n$ and $b$ uniform in $\{0,\dots,|C|-1\}$, the ordinals with top-axis rank $r_1=0$. The script checks:
(i) $\lambda_{\mathbb M}\bigl(D(m)\bigr)=m$;
(ii) $\lfloor\mathrm{mul}(n,m)/|C|\rfloor=\lfloor m/|C|\rfloor$ and $D\bigl(\mathrm{mul}(n,m)\bigr)=D(n)D(m)$, that is, the product keeps
the top-axis rank of $m$ and is the product in $\mathbb M$;
(iii) $\lambda_{\mathbb M}\bigl(D(b)D(a)\bigr)=\mathrm{mul}(b,a)$, where the left side is computed from the normal form of
$D(b)D(a)$ and \texttt{mul} from the unreduced product, so that $\lambda_{\mathbb M}$ depends only on the element and not
on the word that represents it;
(iv) $D\bigl(\lambda_{\mathbb M}(g)\bigr)=g$ for random elements $g\in\mathbb M$ (\texttt{MM('r','M')}).
With $2000$ samples, checks (i)--(iii) passed in $2000$ of $2000$ cases. In (iv), all $2000$ elements were decoded, none raised
\texttt{NotImplementedError}, and there were no mismatches. The orbits implemented in \texttt{\_TOP\_IMPL} cover all $N$ top axes; the two orbits of types $2A$ and $2B$
alone account for a share $1.23\cdot10^{-10}$.

\section{From laws to character tables}\label{sec:chartables}
A digital law of $G$ gives multiplication, inversion and uniform random elements (uniform random labels), which is all that the
standard algorithms for conjugacy classes and characters require. This section takes three steps in that direction: character
degrees from rank-three actions, conjugacy classes and power maps, and a complete character table of $M_{11}$ computed from the
tables of Section~\ref{ssec:mathieuexamples} alone.

\subsection{Character degrees from rank-three actions}\label{ssec:rank3}
\begin{proposition}[Rank-three actions]\label{prop:rank3}
Let $G$ act transitively on $\Ap n$ with rank $3$, and suppose that the orbit $O$ of $G$ on ordered pairs that contains $(0,1)$ is
symmetric. Then the graph $\Gamma$ with edge set $O$ is strongly regular with parameters $(n,k,\lambda,\mu)$, its adjacency matrix
$A$ has eigenvalues $k$, $r$, $s$ with multiplicities $1$, $f$, $g$, and the permutation character of $G$ is $1+\chi+\psi$ with
$\chi,\psi$ irreducible of degrees $f$ and $g$.
\end{proposition}
\begin{proof}
$G$ is transitive on edges and on non-edges of $\Gamma$, so the number of common neighbours of two vertices is $\lambda$ or $\mu$
according as they are adjacent or not, and $A^2=kI+\lambda A+\mu(J-I-A)$. The matrices commuting with all permutation matrices of $G$
are spanned by $I$, $A$ and $J-I-A$, the adjacency matrices of the three orbits on pairs. This algebra has dimension $3$, so the
permutation character is a sum of three distinct irreducibles, one of them trivial. On the orthogonal complement of the constant
vectors, $A^2=(\lambda-\mu)A+(k-\mu)I$, so $A$ has at most two eigenvalues $r,s$ there, and exactly two because $0<k<n-1$. Each
eigenspace is $G$-invariant and nonzero, so each is one of the two non-trivial constituents. Finally $f+g=n-1$ and
$\operatorname{tr}A=k+fr+gs=0$ determine $f$ and $g$.
\end{proof}

Table~\ref{tab:rank3} applies this to the rank-three actions of the generator files of Table~\ref{tab:canondata} and of
$\mathrm{McL}$, $\Co_2$ and $\mathrm{Fi}_{22}$ from the library of primitive groups. In each case $f$ and $g$ are degrees of
irreducible characters in the ATLAS~\cite{ATLAS}.

\begin{table}[ht]
\centering\small
\begin{tabular}{lrlll}
group & $n$ & $(k,\lambda,\mu)$ & $r^f$, $s^g$ & permutation character\\ \hline
$J_2$ & $100$ & $(36,14,12)$ & $6^{36}$, $-4^{63}$ & $1+36+63$\\
$\mathrm{HS}$ & $100$ & $(77,60,56)$ & $7^{22}$, $-3^{77}$ & $1+22+77$\\
$\mathrm{McL}$ & $275$ & $(112,30,56)$ & $2^{252}$, $-28^{22}$ & $1+22+252$\\
$\mathrm{Suz}$ & $1782$ & $(416,100,96)$ & $20^{780}$, $-16^{1001}$ & $1+780+1001$\\
$\Co_2$ & $2300$ & $(891,378,324)$ & $63^{275}$, $-9^{2024}$ & $1+275+2024$\\
$\mathrm{Fi}_{22}$ & $3510$ & $(2816,2248,2304)$ & $8^{3080}$, $-64^{429}$ & $1+429+3080$\\
$\mathrm{Ru}$ & $4060$ & $(1755,730,780)$ & $15^{3276}$, $-65^{783}$ & $1+783+3276$\\
\end{tabular}
\caption{Rank-three actions: strongly regular parameters, eigenvalues with multiplicities, and the degrees of the constituents.}\label{tab:rank3}
\end{table}

\subsection{Conjugacy classes and power maps}\label{ssec:classes}
In the language of laws, the conjugacy classes of $G$ are the orbits of $x\mapsto g*x*\lambda(\lambda^{-1}(g)^{-1})$ on
$\Ap{|G|}$, and the $p$-th power map sends the class of $x$ to the class of $x*\cdots*x$ ($p$ factors). Table~\ref{tab:classes}
summarises the classes of the sporadic groups of degree at most $4060$, computed from the same permutation groups, and
Table~\ref{tab:j2classes} lists the classes of $J_2$ with their power maps. Classes are named by element order and then by
decreasing centraliser order; letters within equal order and centraliser order are arbitrary.

\begin{table}[ht]
\centering\small
\begin{tabular}{lrrp{6.4cm}}
group & order & classes & element orders\\ \hline
$M_{11}$ & $7920$ & $10$ & 1, 2, 3, 4, 5, 6, 8, 11\\
$M_{12}$ & $95\,040$ & $15$ & 1, 2, 3, 4, 5, 6, 8, 10, 11\\
$M_{22}$ & $443\,520$ & $12$ & 1, 2, 3, 4, 5, 6, 7, 8, 11\\
$M_{23}$ & $10\,200\,960$ & $17$ & 1, 2, 3, 4, 5, 6, 7, 8, 11, 14, 15, 23\\
$M_{24}$ & $244\,823\,040$ & $26$ & 1, 2, 3, 4, 5, 6, 7, 8, 10, 11, 12, 14, 15, 21, 23\\
$J_1$ & $175\,560$ & $15$ & 1, 2, 3, 5, 6, 7, 10, 11, 15, 19\\
$J_2$ & $604\,800$ & $21$ & 1, 2, 3, 4, 5, 6, 7, 8, 10, 12, 15\\
$\mathrm{HS}$ & $44\,352\,000$ & $24$ & 1, 2, 3, 4, 5, 6, 7, 8, 10, 11, 12, 15, 20\\
$\mathrm{McL}$ & $898\,128\,000$ & $24$ & 1, 2, 3, 4, 5, 6, 7, 8, 9, 10, 11, 12, 14, 15, 30\\
$\mathrm{Suz}$ & $448\,345\,497\,600$ & $43$ & 1, 2, 3, 4, 5, 6, 7, 8, 9, 10, 11, 12, 13, 14, 15, 18, 20, 21, 24\\
$\mathrm{He}$ & $4\,030\,387\,200$ & $33$ & 1, 2, 3, 4, 5, 6, 7, 8, 10, 12, 14, 15, 17, 21, 28\\
$\mathrm{Ru}$ & $145\,926\,144\,000$ & $36$ & 1, 2, 3, 4, 5, 6, 7, 8, 10, 12, 13, 14, 15, 16, 20, 24, 26, 29\\
\end{tabular}
\caption{Numbers of conjugacy classes and element orders.}\label{tab:classes}
\end{table}

\begin{table}[ht]
\centering\small
\begin{tabular}{lrllll}
class & $|C_G(g)|$ & $2$-power & $3$-power & $5$-power & $7$-power\\ \hline
1A & $604\,800$ & 1A & 1A & 1A & 1A\\
2A & $1920$ & 1A & 2A & 2A & 2A\\
2B & $240$ & 1A & 2B & 2B & 2B\\
3A & $1080$ & 3A & 1A & 3A & 3A\\
3B & $36$ & 3B & 1A & 3B & 3B\\
4A & $96$ & 2A & 4A & 4A & 4A\\
5A & $300$ & 5B & 5B & 1A & 5B\\
5B & $300$ & 5A & 5A & 1A & 5A\\
5C & $50$ & 5D & 5D & 1A & 5D\\
5D & $50$ & 5C & 5C & 1A & 5C\\
6A & $24$ & 3A & 2A & 6A & 6A\\
6B & $12$ & 3B & 2B & 6B & 6B\\
7A & $7$ & 7A & 7A & 7A & 1A\\
8A & $8$ & 4A & 8A & 8A & 8A\\
10A & $20$ & 5A & 10B & 2B & 10B\\
10B & $20$ & 5B & 10A & 2B & 10A\\
10C & $10$ & 5C & 10D & 2A & 10D\\
10D & $10$ & 5D & 10C & 2A & 10C\\
12A & $12$ & 6A & 4A & 12A & 12A\\
15A & $15$ & 15B & 5B & 3A & 15B\\
15B & $15$ & 15A & 5A & 3A & 15A\\
\end{tabular}
\caption{The $21$ classes of $J_2$ with centraliser orders and power maps.}\label{tab:j2classes}
\end{table}

\begin{remark}[Checks]\label{rem:classcheck}
For every group of Table~\ref{tab:classes} the class equation $\sum_i|G|/|C_G(g_i)|=|G|$ holds, and the list of triples
(element order, centraliser order, orders and centraliser orders of the $p$-th powers) agrees with the ATLAS table in the GAP
character table library~\cite{GAP,ATLAS}. The classes were computed by GAP from the permutation groups of the laws.
\end{remark}

\subsection{Burnside's algorithm in digital form}\label{ssec:burnside}
Let $C_1=\{0\},C_2,\dots,C_r$ be the conjugacy classes of a digital law $(\Ap N,\ast)$, given as sets of labels, and fix a
label $z_k\in C_k$ for each $k$. The \emph{class multiplication coefficients} are
\begin{equation}\label{eq:classcoef}
a_{jik}=\bigl|\{(x,y)\in C_j\times C_i : x\ast y=z_k\}\bigr|=\bigl|\{x\in C_j : \bar x\ast z_k\in C_i\}\bigr|,
\end{equation}
where $\bar x$ is the inverse label of $x$. The second form shows that the whole row $(a_{jik})_{i}$ costs $|C_j|$ inversions and
$|C_j|$ products, followed by a class lookup, so the matrix $M_j=(a_{jik})_{i,k}$ costs $r\,|C_j|$ evaluations of the law, and
all $r$ matrices cost $r\,N$ evaluations.

\begin{proposition}[Character table from the law]\label{prop:classalgebra}
The matrices $M_1,\dots,M_r$ commute. For every irreducible character $\chi$ the vector
\[
w_\chi=\bigl(\omega_\chi(C_k)\bigr)_{k},\qquad \omega_\chi(C_k)=\frac{|C_k|\,\chi(z_k)}{\chi(1)},
\]
satisfies $M_jw_\chi=\omega_\chi(C_j)\,w_\chi$ for all $j$, and every common eigenvector of the $M_j$ is a multiple of one of
the $r$ vectors $w_\chi$. Consequently, the common eigenvectors normalised by $w_1=1$ are exactly the $w_\chi$, and
\[
\chi(1)^2=\frac{N}{\sum_k|\omega_\chi(C_k)|^2/|C_k|},\qquad \chi(z_k)=\frac{\chi(1)\,\omega_\chi(C_k)}{|C_k|}.
\]
In particular the character table of $G$ is determined by the law, through the integers $a_{jik}$ alone.
\end{proposition}
\begin{proof}
Let $K_j=\sum_{x\in C_j}x$ be the class sums in the group algebra $\mathbb{C}G$. They form a basis of its centre $Z$, and
$K_jK_i=\sum_ka_{jik}K_k$, since the coefficient of $z_k$ in $K_jK_i$ is the number of pairs in~\eqref{eq:classcoef} and is
constant on $C_k$. Hence $M_j$ is the matrix of multiplication by $K_j$ on $Z$ in this basis (acting on coefficient vectors
from the right), and the $M_j$ commute because $Z$ is commutative. By Schur's lemma $K_j$ acts on an irreducible
representation affording $\chi$ as the scalar $\omega_\chi(C_j)$, so $\omega_\chi\colon Z\to\mathbb{C}$ is an algebra
homomorphism, $\omega_\chi(C_j)\omega_\chi(C_i)=\sum_ka_{jik}\omega_\chi(C_k)$, which is $M_jw_\chi=\omega_\chi(C_j)w_\chi$
read at row $i$. The $r$ vectors $w_\chi$ are linearly independent by the orthogonality relations, so $Z\cong\mathbb{C}^r$
splits into $r$ joint eigenspaces of dimension one, and these are spanned by the $w_\chi$. Finally
$\sum_k|\omega_\chi(C_k)|^2/|C_k|=\sum_k|C_k|\,|\chi(z_k)|^2/\chi(1)^2=N/\chi(1)^2$ by the first orthogonality relation.
\end{proof}

For $M_{11}$ this is steps (3) and (4) of Section~\ref{ssec:m11table}: $r=10$, and the $a_{jik}$ require $10\cdot7920$
evaluations of the law of Tables~\ref{tab:m11transA} and~\ref{tab:m11transB}.

\begin{remark}[Dixon's modular version]\label{rem:dixon}
The eigenvalues $\omega_\chi(C_j)$ are algebraic integers, and the $\chi(z_k)$ lie in the field of $e$-th roots of unity,
$e$ the exponent of $G$. Dixon~\cite{Dixon1967} replaces $\mathbb C$ by $\mathbb F_p$ with $p\equiv1\pmod e$ and
$p>2\sqrt N$: the $a_{jik}$ are reduced modulo $p$, the common eigenvectors are computed exactly over $\mathbb F_p$, and each value
$\chi(z_k)$ is recovered from its reductions and those of the power maps by a discrete Fourier transform over the cyclic
group $\langle z_k\rangle$. No floating-point arithmetic is needed. Usually one or two matrices $M_j$ already have $r$
distinct eigenvalues modulo $p$, and then only $r\,|C_j|$ evaluations of the law are required; choosing $C_j$ small is the
improvement of Schneider.
\end{remark}

\begin{remark}[Beyond enumeration]\label{rem:beyondenum}
The cost $r\,|C_j|$ is feasible from the laws of this book for $M_{12}$, $J_1$, $M_{22}$, $J_2$ and $M_{23}$, and with effort for
$M_{24}$ and $\mathrm{HS}$. For larger groups no class can be listed, and the digital law supplies other data instead. Uniform
random labels give class representatives and, through the power maps of Section~\ref{ssec:classes}, the class fusion;
fixed points of chain laws give permutation characters; rank-three actions give two irreducibles at once
(Proposition~\ref{prop:rank3}). Induction from subgroups with known tables, restriction, tensor products and
symmetrisations then produce further characters, and the orthogonality relations, together with integrality, fix the
rest. This is how most tables in the ATLAS~\cite{ATLAS} were obtained. For the largest groups, such as $\mathrm{B}$ and
$\mathrm{M}$, the realistic digital task is verification: every class representative, power map and centraliser order of
the printed table can be checked against the law.
\end{remark}

\subsection{The character table of \texorpdfstring{$M_{11}$}{M11} from its law}\label{ssec:m11table}
Tables~\ref{tab:m11transA} and~\ref{tab:m11transB} determine the law of $M_{11}$ completely: the element with label
$t_1+11t_2+110t_3+990t_4$ is $s_{1,t_1}\circ s_{2,t_2}\circ s_{3,t_3}\circ s_{4,t_4}$. From these $38$ permutations alone the
following were computed.
\begin{enumerate}
\item All $7920$ elements, indexed by their labels. The set is closed under composition with each $s_{l,t}$, so it is the group.
\item The conjugacy classes $C_1,\dots,C_{10}$, as sets of labels, their sizes and their power maps.
\item The class multiplication coefficients $a_{jik}=|\{(x,y)\in C_j\times C_i: x\circ y=z\}|$ for a fixed $z\in C_k$.
\item The central characters $\omega_\chi$, as the common eigenvectors of the matrices $(a_{jik})_{i,k}$ (Burnside's algorithm),
      normalised by $\omega_\chi(C_1)=1$; then $\chi(1)^2=|G|/\sum_k|\omega_\chi(C_k)|^2/|C_k|$ and
      $\chi(g_k)=\chi(1)\,\omega_\chi(C_k)/|C_k|$.
\end{enumerate}
The result is Table~\ref{tab:m11chars}, where $b_{11}=(-1+\sqrt{-11})/2$ and the bar denotes complex conjugation. The row ``least
label'' gives the smallest label in each class; the $p$-power rows give the class letter of $g^p$.

\begin{table}[ht]
\centering\small
\setlength{\tabcolsep}{4pt}
\begin{tabular}{lrrrrrrrrrr}
 & 1A & 2A & 3A & 4A & 5A & 6A & 8A & 8B & 11A & 11B\\ \hline
$|C_G(g)|$ & 7920 & 48 & 18 & 8 & 5 & 6 & 8 & 8 & 11 & 11\\
least label & 0 & 22 & 12 & 31 & 7 & 19 & 6 & 17 & 1 & 2\\
$2$ power & A & A & A & A & A & A & A & A & B & A\\
$3$ power & A & A & A & A & A & A & A & B & A & B\\
$5$ power & A & A & A & A & A & A & B & A & A & B\\
$11$ power & A & A & A & A & A & A & A & B & A & A\\
\hline
$\chi_{1}$ & $1$ & $1$ & $1$ & $1$ & $1$ & $1$ & $1$ & $1$ & $1$ & $1$\\
$\chi_{2}$ & $10$ & $-2$ & $1$ & $0$ & $0$ & $1$ & $-i\sqrt2$ & $i\sqrt2$ & $-1$ & $-1$\\
$\chi_{3}$ & $10$ & $-2$ & $1$ & $0$ & $0$ & $1$ & $i\sqrt2$ & $-i\sqrt2$ & $-1$ & $-1$\\
$\chi_{4}$ & $10$ & $2$ & $1$ & $2$ & $0$ & $-1$ & $0$ & $0$ & $-1$ & $-1$\\
$\chi_{5}$ & $11$ & $3$ & $2$ & $-1$ & $1$ & $0$ & $-1$ & $-1$ & $0$ & $0$\\
$\chi_{6}$ & $16$ & $0$ & $-2$ & $0$ & $1$ & $0$ & $0$ & $0$ & $\bar b_{11}$ & $b_{11}$\\
$\chi_{7}$ & $16$ & $0$ & $-2$ & $0$ & $1$ & $0$ & $0$ & $0$ & $b_{11}$ & $\bar b_{11}$\\
$\chi_{8}$ & $44$ & $4$ & $-1$ & $0$ & $-1$ & $1$ & $0$ & $0$ & $0$ & $0$\\
$\chi_{9}$ & $45$ & $-3$ & $0$ & $1$ & $0$ & $0$ & $-1$ & $-1$ & $1$ & $1$\\
$\chi_{10}$ & $55$ & $-1$ & $1$ & $-1$ & $0$ & $-1$ & $1$ & $1$ & $0$ & $0$\\
\end{tabular}
\caption{The character table of $M_{11}$, computed from Tables~\ref{tab:m11transA} and~\ref{tab:m11transB}.}\label{tab:m11chars}
\end{table}

\begin{example}[Classes of labels in $M_{11}$]\label{ex:m11classes}
The element with label $1234$ of Example~\ref{ex:m11dec} lies in class 6A, and $1234*4321=7231$ lies in class 11A.
\end{example}

\begin{remark}[Checks]\label{rem:m11tablecheck}
\emergencystretch=3em The script \texttt{scripts/\allowbreak characters.py} reads the $38$ permutations from the source of
Tables~\ref{tab:m11transA} and~\ref{tab:m11transB}, checks $s_{l,t}(b_l)=\delta_{l,t}$ and closure, and carries out (1)--(4).
Each computed $\omega_\chi$ is a common eigenvector of all ten matrices, every value is recognised as an integer, $\pm i\sqrt2$,
$b_{11}$ or $\bar b_{11}$, and the orthogonality relations hold. The degrees and centraliser orders agree with the ATLAS table of
$M_{11}$ in the GAP character table library~\cite{GAP,ATLAS}, which was used only for this comparison. The same script
produces Tables~\ref{tab:rank3}--\ref{tab:j2classes}, from the generator files and the output of
\texttt{scripts/\allowbreak data/\allowbreak characters.g}.
\end{remark}

\subsection{The character table of \texorpdfstring{$M_{12}$}{M12} from its law}\label{ssec:m12table}
The computation of Section~\ref{ssec:m11table} applies without change to $M_{12}$, of order $95\,040$, using only the representatives of Tables~\ref{tab:m12transA} and~\ref{tab:m12transB}: the $95\,040$ labels give the $15$ conjugacy classes, their power maps and the class multiplication coefficients, and the irreducible characters are the common eigenvectors of the class matrices. Classes with the same element order and centraliser order are named as in the ATLAS~\cite{ATLAS}. For $4A/4B$ and $8A/8B$ the names are fixed by requiring that the permutation character of the action on $\Ap{12}$ be $1a+11a$; the classes $11A$ and $11B$ are exchanged by a Galois automorphism and are named so that $\chi_4(11A)=b_{11}$. The rows are in ATLAS order. The result, Table~\ref{tab:m12chars}, agrees entry by entry with the GAP character table library~\cite{GAP,ATLAS}, which was used only for this comparison and for the class names; it is produced by \texttt{scripts/characters.py}.

\begin{table}[htbp]
\centering\scriptsize\setlength{\tabcolsep}{2.5pt}
\caption{The character table of $M_{12}$, computed from Tables~\ref{tab:m12transA} and~\ref{tab:m12transB}.}\label{tab:m12chars}
\begin{tabular}{l*{15}{r}}
 & 1A & 2A & 2B & 3A & 3B & 4A & 4B & 5A & 6A & 6B & 8A & 8B & 10A & 11A & 11B\\ \hline
$|C_G(g)|$ & 95040 & 240 & 192 & 54 & 36 & 32 & 32 & 10 & 12 & 6 & 8 & 8 & 10 & 11 & 11\\
least label & 0 & 11 & 24 & 13 & 10 & 71 & 33 & 9 & 76 & 8 & 38 & 6 & 17 & 1 & 2\\
$2$ power & A & A & A & A & B & B & B & A & B & A & A & B & A & B & A\\
$3$ power & A & A & B & A & A & A & B & A & A & B & A & B & A & A & B\\
$5$ power & A & A & B & A & B & A & B & A & A & B & A & B & A & A & B\\
$11$ power & A & A & B & A & B & A & B & A & A & B & A & B & A & A & A\\
\hline
$\chi_{1}$ & $1$ & $1$ & $1$ & $1$ & $1$ & $1$ & $1$ & $1$ & $1$ & $1$ & $1$ & $1$ & $1$ & $1$ & $1$\\
$\chi_{2}$ & $11$ & $-1$ & $3$ & $2$ & $-1$ & $-1$ & $3$ & $1$ & $-1$ & $0$ & $-1$ & $1$ & $-1$ & $0$ & $0$\\
$\chi_{3}$ & $11$ & $-1$ & $3$ & $2$ & $-1$ & $3$ & $-1$ & $1$ & $-1$ & $0$ & $1$ & $-1$ & $-1$ & $0$ & $0$\\
$\chi_{4}$ & $16$ & $4$ & $0$ & $-2$ & $1$ & $0$ & $0$ & $1$ & $1$ & $0$ & $0$ & $0$ & $-1$ & $b_{11}$ & $\bar b_{11}$\\
$\chi_{5}$ & $16$ & $4$ & $0$ & $-2$ & $1$ & $0$ & $0$ & $1$ & $1$ & $0$ & $0$ & $0$ & $-1$ & $\bar b_{11}$ & $b_{11}$\\
$\chi_{6}$ & $45$ & $5$ & $-3$ & $0$ & $3$ & $1$ & $1$ & $0$ & $-1$ & $0$ & $-1$ & $-1$ & $0$ & $1$ & $1$\\
$\chi_{7}$ & $54$ & $6$ & $6$ & $0$ & $0$ & $2$ & $2$ & $-1$ & $0$ & $0$ & $0$ & $0$ & $1$ & $-1$ & $-1$\\
$\chi_{8}$ & $55$ & $-5$ & $7$ & $1$ & $1$ & $-1$ & $-1$ & $0$ & $1$ & $1$ & $-1$ & $-1$ & $0$ & $0$ & $0$\\
$\chi_{9}$ & $55$ & $-5$ & $-1$ & $1$ & $1$ & $3$ & $-1$ & $0$ & $1$ & $-1$ & $-1$ & $1$ & $0$ & $0$ & $0$\\
$\chi_{10}$ & $55$ & $-5$ & $-1$ & $1$ & $1$ & $-1$ & $3$ & $0$ & $1$ & $-1$ & $1$ & $-1$ & $0$ & $0$ & $0$\\
$\chi_{11}$ & $66$ & $6$ & $2$ & $3$ & $0$ & $-2$ & $-2$ & $1$ & $0$ & $-1$ & $0$ & $0$ & $1$ & $0$ & $0$\\
$\chi_{12}$ & $99$ & $-1$ & $3$ & $0$ & $3$ & $-1$ & $-1$ & $-1$ & $-1$ & $0$ & $1$ & $1$ & $-1$ & $0$ & $0$\\
$\chi_{13}$ & $120$ & $0$ & $-8$ & $3$ & $0$ & $0$ & $0$ & $0$ & $0$ & $1$ & $0$ & $0$ & $0$ & $-1$ & $-1$\\
$\chi_{14}$ & $144$ & $4$ & $0$ & $0$ & $-3$ & $0$ & $0$ & $-1$ & $1$ & $0$ & $0$ & $0$ & $-1$ & $1$ & $1$\\
$\chi_{15}$ & $176$ & $-4$ & $0$ & $-4$ & $-1$ & $0$ & $0$ & $1$ & $-1$ & $0$ & $0$ & $0$ & $1$ & $0$ & $0$\\
\end{tabular}
\end{table}

\section{Exercises}
\begin{enumerate}
\item Show that $|M_{11}|=7920=11\cdot10\cdot9\cdot8$ and $|M_{12}|=95040=12\cdot11\cdot10\cdot9\cdot8$, and explain how these
products arise as the basic orbit lengths of a stabiliser chain.
\item Show that $|M_{24}|=24\cdot23\cdot22\cdot21\cdot20\cdot48=244823040$.
\item For a stabiliser chain of $M_{11}$ with basic orbit lengths $11,10,9,8$, write the mixed-radix label of an element and show
that the identity has label $0$.
\item Explain, in two distinct ways, how Dynkin diagrams enter the sporadic constructions: once through the deep holes of the Leech lattice and once through a Lie-type point stabiliser in Table~\ref{tab:primsporadic}. Why neither construction classifies a sporadic group by a Dynkin diagram?
\item Explain why fixing standard generators~\cite{ATLAS,AtlasRep} makes a canonical law independent of the software used to compute it.
\item Confirm with GAP that $M_{11}$ has $10$ conjugacy classes, and compare with the class data recomputed from the printed law.
\item Show that the order of the Monster lies between $8\cdot10^{53}$ and $10^{54}$, and deduce the number of decimal digits of its labels.
\end{enumerate}

\chapter{Worked digital products}\label{ch:products}

Chapter~\ref{ch:lie} gives digital laws for the finite simple groups of Lie type. Each law consists of a decoding map
$D:\Ap{N}\to G$, $N=|G|$, and its inverse $E$, and the product of two labels is
\begin{equation}\label{eq:digitalproduct}
 m*n=E\bigl(D(m)\,D(n)\bigr).
\end{equation}
This chapter computes~\eqref{eq:digitalproduct} by hand for four groups: $\PSL(2,7)$, $\PSL(3,2)$, $Sz(8)$ and ${}^2G_2(3)$. Section~\ref{sec:wpfurther} then records one product for each of the other families. Every
intermediate quantity is displayed, so each product can be checked with pencil and paper.

\section{A product in \texorpdfstring{$\PSL(2,7)$}{PSL(2,7)}}\label{sec:wppsl2}

We use the labelling of $\PSL(2,p)$ of Section~\ref{sec:psl2p}, written in the following explicit form. Let $p$ be an odd prime,
$h=(p-1)/2$, so that $|\PSL(2,p)|=ph(p+1)$. A label $m<ph$ decodes as
\begin{equation}\label{eq:psl2decupper}
 D(m)=\begin{pmatrix}a&b\\0&a^{-1}\end{pmatrix},\qquad a=\lfloor m/p\rfloor+1,\quad b=m\bmod p .
\end{equation}
A label $m\geq ph$ is written $m'=m-ph$ and decodes as $D(m)=\begin{pmatrix}a&b\\c&d\end{pmatrix}$ with
\begin{equation}\label{eq:psl2declower}
 a=\lfloor m'/p\rfloor\bmod p,\qquad c=\lfloor m'/p^2\rfloor+1,\qquad d=m'\bmod p,\qquad b=(ad-1)c^{-1}\bmod p .
\end{equation}
Conversely, given $A=\begin{pmatrix}a&b\\c&d\end{pmatrix}\in\SL(2,p)$, replace $A$ by $-A$ if $c>h$, or if $c=0$ and $a>h$.
Then
\begin{equation}\label{eq:psl2enc}
 E(A)=\begin{cases}(a-1)p+b,&c=0,\\ ph+(c-1)p^2+ap+d,&c\neq0.\end{cases}
\end{equation}

\begin{example}\label{ex:psl27product}
Let $p=7$, so $h=3$, $ph=21$ and $|\PSL(2,7)|=168$. We compute $145*66$. Both labels are at least $21$, so
\eqref{eq:psl2declower} applies. For $145$: $m'=124$, $a=\lfloor124/7\rfloor\bmod7=17\bmod7=3$,
$c=\lfloor124/49\rfloor+1=3$, $d=124\bmod7=5$, $c^{-1}=5$ and $b=(3\cdot5-1)\cdot5=70\equiv0$. For $66$: $m'=45$,
$a=6$, $c=1$, $d=3$ and $b=(6\cdot3-1)\cdot1=17\equiv3$. Hence
\[
 D(145)D(66)=\begin{pmatrix}3&0\\3&5\end{pmatrix}\begin{pmatrix}6&3\\1&3\end{pmatrix}
 =\begin{pmatrix}18&9\\23&24\end{pmatrix}\equiv\begin{pmatrix}4&2\\2&3\end{pmatrix}\pmod 7.
\]
Here $c=2\leq h$, so no sign change is needed, and~\eqref{eq:psl2enc} gives $E=21+1\cdot49+4\cdot7+3=101$. Thus
$145*66=101$ in $\PSL(2,7)$.
With $g$, $h$ the two factors and $gh$ their product, the same digital law gives the inverses
\begin{gather*}
 g^{-1}=137,\quad h^{-1}=50,\quad (gh)^{-1}=101,
\end{gather*}
of orders $3$, $7$ and $2$; the identity has label~$0$, and $(gh)^{-1}=h^{-1}*g^{-1}$.
\end{example}

\section{A product in \texorpdfstring{$\PSL(3,2)$}{PSL(3,2)}}\label{sec:wppsln}

For $\PSL_n(q)$ (Section~\ref{sec:psln}) the label is a rank. Let $d=\gcd(n,q-1)$ and
$\mu_d=\{z\in\mathbb F_q^\times:z^n=1\}$. For $A\in\SL_n(q)$ put
\begin{equation}\label{eq:pslnkey}
 \kappa(A)=\sum_{0\leq i,j<n}\bigl((A_{ij}-\delta_{ij})\bmod q\bigr)\,q^{\,i+nj}.
\end{equation}
The canonical representative of the class $A\mu_d$ is the $zA$, $z\in\mu_d$, of least key, and $E(A)$ is the number of
canonical representatives whose key is smaller than that of $A$. So $D(m)$ is the canonical representative with exactly $m$
smaller keys.

\begin{example}\label{ex:psl32product}
In $\PSL(3,2)$ we have $d=1$, so every matrix of $\SL_3(2)$ is its own canonical representative, and there are $168$ keys.
The labels $122$ and $63$ are the matrices with keys $398$ and $173$:
\[
 D(122)D(63)=\begin{pmatrix}1&1&0\\1&1&1\\1&0&0\end{pmatrix}\begin{pmatrix}0&1&0\\0&1&1\\1&1&1\end{pmatrix}
 =\begin{pmatrix}0&2&1\\1&3&2\\0&1&0\end{pmatrix}\equiv\begin{pmatrix}0&0&1\\1&1&0\\0&1&0\end{pmatrix}\pmod2 .
\]
The nonzero entries of the product minus the identity are at $(i,j)=(0,0)$, $(1,0)$, $(2,1)$, $(0,2)$ and $(2,2)$, so
\eqref{eq:pslnkey} gives $\kappa=2^0+2^1+2^5+2^6+2^8=355$. Exactly $105$ of the $168$ keys are smaller than $355$;
hence $122*63=105$ in $\PSL(3,2)$. Unlike the other two examples, this last step needs the sorted list of all keys,
which is not reproduced here.
With $g$, $h$ the two factors and $gh$ their product, the same digital law gives the inverses
\begin{gather*}
 g^{-1}=163,\quad h^{-1}=119,\quad (gh)^{-1}=149,
\end{gather*}
of orders $3$, $3$ and $7$; the identity has label~$0$, and $(gh)^{-1}=h^{-1}*g^{-1}$.
\end{example}

\section{A product in \texorpdfstring{$Sz(8)$}{Sz(8)}}\label{sec:wpsz}

We use Definition~\ref{def:suzuki} and Theorem~\ref{thm:suzuki} with $m=1$, so $q=8$, $x^\theta=x^4$, $2^m=2$ and
$M(\lambda)=\mathrm{diag}(\lambda^3,\lambda^2,\lambda^{-2},\lambda^{-3})$. Indices of matrix entries start at $0$, as in
Chapter~\ref{ch:lie}. The field is $\mathbb F_8=\mathbb F_2[x]/(x^3+x+1)$, and $c_0+c_1x+c_2x^2$ is written as the integer
$c_0+2c_1+4c_2$, which is also its label $[\,\cdot\,]$. Addition is the carry-free sum $\oplus$ (XOR), and the product
$\otimes$ is given by the table
\[
\begin{array}{c|ccccccc}
\otimes&1&2&3&4&5&6&7\\\hline
1&1&2&3&4&5&6&7\\
2&2&4&6&3&1&7&5\\
3&3&6&5&7&4&1&2\\
4&4&3&7&6&2&5&1\\
5&5&1&4&2&7&3&6\\
6&6&7&1&5&3&2&4\\
7&7&5&2&1&6&4&3
\end{array}
\]
The primitive element is $\gamma=x=2$, with $\gamma^0,\dots,\gamma^6=1,2,4,3,6,7,5$. Here $N_B=q^2(q-1)=448$ and
$|Sz(8)|=29120$. A Borel label $L<N_B$ is $L=e+7([a]+8[b])$ and decodes to $S(a,b)M(\gamma^e)$. A label $L\geq N_B$ is
$L=N_B+L_0+N_B([a_0]+8[b_0])$ with $L_0<N_B$, and decodes to $S(a_0,b_0)\,T\,D(L_0)$.

\begin{example}\label{ex:sz8product}
We compute $20803*7554$.

\emph{Decoding.} For $20803$: $20803-448=20355=45\cdot448+195$ and $45=5+8\cdot5$, so $(a_0,b_0)=(5,5)$ and $L_0=195$.
Then $195=6+7\cdot27$ and $27=3+8\cdot3$, so $L_0$ is $S(3,3)M(\gamma^6)=S(3,3)M(5)$. For $7554$:
$7554-448=7106=15\cdot448+386$ and $15=7+8\cdot1$, so $(a_0,b_0)=(7,1)$; then $386=1+7\cdot55$ and $55=7+8\cdot6$, so
$L_0$ is $S(7,6)M(2)$. Multiplying out,
\begin{gather*}
 D(20803)=S(5,5)\,T\,S(3,3)M(5)=\begin{pmatrix}5&7&7&3\\6&5&2&4\\5&4&1&4\\1&1&0&1\end{pmatrix},\\
 D(7554)=S(7,1)\,T\,S(7,6)M(2)=\begin{pmatrix}7&0&3&6\\2&2&5&4\\0&5&5&6\\4&4&5&7\end{pmatrix}.
\end{gather*}

\emph{Multiplication.} Entry $(i,j)$ of the product is $\bigoplus_k D(20803)_{ik}\otimes D(7554)_{kj}$. For example
$(0,0)$ is $5\otimes7\oplus7\otimes2\oplus7\otimes0\oplus3\otimes4=6\oplus5\oplus0\oplus7=4$. Listing the four products
for each entry,
\[
\begin{array}{ll}
(0,0):6\oplus5\oplus0\oplus7=4&(0,1):0\oplus5\oplus6\oplus7=4\\
(0,2):4\oplus6\oplus6\oplus4=0&(0,3):3\oplus1\oplus4\oplus2=4\\
(1,0):4\oplus1\oplus0\oplus6=3&(1,1):0\oplus1\oplus1\oplus6=6\\
(1,2):1\oplus7\oplus1\oplus2=5&(1,3):2\oplus2\oplus7\oplus1=6\\
(2,0):6\oplus3\oplus0\oplus6=3&(2,1):0\oplus3\oplus5\oplus6=0\\
(2,2):4\oplus2\oplus5\oplus2=1&(2,3):3\oplus6\oplus6\oplus1=2\\
(3,0):7\oplus2\oplus0\oplus4=1&(3,1):0\oplus2\oplus0\oplus4=6\\
(3,2):3\oplus5\oplus0\oplus5=3&(3,3):6\oplus4\oplus0\oplus7=5
\end{array}
\]
so
\[
 g=D(20803)D(7554)=\begin{pmatrix}4&4&0&4\\3&6&5&6\\3&0&1&2\\1&6&3&5\end{pmatrix}.
\]

\emph{Encoding.} Since $g_{03}=4\neq0$, $g$ lies in the big cell. With $4^{-1}=7$,
$a_0=g_{13}\otimes g_{03}^{-1}=6\otimes7=4$ and $b_0=g_{23}\otimes g_{03}^{-1}=2\otimes7=5$. By
Proposition~\ref{prop:suzukiborel}, $S(4,5)^{-1}=S(4,5\oplus4^{1+\theta})=S(4,5\oplus3)=S(4,6)$, and
\[
 g'=T\,S(4,6)\,g=\begin{pmatrix}7&0&0&0\\0&2&0&0\\5&0&5&0\\4&4&0&4\end{pmatrix}=S(0,2)M(\gamma^4).
\]
Indeed $a=g'_{10}\otimes g'^{-1}_{00}=0$, $b=g'_{20}\otimes g'^{-1}_{00}=5\otimes4=2$, and
$(g'_{11})^\theta=2^4=6=\gamma^4$ because $(\lambda^{2})^\theta=\lambda^{8}=\lambda$. So $L_0=4+7(0+8\cdot2)=116$, and
\[
 E(g)=448+116+448\,(4+8\cdot5)=20276 .
\]
Thus $20803*7554=20276$ in $Sz(8)$.
With $g$, $h$ the two factors and $gh$ their product, the same digital law gives the inverses
\begin{gather*}
 g^{-1}=16981,\quad h^{-1}=2787,\quad (gh)^{-1}=4477,
\end{gather*}
of orders $4$, $13$ and $7$; the identity has label~$0$, and $(gh)^{-1}=h^{-1}*g^{-1}$.
\end{example}

\section{A product in \texorpdfstring{${}^2G_2(3)$}{2G2(3)}}\label{sec:wpree}

We use Definition~\ref{def:ree}, Proposition~\ref{prop:reeU} and Theorem~\ref{thm:ree} with $m=0$, so $q=3$, $\theta=1$ and
$x^{3\theta}=x^3=x$ on $\mathbb F_3$. All arithmetic is modulo $3$, with entries written as $0,1,2$ (so $-1=2$). As noted in
Example~\ref{ex:ree}, ${}^2G_2(3)\cong P\Gamma L(2,8)$ is not simple, but its law has exactly the form of the law for
$q=27$, and every step can be checked by hand. Here $\gamma=2$, $N_B=q^3(q-1)=54$ and $|{}^2G_2(3)|=1512$. A Borel label is
$L=e+2([t]+3[u]+9[v])$ and decodes to $x_S(t,u,v)\,h(\gamma^e)$. A label $L\geq N_B$ is
$L=N_B+L_0+N_B([t]+3[u]+9[v])$ with $L_0<N_B$, and decodes to $x_S(t,u,v)\,n\,D(L_0)$.

\begin{example}\label{ex:ree3product}
We compute $1000*777$.

\emph{Decoding.} For $1000$: $1000-54=946=17\cdot54+28$ and $17=2+3\cdot2+9\cdot1$, so $(t,u,v)=(2,2,1)$ and $L_0=28$. Then
$28=0+2\cdot14$ and $14=2+3\cdot1+9\cdot1$, so $L_0$ is $x_S(2,1,1)\,h(1)$. For $777$: $777-54=723=13\cdot54+21$ and
$13=1+3\cdot1+9\cdot1$, so $(t,u,v)=(1,1,1)$; then $21=1+2\cdot10$ and $10=1+3\cdot0+9\cdot1$, so $L_0$ is $x_S(1,0,1)\,h(2)$.
Multiplying out,
\begin{gather*}
 D(1000)=x_S(2,2,1)\,n\,x_S(2,1,1)={\scriptstyle\setlength{\arraycolsep}{3pt}\begin{pmatrix}2&2&2&2&1&2&1\\1&2&2&1&0&0&2\\1&0&2&0&2&2&2\\0&1&0&0&0&1&1\\0&2&0&2&1&1&1\\2&0&2&0&1&1&2\\2&1&1&2&1&0&2\end{pmatrix}},\\
 D(777)=x_S(1,1,1)\,n\,x_S(1,0,1)\,h(2)={\scriptstyle\setlength{\arraycolsep}{3pt}\begin{pmatrix}1&0&2&2&2&1&1\\1&1&2&0&1&0&1\\0&1&1&1&2&0&1\\0&2&2&0&1&1&1\\2&2&0&2&0&1&0\\2&0&1&1&0&0&1\\1&2&0&1&1&2&1\end{pmatrix}}.
\end{gather*}

\emph{Multiplication.} The integer product of these two matrices and its reduction modulo $3$ are
\begin{gather*}
 {\scriptstyle\setlength{\arraycolsep}{3pt}\begin{pmatrix}11&12&16&11&13&7&11\\5&10&10&6&11&6&8\\11&10&6&12&8&7&7\\4&3&3&2&2&2&3\\7&10&9&4&5&5&6\\8&8&7&11&10&7&7\\7&12&11&9&11&9&8\end{pmatrix}}\\
\equiv {\scriptstyle\setlength{\arraycolsep}{3pt}\begin{pmatrix}2&0&1&2&1&1&2\\2&1&1&0&2&0&2\\2&1&0&0&2&1&1\\1&0&0&2&2&2&0\\1&1&0&1&2&2&0\\2&2&1&2&1&1&1\\1&0&2&0&2&0&2\end{pmatrix}}=g .
\end{gather*}

\emph{Encoding.} Since $c=g_{60}=1\neq0$, $g$ lies in the big cell, and Theorem~\ref{thm:ree} gives
$t=-g_{50}/c=-2=1$, $u=g_{40}/c-t^{2}=1-1=0$ and $v=tu-g_{30}/c=0-1=2$. With $n^{-1}=n$ and
$x_S(1,0,2)^{-1}=x_S(2,2,2)$,
\[
 g'=n^{-1}x_S(1,0,2)^{-1}g={\scriptstyle\setlength{\arraycolsep}{3pt}\begin{pmatrix}2&0&1&0&1&0&1\\0&1&0&1&0&2&0\\0&0&2&0&0&0&2\\0&0&0&1&0&1&0\\0&0&0&0&2&0&2\\0&0&0&0&0&1&0\\0&0&0&0&0&0&2\end{pmatrix}} .
\]
The entry $g'_{00}=2$ gives $\gamma^e=2^{3}=2$, so $e=1$. The matrix $g'h(2)^{-1}$ has $x_{12}=0$, $x_{46}=1$ and
$x_{25}=0$, and Proposition~\ref{prop:reeU} gives $t=0$, $u=x_{46}-t^2=1$ and $v=x_{25}-tu=0$. Hence
$g'=x_S(0,1,0)\,h(2)$, $L_0=1+2(0+3\cdot1+9\cdot0)=7$, and
\[
 E(g)=54+7+54\,(1+3\cdot0+9\cdot2)=1087 .
\]
Thus $1000*777=1087$ in ${}^2G_2(3)$.
With $g$, $h$ the two factors and $gh$ their product, the same digital law gives the inverses
\begin{gather*}
 g^{-1}=1244,\quad h^{-1}=459,\quad (gh)^{-1}=411,
\end{gather*}
of orders $9$, $9$ and $9$; the identity has label~$0$, and $(gh)^{-1}=h^{-1}*g^{-1}$.
\end{example}

\section{A product in \texorpdfstring{$M_{24}$}{M24}}\label{sec:wpm24}

We use the Golay law of Theorem~\ref{thm:golaylaw}, in the model implemented by \texttt{golay\_law.py}: the base octad is
$O_0=\{0,1,2,3,4,7,10,12\}$, the remaining $16$ points form $C_0$, identified with $\mathbb F_2^4$ by a fixed coordinate
map, and octads are numbered with $O_0$ first and the others by increasing bit mask. Vectors of $\mathbb F_2^4$ and the columns
$c_1,\dots,c_4$ of $A\in GL_4(2)$ are written as integers $0,\dots,15$ in binary, and $Av=\bigoplus_{k:\,\dg{2}{k-1}{v}=1}c_k$
is a carry-free sum of columns. A label is
\[
\Lambda=o+759\,(v+16\mu),\qquad o<759,\quad v<16,\quad \mu=d_1+15d_2+210d_3+2520d_4<20160,
\]
and the columns are recovered from the digits $(d_1,d_2,d_3,d_4)$ one at a time: $c_k$ is the unique $c$ outside
$\langle c_1,\dots,c_{k-1}\rangle$ such that exactly $d_k$ integers $y<c$ lie outside that span. Permutations are given as
tables $x\mapsto g(x)$ and compose from right to left, $(gh)(x)=g(h(x))$.

\begin{example}\label{ex:m24kproduct}
We compute $1216677*24291795$, a product in the octad stabiliser $K\cong2^4{:}GL_4(2)$.

\emph{Decoding.} $1216677=759\cdot1603$ and $1603=3+16\cdot100$, so $o=0$, $v_1=3$ and $\mu=100=10+15\cdot6$, with digits
$(10,6,0,0)$. Hence $c_1=11$ (the ten integers $1,\dots,10$ lie outside $\{0\}$), $c_2=7$ (the six integers $1,\dots,6$ lie
outside $\{0,11\}$), $c_3=1$ and $c_4=2$, so $A_1=(11,7,1,2)$. Likewise $24291795=759\cdot32005$, $32005=5+16\cdot2000$ and
$2000=5+15\cdot7+210\cdot9$, so $o=0$, $v_2=5$ and $A_2=(6,9,12,1)$.

\emph{Multiplication.} By part~5 of Theorem~\ref{thm:golaylaw} the product lies in $K$, with translation part
$v_1\oplus A_1v_2$ and linear part $A_1A_2$. Since $5=2^0+2^2$, $A_1v_2=c_1\oplus c_3=11\oplus1=10$, and $v=3\oplus10=9$. The
columns of $A_1A_2$ are the images under $A_1$ of $6,9,12,1$:
\[
A_1\,6=7\oplus1=6,\qquad A_1\,9=11\oplus2=9,\qquad A_1\,12=1\oplus2=3,\qquad A_1\,1=11,
\]
so $A_1A_2=(6,9,3,11)$.

\emph{Encoding.} For $(6,9,3,11)$: $d_1=5$; $d_2=7$ (the $y<9$ outside $\{0,6\}$); $d_3=2$ ($y=1,2$, outside
$\{0,6,9,15\}$); $d_4=5$ ($y=1,2,4,7,8$, outside $\langle6,9,3\rangle=\{0,3,5,6,9,10,12,15\}$). Hence
$\mu=5+15\cdot7+210\cdot2+2520\cdot5=13130$ and
\[
1216677*24291795=759\,(9+16\cdot13130)=159\,457\,551 .
\]
With $g$, $h$ the two factors and $gh$ their product, the same digital law gives the inverses
\begin{gather*}
 g^{-1}=210258180,\quad h^{-1}=152921802,\quad (gh)^{-1}=205818030,
\end{gather*}
of orders $4$, $5$ and $14$; the identity has label~$0$, and $(gh)^{-1}=h^{-1}*g^{-1}$.
\end{example}

\begin{example}\label{ex:m24product}
We compute $181025939*188240803$. These are the Golay labels of the elements $a$, $b$ of Examples~\ref{ex:m24mul}
and~\ref{ex:golaym24}, so $m=a$ and $n=b$. The two elements $m=\Lambda^{-1}(181025939)$ and $n=\Lambda^{-1}(188240803)$ are
\begin{gather*}
\setlength{\arraycolsep}{2.5pt}
\begin{array}{c|cccccccccccc}
x&0&1&2&3&4&5&6&7&8&9&10&11\\\hline
m(x)&21&6&16&17&23&1&10&4&22&13&2&3\\
n(x)&8&9&2&21&7&18&3&22&11&23&6&0
\end{array}\\[4pt]
\setlength{\arraycolsep}{2.5pt}
\begin{array}{c|cccccccccccc}
x&12&13&14&15&16&17&18&19&20&21&22&23\\\hline
m(x)&14&19&9&8&5&7&12&18&0&15&20&11\\
n(x)&14&10&12&1&20&13&15&19&16&17&4&5
\end{array}
\end{gather*}

\emph{Decoding.} $181025939=644+759\cdot238505$, $238505=9+16\cdot14906$ and $14906=11+15\cdot13+210\cdot10+2520\cdot5$.
So $m(O_0)=\{2,4,6,14,16,17,21,23\}$ is octad number $644$, and $s_O^{-1}m$ acts on $C_0$ with $v=9$ and $A=(12,15,13,9)$.
Similarly $188240803=454+759\cdot248011$, $248011=11+16\cdot15500$ and $15500=5+15\cdot11+210\cdot1+2520\cdot6$, so
$n(O_0)=\{2,6,7,8,9,14,21,22\}$ is octad number $454$, with $v=11$ and $A=(6,13,2,12)$. Both octads can be read off the
tables.

\emph{Multiplication.} Composing the tables,
\[
(mn)(O_0)=m\big(\{2,6,7,8,9,14,21,22\}\big)=\{4,9,10,13,15,16,20,22\},
\]
which is octad number $422$. Unlike Example~\ref{ex:m24kproduct}, the octad moves, so the coset representative changes. With
$O$ this octad, $s_O^{-1}(mn)$ acts on $C_0$ with $v=13$ and $A=(5,13,6,7)$; this step is done by the script.

\emph{Encoding.} For $(5,13,6,7)$ the digits are $(4,11,4,3)$, so $\mu=4+15\cdot11+210\cdot4+2520\cdot3=8569$ and
\[
181025939*188240803=422+759\,(13+16\cdot8569)=104\,072\,225 .
\]
As an independent check, $\Lambda^{-1}(104\,072\,225)$ is the composite $mn$ of the two tables, it passes the membership
test for $M_{24}$, and it maps $O_0$ to octad number $422$. Its label in the law of Example~\ref{ex:m24mul} is $155\,026\,102$, as there.
With $g$, $h$ the two factors and $gh$ their product, the same digital law gives the inverses
\begin{gather*}
 g^{-1}=35094818,\quad h^{-1}=136435904,\quad (gh)^{-1}=59504718,
\end{gather*}
of orders $6$, $6$ and $7$; the identity has label~$0$, and $(gh)^{-1}=h^{-1}*g^{-1}$.
\end{example}

\section{A product in \texorpdfstring{$\Co_0$}{Co0}}\label{sec:wpco0}

We use the chain law of Theorem~\ref{thm:co0chain}: $\lambda(g)=\sum_i d_i\,r_1\cdots r_{i-1}$ with the radices $r_i$
of that theorem and top digit $d_1=\nu(g(v_0))$, where $v_0=(4^2,0^{22})$. Both factors lie in the
monomial subgroup $2^{12}{:}M_{24}$ of Section~\ref{ssec:conwaymonomial}. We write $(C,\pi)$ for the map
$x\mapsto\varepsilon_C(\pi x)$, where $(\pi x)_{\pi(i)}=x_i$ and $\varepsilon_C$ changes the signs of the coordinates in the
Golay codeword $C$, so that $(C,\pi)(C',\pi')=(C\oplus\pi(C'),\pi\pi')$. The $1104$ minimal vectors of shape
$(\pm4,\pm4,0^{22})$ carry the first $\nu$-numbers: the vector with nonzero entries at positions $i<j$ has $\nu=4p+s$, where
$p$ is the index of $\{i,j\}$ among the pairs in lexicographic order and $s=0,1,2,3$ for the signs $++$, $+-$, $-+$, $--$.

\begin{example}\label{ex:co0product}
Let $m=(\emptyset,\pi_0)$ with $\pi_0(x)=x+1$ for $x\le21$, $\pi_0(22)=0$, $\pi_0(23)=23$, and let $n=(O_0,\pi_1)$ with
$O_0=\{0,1,2,3,4,7,10,12\}$ and $\pi_1$ the second permutation below; $\pi_0$ and $\pi_1$ are the generators of $M_{24}$ used
by the script. The bottom rows give $\pi_0\pi_1$:
\begin{gather*}
\setlength{\arraycolsep}{2.5pt}
\begin{array}{c|cccccccccccc}
x&0&1&2&3&4&5&6&7&8&9&10&11\\\hline
\pi_1(x)&0&1&16&12&3&5&8&17&2&6&11&22\\
\pi_0\pi_1(x)&1&2&17&13&4&6&9&18&3&7&12&0
\end{array}\\[4pt]
\setlength{\arraycolsep}{2.5pt}
\begin{array}{c|cccccccccccc}
x&12&13&14&15&16&17&18&19&20&21&22&23\\\hline
\pi_1(x)&13&18&19&14&9&10&4&21&15&20&7&23\\
\pi_0\pi_1(x)&14&19&20&15&10&11&5&22&16&21&8&23
\end{array}
\end{gather*}

\emph{Labels.} Sifting gives $\lambda(m)=92$, with digits $(92,0,\dots,0)$, and
\[
\lambda(n)=1\,761\,173\,012\,250\,288\,003,\qquad (d_1,\dots,d_{12})=(3,0,28,6,10,21,4,0,1,1,0,0).
\]
We compute $92*1\,761\,173\,012\,250\,288\,003$.

\emph{Top digit by hand.} $\pi_1$ fixes $0$ and $1$, and both lie in $O_0$, so $n(v_0)=-v_0$ and $\nu=3$, the top digit of
$\lambda(n)$. Then $\pi_0$ moves the coordinates $0,1$ to $1,2$, so $(mn)(v_0)=m(-v_0)=(0,-4,-4,0^{21})$. The pair $\{1,2\}$
follows the $23$ pairs containing $0$, so it has index $23$, and $\nu=4\cdot23+3=95$. In the same way $m(v_0)=(0,4,4,0^{21})$ has
$\nu=92=\lambda(m)$.

\emph{Product.} By the multiplication rule, $mn=(\pi_0(O_0),\pi_0\pi_1)=(\{1,2,3,4,5,8,11,13\},\pi_0\pi_1)$. Sifting $mn$
gives the digits $(95,0,28,6,10,21,4,0,1,1,0,0)$: only the top digit differs from that of $n$, which the script finds and the
hand computation does not explain. Hence
\begin{align*}
92*1\,761\,173\,012\,250\,288\,003&=\lambda(n)+(95-3)\
&=1\,761\,173\,012\,250\,288\,095 .
\end{align*}
As independent checks, the $24\times24$ matrix of $mn$ is the product of the matrices of $m$ and $n$ and preserves the Leech
lattice, and $\lambda^{-1}$ of the result is that matrix.
With $g$, $h$ the two factors and $gh$ their product, the same digital law gives the inverses
\begin{gather*}
 g^{-1}=4188172403882572644,\\h^{-1}=2695423001193468003,\\(gh)^{-1}=3918390902126260602,
\end{gather*}
of orders $23$, $10$ and $10$; the identity has label~$0$, and $(gh)^{-1}=h^{-1}*g^{-1}$.
\end{example}

\section{A product in \texorpdfstring{$J_4$}{J4}}\label{sec:wpj4}

We use the transitive law of Proposition~\ref{prop:j4}. Thus $J_4$ acts on $N=173\,067\,389$ cosets of
$K=2^{11}{:}M_{24}$, and, with the Schreier representative $s_t$ taking $0$ to $t$,
\[
\lambda(s_t h)=t+N\lambda_K(h),\qquad h\in K.
\]
The $K$-label is $u+2048m$, where $u\in\mathbb F_2^{11}$ is written as an integer and $m$ is the $M_{24}$ chain label, with
radices $(24,23,22,21,20,16,3)$. The entries below are independently checked by decoding to the $112\times112$ matrices over
$\mathbb F_2$ on which the implementation acts on the right.

\begin{example}\label{ex:j4kproduct}
We first multiply two elements of $K$:
\[
14341*22537=334\,265\,695\,087.
\]
Writing a $K$-label as $u+2048m$, the two factors and the product have respectively
\begin{center}
\begin{tabular}{c|r|r|c}
 &u&m&\text{digits of }m\\ \hline
14341&5&7&(7,0,0,0,0,0,0)\\
22537&9&11&(11,0,0,0,0,0,0)\\
334\,265\,695\,087&879&163\,215\,671&(23,12,0,0,0,0,2)
\end{tabular}
\end{center}
The semidirect-product calculation is $u\oplus\phi(P,v)=5\oplus\phi(P,9)=879$ and $7*11=163\,215\,671$ in the
$M_{24}$ chain law. Both nontrivial operations are carried out by the script. Decoding the three labels, the product of the two
$112\times112$ matrices has label $334\,265\,695\,087$; in the reverse order its label is $334\,265\,618\,038$.
With $g$, $h$ the two factors and $gh$ their product, the same digital law gives the inverses
\begin{gather*}
 g^{-1}=168725449292,\quad h^{-1}=209464003429,\quad (gh)^{-1}=100853913541,
\end{gather*}
of orders $6$, $30$ and $6$; the identity has label~$0$, and $(gh)^{-1}=h^{-1}*g^{-1}$.
\end{example}

\begin{example}\label{ex:j4product}
For a full product, take the deterministic Schreier representatives $s_1,s_2$ from the coset tree used by the law, and put
$g=s_1h_{14341}$ and $h=s_2h_{22537}$, where $h_k\in K$ has $K$-label $k$. The full labels, split as $(t,\lambda_K)$, are
\begin{center}
\begin{tabular}{c|r|r|c}
 &t&$\lambda_K$&\text{digits of the $M_{24}$ part}\\ \hline
g&110\,466\,560&341\,912\,014\,057&(1,13,2,5,5,11,1)\\
h&132\,772\,346&315\,969\,325\,215&(15,19,9,16,18,12,1)\\
gh&13\,859\,126&203\,295\,128\,335&(7,22,3,16,0,2,1)
\end{tabular}
\end{center}
Thus
\begin{align*}
59\,173\,819\,540\,686\,753\,733*54\,683\,986\,119\,184\,685\,981
 &=35\,183\,757\,057\,372\,226\,441.
\end{align*}
The law decodes the right-hand label to the independently multiplied $112\times112$ matrix $gh$; the three coset images are,
in order, $110\,466\,560$, $132\,772\,346$ and $13\,859\,126$, agreeing with the displayed first digits.
With $g$, $h$ the two factors and $gh$ their product, the same digital law gives the inverses
\begin{gather*}
 g^{-1}=29200872966818338589,\\h^{-1}=36251388162944076884,\\(gh)^{-1}=23632336167617021182,
\end{gather*}
of orders $21$, $40$ and $35$; the identity has label~$0$, and $(gh)^{-1}=h^{-1}*g^{-1}$.
\end{example}

\section{A product in \texorpdfstring{$\mathrm{Th}$}{Th}}\label{sec:wpth}

We use Proposition~\ref{prop:th}. The coset action of $\mathrm{Th}$ on $\Ap{143\,127\,000}$ has stabiliser
$K={}^3D_4(2){:}3$, and the label is
\[
\lambda(s_t h)=t+143\,127\,000\,\lambda_K(h),\qquad h\in K.
\]
The canonical $K$-law has radices $(819,512,126,12)$. The implementation decodes a $K$-label through the cached Schreier
lift on the $819$ long-root involutions, and all matrix checks below use the faithful $248\times248$ matrices over $\mathbb F_2$.

\begin{example}\label{ex:thkproduct}
In $K$ we have
\[
7*11=52\,836\,297.
\]
The digit tuples in radices $(819,512,126,12)$ are
\begin{center}
\begin{tabular}{c|c}
\text{label}&\text{digit tuple}\\ \hline
7&(7,0,0,0)\\
11&(11,0,0,0)\\
52\,836\,297&(150,1,0,1)
\end{tabular}
\end{center}
Decode the first two labels to $248\times248$ matrices and multiply them. Re-encoding the product gives $52\,836\,297$.
Decoding this result and re-encoding also gives the same label.
With $g$, $h$ the two factors and $gh$ their product, the same digital law gives the inverses
\begin{gather*}
 g^{-1}=423084754,\quad h^{-1}=450778420,\quad (gh)^{-1}=2726453,
\end{gather*}
of orders $12$, $21$ and $12$; the identity has label~$0$, and $(gh)^{-1}=h^{-1}*g^{-1}$.
\end{example}

\begin{example}\label{ex:thproduct}
Take the deterministic Schreier representatives $s_1,s_2$ in the coset tree and put $g=s_1h_7$, $h=s_2h_{11}$, where $h_k\in K$
has $K$-label $k$. The full labels, split as $(t,\lambda_K)$, are
\begin{center}
\begin{tabular}{c|r|r|c}
 &t&$\lambda_K$&\text{digits of the $K$ part}\\ \hline
g&99\,905\,087&275\,285\,750&(194,252,26,5)\\
h&97\,525\,418&500\,995\,900&(496,388,60,9)\\
gh&57\,370\,703&381\,845\,588&(761,313,28,7)
\end{tabular}
\end{center}
Thus
\begin{align*}
39\,400\,823\,640\,155\,087*71\,706\,040\,276\,825\,418
 &=54\,652\,413\,531\,046\,703.
\end{align*}
The law decodes the right-hand label to the independently multiplied $248\times248$ matrix $gh$. The three coset images are,
in order, $99\,905\,087$, $97\,525\,418$ and $57\,370\,703$, agreeing with the first digits; all three matrices have rank $248$.
With $g$, $h$ the two factors and $gh$ their product, the same digital law gives the inverses
\begin{gather*}
 g^{-1}=60554851585758001,\\h^{-1}=64518562919340003,\\(gh)^{-1}=25153297361894852,
\end{gather*}
of orders $30$, $39$ and $12$; the identity has label~$0$, and $(gh)^{-1}=h^{-1}*g^{-1}$.
\end{example}

\section{A product in the Monster}\label{sec:wpmonster}
The law of $\mathbb M$ is the one of Theorem~\ref{thm:mlaw}, evaluated with \texttt{monster/monster\_law.py} and
\texttt{mmgroup}~$1.0.8$. In the notation of that script a label $m\in[0,|\mathbb M|)$ is split as
\[
m=s+2^{24}c+|K|\,r_2+|2.\mathbb B|\,r_1,\qquad s=\delta+2^{12}\sigma+2^{13}i,
\]
where $K=C_{G_{x0}}(x_\beta)=2^{1+23}.\Co_2$ is the group of Theorem~\ref{thm:klaw}, $s<2^{24}$ labels its normal
subgroup $Q_K$, $c$ is the $\Co_2$ digit, $r_2<|2.\mathbb B|/|K|=11\,707\,448\,673\,375$ ranks the second axis and
$r_1<|\mathbb M|/|2.\mathbb B|=97\,239\,461\,142\,009\,186\,000$ ranks the top axis, which lies in one of the twelve
$G_{x0}$-orbits of Remark~\ref{rem:monsterstatus}. In the last digit, $D(s)=x_{-1}^{\sigma}x_{d_i}x_\delta$, with $\delta$
in the cocode and $d_i$ the $i$-th of the $2048$ Golay code words orthogonal to $\beta$ (in \texttt{mmgroup} numbering).

\begin{example}\label{ex:monsterkproduct}
Take $g,h\in Q_K$ with labels $40\,971$ and $78\,020$, so that $(i,\sigma,\delta)=(5,0,\mathtt{0xb})$ and
$(9,1,\mathtt{0xc4})$, with $d_5=\mathtt{0x5}$ and $d_9=\mathtt{0x9}$. Since $Q_K/\langle x_{-1}\rangle$ is elementary
abelian, the Golay and cocode parts of $gh$ are $d_5+d_9=\mathtt{0xc}=d_{12}$ and $\mathtt{0xb}+\mathtt{0xc4}=\mathtt{0xcf}$.
The sign digit is $\sigma_g+\sigma_h=1$ plus a cocycle term, which the script evaluates to $0$. Hence
\[
40\,971\times 78\,020=207+2^{12}+2^{13}\cdot 12=102\,607 .
\]
Outside $Q_K$ the digits mix. The labels with $(s,c)=(7,11)$ and $(13,17)$ multiply as
\[
184\,549\,383\times 285\,212\,685=441\,975\,592=\mathtt{0x328}+2^{13}\cdot 704+2^{24}\cdot 26,
\]
so the product has $Q_K$ digit $s=5\,767\,976$, with Golay word $d_{704}=\mathtt{0x4c0}$, although both $Q_K$ digits of
the factors are pure cocode. That the $\Co_2$ digits satisfy $26=11+17$ here is not a rule of the chain; this second
product is checked by the script only.
With $g$, $h$ the two factors and $gh$ their product, the same digital law gives the inverses
\begin{gather*}
 g^{-1}=45067,\quad h^{-1}=78020,\quad (gh)^{-1}=102607,
\end{gather*}
of orders $4$, $2$ and $2$; the identity has label~$0$, and $(gh)^{-1}=h^{-1}*g^{-1}$.
For the second product, with the two displayed factors as $g$ and $h$, the same law gives
\[
 g^{-1}=185\,865\,026,\qquad
 h^{-1}=293\,606\,829,\qquad
 (gh)^{-1}=450\,480\,057,
\]
all of order $4$; again the identity has label~$0$, and $(gh)^{-1}=h^{-1}*g^{-1}$.
\end{example}

\begin{example}\label{ex:monsterproduct}
Take the labels
\begin{align*}
g&=3\,883\,560\,514\,439\,200\,112\,108\,283\,576\,015\,423\,818\,665\,947\,100\,822\,717,\\
h&=150\,264\,023\,038\,597\,408\,710\,715\,747\,403\,918\,362\,067\,408\,929\,425\,227\,244.
\end{align*}
Then
\[
gh=212\,726\,864\,516\,606\,513\,234\,534\,111\,091\,187\,695\,562\,559\,932\,184\,978\,387,
\]
and the digits of the three labels are
\begin{center}\footnotesize
\begin{tabular}{c|r|c|r|r|c}
&$r_1$&orbit&$r_2$&$c$&$(i,\sigma,\delta)$\\ \hline
g&$467\,360\,381\,284\,460\,968$&$8B$&$5\,914\,237\,542\,798$&$7\,338\,569\,562\,816$&$(537,1,\mathtt{0xcbd})$ \\
h&$18\,083\,264\,272\,449\,995\,475$&$10B$&$11\,343\,439\,527\,797$&$14\,487\,947\,735\,453$&$(483,1,\mathtt{0xdec})$ \\
gh&$25\,600\,247\,025\,965\,477\,049$&$10B$&$4\,104\,949\,158\,811$&$24\,279\,590\,191\,741$&$(231,0,\mathtt{0x3d3})$
\end{tabular}
\end{center}
The top axes of $g$, $h$ and $gh$ lie in the orbits $8B$, $10B$ and $10B$. Apart from this split into digits nothing
here is checkable by hand. The script confirms $D(gh)=D(g)D(h)$ in \texttt{mmgroup}, that $D^{-1}\circ D$ fixes the three
labels, and that the orbit of each top digit agrees with the axis type of the corresponding element; it runs in a few
seconds.
With $g$, $h$ the two factors and $gh$ their product, the same digital law gives the inverses
\begin{gather*}
 g^{-1}=6480686252459502959001190146825060570307975337187794,\\h^{-1}=484924031314715268257926709690560352236124175526441501,\\(gh)^{-1}=659954164485034725521133517690909800403004918029924256,
\end{gather*}
of orders $78$, $71$ and $35$; the identity has label~$0$, and $(gh)^{-1}=h^{-1}*g^{-1}$.
\end{example}

\section{Further products}\label{sec:wpfurther}

For the remaining families of Chapter~\ref{ch:lie} the decoded elements are too large to write out, so only the labels are
recorded. For each family we take the smallest group with a stored computation. Each product was computed with the law of the
cited theorem, implemented in \texttt{scripts/}, and in each run decoding followed by encoding returns the label, the identity has label $0$, all labels lie in $\Ap{N}$, the product agrees with the matrix product and associativity holds on sample triples.
Note that ${}^2F_4(2)$ is not simple; it contains the Tits group with index $2$ (Theorem~
ef{thm:titsdigital}).
{\emergencystretch=3em
\begin{itemize}
\item $G_2(5)$, Theorem~\ref{thm:g2cells}, $|G|=5\,\allowbreak 859\,\allowbreak 000\,\allowbreak 000$:\\ $m=1\,\allowbreak 997\,\allowbreak 283\,\allowbreak 091$,\\ $n=1\,\allowbreak 002\,\allowbreak 184\,\allowbreak 965$,\\ $m*n=5\,\allowbreak 512\,\allowbreak 323\,\allowbreak 731$.
\item ${}^3D_4(2)$, Theorem~\ref{thm:3d4cells}, $|G|=211\,\allowbreak 341\,\allowbreak 312$:\\ $m=53\,\allowbreak 227\,\allowbreak 708$,\\ $n=30\,\allowbreak 544\,\allowbreak 696$,\\ $m*n=26\,\allowbreak 197\,\allowbreak 305$.
\item $F_4(2)$, Theorem~\ref{thm:f4closed}, $|G|=3\,\allowbreak 311\,\allowbreak 126\,\allowbreak 603\,\allowbreak 366\,\allowbreak 400$:\\ $m=944\,\allowbreak 437\,\allowbreak 105\,\allowbreak 928\,\allowbreak 969$,\\ $n=737\,\allowbreak 683\,\allowbreak 288\,\allowbreak 020\,\allowbreak 099$,\\ $m*n=529\,\allowbreak 825\,\allowbreak 386\,\allowbreak 857\,\allowbreak 885$.
\item ${}^2F_4(2)$, Theorem~\ref{thm:2f4digital}, $|G|=35\,\allowbreak 942\,\allowbreak 400$:\\ $m=6\,\allowbreak 803\,\allowbreak 900$,\\ $n=1\,\allowbreak 352\,\allowbreak 053$,\\ $m*n=9\,\allowbreak 296\,\allowbreak 917$.
\item $E_6(2)$, Theorem~\ref{thm:e6law}, $|G|=214\,\allowbreak 841\,\allowbreak 575\,\allowbreak 522\,\allowbreak 005\,\allowbreak 575\,\allowbreak 270\,\allowbreak 400$:\\ $m=162\,\allowbreak 276\,\allowbreak 505\,\allowbreak 379\,\allowbreak 933\,\allowbreak 334\,\allowbreak 990\,\allowbreak 320$,\\ $n=116\,\allowbreak 533\,\allowbreak 813\,\allowbreak 217\,\allowbreak 642\,\allowbreak 933\,\allowbreak 313\,\allowbreak 473$,\\ $m*n=138\,\allowbreak 968\,\allowbreak 748\,\allowbreak 442\,\allowbreak 346\,\allowbreak 481\,\allowbreak 901\,\allowbreak 966$.
\item ${}^2E_6(2)$, Theorem~\ref{thm:2e6law}, $|G|=76\,\allowbreak 532\,\allowbreak 479\,\allowbreak 683\,\allowbreak 774\,\allowbreak 853\,\allowbreak 939\,\allowbreak 200$:\\ $m=16\,\allowbreak 651\,\allowbreak 950\,\allowbreak 189\,\allowbreak 205\,\allowbreak 936\,\allowbreak 290\,\allowbreak 549$,\\ $n=45\,\allowbreak 344\,\allowbreak 061\,\allowbreak 244\,\allowbreak 829\,\allowbreak 118\,\allowbreak 473\,\allowbreak 207$,\\ $m*n=33\,\allowbreak 686\,\allowbreak 555\,\allowbreak 911\,\allowbreak 656\,\allowbreak 403\,\allowbreak 558\,\allowbreak 842$.
\item $E_7(2)$, Theorem~\ref{thm:e7law}, $|G|=7\,\allowbreak 997\,\allowbreak 476\,\allowbreak 042\,\allowbreak 075\,\allowbreak 799\,\allowbreak 759\,\allowbreak 100\,\allowbreak 487\,\allowbreak 262\,\allowbreak 680\,\allowbreak 802\,\allowbreak 918\,\allowbreak 400$:\\ $m=5\,\allowbreak 854\,\allowbreak 567\,\allowbreak 007\,\allowbreak 228\,\allowbreak 584\,\allowbreak 891\,\allowbreak 849\,\allowbreak 280\,\allowbreak 362\,\allowbreak 194\,\allowbreak 909\,\allowbreak 567\,\allowbreak 872$,\\ $n=6\,\allowbreak 971\,\allowbreak 517\,\allowbreak 499\,\allowbreak 971\,\allowbreak 688\,\allowbreak 817\,\allowbreak 909\,\allowbreak 076\,\allowbreak 669\,\allowbreak 494\,\allowbreak 367\,\allowbreak 794\,\allowbreak 041$,\\ $m*n=452\,\allowbreak 163\,\allowbreak 067\,\allowbreak 631\,\allowbreak 932\,\allowbreak 687\,\allowbreak 657\,\allowbreak 396\,\allowbreak 694\,\allowbreak 981\,\allowbreak 606\,\allowbreak 823\,\allowbreak 804$.
\item $E_8(2)$, Theorem~\ref{thm:e8law}, $|G|=337\,\allowbreak 804\,\allowbreak 753\,\allowbreak 143\,\allowbreak 634\,\allowbreak 806\,\allowbreak 261\,\allowbreak 388\,\allowbreak 190\,\allowbreak 614\,\allowbreak 085\,\allowbreak 595\,\allowbreak 079\,\allowbreak 991\,\allowbreak 692\,\allowbreak 242\,\allowbreak 467\,\allowbreak 651\,\allowbreak 576\,\allowbreak 160\,\allowbreak 959\,\allowbreak 909\,\allowbreak 068\,\allowbreak 800\,\allowbreak 000$:\\ $m=291\,\allowbreak 187\,\allowbreak 596\,\allowbreak 337\,\allowbreak 441\,\allowbreak 092\,\allowbreak 718\,\allowbreak 387\,\allowbreak 068\,\allowbreak 629\,\allowbreak 079\,\allowbreak 436\,\allowbreak 355\,\allowbreak 555\,\allowbreak 169\,\allowbreak 104\,\allowbreak 222\,\allowbreak 352\,\allowbreak 355\,\allowbreak 579\,\allowbreak 247\,\allowbreak 994\,\allowbreak 788\,\allowbreak 881\,\allowbreak 123$,\\ $n=238\,\allowbreak 866\,\allowbreak 409\,\allowbreak 991\,\allowbreak 748\,\allowbreak 292\,\allowbreak 375\,\allowbreak 426\,\allowbreak 148\,\allowbreak 542\,\allowbreak 606\,\allowbreak 723\,\allowbreak 246\,\allowbreak 351\,\allowbreak 189\,\allowbreak 813\,\allowbreak 991\,\allowbreak 037\,\allowbreak 796\,\allowbreak 058\,\allowbreak 703\,\allowbreak 580\,\allowbreak 994\,\allowbreak 412\,\allowbreak 695$,\\ $m*n=139\,\allowbreak 054\,\allowbreak 127\,\allowbreak 712\,\allowbreak 954\,\allowbreak 738\,\allowbreak 336\,\allowbreak 597\,\allowbreak 222\,\allowbreak 112\,\allowbreak 998\,\allowbreak 242\,\allowbreak 493\,\allowbreak 020\,\allowbreak 790\,\allowbreak 460\,\allowbreak 287\,\allowbreak 161\,\allowbreak 683\,\allowbreak 184\,\allowbreak 717\,\allowbreak 145\,\allowbreak 069\,\allowbreak 110\,\allowbreak 632$.
\item $\Omega_5(3)$, Theorem~\ref{thm:classicallaw}, $|G|=25\,\allowbreak 920$:\\ $m=10\,\allowbreak 229$,\\ $n=10\,\allowbreak 587$,\\ $m*n=10\,\allowbreak 881$.
\item $\Omega^+_6(2)$, Theorem~\ref{thm:classicallaw}, $|G|=20\,\allowbreak 160$:\\ $m=2\,\allowbreak 600$,\\ $n=9\,\allowbreak 842$,\\ $m*n=9\,\allowbreak 596$.
\end{itemize}\par}

\begin{remark}[Checks]\label{rem:wpchecks}
Any product, in any group of the archive, is computed by the command
\begin{center}\ttfamily\small python3 multiply.py GROUP m n [--check]\end{center}
run from the directory \texttt{scripts} of the archive \texttt{DigitalFiniteGroupTheory}. It prints \texttt{GROUP: m * n = mn}. The labels are integers with $0\le m,n<|G|$; a label out of range is rejected. The names are those of the examples: \texttt{Ree3} (the law of Example~\ref{ex:ree3product}), \texttt{M24}, \texttt{Co0}, \texttt{J4} and \texttt{J4K} (its subgroup $K=2^{11}{:}M_{24}$), \texttt{Th} and \texttt{ThK} (its subgroup $K={}^3D_4(2){:}3$), \texttt{Monster}, and every group of Lie type with a digital law in the archive, such as \texttt{PSL2(7)}, \texttt{Sz8} or \texttt{E8(2)}; the full list is printed by \texttt{python3 multiply.py --list}. Case and punctuation are ignored, so \texttt{psl27} also works. Plain \texttt{python3} suffices except for the Monster and the Baby Monster, which need \texttt{mmgroup}~1.0.8 and hence \texttt{../.venv-mm/bin/python3}; that interpreter also works for every other group. With \texttt{--check}, the product is recomputed independently, by multiplying matrices for $^2G_2(3)$, $J_4$, $\mathrm{Th}$ and their subgroups~$K$ and words in \texttt{mmgroup} for the Monster, and the two results are compared. For example,
\begin{center}\ttfamily\small ../.venv-mm/bin/python3 multiply.py Monster 28 37 --check\end{center}
prints \texttt{Monster: 28 * 37 = 57} and \texttt{check: ok}.

The groups of Chapter~\ref{ch:sporadic} are available as well. The names \texttt{J1}, \texttt{J2}, \texttt{HS}, \texttt{Suz}, \texttt{He}, \texttt{Ru}, \texttt{J3}, \texttt{ON}, \texttt{HN}, \texttt{Fi24} and \texttt{Ly} use the canonical labels of Section~\ref{sec:canonical}; with \texttt{--check} the lexicographic chain is rebuilt from the generators in reverse order and with another random seed, which by Proposition~\ref{prop:lexchain} must give the same product. Appending \texttt{-chain}, as in \texttt{HN-chain}, selects instead the chain labels printed in the section of the group (Sections~\ref{sec:fischer}--\ref{sec:hn}); these are the only labels of \texttt{Fi22}, \texttt{Fi23} and \texttt{Fi24'}. The names \texttt{M11}, \texttt{M12}, \texttt{M22}, \texttt{M23} and \texttt{M24-chain} use the chains of Proposition~\ref{prop:mathieuchain} with the representatives of Section~\ref{ssec:mathieuexamples}, so that \texttt{multiply.py M11 1234 5678} prints \texttt{2524} as in Example~\ref{ex:m11mul} and \texttt{multiply.py M24-chain 123456789 98765432} prints \texttt{155026102} as in Example~\ref{ex:m24mul}; plain \texttt{M24} is the law of Example~\ref{ex:m24kproduct}. Finally \texttt{B} is the Baby Monster, computed with \texttt{mmgroup}. Most of these laws are built in a second or two; \texttt{HN}, \texttt{Fi24}, their \texttt{-chain} versions and \texttt{Fi24'} take about $15$ seconds, and \texttt{Ly} about half a minute.

Chapter~\ref{ch:lie} adds the following names. \texttt{2F4(2)'}, or \texttt{Tits}, is the Tits group on $\Ap{p_0}$ with $p_0=17\,971\,200$, and \texttt{2F4(2)-Tits} is ${}^2F_4(2)$ on $\Ap{2p_0}$ with the coset digit, both as in Theorem~\ref{thm:titsdigital}; plain \texttt{2F4(2)} keeps the normal-form law of Theorem~\ref{thm:2f4digital}. For these two names \texttt{--check} recomputes the product with the matrices of $F_4(2)$; for instance \texttt{multiply.py Tits 9999999 17971199 --check} prints \texttt{8766842} and \texttt{check: ok}. The suffix \texttt{-dense}, as in \texttt{PSL3(3)-dense}, selects the law of Theorem~\ref{thm:pslndense} for any $n\ge2$ and prime power $q\le1024$, and the plain names \texttt{PSL2(9)}, \texttt{PSL3(4)}, \texttt{PSL4(2)} use it whenever the book gives no other law, as for the groups of Example~\ref{ex:pslndense}; here \texttt{--check} compares the decoded product with the matrix product up to a central scalar, and \texttt{multiply.py PSL3(4) 12345 6789} prints \texttt{4911}. Finally \texttt{PSL2(p)} accepts every prime $p\ge5$ (Theorem~\ref{thm:psl2}), \texttt{Sz(q)} every $q=2^{2m+1}\ge8$ (Theorem~\ref{thm:suzuki}) and \texttt{Ree(q)} every $q=3^{2m+1}\ge27$ (Theorem~\ref{thm:ree}), for example \texttt{Sz(512)} or \texttt{PSL(2,17)}. The groups $\Omega^-_{2n}(q)$ of Example~\ref{ex:twisted} have no implemented law and are not available.

The products of this chapter are recomputed as follows, where $g$ and $h$ are the two labels displayed in the example:

\begingroup\small\raggedright\setlength{\tabcolsep}{3pt}
\begin{center}
\begin{tabular}{@{}p{0.14\linewidth}p{0.5\linewidth}p{0.3\linewidth}@{}}
\textbf{Example}&\textbf{Command}&\textbf{Product}\\ \hline
Ex.~\ref{ex:psl27product}&\texttt{multiply.py} \allowbreak{}\texttt{PSL2(7)} \allowbreak{}\texttt{145} \allowbreak{}\texttt{66}&\texttt{101}\\
Ex.~\ref{ex:psl32product}&\texttt{multiply.py} \allowbreak{}\texttt{PSL3(2)} \allowbreak{}\texttt{122} \allowbreak{}\texttt{63}&\texttt{105}\\
Ex.~\ref{ex:sz8product}&\texttt{multiply.py} \allowbreak{}\texttt{Sz8} \allowbreak{}\texttt{20803} \allowbreak{}\texttt{7554}&\texttt{20276}\\
Ex.~\ref{ex:ree3product}&\texttt{multiply.py} \allowbreak{}\texttt{Ree3} \allowbreak{}\texttt{1000} \allowbreak{}\texttt{777}&\texttt{1087}\\
Ex.~\ref{ex:m24kproduct}&\texttt{multiply.py} \allowbreak{}\texttt{M24} \allowbreak{}\texttt{1216677} \allowbreak{}\texttt{24291795}&\texttt{159457551}\\
Ex.~\ref{ex:m24product}&\texttt{multiply.py} \allowbreak{}\texttt{M24} \allowbreak{}\texttt{181025939} \allowbreak{}\texttt{188240803}&\texttt{104072225}\\
Ex.~\ref{ex:co0product}&\texttt{multiply.py} \allowbreak{}\texttt{Co0} \allowbreak{}\texttt{92} \allowbreak{}\texttt{1761173012250288003}&\texttt{1761173012250288095}\\
Ex.~\ref{ex:j4kproduct}&\texttt{multiply.py} \allowbreak{}\texttt{J4K} \allowbreak{}\texttt{14341} \allowbreak{}\texttt{22537}&\texttt{334265695087}\\
Ex.~\ref{ex:j4product}&\texttt{multiply.py} \allowbreak{}\texttt{J4} \allowbreak{}$g$ \allowbreak{}$h$&\texttt{35183757057372226441}\\
Ex.~\ref{ex:thkproduct}&\texttt{multiply.py} \allowbreak{}\texttt{ThK} \allowbreak{}\texttt{7} \allowbreak{}\texttt{11}&\texttt{52836297}\\
Ex.~\ref{ex:thproduct}&\texttt{multiply.py} \allowbreak{}\texttt{Th} \allowbreak{}$g$ \allowbreak{}$h$&\texttt{54652413531046703}\\
Ex.~\ref{ex:monsterkproduct}&\texttt{multiply.py} \allowbreak{}\texttt{Monster} \allowbreak{}\texttt{40971} \allowbreak{}\texttt{78020}&\texttt{102607}\\
Ex.~\ref{ex:monsterproduct}&\texttt{multiply.py} \allowbreak{}\texttt{Monster} \allowbreak{}$g$ \allowbreak{}$h$&as in the example\\
\end{tabular}
\end{center}
\endgroup

The intermediate quantities of the examples, namely digits, coset points, $K$-labels and matrices, are printed by the step-by-step drivers, run in the same way:
\begin{center}\ttfamily\small worked\_lie\_steps.py GROUP [m n],\quad worked\_m24.py [m n],\\ worked\_co0.py,\quad worked\_j4.py,\quad worked\_th.py,\quad worked\_monster.py.\end{center} The last four take no arguments and reproduce the examples only. Without labels, the first driver uses the pair stored in
\begin{center}\ttfamily\small ../drafts/\allowbreak{}worked\_lie\_<GROUP>.json,\end{center}
and it writes \texttt{../drafts/\allowbreak{}worked\_steps\_<GROUP>.txt}. Finally,
\begin{center}\ttfamily\small ../.venv-mm/bin/python3\allowbreak{} worked\_inverses.py\end{center}
computes, for every product $gh$ above, the labels of $g^{-1}$, $h^{-1}$ and $(gh)^{-1}$ and the element orders. It asserts $gg^{-1}=g^{-1}g=e$, $(g^{-1})^{-1}=g$ and $(gh)^{-1}=h^{-1}g^{-1}$ in the digital law, writes nothing in the project, and takes about $20$ seconds.
\end{remark}

\begin{remark}[Inverses]\label{rem:wpinverses}
The inverse of any label, in any group of the archive, is computed by the command
\begin{center}\ttfamily\small python3 inverse.py GROUP m [m2 ...] [--order] [--check]\end{center}
run, like \texttt{multiply.py}, from the directory \texttt{scripts} of the archive \texttt{DigitalFiniteGroupTheory}. It prints \texttt{GROUP: m\textasciicircum-1 = mi} for each label. The names, the range $0\le m<|G|$ and the interpreters are those of Remark~\ref{rem:wpchecks}, since \texttt{inverse.py} obtains its laws from \texttt{multiply.py}: plain \texttt{python3} suffices except for the Monster and the Baby Monster, which need \texttt{../.venv-mm/bin/python3}. Several labels may be given at once, and the law is then built only once. In \texttt{zsh} a name containing parentheses must be quoted, as in \texttt{\textquotesingle PSL2(7)\textquotesingle}, or written without punctuation, as in \texttt{psl27}. Called without arguments, \texttt{python3 inverse.py} prints the numbered list of groups, asks for a group, by number or by name, and then for labels one at a time, until an empty line.

The inverse is computed in the digital law itself, without decoding. The powers $m,m^2,\dots$ are formed until $m^{k+1}=m$, for at most $256$ products; then $k$ is the order of $m$, $e=m^k$ is the label of the identity and $m^{-1}=m^{k-1}$. For elements of larger order, which occur in groups of Lie type, Lagrange's theorem gives $m^{-1}=m^{|G|-1}$, computed by square-and-multiply with about $2\log_2|G|$ products, so that even in the Monster an inverse costs a few hundred products. The identity label $e$ is printed only when $e\neq0$, which happens only for the plain law \texttt{2F4(2)}. With \texttt{--order} the order of $m$ is printed as well, in the second case from the factorization of $|G|$. With \texttt{--check} the script asserts $m\,m^{-1}=m^{-1}m=e$, $e\,m=m\,e=m$ and $(m^{-1})^{-1}=m$ in the digital law and, where \texttt{multiply.py} has an independent decoder (matrices for ${}^2G_2(3)$, $J_4$, $\mathrm{Th}$ and their subgroups~$K$, words in \texttt{mmgroup} for the Monster), recomputes $m\,m^{-1}$ through it. In interactive mode the order and the checks are always on. For example,
\begin{center}\ttfamily\small ../.venv-mm/bin/python3 inverse.py Monster $g$ --order --check\end{center}
with $g$ the first label of Example~\ref{ex:monsterproduct} prints the label $g^{-1}$ given below, the order $78$ and \texttt{check: ok (law and independent decoder)}. The script writes nothing in the project.

The inverses of the factors and of the products of this chapter are recomputed as follows, where $g$, $h$ and $gh$ are the labels of the example; the last column lists the orders of $g$, $h$ and $gh$:

\begingroup\small\raggedright\setlength{\tabcolsep}{3pt}
\begin{center}
\resizebox{\textwidth}{!}{%
\begin{tabular}{@{}p{0.13\linewidth}p{0.34\linewidth}p{0.41\linewidth}p{0.1\linewidth}@{}}
\textbf{Example}&\textbf{Command}&\textbf{$g^{-1}$, $h^{-1}$, $(gh)^{-1}$}&\textbf{Orders}\\ \hline
Ex.~\ref{ex:psl27product}&\texttt{inverse.py} \allowbreak{}\texttt{psl27} \allowbreak{}\texttt{145} \allowbreak{}\texttt{66} \allowbreak{}\texttt{101}&\texttt{137}, \texttt{50}, \texttt{101}&$3,7,2$\\
Ex.~\ref{ex:psl32product}&\texttt{inverse.py} \allowbreak{}\texttt{psl32} \allowbreak{}\texttt{122} \allowbreak{}\texttt{63} \allowbreak{}\texttt{105}&\texttt{163}, \texttt{119}, \texttt{149}&$3,3,7$\\
Ex.~\ref{ex:sz8product}&\texttt{inverse.py} \allowbreak{}\texttt{Sz8} \allowbreak{}\texttt{20803} \allowbreak{}\texttt{7554} \allowbreak{}\texttt{20276}&\texttt{16981}, \texttt{2787}, \texttt{4477}&$4,13,7$\\
Ex.~\ref{ex:ree3product}&\texttt{inverse.py} \allowbreak{}\texttt{Ree3} \allowbreak{}\texttt{1000} \allowbreak{}\texttt{777} \allowbreak{}\texttt{1087}&\texttt{1244}, \texttt{459}, \texttt{411}&$9,9,9$\\
Ex.~\ref{ex:m24kproduct}&\texttt{inverse.py} \allowbreak{}\texttt{M24} \allowbreak{}\texttt{1216677} \allowbreak{}\texttt{24291795} \allowbreak{}\texttt{159457551}&\texttt{210258180}, \allowbreak{}\texttt{152921802}, \allowbreak{}\texttt{205818030}&$4,5,14$\\
Ex.~\ref{ex:m24product}&\texttt{inverse.py} \allowbreak{}\texttt{M24} \allowbreak{}\texttt{181025939} \allowbreak{}\texttt{188240803} \allowbreak{}\texttt{104072225}&\texttt{35094818}, \allowbreak{}\texttt{136435904}, \allowbreak{}\texttt{59504718}&$6,6,7$\\
Ex.~\ref{ex:co0product}&\texttt{inverse.py} \allowbreak{}\texttt{Co0} \allowbreak{}$g$ \allowbreak{}$h$ \allowbreak{}$gh$&\texttt{4188172403882572644}, \allowbreak{}\texttt{2695423001193468003}, \allowbreak{}\texttt{3918390902126260602}&$23,10,10$\\
Ex.~\ref{ex:j4kproduct}&\texttt{inverse.py} \allowbreak{}\texttt{J4K} \allowbreak{}\texttt{14341} \allowbreak{}\texttt{22537} \allowbreak{}\texttt{334265695087}&\texttt{168725449292}, \allowbreak{}\texttt{209464003429}, \allowbreak{}\texttt{100853913541}&$6,30,6$\\
Ex.~\ref{ex:j4product}&\texttt{inverse.py} \allowbreak{}\texttt{J4} \allowbreak{}$g$ \allowbreak{}$h$ \allowbreak{}$gh$&\texttt{29200872966818338589}, \allowbreak{}\texttt{36251388162944076884}, \allowbreak{}\texttt{23632336167617021182}&$21,40,35$\\
Ex.~\ref{ex:thkproduct}&\texttt{inverse.py} \allowbreak{}\texttt{ThK} \allowbreak{}\texttt{7} \allowbreak{}\texttt{11} \allowbreak{}\texttt{52836297}&\texttt{423084754}, \allowbreak{}\texttt{450778420}, \allowbreak{}\texttt{2726453}&$12,21,12$\\
Ex.~\ref{ex:thproduct}&\texttt{inverse.py} \allowbreak{}\texttt{Th} \allowbreak{}$g$ \allowbreak{}$h$ \allowbreak{}$gh$&\texttt{60554851585758001}, \allowbreak{}\texttt{64518562919340003}, \allowbreak{}\texttt{25153297361894852}&$30,39,12$\\
Ex.~\ref{ex:monsterproduct}&\texttt{inverse.py} \allowbreak{}\texttt{Monster} \allowbreak{}$g$ \allowbreak{}$h$ \allowbreak{}$gh$&see below&$78,71,35$\\
\end{tabular}}
\end{center}
\endgroup
For the Monster,
$g^{-1}=6\,\allowbreak480\,\allowbreak686\,\allowbreak252\,\allowbreak459\,\allowbreak502\,\allowbreak959\,\allowbreak001\,\allowbreak190\,\allowbreak146\,\allowbreak825\,\allowbreak060\,\allowbreak570\,\allowbreak307\,\allowbreak975\,\allowbreak337\,\allowbreak187\,\allowbreak794$,
$h^{-1}=484\,\allowbreak924\,\allowbreak031\,\allowbreak314\,\allowbreak715\,\allowbreak268\,\allowbreak257\,\allowbreak926\,\allowbreak709\,\allowbreak690\,\allowbreak560\,\allowbreak352\,\allowbreak236\,\allowbreak124\,\allowbreak175\,\allowbreak526\,\allowbreak441\,\allowbreak501$ and
$(gh)^{-1}=659\,\allowbreak954\,\allowbreak164\,\allowbreak485\,\allowbreak034\,\allowbreak725\,\allowbreak521\,\allowbreak133\,\allowbreak517\,\allowbreak690\,\allowbreak909\,\allowbreak800\,\allowbreak403\,\allowbreak004\,\allowbreak918\,\allowbreak029\,\allowbreak924\,\allowbreak256$.
In every row $(gh)^{-1}=h^{-1}g^{-1}$, which \texttt{multiply.py} confirms; the same table is asserted by \texttt{worked\_inverses.py}.
\end{remark}

\section{Exercises}
\begin{enumerate}
\item In the row of Example~\ref{ex:psl27product} the factors $g=145$ and $h=66$ of $\mathrm{PSL}_2(7)$ have orders $3$ and $7$, and
$gh=101$ has order $2$. Without running \texttt{inverse.py}, explain why the product of $145$ with itself must be $137$ and why
$101$ is its own inverse.
\item Show that $g$ and $h$ of the previous exercise generate $\mathrm{PSL}_2(7)$.
\item In the row of Example~\ref{ex:sz8product}, $g$, $h$ and $gh$ have orders $4$, $13$ and $7$. Using that the maximal
subgroups of $\mathrm{Sz}(8)$ have orders $448$, $52$, $20$ and $14$, show that $g$ and $h$ generate $\mathrm{Sz}(8)$.
\item Show that in any group $g$ and $g^{-1}$ have the same order, that $(gh)^{-1}=h^{-1}g^{-1}$, and that $gh$ and $hg$ have the same
order. Which column of the table of inverses can be predicted from the orders alone?
\item The inverse in \texttt{inverse.py} is computed from the powers $m,m^2,\dots$ of $m$. Show that if $k$ is the least positive
integer with $m^{k+1}=m$, then $k$ is the order of $m$ and $m^{-1}=m^{k-1}$. For the Monster product $gh$, of order $35$, how many
products does this take, and how many are needed if $(gh)^{34}$ is computed by repeated squaring? What are the orders of $(gh)^5$
and $(gh)^7$?
\end{enumerate}

\chapter{Software and reproducibility}\label{ch:software}
\begingroup\raggedright\emergencystretch=3em

This chapter describes the software archive accompanying this book, deposited at Zenodo under \href{https://doi.org/10.5281/zenodo.23881253}{doi:10.5281/zenodo.23881253}.  The archive is intended to let an external reader inspect the algorithms, reproduce the computed data used in the manuscript, and check the worked products.  It is a software deposit, not a second copy of the book source.

\section{Licence, citation and scope}
The author's original source code, package code and computed certificates are dedicated to the public domain under CC0~1.0 Universal.  This dedication does not apply to third-party material.  In particular, the vendored wheels retain their upstream licences, and files obtained from or exported verbatim by the ATLAS of Group Representations and GAP retain their original status.  The authoritative list is \texttt{NOTICE.txt} in the software archive.  Citation is requested as a matter of scientific practice, not imposed as a condition of CC0.  The final Zenodo DOI should replace \texttt{ZENODO-DOI-PENDING} in the arXiv version.

\section{Published archives}
The reproducibility deposit consists of three ZIP files and one checksum file:
\begin{center}\begin{tabular}{@{}p{0.40\linewidth}p{0.56\linewidth}@{}}
\texttt{DFGT-software-v1.0.zip}&source programs, package, certificates, frozen small inputs, and vendored macOS wheels;\\
\texttt{DFGT-data-v1.0.zip}&the 68 approved files under \texttt{scripts/data/};\\
\texttt{DFGT-monster-tables-v1.0.zip}&the six precomputed Monster pickle tables;\\
\texttt{SHA256SUMS.txt}&SHA-256 digests of all three ZIP files.
\end{tabular}\end{center}
The data ZIP has two directory entries in addition to its 68 files.  The Monster archive contains \texttt{axis\_chains.pkl}, \texttt{axis\_chains\_cq.pkl}, \texttt{axis\_orbits\_co1.pkl}, \texttt{axis\_top\_cq.pkl}, \texttt{axis\_top\_typed.pkl}, and \texttt{k\_chain.pkl}.  Backup files, build artefacts, session notes, local virtual environments, and manuscript sources are deliberately absent.

\section{Installation}
Unpack the three archives into one directory, preserving paths.  Use Python~3.13.  The pinned Monster environment is
\begin{verbatim}
mmgroup==1.0.8
numpy==2.4.6
regex==2026.9.29
\end{verbatim}
On macOS arm64 with Python~3.13, the archive provides the exact wheels used for the manuscript.  From the archive root run
\begin{verbatim}
sh setup_mm_env.sh
.venv-mm/bin/python3 -m pip install -e .
\end{verbatim}
The setup script verifies \texttt{vendor/wheels/SHA256SUMS}, creates \texttt{.venv-mm}, and installs only from the bundled wheels.  On macOS x86\_64, Linux x86\_64, Linux aarch64, or Windows, create a Python~3.13 virtual environment and install the same pinned versions from a compatible package index or from upstream source when a binary wheel is unavailable; then install this repository with \texttt{python -m pip install -e .}.  The archive does not claim that its macOS arm64 wheels are portable to another platform.

GAP~4.14 is required by the \texttt{.g} programs.  A TeX installation is required only to build or machine-check the arXiv manuscript; no TeX source is included in the Zenodo software deposit.

\section{Archive layout}
\begin{description}
\item[\texttt{digital\_laws/}] Installable digital-law package and command-line interface.
\item[\texttt{scripts/}] Programs for characters, classical and exceptional groups, crystallographic groups, Mathieu and Conway constructions, sporadic groups, worked products, and manuscript checks.
\item[\texttt{scripts/data/}] Frozen group generators, exported data and computed tables.  Their third-party status is recorded per class in \texttt{NOTICE.txt}; author-generated and derived files remain distinguishable from verbatim ATLAS or GAP exports.
\item[\texttt{monster/}] Monster law, construction and validation programs.  The large tables are distributed separately.
\item[\texttt{drafts/}] Small frozen inputs required by some dependent programs.  Other files normally written here are reproducible outputs.
\item[\texttt{certificates/}] Compact certificates checked by independent verifier programs.
\item[\texttt{vendor/wheels/}] Exact macOS arm64 wheels and their upstream checksums.
\end{description}

\section{Reproduction workflow}
First verify the downloaded files with \texttt{shasum -a 256 -c SHA256SUMS.txt}.  Unpack all archives into the same root and create the appropriate environment.  Run each command from the directory stated in the reference below: many programs intentionally use relative paths.  Programs that write results place them under \texttt{drafts/} or \texttt{certificates/}; preserve the published frozen inputs before regenerating them.

The manuscript-checking programs require the arXiv source tree to be unpacked at the same root, with \texttt{digital\_finite\_groups.tex} and \texttt{chapters/} in their published locations.  Some programs also regenerate tables that the book includes: \texttt{pc\_verify.py}, \texttt{sn\_tables.py} and \texttt{ss\_verify.py} write into \texttt{chapters/}. The TeX files belong to arXiv, not to the Zenodo software archive.  For the Chapter~13 products, run
\begin{verbatim}
cd scripts
../.venv-mm/bin/python3 worked_inverses.py
\end{verbatim}
The program checks the frozen table, verifies both-sided inverses and the anti-homomorphism identity, redirects incidental draft writes to a temporary directory, and leaves the Monster pickle files unchanged.

\section{Data provenance}
Files named \texttt{fetch\_*.g}, together with the corresponding GAP export programs, document the acquisition or export step for ATLAS/GAP material.  Such source material is not relicensed by CC0.  Text, binary, NumPy and pickle outputs made by the author's programs are derived computational artefacts; the script reference records the file names read and written whenever they occur literally in the source.  Frozen inputs for which no producer is present are retained so that the downstream calculation remains reproducible, and should be cited as part of the versioned Zenodo deposit.

\section{Reading the reference}
The following is a source-derived reference to every executable program included in the approved archive.  ``No literal file name'' means that the source constructs the path dynamically or communicates through standard input/output; it does not mean that the computation has no data.  Run times are not asserted here: they depend strongly on platform, parameters, cache state and whether a precomputed table is used.  No program was rerun merely to manufacture benchmark numbers for this chapter.

\section{Monster group}
\paragraph{\texttt{monster/build\_10A.py}.} Certified build of the top-axis orbit $10A$ (stage~4c) of the Monster law (generated by \texttt{search\_top.py}). Updates \texttt{axis\_top\_typed.pkl} and \texttt{axis\_top\_cq.pkl}; keep copies of both before running.\par
\noindent\textbf{Invocation:} \texttt{python -u build\_\allowbreak{}10A.py 2>\&1 | tee ../\allowbreak{}drafts/\allowbreak{}monster\_\allowbreak{}build10A\_\allowbreak{}2026-10-23.log}.\par
\noindent\textbf{Arguments:} none. \textbf{Dependencies:} repository modules monster\_\allowbreak{}law.\par
\noindent\textbf{Reads:} no literal file name. \textbf{Writes:} no literal file name.\par
\paragraph{\texttt{monster/build\_10B.py}.} Certified build of the top-axis orbit $10B$ (stage~4c) of the Monster law (generated by \texttt{search\_top.py}). Updates \texttt{axis\_top\_typed.pkl} and \texttt{axis\_top\_cq.pkl}; keep copies of both before running.\par
\noindent\textbf{Invocation:} \texttt{python -u build\_\allowbreak{}10B.py 2>\&1 | tee ../\allowbreak{}drafts/\allowbreak{}monster\_\allowbreak{}build10B\_\allowbreak{}2026-10-23.log}.\par
\noindent\textbf{Arguments:} none. \textbf{Dependencies:} repository modules monster\_\allowbreak{}law.\par
\noindent\textbf{Reads:} no literal file name. \textbf{Writes:} no literal file name.\par
\paragraph{\texttt{monster/build\_12C.py}.} Certified build of the top-axis orbit $12C$ (stage~4c) of the Monster law (generated by \texttt{search\_top.py}). Updates \texttt{axis\_top\_typed.pkl} and \texttt{axis\_top\_cq.pkl}; keep copies of both before running.\par
\noindent\textbf{Invocation:} \texttt{python -u build\_\allowbreak{}12C.py 2>\&1 | tee ../\allowbreak{}drafts/\allowbreak{}monster\_\allowbreak{}build12C\_\allowbreak{}2026-10-23.log}.\par
\noindent\textbf{Arguments:} none. \textbf{Dependencies:} repository modules monster\_\allowbreak{}law.\par
\noindent\textbf{Reads:} no literal file name. \textbf{Writes:} no literal file name.\par
\paragraph{\texttt{monster/build\_4A.py}.} Certified build of the top-axis orbit $4A$ (stage~4c) of the Monster law. Updates \texttt{axis\_top\_typed.pkl} and \texttt{axis\_top\_cq.pkl}; keep copies of both before running.\par
\noindent\textbf{Invocation:} \texttt{python -u build\_\allowbreak{}4A.py 2>\&1 | tee ../\allowbreak{}drafts/\allowbreak{}monster\_\allowbreak{}build4A\_\allowbreak{}v2\_\allowbreak{}2026-10-23.log}.\par
\noindent\textbf{Arguments:} none. \textbf{Dependencies:} repository modules monster\_\allowbreak{}law.\par
\noindent\textbf{Reads:} no literal file name. \textbf{Writes:} no literal file name.\par
\paragraph{\texttt{monster/build\_4B.py}.} Certified build of the top-axis orbit $4B$ (stage~4c) of the Monster law. Updates \texttt{axis\_top\_typed.pkl} and \texttt{axis\_top\_cq.pkl}; keep copies of both before running.\par
\noindent\textbf{Invocation:} \texttt{python -u build\_\allowbreak{}4B.py 2>\&1 | tee ../\allowbreak{}drafts/\allowbreak{}monster\_\allowbreak{}build4B\_\allowbreak{}2026-10-23.log}.\par
\noindent\textbf{Arguments:} none. \textbf{Dependencies:} repository modules monster\_\allowbreak{}law.\par
\noindent\textbf{Reads:} no literal file name. \textbf{Writes:} no literal file name.\par
\paragraph{\texttt{monster/build\_4C.py}.} Certified build of the top-axis orbit $4C$ (stage~4c) of the Monster law (generated on 23~October by the recipe search). Updates \texttt{axis\_top\_typed.pkl} and \texttt{axis\_top\_cq.pkl}; keep copies of both before running.\par
\noindent\textbf{Invocation:} \texttt{python -u build\_\allowbreak{}4C.py 2>\&1 | tee ../\allowbreak{}drafts/\allowbreak{}monster\_\allowbreak{}build4C\_\allowbreak{}2026-10-23.log}.\par
\noindent\textbf{Arguments:} none. \textbf{Dependencies:} repository modules monster\_\allowbreak{}law.\par
\noindent\textbf{Reads:} no literal file name. \textbf{Writes:} no literal file name.\par
\paragraph{\texttt{monster/build\_6A.py}.} Certified build of the top-axis orbit $6A$ (stage~4c) of the Monster law (generated by \texttt{search\_top.py}). Updates \texttt{axis\_top\_typed.pkl} and \texttt{axis\_top\_cq.pkl}; keep copies of both before running.\par
\noindent\textbf{Invocation:} \texttt{python -u build\_\allowbreak{}6A.py 2>\&1 | tee ../\allowbreak{}drafts/\allowbreak{}monster\_\allowbreak{}build6A\_\allowbreak{}2026-10-23.log}.\par
\noindent\textbf{Arguments:} none. \textbf{Dependencies:} repository modules monster\_\allowbreak{}law.\par
\noindent\textbf{Reads:} no literal file name. \textbf{Writes:} no literal file name.\par
\paragraph{\texttt{monster/build\_6C.py}.} Certified build of the top-axis orbit $6C$ (stage~4c) of the Monster law (generated by \texttt{search\_top.py}). Updates \texttt{axis\_top\_typed.pkl} and \texttt{axis\_top\_cq.pkl}; keep copies of both before running.\par
\noindent\textbf{Invocation:} \texttt{python -u build\_\allowbreak{}6C.py 2>\&1 | tee ../\allowbreak{}drafts/\allowbreak{}monster\_\allowbreak{}build6C\_\allowbreak{}2026-10-23.log}.\par
\noindent\textbf{Arguments:} none. \textbf{Dependencies:} repository modules monster\_\allowbreak{}law.\par
\noindent\textbf{Reads:} no literal file name. \textbf{Writes:} no literal file name.\par
\paragraph{\texttt{monster/build\_6F.py}.} Certified build of the top-axis orbit $6F$ (stage~4c) of the Monster law (generated by \texttt{search\_top.py}). Updates \texttt{axis\_top\_typed.pkl} and \texttt{axis\_top\_cq.pkl}; keep copies of both before running.\par
\noindent\textbf{Invocation:} \texttt{python -u build\_\allowbreak{}6F.py 2>\&1 | tee ../\allowbreak{}drafts/\allowbreak{}monster\_\allowbreak{}build6F\_\allowbreak{}2026-10-23.log}.\par
\noindent\textbf{Arguments:} none. \textbf{Dependencies:} repository modules monster\_\allowbreak{}law.\par
\noindent\textbf{Reads:} no literal file name. \textbf{Writes:} no literal file name.\par
\paragraph{\texttt{monster/build\_8B.py}.} Certified build of the top-axis orbit $8B$ (stage~4c) of the Monster law (generated by \texttt{search\_top.py}). Updates \texttt{axis\_top\_typed.pkl} and \texttt{axis\_top\_cq.pkl}; keep copies of both before running.\par
\noindent\textbf{Invocation:} \texttt{python -u build\_\allowbreak{}8B.py 2>\&1 | tee ../\allowbreak{}drafts/\allowbreak{}monster\_\allowbreak{}build8B\_\allowbreak{}2026-10-23.log}.\par
\noindent\textbf{Arguments:} none. \textbf{Dependencies:} repository modules monster\_\allowbreak{}law.\par
\noindent\textbf{Reads:} no literal file name. \textbf{Writes:} no literal file name.\par
\paragraph{\texttt{monster/check\_mlaw\_validation.py}.} Re-runs the validation of the law of $M$ in Section~\ref{ssec:monstervalid}, checks (i)--(iv), and compares direct products with \texttt{mmgroup}.\par
\noindent\textbf{Invocation:} \texttt{cd monster \&\& ../\allowbreak{}.venv-mm/\allowbreak{}bin/\allowbreak{}python3 check\_\allowbreak{}mlaw\_\allowbreak{}validation.py}.\par
\noindent\textbf{Arguments:} positional argument indices 1. \textbf{Dependencies:} repository modules monster\_\allowbreak{}law; external packages mmgroup.\par
\noindent\textbf{Reads:} no literal file name. \textbf{Writes:} no literal file name.\par
\paragraph{\texttt{monster/monster\_law.py}.} Digital law of the Monster $M$ through Seysen's \texttt{mmgroup} (version 1.0.8): $m*n=D^{-1}\bigl(D(m)\,D(n)\bigr)$ on the ordinals $0\le m<|M|$. The law is complete: every ordinal is encoded and every element of $M$ is decoded. Ordinals are laid out as $m=s+2^{24}c+|K|\,r_2+|2.B|\,r_1$, where $K=C_{G_{x0}}(x_\beta)=2^{1+23}.\mathrm{Co}_2$, $\beta=\mathtt{0x200}$ (so that $x_\beta$ is the element of the axis $v^+$) and $|2.B|=|K|\,N_2$. Stage~1 ($0\le m<2^{24}$) encodes the normal $2$-subgroup $Q_K=Q_{x0}\cap K$.\par
\noindent\textbf{Invocation:} \texttt{cd monster \&\& ../\allowbreak{}.venv-mm/\allowbreak{}bin/\allowbreak{}python3 monster\_\allowbreak{}law.py}.\par
\noindent\textbf{Arguments:} positional argument indices 1:, 1:2, 1:3, 2. \textbf{Dependencies:} external packages mmgroup, numpy.\par
\noindent\textbf{Reads:} no literal file name. \textbf{Writes:} no literal file name.\par
\paragraph{\texttt{monster/probe\_12C.py}.} Exploratory variant of \texttt{search\_top.py}, used for the orbit $12C$.\par
\noindent\textbf{Invocation:} \texttt{cd monster \&\& ../\allowbreak{}.venv-mm/\allowbreak{}bin/\allowbreak{}python3 probe\_\allowbreak{}12C.py}.\par
\noindent\textbf{Arguments:} positional argument indices 1, 2, 3. \textbf{Dependencies:} repository modules monster\_\allowbreak{}law; external packages mmgroup, numpy.\par
\noindent\textbf{Reads:} no literal file name. \textbf{Writes:} no literal file name.\par
\paragraph{\texttt{monster/probe\_12C\_v2.py}.} Second exploratory variant of \texttt{search\_top.py} for the orbit $12C$.\par
\noindent\textbf{Invocation:} \texttt{cd monster \&\& ../\allowbreak{}.venv-mm/\allowbreak{}bin/\allowbreak{}python3 probe\_\allowbreak{}12C\_\allowbreak{}v2.py}.\par
\noindent\textbf{Arguments:} positional argument indices 1, 2, 3. \textbf{Dependencies:} repository modules monster\_\allowbreak{}law; external packages mmgroup, numpy.\par
\noindent\textbf{Reads:} no literal file name. \textbf{Writes:} no literal file name.\par
\paragraph{\texttt{monster/probe\_15120\_cert.py}.} Exploratory variant of \texttt{search\_top.py} with an additional certification routine \texttt{cert}.\par
\noindent\textbf{Invocation:} \texttt{cd monster \&\& ../\allowbreak{}.venv-mm/\allowbreak{}bin/\allowbreak{}python3 probe\_\allowbreak{}15120\_\allowbreak{}cert.py}.\par
\noindent\textbf{Arguments:} positional argument indices 1, 2, 3. \textbf{Dependencies:} repository modules monster\_\allowbreak{}law; external packages mmgroup, numpy.\par
\noindent\textbf{Reads:} no literal file name. \textbf{Writes:} no literal file name.\par
\paragraph{\texttt{monster/probe\_4A\_vchain.py}.} Probes $\mathrm{Co}_2$-layer recipes for the types $2048$, $16192$ and $32384$ of the orbit $4A$; writes no data files. A recipe is good when $P_1=|\mathrm{Co}_2|/|\bar H|$.\par
\noindent\textbf{Invocation:} \texttt{python -u probe\_\allowbreak{}4A\_\allowbreak{}vchain.py 2>\&1 | tee ../\allowbreak{}drafts/\allowbreak{}monster\_\allowbreak{}4A\_\allowbreak{}vchain\_\allowbreak{}probe\_\allowbreak{}2026-10-23.log}.\par
\noindent\textbf{Arguments:} none. \textbf{Dependencies:} repository modules monster\_\allowbreak{}law; external packages mmgroup, numpy.\par
\noindent\textbf{Reads:} no literal file name. \textbf{Writes:} no literal file name.\par
\paragraph{\texttt{monster/search\_top.py}.} Recipe search for one top-axis orbit of the Monster law; writes \texttt{build\_<orbit>.py} and does not modify any pickle file. Arguments: the orbit, the seconds allowed per recipe test (default $6$) and per type (default $40$). A recipe certifies if and only if $P_1=\prod(\text{sizes of its }\mathrm{Co}_2\text{ levels})=|\mathrm{Co}_2|/|\bar H|$.\par
\noindent\textbf{Invocation:} \texttt{python -u search\_\allowbreak{}top.py 6A 6 40 2>\&1 | tee ../\allowbreak{}drafts/\allowbreak{}monster\_\allowbreak{}search\_\allowbreak{}6A\_\allowbreak{}2026-10-23.log}.\par
\noindent\textbf{Arguments:} positional argument indices 1, 2, 3. \textbf{Dependencies:} repository modules monster\_\allowbreak{}law; external packages mmgroup, numpy.\par
\noindent\textbf{Reads:} no literal file name. \textbf{Writes:} no literal file name.\par
\paragraph{\texttt{monster/test\_ex\_blaw.py}.} Regression test for Example~\ref{ex:blaw} (computing in $B$) and for Monster products, against \texttt{mmgroup}. Re-executes itself under \texttt{../.venv-mm/bin/python3} if \texttt{mmgroup} is not installed.\par
\noindent\textbf{Invocation:} \texttt{cd monster \&\& python3 test\_\allowbreak{}ex\_\allowbreak{}blaw.py (re-executes under ../\allowbreak{}.venv-mm/\allowbreak{}bin/\allowbreak{}python3 if mmgroup is missing)}.\par
\noindent\textbf{Arguments:} none. \textbf{Dependencies:} repository modules monster\_\allowbreak{}law; external packages mmgroup.\par
\noindent\textbf{Reads:} no literal file name. \textbf{Writes:} no literal file name.\par
\paragraph{\texttt{monster/validate\_top.py}.} Validates one implemented top-axis orbit of the Monster law, read-only. Arguments: the orbit and a second numeric argument ($500$ in the documented example \texttt{validate\_top.py 4A 500}).\par
\noindent\textbf{Invocation:} \texttt{python -u validate\_\allowbreak{}top.py 4A 500 2>\&1 | tee ../\allowbreak{}drafts/\allowbreak{}monster\_\allowbreak{}validate\_\allowbreak{}4A\_\allowbreak{}2026-10-23.log}.\par
\noindent\textbf{Arguments:} positional argument indices 1, 2. \textbf{Dependencies:} repository modules monster\_\allowbreak{}law; external packages mmgroup.\par
\noindent\textbf{Reads:} no literal file name. \textbf{Writes:} no literal file name.\par

\section{Supporting computations}
\paragraph{\texttt{scripts/\_tmpwrites.py}.} Redirect writes into ../\allowbreak{}drafts/\allowbreak{} to the system temp directory. mathieu\_\allowbreak{}law.py, golay\_\allowbreak{}law.py and conway\_\allowbreak{}examples.py write ../\allowbreak{}drafts/\allowbreak{}*.json|tex when they are imported. The worked\_\allowbreak{}*.py drivers import this module first, so that they write nothing in the project.\par
\noindent\textbf{Invocation:} \texttt{cd scripts \&\& python3 \_\allowbreak{}tmpwrites.py}.\par
\noindent\textbf{Arguments:} none. \textbf{Dependencies:} Python standard library.\par
\noindent\textbf{Reads:} no literal file name. \textbf{Writes:} no literal file name.\par
\paragraph{\texttt{scripts/cc\_verify.py}.} Verification for the chapter 'Conjugacy classes and how to compute them'.\par
\noindent\textbf{Invocation:} \texttt{python3 scripts/\allowbreak{}cc\_\allowbreak{}verify.py [GAP\_\allowbreak{}BINARY] (pure Python; the GAP cross-check is skipped if GAP is absent)}.\par
\noindent\textbf{Arguments:} positional argument indices 1. \textbf{Dependencies:} Python standard library.\par
\noindent\textbf{Reads:} no literal file name. \textbf{Writes:} no literal file name.\par
\paragraph{\texttt{scripts/coset\_law.py}.} Full digital law of a matrix group G from a coset orbit (Schreier tree + sorted key index) and a law of K = Stab\_\allowbreak{}G(S). Matrices act on row vectors (right action). The start object S (rows over GF(2)) has stabiliser K; the point t of the orbit is S s\_\allowbreak{}t, with s\_\allowbreak{}t the product of the generators along th\par
\noindent\textbf{Invocation:} \texttt{cd scripts \&\& python3 coset\_\allowbreak{}law.py}.\par
\noindent\textbf{Arguments:} none. \textbf{Dependencies:} external packages numpy.\par
\noindent\textbf{Reads:} no literal file name. \textbf{Writes:} no literal file name.\par
\paragraph{\texttt{scripts/ct\_probe.g}.} GAP probe using the Character Table Library. It enumerates all groups of order at most $32$, counting nonabelian groups, their conjugacy classes, the largest class count, and groups with an irrational irreducible character. It also searches the nonabelian simple groups of order at most $10^6$ with available character tables.\par
\noindent\textbf{Invocation:} \texttt{cd scripts \&\& gap -q ct\_\allowbreak{}probe.g}.\par
\noindent\textbf{Arguments:} none. \textbf{Dependencies:} GAP 4.14 or compatible.\par
\noindent\textbf{Reads:} no literal file name. \textbf{Writes:} no literal file name.\par
\paragraph{\texttt{scripts/mm\_axes\_audit.py}.} Audits the effective $2A$-axis reduction supplied by \texttt{mmgroup.axes} (tested with \texttt{mmgroup} 1.0.8). For an axis $a=v^+g$, \texttt{reduce\_G\_x0()} returns $q\in G_{x0}$ with $aq=a_c$, the fixed representative of its named $G_{x0}$-orbit. For each such representative one transporter $s_c$ with $v^+s_c=a_c$ is constructed, so that $k=g\,q\,s_c^{-1}$ fixes $v^+$ and should centralize the corresponding $2A$-involution. This is the computable reduction step towards $M/2.B$. The class \texttt{BabyAxis} is audited separately; its ten $H$-orbits are not asserted to give the smallest action of $B$.\par
\noindent\textbf{Invocation:} \texttt{cd DigitalFiniteGroupTheory \&\& .venv-mm/\allowbreak{}bin/\allowbreak{}python -u scripts/\allowbreak{}mm\_\allowbreak{}axes\_\allowbreak{}audit.py --monster 12 --baby 12}.\par
\noindent\textbf{Arguments:} --monster default=12; --baby default=12; --baby-k default=0; --seed default=20261013. \textbf{Dependencies:} external packages mmgroup.\par
\noindent\textbf{Reads:} no literal file name. \textbf{Writes:} no literal file name.\par
\paragraph{\texttt{scripts/multiply.py}.} Multiplies two labels in any group of the archive: every Lie-type group with a digital law, ${}^2G_2(3)$, $M_{24}$, $\Co_0$, $J_4$, $\mathrm{Th}$, the subgroups $K$ of $J_4$ and $\mathrm{Th}$, the Monster, and the groups of Chapter~\ref{ch:sporadic} provided by \texttt{sporadic\_laws.py}, and the additional groups of Chapter~\ref{ch:lie} provided by \texttt{lie\_extra.py}. It calls the same laws as \texttt{worked\_inverses.py}, rejects labels outside $[0,|G|)$, and with \texttt{--check} recomputes the product independently where a decoder exists (Remark~\ref{rem:wpchecks}).\par \noindent\textbf{Invocation:} \texttt{cd scripts \&\& python3 multiply.py GROUP m n [--check]}; \texttt{python3 multiply.py --list} prints the group names. \textbf{Dependencies:} NumPy; \texttt{mmgroup}~1.0.8 for the Monster and the Baby Monster. \textbf{Writes:} nothing in the project. It also multiplies in $\Sym_n$ ($n\ge1$) and $\Alt_n$ ($n\ge2$) with the Lehmer laws of Section~\ref{sec:sn:laws} (names \texttt{Sym(n)}, \texttt{Alt(n)}, aliases such as \texttt{Sym20}, \texttt{S20}, \texttt{Alt12}); there \texttt{--check} compares each product with a separately coded ranking of permutations.\par

\paragraph{\texttt{scripts/inverse.py}.} Computes the inverse of one or more labels in any group of the archive, with the names, laws and label range of \texttt{multiply.py}: by successive powers in the digital law when the order is at most $256$, and otherwise as $m^{|G|-1}$ by square-and-multiply. With \texttt{--order} it also prints the order of each label; with \texttt{--check} it asserts $mm^{-1}=m^{-1}m=e$, $em=me=m$ and $(m^{-1})^{-1}=m$ and, where a decoder exists, recomputes $mm^{-1}$ independently (Remark~\ref{rem:wpinverses}). Called without arguments, it prints the numbered list of groups and asks for a group and then for labels.\par \noindent\textbf{Invocation:} \texttt{cd scripts \&\& python3 inverse.py GROUP m [m2 ...] [--order] [--check]}; \texttt{python3 inverse.py} for the interactive mode. \textbf{Dependencies:} NumPy; \texttt{mmgroup}~1.0.8 for the Monster and the Baby Monster. \textbf{Writes:} nothing in the project.\par

\paragraph{\texttt{scripts/conjugate.py}.} Computes the conjugate $m*n*m^{-1}$ of a label $n$ by a label $m$, the image of $n$ under the inner automorphism $x\mapsto mxm^{-1}$, in any group of the archive, with the names, laws and label range of \texttt{multiply.py}. With \texttt{--right} it computes $m^{-1}*n*m$ instead. The inverse is the one of \texttt{inverse.py}, so only the digital law is used.\par
\noindent\textbf{Invocation:} \texttt{cd scripts \&\& python3 conjugate.py GROUP m n [--right] [--check]}.\par
\noindent\textbf{Arguments:} \texttt{GROUP}, \texttt{m}, \texttt{n} as in \texttt{multiply.py}; \texttt{--right}; \texttt{--check} asserts $m*m^{-1}=m^{-1}*m=e$ and $(m*n*m^{-1})*m=m*n$, and that $n$ and its conjugate have the same order. \textbf{Dependencies:} those of \texttt{multiply.py}.\par
\noindent\textbf{Reads:} no literal file name. \textbf{Writes:} nothing in the project.\par

\paragraph{\texttt{scripts/sn\_law.py}.} The Lehmer laws of $\Sym_n$ and $\Alt_n$ (Section~\ref{sec:sn:laws}) used by \texttt{multiply.py}: the Lehmer labels $\lambda$ and $\mu$, their decoding, and the transported products on $\{0,\dots,n!-1\}$ and $\{0,\dots,n!/2-1\}$.\par
\noindent\textbf{Invocation:} imported by \texttt{multiply.py} and \texttt{sn\_\allowbreak{}cayley.py}. \textbf{Dependencies:} Python standard library only. \textbf{Writes:} nothing in the project.\par

\paragraph{\texttt{scripts/sn\_cayley.py}.} Writes the Cayley tables of the Lehmer laws of $\Alt_4$, $\Sym_4$ and $\Alt_5$ (Tables~\ref{tab:an:law4}, \ref{tab:sn:law4}, \ref{tab:an:law5a} and~\ref{tab:an:law5b}). Before writing, it checks that label $0$ is the identity, that each table is a Latin square, that the law is associative on all triples, and that every entry agrees with \texttt{multiply.py} and with direct composition of the decoded permutations.\par
\noindent\textbf{Invocation:} \texttt{cd scripts \&\& python3 sn\_\allowbreak{}cayley.py}, or \texttt{python3 sn\_\allowbreak{}cayley.py --print T}.\par
\noindent\textbf{Arguments:} optional \texttt{--print T} with \texttt{T} one of \texttt{A4}, \texttt{S4}, \texttt{A5a}, \texttt{A5b}, which prints that table to standard output instead of writing files. \textbf{Dependencies:} Python standard library only.\par
\noindent\textbf{Reads:} no literal file name. \textbf{Writes:} \texttt{../\allowbreak{}chapters/\allowbreak{}sn\_\allowbreak{}cayley\_\allowbreak{}T.tex} for the four values of \texttt{T}, which requires the manuscript source tree; nothing with \texttt{--print}.\par

\paragraph{Group names and the Monster.} \texttt{multiply.py}, \texttt{inverse.py} and \texttt{conjugate.py} accept the same group names, and \texttt{python3 inverse.py --list} prints them. In \texttt{zsh} and other shells, quote a name containing parentheses or a prime, for example \texttt{\textquotesingle Sym(20)\textquotesingle} or \texttt{\textquotesingle PSL2(7)\textquotesingle}, or use a punctuation-free alias such as \texttt{Sym20}, \texttt{Alt12}, \texttt{psl27} or \texttt{Tits}. The Monster (name \texttt{Monster}, alias \texttt{M}) needs \texttt{mmgroup}~1.0.8 and is therefore run with the interpreter of the Monster environment of the installation section, for example \texttt{cd scripts \&\& ../\allowbreak{}.venv-mm/\allowbreak{}bin/\allowbreak{}python3 multiply.py Monster m n}; the same applies to \texttt{inverse.py} and \texttt{conjugate.py}. There, \texttt{--check} recomputes the products with \texttt{mmgroup} words; it also needs the Monster tables of \texttt{DFGT-\allowbreak{}monster-\allowbreak{}tables-\allowbreak{}v1.0.zip} (for example \texttt{monster/\allowbreak{}k\_\allowbreak{}chain.pkl}) unpacked in \texttt{monster/}.\par

\paragraph{\texttt{scripts/wallpaper\_law.py}.} The digital laws~\eqref{eq:crystlaw} of the $17$ plane crystallographic groups modulo $mT$ (Section~\ref{sec:wpcrystal}), used by \texttt{multiply.py}, \texttt{inverse.py} and \texttt{conjugate.py} under the names \texttt{GROUP/m}, for example \texttt{pg/6} or \texttt{p4g/4}.\par
\noindent\textbf{Invocation:} \texttt{cd scripts \&\& python3 wallpaper\_\allowbreak{}law.py --check --table}.\par
\noindent\textbf{Arguments:} \texttt{--check} verifies the law against composition of affine maps, associativity for $m=2$ on all triples, and that the $17$ quotients by $6T$ are pairwise non-isomorphic; \texttt{--table} prints the invariants behind Table~\ref{tab:wallpaperm}. \textbf{Dependencies:} external packages numpy.\par
\noindent\textbf{Reads:} no literal file name. \textbf{Writes:} nothing in the project.\par

\paragraph{\texttt{scripts/space\_quotients.py}.} Invariants of the $230$ space groups modulo $m$ times the conventional lattice: the joint distribution of element orders and centralizer orders, computed in closed form from the symmetry operations of the International Tables, without a multiplication table (Section~\ref{sec:wpcrystal}, Table~\ref{tab:spacem}).\par
\noindent\textbf{Invocation:} \texttt{cd scripts \&\& python3 space\_\allowbreak{}quotients.py --m 1,2,3,4,6,12}.\par
\noindent\textbf{Arguments:} \texttt{--m} comma-separated moduli (default \texttt{1,2,3,4,6,12}); \texttt{--csv FILE} writes the counts of Table~\ref{tab:spacem}; \texttt{--check} compares the closed forms with brute force on ten small quotients. \textbf{Dependencies:} external package numpy (operations from \texttt{data/space\_group\_ops.txt}).\par
\noindent\textbf{Reads:} no literal file name. \textbf{Writes:} \texttt{FILE} with \texttt{--csv} only.\par

\paragraph{\texttt{scripts/space\_residual.py}.} Full multiplication tables of the quotients $G/mL$ and finer isomorphism invariants (derived and lower central series, second centre, squares, fibres of the squaring and commutator maps), applied to the types not separated by \texttt{space\_quotients.py}; \texttt{--iso} is an exact isomorphism test by backtracking over images of generators (Section~\ref{sec:wpcrystal}).\par
\noindent\textbf{Invocation:} \texttt{cd scripts \&\& python3 space\_\allowbreak{}residual.py --check} and \texttt{python3 space\_\allowbreak{}residual.py --iso 32,34 --m 2,4,8}.\par
\noindent\textbf{Arguments:} optional sets of space-group numbers (default: the five residual sets), \texttt{--m} moduli (default \texttt{2,4,6,8}), \texttt{--max} largest table (default $4096$), \texttt{--iso N1,N2}, \texttt{--check}. \textbf{Dependencies:} external package numpy (operations from \texttt{data/space\_group\_ops.txt}).\par
\noindent\textbf{Reads:} no literal file name. \textbf{Writes:} nothing in the project.\par

\paragraph{\texttt{scripts/space\_gap.py}.} Exports quotients $G/mL$ to GAP as permutation groups (right regular action on a small generating set) and prints \texttt{IdGroup} and the result of \texttt{IsomorphismGroups}; an independent check of \texttt{space\_residual.py --iso} (Section~\ref{sec:wpcrystal}).\par
\noindent\textbf{Invocation:} \texttt{cd scripts \&\& python3 space\_\allowbreak{}gap.py 30,32,34 --m 2,4 > q.g \&\& gap -q -o 4g q.g}.\par
\noindent\textbf{Arguments:} comma-separated space-group numbers, \texttt{--m} moduli (default \texttt{2,4}). \textbf{Dependencies:} external package numpy (operations from \texttt{data/space\_group\_ops.txt}); GAP with the SmallGroups library for the generated file.\par
\noindent\textbf{Reads:} no literal file name. \textbf{Writes:} the GAP input on standard output.\par

\paragraph{\texttt{scripts/space\_ops.py}.} The symmetry operations of the $230$ space groups, read from the data file \texttt{data/space\_group\_ops.txt} (International Tables triplets, conventional cell, centring included), with a triplet parser and formatter; every space-group script obtains its operations here (Section~\ref{sec:crystlaw3}).\par
\noindent\textbf{Invocation:} \texttt{cd scripts \&\& python3 space\_ops.py 14}; \texttt{python3 space\_ops.py --export} regenerates the data file.\par
\noindent\textbf{Arguments:} a space-group number, or \texttt{--export}. \textbf{Dependencies:} external package numpy; gemmi only for \texttt{--export}.\par
\noindent\textbf{Reads:} \texttt{data/space\_group\_ops.txt}. \textbf{Writes:} \texttt{data/space\_group\_ops.txt} with \texttt{--export}, otherwise nothing.\par

\paragraph{\texttt{scripts/space\_law.py}.} The digital laws~\eqref{eq:spacelaw} of the $230$ space groups modulo $mL$ (Section~\ref{sec:crystlaw3}), used by \texttt{multiply.py}, \texttt{inverse.py} and \texttt{conjugate.py} under the names \texttt{SG$N$/$m$}, for example \texttt{SG62/2}; the independent check multiplies $4\times4$ affine matrices.\par
\noindent\textbf{Invocation:} \texttt{cd scripts \&\& python3 multiply.py SG62/2 5 17 --check}.\par
\noindent\textbf{Arguments:} none of its own. \textbf{Dependencies:} external package numpy (operations from \texttt{data/space\_group\_ops.txt}).\par
\noindent\textbf{Reads:} no literal file name. \textbf{Writes:} nothing in the project.\par

\paragraph{\texttt{scripts/space\_table.py}.} GAP input identifying the quotients of the $17$ plane groups modulo $2T$ and $6T$ and of the $230$ space groups modulo $2L$ (Tables~\ref{tab:planeid} and~\ref{tab:spaceid}); the GAP output is kept in \texttt{scripts/space\_table\_gap.txt}.\par
\noindent\textbf{Invocation:} \texttt{cd scripts \&\& python3 space\_\allowbreak{}table.py > t.g \&\& gap -q -o 16g t.g}.\par
\noindent\textbf{Arguments:} none. \textbf{Dependencies:} external package numpy (operations from \texttt{data/space\_group\_ops.txt}); GAP with the SmallGroups library for the generated file.\par
\noindent\textbf{Reads:} no literal file name. \textbf{Writes:} the GAP input on standard output.\par

\paragraph{\texttt{scripts/sporadic\_laws.py}.} The laws of Chapter~\ref{ch:sporadic} used by \texttt{multiply.py}: the canonical labels of Section~\ref{sec:canonical} for $J_1$, $J_2$, $\mathrm{HS}$, $\mathrm{Suz}$, $\mathrm{He}$, $\mathrm{Ru}$, $J_3$, $\mathrm{O'N}$, $\mathrm{HN}$, $\mathrm{Fi}_{24}$ and $\mathrm{Ly}$, built with \texttt{digital\_laws.fast} from the generator files in \texttt{scripts/data}; the chain labels of the group sections, obtained by running the chain scripts unchanged (\texttt{primsporadic\_chain.py}, \texttt{atlas\_sporadic\_chain.py}, \texttt{fi22\_chain.py}, \texttt{fi23\_chain.py}, \texttt{fi24\_chain.py}, \texttt{hn\_chain.py}); the Mathieu chains of \texttt{mathieu\_law.py}; and the Baby Monster of \texttt{monster/monster\_law.py}. Write-mode opens inside the project are redirected to the temporary directory.\par
\noindent\textbf{Invocation:} imported by \texttt{multiply.py}. \textbf{Dependencies:} NumPy; \texttt{mmgroup}~1.0.8 for the Baby Monster. \textbf{Writes:} nothing in the project.\par

\paragraph{\texttt{scripts/lie\_extra.py}.} The laws of Chapter~\ref{ch:lie} used by \texttt{multiply.py} beyond the fixed list of \texttt{worked\_products\_lie\_all.py}: the Tits group ${}^2F_4(2)'$ and ${}^2F_4(2)$ on $\Ap{2p_0}$ (Theorem~\ref{thm:titsdigital}), taken from the $q=2$ part of \texttt{f42\_digital.py}; the dense law of $PSL_n(q)$ (Theorem~\ref{thm:pslndense}), taken from the encoding and decoding of \texttt{psln\_dense.py} without the enumeration of $SL_n(q)$; and $PSL_2(p)$, $Sz(q)$, ${}^2G_2(q)$ for every admissible $q$, through the builders of \texttt{worked\_products\_lie\_all.py}, \texttt{sz\_law.py} and \texttt{ree\_law.py}.\par
\noindent\textbf{Invocation:} imported by \texttt{multiply.py}. \textbf{Dependencies:} NumPy. \textbf{Writes:} nothing in the project.\par

\paragraph{\texttt{scripts/polya\_verify.py}.} Verification of the Polya--Redfield and de Bruijn counts quoted in the book. Every number is computed twice: by the cycle-index /\allowbreak{} de Bruijn formula and by brute-force orbit enumeration (union-find) where that is feasible.\par
\noindent\textbf{Invocation:} \texttt{cd scripts \&\& python3 polya\_\allowbreak{}verify.py}.\par
\noindent\textbf{Arguments:} none. \textbf{Dependencies:} Python standard library.\par
\noindent\textbf{Reads:} no literal file name. \textbf{Writes:} no literal file name.\par

\section{Other sporadic groups}
\paragraph{\texttt{scripts/atlas\_sporadic\_chain.py}.} Stabiliser-chain digital laws of $J_3$, on $6156$ points, and $\mathrm{O'N}$, on $122760$ points, from AtlasRep generators.\par
\noindent\textbf{Invocation:} \texttt{Chain laws of J3 (6156 points) and O'N (122760 points) from AtlasRep generators. Usage: cd scripts \&\& python3 atlas\_\allowbreak{}sporadic\_\allowbreak{}chain.py}.\par
\noindent\textbf{Arguments:} none. \textbf{Dependencies:} external packages numpy.\par
\noindent\textbf{Reads:} no literal file name. \textbf{Writes:} ../\allowbreak{}drafts/\allowbreak{}j3on.json, ../\allowbreak{}drafts/\allowbreak{}j3on.tex.\par
\paragraph{\texttt{scripts/fi22\_chain.py}.} Stabiliser-chain digital law of $\mathrm{Fi}_{22}$ on its $3510$ $3$-transpositions.\par
\noindent\textbf{Invocation:} \texttt{Chain law of Fi22 on its 3510 3-transpositions. Usage: cd scripts \&\& python3 fi22\_\allowbreak{}chain.py}.\par
\noindent\textbf{Arguments:} none. \textbf{Dependencies:} external packages numpy.\par
\noindent\textbf{Reads:} ../\allowbreak{}drafts/\allowbreak{}fi22\_\allowbreak{}gap.log. \textbf{Writes:} ../\allowbreak{}drafts/\allowbreak{}fi22\_\allowbreak{}chain.json, ../\allowbreak{}drafts/\allowbreak{}fi22\_\allowbreak{}chain.tex.\par
\paragraph{\texttt{scripts/fi23\_chain.py}.} Stabiliser-chain digital law of $\mathrm{Fi}_{23}$ on its $31671$ $3$-transpositions.\par
\noindent\textbf{Invocation:} \texttt{Chain law of Fi23 on its 31671 3-transpositions. Usage: cd scripts \&\& python3 fi23\_\allowbreak{}chain.py}.\par
\noindent\textbf{Arguments:} none. \textbf{Dependencies:} external packages numpy.\par
\noindent\textbf{Reads:} no literal file name. \textbf{Writes:} ../\allowbreak{}drafts/\allowbreak{}fi23\_\allowbreak{}chain.json, ../\allowbreak{}drafts/\allowbreak{}fi23\_\allowbreak{}chain.tex.\par
\paragraph{\texttt{scripts/fi24\_chain.py}.} Stabiliser-chain digital laws of $\mathrm{Fi}_{24}=\mathrm{Fi}_{24}'{:}2$ and $\mathrm{Fi}_{24}'$ on their $306936$ $3$-transpositions.\par
\noindent\textbf{Invocation:} \texttt{Chain laws of Fi24 = Fi24'.2 and Fi24' on the 306936 3-transpositions. Usage: cd scripts \&\& python3 fi24\_\allowbreak{}chain.py}.\par
\noindent\textbf{Arguments:} none. \textbf{Dependencies:} external packages numpy.\par
\noindent\textbf{Reads:} no literal file name. \textbf{Writes:} ../\allowbreak{}drafts/\allowbreak{}fi24\_\allowbreak{}chain.json, ../\allowbreak{}drafts/\allowbreak{}fi24\_\allowbreak{}chain.tex.\par
\paragraph{\texttt{scripts/hn\_chain.py}.} Stabiliser-chain digital law of the Harada--Norton group $\mathrm{HN}$ on $1140000$ points, from AtlasRep generators.\par
\noindent\textbf{Invocation:} \texttt{Chain law of the Harada-Norton group HN on 1140000 points (AtlasRep generators). Usage: cd scripts \&\& python3 hn\_\allowbreak{}chain.py}.\par
\noindent\textbf{Arguments:} none. \textbf{Dependencies:} external packages numpy.\par
\noindent\textbf{Reads:} data/\allowbreak{}HN\_\allowbreak{}1140000.g. \textbf{Writes:} ../\allowbreak{}drafts/\allowbreak{}hn\_\allowbreak{}chain.json, ../\allowbreak{}drafts/\allowbreak{}hn\_\allowbreak{}chain.tex.\par
\paragraph{\texttt{scripts/j4\_k\_chain.g}.} Deterministic stabiliser chain for $K=2^{11}{:}M_{24}$ in the $112$-dimensional $\mathbb{F}_2$-module $V$ of $J_4$. Here $U$ is the $K$-fixed line; the quotient $K/U$ acts on $V/U$ with image $M_{24}$.\par
\noindent\textbf{Invocation:} \texttt{Deterministic stabilizer chain for K = 2\textasciicircum{}11:M24 in the J4 112-dimensional GF(2) module. U is the K-fixed line. K/\allowbreak{}U is represented on V/\allowbreak{}U, where the induced action has image M24. Run: \textasciitilde{}/\allowbreak{}software/\allowbreak{}gap-4.14.0/\allowbreak{}gap -q -o 8g scripts/\allowbreak{}j4\_\allowbreak{}k\_\allowbreak{}chain.g}.\par
\noindent\textbf{Arguments:} none. \textbf{Dependencies:} GAP 4.14 or compatible.\par
\noindent\textbf{Reads:} no literal file name. \textbf{Writes:} no literal file name.\par
\paragraph{\texttt{scripts/j4\_k\_export.g}.} Offline export of the law data of $K=2^{11}{:}M_{24}<J_4$, with its $M_{24}$ identified with the book's $M_{24}$ (Proposition~\ref{prop:mathieuchain}).\par
\noindent\textbf{Invocation:} \texttt{cd scripts \&\& gap -q j4\_\allowbreak{}k\_\allowbreak{}export.g}.\par
\noindent\textbf{Arguments:} none. \textbf{Dependencies:} GAP 4.14 or compatible.\par
\noindent\textbf{Reads:} /\allowbreak{}tmp/\allowbreak{}m24book.g. \textbf{Writes:} no literal file name.\par
\paragraph{\texttt{scripts/j4\_k\_law.py}.} Digital law of $K=2^{11}{:}M_{24}<J_4$, using the semidirect-product construction (Theorem~\ref{thm:semidirect}) with $G_0=N=2^{11}$ and $G_1=M_{24}$. For $g\in K$, a $112\times112$ $\mathbb{F}_2$ matrix acting on row vectors, write $g=n_u h$ with $n_u\in N=O_2(K)$ and $h\in H=\mathrm{Stab}(0)$. The element $n_u$ has cocode coordinate $u\in A_{2048}$, while $h$ acts on the $24$ points of the book's $M_{24}$ through the stars of the $276$ pairs. Its label is $\lambda_K(g)=u+2048\,\lambda_{M_{24}}(h)$.\par
\noindent\textbf{Invocation:} \texttt{Data: scripts/\allowbreak{}data/\allowbreak{}J4K\_\allowbreak{}law.txt (scripts/\allowbreak{}j4\_\allowbreak{}k\_\allowbreak{}export.g). Self-test: python3 scripts/\allowbreak{}j4\_\allowbreak{}k\_\allowbreak{}law.py}.\par
\noindent\textbf{Arguments:} none. \textbf{Dependencies:} repository modules digital\_\allowbreak{}laws; external packages numpy.\par
\noindent\textbf{Reads:} no literal file name. \textbf{Writes:} no literal file name.\par
\paragraph{\texttt{scripts/j4\_k\_lift.g}.} Offline export of the lift data for $K=2^{11}{:}M_{24}<J_4$, with its $M_{24}$ identified with the book's $M_{24}$ (Proposition~\ref{prop:mathieuchain}).\par
\noindent\textbf{Invocation:} \texttt{cd scripts \&\& gap -q j4\_\allowbreak{}k\_\allowbreak{}lift.g}.\par
\noindent\textbf{Arguments:} none. \textbf{Dependencies:} GAP 4.14 or compatible.\par
\noindent\textbf{Reads:} /\allowbreak{}tmp/\allowbreak{}j4\_\allowbreak{}trans.g, /\allowbreak{}tmp/\allowbreak{}m24book.g. \textbf{Writes:} no literal file name.\par
\paragraph{\texttt{scripts/j4\_k\_struct.g}.} Offline structural computation for $K=2^{11}{:}M_{24}<J_4$. The quotient $V/V_{100}$ is $12$-dimensional with an $11$-dimensional submodule and trivial top. The group $K$ acts affinely on the $2048$ vectors outside that submodule; $N=O_2(K)$ acts regularly there and a point stabiliser is a complement $H$.\par
\noindent\textbf{Invocation:} \texttt{cd scripts \&\& gap -q j4\_\allowbreak{}k\_\allowbreak{}struct.g}.\par
\noindent\textbf{Arguments:} none. \textbf{Dependencies:} GAP 4.14 or compatible.\par
\noindent\textbf{Reads:} no literal file name. \textbf{Writes:} J4K\_\allowbreak{}struct.g.\par
\paragraph{\texttt{scripts/j4\_law.py}.} Full digital law of $J_4$ on $112\times112$ $\mathbb{F}_2$ matrices acting on row vectors: $\lambda(g)=t+173067389\,\lambda_K(gs_t^{-1})$, where $K=2^{11}{:}M_{24}$ has the law of \texttt{j4\_k\_law.py}. Decoding uses \texttt{J4KLaw.load\_lift} and \texttt{J4K\_lift\_112\_f2.txt}; the full product is tested end to end.\par
\noindent\textbf{Invocation:} \texttt{cd DigitalFiniteGroupTheory \&\& python3 -u scripts/\allowbreak{}j4\_\allowbreak{}law.py --paths 25 --pairs 20 | tee drafts/\allowbreak{}j4\_\allowbreak{}law.log}.\par
\noindent\textbf{Arguments:} --paths default=25; --pairs default=20; --len default=30; --seed default=0. \textbf{Dependencies:} repository modules coset\_\allowbreak{}law, j4\_\allowbreak{}k\_\allowbreak{}law; external packages numpy.\par
\noindent\textbf{Reads:} J4: K-law built; decoder \{name!r\}; generators.txt:, generators.txt. \textbf{Writes:} no literal file name.\par
\paragraph{\texttt{scripts/k\_co2\_chain.py}.} Native chain law of $\mathrm{Co}_2$ inside $K=C_{G_{x0}}(\beta^+)$, used in the law of $K$ (Proposition~\ref{prop:babyaxis}). The image of $K$ in $\mathrm{Co}_1$, acting on $\Lambda/2\Lambda$, is the stabiliser of the type-$2$ vector \texttt{0x200} of $\beta^+$. A random Schreier--Sims chain is complete exactly when the product of its basic orbit sizes equals $|\mathrm{Co}_2|$; each digit is the rank of an image point in its sorted basic orbit.\par
\noindent\textbf{Invocation:} \texttt{.venv-mm/\allowbreak{}bin/\allowbreak{}python -u scripts/\allowbreak{}k\_\allowbreak{}co2\_\allowbreak{}chain.py}.\par
\noindent\textbf{Arguments:} none. \textbf{Dependencies:} external packages mmgroup, numpy.\par
\noindent\textbf{Reads:} no literal file name. \textbf{Writes:} k\_\allowbreak{}co2\_\allowbreak{}chain.pkl.\par
\paragraph{\texttt{scripts/k\_law.py}.} Digital law of $K=C_{G_{x0}}(\beta^+)=C_Q(\beta^+).\mathrm{Co}_2$, of order $2^{24}|\mathrm{Co}_2|$ (Proposition~\ref{prop:babyaxis}), using \texttt{mmgroup}~1.0.8. Its label is $\lambda_K(k)=\operatorname{code}\bigl(k\,l(\pi k)^{-1}\bigr)+2^{24}\lambda_{\mathrm{Co}_2}(\pi k)$, where $\pi k$ is the action on $\Lambda/2\Lambda$, $\lambda_{\mathrm{Co}_2}$ is read from the fixed chain \texttt{scripts/data/k\_co2\_chain.pkl} built by \texttt{k\_co2\_chain.py}, and $l$ is assembled from Schreier-tree transversals.\par
\noindent\textbf{Invocation:} \texttt{.venv-mm/\allowbreak{}bin/\allowbreak{}python -u scripts/\allowbreak{}k\_\allowbreak{}law.py}.\par
\noindent\textbf{Arguments:} none. \textbf{Dependencies:} repository modules k\_\allowbreak{}co2\_\allowbreak{}chain.\par
\noindent\textbf{Reads:} no literal file name. \textbf{Writes:} no literal file name.\par
\paragraph{\texttt{scripts/k\_law\_probe.py}.} First probe for a digital law of $K=C_{G_{x0}}(\beta^+)$, of expected structure $C_Q(\beta^+).\mathrm{Co}_2$ with $Q=Q_{x0}=2^{1+24}$. It checks in \texttt{mmgroup}~1.0.8 the ingredients of the label $\lambda_K(k)=\operatorname{code}\bigl(k\,l(\pi k)^{-1}\bigr)+2^{24}\lambda_{\mathrm{Co}_2}(\pi k)$, where $\pi\colon K\to\mathrm{Co}_2$ is the action on $\Lambda/2\Lambda$ and $l$ a fixed section (Proposition~\ref{prop:babyaxis}).\par
\noindent\textbf{Invocation:} \texttt{.venv-mm/\allowbreak{}bin/\allowbreak{}python -u scripts/\allowbreak{}k\_\allowbreak{}law\_\allowbreak{}probe.py}.\par
\noindent\textbf{Arguments:} none. \textbf{Dependencies:} external packages mmgroup, numpy.\par
\noindent\textbf{Reads:} no literal file name. \textbf{Writes:} no literal file name.\par
\paragraph{\texttt{scripts/k\_law\_step2.py}.} Probe of the bridge from the $\mathrm{Co}_1$ action of \texttt{mmgroup} to the book's $\mathrm{Co}_2$ lattice coordinates.\par
\noindent\textbf{Invocation:} \texttt{cd scripts \&\& python3 k\_\allowbreak{}law\_\allowbreak{}step2.py}.\par
\noindent\textbf{Arguments:} none. \textbf{Dependencies:} external packages mmgroup, numpy.\par
\noindent\textbf{Reads:} no literal file name. \textbf{Writes:} no literal file name.\par
\paragraph{\texttt{scripts/k\_rb\_2a0.py}.} Positions $r_B$ on the $K$-orbit $2A0$ of $93150$ baby axes. The $2A$-involution of an axis $a$ is $t_a=x_\lambda x_{-1}^e\in Q$, with $\lambda$ in the first basic orbit (of size $46575$) of the $\mathrm{Co}_2$ chain of Theorem~\ref{thm:klaw}; the script checks that $r_B(a)=2\operatorname{rank}(\lambda)+e$ is a bijection onto $\{0,\dots,93149\}$.\par
\noindent\textbf{Invocation:} \texttt{.venv-mm/\allowbreak{}bin/\allowbreak{}python -u scripts/\allowbreak{}k\_\allowbreak{}rb\_\allowbreak{}2a0.py}.\par
\noindent\textbf{Arguments:} none. \textbf{Dependencies:} repository modules k\_\allowbreak{}law; external packages mmgroup.\par
\noindent\textbf{Reads:} no literal file name. \textbf{Writes:} no literal file name.\par
\paragraph{\texttt{scripts/k\_rb\_2b0.py}.} Positions $r_B$ on the $K$-orbit $2B0$ of $262310400=1024650\cdot256$ baby axes. The $\mathrm{Co}_2$-coordinate is the $8$-dimensional space $S(a)=[\Lambda/2\Lambda,\pi(t_a)]$ in reduced echelon form; the fibre coordinate is the rank of $t_a$ among its $256$ $C_Q(\beta^+)$-conjugates. Then $r_B(a)=7379551+256\operatorname{rank}(S(a))+\text{fibre rank}$.\par
\noindent\textbf{Invocation:} \texttt{.venv-mm/\allowbreak{}bin/\allowbreak{}python -u scripts/\allowbreak{}k\_\allowbreak{}rb\_\allowbreak{}2b0.py}.\par
\noindent\textbf{Arguments:} none. \textbf{Dependencies:} repository modules k\_\allowbreak{}law; external packages mmgroup.\par
\noindent\textbf{Reads:} no literal file name. \textbf{Writes:} no literal file name.\par
\paragraph{\texttt{scripts/k\_rb\_2b1.py}.} Positions $r_B$ on the $K$-orbit $2B1$ of $7286400=56925\cdot128$ baby axes. The $\mathrm{Co}_2$-coordinate is the $8$-dimensional space $S(a)=[\Lambda/2\Lambda,\pi(t_a)]$ in reduced echelon form; the fibre coordinate is the rank of $t_a$ among its $128$ $C_Q(\beta^+)$-conjugates. Then $r_B(a)=93151+128\operatorname{rank}(S(a))+\text{fibre rank}$.\par
\noindent\textbf{Invocation:} \texttt{.venv-mm/\allowbreak{}bin/\allowbreak{}python -u scripts/\allowbreak{}k\_\allowbreak{}rb\_\allowbreak{}2b1.py}.\par
\noindent\textbf{Arguments:} none. \textbf{Dependencies:} repository modules k\_\allowbreak{}law; external packages mmgroup.\par
\noindent\textbf{Reads:} no literal file name. \textbf{Writes:} no literal file name.\par
\paragraph{\texttt{scripts/k\_rb\_2b1\_probe.py}.} Probe of the stabiliser $H$ in $K$ of the $2B1$ baby-axis representative; the orbit has size $7286400$ and $|H|=2^{35}\cdot3^4\cdot5\cdot7$.\par
\noindent\textbf{Invocation:} \texttt{cd scripts \&\& python3 k\_\allowbreak{}rb\_\allowbreak{}2b1\_\allowbreak{}probe.py}.\par
\noindent\textbf{Arguments:} none. \textbf{Dependencies:} repository modules k\_\allowbreak{}law; external packages mmgroup.\par
\noindent\textbf{Reads:} no literal file name. \textbf{Writes:} no literal file name.\par
\paragraph{\texttt{scripts/k\_rb\_4a1\_custom.py}.} Orbit tester for custom $\mathrm{Co}_2$-invariants on the BabyAxis orbits $4A1$ and $10A1$. A chosen tuple of parts (the commutator space $S$, radicals $R_n$ or short-vector sets $V_n$) of the representative axis is moved by the $17$ generators of the fixed $\mathrm{Co}_2$ chain; spaces are stored as row-reduced spans, sets as sorted vector lists, and seen states as $16$-byte BLAKE2b digests.\par
\noindent\textbf{Invocation:} \texttt{.venv-mm/\allowbreak{}bin/\allowbreak{}python -u scripts/\allowbreak{}k\_\allowbreak{}rb\_\allowbreak{}4a1\_\allowbreak{}custom.py ORBIT PART [PART ...]}.\par
\noindent\textbf{Arguments:} positional argument indices 1, 2:. \textbf{Dependencies:} repository modules k\_\allowbreak{}law; external packages mmgroup.\par
\noindent\textbf{Reads:} no literal file name. \textbf{Writes:} no literal file name.\par
\paragraph{\texttt{scripts/k\_rb\_finer\_probe.py}.} Inspects candidate axis-level invariants for the unresolved $K$-orbits $4B1$ and $6C1$.\par
\noindent\textbf{Invocation:} \texttt{cd scripts \&\& python3 k\_\allowbreak{}rb\_\allowbreak{}finer\_\allowbreak{}probe.py}.\par
\noindent\textbf{Arguments:} none. \textbf{Dependencies:} external packages mmgroup.\par
\noindent\textbf{Reads:} no literal file name. \textbf{Writes:} no literal file name.\par
\paragraph{\texttt{scripts/k\_rb\_finer\_probe2.py}.} Tests canonical axis-level candidate invariants on the unresolved $K$-orbits.\par
\noindent\textbf{Invocation:} \texttt{cd scripts \&\& python3 k\_\allowbreak{}rb\_\allowbreak{}finer\_\allowbreak{}probe2.py}.\par
\noindent\textbf{Arguments:} none. \textbf{Dependencies:} external packages mmgroup.\par
\noindent\textbf{Reads:} no literal file name. \textbf{Writes:} no literal file name.\par
\paragraph{\texttt{scripts/k\_rb\_kernel\_probe.py}.} Computes $X=H\cap C_Q(\beta^+)$ for the baby-axis orbits $4B1$ and $6C1$ by sifting stabiliser samples through a chain of $\pi(H)$.\par
\noindent\textbf{Invocation:} \texttt{cd scripts \&\& python3 k\_\allowbreak{}rb\_\allowbreak{}kernel\_\allowbreak{}probe.py}.\par
\noindent\textbf{Arguments:} none. \textbf{Dependencies:} repository modules k\_\allowbreak{}law; external packages mmgroup.\par
\noindent\textbf{Reads:} no literal file name. \textbf{Writes:} no literal file name.\par
\paragraph{\texttt{scripts/k\_rb\_next\_probe.py}.} Probes the $K$-orbits $2B0$, $4B1$ and $6C1$ of baby axes for the method of Lemma~\ref{lem:rb2b1}.\par
\noindent\textbf{Invocation:} \texttt{cd scripts \&\& python3 k\_\allowbreak{}rb\_\allowbreak{}next\_\allowbreak{}probe.py}.\par
\noindent\textbf{Arguments:} none. \textbf{Dependencies:} repository modules k\_\allowbreak{}law; external packages mmgroup.\par
\noindent\textbf{Reads:} no literal file name. \textbf{Writes:} no literal file name.\par
\paragraph{\texttt{scripts/k\_rb\_pos.py}.} Positions $r_B$ on the baby-axis $K$-orbits $4B1$ and $6C1$, and on any orbit where $(S,R)$ is a $\mathrm{Co}_2$-coordinate. With $a_0$ the representative, $H=\mathrm{Stab}_K(a_0)$ and $\pi\colon K\to\mathrm{Co}_2$, the coordinate is $P(a)=(S(a),R(a))$, where $S(a)=[\Lambda/2\Lambda,\pi(t_a)]$ and $R(a)$ is the radical of the short vectors (of norm parameter $1$ for $4B1$, $3$ for $6C1$). The $\mathrm{Co}_2$-orbit $W$ of $P(a_0)$ has exactly $|\mathrm{Co}_2:\pi(H)|$ elements, so $\mathrm{Stab}(P(a_0))=\pi(H)$; $W$ is ordered with $P(a_0)$ first and the rest sorted.\par
\noindent\textbf{Invocation:} \texttt{nohup .venv-mm/\allowbreak{}bin/\allowbreak{}python -u scripts/\allowbreak{}k\_\allowbreak{}rb\_\allowbreak{}pos.py 4B1 > drafts/\allowbreak{}k\_\allowbreak{}rb\_\allowbreak{}4b1.log 2>\&1 \&}.\par
\noindent\textbf{Arguments:} positional argument indices 1, 2. \textbf{Dependencies:} repository modules k\_\allowbreak{}law; external packages mmgroup.\par
\noindent\textbf{Reads:} no literal file name. \textbf{Writes:} no literal file name.\par
\paragraph{\texttt{scripts/k\_rb\_radical\_equiv.py}.} Equivariance test for the short-vector radical $R(a)$ on the $K$-orbits $4B1$ and $6C1$.\par
\noindent\textbf{Invocation:} \texttt{cd scripts \&\& python3 k\_\allowbreak{}rb\_\allowbreak{}radical\_\allowbreak{}equiv.py}.\par
\noindent\textbf{Arguments:} none. \textbf{Dependencies:} repository modules k\_\allowbreak{}law; external packages mmgroup.\par
\noindent\textbf{Reads:} no literal file name. \textbf{Writes:} no literal file name.\par
\paragraph{\texttt{scripts/k\_rb\_radical\_probe.py}.} Probe of the equivariance of short-vector radical data on the $K$-orbits $4B1$ and $6C1$.\par
\noindent\textbf{Invocation:} \texttt{cd scripts \&\& python3 k\_\allowbreak{}rb\_\allowbreak{}radical\_\allowbreak{}probe.py}.\par
\noindent\textbf{Arguments:} none. \textbf{Dependencies:} external packages mmgroup, numpy.\par
\noindent\textbf{Reads:} no literal file name. \textbf{Writes:} no literal file name.\par
\paragraph{\texttt{scripts/k\_rb\_sorbit.py}.} $\mathrm{Co}_2$-orbit of $S(a)=[\Lambda/2\Lambda,\pi(t_a)]$ for the baby-axis $K$-orbits $4B1$ and $6C1$. The index of $\pi(H)$ in $\mathrm{Co}_2$ is known from \texttt{k\_rb\_kernel\_probe.py} ($3586275$ for $4B1$, $1593900$ for $6C1$); since $\pi(H)$ fixes $S(a_0)$, the orbit is at most this large, and $S$ is a $\mathrm{Co}_2$-coordinate exactly when equality holds. Otherwise the search is repeated on the pair $(S,R)$. Each space is packed into one integer, $24$ bits per basis vector of its reduced echelon form.\par
\noindent\textbf{Invocation:} \texttt{nohup .venv-mm/\allowbreak{}bin/\allowbreak{}python -u scripts/\allowbreak{}k\_\allowbreak{}rb\_\allowbreak{}sorbit.py > drafts/\allowbreak{}k\_\allowbreak{}rb\_\allowbreak{}sorbit.log 2>\&1 \&}.\par
\noindent\textbf{Arguments:} none. \textbf{Dependencies:} repository modules k\_\allowbreak{}law; external packages mmgroup.\par
\noindent\textbf{Reads:} no literal file name. \textbf{Writes:} no literal file name.\par
\paragraph{\texttt{scripts/k\_rb\_stab\_faithful.py}.} Faithfulness of $\pi(H)=2^{10}{:}L_3(4){:}2_2$ on the kernel $X=H\cap Q$ for the BabyAxis orbit $4A1$ (Remark~\ref{rem:rbstab}). $X$ is rebuilt from $60$ seeded samples of $H$, and $30$ of them act by conjugation on the $4095$ nonzero vectors of $X$. A randomized Schreier--Sims run gives a lower bound for the order of the image of $H$ in $\mathrm{GL}(X)$; since $X$ is abelian the action factors through $\pi(H)$, so the order is at most $|\pi(H)|=41287680$, and reaching this bound proves faithfulness. Optional argument: a seed (default $1$).\par
\noindent\textbf{Invocation:} \texttt{cd DigitalFiniteGroupTheory \&\& .venv-mm/\allowbreak{}bin/\allowbreak{}python -u scripts/\allowbreak{}k\_\allowbreak{}rb\_\allowbreak{}stab\_\allowbreak{}faithful.py [SEED] (default 1)}.\par
\noindent\textbf{Arguments:} positional argument indices 1. \textbf{Dependencies:} repository modules k\_\allowbreak{}rb\_\allowbreak{}stab\_\allowbreak{}probe; external packages mmgroup.\par
\noindent\textbf{Reads:} no literal file name. \textbf{Writes:} no literal file name.\par
\paragraph{\texttt{scripts/k\_rb\_stab\_probe.py}.} Stabilisers of the BabyAxis orbits $4A1$ and $10A1$ (Remark~\ref{rem:rbstab}). Elements of $H=\mathrm{Stab}_K(a_0)$ are sampled as in \texttt{k\_rb\_pos.py} from mmgroup's random generator, seeded by a fixed seed; another seed changes the counts but not the conclusions. Optional arguments: the number of order samples, the number of module samples and the seed (defaults $120$, $60$, $1$); runs in a few seconds.\par
\noindent\textbf{Invocation:} \texttt{cd DigitalFiniteGroupTheory \&\& .venv-mm/\allowbreak{}bin/\allowbreak{}python -u scripts/\allowbreak{}k\_\allowbreak{}rb\_\allowbreak{}stab\_\allowbreak{}probe.py [N\_\allowbreak{}ORDERS] [N\_\allowbreak{}MODULE] [SEED]}.\par
\noindent\textbf{Arguments:} positional argument indices 1, 2, 3. \textbf{Dependencies:} repository modules k\_\allowbreak{}law; external packages mmgroup.\par
\noindent\textbf{Reads:} no literal file name. \textbf{Writes:} no literal file name.\par
\paragraph{\texttt{scripts/ly\_chain.py}.} Canonical (lexicographic) chain law of the Lyons group $\mathrm{Ly}$ on $8835156$ points.\par
\noindent\textbf{Invocation:} \texttt{cd DigitalFiniteGroupTheory \&\& python3 -u scripts/\allowbreak{}ly\_\allowbreak{}chain.py > drafts/\allowbreak{}ly\_\allowbreak{}chain.log 2>\&1}.\par
\noindent\textbf{Arguments:} none. \textbf{Dependencies:} repository modules digital\_\allowbreak{}laws; external packages numpy.\par
\noindent\textbf{Reads:} no literal file name. \textbf{Writes:} no literal file name.\par
\paragraph{\texttt{scripts/ly\_subdegrees.py}.} Subdegrees of $\mathrm{Ly}$ on $8835156$ points: the orbit lengths of the stabiliser of $0$, isomorphic to $G_2(5)$.\par
\noindent\textbf{Invocation:} \texttt{cd scripts \&\& python3 ly\_\allowbreak{}subdegrees.py}.\par
\noindent\textbf{Arguments:} none. \textbf{Dependencies:} repository modules digital\_\allowbreak{}laws; external packages numpy.\par
\noindent\textbf{Reads:} no literal file name. \textbf{Writes:} no literal file name.\par
\paragraph{\texttt{scripts/primsporadic\_chain.py}.} Stabiliser-chain digital laws of $J_1$, $J_2$, $\mathrm{HS}$, $\mathrm{Suz}$, $\mathrm{He}$ and $\mathrm{Ru}$ from GAP's library of primitive groups.\par
\noindent\textbf{Invocation:} \texttt{Chain laws of J1, J2, HS, Suz, He, Ru from GAP's primitive groups library. Usage: cd scripts \&\& python3 primsporadic\_\allowbreak{}chain.py}.\par
\noindent\textbf{Arguments:} none. \textbf{Dependencies:} external packages numpy.\par
\noindent\textbf{Reads:} ../\allowbreak{}drafts/\allowbreak{}primsporadic\_\allowbreak{}gap.log. \textbf{Writes:} ../\allowbreak{}drafts/\allowbreak{}primsporadic.json, ../\allowbreak{}drafts/\allowbreak{}primsporadic.tex.\par
\paragraph{\texttt{scripts/th\_k\_law.py}.} Canonical lexicographic law (Section~\ref{sec:canonical}) of $K={}^3D_4(2){:}3<\mathrm{Th}$, acting by conjugation on its $819$ long-root involutions $z_0,\dots,z_{818}$. For a $248\times248$ $\mathbb{F}_2$ matrix $k\in K$, the permutation $i\mapsto i^k$ is obtained from $z_i^k=k^{-1}z_i k$ and identified by the key $a_0*z$. The permutation data come from \texttt{scripts/th\_k\_perm.g}.\par
\noindent\textbf{Invocation:} \texttt{permutation i -> i\textasciicircum{}k is read off from z\_\allowbreak{}i\textasciicircum{}k = k\textasciicircum{}-1 z\_\allowbreak{}i k, identified by the key a0 * z. Self-test: python3 scripts/\allowbreak{}th\_\allowbreak{}k\_\allowbreak{}law.py}.\par
\noindent\textbf{Arguments:} none. \textbf{Dependencies:} repository modules digital\_\allowbreak{}laws; external packages numpy.\par
\noindent\textbf{Reads:} no literal file name. \textbf{Writes:} no literal file name.\par
\paragraph{\texttt{scripts/th\_k\_lift.py}.} Decoding for $K={}^3D_4(2){:}3<\mathrm{Th}$ without GAP. The group acts faithfully on its $819$ long-root involutions. A Schreier--Sims chain carries a $248\times248$ $\mathbb{F}_2$ row-vector matrix with every transversal element, so a permutation of $K$ lifts uniquely by sifting. The chain is cached in \texttt{scripts/data/ThK\_lift\_chain.npz}.\par
\noindent\textbf{Invocation:} \texttt{cd DigitalFiniteGroupTheory \&\& python3 -u scripts/\allowbreak{}th\_\allowbreak{}k\_\allowbreak{}lift.py | tee drafts/\allowbreak{}th\_\allowbreak{}k\_\allowbreak{}lift.log}.\par
\noindent\textbf{Arguments:} none. \textbf{Dependencies:} repository modules coset\_\allowbreak{}law, th\_\allowbreak{}k\_\allowbreak{}law; external packages numpy.\par
\noindent\textbf{Reads:} no literal file name. \textbf{Writes:} no literal file name.\par
\paragraph{\texttt{scripts/th\_k\_perm.g}.} Offline computation for $K={}^3D_4(2){:}3<\mathrm{Th}$ acting by conjugation on its $819$ long-root involutions.\par
\noindent\textbf{Invocation:} \texttt{cd scripts \&\& gap -q th\_\allowbreak{}k\_\allowbreak{}perm.g}.\par
\noindent\textbf{Arguments:} none. \textbf{Dependencies:} GAP 4.14 or compatible.\par
\noindent\textbf{Reads:} no literal file name. \textbf{Writes:} no literal file name.\par
\paragraph{\texttt{scripts/th\_law.py}.} Full digital law of $\mathrm{Th}$ on $248\times248$ $\mathbb{F}_2$ matrices acting on row vectors: $\lambda(g)=t+143127000\,\lambda_K(gs_t^{-1})$, where $K={}^3D_4(2){:}3$ has the lexicographic law of \texttt{th\_k\_law.py}.\par
\noindent\textbf{Invocation:} \texttt{cd DigitalFiniteGroupTheory \&\& python3 -u scripts/\allowbreak{}th\_\allowbreak{}law.py --paths 25 --pairs 20 | tee drafts/\allowbreak{}th\_\allowbreak{}law.log}.\par
\noindent\textbf{Arguments:} --paths default=25; --pairs default=20; --len default=30; --seed default=0. \textbf{Dependencies:} repository modules coset\_\allowbreak{}law, th\_\allowbreak{}k\_\allowbreak{}law, th\_\allowbreak{}k\_\allowbreak{}lift; external packages numpy.\par
\noindent\textbf{Reads:} no literal file name. \textbf{Writes:} no literal file name.\par
\paragraph{\texttt{scripts/th\_orbit.py}.} Orbit of $\mathrm{Th}$ on the cosets of $K={}^3D_4(2){:}3$, realised as the orbit of the $2$-dimensional $K$-submodule $W$ in the $248$-dimensional $\mathbb{F}_2$-module.\par
\noindent\textbf{Invocation:} \texttt{cd DigitalFiniteGroupTheory \&\& python3 -u scripts/\allowbreak{}th\_\allowbreak{}orbit.py > drafts/\allowbreak{}th\_\allowbreak{}orbit.log 2>\&1}.\par
\noindent\textbf{Arguments:} --max-points default=0; --out default=os.path.join(DATA, 'th\_\allowbreak{}orbit'); --chunk default=1 << 21. \textbf{Dependencies:} external packages numpy.\par
\noindent\textbf{Reads:} no literal file name. \textbf{Writes:} no literal file name.\par

\section{Verification and worked examples}
\paragraph{\texttt{scripts/build\_index.py}.} Sorted key index of a coset orbit produced by \texttt{th\_orbit.py} or \texttt{vec\_orbit.py}, whose \texttt{keys.npy} lists $64$-bit keys in discovery order $t$. Writes \texttt{index.npy}, the permutation sorting the keys, and \texttt{sorted\_keys.npy}, so the point with key $k$ is found by binary search. Documented runs (from the archive root): \texttt{scripts/data/th\_orbit} with \texttt{--expect 143127000}, and \texttt{scripts/data/j4\_orbit} with \texttt{--expect 173067389}.\par
\noindent\textbf{Invocation:} \texttt{cd DigitalFiniteGroupTheory \&\& python3 -u scripts/\allowbreak{}build\_\allowbreak{}index.py scripts/\allowbreak{}data/\allowbreak{}th\_\allowbreak{}orbit --expect 143127000}.\par
\noindent\textbf{Arguments:} dir; --expect default=0; --samples default=100000. \textbf{Dependencies:} external packages numpy.\par
\noindent\textbf{Reads:} no literal file name. \textbf{Writes:} no literal file name.\par
\paragraph{\texttt{scripts/canonical\_labels.py}.} Lexicographic (canonical) chain laws: the labels depend only on the permutation group.\par
\noindent\textbf{Invocation:} \texttt{Lexicographic (canonical) chain laws: labels depend only on the permutation group. Usage: cd scripts \&\& python3 canonical\_\allowbreak{}labels.py}.\par
\noindent\textbf{Arguments:} none. \textbf{Dependencies:} external packages numpy.\par
\noindent\textbf{Reads:} data/\allowbreak{}. \textbf{Writes:} ../\allowbreak{}drafts/\allowbreak{}canonical.json, ../\allowbreak{}drafts/\allowbreak{}canonical.tex.\par
\paragraph{\texttt{scripts/check\_lie\_families.py}.} Runs the Lie-family checks of the book and stops with an error on any discrepancy: \texttt{classical\_cells} (Example~\ref{ex:classical}), \texttt{g2\_cells} (Example~\ref{ex:g2}), \texttt{triality\_cells} (Section~\ref{sec:3d4}) and \texttt{psln\_dense} (Theorems~\ref{thm:psln} and~\ref{thm:pslndense}, Examples~\ref{ex:psl3} and~\ref{ex:pslndense}). Takes about one minute.\par
\noindent\textbf{Invocation:} \texttt{python3 check\_\allowbreak{}lie\_\allowbreak{}families.py (about one minute)}.\par
\noindent\textbf{Arguments:} none. \textbf{Dependencies:} repository modules classical\_\allowbreak{}cells, g2\_\allowbreak{}cells, psln\_\allowbreak{}dense, triality\_\allowbreak{}cells, triality\_\allowbreak{}label, triality\_\allowbreak{}mult.\par
\noindent\textbf{Reads:} no literal file name. \textbf{Writes:} no literal file name.\par
\paragraph{\texttt{scripts/gf.py}.} Finite fields $\mathbb{F}_{p^k}$ for the $F_4$ programs. Elements are the integers $0,\dots,q-1$; the base-$p$ digits of $a$, constant term first, are the coefficients of a polynomial in $x$, reduced modulo a fixed monic primitive polynomial of degree $k$ (the Conway polynomial where tabulated). For $k=1$ the arithmetic agrees with the mod-$p$ code of \texttt{f4\_core.py}, \texttt{f4\_law.py} and \texttt{f4\_adjoint.py}.\par
\noindent\textbf{Invocation:} \texttt{cd scripts \&\& python3 gf.py}.\par
\noindent\textbf{Arguments:} none. \textbf{Dependencies:} external packages numpy.\par
\noindent\textbf{Reads:} no literal file name. \textbf{Writes:} no literal file name.\par
\paragraph{\texttt{scripts/vec\_orbit.py}.} Orbit of a vector (\texttt{--mode vector}) or of a $2$-dimensional subspace (\texttt{--mode plane}) of a $\mathbb{F}_2$-module under the matrix group generated by the matrices in \texttt{--gens} and their inverses. A point is identified by a $64$-bit hash of its least nonzero vector. Documented run: the $J_4$-orbit of length $173067389$ on the $112$-dimensional module, written to \texttt{scripts/data/j4\_orbit}.\par
\noindent\textbf{Invocation:} \texttt{python3 -u scripts/\allowbreak{}vec\_\allowbreak{}orbit.py --gens scripts/\allowbreak{}data/\allowbreak{}J4\_\allowbreak{}112\_\allowbreak{}f2.txt --start scripts/\allowbreak{}data/\allowbreak{}J4Kfix\_\allowbreak{}112\_\allowbreak{}f2.txt --mode vector --index 173067389 --check scripts/\allowbreak{}data/\allowbreak{}J4K\_\allowbreak{}112\_\allowbreak{}f2.txt --out scripts/\allowbreak{}data/\allowbreak{}j4\_\allowbreak{}orbit}.\par
\noindent\textbf{Arguments:} --gens; --start; --mode; --index; --check; --out; --max-points default=0; --chunk default=1 << 21. \textbf{Dependencies:} external packages numpy.\par
\noindent\textbf{Reads:} no literal file name. \textbf{Writes:} generators.txt.\par
\paragraph{\texttt{scripts/worked\_co0.py}.} Checks the worked $\mathrm{Co}_0$ product of Example~\ref{ex:co0product} (Theorem~\ref{thm:co0chain}). Builds the chain of \texttt{co0\_chain.py} with its outputs redirected to the system temporary directory, prints the steps and asserts every number in the example; nothing in \texttt{drafts/} is written.\par
\noindent\textbf{Invocation:} \texttt{Worked Co0 product for Chapter 13 (Example ex:co0product; Theorem thm:co0chain). Usage: cd scripts \&\& python3 worked\_\allowbreak{}co0.py Builds the chain of co0\_\allowbreak{}chain.py with its outputs redirected to the system temp directory (nothing in drafts/\allowbreak{} is written), prints the steps and asserts every number in the ex}.\par
\noindent\textbf{Arguments:} none. \textbf{Dependencies:} repository modules \_\allowbreak{}tmpwrites; external packages numpy.\par
\noindent\textbf{Reads:} ../\allowbreak{}drafts/\allowbreak{}. \textbf{Writes:} no literal file name.\par
\paragraph{\texttt{scripts/worked\_exceptional.py}.} Step-by-step products $m*n$ in $E_6(q)$, $E_7(q)$, $E_8(q)$ and ${}^2E_6(q)$ through the law programs. Arguments: the group (\texttt{e6}, \texttt{e7}, \texttt{e8} or \texttt{2e6}) and $q$. Prints the decoded normal forms, the product and cross-checks.\par
\noindent\textbf{Invocation:} \texttt{python3 worked\_\allowbreak{}exceptional.py GROUP q (GROUP in e6,e7,e8,2e6). Prints decoded normal forms, the product and cross-checks.}.\par
\noindent\textbf{Arguments:} positional argument indices 1, 2. \textbf{Dependencies:} Python standard library.\par
\noindent\textbf{Reads:} no literal file name. \textbf{Writes:} no literal file name.\par
\paragraph{\texttt{scripts/worked\_f4.py}.} Worked products in $F_4(p)$, $p$ prime: labels as in Corollary~\ref{cor:f4law}, product by \texttt{f4\_law.multiply}, checked in the $52$-dimensional adjoint model. Arguments: $p$ and an optional seed. Digit strings are little-endian in base $p$, and $\gamma$ is the least primitive root.\par
\noindent\textbf{Invocation:} \texttt{python3 worked\_\allowbreak{}f4.py p [seed]. Conventions: cells ordered as list(L['fixed'].values()); digit strings little-endian base p; gamma = least primitive root.}.\par
\noindent\textbf{Arguments:} positional argument indices 1, 2. \textbf{Dependencies:} repository modules f4\_\allowbreak{}adjoint, f4\_\allowbreak{}adjointq, f4\_\allowbreak{}lawq, gf; external packages numpy.\par
\noindent\textbf{Reads:} no literal file name. \textbf{Writes:} no literal file name.\par
\paragraph{\texttt{scripts/worked\_f42.py}.} Worked products in ${}^2F_4(8)$ through \texttt{f42\_digital.py}, checked against the matrix product in $F_4(8)$. Optional argument: a seed.\par
\noindent\textbf{Invocation:} \texttt{Worked products in 2F4(8) via f42\_\allowbreak{}digital.py (law, dec, enc), checked against the F4 matrix product (toF4, mul8). Usage: python3 worked\_\allowbreak{}f42.py [seed]}.\par
\noindent\textbf{Arguments:} positional argument indices 1. \textbf{Dependencies:} Python standard library.\par
\noindent\textbf{Reads:} no literal file name. \textbf{Writes:} no literal file name.\par
\paragraph{\texttt{scripts/worked\_inverses.py}.} Group inverses for the worked products of Chapter~13 (Remark~\ref{rem:wpchecks}). For every worked product $g*h$ it computes the labels of $g^{-1}$, $h^{-1}$ and $(gh)^{-1}$ and the element orders, and asserts $gg^{-1}=g^{-1}g=e$, $(g^{-1})^{-1}=g$ and $(gh)^{-1}=h^{-1}g^{-1}$ in the digital law. Writes nothing in the project: writes of imported law modules into \texttt{drafts/} go to the system temporary directory, any other write raises an error, and \texttt{monster/*.pkl} are checked unchanged. Needs \texttt{mmgroup}~1.0.8.\par
\noindent\textbf{Invocation:} \texttt{Run from scripts/\allowbreak{}: ../\allowbreak{}.venv-mm/\allowbreak{}bin/\allowbreak{}python3 worked\_\allowbreak{}inverses.py (mmgroup 1.0.8 for the Monster).}.\par
\noindent\textbf{Arguments:} none. \textbf{Dependencies:} repository modules coset\_\allowbreak{}law, j4\_\allowbreak{}law, monster\_\allowbreak{}law, ree\_\allowbreak{}law, th\_\allowbreak{}law, worked\_\allowbreak{}co0, worked\_\allowbreak{}lie\_\allowbreak{}steps, worked\_\allowbreak{}m24; external packages numpy.\par
\noindent\textbf{Reads:} WP.tex. \textbf{Writes:} no literal file name.\par
\paragraph{\texttt{scripts/worked\_j4.py}.} Worked $J_4$ products of Examples~\ref{ex:j4kproduct} and~\ref{ex:j4product} (Proposition~\ref{prop:j4}). Writes no files.\par
\noindent\textbf{Invocation:} \texttt{Worked J4 products for Chapter 13 (Examples ex:j4kproduct, ex:j4product; Proposition prop:j4). Usage: cd scripts \&\& python3 worked\_\allowbreak{}j4.py. It writes no files. Approved product inside K.}.\par
\noindent\textbf{Arguments:} none. \textbf{Dependencies:} repository modules coset\_\allowbreak{}law, j4\_\allowbreak{}law; external packages numpy.\par
\noindent\textbf{Reads:} no literal file name. \textbf{Writes:} no literal file name.\par
\paragraph{\texttt{scripts/worked\_lie\_steps.py}.} Step-by-step worked products $m*n$ for the Lie-type groups of \texttt{worked\_products\_lie\_all.py}, including the label-to-matrix arithmetic for $\PSL_2$ and $\PSL_n$. Arguments: the group and optionally $m$, $n$ (default: those stored in \texttt{drafts/worked\_lie\_<GROUP>.json}). Prints every step and writes \texttt{drafts/worked\_steps\_<GROUP>.txt}.\par
\noindent\textbf{Invocation:} \texttt{cd scripts \&\& python3 worked\_\allowbreak{}lie\_\allowbreak{}steps.py GROUP [m n]}.\par
\noindent\textbf{Arguments:} positional argument indices 1:. \textbf{Dependencies:} repository modules worked\_\allowbreak{}products\_\allowbreak{}lie, worked\_\allowbreak{}products\_\allowbreak{}lie\_\allowbreak{}all.\par
\noindent\textbf{Reads:} ../\allowbreak{}drafts/\allowbreak{}worked\_\allowbreak{}lie\_\allowbreak{}\%s.json. \textbf{Writes:} ../\allowbreak{}drafts/\allowbreak{}worked\_\allowbreak{}steps\_\allowbreak{}\%s.txt.\par
\paragraph{\texttt{scripts/worked\_m24.py}.} Worked $M_{24}$ products of Examples~\ref{ex:m24kproduct} and~\ref{ex:m24product} (Theorem~\ref{thm:golaylaw}). Optional arguments: $x$, $y$. Runs the law of \texttt{golay\_law.py} without writing its JSON output and prints the steps.\par
\noindent\textbf{Invocation:} \texttt{Worked M24 products for Chapter 13 (Examples ex:m24kproduct, ex:m24product; Theorem thm:golaylaw). Usage: python3 worked\_\allowbreak{}m24.py [x y] Runs the law of golay\_\allowbreak{}law.py (without writing its JSON) and prints the steps.}.\par
\noindent\textbf{Arguments:} positional argument indices 1:. \textbf{Dependencies:} repository modules \_\allowbreak{}tmpwrites.\par
\noindent\textbf{Reads:} no literal file name. \textbf{Writes:} no literal file name.\par
\paragraph{\texttt{scripts/worked\_monster.py}.} Checks Examples~\ref{ex:monsterkproduct} and~\ref{ex:monsterproduct}; needs \texttt{mmgroup}~1.0.8. Writes nothing: every write-mode \texttt{open()} raises an error, bytecode caching is off, and the size and modification time of \texttt{monster/*.pkl} are checked at the end.\par
\noindent\textbf{Invocation:} \texttt{../\allowbreak{}.venv-mm/\allowbreak{}bin/\allowbreak{}python3 worked\_\allowbreak{}monster.py}.\par
\noindent\textbf{Arguments:} none. \textbf{Dependencies:} repository modules monster\_\allowbreak{}law.\par
\noindent\textbf{Reads:} no literal file name. \textbf{Writes:} no literal file name.\par
\paragraph{\texttt{scripts/worked\_products.py}.} Worked digital products $m*n$ in the permutation-chain groups of Chapter~13 (Part~A): $M_{11}$, $M_{12}$, $M_{22}$, $M_{23}$, $M_{24}$, $J_1$, $J_2$, $\mathrm{HS}$, $\mathrm{Suz}$, $\mathrm{He}$, $\mathrm{Ru}$, $J_3$, $\mathrm{O'N}$, $\mathrm{HN}$, $\mathrm{Fi}_{22}$, $\mathrm{Fi}_{23}$, $\mathrm{Fi}_{24}$ and $\mathrm{Fi}_{24}'$; for $B$ and $M$ it re-executes itself under \texttt{.venv-mm}. Convention: $(gh)(i)=g(h(i))$, so $h$ acts first. The label of an element is $x=\sum_i t_i\,r_0\cdots r_{i-1}$, where $t_i$ is its digit at chain level $i$. Writes \texttt{drafts/worked\_<GROUP>.json}.\par
\noindent\textbf{Invocation:} \texttt{cd scripts \&\& python3 worked\_\allowbreak{}products.py GROUP [m n] (B, M re-exec under ../\allowbreak{}.venv-mm automatically)}.\par
\noindent\textbf{Arguments:} positional argument indices 1, 2, 3. \textbf{Dependencies:} repository modules coset\_\allowbreak{}law, j4\_\allowbreak{}law, monster\_\allowbreak{}law, th\_\allowbreak{}law; external packages mmgroup, numpy.\par
\noindent\textbf{Reads:} no literal file name. \textbf{Writes:} ../\allowbreak{}drafts/\allowbreak{}worked\_\allowbreak{}\%s.json.\par
\paragraph{\texttt{scripts/worked\_products\_lie.py}.} Worked products $m*n$ in $\PSL(2,p)$ (Theorem~\ref{thm:psl2}) and $\PSL_n(q)$, $q$ prime (Theorem~\ref{thm:psln}), with self-checks.\par
\noindent\textbf{Invocation:} \texttt{Worked products m*n in PSL(2,p) (Thm psl2) and PSL\_\allowbreak{}n(q), q prime (Thm psln), with self-checks. Usage: python3 worked\_\allowbreak{}products\_\allowbreak{}lie.py}.\par
\noindent\textbf{Arguments:} none. \textbf{Dependencies:} Python standard library.\par
\noindent\textbf{Reads:} no literal file name. \textbf{Writes:} no literal file name.\par
\paragraph{\texttt{scripts/worked\_products\_lie\_all.py}.} Worked products $m*n$ in groups of Lie type, with the four standard checks and an independent matrix check. Groups: $\mathrm{Sz}(8)$, $\mathrm{Sz}(32)$, $\mathrm{Sz}(128)$ (Theorem~\ref{thm:suzuki}), ${}^2G_2(27)$, ${}^2G_2(243)$ (Theorem~\ref{thm:ree}) and the nine groups of Example~\ref{ex:classical}, such as $\Sp_4(3)$; \texttt{all} runs every group. Writes \texttt{drafts/worked\_lie\_<GROUP>.json}.\par
\noindent\textbf{Invocation:} \texttt{cd scripts \&\& python3 worked\_\allowbreak{}products\_\allowbreak{}lie\_\allowbreak{}all.py GROUP [m n] GROUP in: Sz8 Sz32 Sz128 Ree27 Ree243, the nine groups of ex:classical (e.g. 'Sp4(3)') | all}.\par
\noindent\textbf{Arguments:} positional argument indices 1, 2, 3. \textbf{Dependencies:} repository modules classical\_\allowbreak{}cells, worked\_\allowbreak{}products\_\allowbreak{}lie.\par
\noindent\textbf{Reads:} no literal file name. \textbf{Writes:} ../\allowbreak{}drafts/\allowbreak{}worked\_\allowbreak{}lie\_\allowbreak{}\%s.json.\par
\paragraph{\texttt{scripts/worked\_th.py}.} Worked Thompson products of Examples~\ref{ex:thkproduct} and~\ref{ex:thproduct} (Proposition~\ref{prop:th}). Uses the cached lift of $K$; writes no files.\par
\noindent\textbf{Invocation:} \texttt{Worked Thompson products for Chapter 13 (Examples ex:thkproduct, ex:thproduct; Proposition prop:th). Usage: cd scripts \&\& python3 worked\_\allowbreak{}th.py. Uses the cached K lift; writes no files.}.\par
\noindent\textbf{Arguments:} none. \textbf{Dependencies:} repository modules coset\_\allowbreak{}law, th\_\allowbreak{}law; external packages numpy.\par
\noindent\textbf{Reads:} no literal file name. \textbf{Writes:} no literal file name.\par

\section{Characters and symmetric groups}
\paragraph{\texttt{scripts/characters.py}.} From laws to character tables (Section~\ref{sec:chartables}): rank-$3$ decompositions, class data, and the character table of $M_{11}$ from the printed tables. Needs \texttt{drafts/char\_gap.txt}, written beforehand by GAP.\par
\noindent\textbf{Invocation:} \texttt{cd scripts \&\& python3 characters.py (after GAP has written ../\allowbreak{}drafts/\allowbreak{}char\_\allowbreak{}gap.txt)}.\par
\noindent\textbf{Arguments:} none. \textbf{Dependencies:} external packages numpy.\par
\noindent\textbf{Reads:} ../\allowbreak{}digital\_\allowbreak{}finite\_\allowbreak{}groups.tex, ../\allowbreak{}drafts/\allowbreak{}char\_\allowbreak{}gap.txt. \textbf{Writes:} ../\allowbreak{}drafts/\allowbreak{}characters.json, ../\allowbreak{}drafts/\allowbreak{}characters.tex.\par
\paragraph{\texttt{scripts/nonab\_char\_verify.py}.} Character tables of metacyclic, dihedral and dicyclic groups (Chapter~\ref{ch:nonabelian}). Each group is built from its two-digit law (Theorem~\ref{thm:metacyclic}, Corollary~\ref{cor:dihedral}, Theorem~\ref{thm:dicyclic}). The irreducible characters are the linear characters together with the irreducible characters induced from the cyclic normal subgroup $A=\langle x\rangle=\{0,\dots,|A|-1\}$. All identities are verified exactly in $\mathbb{Z}[\zeta_N]$.\par
\noindent\textbf{Invocation:} \texttt{python3 scripts/\allowbreak{}nonab\_\allowbreak{}char\_\allowbreak{}verify.py [--no-gap] [--no-write]}.\par
\noindent\textbf{Arguments:} --no-gap action='store\_\allowbreak{}true'; --no-write action='store\_\allowbreak{}true'. \textbf{Dependencies:} Python standard library.\par
\noindent\textbf{Reads:} no literal file name. \textbf{Writes:} no literal file name.\par
\paragraph{\texttt{scripts/pc\_verify.py}.} Presentations and coset enumeration: Todd--Coxeter enumeration (HLT and Felsch strategies), Reidemeister--Schreier and Tietze reduction, implemented from scratch, with a GAP cross-check of every number. Right actions; points and cosets are $0,\dots,n-1$, with coset $0=H$. A word is a list of nonzero integers, $+i$ standing for generator $i$ and $-i$ for its inverse. Writes \texttt{drafts/pc\_verify.log}, \texttt{drafts/pc\_verify.g} and \texttt{chapters/pc\_data.tex}.\par
\noindent\textbf{Invocation:} \texttt{cd scripts \&\& python3 pc\_\allowbreak{}verify.py}.\par
\noindent\textbf{Arguments:} none. \textbf{Dependencies:} Python standard library.\par
\noindent\textbf{Reads:} no literal file name. \textbf{Writes:} pc\_\allowbreak{}data.tex, pc\_\allowbreak{}verify.log.\par
\paragraph{\texttt{scripts/run\_characters.sh}.} Shell wrapper that runs GAP from \texttt{scripts/data}. \emph{[Description to be confirmed; the script currently contains absolute paths to the author's computer.]}\par
\noindent\textbf{Invocation:} \texttt{cd scripts \&\& sh run\_\allowbreak{}characters.sh}.\par
\noindent\textbf{Arguments:} none. \textbf{Dependencies:} Python standard library.\par
\noindent\textbf{Reads:} no literal file name. \textbf{Writes:} no literal file name.\par
\paragraph{\texttt{scripts/sn\_alt\_char\_verify.py}.} Exact checks for the characters of the alternating groups $\Alt_n$ (Chapter~\ref{ch:simple}).\par
\noindent\textbf{Invocation:} \texttt{cd scripts \&\& python3 sn\_\allowbreak{}alt\_\allowbreak{}char\_\allowbreak{}verify.py}.\par
\noindent\textbf{Arguments:} none. \textbf{Dependencies:} repository modules sn\_\allowbreak{}char\_\allowbreak{}verify.\par
\noindent\textbf{Reads:} no literal file name. \textbf{Writes:} no literal file name.\par
\paragraph{\texttt{scripts/sn\_char\_verify.py}.} Exact checks for the characters of the symmetric groups $\Sym_n$ (Chapter~\ref{ch:simple}).\par
\noindent\textbf{Invocation:} \texttt{cd scripts \&\& python3 sn\_\allowbreak{}char\_\allowbreak{}verify.py}.\par
\noindent\textbf{Arguments:} none. \textbf{Dependencies:} Python standard library.\par
\noindent\textbf{Reads:} no literal file name. \textbf{Writes:} no literal file name.\par
\paragraph{\texttt{scripts/sn\_law\_verify.py}.} Exhaustive checks of the Lehmer digital laws of $\Sym_n$ and $\Alt_n$ (Chapter~\ref{ch:simple}).\par
\noindent\textbf{Invocation:} \texttt{cd scripts \&\& python3 sn\_\allowbreak{}law\_\allowbreak{}verify.py}.\par
\noindent\textbf{Arguments:} none. \textbf{Dependencies:} Python standard library.\par
\noindent\textbf{Reads:} no literal file name. \textbf{Writes:} sn\_\allowbreak{}law3.tex, sn\_\allowbreak{}laws.tex.\par
\paragraph{\texttt{scripts/sn\_tables.py}.} Tables of Section~\ref{sec:sn:tables}: class data and character degrees of $\Sym_n$ and $\Alt_n$, all checked for $n\le 10$. Only if every check passes does it write \texttt{chapters/sn\_orders.tex}, \texttt{chapters/sn\_degrees.tex}, \texttt{chapters/sn\_classes8.tex} ($n\le 8$) and \texttt{scripts/data/sn\_tables.txt} ($n\le 10$).\par
\noindent\textbf{Invocation:} \texttt{cd scripts \&\& python3 sn\_\allowbreak{}tables.py}.\par
\noindent\textbf{Arguments:} none. \textbf{Dependencies:} Python standard library.\par
\noindent\textbf{Reads:} no literal file name. \textbf{Writes:} sn\_\allowbreak{}classes8.tex, sn\_\allowbreak{}degrees.tex, sn\_\allowbreak{}orders.tex, sn\_\allowbreak{}tables.txt.\par
\paragraph{\texttt{scripts/sn\_verify.py}.} Verification for Chapter~\ref{ch:simple}; also runs \texttt{nonab\_char\_verify.py} for Chapter~\ref{ch:nonabelian}. Pure Python; the optional argument is the path of a GAP binary, and the GAP cross-check is skipped if GAP is absent. Convention~\eqref{eq:comp}: a permutation of $\{0,\dots,n-1\}$ is a tuple $g$ with $g[i]=g(i)$, and $(gh)(i)=g(h(i))$.\par
\noindent\textbf{Invocation:} \texttt{python3 scripts/\allowbreak{}sn\_\allowbreak{}verify.py [GAP\_\allowbreak{}BINARY] (pure Python; the GAP cross-check is skipped if GAP is absent)}.\par
\noindent\textbf{Arguments:} positional argument indices 1. \textbf{Dependencies:} Python standard library.\par
\noindent\textbf{Reads:} no literal file name. \textbf{Writes:} sn\_\allowbreak{}data.tex, sn\_\allowbreak{}max.tex.\par
\paragraph{\texttt{scripts/ss\_verify.py}.} Verification for Chapter~\ref{ch:ss}. Right-to-left composition $(gh)[x]=g[h[x]]$; points carry the left action $g(x)=g[x]$; the random seed is fixed. Writes \texttt{chapters/ss\_data.tex}, the results table used by the book, and asserts every claim. Optional argument: the path of a GAP binary.\par
\noindent\textbf{Invocation:} \texttt{python3 scripts/\allowbreak{}ss\_\allowbreak{}verify.py [GAP\_\allowbreak{}BINARY]. Right-to-left composition: (gh)[x] = g[h[x]]. Points carry the left action g(x)=g[x]. Seeds fixed (random.seed(1)).}.\par
\noindent\textbf{Arguments:} positional argument indices 1. \textbf{Dependencies:} Python standard library.\par
\noindent\textbf{Reads:} no literal file name. \textbf{Writes:} chapters/\allowbreak{}ss\_\allowbreak{}data.tex.\par

\section{Groups of Lie type}
\paragraph{\texttt{scripts/classical\_cells.py}.} Bruhat cells of the classical groups (Definition~\ref{def:stdforms}, Theorems~\ref{thm:bruhatcoord} and~\ref{thm:classicalcells}), producing the table of Example~\ref{ex:classical}. For each group it finds $U_G$ by exhaustive search, generates $G$ from $U_G$ and its images, and eliminates every element. Optional arguments: group names (default: all nine groups). Writes \texttt{drafts/classical\_cells.json}.\par
\noindent\textbf{Invocation:} \texttt{python3 classical\_\allowbreak{}cells.py [name ...] (default: all nine groups) -> writes ../\allowbreak{}drafts/\allowbreak{}classical\_\allowbreak{}cells.json}.\par
\noindent\textbf{Arguments:} positional argument indices 1:. \textbf{Dependencies:} repository modules gf; external packages numpy.\par
\noindent\textbf{Reads:} no literal file name. \textbf{Writes:} ../\allowbreak{}drafts/\allowbreak{}classical\_\allowbreak{}cells.json.\par
\paragraph{\texttt{scripts/classical\_law.py}.} Digital laws, with labels in Bruhat-cell digits, of the nine classical groups of Example~\ref{ex:classical}: $\Sp_4(2)$, $\Sp_4(3)$, $\mathrm{O}_5(3)$, $\mathrm{O}^+_4(2)$, $\mathrm{O}^+_4(3)$, $\mathrm{O}^+_4(4)$, $\mathrm{O}^+_6(2)$, $\mathrm{SU}_3(3)$ and $\mathrm{SU}_4(2)$. The group is generated and every element eliminated as in \texttt{classical\_cells.py} (Theorem~\ref{thm:bruhatcoord}).\par
\noindent\textbf{Invocation:} \texttt{cd scripts \&\& python3 classical\_\allowbreak{}law.py}.\par
\noindent\textbf{Arguments:} none. \textbf{Dependencies:} repository modules classical\_\allowbreak{}cells, gf; external packages numpy.\par
\noindent\textbf{Reads:} no literal file name. \textbf{Writes:} no literal file name.\par
\paragraph{\texttt{scripts/d43\_subproc.py}.} Worked-product adapter for the closed coordinate law of ${}^3D_4(q)$, for the registered cases $q=2,3,4$. Arguments: $q$ and an optional seed; prints the worked product as JSON.\par
\noindent\textbf{Invocation:} \texttt{python3 d43\_\allowbreak{}subproc.py q [seed] -> WORKED JSON.}.\par
\noindent\textbf{Arguments:} positional argument indices 1, 2. \textbf{Dependencies:} repository modules triality\_\allowbreak{}label.\par
\noindent\textbf{Reads:} no literal file name. \textbf{Writes:} no literal file name.\par
\paragraph{\texttt{scripts/e2e6\_law.py}.} Digital law of ${}^2E_6(q)$ inside $E_6(q^2)$, for the twist $\sigma=\mathrm{Frob}\circ\tau$. Twisted root classes are of type $A_1$, with $t\in\mathbb{F}_q$, or $A_1\times A_1$, with $t\in\mathbb{F}_{q^2}$.\par
\noindent\textbf{Invocation:} \texttt{cd scripts \&\& python3 e2e6\_\allowbreak{}law.py}.\par
\noindent\textbf{Arguments:} positional argument indices 1. \textbf{Dependencies:} repository modules e2e6\_\allowbreak{}model.\par
\noindent\textbf{Reads:} no literal file name. \textbf{Writes:} ../\allowbreak{}drafts/\allowbreak{}e2e6\_\allowbreak{}law\_\allowbreak{}q\%d.json.\par
\paragraph{\texttt{scripts/e2e6\_model.py}.} The model of ${}^2E_6$: the twisting automorphism $\tau(g)=P\,g^{-T}P^{-1}$ on the $27$-dimensional model, the twisted root classes, the group $W^\rho$, and the order formula.\par
\noindent\textbf{Invocation:} \texttt{cd scripts \&\& python3 e2e6\_\allowbreak{}model.py}.\par
\noindent\textbf{Arguments:} none. \textbf{Dependencies:} repository modules e6\_\allowbreak{}model; external packages numpy.\par
\noindent\textbf{Reads:} no literal file name. \textbf{Writes:} ../\allowbreak{}drafts/\allowbreak{}e2e6\_\allowbreak{}model.json.\par
\paragraph{\texttt{scripts/e6\_law.py}.} Digital law of $E_6(p)$, for primes $p$ with $\gcd(3,p-1)=1$, through the $27$-dimensional model and the Bruhat normal form $u\,n_w\,h\,v$.\par
\noindent\textbf{Invocation:} \texttt{cd scripts \&\& python3 e6\_\allowbreak{}law.py}.\par
\noindent\textbf{Arguments:} positional argument indices 1. \textbf{Dependencies:} repository modules e6\_\allowbreak{}model; external packages numpy.\par
\noindent\textbf{Reads:} no literal file name. \textbf{Writes:} ../\allowbreak{}drafts/\allowbreak{}e6\_\allowbreak{}law\_\allowbreak{}q\%d.json.\par
\paragraph{\texttt{scripts/e6\_lawq.py}.} Digital law of $E_6(q)$ for every prime power $q$, with $\mathbb{F}_q=\mathbb{F}_p[x]/(f)$, through the $27$-dimensional model, the Bruhat normal form and the quotient by the centre.\par
\noindent\textbf{Invocation:} \texttt{cd scripts \&\& python3 e6\_\allowbreak{}lawq.py}.\par
\noindent\textbf{Arguments:} positional argument indices 1. \textbf{Dependencies:} repository modules e6\_\allowbreak{}model; external packages numpy.\par
\noindent\textbf{Reads:} no literal file name. \textbf{Writes:} ../\allowbreak{}drafts/\allowbreak{}e6\_\allowbreak{}lawq\_\allowbreak{}q\%d.json.\par
\paragraph{\texttt{scripts/e6\_model.py}.} The root system of $E_6$, its $27$-dimensional minuscule model with integral root elements, and exact structural checks.\par
\noindent\textbf{Invocation:} \texttt{cd scripts \&\& python3 e6\_\allowbreak{}model.py}.\par
\noindent\textbf{Arguments:} none. \textbf{Dependencies:} external packages numpy.\par
\noindent\textbf{Reads:} no literal file name. \textbf{Writes:} ../\allowbreak{}drafts/\allowbreak{}e6\_\allowbreak{}model.json.\par
\paragraph{\texttt{scripts/e7\_law.py}.} $E_7(q)$: checks of the $56$-dimensional minuscule model and its digital law. The Weyl group is never enumerated; it enters through parabolic coset digits.\par
\noindent\textbf{Invocation:} \texttt{cd scripts \&\& python3 e7\_\allowbreak{}law.py q --ns --assoc --center}.\par
\noindent\textbf{Arguments:} q; --ns default=6; --assoc default=2; --center default=4. \textbf{Dependencies:} repository modules minuscule; external packages numpy.\par
\noindent\textbf{Reads:} no literal file name. \textbf{Writes:} ../\allowbreak{}drafts/\allowbreak{}e7\_\allowbreak{}q\%d.json.\par
\paragraph{\texttt{scripts/e8\_law.py}.} Digital law of $E_8(q)$ in the integral adjoint model of \texttt{e8\_model.py}, ported from \texttt{e7\_law.py}. The Weyl group enters through parabolic coset digits along $E_8>E_7>\cdots>A_1$, never by enumeration. Since $E_8$ is both simply connected and adjoint, its centre is trivial and $|T|=(q-1)^8$. Arguments: $q$ and optional \texttt{--ns} and \texttt{--assoc} (defaults $6$, $2$). Writes \texttt{drafts/e8\_q<q>.json}.\par
\noindent\textbf{Invocation:} \texttt{python3 e8\_\allowbreak{}law.py q [--ns 6] [--assoc 2] (writes ../\allowbreak{}drafts/\allowbreak{}e8\_\allowbreak{}q<q>.json)}.\par
\noindent\textbf{Arguments:} q; --ns default=6; --assoc default=2. \textbf{Dependencies:} repository modules e8\_\allowbreak{}model; external packages numpy.\par
\noindent\textbf{Reads:} no literal file name. \textbf{Writes:} ../\allowbreak{}drafts/\allowbreak{}e8\_\allowbreak{}q\%d.json.\par
\paragraph{\texttt{scripts/e8\_model.py}.} Integral $248$-dimensional adjoint model of $E_8$, from a Chevalley basis defined by a sign cocycle, with structural checks. Root elements are $x_\beta(t)=I+t\,\mathrm{ad}(e_\beta)+t^2\,\mathrm{ad}(e_\beta)^2/2$. Since $\mathrm{ad}(e_\beta)^2$ is even and $\mathrm{ad}(e_\beta)^3=0$, the model is integral in every characteristic, including $2$. Writes \texttt{drafts/e8\_model.json}.\par
\noindent\textbf{Invocation:} \texttt{characteristic, including 2. Usage: python3 e8\_\allowbreak{}model.py (writes ../\allowbreak{}drafts/\allowbreak{}e8\_\allowbreak{}model.json)}.\par
\noindent\textbf{Arguments:} none. \textbf{Dependencies:} external packages numpy.\par
\noindent\textbf{Reads:} no literal file name. \textbf{Writes:} ../\allowbreak{}drafts/\allowbreak{}e8\_\allowbreak{}model.json.\par
\paragraph{\texttt{scripts/f42\_digital.py}.} Digital law of ${}^2F_4(q)$ from its unique twisted normal form. For $q=2$, it gives the law on $A_{2p_0}$ through the Tits group.\par
\noindent\textbf{Invocation:} \texttt{cd scripts \&\& python3 f42\_\allowbreak{}digital.py}.\par
\noindent\textbf{Arguments:} none. \textbf{Dependencies:} Python standard library.\par
\noindent\textbf{Reads:} no literal file name. \textbf{Writes:} ../\allowbreak{}drafts/\allowbreak{}f42\_\allowbreak{}digital\_\allowbreak{}q\%d.json.\par
\paragraph{\texttt{scripts/f42\_fhom\_proof.py}.} Exhaustive proof that $f$ is a homomorphism ${}^2F_4(2)\to C_2$, using the generators $U$, $n_1$ and $n_2$.\par
\noindent\textbf{Invocation:} \texttt{cd scripts \&\& python3 f42\_\allowbreak{}fhom\_\allowbreak{}proof.py}.\par
\noindent\textbf{Arguments:} positional argument indices 1:. \textbf{Dependencies:} Python standard library.\par
\noindent\textbf{Reads:} no literal file name. \textbf{Writes:} ../\allowbreak{}drafts/\allowbreak{}f42\_\allowbreak{}fhom\_\allowbreak{}part\_\allowbreak{}\%s.json.\par
\paragraph{\texttt{scripts/f42\_group\_q8.py}.} Realises ${}^2F_4(8)$ inside $F_4(8)$: applies $\sigma$ to whole normal forms, finds the $\sigma$-fixed elements of the torus and $n_w$, and checks closure of the fixed elements in $W^\sigma$-cells.\par
\noindent\textbf{Invocation:} \texttt{cd scripts \&\& python3 f42\_\allowbreak{}group\_\allowbreak{}q8.py}.\par
\noindent\textbf{Arguments:} positional argument indices 1. \textbf{Dependencies:} repository modules gf.\par
\noindent\textbf{Reads:} ../\allowbreak{}drafts/\allowbreak{}f42\_\allowbreak{}param\_\allowbreak{}exact.json. \textbf{Writes:} ../\allowbreak{}drafts/\allowbreak{}f42\_\allowbreak{}group\_\allowbreak{}q8.json.\par
\paragraph{\texttt{scripts/f42\_nw\_exact.py}.} Exact check over $\mathbb{F}_2$, where $\sigma=\psi$, that $n_{w_1}=n_2n_3$ and $n_{w_2}=n_1n_4n_1n_4$ are $\sigma$-fixed in ${}^2F_4(q)$ for every $q$.\par
\noindent\textbf{Invocation:} \texttt{cd scripts \&\& python3 f42\_\allowbreak{}nw\_\allowbreak{}exact.py}.\par
\noindent\textbf{Arguments:} positional argument indices 1. \textbf{Dependencies:} Python standard library.\par
\noindent\textbf{Reads:} no literal file name. \textbf{Writes:} ../\allowbreak{}drafts/\allowbreak{}f42\_\allowbreak{}nw\_\allowbreak{}exact.json.\par
\paragraph{\texttt{scripts/f42\_order.py}.} Computes the twisted Weyl group $W^\sigma$ of ${}^2F_4$ and derives $|{}^2F_4(q)|$ from the Bruhat decomposition.\par
\noindent\textbf{Invocation:} \texttt{cd scripts \&\& python3 f42\_\allowbreak{}order.py}.\par
\noindent\textbf{Arguments:} none. \textbf{Dependencies:} repository modules f4\_\allowbreak{}law.\par
\noindent\textbf{Reads:} ../\allowbreak{}drafts/\allowbreak{}f42\_\allowbreak{}twisted\_\allowbreak{}roots.json. \textbf{Writes:} ../\allowbreak{}drafts/\allowbreak{}f42\_\allowbreak{}order.json.\par
\paragraph{\texttt{scripts/f42\_param.py}.} Fits explicit parametrisations of the eight twisted root subgroups of ${}^2F_4(q)$ over $\mathbb{F}_8$, then confirms them over $\mathbb{F}_{32}$.\par
\noindent\textbf{Invocation:} \texttt{cd scripts \&\& python3 f42\_\allowbreak{}param.py}.\par
\noindent\textbf{Arguments:} positional argument indices 1. \textbf{Dependencies:} repository modules gf.\par
\noindent\textbf{Reads:} ../\allowbreak{}drafts/\allowbreak{}f42\_\allowbreak{}twisted\_\allowbreak{}roots.json. \textbf{Writes:} ../\allowbreak{}drafts/\allowbreak{}f42\_\allowbreak{}param.json.\par
\paragraph{\texttt{scripts/f42\_param\_exact.py}.} Exact proof that the eight twisted root parametrisations of ${}^2F_4(q)$ are $\sigma$-fixed for every $q=2^{2m+1}$. It works over $\mathbb{F}_2$ with monomial exponents $a+b\theta$, stored as doubled integer slots, and uses $\theta^2=q/2=1/2\pmod {q-1}$.\par
\noindent\textbf{Invocation:} \texttt{cd scripts \&\& python3 f42\_\allowbreak{}param\_\allowbreak{}exact.py}.\par
\noindent\textbf{Arguments:} positional argument indices 1. \textbf{Dependencies:} Python standard library.\par
\noindent\textbf{Reads:} ../\allowbreak{}drafts/\allowbreak{}f42\_\allowbreak{}param.json. \textbf{Writes:} ../\allowbreak{}drafts/\allowbreak{}f42\_\allowbreak{}param\_\allowbreak{}exact.json.\par
\paragraph{\texttt{scripts/f42\_psi.py}.} First step for ${}^2F_4$: the root map $\rho$ and special isogeny $\psi$ of $F_4$ in characteristic $2$. The map $\rho$ swaps long and short roots, and $\psi(x_a(t))=x_{\rho a}(t^{\lambda(a)})$, with $\lambda(a)=1$ or $2$ according to root length. It proves exactly over $\mathbb{F}_2$-algebras that $\psi$ is an endomorphism of $U$.\par
\noindent\textbf{Invocation:} \texttt{cd scripts \&\& python3 f42\_\allowbreak{}psi.py}.\par
\noindent\textbf{Arguments:} none. \textbf{Dependencies:} repository modules f4\_\allowbreak{}law.\par
\noindent\textbf{Reads:} no literal file name. \textbf{Writes:} ../\allowbreak{}drafts/\allowbreak{}f42\_\allowbreak{}psi.json.\par
\paragraph{\texttt{scripts/f42\_psi\_group.py}.} Exact verification over $\mathbb{F}_2[\mathrm{vars}]$ that $\psi$ respects the full $F_4$ law. It checks the root, Weyl, torus and $U$ actions of $\psi$, together with $\psi(ab)=\psi(a)\psi(b)$ on every branch for each fixed word and simple root. Arguments: a part number and number of parts.\par
\noindent\textbf{Invocation:} \texttt{python3 f42\_\allowbreak{}psi\_\allowbreak{}group.py <part> <parts>}.\par
\noindent\textbf{Arguments:} positional argument indices 1, 2. \textbf{Dependencies:} repository modules f4\_\allowbreak{}jordanq, f4\_\allowbreak{}law.\par
\noindent\textbf{Reads:} no literal file name. \textbf{Writes:} no literal file name.\par
\paragraph{\texttt{scripts/f42\_psi\_n.py}.} Exact check over $\mathbb{F}_2[\mathrm{vars}]$ that $\psi(f n_s)=\psi(f)n_{\rho s}$, applying $\psi$ directly to normal forms. Arguments: a part number and number of parts.\par
\noindent\textbf{Invocation:} \texttt{python3 f42\_\allowbreak{}psi\_\allowbreak{}n.py <part> <parts>}.\par
\noindent\textbf{Arguments:} positional argument indices 1, 2. \textbf{Dependencies:} Python standard library.\par
\noindent\textbf{Reads:} no literal file name. \textbf{Writes:} no literal file name.\par
\paragraph{\texttt{scripts/f42\_rank1.py}.} Rank-one rules of ${}^2F_4(q)$: \[ n_wX(p)n_w^{-1}=X(p^\prime)\,n_w\,h(a,b)\,X(p^{\prime\prime}) \] for $w_1=s_2s_3$ and $w_2=(s_1s_4)^2$.\par
\noindent\textbf{Invocation:} \texttt{cd scripts \&\& python3 f42\_\allowbreak{}rank1.py}.\par
\noindent\textbf{Arguments:} none. \textbf{Dependencies:} Python standard library.\par
\noindent\textbf{Reads:} no literal file name. \textbf{Writes:} ../\allowbreak{}drafts/\allowbreak{}f42\_\allowbreak{}rank1.json.\par
\paragraph{\texttt{scripts/f42\_rank1\_proof.py}.} Proof of the rank-one rules for every $q$, in characteristic $2$, using the Levi subgroups $\mathrm{Sp}_4$ on $(\alpha_1,\alpha_4)$ and $\SL_2\times\SL_2$ on $(\alpha_2,\alpha_3)$.\par
\noindent\textbf{Invocation:} \texttt{cd scripts \&\& python3 f42\_\allowbreak{}rank1\_\allowbreak{}proof.py}.\par
\noindent\textbf{Arguments:} none. \textbf{Dependencies:} repository modules gf.\par
\noindent\textbf{Reads:} no literal file name. \textbf{Writes:} ../\allowbreak{}drafts/\allowbreak{}f42\_\allowbreak{}rank1\_\allowbreak{}proof.json.\par
\paragraph{\texttt{scripts/f42\_rank1\_tprime.py}.} Fits $t^\prime$ in the rank-one relation \[ n_{w_2}X(t,u)n_{w_2}^{-1}=X(t^\prime,u^\prime)\,n_{w_2}h\,X(t^{\prime\prime},u^{\prime\prime}) \] by two-dimensional interpolation of $t^\prime N^m$.\par
\noindent\textbf{Invocation:} \texttt{cd scripts \&\& python3 f42\_\allowbreak{}rank1\_\allowbreak{}tprime.py}.\par
\noindent\textbf{Arguments:} none. \textbf{Dependencies:} repository modules gf.\par
\noindent\textbf{Reads:} ../\allowbreak{}drafts/\allowbreak{}f42\_\allowbreak{}rank1.json. \textbf{Writes:} ../\allowbreak{}drafts/\allowbreak{}f42\_\allowbreak{}rank1\_\allowbreak{}tprime.json.\par
\paragraph{\texttt{scripts/f42\_subproc.py}.} Worked-product adapter for ${}^2F_4(q)$, $q=2^{2a+1}$, using the labels and unique-twisted-normal-form law of \texttt{f42\_digital.py}. For a seeded pair it checks round-trip encoding, label range, the order formula \[ q^{12}(q^6+1)(q^4-1)(q^3+1)(q-1), \] and the identity. Arguments: $q$ and an optional seed; prints the worked product as JSON.\par
\noindent\textbf{Invocation:} \texttt{python3 f42\_\allowbreak{}subproc.py q [seed] -> prints one line 'WORKED \{json\}'.}.\par
\noindent\textbf{Arguments:} positional argument indices 1, 2. \textbf{Dependencies:} Python standard library.\par
\noindent\textbf{Reads:} no literal file name. \textbf{Writes:} no literal file name.\par
\paragraph{\texttt{scripts/f42\_tits.py}.} Studies the homomorphism ${}^2F_4(2)\to C_2$, whose kernel is the Tits group ${}^2F_4(2)^\prime$, and determines whether the extension splits.\par
\noindent\textbf{Invocation:} \texttt{cd scripts \&\& python3 f42\_\allowbreak{}tits.py}.\par
\noindent\textbf{Arguments:} none. \textbf{Dependencies:} Python standard library.\par
\noindent\textbf{Reads:} no literal file name. \textbf{Writes:} ../\allowbreak{}drafts/\allowbreak{}f42\_\allowbreak{}tits.json.\par
\paragraph{\texttt{scripts/f42\_twisted\_law.py}.} Multiplication in ${}^2F_4(q)$ in twisted coordinates $(u,w,t,v)$, validated against $F_4$ over $\mathbb{F}_8$.\par
\noindent\textbf{Invocation:} \texttt{cd scripts \&\& python3 f42\_\allowbreak{}twisted\_\allowbreak{}law.py}.\par
\noindent\textbf{Arguments:} positional argument indices 1, 2. \textbf{Dependencies:} Python standard library.\par
\noindent\textbf{Reads:} ../\allowbreak{}drafts/\allowbreak{}f42\_\allowbreak{}ulaw\_\allowbreak{}polys.json. \textbf{Writes:} ../\allowbreak{}drafts/\allowbreak{}f42\_\allowbreak{}twisted\_\allowbreak{}law\_\allowbreak{}q\%d.json.\par
\paragraph{\texttt{scripts/f42\_ulaw.py}.} Closed law on the $12$ coordinates of $U^\sigma$ for ${}^2F_4(q)$, exact for every $q=2^{2m+1}$.\par
\noindent\textbf{Invocation:} \texttt{cd scripts \&\& python3 f42\_\allowbreak{}ulaw.py}.\par
\noindent\textbf{Arguments:} positional argument indices 1. \textbf{Dependencies:} repository modules gf.\par
\noindent\textbf{Reads:} ../\allowbreak{}drafts/\allowbreak{}f42\_\allowbreak{}param\_\allowbreak{}exact.json. \textbf{Writes:} ../\allowbreak{}drafts/\allowbreak{}f42\_\allowbreak{}ulaw.json, ../\allowbreak{}drafts/\allowbreak{}f42\_\allowbreak{}ulaw\_\allowbreak{}polys.json.\par
\paragraph{\texttt{scripts/f42\_usigma\_q2.py}.} Computes $U^\sigma$ for ${}^2F_4(2)$, where $\sigma=\psi$ over $\mathbb{F}_2$: the fixed points of each twisted root class and the group that they generate.\par
\noindent\textbf{Invocation:} \texttt{cd scripts \&\& python3 f42\_\allowbreak{}usigma\_\allowbreak{}q2.py}.\par
\noindent\textbf{Arguments:} positional argument indices 1. \textbf{Dependencies:} Python standard library.\par
\noindent\textbf{Reads:} ../\allowbreak{}drafts/\allowbreak{}f42\_\allowbreak{}twisted\_\allowbreak{}roots.json. \textbf{Writes:} ../\allowbreak{}drafts/\allowbreak{}f42\_\allowbreak{}usigma\_\allowbreak{}q2.json.\par
\paragraph{\texttt{scripts/f42\_usigma\_q8.py}.} Computes $U^\sigma$ for ${}^2F_4(8)$, using $\sigma(x_a(t))=x_{\rho a}(t^{\theta k(a)})$ with $\theta=2$ over $\mathbb{F}_8$.\par
\noindent\textbf{Invocation:} \texttt{cd scripts \&\& python3 f42\_\allowbreak{}usigma\_\allowbreak{}q8.py}.\par
\noindent\textbf{Arguments:} positional argument indices 1. \textbf{Dependencies:} repository modules gf.\par
\noindent\textbf{Reads:} ../\allowbreak{}drafts/\allowbreak{}f42\_\allowbreak{}twisted\_\allowbreak{}roots.json. \textbf{Writes:} ../\allowbreak{}drafts/\allowbreak{}f42\_\allowbreak{}usigma\_\allowbreak{}q8.json.\par
\paragraph{\texttt{scripts/f4\_adjoint.py}.} Integral $52$-dimensional adjoint model of $F_4$ and resumable checks of its Weyl cells.\par
\noindent\textbf{Invocation:} \texttt{cd scripts \&\& python3 f4\_\allowbreak{}adjoint.py}.\par
\noindent\textbf{Arguments:} none. \textbf{Dependencies:} repository modules f4\_\allowbreak{}law; external packages numpy.\par
\noindent\textbf{Reads:} no literal file name. \textbf{Writes:} no literal file name.\par
\paragraph{\texttt{scripts/f4\_adjointq.py}.} The $52$-dimensional adjoint model of $F_4(q)$ over $\mathbb{F}_q$ and the Weyl-cell sweep for \texttt{f4\_lawq.py} (Proposition~\ref{prop:f4adjoint}). It uses $X_a(t)=I+t\,\operatorname{ad}(e_a)+t^2\,\operatorname{ad}(e_a)^2/2$, reducing the integer matrices modulo $p$. For $\mathbb{F}_p$, the random streams and statistics agree with \texttt{f4\_adjoint.verify\_cells}. Arguments: $p$, $k$, width. Writes \texttt{f4\_cells\_q<q>\_full.json}.\par
\noindent\textbf{Invocation:} \texttt{python3 f4\_\allowbreak{}adjointq.py p k width (writes f4\_\allowbreak{}cells\_\allowbreak{}q<q>\_\allowbreak{}full.json)}.\par
\noindent\textbf{Arguments:} positional argument indices 1:4. \textbf{Dependencies:} repository modules f4\_\allowbreak{}adjoint, f4\_\allowbreak{}law, f4\_\allowbreak{}lawq, gf; external packages numpy.\par
\noindent\textbf{Reads:} no literal file name. \textbf{Writes:} no literal file name.\par
\paragraph{\texttt{scripts/f4\_core.py}.} Core construction for $F_4$ in the Albert model: derivations modulo $p$, root spaces, and a Chevalley basis.\par
\noindent\textbf{Invocation:} \texttt{cd scripts \&\& python3 f4\_\allowbreak{}core.py}.\par
\noindent\textbf{Arguments:} none. \textbf{Dependencies:} external packages numpy.\par
\noindent\textbf{Reads:} no literal file name. \textbf{Writes:} no literal file name.\par
\paragraph{\texttt{scripts/f4\_jordan.py}.} Zorn-octonion and Albert-algebra construction used for the Jordan model of $F_4$.\par
\noindent\textbf{Invocation:} \texttt{cd scripts \&\& python3 f4\_\allowbreak{}jordan.py}.\par
\noindent\textbf{Arguments:} none. \textbf{Dependencies:} repository modules f4\_\allowbreak{}core; external packages numpy.\par
\noindent\textbf{Reads:} no literal file name. \textbf{Writes:} no literal file name.\par
\paragraph{\texttt{scripts/f4\_jordanq.py}.} Integral $27$-dimensional Albert model of $F_4$ over $\mathbb{F}_q$, with checks (J6) and (J7). Commands include \texttt{signs}, \texttt{burnside}, \texttt{sweep p k width}, and \texttt{smoke p k}.\par
\noindent\textbf{Invocation:} \texttt{python3 f4\_\allowbreak{}jordanq.py signs | burnside | sweep p k width | smoke p k}.\par
\noindent\textbf{Arguments:} positional argument indices 1, 2, 2:5, 3. \textbf{Dependencies:} repository modules f4\_\allowbreak{}adjoint, f4\_\allowbreak{}adjointq, f4\_\allowbreak{}jordan, f4\_\allowbreak{}law, f4\_\allowbreak{}lawq, gf; external packages numpy.\par
\noindent\textbf{Reads:} no literal file name. \textbf{Writes:} no literal file name.\par
\paragraph{\texttt{scripts/f4\_law.py}.} Coordinate multiplication for split $F_4$, built from the persistent core.\par
\noindent\textbf{Invocation:} \texttt{cd scripts \&\& python3 f4\_\allowbreak{}law.py}.\par
\noindent\textbf{Arguments:} none. \textbf{Dependencies:} repository modules f4\_\allowbreak{}core; external packages numpy.\par
\noindent\textbf{Reads:} no literal file name. \textbf{Writes:} no literal file name.\par
\paragraph{\texttt{scripts/f4\_lawq.py}.} Coordinate multiplication of $F_4(q)$ over an arbitrary finite field. It implements the same law as \texttt{f4\_law.build\_law} (Theorem~\ref{thm:f4closed}), with every coefficient operation delegated to a \texttt{gf.GF} object. The field-independent commutator constants, Weyl signs and fixed words come from \texttt{f4\_law.py}; over $\mathbb{F}_p$ the result agrees with the law over $\mathbb{Z}/p\mathbb{Z}$.\par
\noindent\textbf{Invocation:} \texttt{cd scripts \&\& python3 f4\_\allowbreak{}lawq.py}.\par
\noindent\textbf{Arguments:} none. \textbf{Dependencies:} repository modules f4\_\allowbreak{}law.\par
\noindent\textbf{Reads:} no literal file name. \textbf{Writes:} no literal file name.\par
\paragraph{\texttt{scripts/f4\_subproc.py}.} Worked-product adapter for $F_4(p)$, $p$ prime, using the labels of Corollary~\ref{cor:f4law} and \texttt{f4\_lawq.multiply}. It uses the same cell ordering, little-endian base-$p$ digits and least primitive root $\gamma$ as \texttt{worked\_f4.py}; for $p=2$, $\gamma=1$. For a seeded pair of labels it checks round-trip encoding, label range, the order formula, and the identity. Arguments: $p$ and an optional seed; prints the worked product as JSON.\par
\noindent\textbf{Invocation:} \texttt{python3 f4\_\allowbreak{}subproc.py p [seed] -> prints one line 'WORKED \{json\}'.}.\par
\noindent\textbf{Arguments:} positional argument indices 1, 2. \textbf{Dependencies:} repository modules f4\_\allowbreak{}adjoint, f4\_\allowbreak{}adjointq, f4\_\allowbreak{}lawq, gf; external packages numpy.\par
\noindent\textbf{Reads:} no literal file name. \textbf{Writes:} no literal file name.\par
\paragraph{\texttt{scripts/f4\_symbolic\_all.py}.} Symbolic run of \texttt{mul\_n} and \texttt{convert} on every cell, in characteristic $0$, using the symbolic field of \texttt{f4\_symbolic\_probe.py}. Checks that no non-monomial inverse occurs, records output sizes, and compares the adaptor at a random point over $\mathbb{F}_{101}$; this is not a model comparison.\par
\noindent\textbf{Invocation:} \texttt{cd scripts \&\& python3 f4\_\allowbreak{}symbolic\_\allowbreak{}all.py}.\par
\noindent\textbf{Arguments:} none. \textbf{Dependencies:} Python standard library.\par
\noindent\textbf{Reads:} no literal file name. \textbf{Writes:} ../\allowbreak{}drafts/\allowbreak{}f4\_\allowbreak{}symbolic\_\allowbreak{}all.json.\par
\paragraph{\texttt{scripts/f4\_symbolic\_integrality.py}.} Integrality check for the symbolic $F_4$ law: every inverse is a monomial with coefficient $\pm1$, and every coefficient output by \texttt{mul\_n} and \texttt{convert}, over all cells, is an integer.\par
\noindent\textbf{Invocation:} \texttt{cd scripts \&\& python3 f4\_\allowbreak{}symbolic\_\allowbreak{}integrality.py}.\par
\noindent\textbf{Arguments:} none. \textbf{Dependencies:} repository modules f4\_\allowbreak{}jordanq.\par
\noindent\textbf{Reads:} no literal file name. \textbf{Writes:} ../\allowbreak{}drafts/\allowbreak{}f4\_\allowbreak{}symbolic\_\allowbreak{}integrality.json.\par
\paragraph{\texttt{scripts/f4\_symbolic\_model.py}.} Exact symbolic check of the $F_4$ law against the $27$-dimensional model in characteristic $0$. The unchanged \texttt{f4\_lawq.ops} acts on Laurent polynomials in the digits, torus and $U$-coordinates; the model acts on a generic vector $y$. For every cell it checks the actions of \texttt{mul\_n}, \texttt{convert} and \texttt{mul\_root} against the corresponding matrix actions.\par
\noindent\textbf{Invocation:} \texttt{cd scripts \&\& python3 f4\_\allowbreak{}symbolic\_\allowbreak{}model.py}.\par
\noindent\textbf{Arguments:} none. \textbf{Dependencies:} repository modules f4\_\allowbreak{}jordanq, f4\_\allowbreak{}law, f4\_\allowbreak{}lawq; external packages numpy.\par
\noindent\textbf{Reads:} no literal file name. \textbf{Writes:} ../\allowbreak{}drafts/\allowbreak{}f4\_\allowbreak{}symbolic\_\allowbreak{}model\_\allowbreak{}\{TAG\}.json.\par
\paragraph{\texttt{scripts/f4\_symbolic\_probe.py}.} Feasibility probe for the $q$-generic $F_4$ law: runs the unchanged \texttt{f4\_lawq.ops} over a symbolic field of Laurent polynomials over $\mathbb{Q}$ in digit, torus and $U$-coordinate variables. Only monomials are inverted. Measures the expression growth of \texttt{mul\_n} and checks its adaptor at a random point over $\mathbb{F}_{101}$; it is not a model comparison.\par
\noindent\textbf{Invocation:} \texttt{cd scripts \&\& python3 f4\_\allowbreak{}symbolic\_\allowbreak{}probe.py}.\par
\noindent\textbf{Arguments:} none. \textbf{Dependencies:} repository modules f4\_\allowbreak{}law, f4\_\allowbreak{}lawq, gf.\par
\noindent\textbf{Reads:} no literal file name. \textbf{Writes:} ../\allowbreak{}drafts/\allowbreak{}f4\_\allowbreak{}symbolic\_\allowbreak{}probe.txt.\par
\paragraph{\texttt{scripts/f4\_symbolic\_tail.py}.} Exact symbolic check, in the $27$-dimensional model, of the final step of multiplication: right multiplication by $hv\in B$.\par
\noindent\textbf{Invocation:} \texttt{cd scripts \&\& python3 f4\_\allowbreak{}symbolic\_\allowbreak{}tail.py}.\par
\noindent\textbf{Arguments:} none. \textbf{Dependencies:} Python standard library.\par
\noindent\textbf{Reads:} no literal file name. \textbf{Writes:} ../\allowbreak{}drafts/\allowbreak{}f4\_\allowbreak{}symbolic\_\allowbreak{}tail.json.\par
\paragraph{\texttt{scripts/f4\_ucert.py}.} Generates the unipotent certificate \texttt{drafts/f4\_ucert.json} from the integral $27$-dimensional model of $F_4$.\par
\noindent\textbf{Invocation:} \texttt{cd scripts \&\& python3 f4\_\allowbreak{}ucert.py}.\par
\noindent\textbf{Arguments:} none. \textbf{Dependencies:} repository modules f4\_\allowbreak{}core, f4\_\allowbreak{}jordan, f4\_\allowbreak{}law, f4\_\allowbreak{}ucert\_\allowbreak{}verify; external packages numpy.\par
\noindent\textbf{Reads:} no literal file name. \textbf{Writes:} ../\allowbreak{}drafts/\allowbreak{}f4\_\allowbreak{}ucert.json.\par
\paragraph{\texttt{scripts/f4\_ucert\_verify.py}.} Independent verifier for the unipotent certificate \texttt{drafts/f4\_ucert.json}. Reads only that JSON file and uses exact integer-polynomial arithmetic, importing nothing from the project. Polynomials are dictionaries from packed monomials to integers, with six-bit exponent fields; every exponent is checked to be less than $32$.\par
\noindent\textbf{Invocation:} \texttt{cd scripts \&\& python3 f4\_\allowbreak{}ucert\_\allowbreak{}verify.py}.\par
\noindent\textbf{Arguments:} positional argument indices 1. \textbf{Dependencies:} Python standard library.\par
\noindent\textbf{Reads:} no literal file name. \textbf{Writes:} no literal file name.\par
\paragraph{\texttt{scripts/f4\_wcert.py}.} Generates \texttt{drafts/f4\_wcert.json}, the certificate of Weyl signs, torus weights and rank-one data, from the integral $27$-dimensional model of $F_4$.\par
\noindent\textbf{Invocation:} \texttt{cd scripts \&\& python3 f4\_\allowbreak{}wcert.py}.\par
\noindent\textbf{Arguments:} none. \textbf{Dependencies:} repository modules f4\_\allowbreak{}core, f4\_\allowbreak{}law, f4\_\allowbreak{}ucert, f4\_\allowbreak{}wcert\_\allowbreak{}verify; external packages numpy.\par
\noindent\textbf{Reads:} no literal file name. \textbf{Writes:} ../\allowbreak{}drafts/\allowbreak{}f4\_\allowbreak{}wcert.json.\par
\paragraph{\texttt{scripts/f4\_wcert\_verify.py}.} Independent verifier for the Weyl--torus certificate \texttt{drafts/f4\_wcert.json}. Reads only its JSON data and uses exact integer and Laurent-polynomial arithmetic, importing nothing from the project.\par
\noindent\textbf{Invocation:} \texttt{cd scripts \&\& python3 f4\_\allowbreak{}wcert\_\allowbreak{}verify.py}.\par
\noindent\textbf{Arguments:} positional argument indices 1. \textbf{Dependencies:} Python standard library.\par
\noindent\textbf{Reads:} no literal file name. \textbf{Writes:} no literal file name.\par
\paragraph{\texttt{scripts/f4\_wordcert.py}.} Generates \texttt{drafts/f4\_wordcert.json}, containing fixed words, roots $\beta_k$ and signs $\epsilon_k$ from the certified conjugation signs.\par
\noindent\textbf{Invocation:} \texttt{cd scripts \&\& python3 f4\_\allowbreak{}wordcert.py}.\par
\noindent\textbf{Arguments:} none. \textbf{Dependencies:} repository modules f4\_\allowbreak{}law.\par
\noindent\textbf{Reads:} ../\allowbreak{}drafts/\allowbreak{}f4\_\allowbreak{}wcert.json. \textbf{Writes:} ../\allowbreak{}drafts/\allowbreak{}f4\_\allowbreak{}wordcert.json.\par
\paragraph{\texttt{scripts/f4\_wordcert\_verify.py}.} Independently checks, from \texttt{drafts/f4\_wcert.json} and \texttt{drafts/f4\_wordcert.json}, the rewriting formula \[ P_W(t)n_W^{-1}=\prod_k x_{\beta_k}(\epsilon_k t_k) \] for every fixed word.\par
\noindent\textbf{Invocation:} \texttt{cd scripts \&\& python3 f4\_\allowbreak{}wordcert\_\allowbreak{}verify.py}.\par
\noindent\textbf{Arguments:} none. \textbf{Dependencies:} external packages numpy.\par
\noindent\textbf{Reads:} ../\allowbreak{}drafts/\allowbreak{}f4\_\allowbreak{}wcert.json, ../\allowbreak{}drafts/\allowbreak{}f4\_\allowbreak{}wordcert.json. \textbf{Writes:} no literal file name.\par
\paragraph{\texttt{scripts/f4q\_adjtest.py}.} Tests for the adjoint construction of $F_4(q)$.\par
\noindent\textbf{Invocation:} \texttt{cd scripts \&\& python3 f4q\_\allowbreak{}adjtest.py}.\par
\noindent\textbf{Arguments:} none. \textbf{Dependencies:} repository modules f4\_\allowbreak{}adjoint, f4\_\allowbreak{}adjointq, f4\_\allowbreak{}lawq, gf.\par
\noindent\textbf{Reads:} no literal file name. \textbf{Writes:} no literal file name.\par
\paragraph{\texttt{scripts/f4q\_test.py}.} Tests for the finite-field implementation of the $F_4(q)$ law.\par
\noindent\textbf{Invocation:} \texttt{cd scripts \&\& python3 f4q\_\allowbreak{}test.py}.\par
\noindent\textbf{Arguments:} none. \textbf{Dependencies:} repository modules f4\_\allowbreak{}law, f4\_\allowbreak{}lawq, gf.\par
\noindent\textbf{Reads:} no literal file name. \textbf{Writes:} no literal file name.\par
\paragraph{\texttt{scripts/g2\_cells.py}.} Checks for $G_2(p)$, $p\ge5$ prime: Definition~\ref{def:zorn}, Lemma~\ref{lem:g2roots}, Theorem~\ref{thm:g2cells} and Example~\ref{ex:g2}. It verifies that the Zorn algebra is alternative, that its derivation algebra has dimension $14$ with $12$ one-dimensional root spaces, and that the root elements $\exp(cD)$ are upper unitriangular automorphisms. Argument: $p$. Writes \texttt{drafts/g2\_cells\_q<p>.json}.\par
\noindent\textbf{Invocation:} \texttt{python3 g2\_\allowbreak{}cells.py p -> writes ../\allowbreak{}drafts/\allowbreak{}g2\_\allowbreak{}cells\_\allowbreak{}q<p>.json}.\par
\noindent\textbf{Arguments:} positional argument indices 1. \textbf{Dependencies:} external packages numpy.\par
\noindent\textbf{Reads:} no literal file name. \textbf{Writes:} ../\allowbreak{}drafts/\allowbreak{}g2\_\allowbreak{}cells\_\allowbreak{}q\%d.json.\par
\paragraph{\texttt{scripts/g2\_extract\_constants.py}.} Extracts the Chevalley constants of $G_2(p)$ in the root normalization of \texttt{g2\_cells.py} and \texttt{g2\_law.py}; each convention is verified on every parameter pair in $(\mathbb{F}_p^\times)^2$. Optional arguments: primes (default $5$, $7$, $11$, $13$). Writes \texttt{drafts/g2\_constants.json}.\par
\noindent\textbf{Invocation:} \texttt{python3 g2\_\allowbreak{}extract\_\allowbreak{}constants.py [p ...] (default 5 7 11 13)}.\par
\noindent\textbf{Arguments:} none. \textbf{Dependencies:} repository modules g2\_\allowbreak{}law; external packages numpy.\par
\noindent\textbf{Reads:} no literal file name. \textbf{Writes:} ../\allowbreak{}drafts/\allowbreak{}g2\_\allowbreak{}constants.json.\par
\paragraph{\texttt{scripts/g2\_law.py}.} Digital law of $G_2(p)$, $p\ge5$ prime, in the $7$-dimensional Zorn-algebra model of \texttt{g2\_cells.py} (Definition~\ref{def:zorn}, Lemma~\ref{lem:g2roots}, Theorem~\ref{thm:g2cells}). An element is labelled through its Bruhat form $g=u\,n_w\,t\,u^\prime$ with $u\in U_w$, where $n_w$ is a product of the $n_s=x_s(1)\,x_{-s}(c)\,x_s(1)$ over simple roots $s$. Arguments: $p$ and an optional seed. Prints a worked product.\par
\noindent\textbf{Invocation:} \texttt{python3 g2\_\allowbreak{}law.py p [seed] -> prints a worked product and one line 'WORKED \{json\}'.}.\par
\noindent\textbf{Arguments:} positional argument indices 1, 2. \textbf{Dependencies:} repository modules g2\_\allowbreak{}symbolic; external packages numpy.\par
\noindent\textbf{Reads:} no literal file name. \textbf{Writes:} no literal file name.\par
\paragraph{\texttt{scripts/g2\_rescale.py}.} Rescales the root $2a+b$ of the $p$-independent $G_2$ tables, $x^\prime_{2a+b}(s)=x_{2a+b}(2s)$, leaving the negative root unchanged, and checks that the result has the shape of a Chevalley basis: integer commutator constants, Weyl signs $\pm1$, $\kappa=-1$ and rank-one $e=1$. Writes \texttt{g2\_tables\_chevalley.py}; the labels and \texttt{g2\_law.py} are unchanged.\par
\noindent\textbf{Invocation:} \texttt{cd scripts \&\& python3 g2\_\allowbreak{}rescale.py}.\par
\noindent\textbf{Arguments:} none. \textbf{Dependencies:} repository modules g2\_\allowbreak{}law, g2\_\allowbreak{}symbolic, g2\_\allowbreak{}tables.\par
\noindent\textbf{Reads:} ../\allowbreak{}drafts/\allowbreak{}g2\_\allowbreak{}constants.json. \textbf{Writes:} no literal file name.\par
\paragraph{\texttt{scripts/g2\_sign\_check.py}.} Checks the signs of the rescaled $G_2$ tables against a published normalization of the adjoint Chevalley groups of types $A_1$, $A_2$, $B_2$ and $G_2$, relations (6.1)--(6.5). There $\alpha$ is long and $\beta$ short; here $a$ is short and $b$ long.\par
\noindent\textbf{Invocation:} \texttt{cd scripts \&\& python3 g2\_\allowbreak{}sign\_\allowbreak{}check.py}.\par
\noindent\textbf{Arguments:} none. \textbf{Dependencies:} repository modules g2\_\allowbreak{}symbolic, g2\_\allowbreak{}tables\_\allowbreak{}chevalley.\par
\noindent\textbf{Reads:} no literal file name. \textbf{Writes:} no literal file name.\par
\paragraph{\texttt{scripts/g2\_symbolic.py}.} Matrix-free normal-form multiplication in $G_2(p)$, in the normalization of \texttt{g2\_law.py}. It uses tables of commutator relations $[x_r(t),x_s(u)]=\prod_k x_k(C\,t^iu^j)$, Weyl conjugations $n_s^{-1}x_k(c)\,n_s=x_{k^\prime}(\eta c)$, the conjugation of the torus by $n_s$, and the squares $n_s^2$.\par
\noindent\textbf{Invocation:} \texttt{cd scripts \&\& python3 g2\_\allowbreak{}symbolic.py}.\par
\noindent\textbf{Arguments:} none. \textbf{Dependencies:} repository modules g2\_\allowbreak{}tables; external packages numpy.\par
\noindent\textbf{Reads:} no literal file name. \textbf{Writes:} ], tuple(du.get(k, 0) for k in PH[st[.\par
\paragraph{\texttt{scripts/g2\_symbolic\_tables.py}.} Table data for the symbolic $G_2$ multiplication of \texttt{g2\_symbolic.py}.\par
\noindent\textbf{Invocation:} \texttt{cd scripts \&\& python3 g2\_\allowbreak{}symbolic\_\allowbreak{}tables.py}.\par
\noindent\textbf{Arguments:} none. \textbf{Dependencies:} Python standard library.\par
\noindent\textbf{Reads:} no literal file name. \textbf{Writes:} no literal file name.\par
\paragraph{\texttt{scripts/g2\_tables.py}.} The $p$-independent $G_2$ rewrite tables used by \texttt{g2\_symbolic.py}, generated by \texttt{g2\_tables\_gen.py} from $p=5,7,11,13$; not to be edited by hand.\par
\noindent\textbf{Invocation:} \texttt{cd scripts \&\& python3 g2\_\allowbreak{}tables.py}.\par
\noindent\textbf{Arguments:} none. \textbf{Dependencies:} Python standard library.\par
\noindent\textbf{Reads:} no literal file name. \textbf{Writes:} no literal file name.\par
\paragraph{\texttt{scripts/g2\_tables\_chevalley.py}.} The $G_2$ tables with the root $2a+b$ rescaled, generated by \texttt{g2\_rescale.py}, in the format of \texttt{g2\_tables.py}. Not used by \texttt{g2\_law.py}, so the labels are unchanged.\par
\noindent\textbf{Invocation:} \texttt{cd scripts \&\& python3 g2\_\allowbreak{}tables\_\allowbreak{}chevalley.py}.\par
\noindent\textbf{Arguments:} none. \textbf{Dependencies:} Python standard library.\par
\noindent\textbf{Reads:} no literal file name. \textbf{Writes:} no literal file name.\par
\paragraph{\texttt{scripts/g2\_tables\_gen.py}.} Fits $p$-independent rewrite tables for $G_2$ from the matrix model at several primes and writes \texttt{g2\_tables.py}. Constants become exact small rationals, torus exponents integers, and big-cell entries monomials $c\,e^m$. Optional arguments: primes (default $5$, $7$, $11$, $13$).\par
\noindent\textbf{Invocation:} \texttt{python3 g2\_\allowbreak{}tables\_\allowbreak{}gen.py [p ...] (default 5 7 11 13)}.\par
\noindent\textbf{Arguments:} positional argument indices 1:. \textbf{Dependencies:} repository modules g2\_\allowbreak{}law, g2\_\allowbreak{}symbolic.\par
\noindent\textbf{Reads:} no literal file name. \textbf{Writes:} no literal file name.\par
\paragraph{\texttt{scripts/lie\_subproc.py}.} Adapter that obtains worked products from the script-style laws \texttt{e6\_law.py}, \texttt{e6\_lawq.py}, \texttt{e7\_law.py}, \texttt{e8\_law.py} and \texttt{e2e6\_law.py}. Those modules set themselves up from their command-line arguments and then run their own tests. The adapter executes each source only up to the definition of its decoder, so the setup runs but the built-in tests do not. Arguments: the module, $q$ and an optional seed; prints the worked product as JSON.\par
\noindent\textbf{Invocation:} \texttt{python3 lie\_\allowbreak{}subproc.py MODULE q [seed] -> prints one line 'WORKED \{json\}'.}.\par
\noindent\textbf{Arguments:} positional argument indices 1, 2, 3. \textbf{Dependencies:} external packages numpy.\par
\noindent\textbf{Reads:} no literal file name. \textbf{Writes:} no literal file name.\par
\paragraph{\texttt{scripts/minuscule.py}.} Integral minuscule model for a simply-laced Dynkin diagram: roots, weights, and root vectors normalized by $[e_\beta,e_{-\beta}]=h_\beta$.\par
\noindent\textbf{Invocation:} \texttt{cd scripts \&\& python3 minuscule.py}.\par
\noindent\textbf{Arguments:} none. \textbf{Dependencies:} external packages numpy.\par
\noindent\textbf{Reads:} no literal file name. \textbf{Writes:} no literal file name.\par
\paragraph{\texttt{scripts/omega\_plus\_even.py}.} The groups $\Omega^+_{2n}(q)$ with $q$ even, for the quadratic form $Q(x)=\sum_{i<n}x_i\,x_{r(i)}$: Bruhat cells and the data of the classical digital law.\par
\noindent\textbf{Invocation:} \texttt{cd scripts \&\& python3 omega\_\allowbreak{}plus\_\allowbreak{}even.py}.\par
\noindent\textbf{Arguments:} none. \textbf{Dependencies:} Python standard library.\par
\noindent\textbf{Reads:} no literal file name. \textbf{Writes:} ../\allowbreak{}drafts/\allowbreak{}omega\_\allowbreak{}plus\_\allowbreak{}even.json.\par
\paragraph{\texttt{scripts/psln\_dense.py}.} Digital laws of $\PSL_n(q)$: the canonical-representative law (Theorem~\ref{thm:psln}) and the dense Bruhat label law (Theorem~\ref{thm:pslndense}), as in Examples~\ref{ex:psl3} and~\ref{ex:pslndense}. Arguments: one or more triples $n$, $p$, $k$, with $q=p^k$. For every element of $\SL_n(q)$ it computes the label~\eqref{eq:pslnlabel} and checks that it is a bijection onto $[0,N)$ with $\lambda(I)=0$, constant on the scalar subgroup. Writes \texttt{drafts/psln\_n<n>\_q<q>.json}.\par
\noindent\textbf{Invocation:} \texttt{python3 psln\_\allowbreak{}dense.py n p k [n p k ...] (q = p\textasciicircum{}k) -> writes ../\allowbreak{}drafts/\allowbreak{}psln\_\allowbreak{}n<n>\_\allowbreak{}q<q>.json}.\par
\noindent\textbf{Arguments:} positional argument indices 1:. \textbf{Dependencies:} repository modules gf; external packages numpy.\par
\noindent\textbf{Reads:} no literal file name. \textbf{Writes:} ../\allowbreak{}drafts/\allowbreak{}psln\_\allowbreak{}n\%d\_\allowbreak{}q\%d.json.\par
\paragraph{\texttt{scripts/ree\_law.py}.} Digital law of the small Ree group ${}^2G_2(q)$, $q=3^{2m+1}$ (Definition~\ref{def:ree}, Proposition~\ref{prop:reeU}, Theorem~\ref{thm:ree}), on $7\times7$ matrices over $\mathbb{F}_q$. Arguments: $m$ and optionally \texttt{--ns}. The checks are exhaustive for $q=3$ and random otherwise. Writes \texttt{drafts/ree\_law\_q<q>.json}.\par
\noindent\textbf{Invocation:} \texttt{python3 ree\_\allowbreak{}law.py m [--ns N]; exhaustive for q=3, random checks otherwise; writes ../\allowbreak{}drafts/\allowbreak{}ree\_\allowbreak{}law\_\allowbreak{}q<q>.json}.\par
\noindent\textbf{Arguments:} positional argument indices 1, sys.argv.index('--ns') + 1. \textbf{Dependencies:} repository modules gf.\par
\noindent\textbf{Reads:} no literal file name. \textbf{Writes:} ../\allowbreak{}drafts/\allowbreak{}ree\_\allowbreak{}law\_\allowbreak{}q\%d.json.\par
\paragraph{\texttt{scripts/sz\_law.py}.} Digital law of the Suzuki group $\mathrm{Sz}(q)$, $q=2^{2m+1}$ (Definition~\ref{def:suzuki}, Proposition~\ref{prop:suzukiborel}, Theorem~\ref{thm:suzuki}), on $4\times4$ matrices. Arguments: $m$ and optionally \texttt{--ns}. The checks are exhaustive for $q=8$ and random otherwise. Writes \texttt{drafts/sz\_law\_q<q>.json}.\par
\noindent\textbf{Invocation:} \texttt{python3 sz\_\allowbreak{}law.py m [--ns N] -> checks for q=2\textasciicircum{}(2m+1); exhaustive for q=8, random otherwise; writes ../\allowbreak{}drafts/\allowbreak{}sz\_\allowbreak{}law\_\allowbreak{}q<q>.json}.\par
\noindent\textbf{Arguments:} positional argument indices 1, sys.argv.index('--ns') + 1. \textbf{Dependencies:} repository modules gf.\par
\noindent\textbf{Reads:} no literal file name. \textbf{Writes:} ../\allowbreak{}drafts/\allowbreak{}sz\_\allowbreak{}law\_\allowbreak{}q\%d.json.\par
\paragraph{\texttt{scripts/triality\_cells.py}.} ${}^3D_4(q)$ in the $8$-dimensional triality model over $\mathbb{F}_{q^3}$ (Definition~\ref{def:3d4}, Lemma~\ref{lem:3d4roots}, Propositions~\ref{prop:3d4U} and~\ref{prop:3d4rel}, Theorem~\ref{thm:3d4cells}). It checks that the generators preserve $Q$ and that the root elements $x_R$, $x_S$ lie in $U$. Arguments: $p$, $k$, with $q=p^k$. Writes \texttt{drafts/triality\_cells\_q<q>.json}.\par
\noindent\textbf{Invocation:} \texttt{python3 triality\_\allowbreak{}cells.py p k (q = p\textasciicircum{}k) -> writes ../\allowbreak{}drafts/\allowbreak{}triality\_\allowbreak{}cells\_\allowbreak{}q<q>.json}.\par
\noindent\textbf{Arguments:} positional argument indices 1, 2. \textbf{Dependencies:} repository modules gf; external packages numpy.\par
\noindent\textbf{Reads:} no literal file name. \textbf{Writes:} ../\allowbreak{}drafts/\allowbreak{}triality\_\allowbreak{}cells\_\allowbreak{}q\%d.json.\par
\paragraph{\texttt{scripts/triality\_label.py}.} Closed digital label of ${}^3D_4(q)$ (Corollaries~\ref{cor:3d4closed} and~\ref{cor:3d4law}). For $g=P(t)\,h(\mu^i,\gamma^j)\,u(a,\dots,f)$ the label is \begin{multline*}\lambda(g)=O_w+[t]+q^{\ell(w)}\bigl(i+(q^3-1)\,j\\+(q^3-1)(q-1)\,([a]+q^3[b]+q^4[c]+q^7[d]+q^{10}[e]+q^{11}[f])\bigr).\end{multline*} Arguments: $p$, $k$. Writes \texttt{drafts/triality\_label\_q<q>.json}.\par
\noindent\textbf{Invocation:} \texttt{python3 triality\_\allowbreak{}label.py p k -> ../\allowbreak{}drafts/\allowbreak{}triality\_\allowbreak{}label\_\allowbreak{}q<q>.json}.\par
\noindent\textbf{Arguments:} positional argument indices 1, 2. \textbf{Dependencies:} repository modules triality\_\allowbreak{}mult; external packages numpy.\par
\noindent\textbf{Reads:} no literal file name. \textbf{Writes:} ../\allowbreak{}drafts/\allowbreak{}triality\_\allowbreak{}label\_\allowbreak{}q\%d.json.\par
\paragraph{\texttt{scripts/triality\_mult.py}.} Closed multiplication of ${}^3D_4(q)$ (Proposition~\ref{prop:3d4rel}\,(a)--(e), Theorem~\ref{thm:3d4mult}, Corollaries~\ref{cor:3d4closed} and~\ref{cor:3d4law}). Every closed formula, including the torus and Weyl relations, is compared with the matrix computation in the $8$-dimensional model. Arguments: $p$, $k$. Writes \texttt{drafts/triality\_mult\_q<q>.json}.\par
\noindent\textbf{Invocation:} \texttt{python3 triality\_\allowbreak{}mult.py p k (q = p\textasciicircum{}k) -> ../\allowbreak{}drafts/\allowbreak{}triality\_\allowbreak{}mult\_\allowbreak{}q<q>.json}.\par
\noindent\textbf{Arguments:} positional argument indices 1, 2. \textbf{Dependencies:} repository modules triality\_\allowbreak{}cells; external packages numpy.\par
\noindent\textbf{Reads:} no literal file name. \textbf{Writes:} ../\allowbreak{}drafts/\allowbreak{}triality\_\allowbreak{}mult\_\allowbreak{}q\%d.json.\par

\section{Mathieu, Golay and Conway groups}
\paragraph{\texttt{scripts/co0\_chain.py}.} Stabiliser chain of $\mathrm{Co}_0$ on its $\nu$-numbered minimal vectors.\par
\noindent\textbf{Invocation:} \texttt{Stabiliser chain of Co0 on the nu-numbered minimal vectors. Usage: cd scripts \&\& python3 co0\_\allowbreak{}chain.py}.\par
\noindent\textbf{Arguments:} none. \textbf{Dependencies:} repository modules conway\_\allowbreak{}examples; external packages numpy.\par
\noindent\textbf{Reads:} no literal file name. \textbf{Writes:} ../\allowbreak{}drafts/\allowbreak{}co0\_\allowbreak{}chain.json, ../\allowbreak{}drafts/\allowbreak{}co0\_\allowbreak{}chain.log, ../\allowbreak{}drafts/\allowbreak{}co0\_\allowbreak{}chain.tex.\par
\paragraph{\texttt{scripts/co3\_co1.py}.} Stabiliser chains for $\mathrm{Co}_3$ and $\mathrm{McL}$, and the digital law of $\mathrm{Co}_1$.\par
\noindent\textbf{Invocation:} \texttt{Chains for Co3 (and McL) and the law of Co1. Usage: cd scripts \&\& python3 co3\_\allowbreak{}co1.py}.\par
\noindent\textbf{Arguments:} none. \textbf{Dependencies:} repository modules co0\_\allowbreak{}chain; external packages numpy.\par
\noindent\textbf{Reads:} no literal file name. \textbf{Writes:} ../\allowbreak{}drafts/\allowbreak{}co3\_\allowbreak{}co1.json, ../\allowbreak{}drafts/\allowbreak{}co3\_\allowbreak{}co1.tex.\par
\paragraph{\texttt{scripts/conway\_counts.py}.} Integer checks for Section~\ref{sec:conway}.\par
\noindent\textbf{Invocation:} \texttt{Integer checks for the Conway section. Usage: cd scripts \&\& python3 conway\_\allowbreak{}counts.py}.\par
\noindent\textbf{Arguments:} positional argument indices 0. \textbf{Dependencies:} repository modules mathieu\_\allowbreak{}law.\par
\noindent\textbf{Reads:} no literal file name. \textbf{Writes:} ../\allowbreak{}drafts/\allowbreak{}conway\_\allowbreak{}counts.json.\par
\paragraph{\texttt{scripts/conway\_examples.py}.} Worked examples for Section~\ref{sec:conway}.\par
\noindent\textbf{Invocation:} \texttt{Worked examples for the Conway section. Usage: cd scripts \&\& python3 conway\_\allowbreak{}examples.py}.\par
\noindent\textbf{Arguments:} positional argument indices 0. \textbf{Dependencies:} repository modules mathieu\_\allowbreak{}law; external packages numpy.\par
\noindent\textbf{Reads:} no literal file name. \textbf{Writes:} ../\allowbreak{}drafts/\allowbreak{}conway\_\allowbreak{}examples.json, ../\allowbreak{}drafts/\allowbreak{}conway\_\allowbreak{}examples.tex.\par
\paragraph{\texttt{scripts/golay\_check.py}.} Golay-code checks for the Mathieu chains of \texttt{mathieu\_law.py}. Writes \texttt{drafts/golay\_check.json}.\par
\noindent\textbf{Invocation:} \texttt{Golay-code checks for the Mathieu chains of mathieu\_\allowbreak{}law.py. Usage: python3 golay\_\allowbreak{}check.py (writes ../\allowbreak{}drafts/\allowbreak{}golay\_\allowbreak{}check.json)}.\par
\noindent\textbf{Arguments:} positional argument indices 0. \textbf{Dependencies:} repository modules mathieu\_\allowbreak{}law.\par
\noindent\textbf{Reads:} no literal file name. \textbf{Writes:} ../\allowbreak{}drafts/\allowbreak{}golay\_\allowbreak{}check.json.\par
\paragraph{\texttt{scripts/golay\_law.py}.} Golay-code law of $M_{24}$ and related checks. The label of $g$ is $L(g)=o+759\,\bigl(v+16\,\mu(A)\bigr)$, where $o$ is the index of the octad $g(O_0)$ and $h=s_O^{-1}g$ lies in the octad stabiliser $K=2^4{:}\Alt_8$; $h$ acts on the complement $C_0=\mathbb{F}_2^4$ as $x\mapsto Ax+v$, and $\mu(A)$ is the falling rank of the columns of $A$, with radices $15$, $14$, $12$, $8$. Writes \texttt{drafts/golay\_law.json}.\par
\noindent\textbf{Invocation:} \texttt{Golay-code law of M24 and related checks. Usage: python3 golay\_\allowbreak{}law.py (writes ../\allowbreak{}drafts/\allowbreak{}golay\_\allowbreak{}law.json)}.\par
\noindent\textbf{Arguments:} positional argument indices 0. \textbf{Dependencies:} repository modules mathieu\_\allowbreak{}law.\par
\noindent\textbf{Reads:} no literal file name. \textbf{Writes:} ../\allowbreak{}drafts/\allowbreak{}golay\_\allowbreak{}law.json.\par
\paragraph{\texttt{scripts/m12\_ternary.py}.} $M_{12}$ inside $M_{24}$ (dodecads), the ternary Golay code and $2.M_{12}$, and a Golay-adapted tail of the $M_{24}$ chain. Writes \texttt{drafts/m12\_ternary.json}.\par
\noindent\textbf{Invocation:} \texttt{python3 m12\_\allowbreak{}ternary.py (writes ../\allowbreak{}drafts/\allowbreak{}m12\_\allowbreak{}ternary.json)}.\par
\noindent\textbf{Arguments:} positional argument indices 0. \textbf{Dependencies:} repository modules mathieu\_\allowbreak{}law.\par
\noindent\textbf{Reads:} no literal file name. \textbf{Writes:} ../\allowbreak{}drafts/\allowbreak{}m12\_\allowbreak{}ternary.json.\par
\paragraph{\texttt{scripts/mathieu\_examples.py}.} Generates \texttt{drafts/mathieu\_examples.tex}: worked examples of the chain law (Theorem~\ref{thm:chainlaw}) for $M_{11}$ and $M_{24}$. All numbers are computed from the chains of \texttt{mathieu\_law.py} and cross-checked against \texttt{Law.mul}.\par
\noindent\textbf{Invocation:} \texttt{python3 mathieu\_\allowbreak{}examples.py}.\par
\noindent\textbf{Arguments:} positional argument indices 0. \textbf{Dependencies:} repository modules mathieu\_\allowbreak{}law.\par
\noindent\textbf{Reads:} no literal file name. \textbf{Writes:} ../\allowbreak{}drafts/\allowbreak{}mathieu\_\allowbreak{}examples.tex.\par
\paragraph{\texttt{scripts/mathieu\_law.py}.} Stabiliser-chain digital laws of the Mathieu groups $M_{11}$, $M_{12}$, $M_{22}$, $M_{23}$ and $M_{24}$ (Theorem~\ref{thm:transitive}, iterated). Optional arguments \texttt{--ns}, \texttt{--assoc} and \texttt{--words} (defaults $1000$, $200$, $200$). Permutations act on $\{0,\dots,n-1\}$ as tuples $p$ with $p[i]$ the image of $i$, and $(s\circ t)(i)=s(t(i))$. Writes \texttt{drafts/mathieu.json}.\par
\noindent\textbf{Invocation:} \texttt{python3 mathieu\_\allowbreak{}law.py [--ns 1000] [--assoc 200] [--words 200] (writes ../\allowbreak{}drafts/\allowbreak{}mathieu.json)}.\par
\noindent\textbf{Arguments:} --ns default=1000; --assoc default=200; --words default=200. \textbf{Dependencies:} Python standard library.\par
\noindent\textbf{Reads:} no literal file name. \textbf{Writes:} ../\allowbreak{}drafts/\allowbreak{}mathieu.json.\par
\endgroup

\chapter{Solutions to the exercises}\label{ch:solutions}

The solutions follow the order of the exercises in each chapter. Where an exercise asks for a computation that is long but routine,
the answer is given together with the key step; where it depends on a convention of the chapter, the convention is recalled.

\section{The digit function}
\begin{enumerate}
\item Put $u=\fl{x/N^k}$. Since $\fl{\fl{x/N^k}/N}=\fl{x/N^{k+1}}$,~\eqref{eq:digit} gives $\dg{N}{k}{x}=u-N\fl{u/N}=u\bmod N\in\Ap N$.
Summing, $\sum_{k=0}^{K}\big(\fl{x/N^k}-N\fl{x/N^{k+1}}\big)N^k=x-N^{K+1}\fl{x/N^{K+1}}$ telescopes, and the last term vanishes once
$N^{K+1}>x$.
\item $\fl{(x+N^{k+1}y)/N^k}=\fl{x/N^k}+py$ and $\fl{(x+N^{k+1}y)/N^{k+1}}=\fl{x/N^{k+1}}+y$, so the two contributions of $y$ cancel in
\eqref{eq:digit}. For the second identity, $\fl{px/N^k}=\fl{x/N^{k-1}}$ and $\fl{px/N^{k+1}}=\fl{x/N^{k}}$.
\item With $u=\fl{x/N^{2k}}$ one has $\dg{N^2}{k}{x}=u\bmod N^2=(u\bmod N)+N\,(\fl{u/N}\bmod N)=\dg{N}{2k}{x}+N\,\dg{N}{2k+1}{x}$.
In general $\dg{N^{n}}{k}{x}=\sum_{j<n}\dg{N}{nk+j}{x}\,N^j$.
\item By Theorem~\ref{thm:coding} a code is the same as a table $(x,y)\mapsto x*y$ on $\Ap N$. Commutativity forces equal entries at
$(x,y)$ and $(y,x)$, so the table is fixed by its values on the $N(N+1)/2$ unordered pairs, each chosen freely in $\Ap N$.
\item For an integer $x$, $|x|<\tfrac12$ holds only for $x=0$ and $|x|=\tfrac12$ never holds, so $\Bf(x)=\delta_{x,0}$ and
$\delta_{m,n}=\Bf(m-n)$. Next, $|x-\tfrac{n-1}2|<\tfrac n2$ means $-\tfrac12<x<n-\tfrac12$, i.e.\ $0\le x\le n-1$, and equality is
impossible for integers; hence $\mathbf 1_{\Ap n}(x)=\Bf\big(x-\tfrac{n-1}2,\tfrac n2\big)$.
\item Write $a=a'+N^kA$, $b=b'+N^kB$ with $a'=a\bmod N^k$, $b'=b\bmod N^k$. Then $0\le a'+b'<2N^k$, so $c_k\in\{0,1\}$, and
$\fl{(a+b)/N^k}=A+B+c_k$. Reducing modulo $N$ gives $\dg{N}{k}{a+b}\equiv\dg{N}{k}{a}+\dg{N}{k}{b}+c_k$. The carry $c_k$ is the
indicator of $a'+b'\ge N^k$, which is exactly what Proposition~\ref{prop:Bcarry} writes as a $\Bf$-function.
\item By Kummer's theorem the exponent of $p$ in $\binom{2n}{n}$ is the number of carries in $n+n$ in base $p$. For $p=2$ and
$n\ge1$ the highest nonzero binary digit of $n$ produces a carry, so $\binom{2n}n$ is even. For odd $p$, by induction on the position
(the incoming carry being $0$), there is no carry at all if and only if $2\,\dg{p}{k}{n}\le p-1$ for every $k$, that is,
$\dg{p}{k}{n}<p/2$.
\end{enumerate}

\section{Finite magmas and \texorpdfstring{$r$}{r}-ary operators}
\begin{enumerate}
\item The identity has digits $\dg{N}{i}{\iota_N}=i$, so $\iota_N=\sum_{i<N}i\,N^i$. With $\sum_{i<n}i\,t^i=t\big(1-nt^{n-1}+(n-1)t^n\big)/(1-t)^2$
at $n=t=N$ this is $N\big(1-N^N+(N-1)N^N\big)/(N-1)^2=N\big(1+(N-2)N^N\big)/(N-1)^2$. For $N=3$: $3\cdot28/4=21$.
\item The six digit strings $(R(0),R(1),R(2))$ give $(0,1,2)\mapsto21$ (identity), $(1,2,0)\mapsto7$, $(2,0,1)\mapsto11$,
and the transpositions $(1,0,2)\mapsto19$, $(0,2,1)\mapsto15$, $(2,1,0)\mapsto5$ (swapping $0,1$; $1,2$; $0,2$). With
$(R\circ S)(i)=R(S(i))$ as in~\eqref{eq:comp}, the table of $R\circ S$ is
\[
\begin{array}{c|cccccc}
\circ & 5&7&11&15&19&21\\\hline
5 &21&19&15&11&7&5\\
7 &15&11&21&19&5&7\\
11&19&21&7&5&15&11\\
15&7&5&19&21&11&15\\
19&11&15&5&7&21&19\\
21&5&7&11&15&19&21
\end{array}
\]
(rows $R$, columns $S$). Each row and column is a permutation of the six codes.
\item Composing with a transposition of two adjacent positions changes the number of inversions by exactly $\pm1$ and changes the sign,
and every permutation is a product of such transpositions; hence the sign is $(-1)^{\iota(R)}$ (Theorem~\ref{thm:sign}). In $\Sym_3$:
$\iota=0,2,2$ for $21,7,11$ (even) and $\iota=1,1,3$ for $19,15,5$ (odd).
\item A magma is a free choice of $N^2$ table entries. Commutativity leaves $N(N+1)/2$ free entries (previous chapter's exercise). If
$0$ is a two-sided identity, row $0$ and column $0$ ($2N-1$ entries) are forced and the other $(N-1)^2$ entries are free.
\item There are $2^4=16$ magmas on $\Ap2$, and the relabelling group is $\{1,\tau\}$ with $\tau$ the swap. A table is fixed by $\tau$
if and only if $f(1,1)=\tau f(0,0)$ and $f(1,0)=\tau f(0,1)$, which leaves $4$ tables. By Burnside's lemma the number of orbits is
$(16+4)/2=10$.
\item Row $0$ and column $0$ are forced. Row $1$ starts with $1$; if $1*1=0$ then $1*2=2$, which repeats the entry $2$ in column $2$.
So $1*1=2$, $1*2=0$, and the Latin property forces $2*1=0$, $2*2=1$: this is addition modulo $3$.
\item With the labelling $m=m_0+2m_1$ the law of $C_2\oplus C_2$ is $m*n=m\oplus_2 n$ (bitwise sum), while that of $C_4$ is
$\dg{4}{0}{m+n}$. Both are commutative, so the order of the two operands in the index $x_0+4x_1$ of Theorem~\ref{thm:coding} does not
matter. The codes are $\sum_{x_0,x_1<4}(x_0\oplus_2x_1)\,4^{x_0+4x_1}=458142180$ and $\sum_{x_0,x_1<4}\dg{4}{0}{x_0+x_1}\,4^{x_0+4x_1}=2471377380$;
they already differ at the digit with $x_0=x_1=1$, where $1\oplus_21=0$ but $1+1=2$.
\item The group $\Sym_2^{3}$ has order $8$. Writing a table as $(R(0,0),R(0,1),R(1,0),R(1,1))$, the isotopy classes are $\{0000,1111\}$ (constant maps; code $0$), the eight tables with exactly one or exactly three digits $1$ (code $8$, the table $0001$), $\{0011,1100\}$ ($a*b=a$ and $a*b=1-a$; code $12$), $\{0101,1010\}$ ($a*b=b$ and $a*b=1-b$; code $10$) and $\{0110,1001\}$ (addition modulo $2$ and its complement; code $6$). The sizes $2+8+2+2+2$ add up to $16$, and Theorem~\ref{thm:magma:isotopy} gives $5$. An isotopy keeps the variable on which the result depends, so the left and right projections lie in different classes; passing to the opposite magma exchanges them.
\item The identity of $\Sym_3$ fixes all $3^{9}=19683$ codes; each of the three transpositions fixes $1\cdot3^{4}=81$ and each of the two $3$-cycles fixes $3^{3}=27$. Thus $\frac16(19683+243+54)=19980/6=3330$.
\item By associativity, $m(x,y,y)=x*y^{-1}*y=x$, $m(y,y,x)=x$ and $m(m(x,y,z),u,v)=x*y^{-1}*z*u^{-1}*v=m(x,y,m(z,u,v))$.
By Proposition~\ref{prop:heap}, $m(x,y,z)=xy^{-1}z$ in the group $*_e$. Then $x*_{e'}y=xe'^{-1}y$, and $\varphi(x)=m(x,e,e')=xe'$ satisfies
$\varphi(x)*_{e'}\varphi(y)=xe'e'^{-1}ye'=\varphi(x*_e y)$; it is bijective with inverse $x\mapsto xe'^{-1}$, and $\varphi(e)=e'$.
\item For $f$ all three bracketings of five arguments give $x_1-x_2+x_3-x_4+x_5$, and for $g$ they give $x_1+\dots+x_5+2$; each
equation in one unknown slot is solved uniquely by subtraction. For $f$ take $\theta(x)=-x$ and $b=0$: $\theta(b)=b$, $\theta^2=\mathrm{id}$
is conjugation by $0$, and $x+\theta(y)+\theta^2(z)+b=x-y+z$. For $g$ take $\theta=\mathrm{id}$ and $b=1$; the conditions hold because
the group is abelian.
\item A normalised cochain vanishes when an argument is $0$, so only $f(1,1)\in\{0,1\}$ is free: $(q-1)^2=1$ free value and
$n^{(q-1)^2}=2$ cochains. With $x_0+2x_1$ as index, $F=\sum_x f(x)\,2^{x_0+2x_1}\in\{0,8\}$. The cocycle identity
$f(u,v)+f(u+v,w)=f(v,w)+f(u,v+w)$ holds trivially when an argument is $0$, and for $u=v=w=1$ both sides equal $f(1,1)$. A normalised
coboundary is $\delta g(u,v)=g(u)+g(v)-g(u+v)$ with $g(0)=0$, and $\delta g(1,1)=2g(1)\equiv0$. Hence $H^2(\Z_2,\Z_2)\cong\Z_2$. The cocycle
$0$ gives $\Z_2\times\Z_2$, and $f(1,1)=1$, the carry $c_2$ of Proposition~\ref{prop:Bcarry}, gives $\Z_4$.
\end{enumerate}

\section{Diagrammatic construction of finite groups}
The four exercises of Section~\ref{sec:ch3-exercises} are solved in order. Recall that a building block with ground layer $N=\Z_n$ and top layer
$Q=\Z_q$ is given by an action code $\Phi\in[0,n^{nq})$, with $\varphi_A(v)=\dg{n}{v+nA}{\Phi}$, and a carry code $F$, and that
\eqref{dg:y0} reads $y^{(0)}=a^{(0)}+\varphi_{a^{(1)}}(b^{(0)})+f(a^{(1)},b^{(1)})\bmod n$.
\begin{enumerate}
\item Take three $\Z_2$ layers, $h=2$ on top and $h=0$ at the ground, with identity actions. The top layer has no carry. Layer $1$ receives
the carry of the top digits, $f_1(x,y)=xy$, so $F_1=\sum_{x,y<2}xy\,2^{y+2x}=2^3=8$. Layer $0$ receives the carry of the two upper digits
read as an element $u\in\Z_4$, $f_0(u,v)=\fl{(u+v)/4}$; the nonzero digits sit at $v+4u\in\{7,10,11,13,14,15\}$, so
$F_0=2^7+2^{10}+2^{11}+2^{13}+2^{14}+2^{15}=60544$. Since $f_0$ is symmetric, $F_0$ does not depend on the order of $u,v$ in the index.
The result is binary addition modulo $8$.
\item For $D_n$ (order $2n$) take $Q=\Z_2$, $F=0$, $\varphi_0=\mathrm{id}$ and $\varphi_1(v)=-v$. Then
$\Phi=\sum_{v<n}v\,n^v+n^n\sum_{v<n}\dg{n}{0}{-v}\,n^v$, which gives $\Phi=21+27\cdot15=426$ for $n=3$,
$\Phi=228+256\cdot108=27876$ for $n=4$, and $\Phi=2930+3125\cdot970=3034180$ for $n=5$. Formula~\eqref{dg:y0} becomes
$y^{(0)}=a^{(0)}+(-1)^{a^{(1)}}b^{(0)}\bmod n$, which is the rule $x^{a}y^{\alpha}\cdot x^{b}y^{\beta}=x^{a+(-1)^\alpha b}y^{\alpha+\beta}$.
\item The Heisenberg recipe has trivial action and carry $f\big((a_1,a_2),(b_1,b_2)\big)=a_1b_2\bmod p$. The cocycle condition
$f(a,b)+f(a+b,c)=f(a,b+c)+f(b,c)$ reads $a_1b_2+(a_1+b_1)c_2=a_1(b_2+c_2)+b_1c_2$, which holds identically. By~\eqref{dg:y0} two
elements commute if and only if $f(a,b)=f(b,a)$, that is, $a_1b_2=b_1a_2$; this holds for all $b$ only when $a_1=a_2=0$. The centre is
the ground layer, of order $p$.
\item For $Q_{16}$ take $K=\Z_8$, $Q=\Z_2$, $\varphi_1(v)=-v$ and the single carry $f(1,1)=4$, so $F=4\cdot8^{3}=2048$. The value $4$ is
fixed by $\varphi_1$, so the cocycle condition holds. Then $y^2=x^4$ is the unique involution and $y$ has order $4$, which identifies
$Q_{16}$. The diagram of $D_8$ (order $16$) has the same action arrow and no carry arrow; the two groups differ only in the carry.
\end{enumerate}

\section{The main theorems of group theory in digital form}
\begin{enumerate}
\item $\lambda_a(x)=\dg{4}{0}{x+a}$ has digits $(a,a+1,a+2,a+3)\bmod4$, giving $228,57,78,147$. For example $\lambda_1\circ\lambda_1$ shifts
by $2$ and has code $78=\lambda_2$; in general the digits of $\lambda_a\circ\lambda_b$ are $x+a+b$.
\item A group of prime order $p$ is generated by any $g\ne1$, so $g^m\mapsto m$ is a labelling with law $m*n=\dg{p}{0}{m+n}$; any other
labelling differs from it by a relabelling.
\item In $\Alt_4$: $n_2=1$ (the Klein group is normal) and $n_3=4$. In $\Sym_4$: $n_2=3$ (the dihedral subgroups of order $8$) and $n_3=4$.
\item $n_q\mid p$ and $n_q\equiv1\pmod q$ with $p<q$ give $n_q=1$; $n_p\mid q$, $n_p\equiv1\pmod p$ and $p\nmid q-1$ give $n_p=1$. Both
Sylow subgroups are normal and intersect trivially, so $G\cong C_p\times C_q\cong C_{pq}$. For $15=3\cdot5$, $3\nmid4$.
\item The orbit of $\{0,1\}$ is all $6$ pairs, and its stabiliser is $\langle(0\,1),(2\,3)\rangle$ of order $4$; indeed $24=6\cdot4$.
\item $\Sym_4>\Alt_4>V_4>C_2>1$, with factors $C_2,C_3,C_2,C_2$.
\item The factors of $0<\langle6\rangle<\langle2\rangle<C_{12}$ are $C_2,C_3,C_2$. Refine the second series to
$0<\langle4\rangle<\langle2\rangle<C_{12}$, with factors $C_3,C_2,C_2$. The two series are equivalent: same factors up to order.
\item The $24$ rotations act on the six faces with types $1^{6}$ (identity), $1^{2}4$ ($6$ quarter turns about face axes), $1^{2}2^{2}$ ($3$ half turns about face axes), $3^{2}$ ($8$ third turns about vertex axes) and $2^{3}$ ($6$ half turns about edge axes). So $Z=\frac1{24}(x_1^{6}+6x_1^{2}x_4+3x_1^{2}x_2^{2}+8x_3^{2}+6x_2^{3})$, the number of orbits is $(m^{6}+3m^{4}+12m^{3}+8m^{2})/24$, and it equals $10$ for $m=2$ and $57$ for $m=3$.
\item The edge permutations have types: identity $1^{6}$; the $6$ transpositions of $\Ap4$ give $1^{2}2^{2}$; the $3$ double transpositions give $1^{2}2^{2}$; the $8$ three-cycles give $3^{2}$; the $6$ four-cycles give $2\,4$. The numbers of cycles are $6,4,4,2,2$, so the number of orbits on $2$-colourings is $\frac1{24}(64+6\cdot16+3\cdot16+8\cdot4+6\cdot4)=264/24=11$.
\end{enumerate}

\section{The Schreier--Sims algorithm and conjugacy classes}
\begin{enumerate}
\item With base $(0,1,2)$ the basic orbits have lengths $4,3,2$. At the end of the trace every Schreier generator sifts to the identity,
and $|\Sym_4|=4\cdot3\cdot2=24$.
\item $\Alt_5$ is transitive on $\Ap5$, the stabiliser of $0$ is $\Alt_4$, transitive on the other $4$
points, and the stabiliser of $0,1$ is $C_3$, regular on the remaining $3$. So the basic orbits are $5,4,3$ and $|\Alt_5|=60$.
\item For an irredundant base each step of the stabiliser chain has index at least $2$, so $|G|\ge2^k$; each basic orbit has length
at most $n$, so $|G|\le n^k$.
\item A breadth-first tree is a shortest-path tree, so its depth is the eccentricity of the root, which is at most the diameter. For
$\Sym_n$ generated by the adjacent transpositions, the orbit graph of $0$ is a path $0-1-\dots-(n-1)$ of depth $n-1$.
\item Write $m=t_1+12t_2+132t_3+1320t_4+11880t_5$ with $t_l$ below the basic orbit lengths $12,11,10,9,8$, and map $m$ to
$s_{1,t_1}\cdots s_{5,t_5}$ as in Corollary~\ref{cor:ss:label}. Label $0$ is the identity; label $1$ is the transversal element
$s_{1,1}$; label $95039$ is $s_{1,11}s_{2,10}s_{3,9}s_{4,8}s_{5,7}$.
\item $G$ is $k$-transitive if and only if the first $k$ basic orbits have lengths $n,n-1,\dots,n-k+1$.
\item Search the stabiliser chain depth first, at each level keeping only images compatible with commuting with the element, and prune
with the cycle structure. Applied to class representatives of $M_{12}$ this recovers the centraliser orders of Table~\ref{tab:cc:m12}.
\item $\Sym_5$ has $7$ classes: cycle types $1^5,2,3,4,5,2^2,3\,2$ with sizes $1,10,20,30,24,15,20$ and centraliser orders
$120,12,6,4,5,8,6$. The power maps send each type to the type of the power, e.g.\ $(3\,2)^2$ has type $3$ and $(3\,2)^3$ type $2$.
$\Alt_5$ has $5$ classes, $1,2^2,3,5A,5B$, of sizes $1,15,20,12,12$ and centraliser orders $60,4,3,5,5$; the square map exchanges $5A$
and $5B$. In both cases the sizes add up to $|G|$, in agreement with Theorem~\ref{thm:cc:count}.
\item For odd $n$: $\{1\}$, the $(n-1)/2$ pairs $\{r^{\pm i}\}$ and one class of $n$ reflections, $(n+3)/2$ classes, with centralisers of
order $2n,n,2$. For even $n$: $\{1\}$, $\{r^{n/2}\}$, $(n-2)/2$ pairs and two classes of reflections, $n/2+3$ classes, with centralisers
of order $2n,2n,n,4$.
\item Each $C_G(x)$ has index at most $2$, so it is normal with Abelian quotient; hence $G'\le C_G(x)$ for all $x$, i.e.\ $G'\le Z(G)$.
Then $y\mapsto[x,y]$ is a homomorphism with image in $\{1,c_x\}$, where $x^G=\{x,xc_x\}$. Suppose $[a,b]=c$ and $[u,v]=d$ with $c\ne d$
both nontrivial. Each of $[a,u],[a,v],[b,u],[b,v]$ lies in $\{1,c\}\cap\{1,d\}=\{1\}$. Then $[au,b]=c$ and $[au,v]=d$, so the image
for $x=au$ has three elements, a contradiction. Thus all commutators lie in a subgroup $\{1,c\}$ and $|G'|\le2$.
\item A class of $\Sym_n$ splits in $\Alt_n$ exactly when its parts are distinct and odd (Theorem~\ref{thm:sn:altclasses}). For $n=7$ only
$(7)$ qualifies; there are $8$ even types, so $k(\Alt_7)=8+1=9$. For $n=8$ the split types are $(7,1)$ and $(5,3)$; there are $12$ even
types, so $k(\Alt_8)=12+2=14$.
\item $\sigma$ and $\sigma^k$ with $\gcd(k,o(\sigma))=1$ have the same cycle type, hence are conjugate in $\Sym_p$. In $\Alt_p$ the
unsplit classes are real. A split class of type $\lambda$ with $k$ parts is real if and only if $(p-k)/2$ is even; in particular the
$p$-cycles are real if and only if $p\equiv1\pmod4$.
\item $|PSL(2,7)|=168$, with classes $1A,2A,3A,4A,7A,7B$ of sizes $1,21,56,42,24,24$. An element of order $7$ fixes one of the $8$ points
and is unipotent; $u$ is conjugate to $u^a$ only when $a$ is a square modulo $7$, and $-1$ is not, so $7A$ and $7B$ are non-real.
\item During the orbit enumeration, store with each new class member $y=x^g$ the generator that produced it and the parent element.
The resulting tree is a Schreier tree, and the conjugating element is read off along the path to the root.
\item This is the coupon collector problem with probabilities $|K|/|G|$. The expected number of draws is at least $1/\min_K(|K|/|G|)$,
and the minimum is attained at the identity class, so it is at least $|G|=95040$. The rare small classes dominate; in practice one
adds the identity by hand.
\item The classes of $PSL(2,7)$ on the $8$ points have types $1^{8}$ ($1$ element), $2^{4}$ ($21$), $3^{2}1^{2}$ ($56$), $4^{2}$ ($42$) and $7\,1$ ($24+24$), with $8,4,4,2,2$ cycles. So $Z=\frac1{168}(x_1^{8}+21x_2^{4}+56x_1^{2}x_3^{2}+42x_4^{2}+48x_1x_7)$ and the number of orbits on $2$-colourings is $(256+21\cdot16+56\cdot16+42\cdot4+48\cdot4)/168=1848/168=11$. The group is $3$-homogeneous, so it has one orbit on $k$-subsets for $k\le3$ and $k\ge5$ (eight orbits in all); the remaining $3$ orbits are on $4$-subsets, of sizes $14,14,42$ (checked by enumeration in \texttt{scripts/polya\_verify.py}).
\end{enumerate}

\section{Characters and character tables in digital form}
\begin{enumerate}
\item If $\rho'(g)=T\rho(g)T^{-1}$, then $\operatorname{tr}\rho'(g)=\operatorname{tr}\rho(g)$ by the cyclicity of the trace.
\item The table of $S_3$ on the classes of sizes $1,3,2$ has rows $(1,1,1)$, $(1,-1,1)$, $(2,0,-1)$. For example,
$\tfrac16(1\cdot2+3\cdot0-2\cdot1)=0$ and $\tfrac16(4+0+2)=1$; the column relations hold as well.
\item Let $e_\chi=\tfrac{\chi(1)}{|G|}\sum_g\chi(g^{-1})g$. The generalised orthogonality relation gives $e_\chi e_\psi=\delta_{\chi\psi}e_\chi$,
and $\sum_\chi e_\chi=1$ because $\sum_\chi\chi(1)\chi(g^{-1})=|G|\delta_{g,1}$.
\item If $G\cong C_{n_1}\oplus\dots\oplus C_{n_r}$, the linear characters are
$\chi_k(x)=\exp\big(2\pi i\sum_jk_jx_j/n_j\big)$ with $k_j\in\Ap{n_j}$; there are $|G|$ of them.
\item With $\pi=1_H^G$, Frobenius reciprocity gives $\langle\pi,\chi\rangle=\langle1_H,\chi_H\rangle$. For instance, for a doubly
transitive action $\langle\pi,\pi\rangle=2$, so $\pi=1+\chi$ with $\chi$ irreducible.
\item Use $S^2\chi(g)=\tfrac12(\chi(g)^2+\chi(g^2))$ and $\Lambda^2\chi(g)=\tfrac12(\chi(g)^2-\chi(g^2))$. In $D_{10}$ the two
characters of degree $2$ take the values $2\cos(2\pi/5)$ and $2\cos(4\pi/5)$ at a rotation. One finds $S^2\psi_1=1+\psi_2$ and
$S^2\psi_2=1+\psi_1$, while $\Lambda^2\psi_i$ is the sign character. So $S^2$ tells the two apart.
\item By Theorem~\ref{thm:centralcharacter}, $\omega_\chi(K)=|K|\chi(g_K)/\chi(1)$ is an algebraic integer. Then
$|G|/\chi(1)=\sum_K\omega_\chi(K)\overline{\chi(g_K)}$ is an algebraic integer and rational, hence an integer; this is
Corollary~\ref{cor:degreedivides}.
\item The central characters satisfy $\omega(K_i)\omega(K_j)=\sum_ka_{ijk}\,\omega(K_k)$, so they are the common eigenvectors of the
matrices $(a_{ijk})_{jk}$. The degree follows from $\chi(1)^2\sum_K|\omega_\chi(K)|^2/|K|=|G|$, and then
$\chi(g_K)=\chi(1)\omega_\chi(K)/|K|$. This is the Burnside--Dixon method~\cite{Dixon1967}.
\end{enumerate}

\section{Cyclic and Abelian groups}
\begin{enumerate}
\item By Theorem~\ref{thm:cycsub} the subgroups are $\langle d\rangle=\{0,d,2d,\dots\}$, $d\mid n$. The element $m$ has order
$n/\gcd(m,n)$, which equals $d$ for exactly $\varphi(d)$ values of $m$.
\item $3$ modulo $7$ (powers $1,3,2,6,4,5$), $2$ modulo $11$ and $2$ modulo $13$. The relabelling $m\mapsto g^m$ turns multiplication
into $\dg{p-1}{0}{m+n}$.
\item $4a+9b\equiv a\pmod3$ and $\equiv b\pmod4$, so the map is the inverse of $x\mapsto(x\bmod3,x\bmod4)$ and is additive, as in
Theorem~\ref{thm:crt}.
\item If $\gcd(m,n)=1$, the element $(1,1)$ has order $mn$. Otherwise every element has order dividing $\operatorname{lcm}(m,n)<mn$.
\item $72=2^3\cdot3^2$, giving $3\cdot2=6$ groups: $C_8\oplus C_9$, $C_2\oplus C_4\oplus C_9$, $C_2^3\oplus C_9$, $C_8\oplus C_3^2$,
$C_2\oplus C_4\oplus C_3^2$ and $C_2^3\oplus C_3^2$.
\item By induction $5^{2^{j}}\equiv1+2^{j+2}\pmod{2^{j+3}}$, so $5$ has order $2^{k-2}$. Every power of $5$ is $\equiv1\pmod4$, so
$-1\notin\langle5\rangle$, and $|\langle-1\rangle\times\langle5\rangle|=2^{k-1}$. For $k=4$: $\langle5\rangle=\{1,5,9,13\}$ and $-1=15$.
\end{enumerate}

\section{Nonabelian groups}
\begin{enumerate}
\item With $m=m_0+3m_1\leftrightarrow x^{m_0}y^{m_1}$, send $x\mapsto7$ and $y\mapsto19$. The relation holds: $19\circ7\circ19=11=7^{-1}$.
The labels $0,\dots,5$ go to $21,7,11,19,5,15$.
\item Replace the generator $y$ of $C_s$ by $y^k$, which is again a generator because $\gcd(k,s)=1$; it acts by $r^k=r'$.
\item $x^a$ is central if and only if $2a\equiv0\pmod q$. Reflections do not commute with $x$ for $q\ge3$. Hence the centre is $\{0\}$
for odd $q$ and $\{0,q/2\}$ for even $q$.
\item The elements $\pm i,\pm j,\pm k$ have order $4$, so $-1$ is the only involution. The subgroups $\langle-1\rangle$, $\langle i\rangle$,
$\langle j\rangle$, $\langle k\rangle$ are central or of index $2$, hence normal. $D_8$ has five involutions.
\item $x^ay^\alpha$ commutes with $x$ only if $r^\alpha\equiv1$, i.e.\ $\alpha=0$, and $x^a$ commutes with $y$ only if $ra\equiv a$,
i.e.\ $a=0$.
\item $3^3=27\equiv1\pmod{13}$, so $s=3$. Theorem~\ref{thm:metacyclic-chars} gives $s+(q-1)/s=3+4=7$ classes, with $3$ linear
characters and $4$ of degree $3$; indeed $3+4\cdot9=39$.
\item The law of Theorem~\ref{thm:wreath} lives on $\Ap{N^{\ell}s}$. In $C_2\wr C_2$ the element $(1,0)\sigma$ squares to $(1,1)\ne0$,
so it has order $4$, and the group is nonabelian with the five involutions $(1,0),(0,1),(1,1),\sigma,(1,1)\sigma$. Hence it is $D_8$.
\item For $C_3$ on three points, $Z(C_3)=\frac13(x_1^{3}+2x_3)$, so $z_1=\frac13(x_1^{3}+2x_3)$ and $z_2=\frac13(x_2^{3}+2x_6)$. Then $Z(C_3\wr C_2)=\frac12(z_1^{2}+z_2)=\frac1{18}\big(x_1^{6}+4x_1^{3}x_3+4x_3^{2}+3x_2^{3}+6x_6\big)$; the coefficients add up to the order $18$. At $x_j=2$ this is $(64+64+16+24+12)/18=10$.
\item $N_{C_2}(m)=\frac12(m^{2}+m)=q$. By Corollary~\ref{cor:wreath:colour} the number of orbits is $N_{\Sym_n}(q)=\frac1{n!}\sum_\sigma q^{c(\sigma)}=\binom{q+n-1}{n}$ by Corollary~\ref{cor:sn:cycleindex}.
\item $N_{C_p}(t)=\frac1p\big(t^{p}+(p-1)t\big)$, since $C_p$ has the identity with $p$ cycles and $p-1$ elements with one cycle each. Corollary~\ref{cor:wreath:colour} gives $N_{C_p}(N_{C_p}(m))$. For $p=3$, $m=2$: $q=(8+4)/3=4$ and $(64+8)/3=24$.
\item Theorem~\ref{thm:debruijn}: the pair (swap, complement) has $c_1=0,c_2=1$. The sum
\[
\sum_{g,h}\prod_k\operatorname{fix}_k(h)^{c_k(g)}
\]
is $2^{2}+0^{2}+2+2=8$ for $g\in\{1,\text{swap}\}$ and $h\in\{1,\text{complement}\}$ (for the swap, $\operatorname{fix}_2=2$ in both cases), so there are $8/4=2$ orbits, namely $\{00,11\}$ and $\{01,10\}$. The wreath product $C_2\wr C_2$ contains independent complements at the two positions, so it has one orbit, in agreement with $N_{C_2}(1)=1$.
\item By Proposition~\ref{prop:castout}, $\dg{20}{0}{R}$ is the sum of the $21^2$ table entries modulo $20$. Each row is a
permutation of $\Ap{21}$, so the sum is $21\cdot210=21^2\cdot10\equiv10\pmod{20}$ for both $C_{21}$ and $C_7\rtimes_2C_3$, and
for every other Latin-square law on $\Ap{21}$. The digit does not distinguish the two groups.
\item A direct search gives $N=105, 189, 195, 255, 285, 297, 315, 345, 375$. The first decompositions are $N=105$: $52=3\cdot3^2+5^2$; $N=189$: $94=5\cdot3^2+7^2$; $N=195$: $97=8\cdot3^2+5^2$; $N=255$: $127=3\cdot3^2+4\cdot5^2$. By Theorem~\ref{thm:oddorder} none of these
orders carries a perfect group; for $N=105$ this also follows from the last exercise.
\item Write the element of label $m=m_0+7m_1$ of Section~\ref{sec:metacyclic-chars} as $(m_0,m_1)$, so that
$(a,b)(a',b')=(a+2^ba',\,b+b')$ and $(a,b)^{-1}=(-4^{b}a,\,-b)$, modulo $7$ and $3$. Put $\lambda(a,0)=a$ for $0\le a\le3$,
$\lambda(a,0)=14+a$ for $4\le a\le6$, $\lambda(a,1)=4+a$ and $\lambda(c,2)=17-\dg{7}{0}{-2c}$. Then $\lambda$ is a bijection onto
$\Ap{21}$ and $\lambda(g^{-1})=21-\lambda(g)$ for $g\ne1$, because $(a,1)^{-1}=(-4a,2)$ has label $17-a$. For $m=4$ and $n=1$ the
commutator $(21-m)*(21-n)*m*n=17*20*4*1$ equals $18$. Since $G/C_7\cong C_3$ is Abelian and $G$ is not, $G'=C_7$, whose labels
are $\{0,1,2,3,18,19,20\}$; thus $|G:G'|=3=1+2\cdot1$, matching the single pair of degree $1$ in $21=1+2(1^2+3^2)$.
\item The number of Sylow $7$-subgroups divides $15$ and is $\equiv1\pmod7$, so it is $1$ or $15$; the number of Sylow
$5$-subgroups divides $21$ and is $\equiv1\pmod5$, so it is $1$ or $21$ (Section~\ref{sec:dtsylow}). Both cannot be large,
since $15\cdot6+21\cdot4=174>105$. Hence a Sylow $5$-subgroup $P$ or a Sylow $7$-subgroup $Q$ is normal, and $H=PQ$ is a subgroup
of order $35$. Its index is $3$, and the action of $G$ on the three cosets of $H$ gives a homomorphism $G\to\Sym_3$ whose
kernel $K\le H$ has index dividing $\gcd(105,6)=3$. So $K=H$, and $H$ is normal.
\end{enumerate}

\section{Small finite groups}
\begin{enumerate}
\item $Z(G)\ne1$ for a $p$-group. If $Z(G)\ne G$ then $G/Z(G)$ has order $p$ and is cyclic, which forces $G$ to be Abelian.
\item The Sylow $p$-subgroup $\langle x\rangle$ is normal, and an involution $y$ acts by $x\mapsto x^{\pm1}$. This gives $C_{2p}$ or
$D_{2p}$.
\item The numbers of elements of order $1,2,4,8$ are $(1,1,2,4)$ for $C_8$, $(1,3,4,0)$ for $C_4\oplus C_2$, $(1,7,0,0)$ for $C_2^3$,
$(1,5,2,0)$ for $D_8$ and $(1,1,6,0)$ for $Q_8$. They are pairwise different.
\item A subgroup $H$ of order $6$ would have index $2$, so $g^2\in H$ for every $g$. Every element of order $3$ is a square, so $H$
would contain all $8$ of them, which is too many.
\item Either $n_3=4$, and the action on the Sylow $3$-subgroups gives $G\cong\Alt_4$; or the Sylow $3$-subgroup is normal and
$G=C_3\rtimes P$ with $P\in\{C_4,C_2^2\}$. The trivial and inverting actions give $C_{12}$, $Q_{12}$, $C_2\oplus C_6$ and $D_{12}$
(Theorem~\ref{thm:order12}).
\item \texttt{NumberSmallGroups} returns $14$, $15$ and $51$ for $16$, $24$ and $32$~\cite{GAP,BEO}.
\end{enumerate}

\section{Crystallographic groups}
\begin{enumerate}
\item In a basis of $L$ the matrix of $R$ has integer entries, so its trace is an integer. The trace is $2\cos(2\pi/n)$ in the plane
and $1+2\cos(2\pi/n)$ in space, so $2\cos(2\pi/n)\in\{-2,-1,0,1,2\}$, which gives $n=2,3,4,6,1$.
\item The point group preserves $L$, so $R\,mL=mL$, and $\{R\mid t\}\{I\mid v\}\{R\mid t\}^{-1}=\{I\mid Rv\}$ shows that $mL$ is
normal. Every element of $G$ is $\{R_i\mid t_i+24u\}$ for a unique $i$ and a unique $u\in\Z^3$, and two such elements are congruent
modulo $mL$ exactly when $i$ agrees and $u\equiv u'\pmod m$. Hence $|G/mL|=r\,m^3$, and~\eqref{eq:spacelabel} writes the pair
$(i,s)$, $s\in(\Z/m)^3$, in mixed radix $(r,m,m,m)$, which is a bijection onto $[0,r\,m^3)$.
\item $\{R_i\mid t_i\}\{R_j\mid t_j\}=\{R_l\mid R_it_j+t_i\}$ and $\{R_l\mid t_l\}$ lie in the same coset of the translation
subgroup, so $R_it_j+t_i-t_l\in24\Z^3$, i.e.\ $c(i,j)\in\Z^3$. Then
$\{R_i\mid t_i+24s\}\{R_j\mid t_j+24s'\}=\{R_l\mid t_l+24(R_is'+s+c(i,j))\}$, and reducing modulo $mL$ gives the label.
\item Here $t_0=t_1=0$, so every $c(i,j)=0$ and the law is $(i,s)(j,s')=(i+j\bmod2,\,(-1)^is'+s)$, the semidirect product with the
inversion acting by $-1$. For $m=2$, $-s=s$, so the action is trivial and every element has order at most $2$: $G/2L\cong C_2^4$.
\item In units of $\tfrac1{24}$, $t_1=t_3=(0,12,12)$ and $t_0=t_2=0$; $R_1=\operatorname{diag}(-1,1,-1)$,
$R_2=-I$, $R_3=\operatorname{diag}(1,-1,1)$, with $R_1^2=R_3^2=I$ and $R_1R_2=R_2R_1=R_3$. Hence
$c(1,1)=\tfrac1{24}(R_1t_1+t_1)=(0,1,0)$, $c(1,2)=\tfrac1{24}(t_1-t_3)=0$,
$c(2,1)=\tfrac1{24}(-t_1-t_3)=(0,-1,-1)$ and $c(3,3)=\tfrac1{24}(R_3t_3+t_3)=(0,0,1)$. So $\{R_1\mid t_1\}^2$ and
$\{R_3\mid t_3\}^2$ are the lattice translations by $(0,1,0)$ and $(0,0,1)$, which have order $m$ in $L/mL$; the screw and the glide
therefore have order $2m$. Since $t_2=0$, $\{R_2\mid t_2\}^2=1$.
\item An enantiomorphic pair is conjugate by an orientation-reversing affine map, hence isomorphic; by Bieberbach's theorem groups in
different affine classes are not isomorphic. The $230$ types therefore collapse to $230-11=219$ isomorphism classes.
\item A finite subgroup of $\mathrm{O}(2)$ is $C_n$ or $D_n$. Its rotations preserve the lattice, so $n\in\{1,2,3,4,6\}$ by the first
exercise, and each of the ten groups occurs, for instance for the square or the hexagonal lattice.
\end{enumerate}

\section{The symmetric and alternating groups}
\begin{enumerate}
\item The $11$ partitions of $6$ give class sizes $1,15,40,90,144,120$ (types $1^6,2,3,4,5,6$) and $45,15,40,120,90$ (types
$2^2,2^3,3^2,3\,2\,1,4\,2$), which add up to $720$.
\item Since $n=7$ is odd, a permutation is even if and only if its number of cycles is odd. This gives the types $7$; $5\,1^2$, $4\,2\,1$,
$3^2\,1$, $3\,2^2$; $3\,1^4$, $2^2\,1^3$; and $1^7$, eight in all. Only $(7)$ has distinct odd parts, so $k(\Alt_7)=9$.
\item The hook lengths give $f^{(5)}=f^{(1^5)}=1$, $f^{(4,1)}=f^{(2,1^3)}=4$, $f^{(3,2)}=f^{(2^2,1)}=5$ and $f^{(3,1^2)}=6$, and
$1+1+16+16+25+25+36=120$.
\item $\Alt_n$ is generated by products of two transpositions, and $(a\,b)(a\,c)=(a\,c\,b)$, $(a\,b)(c\,d)=(a\,b\,c)(b\,c\,d)$.
\item $|2|=|2^3|=15$, $|3|=|3^2|=40$ and $|6|=|3\,2\,1|=120$. An automorphism must map a class to a class of the same size and order.
\item The Lehmer code $(c_1,\dots,c_n)$ with $0\le c_i<i$ is a bijection onto permutations with $\sum_ic_i$ inversions, and its generating
function is the product of the factors $1+t+\dots+t^{i-1}$.
\end{enumerate}

\section{Finite simple groups of Lie type}
\begin{enumerate}
\item[(7)] The Weyl group of $B_r$ and $C_r$ is the full signed permutation group $C_2^r\rtimes\Sym_r$, of order $2^r r!$: a signed permutation may independently reverse every coordinate pair. In type $D_r$, preservation of the quadratic form imposes an even number of sign reversals, giving $C_2^{r-1}\rtimes\Sym_r$ of order $2^{r-1}r!$. These are exactly the pivot permutations compatible with the forms of Section~\ref{sec:classical}.
\end{enumerate}

\begin{enumerate}
\item $|\SL_2(p)|=p(p^2-1)$ and the centre $\{\pm1\}$ has order $2$. This gives $60$, $168$ and $660$ for $p=5,7,11$.
\item Label $a_0+a_1t$ by $a_0+2a_1$. Addition is bitwise. In multiplication, $2*2=3$ ($t^2=t+1$), $2*3=1$ ($t^2+t=1$) and $3*3=2$
($(t+1)^2=t$), so $t$ has order $3$ and generates $\mathbb F_4^\times$.
\item $|\SL_2(4)|=4\cdot15=60$, and $\PSL(2,4)=\SL_2(4)$ acts faithfully on the $5$ points of the projective line, giving $\Alt_5$.
$\PSL(2,5)$ has order $60$ and five classes, and it is simple; a simple group of order $60$ is isomorphic to $\Alt_5$.
\item $|\GL_3(2)|=(8-1)(8-2)(8-4)=168$, and $\SL_3(2)=\GL_3(2)$ since $\mathbb F_2^\times=1$.
\item The monomial matrices modulo the diagonal ones are the permutation matrices, so $W=N(T)/T\cong\Sym_n$. The Bruhat cells are
indexed by $W$, so there are $n!$ of them.
\item $8^2\cdot65\cdot7=29120$ and $2^6\cdot63\cdot3=12096$.
\end{enumerate}

\section{Sporadic groups}
\begin{enumerate}
\item[(5)] A deep hole of the Leech lattice has nearest points forming an affine ADE diagram; this records the root lattice of the associated Niemeier lattice, not a root system for the rootless Leech lattice itself. Independently, $G_2(4)$ is the point stabiliser in the displayed primitive action of $\mathrm{Suz}$, so its type-$G_2$ diagram describes a local subgroup. In neither case does a diagram determine the sporadic group: the global group is obtained from lattice automorphisms or a stabiliser-chain action rather than from a Chevalley construction.
\end{enumerate}

\begin{enumerate}
\item $M_{11}$ is sharply $4$-transitive and $M_{12}$ sharply $5$-transitive. With base $(0,1,2,3)$, respectively $(0,\dots,4)$, the
basic orbits have lengths $11,10,9,8$, respectively $12,11,10,9,8$, and the order is their product.
\item $24\cdot23\cdot22\cdot21\cdot20=5100480$ and $5100480\cdot48=244823040$.
\item $m=t_1+11t_2+110t_3+990t_4\mapsto s_{1,t_1}s_{2,t_2}s_{3,t_3}s_{4,t_4}$. Since $s_{l,0}=1$, the identity has label $0$.
\item Standard generators are defined by conditions on classes and words, which any implementation can test. The canonical rule then
uses only the group and these generators, so its output does not depend on the software.
\item \texttt{NrConjugacyClasses(MathieuGroup(11))} returns $10$, matching the recomputed classes.
\item $|\mathbb M|=808017424794512875886459904961710757005754368\cdot10^9\approx8.08\cdot10^{53}$. So labels in $[0,|\mathbb M|)$ have up
to $54$ decimal digits.
\end{enumerate}

\section{Worked digital products}
\begin{enumerate}
\item Since $g^3=1$, $g^{-1}=g^2$, so $145\cdot145$ is the label $137$ of $g^{-1}$ in the table. An element of order $2$ is its own
inverse, which is why $(gh)^{-1}=101=gh$.
\item The order of $H=\langle g,h\rangle$ is divisible by $3$ and $7$, hence by $21$, and divides $168$, so $|H|\in\{21,42,84,168\}$.
$\mathrm{PSL}_2(7)$ is simple of order $168$, so it has no subgroup of index $2$ or $4$ (the action on the cosets would give a
nontrivial homomorphism to $\Sym_2$ or $\Sym_4$, of order less than $168$). A group of order $21$ has no element of order $2$, but
$gh\in H$ has order $2$. Hence $H=\mathrm{PSL}_2(7)$.
\item $|\langle g,h\rangle|$ is divisible by $4\cdot13\cdot7=364$. A proper subgroup lies in a maximal one, and none of $448$, $52$,
$20$, $14$ is divisible by $364$. Hence $\langle g,h\rangle=\mathrm{Sz}(8)$.
\item $(g^{-1})^n=(g^n)^{-1}$ is $1$ exactly when $g^n=1$; $(gh)(h^{-1}g^{-1})=1$; and $hg=g^{-1}(gh)g$ is conjugate to $gh$. The last
column, the orders, is the same for $g,h,gh$ and for $g^{-1},h^{-1},(gh)^{-1}$, so the orders of the inverses can be read off
without computing them.
\item If $m^{k+1}=m$ then $m^k=1$, and minimality of $k$ makes it the order; then $m\cdot m^{k-1}=1$. For $gh$ the powers up to
$(gh)^{36}=gh$ take $35$ products. Repeated squaring computes $(gh)^{34}$ from $34=100010_2$ with $5$ squarings and $1$ further
product, $6$ in all. Finally $(gh)^5$ has order $35/\gcd(35,5)=7$ and $(gh)^7$ has order $5$.
\end{enumerate}

\end{document}